\documentclass{article}
\usepackage[T1]{fontenc}
\usepackage[latin1]{inputenc}
\usepackage{amsthm, fullpage}
\usepackage{amsmath,enumerate}
\usepackage{amssymb}
\usepackage{amsfonts}
\usepackage{calrsfs}
\usepackage{eucal}
\usepackage[all]{xy}
\usepackage{hyperref}
\providecommand{\scr}{\mathcal}
\usepackage{mathrsfs}
\DeclareMathOperator{\ad}{ad}

\DeclareMathOperator{\Mod}{Mod}
\DeclareMathOperator{\coh}{Coh}
\DeclareMathOperator{\Ker}{Ker}

\DeclareMathOperator{\Frac}{Frac}

\DeclareMathOperator{\Spa}{Spa}
\DeclareMathOperator{\Spv}{Spv}
\DeclareMathOperator{\ind}{Ind}
\DeclareMathOperator{\indcoh}{IndCoh}

\DeclareMathOperator{\supp}{supp}

\newtheorem{prop}{Proposition}[subsection]
\newtheorem{theo}[prop]{Theorem}
\newtheorem{coro}[prop]{Corollary}
\newtheorem{lemm}[prop]{Lemma}
\newtheorem{lem}[prop]{Lemma}
\newtheorem*{lemm*}{Lemma}

\newtheorem{exs}[prop]{Examples} 
\newtheorem{ex}[prop]{Example}

\theoremstyle{definition}

\newtheorem{empt}[prop]{}
\newtheorem{dfn}[prop]{Definition}
\newtheorem{rem}[prop]{Remark}
\newtheorem{ntn}[prop]{Notation}
\newtheorem*{rem*}{Remark}

\theoremstyle{thm}
\newtheorem{thm}[prop]{Theorem}
\newtheorem*{thm*}{Theorem}
\newtheorem*{lem*}{Lemma}
\newtheorem{cor}[prop]{Corollary}
\newtheorem*{cor*}{Corollary}
\newtheorem*{prop*}{Proposition}

\theoremstyle{dfn}
\newtheorem*{dfn*}{Definition}

\numberwithin{equation}{prop}

\newcommand{\riso}{ \overset{\sim}{\longrightarrow}\, }
\newcommand{\liso}{ \overset{\sim}{\longleftarrow}\, }

\newcommand{\Spec}{\mathrm{Spec}\,}
\newcommand{\Spf}{\mathrm{Spf}\,}

\renewcommand{\sp}{\mathrm{sp}}
\renewcommand{\det}{\mathrm{det}}
\newcommand{\gr}{\mathrm{gr}}

\renewcommand{\AA}{{\mathcal{A}}}

\newcommand{\FF}{{\mathcal{F}}}
\newcommand{\B}{{\mathcal{B}}}

\newcommand{\E}{{\mathcal{E}}}
\newcommand{\G}{{\mathcal{G}}}
\renewcommand{\H}{{\mathcal{H}}}
\newcommand{\M}{{\mathcal{M}}}
\newcommand{\D}{{\mathcal{D}}}
\newcommand{\I}{{\mathcal{I}}}

\newcommand{\K}{{\mathcal{K}}}

\newcommand{\bbA}{{\mathbb{A}}}

\newcommand{\bbD}{{\mathbb{D}}}

\newcommand{\bbL}{{\mathbb{L}}}
\newcommand{\bbN}{{\mathbb{N}}}
\newcommand{\bbP}{{\mathbb{P}}}
\newcommand{\bbQ}{{\mathbb{Q}}}
\newcommand{\bbR}{{\mathbb{R}}}
\newcommand{\bbZ}{{\mathbb{Z}}}

\newcommand{\A}{\mathbb{A}}

\newcommand{\DD}{\mathbb{D}}

\renewcommand{\L}{\mathbb{L}}
\renewcommand{\P}{\mathbb{P}}
\newcommand{\R}{\mathbb{R}}
\newcommand{\Q}{\mathbb{Q}}
\newcommand{\Z}{\mathbb{Z}}
\newcommand{\N}{\mathbb{N}}

\newcommand{\cB}{{\mathcal{B}}}
\newcommand{\cC}{{\mathcal{C}}}
\newcommand{\cD}{{\mathcal{D}}}
\newcommand{\cE}{{\mathcal{E}}}
\newcommand{\cF}{{\mathcal{F}}}
\newcommand{\cG}{{\mathcal{G}}}
\newcommand{\cH}{{\mathcal{H}}}
\newcommand{\cI}{{\mathcal{I}}}

\newcommand{\cK}{{\mathcal{K}}}

\newcommand{\cM}{{\mathcal{M}}}
\newcommand{\cN}{{\mathcal{N}}}
\newcommand{\cO}{{\mathcal{O}}}
\newcommand{\cP}{{\mathcal{P}}}

\newcommand{\cV}{{\mathcal{V}}}
\newcommand{\cW}{{\mathcal{W}}}

\newcommand{\fB}{\mathfrak{B}}
\newcommand{\fC}{{\mathfrak{C}}}
\newcommand{\fD}{{\mathfrak{D}}}

\newcommand{\fI}{{\mathfrak{I}}}
\newcommand{\fM}{{\mathfrak{M}}}
\newcommand{\fP}{{\mathfrak{P}}}
\newcommand{\fQ}{{\mathfrak{Q}}}
\newcommand{\fR}{{\mathfrak{R}}}
\newcommand{\fS}{{\mathfrak{S}}}
\newcommand{\fT}{{\mathfrak{T}}}
\newcommand{\fU}{{\mathfrak{U}}}

\newcommand{\fX}{{\mathfrak{X}}}
\newcommand{\fY}{{\mathfrak{Y}}}
\newcommand{\fZ}{{\mathfrak{Z}}}

\newcommand{\fa}{\mathfrak{a}}
\newcommand{\fb}{\mathfrak{b}}

\newcommand{\fm}{\mathfrak{m}}
\newcommand{\fp}{\mathfrak{p}}
\newcommand{\fq}{\mathfrak{q}}
\newcommand{\fr}{\mathfrak{r}}

\renewcommand{\O}{{\mathcal{O}}}

\newcommand{\V}{\mathcal{V}}
\newcommand{\W}{\mathcal{W}}
\newcommand{\T}{{\mathfrak{T}}}
\newcommand{\Y}{\mathfrak{Y}}
\newcommand{\ZZ}{\mathfrak{Z}}
\newcommand{\X}{\mathfrak{X}}

\newcommand{\U}{\mathfrak{U}}
\newcommand{\PP}{{\mathfrak{P}}}

\newcommand{\hdag}{  \phantom{}{^{\dag} }    }

\def\debrom{
\makeatletter
\renewcommand{\theenumi}{(\roman{enumi})}
\renewcommand{\labelenumi}{\theenumi}
\makeatother\begin{enumerate}}
\def\finrom{\end{enumerate}}

\begin{document}

\title{Arithmetic $\D$-modules over Laurent series fields: absolute case}
\author{Daniel Caro}
\date{March 18, 2021}

\maketitle

\tableofcontents

\bigskip

\begin{abstract}
Let $k$ be a perfect field of characteristic $p>0$.
Within Berthelot's theory of arithmetic $\D$-modules,
we construct a $p$-adic formalism of Grothendieck's six operations for quasi-projective schemes over $\Spec k[[t]]$. 
\end{abstract}

\section*{Introduction}
Let $\V$ be a complete discrete valuation ring of mixed characteristic $(0,p)$, 
$\pi$ be a uniformizer,
$k:= \V/\pi \V$ be its residue field and $K$ be its fraction field. 
In order to build a $p$-adic formalism of Grothendieck six operations for $k$-varieties (i.e. separated $k$-schemes of finite type), 
Berthelot introduced an arithmetic avatar of the theory of modules over the differential operators ring. 
The objects appearing in his theory are called arithmetic $\D$-modules or complexes of arithmetic $\D$-modules
(for an introduction,
see \cite{Beintro2}). 

Within Berthelot's arithmetic $\D$-modules theory,
such a $p$-adic formalism was already known in different contexts. Let us describe these known cases. 
With N. Tsuzuki  (see \cite{caro-Tsuzuki}),
the author got such a formalism  for overholonomic $F$-complexes of arithmetic $\D$-modules 
(i.e. complexes together with a Frobenius structure)
over realizable $k$-varieties (i.e. 
$k$-varieties which can be embedded into a proper formal $\V$-scheme). 
Another example was given later (do not focus on the publication date) with 
holonomic $F$-complexes of arithmetic $\D$-modules over quasi-projective varieties (\cite{caro-stab-holo}).
In a wider geometrical context, T. Abe 
established a six functors formalism for admissible stacks, 
namely algebraic stacks of finite type with finite diagonal morphism
(see \cite[2.3]{Abe-Langlands}). The starting point of his work was the case of quasi-projective $k$-varieties. 
Again, some Frobenius structures are involved in his construction. 
Finally, without Frobenius structure, in \cite{caro-unip}, 
we explained how to build 
such a $p$-adic formalism of Grothendieck's six functors, e.g. with quasi-unipotent complexes of arithmetic $\D$-modules (see \cite{caro-unip}).

Recently, Lazda and P\'al have developped in their book \cite{Lazda-Pal-Book} a theory of overconvergent isocrystals on schemes of finite type over 
$\Spec ~ k[[t]]$. Their constructions are inspired by Berthelot's one. 
One main difference is the use of adic spaces instead of Tate's analytic rigid spaces. 
Similarly to Berthelot's category of overconvergent isocrystals, these  overconvergent isocrystals are stable under
tensor products, pull-backs, duality, extraordinary pull-backs, rigid cohomology. 
But, we do not have the stability under push-forwards by a closed immersion. Hence, two of Grothendieck's six operations (the push-forward
and the extraordinary push-forward) are missing. 
In order to obtain succefully a $p$-adic formalism of Grothendieck's six operations valid for schemes of finite type over $\Spec k[[t]]$, 
the purpose of this paper is to extend their work in the theory of arithmetic $\cD$-modules. 
Here, we focus on an ``absolute'' theory, i.e. if $X$ is a $\Spec k[[t]]$-scheme of finite type,
then we introduce a theory of arithmetic $\cD$-modules on $X /\Spec k$ and not on $X /\Spec k [[t]]$, which is the ``relative'' case. 
In order to shorten and simplify the presentation for the author and the reader, we have decided to treat later in another work the relative case. 

\bigskip

Let us clarify the content of the paper. 
Let $\fS := \Spf (\cV)$, $i \geq 0$ be an integer, $S _i := \Spec~ (\cV / \pi  ^{i+1} \cV)$.
For any integer $r\geq 0$, we set 
$\bbD ^r _{\fS}:= 
\Spf \cV  [[ t _1, \dots, t _r]]$
and 
$\bbD ^r _{S _{i}}
:=
\Spec (\cV / \pi  ^{i+1} \cV) [[ t _1, \dots, t _r]]$.
In the first chapter, we give some details and complements on the notion of relative perfect morphisms. 
Moreover, in order to study arithmetic $\cD$-modules in a nice wider context than that of smooth formal $\fS$-schemes or 
smooth $S _{i}$-schemes, we introduce the notion of morphism (locally) of formal finite type.
More precisely, let $Y$ be an $S _i$-scheme of finite type.
An $S_i$-morphism $f\colon X \to Y$ is ``of formal finite type'' is there exists an integer $r\geq 0$ such that 
$X$ is a $Y \times _{S _i}\bbD ^r _{S _{i}}$-scheme of finite type.
When this property is satisfied Zariski locally in $X$, we say that $f$ is ``locally of formal finite type''.
If $f$ is a formally smooth morphism locally of formal finite type then we can check that $f$ is flat
(see \ref{f0formétale-fforméta0}).
Moreover, the usual properties satisfied by étale morphisms extend
to the case of formally étale morphisms of formal finite type (see \ref{EGAIV18.1.2}).
If $\fY$ is a formal $\fS$-scheme of finite type, 
then we get similarly a notion of 
morphism $\fX\to \fY$ locally of formal finite type.

In the second chapter, 
we give an overview with some complements
of the notion of valued triples, analytic triples,
adic spaces and Zariski-Riemann spaces. This is the geometric context of 
Lazda and P\'al in their book \cite{Lazda-Pal-Book}. We will not give some comparison theorem between 
Lazda and P\'al's overconvergent isocrystals and some arithmetic $\cD$-modules. 
However, 
in order to define 
the local cohomology with support in a closed subscheme having locally finite $p$-bases of the constant coefficient
(see definition \ref{ntn-GammaZO-rig} and the remark \ref{rem-why-2dfn-loc-coh}), 
we will need the local cohomology in the context of adic spaces as defined by 
Lazda and P\'al in \cite{Lazda-Pal-Book} (see below the definition \ref{loc-coh&cech-reso}.\ref{loc-coh&cech-reso1}).
To be able to define the local cohomology in a wider context (in the chapter \ref{chapter12}),
we will need the coherence of the constant coefficient whose proof uses the very special case of 
the local cohomology with support in a closed subscheme having locally finite $p$-bases of the constant coefficient.
This is the main reason why we had to use  adic spaces. 

Let  $\fP$ be a separated formal $\fS$-scheme (for the $p$-adic topology) locally of formal finite type and having locally $p$-bases over $\fS$.
The special fiber of $\fP$, the $k$-scheme equal to its reduction modulo $\pi$, is denoted by $P$.
In the third chapter, we give the Berthelot's construction of $\cD ^{\dag} _{\fP/\fS}$,
the sheaf  of differential operators on $\fP/\fS$. The sheaf $\cD ^{\dag} _{\fP/\fS}$ is some kind of weak completion as $\cO _{\fP}$-ring of 
Grothendieck sheaf $\cD _{\fP/\fS}$ of differential operators of finite order. 
We recall that some properties of $\cD ^{\dag} _{\fP/\fS}$ was studied in a wider context by the author and Vauclair in \cite{Caro-Vauclair}.
To sum up, $\cD ^{\dag} _{\fP/\fS}$ behaves as nicely as in the case where $\fP/\fS$ is smooth.
In this paper, we only focus on schemes locally of finite type over the base. 
As explained in the first chapter, in this context relative perfectness behaves better
and satisfies similar to etaleness properties.
We also give the construction of the extraordinary pull-back and the push-forward by a morphism $f$ of schemes (hence the level is fixed) denoted respectively by
$f ^!$ and $f _+$ (or $f ^{!(m)}$ and $f ^{(m)}_+$ if we prefer to indicate the fixed level $m$).
Since these functors depend on the base, we study their behavior with respect to the change of the base (e.g. see \ref{prop-u!=u!/T(m)},
\ref{isooub-u+(m)=u+(m)/T}). 
Finally, we check some  Spencer resolutions and the projection formula.

In the forth chapter, we show that Berthelot's proof of Frobenius descent of the level
is still valid in the context of schemes having locally $p$-bases over the base and locally of formal finite type.
We have tried to be a bit complementary : 
we give sometimes some details not appearing in Berthelot's original proof and 
it is omited to write the proof when this is completely the same word by word.
In the case of the sheaf of differential operators of finite level, the Frobenius descent can 
simply by stated by saying the direct image by Frobenius and the inverse image by Frobenius induce quasi-inverse equivalences of categories.

In the fifth chapter, 
we recall Berthelot's notion of derived categories of inductive systems of arithmetic $\D$-modules on $\fP/\fS$.
Some objects in theses categories will give our coefficients satisfying a six functors formalism. 
Two Berthelot's (too technical to sum-up) notions are fundamental in theses categories :  that of ``quasi-coherence'' and that of ``coherence''. 

In the sixth chapter, 
we study the localization functor outside a divisor $T$ of $P$ and the forgetful functor of a divisor $T$ of $P$. 
We check both functors preserve the quasi-coherence.
Next, we give a coherence stability criterion involving a change of divisors
which is one fundamental property of the theory (see \ref{limTouD}).

In the seventh chapter, 
in the context of inductive systems of arithmetic $\D$-modules on $\fP/\fS$, 
we give the construction of extraordinary inverse images, direct images,
duality and base change.
We check the commutation of both functors with respect to the localization functor outside a divisor.
We prove the commutation of the base change and its commutation with 
tensor products, extraordinary pullbacks, direct images, duality.

Let $f$ be a closed immersion of formal 
$\fS$-schemes of formal finite type having locally $p$-bases.
In the eighth chapter, by proving the fundamental local isomorphism, 
we check that the relative duality isomorphism holds for $f$ and that we have the adjoint paire 
$(f _+, f ^!)$. 
The adjunction morphisms in this case are very explicit. 
Finally, if $X$ is a closed subscheme of $P$ having locally $p$-bases (over $S$), 
we construct by glueing the category of 
coherent arithmetic $\D$-modules over $X/\fS$.
More precisely, first we choose $(X _\alpha)$ an affine open covering of $X$, 
and for each $\alpha$ we choose a formal $\fS$-scheme $\X _\alpha$ having locally $p$-bases which is a lifting of $X _\alpha$.
Then, a coherent arithmetic $\D$-module 
over $X/\fS$ is the data of a family of coherent arithmetic $\D$-module on $\fX _\alpha$
together with glueing isomorphisms
satisfying a cocycle condition.
We check that we have a canonical equivalence of categories between that
of coherent arithmetic $\D$-modules 
over $X/\fS$ and that of coherent arithmetic $\D$-modules on $\fP$
with support in $X$
(see Theorem \ref{prop1}).
This extends Berthelot's theorem of his arithmetic version of Kashiwara theorem
appearing in the classical $\D$-modules theory.

In the ninth chapter, we introduce the notion of convergent isocrystals
in the framework of arithmetic $\cD$-modules.
More precisely, we denote by $\mathrm{MIC} ^{\dag \dag} (\fP/K)$,
the full subcategory of
left $\D ^\dag _{\fP, \Q} $-modules
consisting of 
left $\D ^\dag _{\fP, \Q} $-modules
which are 
$\O _{\fP,\Q} $-coherent.
The category $\mathrm{MIC} ^{\dag \dag} (\fP/K)$
can be seen as an analogue of the category of  convergent isocrystals on $\fP/\fS$.
Moreover, let $u \colon X  \hookrightarrow P $ be a purely of codimension $r$ closed immersion of schemes having locally finite $p$-bases over $\Spec k$.
We set 
$\R \underline{\Gamma} ^\dag _X \O _{\fP,\Q} 
:=\R \sp _* \underline{\Gamma} ^\dag _X ( \O _{\fP _K})$,
where $\sp \colon \fP _K \to \fP$
is  the specialization morphism from $\fP _K$, the adic space associated to $\fP$ (see \ref{ntn-GammaZO-rig})
and where 
$\underline{\Gamma} ^\dag _X$ is Lazda and P\'al's local cohomology.
The complex
$\R \underline{\Gamma} ^\dag _X \O _{\fP,\Q} $ 
is by definition the ``local cohomology with support in $X$
of the constant coefficient $\cO _{\fP,\bbQ}$''.
From the canonical morphism
$\underline{\Gamma} ^\dag _X (\O _{\fP _K}) 
\to 
\O _{\fP _K}$,
we get the morphism
$\R \underline{\Gamma} ^\dag _X \O _{\fP,\Q} 
\to 
\O _{\fP,\Q}$ (this map is a key tool in the proof of the coherence of the contant coefficient described below). 
Suppose there exists 
a finite $p$-basis 
$t _1, \dots, t _d $ of $\fP/\fS$.
Let $T$ be the divisor of $P$ defined by setting
$T := V ( \overline{t} _1\cdots \overline{t} _r)$ with $r \leq d$,
where
$\overline{t} _1, \dots, \overline{t} _r$ are the images of $t _1,\dots, t _r$ in $\Gamma (P ,\O _{P})$.
Then Berthelot's global presentation of 
$\cO _{\fP} (\hdag T) _{\bbQ}$, 
the constant coefficient of $\fP$ with overconvergent singularities along $T$,
is still valid and in particular
$\cO _{\fP} (\hdag T) _{\bbQ}$
 is a coherent 
$\D ^\dag _{\fP, \Q}$-module (see \ref{NCDgencoh}).
This situation arises  for instance  when $(P,T)$ is a strict semi-stable paire. 
This coherence theorem will be extended for any divisor $T$ but when $\fP$ is of finite type over $\bbD ^1 _\fS$ in the eleventh chapter.

In the tenth chapter, we study external tensor products. 
In order to be able to define external tensor products in our context, 
we need to have geometrical objects stable under products. 
To do so, 
we define the category $\scr{C} _{\fS}$  whose objects are
finite type morphisms of formal $\fS$-schemes of the form 
$\fP \to \bbD ^r _{\fS}$ for some integer $r$.
A morphism $f \to g$ of $\scr{C} _{\fS}$, where 
$f\colon \fP \to \bbD ^r _{\fS}$ and
$g\colon \fQ \to \bbD ^s _{\fS}$, consists in 
two morphisms 
$\alpha \colon 
\fP \to \fQ$ and 
$\beta \colon \bbD ^r _{\fS} \to \bbD ^s _{\fS}$
of formal $\fS$-schemes
making commutative the diagram
\begin{equation}
\notag
\xymatrix{
{\fP} 
\ar[r] ^-{\alpha}
\ar[d] ^-{f}
& 
{\fQ} 
\ar[d] ^-{g}
\\ 
{\bbD ^r _{\fS}} 
\ar[r] ^-{\beta}
& 
{\bbD ^s _{\fS}.} 
}
\end{equation}
Let $f\colon \fP \to \bbD ^r _{\fS}$ and
$g\colon \fQ \to \bbD ^s _{\fS}$ be two objects of 
$\scr{C} _{\fS}$.
We put 
\begin{equation}
\notag
\fP \times _{\scr{C} _{\fS}} \fQ : = 
\left ( \fP \times _{\bbD ^r _{\fS}} \bbD ^{r+s} _{\fS} \right ) 
\times _{\bbD ^{r+s} _{\fS}}
\left ( \bbD ^{r+s} _{\fS} \times _{\bbD ^s _{\fS}} \fQ \right )
\riso 
 \fP \times _{\bbD ^r _{\fS}} 
 \bbD ^{r+s} _{\fS} \times _{\bbD ^s _{\fS}} \fQ.
\end{equation}
We denote by 
$f \times  _{\scr{C} _{\fS}} g 
\colon \fP \times _{\scr{C} _{\fS}} \fQ
\to 
 \bbD ^{r+s} _{\fS} $, 
$pr _1 \colon \fP  \times _{\scr{C} _{\fS}} \fQ \to \fP$
and by
$pr_2 \colon \fP \times _{\scr{C} _{\fS}} \fQ \to \fQ$
the canonical projections, 
by 
$pr _{1} ^{r ,s}\colon 
\bbD ^{r+s} _{\fS}
\to 
\bbD ^r _{\fS}$
and
$pr _{2} ^{r ,s}\colon 
\bbD ^{r+s} _{\fS}
\to \bbD ^s _{\fS}$ the canonical morphisms.
Such morphisms of the form 
$pr_1$ or $pr_2$ are called ``projection morphism''. 
We check that $f \times  _{\scr{C} _{\fS}}g $ equipped with the morphisms
$(pr _1 , pr _{1} ^{r ,s}) $
and 
$(pr _2 , pr _{2} ^{r  ,s}) $
of $\scr{C} _{\fS}$
satisfies the universal property of the product
in $\scr{C} _{\fS}$ (see \ref{dfn-CfS2}).
We define in this context of external tensor products 
and we check they satisfy the expected properties such as 
the commutation of external tensor products with pull-backs or push-forwards 
(see \ref{prop-boxtimes}), the commutation of push-forwards with pull-backs by projection morphisms (see \ref{theo-iso-chgtbase2}).
Since this latter property was also checked for pull-backs by a closed immersion, this yields 
the commutation of push-forwards with pull-backs by projective morphisms, which can be called the ``base change isomorphism'' for projective morphisms.
This implies the relative duality isomorphism for projective morphisms and that we have the adjoint paire 
$(f _+, f ^!)$
(see \ref{rel-dual-isom-proj-formal}).
We also extend the relative duality isomorphism 
and the adjunction $(f _+, f ^!)$ for relatively proper complexes and quasi-projective morphisms,
i.e. we replace the properness hypothesis of the morphism $f$
by the properness via $f$ of the support of our complexes
(see \ref{rel-dual-isom}).

Let $f\colon \Y  \to \X  $ be 
a morphism of 
formal  $\fS$-schemes
of formal finite type
and having locally finite $p$-bases  over $\fS$.
We suppose  that the induced morphism 
$f _0 \colon Y \to X$ is a finite, surjective and radicial morphism. 
Then we prove in the eleventh chapter that the functor 
$f _+$ 
from the category of 
coherent left  $\D ^{\dag} _{\Y /\fS ,\Q}$-modules 
to that of 
coherent left  $\D ^{\dag} _{\X /\fS ,\Q}$-modules
is an exact  quasi-inverse equivalence of categories of $f ^* $ (see \ref{256Be2}).
The functors 
$f _+$ and $f ^!$ induce moreover quasi-inverse equivalences of categories between 
$\mathrm{MIC} ^{\dag \dag} (\fX/K)$
and 
$\mathrm{MIC} ^{\dag \dag} (\Y/K)$ 
(see \ref{univhomeo-eqcat-isoc}).
Let $\X $ be a
formal  $\Spf \, \V [[t]]$-scheme of finite type
and having locally finite $p$-bases  over $\fS$.
Let $Z$ be a divisor of $X$.
Adapting Berthelot's original proof, we check that 
$\O _{\X} (\hdag Z) _\Q$,
the constant coefficient on $\fX$
with overconvergent singularities along $Z$,
 is a coherent 
$\D ^\dag _{\X, \Q}$-module (see \ref{coh-cst-div}).
The key point is to use de Jong's desingularisation theorem (see \cite[6.5]{dejong})
which implies that there exist a trait 
$\mathbb{D} ^1 _{S'}=\Spec k '[[u]]$ (with $S ': =\Spec k'$) 
finite over 
$\mathbb{D} ^1 _{S}=\Spec k [[t]]$ such that 
$k [[t]] \to k' [[u]] $ is a morphism of traits, 
a separated $\mathbb{D} ^1 _{S'}$-scheme of finite type $X '$, an alteration of schemes over $\mathbb{D} ^1 _{S}$ (in the sense of \cite[2.20]{dejong})
$g _0 \colon X ' \to X$ 
and an open immersion $j ' \colon X ' \to \overline{X} '$ of $\mathbb{D} ^1 _{S'}$-schemes, with the following properties: 
\begin{enumerate}[(a)]
\item $\overline{X} '$ is an integral projective $\mathbb{D} ^1 _{S'}$-scheme with geometrically irreducible generic fibre, and 
\item  the pair $(\overline{X} ', g _0 ^{-1} (Z ) _\mathrm{red} \cup (\overline{X} ' \setminus j ' (X '))$ is strict semi-stable.
\end{enumerate}
In particular, we get that 
$(X ', g _0  ^{-1} (Z ) _\mathrm{red})$ is a strict semi-stable pair
and there exists a closed immersion of the form
$u _0\colon X ' \hookrightarrow \P ^n _{X}$
whose composition with the projection 
$\P ^n _{X} \to X$ is $g _0 $.
Hence, by universal homeomorphic descent, by using 
the local cohomology with support in a closed subscheme having locally finite $p$-bases 
of the constant coefficient (constructed in chapter nineth), we reduce to the case where 
$(X,Z)$ is a strict semi-stable paire, which was already proved in the nineth chapter.

Let $\fP$ be a formal $\fS$-scheme of formal finite type 
and having locally $p$-bases.
If  $T$ be a divisor of $P$, then 
we do not know if 
$\O _{\fP} (\hdag T) _{\Q}$ is a $\cD ^\dag _{\fP, \Q}$-coherent (because when 
$\fP$ is not a $\bbD ^1 _{\fS}$-scheme of finite type,
de Jong's desingularization theorem can not apply directly).
Hence, in the twelfth  chapter, 
we introduce the notion of ``weak admissible'' subschemes 
$Y$ of $P$. 
Roughly speaking (more precisely, see \ref{dfn-4.3.4bisY}), 
the inclusion $(Y \subset \fP)$ is weak admissible means that 
the constant coefficient on $Y ' $ in $\fP$
is $\cD ^\dag _{\fP, \Q}$-coherent for any subscheme $Y'$ of $Y$ (and this property has to be stable under
pullbacks by smooth projections). The word ``weak'' is added because the stability under duality is not clear (the notion
of admissibility will be define in the next chapter).
We introduce the notion of 
local cohomological functor with strict support over a weak admissible subscheme $Y$ of $P$ that we denote by 
$\R \underline{\Gamma} ^\dag _Y$. 
Next,
we check that expected properties satisfied by local cohomological functors are still valid, e.g. its commutation 
with pushforwards and extraordinary pullbacks. We also check some base change isomorphism (see \ref{theo-iso-chgtbase}).

In the thirteenth chapter, 
we adapt the construction given in  \cite{caro-unip}
of a formalism of Grothendieck six functors. 
We first introduce the notion of ``restricted'' data of absolute coefficients over $\fS$, 
i.e. we attach some coefficients of arithmetic $\cD$-modules 
to formal $\bbD ^1 _{\fT}$-schemes of finite type, where $\fT := \Spf ~ \cW$ 
with $\W$ a complete discrete valued $\cV$-algebra  of mixed characteristic $(0,p)$
with perfect residue field. 
Via Theorem \ref{dfnquprop} and the example \ref{ex-datastableevery},
we explain how to build a restricted data of absolute coefficients 
which contains the constant coefficient, 
which is 
local,
stable by devissages, direct summands, 
local cohomological functors, 
pushforwards, quasi-projective extraordinary pullbacks, base change, tensor products, duality.

The problem of the restricted version
is that we can not define external tensor products since 
formal $\bbD ^1 _{\fS}$-schemes of finite type are not stable under products.
In order to get some data stable under external tensor products, we introduce 
the notion of data of absolute coefficients over objects of $\cC _{\fS}$.
We give some receipt to construct some data of absolute coefficients which
contain the constant coefficient (without overconvergent singularities), 
satisfy $BK ^!$,
are
local,
are stable by devissages, direct summands, 
admissible local cohomological functors, 
pushforwards, extraordinary pullbacks by projections, 
base change, 
admissible external tensor products,
admissible duality (see \ref{theo-V-6operations}).
The notions of stability under admissible local cohomological functors, 
admissible external tensor products,
admissible duality are defined respectively similarly to 
the stability under  local cohomological functors, 
external tensor products,
duality except that roughly speaking the property is checked 
after restricting on admissible subschemes.

Finally, in the last chapter, 
we introduce the notion of ``frames over $\cV[[t]]$''
(see \ref{dfn-framW[[t]]}) as follows :
the objects are 
of the form $(Y,X,\fP)$  
where 
$\fP$ is a quasi-projective smooth  formal $\bbD ^1 _{\fS}$-scheme, 
$X$ is a reduced closed subscheme of the special fiber $P$ of $\fP$ 
and $Y$ is an open subscheme of $X$. 
We define the category of ``couples over $\cV[[t]]$'' 
whose objects are of the form
 $(Y,X)$, where $X$ is a quasi-projective  $\bbD ^1 _{S}$-scheme
 and $Y$ is an open subscheme of $X$.
Such couples 
can be enclosed into a frame over $\cV[[t]]$ of the form
$(Y,X,\fP)$. 
For an enough stable data of absolute coefficients $\fC$, a coefficient of $\fC$ 
over $(Y,X,\fP)$ is a coefficient of $\fC$ over $\fP$ with support in $X$ and having overconvergent singularities along $X \setminus Y$
(i.e. which is isomorphic under its image via $\R \underline{\Gamma} ^\dag _{X \setminus Y}$).
We prove the independence (for such data $\fC$) with respect to the choice of the frame $(Y,X,\fP)$ over $\cV[[t]]$
enclosing $(Y,X)$ 
(\ref{ind-CYW}), and we denote them by $\fC (Y,X/\cV [[t]])$. 
This yields a 
$p$-adic formalism of Grothendieck six operations over couples $(Y,X)/\cV [[t]]$.
Finally, when $X$ is projective over $\bbD ^1 _S$, then 
the category 
$\fC (Y,X/\cV [[t]])$
is  independent (up to canonical equivalence of categories) 
of the choice of such projective $\bbD ^1 _S$-scheme $X$ enclosing $Y$. 
Hence, we denote the corresponding category by 
$\fC (Y/\cV [[t]])$.
This yields a formalism of Grothendieck's six operations over quasi-projective $\bbD ^1 _{S}$-schemes.
 
\subsection*{Acknowledgment}
The author was supported by the IUF.

\section*{Notation}

Let $\V$ be a complete discrete valuation ring of mixed characteristic $(0,p)$, 
$\pi$ a uniformiser, $e$ the ramification index, 
$k$ be its residue field and $K$ its  field of fractions. 
We put $S := \Spec k$, $\mathfrak{S}: = \Spf \V$.  
A formal $\fS$-scheme $\X$ means is a noetherian 
$p$-adic formal scheme endowed with 
a structural morphism of $p$-adic formal schemes
$\X \to \Spf \V$.
We will work mostly with noetherian formal schemes
but we do not follow Grothendieck's terminology of EGA,
i.e. a formal scheme is not necessarily noetherian.

Sheaves will be denoted with calligraphic letters and their global sections with the associated straight letter. 
By default, a module means a left module. 
We denote by a hat the $p$-adic completion and if 
$\E$ is an abelian sheaf of groups, we set $\E _{\Q}:= \E \otimes _{\Z} \Q$. 
Let  $\AA$ be a sheaf of rings.
If $*$ is one of the symboles $+$, $-$, ou $\mathrm{b}$, 
$D ^* ( \AA )$ means the derived category of the complexes of 
(left) $\AA$-modules satisfying the corresponding condition of vanishing of cohomological spacesd. 
When we would like to clarify between right and left,  we will write
$D ^* ( {} ^{\mathrm{l}}\AA )$ or $D ^* ( {} ^{\mathrm{r}}\AA )$.
We denote by $D ^{\mathrm{b}} _{\mathrm{coh}} ( \AA )$
the subcategory of  $D  ( \AA )$
of bounded and coherent complexes.

Formal $\fS$-schemes will be indicated with gothic letters and their reduction modulo $\pi$ 
with the associated roman letter. 
Formal $\fS$-schemes or 
 $S$-schemes will be quasi-compact and separated. 
Finally, when 
$f \colon \X  \to \cP $
is a morphism of formal  $\fS $-schemes,
for any integer $i \in \N$,
we denote by  $ f _{i}   \colon X  _i \to P  _i$
the induced morphism modulo $\pi ^{i+1}$.

If $A$ is $k$-algebra, we denote by 
$A ^{(p)}$ be the $A$-algebra given by the absolute Frobenius
$F _A \colon A \to A ^{(p)} $.
We denote by $A ^p$ the image of the Frobenius homomorphism
$F _A \colon A \to A$. 
Unless otherwise stated, 
we suppose that $k ^{(p)}$ is a finite $k$-algebra 
(see \ref{lem-k[[t]]pbases} to see why we need this hypothesis).

We fix a Cohen algebra $C (k)$ with residue field $k$.
For any $i \in \N$, set $C _i (k):= C (k ) / p ^{i+1}C (k )$,
$V _i := \V / \pi ^{i+1} \V$. 
Let  $S _i= \Spec (V _i)$. 
We denote by $\mathbb{D} ^n _{S _i}:= \Spec ( V _i [[ T _1, \dots, T _n]])$
and by
$\mathbb{D} ^n _{\fS}:= \Spf ~( \cV  [[ T _1, \dots, T _n]])$ the formal $\fS$-scheme.

\section{Complements on formal smoothness}

\subsection{Relatively perfect morphisms}

\begin{empt}
[Around $p$-bases]
\label{aroundpbases}
Let $A \to B$ be an homomorphism of $k$-algebras.
Let $(b _i ) _{i\in I} \in B$ be some elements and
$A [ \underline{T}]= A [ T _i ; \; i \in I] \to B$ be the $A$-algebra homomorphism given by 
$T _i \mapsto b _i$.
\begin{enumerate}[(a)]
\item 
\label{aroundpbases1}
The homomorphism 
$A [ \underline{T}] \to B$
is relatively perfect in the sense of Kato if and only if 
the canonical homomorphism
$\left ( A [ \underline{T}] \right )  ^{(p)}
\otimes _{A [ \underline{T}]}
B
\to 
B ^{(p)}$
is an isomorphism.

\item 
\label{aroundpbases2}
We have the canonical isomorphisms
\begin{equation}
\notag
\left ( A [ \underline{T}] \right )  ^{(p)}
\otimes _{A [ \underline{T}]}
B
\riso 
\left ( A [ \underline{T}] \right )  ^{(p)}
\otimes _{A ^{(p)} [ \underline{T}]}
\left ( A ^{(p)} [ \underline{T}] 
\otimes _{A [ \underline{T}]}
B \right )
\riso 
\left ( A [ \underline{T}] \right )  ^{(p)}
\otimes _{A ^{(p)} [ \underline{T}]}
\left ( A ^{(p)} 
\otimes _{A} 
B \right ).
\end{equation}
This yields that 
the homomorphism 
$A [ \underline{T}] \to B$
is relatively perfect in the sense of Kato if and only if 
the canonical homomorphism
\begin{equation}
\label{rel-perf-prop0}
\left ( A [ \underline{T}] \right )  ^{(p)}
\otimes _{A ^{(p)} [ \underline{T}]}
\left ( A ^{(p)} 
\otimes _{A} 
B \right )
\to 
B ^{(p)}
\end{equation}
is an isomorphism.
The $A ^{(p)} [ \underline{T}]$-module
$\left ( A [ \underline{T}] \right )  ^{(p)}$ is free with the basis
$\prod _{i\in I} T _i ^{k _i}$, for $k _i<p$ for any $i$
and $(k _i) _{i\in I}$ has finite support.
Hence, 
the homomorphism  \ref{rel-perf-prop0} 
is an isomorphism if and only if 
$B ^{(p)}$
is a free 
$A ^{(p)} \otimes _{A} B$-module with the basis $\prod _{i\in I}  b _i ^{k _i}$, for $k _i<p$ for any $i$
and $(k _i) _{i\in I}$ has finite support.
In that case, following Kato's terminology, 
$(b _i) _{i\in I}$ forms a $p$-basis of $B/A$.

\item 
The image of the homomorphism 
$A ^{(p)} \otimes _{A} B \to B ^{(p)}$
is 
$A ^{(p)} [ F _B (B )]$ (which is equal
to $A [ B ^p ]$ if we forget $(p)$).
Recall that following \cite[0.21.1.9]{EGAIV1},
$(b _i) _{i\in I}$ is a $p$-basis of $B/A$
if $B ^{(p)}$ is a free $A ^{(p)} [ F _B (B )]$-module with the basis
$\prod _{i\in I}  b _i ^{k _i}$, for $k _i<p$ for any $i$
and $(k _i) _{i\in I}$ has finite support.

\item 
\label{aroundpbases4}
Hence, 
$(b _i) _{i\in I}$ forms a $p$-basis of $B/A$ in the sense of Kato in \cite[1.3]{Kato-explicity-recip91}
if and only if the homomorphism 
$A ^{(p)} \otimes _{A} B \to B ^{(p)}$ is injective and 
$(b _i) _{i\in I}$ forms a $p$-basis of $B/A$ in the sense of Grothendieck in \cite[0.21.1.9]{EGAIV1}.

\end{enumerate}

\end{empt}

\begin{rem}
\label{rem-aroundpbases}
We keep notation of \ref{aroundpbases}.
When $B/A$ is relatively perfect, then 
$B/A$ is formally étale (this is hidden in \cite[0.21.2.7]{EGAIV4} 
since this corresponds to the case where
the basis is empty, a proof can be found in \cite{Caro-Vauclair}). 
Moreover, following Theorem \cite[0.21.2.7]{EGAIV4},
 if $(b _i) _{i\in I}$ forms a $p$-basis of $B/A$ in the sense of Kato 
then $B/A$ is formally smooth.
When $(b _i) _{i\in I}$ forms a $p$-basis of $B/A$ in the sense of Grothendieck 
then $B/A$ is not necessarily formally smooth, which justifies why we prefer Kato's definition.
\end{rem}

\begin{ex}
\label{ex-ext-fieldpbases}
Let $K/k$ be a field extension of characteristic $p>0$.
Following \cite[0.21.4.2]{EGAIV1}, 
the extension $K/k$ has always a $p$-basis in the sense of Grothendieck. 
The following statement are equivalent.
\begin{enumerate}[(i)]
\item \label{ex-ext-fieldpbases1} The extension $K/k$ is separable. 
\item \label{ex-ext-fieldpbases2} The extension $K/k$ is formally smooth (for the discrete topology). 
\item \label{ex-ext-fieldpbases3} The extension $K/k$ has a $p$-basis in the sense of Kato.
\end{enumerate}
The equivalence $\ref{ex-ext-fieldpbases1} \Leftrightarrow \ref{ex-ext-fieldpbases2}$ is Cohen theorem 
(see \cite[0.19.6.1]{EGAIV1}). 
The implication 
$\ref{ex-ext-fieldpbases3} \Rightarrow \ref{ex-ext-fieldpbases2}$
is \cite[0.21.2.7]{EGAIV1}  (modulo the remark \ref{aroundpbases}.\ref{aroundpbases4}). 
It remains to check  
$\ref{ex-ext-fieldpbases1} \Rightarrow \ref{ex-ext-fieldpbases3}$.
Suppose $K/k$ is separable.
Then 
$k ^{(p)} \otimes _{k} K \to K ^{(p)}$
(i.e. $k \otimes _{k ^p} K ^p \to K $)
is injective. Indeed, if $a \in k ^{(p)} \otimes _{k} K$ is in the kernel of 
$k ^{(p)} \otimes _{k} K \to K ^{(p)}$, 
then $a ^p = 0$. By definition, since $K/k$ is separable then
$ k ^{(p)} \otimes _{k} K$ is reduced. Hence $a =0$ and we are done. 
Following
\ref{aroundpbases}.\ref{aroundpbases4} and \cite[0.21.4.2]{EGAIV1}, 
this yield that the extension $K/k$ has a $p$-basis (in the sense of Kato).
\end{ex}

\begin{lem}
\label{lem-k[[t]]pbases}
Let $A _0$ be a $k$-algebra such that
the absolute Frobenius
$F _{A _0}
\colon 
A _0 \to A _0 ^{(p)}$ is finite. 
Set $A _0[[ \underline{T}]]
:=
A _0[[ T _1, \dots, T _n]]$.
Then $T _1, \dots, T _n$ forms a finite $p$-basis of 
$A _0[[ \underline{T}]]/A _0$.
\end{lem}

\begin{proof}
We have to check that 
$A _0[ \underline{T}]
\to 
A _0[[ \underline{T}]]$
is relatively perfect, i.e. following \ref{aroundpbases}.\ref{aroundpbases2}
that 
the canonical homomorphism
$$\left ( A _0[ \underline{T}] \right ) ^{(p)} 
\otimes _{A _0 ^{(p)}[ \underline{T}]}
\left (A _0 ^{(p)} \otimes _{A _0}  A _0[[ \underline{T}]]  \right ) 
\to 
\left ( A _0[[ \underline{T}]] \right ) ^{(p)}$$
is an isomorphism, i.e. 
$\left ( A _0[[ \underline{T}]] \right ) ^{(p)}$
is a free 
$A _0 ^{(p)} \otimes _{A _0}  A _0[[ \underline{T}]]$-module with the basis 
$\prod _{i=1} ^{n} T _i ^{k _i}$, for $k _i<p$ for any $i$.
Since 
$F _{A _0}$ 
is finite, then 
we have the equality
$A _0 ^{(p)} \otimes _{A _0}  A _0[[ \underline{T}]] 
=
A _0 ^{(p)} [[ \underline{T}]]$.
The induced map 
$A _0 ^{(p)} [[ \underline{T}]]
\to 
\left ( A _0[[ \underline{T}]] \right ) ^{(p)}$
is given by 
$\sum a _k T ^k \to 
\sum a _k T ^{pk}$.
Hence, we conclude by an easy computation.
\end{proof}

\begin{lem}
\label{pbasis-Kato}
Let $Y$ be an $S  _i$-scheme.
Then, the canonical morphism
$\mathbb{D} ^n _Y \times _Y \A ^{n'} _Y
\to \A ^{n+n'} _Y$
is relatively perfect and 
$\Omega _{\mathbb{D} ^n _Y \times _Y \A ^{n'} _Y/Y} $ is $\O _{\mathbb{D} ^n _Y \times _Y \A ^{n'} _Y}$-free of rank $n+n'$. 

\end{lem}
\begin{proof}
Since 
$C _i (k)[T _1, \dots, T _n]$ is noetherian, then
the homomorphism 
$C _i (k)[T _1, \dots, T _n] \to C _i (k)[[ T _1, \dots, T _n]]$
given by the $(T _1,\dots,T _n)$-adic completion is flat.
Following \ref{lem-k[[t]]pbases}, 
$T _1, \dots, T _n$ forms a finite $p$-basis of 
$k[[ T _1, \dots, T _n]]$ over $k$, 
Using \cite[1.6]{Kato-explicity-recip91}, this yields that
$C _i (k)[T _1, \dots, T _n] \to C _i (k)[[ T _1, \dots, T _n]]$
is relatively perfect.
Since $C _i (k) \to V _i$ is finite, 
by applying the functor $V _i \otimes _{C _i (k)}-$, 
we get the relatively perfect homomorphism
$V _i [T _1, \dots, T _n] \to V _i[[ T _1, \dots, T _n]]$.
In other words, 
$\mathbb{D} ^n _{S _i}
\to \A ^{n} _{S _i}$
is relatively perfect. 
Since relatively perfect morphism are closed under base change,
this yields that 
$f \colon \mathbb{D} ^n _Y \times _Y \A ^{n'} _Y
\to \A ^{n+n'} _Y$
is relatively perfect.
In particular, it is formally étale and 
then the canonical morphism
$f ^* 
\Omega _{\A ^{n+n'} _Y/Y} 
\to 
\Omega _{\mathbb{D} ^n _Y \times _Y \A ^{n'} _Y/Y} $
is an isomorphism  (see \cite[17.2.4]{EGAIV4}).
\end{proof}

\begin{empt}
Set 
$\widehat{\Omega} _{\V[[ T _1, \dots, T _n]]/\V}
:=
\underleftarrow{\lim} _i\Omega _{R _i[[ T _1, \dots, T _n]]/R _i}$.
Then,
$d T _1, \dots, d T _n$ is a basis of 
of the free $\V[[ T _1, \dots, T _n]]$-module
$\widehat{\Omega} _{\V[[ T _1, \dots, T _n]]/\V}$.

\end{empt}

Recall the following definitions. 
\begin{dfn}
\label{dfn-pbasispadic}
\begin{enumerate}[(a)]
\item Let $X \to Y$ be a morphism of $V_i$-schemes. 
Let $t _1, \dots, t _d$ be elements of $\Gamma (X, \O _{X})$.
According to definition \cite[1.2]{Kato-explicity-recip91},
we say that  $t _1, \dots, t _d$  form a finite $p$-basis over $Y$ is 
the morphism 
$\X \to \A ^d _{S_i} \times Y$ is relatively perfect, i.e. 
if $X \to \A ^d _{S_i} \times Y$ 
is formally étale and its reduction modulo $p$ is relatively perfect
in the sense that the square given by the respective absolute Frobenius is cartesian.

 \item Let $f \colon \X \to \Y$ be a morphism of $\V$-formal schemes. 
 We say that $f$ is relatively perfect if $f$ is formally étale and 
 $f _0$ 
 is relatively perfect in the sense of Kato. 
 Beware that if this is not clear that 
 if $f $ is relatively perfect then so is $f _e$. 

\item Let $\X \to \Y$ be a morphism of $\V$-formal schemes. 
Let $t _1, \dots, t _d$ be elements of $\Gamma (\X, \O _{\X})$.
We say that  $t _1, \dots, t _d$  form a finite $p$-basis over $\Y$ is 
the morphism 
$\X \to \widehat{\A} ^d _{\V} \times \Y$ is relatively perfect, i.e. 
if $\X \to \widehat{\A} ^d _{\V} \times \Y$ is formally étale and 
$X _0\to \A ^d _{Y _0}$ is relatively perfect.

\end{enumerate}
\end{dfn}

\subsection{Semistable schemes over $k[[t]]/k$}
Let $R:= k[[t]]$.
We denote by $\eta$ (resp. $s$) the generic (resp. closed) point of $\Spec R$.
We recall the following definition. 
\begin{empt}
[Strictly semi-stable over $\Spec R$]
\label{dfn-sss}
Let  $X$ be an integral $\Spec R$-scheme of finite type. 
Let $X _i$, $i\in I$ be the irreducible components of $X_s$. Put $X _J:= \cap _{j \in J} X _j$ (scheme-theoretic intersection), for a nonempty subset $J$ of $I$. 
We recall that $X$ is ``strictly semi-stable over $\Spec R$''
 means that $X/\Spec R$ satisfy the following properties  (see \cite[2.16]{dejong}, and also 40.21.1--2 of the stack project) : 
\begin{enumerate}[(a)]
\item $X _\eta$ is smooth over $\kappa (\eta)$, 
\item $X _s $ is a reduced scheme, i.e. $X _s = \cup _{i\in I} X _i$ scheme-theoretically,
\item for each $i\in I$, $X _i$ is a divisor on $X$, 
\item for each nonempty $J \subset I$, the scheme $X_J$  is smooth over $k (s)$ and has codimension $\sharp J$ in $X$.
\end{enumerate}
\end{empt}

\begin{rem}
\label{locdes-ssv}
Let $X$ be a strictly semi-stable scheme over $\Spec R$.
\begin{enumerate}[(a)]
\item Remark that when $k= \kappa (s)$ is perfect, then conditions $2), 3), 4)$ are equivalent 
to say that $X _s$ is a divisor with strict normal crossing on $X$ (see the definition \cite[2.10]{dejong}).

\item Let $x \in X _s$. Let $X _1, \dots, X _n $ be the irreducible components of $X _s$ containing $x$.
Then there exists an open neighborhood $U$ of $x$ and a smooth morphism
$U \to \Spec R [ t _1, \dots, t _n] / (t - t _1 \cdots t _n)$
such that $X _i \cap U$ is given by $t _i = 0$ and (this is a consequence of the local description of 
\cite[2.16]{dejong} and of \cite[17.5.3]{EGAIV4}).

\end{enumerate}

\end{rem}

\begin{lem}
\label{lem-locdesc-sstvar}
We have the canonical cartesian diagram
\begin{equation}
\label{lem-locdesc-sstvar-diag}
\xymatrix{
{\Spec R [ t _1, \dots, t _n] / (t - t _1 \cdots t _n)} 
\ar[r] ^-{f}
\ar@{}[rd] |{\square}
& 
{\A ^{n} _k} 
\\ 
{V (t)} 
\ar[r] ^-{}
\ar@{^{(}->}[u] ^-{}
& 
{V(T _1 \cdots T _n) } 
\ar@{^{(}->}[u] ^-{u}
}
\end{equation}
where $\A ^{n} _k = \Spec k [ T _1,\dots, T _n]$,
$f$ is the morphism given by $T _i \mapsto t _i \mod t _1 \cdots t _n -t$,
and
$u$ is the closed immersion induced by 
$T _1 \cdots T _n$.
Moreover, $f$ is a relatively perfect morphism.
\end{lem}

\begin{proof}
The cartesianity of the diagram \ref{lem-locdesc-sstvar-diag} is straighforward.
It remains to check that $f$ is relatively perfect.
The morphism 
$k [t, t _1, \dots, t _n] 
\to 
k [T  _1, \dots, T _n]$
given by 
$t \mapsto T  _1 \cdots T _n$
and by 
$t _i \mapsto T _i$, induces the isomorphism
$k [t, t _1, \dots, t _n]/( t _1 \cdots t _n -t)
\riso
k [T  _1, \dots, T _n]$.
Since $k [t] \to k [[t]]$ is relatively perfect, 
since the relative perfectness is stable under base change, 
the canonical morphism
$k [T  _1, \dots, T _n] \liso 
k [t, t _1, \dots, t _n]/( t _1 \cdots t _n -t)
\to 
k [[t]] [t _1, \dots, t _n]/( t _1 \cdots t _n -t)$
is relatively perfect. 
This morphism sends 
$T  _1 \cdots T _n$
to $t$. Hence, we are done.
\end{proof}

\begin{prop}
\label{locdesc-sstvar}
Let $X$ be a semistable scheme over $\Spec R$.
Zariski locally on $X$, there exists a cartesian diagram of the form
\begin{equation}
\notag
\xymatrix{
{X} 
\ar[r] ^-{g}
\ar@{}[rd] |{\square}
& 
{\A ^d _k} 
\\ 
{X _s} 
\ar[r] ^-{}
\ar@{^{(}->}[u] ^-{}
& 
{V(t _1 \cdots t _n) } 
\ar@{^{(}->}[u] ^-{u}
}
\end{equation}
where $g$ is a relatively perfect morphism,
$n \leq d$ are two integers,
$\A ^d _k = \Spec k [ t _1,\dots, t _d]$,
and
$u$ is the closed immersion induced by 
$t _1 \cdots t _n$.
\end{prop}

\begin{proof}
Following the second remark of \ref{locdes-ssv}, 
Zariski locally on $X$, there exists a smooth morphism of the form 
$X \to \Spec R [ t _1, \dots, t _n] / (t - t _1 \cdots t _n)$.
Hence, Zariski locally on $X$, there exists an étale morphism of the form 
$X 
\to 
\Spec  (R [ t _1, \dots, t _n] / (t - t _1 \cdots t _n) ) \times _{\Spec k }\A ^m _k$.
Following  Lemma \ref{lem-locdesc-sstvar}, 
we get a relatively perfect morphism
$\Spec  (R [ t _1, \dots, t _n] / (t - t _1 \cdots t _n) ) \times _{\Spec k }\A ^m _k
\to 
\A ^{n+m} _k$. Hence, we are done.
\end{proof}

\begin{empt}
[Strictly semi-stable pairs over $\Spec R$]
\label{dfn-ssp}
We recall from \cite[6.3]{dejong} that $(X,Z)$ is a {\it strict semi-stable pair over $\Spec R$} if 
\begin{enumerate}[(a)]
\item $X $ is strict semi-stable over $S$ (see \ref{dfn-sss}), 
\item Let $Z _f := \cup _{i\in I} Z _i$ be the union of the irreducible components $Z _i$ of $Z$ which dominate $S$.
For each $J \subset I$, the scheme $Z _J:= \cap _{j\in J}Z _j$ is a disjoint union of strict semi-stable schemes over $S$.
\item $Z$ is a divisor with strict normal crossing on $X$ and $Z= Z _f \cup X _s$. 
\end{enumerate}

\end{empt}

\begin{rem}
\label{ssp-locdesc}
Let $(X,Z)$ be a strict semi-stable pair over $\Spec R$. 
Let $x \in X _s$. Let $X _1, \dots, X _n $ be the irreducible components of $X _s$ containing $x$ 
and $Z _1, \dots, Z _m$ be the irreducible components of $Z _f$ containing $x$.
Then there exist an open neighborhood $U$ of $x$ and a smooth morphism
$U \to \Spec R [ t _1, \dots, t _n, s _1, \dots s _m] / (t - t _1 \cdots t _n)$
such that $X _i \cap U$ is given by $t _i = 0$ and 
$Z _j \cap U$ is given by $s _j =0$ (this is a consequence of the local description of 
\cite[6.4]{dejong} and of \cite[17.5.3]{EGAIV4}.
\end{rem}

\begin{prop}
\label{proplocdesc-sstp}
Let $(X,Z)$ be a strict semi-stable pair over $\Spec R$. 
Zariski locally on $X$, there exist integers $n,m,d$ such that $n+m \leq d$,
there exists a relatively perfect morphism
$f \colon X \to \A ^d _k$ with 
$\A ^d _k = \Spec k [ t _1,\dots, t _d]$
such that 
$X _s = f ^{-1} (V(t _1 \cdots t _n) )$
and $Z _f= f ^{-1} (V(t _{n+1} \cdots t _{n+m}) )$. 
\end{prop}

\begin{proof}
i) Following the remark of \ref{ssp-locdesc}, 
Zariski locally on $X$, 
there exists a smooth morphism
$X \to \Spec R [ t _1, \dots, t _n, s _1, \dots s _m] / (t - t _1 \cdots t _n)$
such that $X _i$ is given by $t _i = 0$ and 
$Z _j $ is given by $s _j =0$.
Hence, 
Zariski locally on $X$, 
there exists an  étale morphism of the form
\begin{equation}
\label{const-rel-perf-h1}
X \to \Spec R [ t _1, \dots, t _n, s _1, \dots s _m, r _1, \dots ,r _l] / (t - t _1 \cdots t _n) 
\end{equation}
such that $X _i $ is given by $t _i = 0$ and 
$Z _j $ is given by $s _j =0$.

ii) Consider  the following canonical cartesian diagram
\begin{equation}
\label{lem-locdesc-sstp-diag}
\xymatrix{
{\Spec R [ t _1, \dots, t _n, s _1, \dots s _m, r _1, \dots ,r _l] / (t - t _1 \cdots t _n)  } 
\ar[r] ^-{h}
\ar@{}[rd] |{\square}
& 
{\A ^{n} _k \times \A ^{m} _k \times \A ^l _k } 
\\ 
{V (s _1\dots s _m)} 
\ar[r] ^-{}
\ar@{^{(}->}[u] ^-{}
& 
{\A ^{n} _k \times V(S _{1} \cdots S _{m}) \times \A ^l _k } 
\ar@{^{(}->}[u] ^-{u}
}
\end{equation}
where $\A ^{n} _k = \Spec k [ T _1,\dots, T _{n}]$,
$\A ^{m} _k = \Spec k [ S _1,\dots, S _{m}]$,
$\A ^{l} _k = \Spec k [ R _1,\dots, R _{l}]$,
$h$ is the morphism given by $T _i \mapsto t _i \mod t _1 \cdots t _n -t$,
$S _j \mapsto s _j \mod t _1 \cdots t _n -t$,
and
$R _k \mapsto r _k \mod t _1 \cdots t _n -t$,
and
$u$ is the closed immersion induced by 
$S _1 \cdots S _m$.
The morphism $h$ of \ref{lem-locdesc-sstp-diag} is induced by base change from the morphism $f$ of 
the diagram \ref{lem-locdesc-sstvar-diag}. 
Hence, since $f$ is relatively perfect, then so is $h$.  
We conclude by taking the composition of $h$ with \ref{const-rel-perf-h1}. 
\end{proof}

\subsection{Formally  smooth morphisms of formal finite type}

\begin{prop}
\label{regularity/formalsm}
Let $X$ be a noetherian formally smooth scheme over $\Spec k$.
Then $X$ is 
geometrically 
regular. 
\end{prop}

\begin{proof}
This is a consequence of  \cite[0.19.3.8]{EGAIV1}
and \cite[0.19.6.5]{EGAIV1}.
\end{proof}

\begin{dfn}
\label{dfn-fft}
\begin{enumerate}[(a)]
\item Let $f \colon X \to Y$ be  an $S_i$-morphism.
We say that the $f$ is 
{\it an $S_i$-morphism locally of formal finite type} if Zariski locally on $X$, 
there exist $n,n' \in \N$,  and a closed $Y$-immersion of the form
$X \hookrightarrow \mathbb{D} ^n _{S _i} \times _{S _i} \A ^{n'} _Y$.
Moreover, the notion of ``locally of formal finite type'' depends on the field $k$ but does not depend on $i$ : 
this means that if $f$ can also be viewed as an $S_{i+1}$-morphism then $f$ is of formal finite type as $S_i$-morphism if and only if 
$f$ is of formal finite type as $S_{i+1}$-morphism.
When $Y=S_i$ and $f$ is the structural morphism of $X$,
we say that $X$ is  an $S_i$-scheme locally of formal finite type. 

\item Let $f \colon X \to Y$ be  an $S_i$-morphism.
We say that the $f$ is 
{\it an $S_i$-morphism of formal finite type} 
if there exists an integer $n \geq 0$, a finite type $Y$-morphism of the form 
$g \colon X \to \mathbb{D} ^n _{S _i} \times _{S _i} Y$.
When $Y=S_i$ and $f$ is the structural morphism of $X$,
we say that $X$ is  an $S_i$-scheme of formal finite type.

\item Let $f \colon \fX \to \fY$ be  a morphism of formal $\fS$-schemes.
We say that the $f$ is 
{\it an $\fS$-morphism locally of formal finite type} if Zariski locally on $\fX$, 
there exist $n,n' \in \bbN$,  and a closed $\fY$-immersion of the form
$\fX \hookrightarrow \mathbb{D} ^n _{\fS} \times _{\fS } \widehat{\A} ^{n'} _\fY$.
When $\fY = \fS$ and $f$ is the structural morphism of $\fX$,
we say that $\X$ is a formal $\fS$-scheme locally of formal finite type.

\item Let $f \colon \fX \to \fY$ be  a morphism of formal $\fS$-schemes.
We say that the $f$ is 
{\it an $\fS$-morphism of formal finite type} if 
there exists an integer $n \geq 0$, a finite type morphism of $\fY$-schemes of the form
$\fX \to \mathbb{D} ^n _{\fS} \times _{\fS } \fY$.
When $\fY = \fS$ and $f$ is the structural morphism of $\fX$,
we say that $\X$ is a formal $\fS$-scheme of formal finite type.

\end{enumerate}

\end{dfn}

\begin{rem}
\label{rem-Dnn'Ynoeth}
Let $f \colon X \to Y$ be  an $S_i$-morphism.

\begin{enumerate}[(i)]
\item 
\label{rem-Dnn'Ynoeth-i)}
Suppose $Y$ is  an $S_i$-scheme of finite type.
Then $f$ is an $S_i$-morphism  locally of formal finite type if and only if $X$ is  an $S_i$-scheme locally of formal finite type. 
In that case, $X$ is noetherian 
(because so are
$\mathbb{D} ^n _{S _i} \times _{S _i} \A ^{n'} _Y$ for any integers $n$ and $n'$)
and $X \times _{S _i}Y$ is also locally of formal finite type and noetherian. 
For short, we say that 
$X$ is  an $Y$-scheme locally of formal finite type.

\item Beware that when $Y$ is not of finite type over $S_i$, then 
$\mathbb{D} ^n _{S _i} \times _{S _i} \A ^{n'} _Y$ is not necessarily noetherian
even if $n'= 0$ and $Y = \mathbb{D} ^r _{S _i}$ for $r \geq 1$.

\item Beware that if $Y$ is an $S _i$-scheme locally of formal finite type and $f$ is an $S _i$-morphism locally of formal finite type then
this is not clear that $X$ is an $S _i$-scheme locally of formal finite type.
\end{enumerate}

\end{rem}

\begin{prop}
\label{EGA17.2.3}
Let 
$f\colon X \to Y$
be   an $S_i$-morphism locally of formal finite type.
\begin{enumerate}[(a)]
\item Then 
$\Omega _{X/Y}$ is an $\O _X$-module of finite type. 
\item If $f$ is moreover formally smooth, then 
$\Omega _{X/Y}$ is an $\O _X$-module locally free of finite type. 
\item 
\label{EGA17.2.3-3}
Let $x _0  \in X$. The following assertions are equivalent
\begin{enumerate}[(a)]
\item 
There exist an open subset $U$ of $X$ containing $x _0$ such that 
$U \to Y$ is formally unramified.
\item $\Omega _{X/Y, x _0}=0$.
\item $\Omega _{X/Y,x _0} \otimes _{\O _{X _{x _0}}} k (x _0)=0$.
\end{enumerate}

\end{enumerate}
\end{prop}

\begin{proof}
Since this is local, we can suppose 
there exist $n,n' \in \N$ and a closed $Y$-immersion of the form
$u\colon X \hookrightarrow \mathbb{D} ^n _{S _i} \times _{S _i} \A ^{n'} _Y$.
Hence, using \cite[16.4.21]{EGAIV4},
we get the epimorphism of $\O _X$-modules
$u ^* \Omega _{\mathbb{D} ^n _{S _i} \times _{S _i}   \A ^{n'} _Y/Y} 
\twoheadrightarrow\Omega _{X/Y}$.
Since $\Omega _{\mathbb{D} ^n  _{S _i} \times _{S _i}   \A ^{n'} _Y/Y} $ is $\O _{\mathbb{D} ^n  _{S _i} \times _{S _i}   \A ^{n'} _Y}$-free of rank $n+n'$ 
(see \ref{pbasis-Kato}),
we conclude the first part. 

The second assertion (resp. third assertion) is a consequence of \cite[17.2.3.(i)]{EGAIV4}
(resp. \cite[17.2.1]{EGAIV4}) and of the first assertion.
\end{proof}

\begin{prop}
\label{EGAIV18.1.1}
Let $u\colon 
Y _0
 \hookrightarrow 
 Y $ 
 be a closed immersion of $S_i$-schemes of finite type.
Let 
$f _0 \colon X _0 \to Y _0$ be a formally smooth locally of formal finite type morphism.
Let $x _0 \in X _0$.

\begin{enumerate}[(a)]
\item There exist an open subset $U _0$ of $X _0$ containing $x _0$, 
and a formally smooth locally of formal finite type $S_i$-morphism  $f \colon U \to Y$ 
together with an isomorphism of the form
$U \times _{Y} Y _0 \riso X _0$.
\item Moreover, 
if $f _0$ is formally étale then so is such $f$.
\end{enumerate}

\end{prop}

\begin{proof}
Using \ref{EGA17.2.3}.\ref{EGA17.2.3-3}, we reduce to check the first assertion. 
We can follow the proof of \cite[18.1.1]{EGAIV4}.
Since this is local, we can suppose $Y _0 = \Spec ( A _0) $,
$Y = \Spec ( A) $, and
there exist $n,n' \in \N$ and a closed $Y _0$-immersion of the form
$\iota\colon X _0 \hookrightarrow 
\mathbb{D} ^n _{S _i} \times _{S _i} \A ^{n'} _{Y _0}$.
Put 
$B _0:= A _0 
\otimes _{V _i}  V _i  [[t _1, \dots, t _n  ]] \otimes _{V _i}  V _i [t _{n+1},\dots, t _{n+n'}] $
and 
$B:= A \otimes _{V _i}  V _i  [[t _1, \dots, t _n  ]] \otimes _{V _i}  V _i [t _{n+1},\dots, t _{n+n'}] $
be respectively the global section of the structural sheaf of 
$\mathbb{D} ^n _{S _i} \times _{S _i} \A ^{n'} _{Y _0}$
and
$\mathbb{D} ^n _{S _i} \times _{S _i} \A ^{n'} _{Y}$.
Let $I _0 $ be an ideal of $B _0$ 
such that
$X _0 = \Spec ( B _0  /I _0)$. 
Set $C _0 := B _0  /I _0$.
Since $C _0/A _0$ is formally smooth, following
\cite[0.20.5.14]{EGAIV1}, the sequence of $C _0$-modules
\begin{equation}
\label{EGAIV18.1.1-ses}
0 \to 
I _0 / I _{0} ^2 
\to 
\Omega _{B _0 /A _0}
\otimes _{B _0}
C _0
\to 
\Omega _{C _0 /A _0}
\to 
0
\end{equation}
is split exact.

Let $\mathfrak{p} _0 $ be a prime ideal of $C _0$,
$\mathfrak{q}  _0 $ (resp. 
$\mathfrak{q}$)
be the induced prime ideal of $B _0$
(resp. $B$).
The $B _0$-module
$\Omega _{B _0 /A _0}$
is free with the basis 
$d t _1,\dots, d t _{n+n'}$ (see \ref{pbasis-Kato}).
We denote by $\partial _1, \dots, \partial _{n+n'}$ the corresponding derivations.
Since the sequence \ref{EGAIV18.1.1-ses} is split exact, 
then there exist
$x _1, \dots, x _m \in I _0$ such that
the images of $x _1, \dots, x _m$ in 
$( I _0 / I ^2  _0 ) _{\mathfrak{p} _0}$
generate 
$( I _0 / I ^2 _0  ) _{\mathfrak{p} _0}$
and such that there exist
$n _1, \dots, n _m \in \{1,\dots,n+n' \}$
satisfying 
$\det (\partial _{n _i} ( x _j))
\not \in \mathfrak{q} _0$.
Since $B _0$ is noetherian (see \ref{rem-Dnn'Ynoeth}.\ref{rem-Dnn'Ynoeth-i)}), $I _0$ is a $B _0$-module of finite type.
Hence, 
since 
$( I _0 / I ^2  _0 ) _{\mathfrak{p} _0}=
 I  _{0,\mathfrak{q} _0}/ I  ^2 _{0,\mathfrak{q} _0}$,
 since
 $ I  _{0,\mathfrak{q} _0}\subset \mathfrak{q} _0 B _{0, \mathfrak{q} _0}$
 then 
using Nakayama lemma, 
the image of $x _1, \dots, x _m$ in 
$ I _{0,\mathfrak{q} _0}$
generates 
$ I  _{0,\mathfrak{q} _0}$.
Let $f _0 \in B _0 \setminus \mathfrak{p} _0$ 
such that the image of $x _1, \dots, x _m$ in $B _{0, f _0}$ 
generate
$ I  _{0, f _0}$.
Let $U _0:= \Spec (B _{0, f _0} /  I  _{0, f _0})$.

Let $f$ and $y _1,\dots, y _m\in B$ whose images in $B _0$ are
$f _0$ and $x _1, \dots, x _m$.
Let $I $ be the ideal of $B$ generated by 
$y _1,\dots, y _m$ and $C := B/I$.
Since
$B _f/I _f \otimes _{A} A _0 \riso 
B _{0, f _0} /I _f B_{0, f _0}
= 
B _{0, f _0} /I  _{0, f _0}$,
since 
$\mathfrak{q} _{0, f _0}$ contains $I _{0, f _0}$
then 
$\mathfrak{q} _f$ contains $I _f$
and 
$\mathfrak{p}:=\mathfrak{q} _f/I _f$
is the prime ideal of 
$B _f/I _f$ corresponding to 
$\mathfrak{p} _{0, f _0}$ via the closed immersion 
$U _0= \Spec (B _{0, f _0} /I  _{0, f _0} )\hookrightarrow \Spec (B _f/I _f )=:U$.
Since
$\det (\partial _{n _i} ( y _j))
\not \in \mathfrak{q} $
(and then 
$\det (\partial _{n _i} ( y _j)) \mod I _f
\not \in \mathfrak{p} $),
since 
the $B _f$-module
$\Omega _{B _f  /A }$
is free,
then 
using 
\cite[0.19.1.12]{EGAIV1},
we check that 
$( I _f  / I _f ^2 ) _{\mathfrak{p}}
\to 
\Omega _{B _f  /A}
\otimes _{B _f}
(C_f) _{\mathfrak{p}}$
is left invertible.
Using
\cite[0.22.6.4]{EGAIV1}, we conclude that 
replacing $f$ by a multiple if necessary, 
we have that $U$ is formally smooth over $Y= \Spec A$.
Hence,  we are done.
\end{proof}

\begin{lem}
\label{EGAIV18.1.2-fullyfaithful}
Let $Y _0 \hookrightarrow Y$ be a nilpotent closed $S_i$-immersion
of schemes of finite type. 
Let $X$ and $X'$ be $Y$-schemes. 
Suppose $X'$ is  formally étale over $Y$. 
Put 
$X _0 := X \times _{Y} Y _0$
and
$X '_0 := X ' \times _{Y} Y _0$.
Then the canonical map
$\mathrm{Hom} _{Y } ( X  , X ' )
\to 
\mathrm{Hom} _{Y _0} ( X _0 , X ' _0)$
is a bijection. 
\end{lem}

\begin{proof}
An element of $\mathrm{Hom} _{Y _0} ( X _0 , X ' _0)$ is equivalent to the data of a morphism
$X _0 \to X'$ making commutative the diagram
\begin{equation}
\notag
\xymatrix{
{X}
\ar[r] ^-{} 
& 
{Y} 
\\ 
{X _0} 
\ar[r] ^-{} 
\ar@{^{(}->}[u] ^-{} 
& 
{X',} 
\ar[u] ^-{}
}
\end{equation}
whose three other morphisms are the structural or canonical ones.
Since $X'/Y$ is formally étale, since the closed immersion
$X _0 \hookrightarrow X$ is nilpotent, 
this yields that 
the canonical map
$\mathrm{Hom} _{Y } ( X  , X ' )
\to 
\mathrm{Hom} _{Y _0} ( X _0 , X ' _0)$
is a bijection.\end{proof}

\begin{thm}
\label{EGAIV18.1.2}
Let $Y _0 \hookrightarrow Y$ be a nilpotent closed immersion
of $S_i$-schemes of finite type. 
Then the functor
$X /Y\mapsto X \times _{Y} Y _0/Y _0$
from the category of $Y$-schemes which are 
locally of formal finite type over $S_i$ and formally étale over $Y$
to the category of 
$Y _0$-schemes which are 
locally of formal finite type over $S_i$ and formally étale over $Y _0$
is an equivalence of categories.
\end{thm}

\begin{proof}
Following \ref{EGAIV18.1.2-fullyfaithful}, the functor is fully faithful. 
By using the full faithfulness, 
the essential surjectivity is local and we conclude using \ref{EGAIV18.1.1}.
\end{proof}

\begin{cor}
\label{lifting-pbasis}
\begin{enumerate}[(a)]
\item Let $Y _0 \hookrightarrow Y$ be a nilpotent closed immersion
of $S_i$-schemes of finite type. 
If $X _0$ is a $Y _0$-scheme locally of formal finite type having finite $p$-basis,
then there exists a (unique up to isomorphisms) 
$Y$-scheme locally of formal finite type $X$ having finite $p$-basis
such that 
$X \times _{Y} Y _0 \riso X _0$.

\item 
\label{lifting-pbasis-p2}
Let $\Y$ be a $\V$-formal scheme of finite type, 
$Y _0 := \Y \times _{\Spf \V} \Spec (S_i)$. 
If $X _0$ is a $Y _0$-scheme locally of formal finite type having finite $p$-basis,
then there exists a (unique up to isomorphisms) 
formal $\Y$-scheme locally of formal finite type $\X$ having finite $p$-basis
and such that $\X \times _{\Y} Y _0 \riso X _0$.

\end{enumerate}
\end{cor}

\begin{proof}
1) Let us consider the first part. 
Let $f _0\colon X _0 \to \A ^d _{Y _0}$ be a relatively perfect morphism.
Then following \ref{EGAIV18.1.2} there exists 
a formally étale morphism
$f \colon X \to \bbA  ^d _Y$ such that 
$X$ is a $Y$-scheme locally of formal finite type and 
the reduction of $f$ modulo $\pi$
is $f _0$

2) Let us consider the second part. 
Let $f _0\colon X _0 \to \A ^d _{Y _0}$ be a relatively perfect morphism.
For any integer $n \geq i$, set
$Y _n := \Y \times _{\Spf \V} \Spec (S_n)$. 
Using \ref{EGAIV18.1.2} iteratively, 
there exists a family of
 formally étale morphism of the form 
$f _{n} \colon X _{n} \to \bbA  ^d _{Y _{n}}$ such that 
$X _{n}$ is a $Y _{n}$-scheme locally of formal finite type and such that 
the reduction of $f _{n}$ modulo $\pi ^{n}$
is $f _{n-1}$.
By taking inductive limits of the family $(f _n) _n$,
this yields a formally étale morphism
$f \colon \X \to \widehat{\A}  ^d _\Y$ such that 
$\X$ is a formal $\fS$-scheme locally of formal finite type and 
the reduction of $f$ modulo $\pi$
is $f _0$.
\end{proof}

\begin{thm}
\label{f0formétale-fforméta0}
Let $Y$ be  an $S_i$-scheme of finite type.
Let $f \colon X \to Y$ be an $S_i$-morphism locally of formal finite type. 
If $f$ is formally smooth then $f$ is flat.
\end{thm}

\begin{proof}
We can use the ideas of the proof of $a) \Rightarrow b)$ of \cite[17.5.1]{EGAIV4}
as follows.
Since this is local, we can suppose
$Y = \Spec ( A) $, and
there exist $n,n' \in \N$ and a closed $Y$-immersion of the form
$\iota\colon X \hookrightarrow 
\mathbb{D} ^n _{S _i} \times _{S _i} \A ^{n'} _{Y}$.
Let  
$B:= A \otimes _{V _i}  V _i  [[t _1, \dots, t _n  ]] \otimes _{V _i}  V _i [t _{n+1},\dots, t _{n+n'}] $
be the global section of the structural sheaf of 
$\mathbb{D} ^n _{S _i} \times _{S _i} \A ^{n'} _{Y}$.
Let $I $ be an ideal   of $B$ 
such that
$X = \Spec ( B/I)$. 
Set $C := B   /I $.
Since $C/A$ is formally smooth, following
\cite[0.20.5.14]{EGAIV1}, the sequence of $C $-modules
$$0 \to 
I  / I ^2 
\to 
\Omega _{B /A}
\otimes _{B}
C 
\to 
\Omega _{C  /A }
\to 
0$$
is split exact.

Let $\mathfrak{p}  $ be a prime ideal of $C $,
$\mathfrak{q}$ (resp. $\fr$)
be the induced prime ideal of $B$ (resp. $A$).
The $B $-module
$\Omega _{B  /A }$
is free with the basis 
$d t _1,\dots, d t _{n+n'}$ (see \ref{pbasis-Kato}).
We denote by $\partial _1, \dots, \partial _{n+n'}$ the corresponding derivations.
Hence, similarly to \cite[0.19.1.12]{EGAIV1}, we check there exist
$x _1, \dots, x _m \in I $ such that
the images of $x _1, \dots, x _m$ in 
$( I  / I ^2   ) _{\mathfrak{p} }$
generate 
$( I  / I ^2   ) _{\mathfrak{p} }$
and such that there exist
$n _1, \dots, n _m \in \{1,\dots,n+n' \}$
satisfying
$\det (\partial _{n _i} ( x _j))
\not \in \mathfrak{q} $.
Since $B $ is noetherian (see \ref{rem-Dnn'Ynoeth}.\ref{rem-Dnn'Ynoeth-i)}), $I $ is a $B $-module of finite type.
Hence, 
since 
$( I  / I ^2   ) _{\mathfrak{p} }=
 I  _{\mathfrak{q} }/ I  ^2 _{\mathfrak{q} }$,
 since
 $ I  _{\mathfrak{q} }\subset \mathfrak{q}  B _{ \mathfrak{q} }$
 then 
using Nakayama lemma, 
the images $x '_1, \dots, x ' _m$ of $x _1, \dots, x _m$ in 
$ I _{\mathfrak{q} }$
generate
$ I  _{\mathfrak{q} }$.

Let $y _1, \dots, y _m$ be the image of 
$x _1, \dots, x _m$ in the maximal ideal $\fm:= \fq B _{\fq} / \fr B _{\fq} $ of 
$B _{\fq} / \fr B _{\fq} $.
Then, $y _1, \dots, y _m$ are linearly independent modulo $\fm ^2$.
Indeed, otherwise there exist 
$w _1,\dots, w _m \in B$ with $w _{j _0} \not \in \fq$ for at least one $j _0$
such that
$\sum _{j=1} ^{m} w _j x _j \in \fq ^2$.
This yields 
$\sum _{j=1} ^{m} w _j \partial _{n _i} (x _j) \in \fq$, for any $i=1,\dots, m$.
Hence, 
we get a contradiction with 
$\det (\partial _{n _i} ( x _j))
\not \in \mathfrak{q} $. 

Since $k (\fr) =A _{\fr} / \fr A _{\fr}$ is a finite $V _i$-module, we get
$k (\fr)  \otimes _{V _i}  V _i  [[t _1, \dots, t _n  ]] 
=
k (\fr)   [[t _1, \dots, t _n  ]] $.
Hence,
$B _{\fr} / \fr B _{\fr} \riso
k (\fr)    [[t _1, \dots, t _n  ]] \otimes _{k (\fr)}  k (\fr) [t _{n+1},\dots, t _{n+n'}] 
=
( k (\fr)    [[t _1, \dots, t _n  ]] )  [t _{n+1},\dots, t _{n+n'}]  $
is a regular Noetherian ring. Since
$B _{\fq} / \fr B _{\fq} $ is a localisation of $B _{\fr} / \fr B _{\fr}$,
then $B _{\fq} / \fr B _{\fq} $ is also a
regular Noetherian ring. This yields that 
$y _1, \dots, y _m$ is a regular sequence of 
$B _{\fq} / \fr B _{\fq} $ (see \cite[0.17.1.7]{EGAIV4}). 
Using \cite[0.10.2.4]{EGAIII1} (to the homomorphism of Noetherian local rings
$A _{\fr} \to B _{\fq}$), 
this yields that 
$ B _{\fq}
\overset{x' _1}{\longrightarrow} 
 B _{\fq}$
 is injective
 and that 
$B _{\fq}/ x '_1 B _{\fq}$ is flat over $A _{\fr}$.
Using again inductively \cite[0.10.2.4]{EGAIII1} (to the homomorphism of Noetherian local rings
$A _{\fr} \to B _{\fq}/ \sum _{i=1} ^{j} x '_i B _{\fq} B _{\fq}$ for $j =1,\dots, m-1$),
we prove that 
$x '_1, \dots, x '_m$ is a regular sequence of 
$B _{\fq}$ and that 
$B _{\fq}/ \sum _{i=1} ^{m} x '_i B _{\fq}$ is flat over $A _{\fr}$.
Since 
$\sum _{i=1} ^{m} x '_i B _{\fq} = I _{\fq}$, this means that 
$C _{\fq}= B _{\fq}/ I _{\fq}$ is flat over $A _{\fr}$.
\end{proof}

\begin{rem}
Let 
$Y $ 
 be an  $S_i$-scheme.
Let  $f \colon X  \to \A ^d _{Y}$ be a relatively perfect morphism. 
\begin{enumerate}[(a)]
\item If $i =0$ and  $Y _0$ is  a noetherian regular $S_0$-scheme,
then  following a result of Gabber (see \cite[1.5]{Kato-explicity-recip91}),
$f $ is flat. 

\item 
If  $Y$ is an $S_i$-scheme of finite type, then since $f$ is formally étale 
$f$ is flat (see  \ref{f0formétale-fforméta0}). 
\end{enumerate}
\end{rem}

\begin{cor}
\label{f0formétale-fforméta0-cor}
Let $\fY$ be a formal $\fS$-scheme of finite type.
Let $f \colon \fX \to \fY$ a morphism locally of formal finite type 
of formal $\fS$-schemes. 
If $f$ is formally smooth then $f$ is flat.
\end{cor}

\begin{proof}
Since $\fX$ and $\fY$ are $p$-adic Noetherian formal schemes, 
by using \cite[Theorem 1 of III.5.2]{bourbaki} (and by using the Krull intersection theorem), 
this is a consequence of \ref{f0formétale-fforméta0}.
\end{proof}

\begin{ex}
The main example of Theorem \ref{f0formétale-fforméta0-cor} 
and 
Corollary \ref{f0formétale-fforméta0-cor}
is when $Y = S_i$ and $\fY =\fS$.
More precisely, 
our main objects will be 
$S _i$-schemes locally of formal finite type
and having locally finite $p$-bases over $S _i$
(resp. formal $\fS$-schemes locally of formal finite type
and having locally finite $p$-bases over $\fS$)
which are then $S _i$-flat (resp. $\fS$-flat).
\end{ex}

 \begin{prop}
 \label{f0formétale-fforméta}
Let $u\colon 
Y _0
 \hookrightarrow 
 Y $ 
 be a nilpotent closed immersion of $S_i$-schemes of finite type.
Let 
$X$ be  an $S_i$-scheme locally of formal finite type
and $f \colon X \to Y$ be  an $S_i$-morphism. 
Let $X _0 := X \times _{Y} Y _0$ and $f _0 \colon X _0 \to Y _0$ be the induced morphism. 
\begin{enumerate}[(a)]
\item If $f _0$ is formally smooth and $f$ is flat then $f$ is formally smooth.
\item If $f _0$ is formally étale and $f$ is flat then $f$ is formally étale.
\end{enumerate}

\end{prop}

\begin{proof}
Let us check the first assertion. 
Since this is local (see \cite[17.1.6]{EGAIV4}), we can suppose $Y _0 = \Spec ( A _0) $,
$Y = \Spec ( A) $, and
there exist $n,n' \in \N$ and a closed $S_i$-immersion of the form
$\iota\colon X \hookrightarrow \mathbb{D} ^n _{S_i} \times _{S_i} \A ^{n'} _{S_i}$.
Put 
$D  := (V _i  [[t _1, \dots, t _n  ]])[t _{n+1},\dots, t _{n+n'}] $
the global section of the structural sheaf of 
$\mathbb{D} ^n _{S_i} \times _{S_i} \A ^{n'} _{S_i}$.
Put
$B _0:= A _0 \otimes _{V _i} D $, 
$B:= A \otimes _{V _i} D$. Let $I$ be the ideal of $B$ such that 
$X = \Spec (B/I)$.
Set $I _0 = I B _0$, 
$C := B /I$,
$C _0 := B _0  /I _0$.
Hence, 
 $X _0 = \Spec ( C_0)$.
Let 
$\mathfrak{p} $ be a prime ideal of $C$. 
Since the closed immersion is nilpotent, it is associated with a prime ideal $\mathfrak{p} _0 $ of $C _0$
such that $k (\mathfrak{p} _0 ) =k (\mathfrak{p} )$.
 Since $A \to C$ is flat, then we get the exact sequence
 $ 0 \to  I \otimes _{A} A _0 \to B \otimes _{A} A _0
 \to C \otimes _{A} A _0 \to 0$.
Since 
$B _0  =B \otimes _{A} A _0$,
and since the canonical morphism 
$I \otimes _{A} A _0
\to 
I \otimes _{B} B _0$
 is an isomorphism, 
then this yields that 
the canonical epimorphism $ I \otimes _{B} B _0 \to I B _0=I _0$
is an isomorphism.
Since $C  \otimes _{B} B _0 \to C _0$ is an isomorphism, 
this yields
 $(I  / I ^2 )
 \otimes _{C} C _0
 \riso
 (I  / I ^2 )
 \otimes _{B} B _0
 \riso
 I _0 / I _{0} ^2 $.
This implies that the canonical morphism
 $(I  / I ^2 )
 \otimes _{C}k (\mathfrak{p} )
\to 
\Omega _{B /A}
\otimes _{B }
k (\mathfrak{p} )
$
 is isomorphic to 
 $I _0 / I _{0} ^2 
  \otimes _{C _0}k (\mathfrak{p} _0 )
\to 
\Omega _{B _0 /A _0}
\otimes _{B _0}
k (\mathfrak{p}  _0)$.
Since $C _0/A _0$ is formally smooth, 
then 
$I _0 / I _{0} ^2 
  \otimes _{C _0}k (\mathfrak{p} _0 )
\to 
\Omega _{B _0 /A _0}
\otimes _{B _0}
k (\mathfrak{p}  _0)$
is injective (see \cite[0.20.5.14]{EGAIV1}).
Hence, 
so is 
$(I  / I ^2 )
 \otimes _{C}k (\mathfrak{p} )
\to 
\Omega _{B /A}
\otimes _{B }
k (\mathfrak{p} )
$
and 
we conclude by using 
\cite[0.22.6.4]{EGAIV1}.

Let us check the second assertion.
From the first part, we already know that $f$ is formally smooth. 
Using \ref{EGA17.2.3}, this yields that 
$\Omega _{X/Y}$ is an $\O _X$-module locally free of finite type. 
Since $f _0$ is formally unramified, then using 
\cite[16.4.5 and 17.2.1]{EGAIV4}
we get 
$\Omega _{X/Y}\otimes _{\O _X} \O _{X _0}
\riso
\Omega _{X _0/Y _0}
=0$. Hence, $\Omega _{X/Y}=0$, i.e. $f$ is formally unramified.\end{proof}

\begin{prop}
\label{lem-flat-ff0}
\begin{enumerate}[(a)]
\item Let $u\colon 
Y _0
 \hookrightarrow 
 Y $ 
 be a nilpotent closed of $S_i$-schemes.
 Let $f \colon X \to X'$ be a morphism of noetherian $Y$-schemes.
Let $X _0 := X \times _{Y} Y _0$,
$X ' _0 := X ' \times _{Y} Y _0$ and $f _0 \colon X _0 \to X'  _0$ be the induced morphism. 
We suppose that $X/Y$ is flat. 
Then $f$ is flat if and only if $f _0$ is flat. 

\item Let $f \colon \X \to \Y$ be a morphism of noetherian formal $\fS$-schemes without $p$-torsion.
Then $f$ is flat if and only if $f _0$ is flat. 
\end{enumerate}
\end{prop}

\begin{proof}
If $f$ is flat, then $f _0$ is always flat. 
The converse of the first statement is a consequence of the flatness criterium by fiber (see \cite[Theorem 11.3.10]{EGAIV3}).
Let us prove the second one.
Since this is local, we can suppose $f$ corresponds to a morphism
$\phi \colon A \to B$ of $p$-adically complete $\V$-algebras without $p$-torsion
such that $A /\pi A \to B /\pi B$ is flat.
Since $A$ and $B$ have no $p$-torsion, 
then we get respectively the last and the last isomorphism
$B  \otimes  _{\V} k  \riso B \otimes ^\L _{\V} k 
\riso B \otimes _{A} ^\L (A/\pi A)$.
Using the Krull intersection theorem, the noetherianity of $A$ and $B$, the separated completeness of $B$, 
we check that for any ideal $\mathfrak{a}$ of $A$, 
$\mathfrak{a} \otimes _A B$ is a $B$-module of finite type and is Hausdorff for the $p$-adic topology. 
Hence, thanks to \cite[Theorem 1 of III.5.2]{bourbaki}, this implies that $A \to B$ is flat.
\end{proof}

\begin{empt}
\label{rem-dfn-pbasispadic}
Let $Y $  be  an $S_i$-scheme of finite type.
Let $g\colon X \to Y$ be a flat  $S_i$-morphism 
locally of formal finite type.

\begin{enumerate}[(i)]

\item \label{rem-dfn-pbasispadic2)}
Suppose given a  $Y _0$-morphism of the form $f _0\colon X _0 \to \A ^d _{Y _0}$ which is  formally étale. 
Choose a $Y$-morphism 
$f \colon X \to \A ^d _{Y}$ which is a lifting of $f _0$.
Then $f$ is formally étale and flat. 
Indeed, since $f _0$ is formally étale and locally of formal finite type, then following \ref{f0formétale-fforméta0},
$f _0$ is flat. 
Since 
$X/Y$ is flat, then by using \ref{lem-flat-ff0}
this yields that $f$ is flat.
Hence, $f$ is formally étale (see \ref{f0formétale-fforméta}).

\item 
 This yields that $X/Y$ has locally finite $p$-bases 
if and only if $X _0/Y _0$ has locally finite $p$-bases.

\end{enumerate}
\end{empt}

\begin{empt}
\label{rem-dfn-pbasispadicbis}
Let $\X$ be a formal $\fS$-scheme locally of formal finite type without $p$-torsion.

\begin{enumerate}[(i)]
\item
Suppose given a relatively perfect morphism of the form $f _0\colon X \to \A ^d _{k}$. 
Choose a lifting 
$f \colon \X \to \widehat{\A} ^d _{\V}$ of $f _0$.
Similarly, we check that such a lifting $f$ is flat and formally étale, 
e.g. the elements $t _1, \dots, t _d$ of $\Gamma (\X, \O _{\X})$ given by $f$ form a finite $p$-basis.

\item This yields that $\X$ has locally finite $p$-bases over $\fS$ if and only if 
$X$ has locally finite $p$-bases over $k$.
\end{enumerate}

\end{empt}

\subsection{Finite $p$-bases and closed immersions}

\begin{lem}
\label{lem-pbasisOmega}
Let 
$X=\Spec A$ be an affine $k$-scheme having 
a finite $p$-basis $b _1, \dots, b _d$ over $k$ (in the sense of Kato).
Let $t _1, \dots, t _d$ be elements of $A$. The following conditions are equivalent.

\begin{enumerate}[(a)]
\item $d {t _1}, \dots, d {t _d}$ generate $\Omega _{X/k}$; 

\item $t _1, \dots, t _d$ form a finite $p$-basis of $X$ over $k$.
\end{enumerate}

\end{lem}

\begin{proof}
1) Suppose $t _1, \dots, t _d$ form a finite $p$-basis of $X$ over $k$.
Then the morphism 
$X \to \A ^d _{k}$ given by $t _1, \dots, t _d$  is relatively perfect and then  formally étale.
We conclude 
by using \cite[17.2.4]{EGAIV4}.  

2) Suppose now $d {t _1}, \dots, d {t _d}$ generate $\Omega _{X/k}$.
By hypothesis, the fact that $b _1, \dots, b _d$ is a finite $p$-basis  over $k$ means two things: 
\begin{enumerate}[(a)]
\item $k ^{(p)} \otimes _{k} A \to A ^{(p)}$ is injective and 
\item  $b _1, \dots, b _n $ form a finite $p$-basis of $A/k$ in the sense of Grothendieck in \cite[0.21.1.9]{EGAIV1} (see \ref{aroundpbases}), i.e., 
$A$ is a free $k [A ^p]$-module with the basis
$\prod _{i=1} ^{d} b _i ^{n _i}$, with $n _i<p$ for any $i$.
\end{enumerate}
Since the statement $(a)$ is satisfied by hypothesis, 
then it is enough to check that 
$t _1, \dots, t _n $ form a finite $p$-basis of $A/k$ in the sense of Grothendieck.
Since 
$d {t _1}, \dots, d {t _d}$ generate $\Omega _{X/k}$, 
then following 
\cite[0.21.1.7]{EGAIV1}, the family
$B:= \{\prod _{i=1} ^{d} t _i ^{n _i} , \text{with $n _i<p$ for any $i$}\}$
generates the $k [A ^p]$-module $A$.
Since $A$ is a free $k [A ^p]$-module whose rank is equal to the cardinal of the family $B$, 
then this family $B$ is a basis (see \cite[Corollary 5 of II.3.3]{bourbaki}), 
which exactly means that
$t _1, \dots, t _n $ form a finite $p$-basis of $A/k$ in the sense of Grothendieck in \cite[0.21.1.9]{EGAIV1}.
\end{proof}

\begin{lemm}
\label{empt-closed-immer-local}
Let $u \colon Z \hookrightarrow X$ be a closed immersion of noetherian $p$-smooth $S_i$-schemes
(resp. be a closed immersion of noetherian $p$-smooth formal $\fS$-schemes).
Let $\I$ be the ideal defining $u$. 
Let $S:= S_i$ (resp. $S := \Spf \V$).

\begin{enumerate}[(a)]
\item The sequence
\begin{equation}
\label{EGAIV.0.20.5.14}
0
\to 
\I /\I ^2 
\to 
u ^* \Omega _{X/S}
\to 
 \Omega _{Z/S}
 \to 0
\end{equation}
is an exact and locally split sequence of locally free $\O _Z$-modules of finite type. 

\item Let $x \in Z$. 
There exists an open affine subset $ U _x \subset X$ and
sections $t _{1},\dots , t _{d} \in \Gamma ( U _x,\O _X)$
such that 
\begin{enumerate}[(a)]
\item $t _{r +1},\dots , t _{d} \in \Gamma ( U _x,\I)$
generate
$\Gamma ( U _x,\I)$, 
\item 
$d \overline{t} _{1},\dots , d \overline{t} _{r}$ is a basis of 
$ \Omega _{Z\cap U _x/S_i}$,
where
$\overline{t} _{1},\dots , \overline{t} _{r} \in\Gamma ( Z \cap U _x,\O _{Z})$
are the image of 
$t _{1},\dots , t _{r}$
via 
$\Gamma ( U _x,\O _{X})
\to
\Gamma ( Z \cap U _x,\O _{Z})$,
\item $d t _{1},\dots , d t _{d} $ is a basis of 
$\Omega _{U _x/S_i}$.
\end{enumerate}

\end{enumerate}

\end{lemm}

\begin{proof}
Since the respective case is a consequence of the non respective one, let us focus on this latter case. 
The fact that the sequence \ref{EGAIV.0.20.5.14} is exact and locally split
is a consequence of 
\cite[0.20.5.14]{EGAIV1} and the fact that 
$Z/S_i$ is formally smooth.
Since $X/S_i$ and $Z/S_i$ are $p$-smooth, 
$ \Omega _{X/S_i}$ and $ \Omega _{Z/S_i}$ are locally free of finite type. 
Hence so is $\I /\I ^2$.
We get from \ref{EGAIV.0.20.5.14} the split exact sequence 
\begin{equation}
\label{EGAIV.0.20.5.14x}
0
\to 
\I _x \otimes _{\O _{X,x}} k (x)
\to 
\Omega _{X/S_i, x} \otimes _{\O _{X,x}} k (x)
\to 
 \Omega _{Z/S_i,x} \otimes _{\O _{Z,x}} k (x)
 \to 0.
\end{equation}
There exists an open affine subset $ U _x \subset X$ such that 
the restriction to $Z \cap U _x$ of the terms of the sequence
\ref{EGAIV.0.20.5.14}
are free $\O _{Z \cap U _x}$-modules.
Let $t _{r +1},\dots , t _{d} \in \Gamma ( U _x,\I)$
which induces a basis of the restriction of $\I /\I ^2 $ over $Z \cap U _x$.
Since $\I _x$ is finitely generated $\O _{X,x}$-module, 
shrinking $U _x$ if necessarily,
by using Nakayamma's lemma we can suppose $t _{r +1},\dots , t _{d}$
generate
$\Gamma ( U _x,\I)$.
Moreover, 
since $Z/S _i$ is $p$-smooth (a fortiori around $x$), 
then shrinking again $U _x$ if necessarily, there exist $t _{1},\dots , t _{r} \in \Gamma ( U _x,\O _{X})$ such that
$d \overline{t} _{1},\dots , d \overline{t} _{r}$ is a basis of 
$ \Omega _{Z\cap U _x/S_i}$,
where
$\overline{t} _{1},\dots , \overline{t} _{r} \in\Gamma ( Z \cap U _x,\O _{Z})$
are the image of 
$t _{1},\dots , t _{r}$
via 
$\Gamma ( U _x,\O _{X})
\to
\Gamma ( Z \cap U _x,\O _{Z})$.
Using \ref{EGAIV.0.20.5.14x}, the image of $d t _{1},\dots , d t _{d} $ in 
$\Omega _{X/S_i, x} \otimes _{\O _{X,x}} k (x)$ is a basis of the $k(x)$-vector space
$\Omega _{X/S_i, x} \otimes _{\O _{X,x}} k (x)$.
Using 
\cite[Corollary of the Proposition 6 of II.3.2]{bourbaki},
the image of $d t _{1},\dots , d t _{d} $ in 
$\Omega _{X/S_i, x}$ is a basis of the free $\O _{X,x}$-module
$\Omega _{X/S_i, x}$.
Hence, 
shrinking $U _x$ if necessarily, 
we get that $d t _{1},\dots , d t _{d} $ is a basis of 
$\Omega _{U _x/S_i}$.
\end{proof}

\begin{prop}
\label{prop-pbasisOmega}
Let $X= \Spec A$ be an affine 
flat $S_i$-scheme 
(resp. $\X= \Spf A$ be an affine  
formal $\fS$-scheme without $p$-torsion)
locally of formal finite type
and having a finite $p$-basis $b _1,\dots, b _d$ over $S_i$ (resp. over $\fS$).
Let $t _1, \dots, t _d$ be elements of $A$,
and
$\overline{t} _1, \dots, \overline{t} _d$ be their image in  $A/\pi A$.
The following conditions are equivalent.
\begin{enumerate}[(a)]
\item $d {t _1}, \dots, d {t _d}$ generate $\Omega _{X/S_i}$ (resp. $\Omega _{\X/\V}$); 

\item $t _1, \dots, t _d$ form a finite $p$-basis of $X$ over $S_i$ (resp. $\X$ over $\fS$) ; 

\item $d \overline{t} _1, \dots, d \overline{t} _d$ generate $\Omega _{X _0/k}$ ;

\item $\overline{t} _1, \dots, \overline{t} _d$ form a finite $p$-basis of $X_0$ over $k$.

\end{enumerate}

\end{prop}

\begin{proof}
Following \ref{lem-pbasisOmega}, we have the equivalence between $1$ and $2$ 
and between $3$ and $4$. 
Using  
\ref{rem-dfn-pbasispadic} (resp. \ref{rem-dfn-pbasispadicbis}), 
we get the equivalence between $2$ and $4$.
\end{proof}

\begin{cor}
\label{cor-closed-immer-local}
Let $u \colon Z \hookrightarrow X$ be a closed immersion of 
flat $S_i$-schemes 
locally of formal finite type
and having locally finite $p$-bases  over $S_i$
(resp. be a closed immersion of formal $\fS$-schemes without $p$-torsion and 
having locally finite $p$-bases over $\fS$).
Let $\I$ be the ideal given by $u$. 
Let $S:= S_i$ (resp. $S := \Spf \V$).

Then, Zariski locally on $X$, there exist
$t _{r +1},\dots , t _{d}\in \Gamma (X,\I)$ generating 
$\Gamma (X,\I)$,
$t _{1},\dots , t _{r}  \in \Gamma ( X,\O _{X})$
such that, denoting by 
$\overline{t} _1, \dots, \overline{t} _d$ 
the image of 
$t _1, \dots, t _d$ on $\Gamma ( X,\O _{X})$,
 the following properties hold :
\begin{enumerate}[(a)]
\item $t _{1},\dots ,t  _{d}$ form a finite $p$-basis of $X$ over $S_i$ (resp. over $\fS$);
\item  $\overline{t} _{1},\dots ,\overline{t} _{r}$ form a finite $p$-basis of $Z$ over $S_i$  (resp. over $\fS$);
\item $\overline{t} _{r +1},\dots ,\overline{t} _{d}$ is a basis of $\I /\I ^2$.
\end{enumerate}

\end{cor}

\begin{proof}
This is a consequence of \ref{empt-closed-immer-local} and \ref{prop-pbasisOmega}.
\end{proof}

\section{Rigid geometry and adic spaces}
We will need later to work with the
direct image by the specialization morphism 
of the constant coefficient
when the boundary is not a divisor
(see \ref{directimagespovcv}).
The purpose will be to use
the local cohomology with support in a closed subscheme having locally finite $p$-bases of the constant coefficient
(see definition \ref{ntn-GammaZO-rig} and the remark \ref{rem-why-2dfn-loc-coh}), 
which use the local cohomology for adic spaces (see below the definition \ref{loc-coh&cech-reso}.\ref{loc-coh&cech-reso1}).
In order to make it clearer and fix the corresponding notations,
first we give below an overview (with further details as in \ref{ZR-adic-affinoid} when it is important in our work) of the notion of valued triples, analytic triples,
adic spaces and Zariski-Riemann spaces.

\subsection{Valuations, valuation rings and $I$-valuative local ring}

\begin{dfn}
[Valuations]
\label{ntn-Gamma-K/V}
We follow in this paper Huber convention of valuations (see \cite[2]{HuberContVal}) that we recall below.
Let $B$ be a commutative ring. 
Let $\Gamma$ be a totally ordered commutative group (whose law is written multiplicatively).
We get a totally ordered commutative group structure on the set 
$\Gamma 
\cup \{ 0\}$ (this is a disjoint union)
by decreeing that 
$0 < \gamma$ and $0 \cdot \gamma = \gamma \cdot 0 =0$ 
for all $\gamma \in\Gamma $
and $0 \cdot 0 = 0$. 
Remark that for any $\alpha , \beta \in \Gamma \cup \{0\}$, 
we have $\alpha \cdot \beta = 0$ if and only if 
$\alpha = 0$ or $\beta =0$.

A ``valuation of $B$ with value in $\Gamma \cup \{0\}$''  
is a mapping $v \colon B \to \Gamma \cup \{0\}$ 
which satisfies the following properties :
\begin{enumerate}[(a)]
\item 
$v  (x+y)
\leq \max
\{
v  (x),
v  (y) 
\} $,
for all $x,y \in B$ ;

\item $v  (xy)
=
v  (x)
v  (y)$,
for all $x,y \in B$ ;

\item $v (0) = 0$ and $v (1) = 1$.

\end{enumerate}
\end{dfn}

\begin{dfn}
[Equivalent valuations]
\label{valuation2valuationring}
Let $B$ be a commutative ring.

\begin{enumerate}[(a)]
\item 
Let $v \colon B \to \Gamma \cup \{0\}$ be a valuation of $B$.
The ``support of $v$'' is the prime ideal $\supp (v) := v ^{-1} (\{0\} $ of $B$.
Let $K _v$ be the fraction field of 
$B / \supp (v)$. 
We get the factorization 
$\overline{v}\colon K _v \to  \Gamma  \cup \{0\} $ of $v$ which is also a valuation.
The ``value group of $v$'' is denoted by $\Gamma _v$ and is by definition the subgroup of $\Gamma$ defined by
$\Gamma _v : = \overline{v} ( K _v ^\times)$.
Remark that $\Gamma _v$ is generated by
$v (B) \cap \Gamma$ and we have the mapping 
$v \colon B \to  \Gamma _v \cup \{0\}$.
Finally, we denote by 
$V _v:= \{ x\in K _v \; ;\; \overline{v} (x) \leq 1\}$ the ``valuation ring of $\overline{v}$''.

\item \label{valuation2valuationring-b}
Let $v$ and $w$ be two valuations of $B$.
We say that $v$ and $w$ are ``equivalent'' if 
the following equivalent
conditions are satisfied
\begin{enumerate}[(i)]
\item There is an isomorphism of ordered monoids 
$f\colon \Gamma _v \cup \{0\} \riso \Gamma _w \cup \{0\}$ with
$w =fo v$ ;
\item \label{valuation2valuationring-b-ii} $\supp (v) = \supp (w)$ and $V _v = V _w$ ; 
\item for all $a,b \in B$, 
$v(a)> v(b)$ if and only if $w(a)> w(b)$.
\end{enumerate}

\item The valuation spectrum
$\Spv B$ is the topological space defined as follows.

\begin{enumerate}[(a)]
\item 
$\Spv (B)$ is the set of 
equivalence classes of valuations of $B$.

\item The topology is the one generated by the subsets of the form
$$ \{v \in \Spv V \; ; \; v (x) \leq v (y) \not = 0 \}$$
for any $x,y \in B$.

\end{enumerate}

\item When $h \colon B \to C$ is an homomorphism of rings, 
then we get $\Spv (h ) \colon \Spv C \to \Spv B$ given by $v \mapsto v \circ h$.

\end{enumerate}

\end{dfn}

\begin{empt}
[Valuations associated to valuation rings, completions]
\label{valuation-affinoid}
Let $V$ be a  valuation ring with fraction field $\cK$
and maximal ideal $\fm _V$.

1) Set $\Gamma  := \cK ^\times  / V ^\times$. 
We endowed  canonically  $\Gamma $ with a totally ordered commutative group structure
(whose law is written multiplicatively) as follows: 
for any $a, b \in \cK ^\times$, 
\begin{equation}
\label{valuation-affinoid-dfn}
a \mod V  ^\times \leq b \mod V  ^\times
\text{ if and only if }
a \in b V 
\end{equation}
(see \cite[0.6.1.9]{FujiwaraKatoBookI} and use \cite[0.6.2.1]{FujiwaraKatoBookI}). 
The canonical map 
$v _{V }
\colon 
\cK
\to 
\Gamma 
\cup \{ 0\}$
defined by 
$v _{V } (a) = a \mod V  ^\times$ if $a \in \cK ^\times$,
and 
$v _{V } (0)=0$
is a valuation.

2) 
We suppose there exists $x \in \fm _V \setminus \{ 0\}$ (hence $V$ is not a field)
such that 
$V$ is $x$-adically separated.
Let $\widehat{V} $ be the $x$-adic completion of 
$V $,
let  $\widehat{\cK} $ be the fraction field of $\widehat{V} $, 
and $\widehat{\Gamma}  := \widehat{\cK} ^\times  / \widehat{V}  ^\times$. 
Since $V $ is a $x$-adically separated valuation ring, 
then $\widehat{V} $ is a valuation ring and we have the canonical isomorphism of totally ordered groups
$\Gamma  \riso \widehat{\Gamma} $ (see 
\cite[0.9.1.1]{FujiwaraKatoBookI}).

\end{empt}

\begin{empt}
[Valuation of the valuation ring associated to a valuation]
\label{valuation2valuationringbis}
Let $B$ be a commutative ring and 
$v \colon B \to \Gamma \cup \{0\}$ be a valuation of $B$.

\begin{enumerate}[(a)]
\item With notation \ref{valuation2valuationring}, 
we get 
$V _v ^\times = \{ x\in K _v \; ;\; \overline{v} (x) =1\}$.
Hence, the valuation $\overline{v}$ induces the isomorphism
of groups 
$K _v ^\times / V _v ^\times 
\riso \Gamma _v$, given by 
$a \mod V _v ^\times 
\mapsto 
\overline{v} (a)$.
This isomorphism preserves the order law.
Hence, $\overline{v}$  and $v _{V _v}$ are equivalent.

\item When $B$ is a field, then $K _v = B$, $v = \overline{v}$
and $v$ is equivalent to $v _{V _v}$ 
where $V _v :=  \{ x\in B \; ;\; v (x) \leq 1\}$. 
\end{enumerate}

\end{empt}

Later in the fundamental bijection of \ref{ZR-adic-affinoid}, 
we will need the following Lemma. 
\begin{lemm}
\label{lemm-VWdom-valequiv}
Let $V$ and $W$ be two valuations rings, 
let $f \colon V \hookrightarrow W$ be an injective map 
such that $W$ dominates $V$.
Let $\cK _V$ (resp. $\cK _W$)
be the fraction field of $V$ (resp. $W$).

\begin{enumerate}[(a)]
\item 
\label{lemm-VWdom-valequiv-1}
Then $f$ induces canonically the homomorphism 
of fields $g \colon \cK _V \to \cK _W$ and 
the injective homomorphism of 
totally ordered commutative groups
$\phi \colon 
\cK _V ^\times  / V ^\times
\to 
\cK _W ^\times  / W ^\times$
making commutative the diagram
\begin{equation}
\notag
\xymatrix{
{\cK _V} 
\ar[d] ^-{g}
\ar[r] ^-{v _V}
& 
{\cK _V ^\times  / V ^\times\cup \{ 0\}} 
\ar[d] ^-{\overline{\phi}}
\\ 
{\cK _W} 
\ar[r] ^-{v _W}
& 
{\cK _W ^\times  / W ^\times\cup \{ 0\},}
}
\end{equation}
where 
$\overline{\phi}$ is the extension of $\phi$ such that
$\overline{\phi} (0)=0$.

\item 
\label{lemm-VWdom-valequiv-2}
The valuations $v _W \circ g$ and $v _V$ are equivalent. 

\item 
\label{lemm-VWdom-valequiv-3}
Let $a \in \fm _V \setminus \{0\}$. 
Suppose $V$ and $W$ are $a$-adically separated. 
Let $\hat{f}\colon \widehat{V} \to \widehat{W}$ be the $a$-adic completion of $f$. 
Then $\hat{f}$ is injective, 
$\widehat{V}$ and 
$\widehat{W}$ are valuation rings such that
$\widehat{W}$
dominates
$\widehat{V}$.
\end{enumerate}
\end{lemm}

\begin{proof}
Since $V \subset g ^{-1}(W)$, then
$\phi$ (which sends 
$a \mod V  ^\times$ 
to 
$g (a) \mod W  ^\times$)
is a homomorphism of ordered commutative groups. 
In fact, we have $ g ^{-1}(W) = V$. 
(Indeed, reductio ad absurdum suppose
$x \in g ^{-1}(W) $ and $x \not \in V$.
Then $x ^{-1} \in V$. Hence $g (x ^{-1})  \in W$
and then $g (x ^{-1})  \in W ^\times$. Since $W$ dominates $V$,
this yields $x ^{-1} \in V ^\times$, which is absurd.)
This yields 
$g ^{-1}(W ^\times) = V ^\times$.
Hence, the homomorphism 
$\phi$ (and then $\overline{\phi}$)
is injective. 
The second statement is a straightforward consequence of the first one.
Let us consider the third one.
We have 
$V \cap g ^{-1}(a ^n W) = a ^n V$, for any integer $n\geq 0$.
(Indeed, $V \cap g ^{-1}(a ^n W) \supset a ^n V$ is obvious.
Moreover, let $x \in V \cap g ^{-1}(a ^n W)$.
Then 
$x /a ^n \in g ^{-1}(W) = V$.
Hence,
$x \in a ^n V$.)
Hence, the homomorphism 
$V / a ^n V \to W / a ^n W$ induced by $f$ is injective. 
This yields that $\widehat{f}$ is injective. The rest of the statement is a consequence of
\cite[0.9.1.1.(1), (2) and (5)]{FujiwaraKatoBookI}.
\end{proof}

\begin{dfn}
\label{dfn-I-adm}
Let $A$ be a commutative ring and $I$ be a  finitely generated (for simplicity) ideal of $A$.

\begin{enumerate}[(a)]
\item An ideal $J$ of $A$ is said to be ``$I$-admissible'' if $J$ is finitely generated 
and there exists an integer $n\geq 1$ such that $I ^n \subset J$. 

\item An ideal $J$ of $A$ is said to be ``invertible'' if 
there exists an $A$-submodule $K$ of $\mathrm{Frac} (A)$ the total field of fractions of $A$
such that $J \cdot K = \mathrm{Frac} (A)$.
We recall the following facts (see \cite[II.5.6]{bourbaki}) :
An invertible ideal is projective of rank $1$.
Conversely, if $J$ is non-degenerate (i.e. $J$ contains at least one non zero divisor) and $J _{\fm}$ is a principally generated ideal of 
$A _\fm$ for any maximal ideal $\fm$ of $A$, then $J$ is invertible.

\item We say that $A$ is ``$I$-valuative'' if any $I$-admissible ideal is invertible
(see \cite[0.8.7]{FujiwaraKatoBookI}).

\item We say that $A$ is an ``$I$-valuative local ring'' if 
$A$ is a local ring which is $I$-valuative. 
\end{enumerate}

\end{dfn}

\begin{rem}
\label{rem-Iadmloc}
Let $A$ be an $I$-valuative local ring, 
where $I$ is a finitely generated ideal. 
Since $I$ is finitely generated, then
$I$ itself is $I$-admissible and then invertible. 
Since $A$ is moreover local, this yields that $I$ is a free $A$-module of rank $1$, 
i.e. $I$ is generated by a non zero divisor element of $A$.
If $x$ is a generator of $I$, then 
the family of morphisms 
$\mathrm{Hom} _A ( I ^n , A)
\to 
A _x$
given by 
$\phi \mapsto \phi (x ^n) / x ^{n}$
 induces the isomorphism
$\underrightarrow{\lim} _{n\geq 0} 
\mathrm{Hom} _A ( I ^n , A)
\riso A _x$.
In particular, 
this yields that the canonical morphism
$A \to \mathrm{Hom} _A ( I ^n , A)$ is injective. 

\end{rem}

We recall below the following Theorem of 
\cite[0.8.7.8]{FujiwaraKatoBookI}.
\begin{thm}
\label{07878FK}
The statements below give the link between $a$-valuative local rings
and $a$-adically separated valuation rings. 
\begin{enumerate}[(a)]
\item \label{07878FK-1}
Let $A$ be an $I$-valuative local ring, 
where $I$ is a non zero proper finitely generated ideal. 
Let $a$ be a generator  of $I$. 
Set $J := \cap _{n\geq 1} I ^n$, 
$V:= A/J$ and $\overline{a}$ the image of $a$ in $V$.
Then 
\begin{enumerate}[(a)]
\item 
\label{07878FK-1a}
$B := \underrightarrow{\lim} _{n\geq 0} 
\mathrm{Hom} _A ( I ^n , A)$ is a local ring 
whose maximal ideal is equal to $J$ ;

\item $V$ is an $\overline{a}$-adically separated valuation ring
for $B/J$, the residue field of $B$ ;

\item 
\label{07878FK-1c}
$A = \{ f \in B \; ; \; \overline{f} \in V\}$, 
where if $b \in B$ then we denote by 
$\overline{b}$ the image of $b$ via the projection
$B \to B/J$.
\end{enumerate}

\item 
\label{07878FK-2}
Conversely, let $B$ be a local ring, 
$K$ be its residue field and
for any $b\in B$,  denote by 
$\overline{b}$ the image of $b$ via the projection
$B \to K$.
Let $a \in B ^\times$
and an $\overline{a}$-adically separated valuation ring $V$ with field of fraction $K$,
let $A := \{ f \in B \; ; \; \overline{f} \in V\}$ be the subring of $B$ defined as in \ref{07878FK-1c} above.
Then $A$ is an $a$-valuative local ring and 
$B = A _{a}$.
\end{enumerate}

\end{thm}

\begin{lem}
\label{valuation-v(A,I)}
Let $A$ be an $I$-valuative local ring, 
where $I$ is a non zero proper finitely generated ideal. 
We can associate canonically from $(A,I)$ a valuation
$v _{(A,I)}$ on 
$B := \underrightarrow{\lim} _{n\geq 0} 
\mathrm{Hom} _A ( I ^n , A)$
such that
$ A _{v _{(A,I)}}:= \{b \in B  \; | \; v _{(A,I)} (b) \leq 1 \} 
= A$,
$\supp v _{(A,I)} =  \cap _{n\in \N } \, I ^n$
and 
$ \{b \in B  \; | \; v _{(A,I)}  (b) < 1 \} = \fm _{A}$,
where $\fm _{A}$ is
the maximal ideal of $A$.
\end{lem}

\begin{proof}
Let $a$ be a generator  of $I$ (see \ref{rem-Iadmloc}) 
and then $B= A _a$.
Set  
$J  = \cap _{n\in \N } \, a ^n A$, 
$K  : =B  / J $, $V := A  /J $.
Following \ref{07878FK},
$B $ is a local ring whose  maximal ideal is equal to $J$,
$V $ is a valuation ring with field of fraction equal to $K$,
and we have the equality 
$A = \{ b \in B  \; | \; b \mod J   \in V \}$.
By setting 
$\Gamma  := K  ^\times / V  ^\times$, 
we get the valuation 
$v _V \colon K \to \Gamma  \cup \{0\}$ (see \ref{valuation-affinoid}).
This yields the valuation
$v _{(A,I)} \colon B  \to \Gamma  \cup \{0\}$
induced by composing $v _V$ with
the homomorphisms of rings
$B \to K $.
We have by construction
$\supp v _{(A,I)} = J$.
Morevoer,using \ref{07878FK}.\ref{07878FK-1c}, 
we get 
$A  = \{b \in B  \; | \; v _{(A,I)} (b) \leq 1 \} $.
Finally, for any $x \in A $, the property
$\overline{x} \in V ^\times$ is equivalent to
$x \in A   ^\times$.
Hence, 
$A  ^\times = \{b \in B  \; | \; v _{(A,I)} (b) =1 \} $.
This yields
$ \{b \in B  \; | \; v _{(A,I)}  (b) < 1 \} = \fm _{A}$.
\end{proof}

\subsection{Huber Adic spaces}
We recall some definitions.

\begin{dfn}
\label{fadic-dfn}
An ``f-adic ring'' is a topological ring 
$B$ that admits an open subring 
$A _0\subset  B$ 
such that the induced topology on 
$A _0$ is an adic topology defined by a finitely generated ideal $I _0$ 
of $A _0$.
In this situation, the subring $A _0$ 
is called a ``ring of definition'', and the ideal $I _0$ is called an ``ideal of definition'' of $B$ (or of $A _0$).

Let $B$ be an f-adic ring. 
A subring $A$ of $B$ that is open, 
integrally closed in $B$ and contained in $B ^{0}$ (the set of power bounded elements of $B$)
is called a ``ring of integral elements of $B$''. 
\end{dfn}

\begin{rem}
\label{rem-B=cupA:In}
Let 
$B$ be a ring, 
$A \subset B$ be a subring, 
and $I\subset A$ be a finitely generated ideal of $A$. 
Following \cite[0.B.1.1]{FujiwaraKatoBookI},
the ring $B$ 
endowed 
with the topological $A$-module structure given by the filtration 
$\{ I ^n\} _{n\geq 0} $ 
is an f-adic ring if and only if it is a topological ring if and only if the following equality holds
\begin{equation}
\label{B=cupA:In}
B=\cup _{n\geq 0} [ A : I ^{n}] .
\end{equation}

\end{rem}

\begin{dfn}
\label{fadic-dfnTate}
Let $B$ be an f-adic ring. 

\begin{enumerate}[(a)]

\item The f-adic ring $B$ is said to be a ``extremal'' 
if it has an ideal of definition $I _0$ such that
$I _0B= B$ (in that case, any ideal of definition $I$ satisfies $I B= B$).

\item The f-adic ring $B$ is said to be a ``Tate ring'' 
if there exists at least one unit of $B$ which is topologically nilpotent.

\item The f-adic ring  (resp. Tate ring) $B$ is said to be ``complete'' if it is separated and complete.
\end{enumerate}

\end{dfn}

\begin{lem}
\label{lem-Tate=ext+princ}
Let $B$ be an f-adic ring. 
Then $B$ is a Tate ring if and only if $B$ is extremal and has a principal ideal of definition.
In that case, the ideal generated by any unit of $B$ which is topologically nilpotent 
is an ideal of definition.
\end{lem}

\begin{proof}
Suppose $B$ is extremal and has a principal ideal of definition $I=(a)$.
Since $I$ is an ideal of definition and $a \in I$, then $a$ is topologically nilpotent. 
Since $B$ is extremal, then $a$ is a unit. 

Conversely, suppose there exists a unit $a$ of $B$ which is topologically nilpotent.
Let $A$ be a ring of definition of $B$ and $I$ be an ideal of definition of $A$.
Since $a$ is topologically nilpotent, there exists an integer $n _0$  such that
$a ^{n _0} \in I$. Hence, replacing $a$ by $a ^{n _0}$ if necessary, 
we can suppose $a \in I$.
Following \ref{B=cupA:In}, there exists a positive integer large enough $n$
such that $a ^{-1} I ^{n} \subset A$.
This implies 
$I ^{n} \subset a A$.
Hence the $a$-adic topology and the $I$-adic topology are identical.
\end{proof}

\begin{ex}
For instance, let $A$ be a flat $\cV$-algebra, 
$A _K: = A \otimes _{\cV} K  \riso A _p$.
Then $A _K$ can be endowed with  a Tate ring structure such that
$A$ is ring of definition  and $pA$ is an ideal of definition.
When $A$ is $p$-adically separated and complet,
then $A _K$ becomes a complete Tate ring.

\end{ex}

\begin{dfn}
An ``affinoid ring'' is a pair $A = ( A ^{\pm},A ^{+})$ consisting of an 
$f$-adic ring $A ^{\pm}$ and 
of a ring of integral elements $A ^+$ of $A ^{\pm} $.
An affinoid ring 
$A = ( A ^{\pm},A ^{+})$
is said to be extremal (resp. Tate)
if $A ^{\pm}$ is extremal (resp. Tate).
\end{dfn}

\begin{dfn}
\label{dfn-HuberSpa}
Let $A = ( A ^{\pm},A ^{+})$ be an affinoid ring. 
The associated adic spectrum
$\Spa A$ is the topological space defined as follows.

\begin{enumerate}[(a)]
\item As a set this is a subset of 
$\Spv (A ^{\pm})$ (see notation \ref{valuation2valuationring}). 
More precisely, the set $\Spa A$ consists of equivalence classes of valuations  
$v \colon A ^{\pm} \to \Gamma \cup \{ 0\}$ of $A ^{\pm}$
that satisfy $v (x) \leq 1$ for $x \in A ^{+}$ and are continuous. 
Here, the valuation $v$ is
``continuous'' means that for any
$\gamma \in \Gamma $,
there exists an open
neighborhood $U$ of $0$ in $A$ such that 
$v(a) < \gamma$
for every $a \in U$.

\item The topology is the one generated by the subsets of the form
$$ \{v \in \Spa A \; ; \; v (x) \leq v (y) \not = 0 \}$$
for any $x,y \in A ^{\pm}$.
\end{enumerate}
\end{dfn}

We have the following example of Tate affinoid ring.
\begin{lem}
\label{valuation-v(A,I)-cont}
Let $A$ be an $I$-valuative local ring, 
where $I$ is a non zero proper finitely generated ideal. 
Let $v _{(A,I)}$ be the valuation
 on 
$B := \underrightarrow{\lim} _{n\geq 0} 
\mathrm{Hom} _A ( I ^n , A)$ associated to $(A,I)$ (see \ref{valuation-v(A,I)}).

\begin{enumerate}[(a)]
\item $(B,A)$ is a Tate affinoid ring such that 
$A$ is a ring of definition and $I$ is an ideal of definition of $B$.

\item We have
$v _{(A,I)}
\in 
\Spa (B,A)$.

\end{enumerate}
\end{lem}

\begin{proof}
1) Let $a$ be a generator  of $I$ (see \ref{rem-Iadmloc}) 
and then $B= A _a$ (see \ref{rem-Iadmloc} and \ref{07878FK}.\ref{07878FK-1a}).
Since we have also 
$B=\cup _{n\geq 0} [ A : I ^{n}]$,
then $B$ can be endowed with a Tate ring structure 
such that $A$ is a ring of definition and $I$ is an ideal of definition 
$B$. 
Set  
$J  = \cap _{n\in \N } \, a ^n A$, 
$K  : =B  / J $, $V := A  /J $,
$\Gamma  := K  ^\times / V  ^\times$.
Since $V$ is a valuation ring then $V$ is integrally closed. 
Using \ref{07878FK}.\ref{07878FK-1c}, this yields that $A$ is integrally closed in $B$.
Hence $(B,A)$ is a Tate affinoid ring.

2) For any $b \in B$, 
we denote by 
$\overline{b}$ the image of $B$ in $K$.
Since for any $x \in A$, $\overline{x}\in V$, then
$v _{(A,I)} (x) \leq 1$. It remains to check that
$v _{(A,I)} \colon B  \to \Gamma  \cup \{0\}$ is continuous. 
Let $\gamma _0 \in \Gamma$. 
Choose $b _0 \in B  \setminus J$ such that
$v _{(A,I)} (b _0) = \gamma _0$ (indeed, $v _{(A,I)}$ is surjective). 
Since $V$ is $\overline{a}$-adically separated and since $\overline{b} _0 \not = 0$, 
there exists an integer $n$ large enough such that 
$\overline{b} _0 \not \in \overline{a} ^n V$, i.e. 
$v _{(A,I)} (b _0) > v _{(A,I)} (a ^n)$ (see \ref{valuation-affinoid-dfn}). 
Since the topology on $A$ is the $a$-adic topology
and $A$ is an open subring of $B$, then
we have checked there exists an open neighborhood $U=a ^n A$ of $0$ in $B$ such that 
$v _{(A,I)} (x) < v _{(A,I)} (b _0)$
for every $x \in U$.
\end{proof}

\begin{dfn}
\label{dfn-Spa}
Let $A = ( A ^{\pm},A ^{+})$ be an affinoid ring,
and $X:= \Spa A$. 
Let $f_0,\dots, f _n \in A ^{\pm}$ such that the ideal $(f _1,\dots, f_n)$ is open
(when $A$ is Tate, this is equivalent to saying $(f _1,\dots, f_n) =A ^{\pm}$).
\begin{enumerate}[(a)]
\item Huber defines the f-adic ring 
$A ^{\pm} (\tfrac{f_1,\dots, f _n}{f _0})$ as follows (see \cite[1]{Huber-gen-rig-an-var}). 
\begin{enumerate}[(i)]
\item As a ring, 
$A ^{\pm} (\tfrac{f_1,\dots, f _n}{f _0}) := A ^{\pm} [\tfrac{1}{f _0}]$.
\item 
$A ^{\pm} (\tfrac{f_1,\dots, f _n}{f _0})$
has the ring of definition 
$A _0 [ \tfrac{f_1}{f _0}, \dots, \tfrac{f _n}{f _0}]$
with the ideal of definition 
$I A_0 [ \tfrac{f_1}{f _0}, \dots, \tfrac{f _n}{f _0}]$,
where $A _0$ is a ring of definition of $A ^{\pm} $ 
with the ideal of definition $I _0$ of $A_0$ (this topology is independent of the choice of $A _0$ and $I _0$).
\end{enumerate}

\item 
Let $A ^{+} (\tfrac{f_1,\dots, f _n}{f _0})$ be
the integral closure of 
$A ^{+} [\tfrac{f_1,\dots, f _n}{f _0}]$ in 
$A ^{\pm} [\tfrac{1}{f _0}]$.
Then $A ^{+} (\tfrac{f_1,\dots, f _n}{f _0})$ is a ring of integral elements of 
$A ^{\pm} (\tfrac{f_1,\dots, f _n}{f _0})$.
We get the affinoid ring 
$A(\tfrac{f_1,\dots, f _n}{f _0}):=
(A ^{\pm} (\tfrac{f_1,\dots, f _n}{f _0}), A ^{+} (\tfrac{f_1,\dots, f _n}{f _0})) $.

\item We denote the completion of the affinoid ring
$A(\tfrac{f_1,\dots, f _n}{f _0})$
by
\begin{equation}
\label{ntn-Spa<>}
A<\tfrac{f_1,\dots, f _n}{f _0}>:=
(A ^{\pm} <\tfrac{f_1,\dots, f _n}{f _0}>, A ^{+} <\tfrac{f_1,\dots, f _n}{f _0}>) .
\end{equation}

\item We define the open subset of $X$ by setting
\begin{equation}
\label{ntnX(f_i/f0)}
X (\tfrac{f_1,\dots, f _n}{f _0})
:=  \{v \in X \; ; \; v (f _i) \leq v (f _0) \not = 0, 
\ 
i= 1,\dots, n\}.
\end{equation}

The open subsets of the form $X (\tfrac{f_1,\dots, f _n}{f _0})$
form a basis of the topology of $X$
and are said to be ``rational''.

\item Following \cite[Lemma 1.5.(ii)]{Huber-gen-rig-an-var}, we have the homeomorphism
\begin{equation}
\label{15(ii)Huber}
X (\tfrac{f_1,\dots, f _n}{f _0})
\riso 
\Spa (A  <\tfrac{f_1,\dots, f _n}{f _0}>) .
\end{equation}

\end{enumerate}

\end{dfn}

\begin{dfn}
[Affinoid adic space]
\label{AffAdSp}
Let $A = ( A ^{\pm},A ^{+})$ be an affinoid ring,
and $X:= \Spa A$. 

\begin{enumerate}[(a)]
\item We define the presheaf $\cO _X$ of complete topological rings on $X$ on the basis of rational open subsets of $X$ by setting
$$\Gamma (X (\tfrac{f_1,\dots, f _n}{f _0}),  \cO _{X} ):= A ^{\pm} <\tfrac{f_1,\dots, f _n}{f _0}>,$$
where $f_0,\dots, f _n \in A ^{\pm}$  are such that the ideal $(f _1,\dots, f_n)$ is open.

\item Let $x \in X$, i.e. 
it corresponds to a continous valuation
$v _{X,x} \colon A ^{\pm} \to  \Gamma _x \cup \{0\}$
such that $v _{X,x}  (a) \leq 1$ for all $a \in A ^+$.
For every rational subset $U$ of $X$  such that
$x\in U$, the valuation $v _{X,x}$
extends uniquely to a continuous valuation 
$v _{U,x}
\colon 
\Gamma (U,\cO _X) 
\to  \Gamma _x \cup \{0\}$.
Then the valuations $v _{U,x}$ define a valuation
$v _x 
\colon 
\cO  _{X,x}
\to  \Gamma _x \cup \{0\}$.
Following
\cite[Proposition 1.6 (i)]{Huber-gen-rig-an-var},
$\cO  _{X,x}$ 
(where $\cO _{X,x}$ denotes the inductive limit 
$\underset{x \in U}{\underrightarrow{\lim}} \Gamma (U, \cO _X)$ 
in the category of rings)
is a local ring whose maximal ideal is equal to the support
$\supp (v _x)$
of $v _x$. 

If $v$ and $w$ are two equivalents valuations then so are
$v _x$ and $w _x$. 
Hence we have the mapping 
$\Spa (A) \to \Spv (\cO  _{X,x})$ given by 
$v \mapsto v _x$.

\item When $A ^{\pm}$ has a noetherian ring of definition, then 
$\cO _X$ is a sheaf of complete topological rings on $X$ (see \cite[Theorem 2.2]{Huber-gen-rig-an-var}).

\item 
The ``affinoid adic space'' associated with $A$ 
is by definition 
$\Spa (A) := (X,  \cO _X, \{v _x \} _{x \in X})$.

Beware that from now 
$\Spa (A) $ will mean an object of ${\bf V}$ (see below \ref{dfn-catV})
and not only its underlying topological space. 

\end{enumerate}

\end{dfn}

\begin{dfn}
\label{dfn-catV}
Following \cite[2]{Huber-gen-rig-an-var},
we have the following definitions.

\begin{enumerate}[(a)]
\item First, we need the category
{\bf V} defined as follows.
The objects are the triples 
$X=
(X, \cO _X, \{v _x \} _{x \in X})$, 
where $X$ is a topological space, $\cO _X$ is a sheaf of
complete topological rings on $X$ 
and 
$v _x
\in \mathrm{Spv} (\cO _{X,x}) $.

The morphisms $X \to Y$ are the pairs
$(\phi,h)$, where 
$\phi \colon X \to Y$ is a continuous
mapping and 
$h\colon 
 \colon \cO _Y \to \phi _* \cO _X$
 is a morphism of sheaves of topological rings
such that, 
for every $ x \in X$, $v _{\phi (x)}$ is equivalent to  $v_x \circ h _x$
(i.e. 
$\mathrm{Spv} (h _x) ( v _x) = v _{\phi (x)}$).

\item An ``affinoid adic space'' (resp. ``analytic affinoid adic space'') 
is an object of 
{\bf V}  which is isomorphic to
the affinoid adic space associated with an affinoid ring (resp. Tate affinoid ring).

\item An ``adic space'' (resp. ``analytic adic space'') 
is an object
$X=(X, \cO _X, \{v _x \} _{x \in X})$ of
{\bf V}  which is locally an affinoid adic space, i.e., every
$x\in  X$ has an open neighbourhood $U \subset X$ such that 
$(U, \cO _X | U, \{v _x \} _{x \in U})$ 
is an affinoid adic space (resp. analytic affinoid adic space). 
A morphism 
$X\to Y$ between adic spaces 
(resp. analytic adic spaces) 
$X$, $Y$ is a
morphism in {\bf V}.
We denote by 
{\bf Ad} 
(resp. {\bf AnAd})
the category of adic spaces
(resp. analytic adic spaces). 

\end{enumerate}

\end{dfn}

\begin{rem}
\label{rem-an-rings=extremal}
Following the remark after \cite[A3.10]{FujiwaraKatoBookI},
in the definition of analytic adic spaces, we can replace Tate affinoid rings 
by extremal affinoid rings without changing the category {\bf AnAd}.
\end{rem}

\begin{empt}
\label{dfncOX+}
Let $X=(X, \cO _X, \{v _x \} _{x \in X})$ be an adic space. 

\begin{enumerate}[(a)]
\item For any open subset $U$ of $X$, for any $x \in U$,
for any $f \in \Gamma (U,  \cO _{X} ) $, we denote by 
$f _x \in \cO _{X,x}$ the image of $f$
via the canonical homomorphism $\Gamma (U,  \cO _{X} ) \to \cO _{X,x}$.

\item We define the subsheaf $\cO ^+ _X$ of rings on $X$  of $\cO _X$ 
by setting 
$$\Gamma (U,  \cO ^+ _{X} ):= \{f \in \cO _X (U) \; ; \; v _{x} (f _x) \leq 1
\text{ for any $x \in U$}\}$$
for any open subset $U$ of $X$.
We can check that the sheaf $\cO ^+ _X$ is an open subsheaf of $\cO _X$
(i.e. for any open subset $U$ of $X$, 
$\cO ^+ _X(U)$ is an open subset of 
$\cO _X(U)$).
Indeed, since this is local we can suppose that $X$ is 
the affinoid adic space associated with an affinoid ring. 
By using \cite[1.5.(ii),(iii) and 1.6.(iv)]{Huber-gen-rig-an-var}, 
if $U$ is a rational open subset of $X$, 
then 
$(\cO  _{X} (U) , \cO ^+ _{X} (U))$ is an affinoid paire and 
we have the isomorphism 
\begin{equation}
\label{U-isoSPa}
U\riso  
\Spa (\cO  _{X} (U) , \cO ^+ _{X} (U)).
\end{equation}
In general, since rational open subsets of $X$ form
a basis of the topology of $X$, then we are done.

\item 
\label{dfncOX+-3}
Let $x \in X$. 
Following
\cite[Proposition 1.6 (i)]{Huber-gen-rig-an-var},
we get that $\cO  _{X,x}$ 
is a local ring whose maximal ideal is equal to the support
$\supp (v _x)$ of $v _x$. 
Following \cite[1.6.(ii)]{Huber-gen-rig-an-var}, we get the equality
$\cO ^+ _{X,x} = \{f \in \cO _{X,x}  \; ; \; v _{x} (f _x) \leq 1\}$.
Moreover, $\cO ^+ _{X,x}$ is a local ring
with maximal ideal 
$\{f \in \cO _{X,x}  \; ; \; v _{x} (f _x) <1\}$.

\end{enumerate}

\end{empt}

\begin{lem}
\label{lem-suppTatering}
Let $B$ be a Tate ring an
$v \colon B \to \Gamma \cup \{0\}$ be a continuous valuation. 
Let $\varpi \in B$ be a unit which is topologically nilpotent.

\begin{enumerate}[(a)]
\item 
\label{lem-suppTatering-1}
Let $A _0$ be a ring of definition of $B$. Then 
$B= A _0 [\frac{1}{\varpi}]$.

\item Let $A$ be the ring of integral elements of $B$ defined by
$A := \{b \in B \; ;\; v (b)\leq 1\}$.
We have the equality
$$\supp (v) 
:= 
v ^{-1} (\{0\}) 
= 
\cap _{n\geq 0} \varpi ^n A.$$ 

\end{enumerate}

\end{lem}

\begin{proof}
1) Since $\varpi$ is topologically nilpotent, replacing $\varpi$ by a power of $\varpi$ if necessary, 
we can suppose $\varpi \in A _0$. 
Since $\varpi$ is a unit of $B$, we get 
the canonical injective homomorphism of rings
$( A _0) _\varpi \to B$.
Set $I := \varpi A _0$. 
Since $B=\cup _{n\geq 0} [ A _0: I ^{n}] $
(see \ref{B=cupA:In}), 
then the homomorphism $( A _0) _\varpi \to B$ is surjective.

2) 
Let $x \in B$ such that $v (x) \not =0$
Since $v$ is continuous, 
since $\varpi$ is topologically nilpotent, 
then there exists an integer large enough $n\geq 0$
such that $\varpi ^n \in  
\{b \in B \; ;\; v (b)< v(x)\}$.
This means
$v(\varpi ^n) <v(x)$. 
But since for any $y \in \varpi ^n A$, we have $v(y) \leq v(\varpi ^n)$, 
then we get
$x\not \in \varpi ^n A$.
Hence,
$x\not \in \cap _{n\geq 0} \varpi ^n A$.
Conversely, let $x \in B$ such that $v (x) =0$. For any integer $n\geq 0$, we have
$x = \varpi ^n ( \varpi ^{-n} x)$.
Since $v(\varpi ^{-n} x)=0$, then $x \in \varpi ^n A$. Hence, we are done. 
\end{proof}

\begin{empt}
\label{anadsp-valuative}
Let $( B,A )$ be a Tate affinoid ring,
$(X,  \cO _X, \{v _x \} _{x \in X})
:=
\Spa ( B,A )
\in 
{\bf AnAd}$.
Let $\cO ^+ _X$ be
the sheaf defined in \ref{dfncOX+}.
Let $A _0 \subset A$ be a ring of definition of $B$ and 
let $\varpi \in A _0$ 
be a unit of $B$ topologically nilpotent.
Let $I:= \varpi A _0$ be the associated ideal of definition of $B$.

a) From \ref{lem-suppTatering}.\ref{lem-suppTatering-1}, we get the first equality which implies the second one:
$B= \underrightarrow{\lim} _{n\geq 0} 
\mathrm{Hom} _{A _0} ( I ^n , A_0)
=
\underrightarrow{\lim} _{n\geq 0} 
\mathrm{Hom} _A ( I ^n , A)
\riso 
A [ \frac{1}{\varpi}]$.

b) 
Let $x \in X$. Set 
$A _x = \cO ^+ _{X,x} $ 
and 
$B _x = \cO  _{X,x} $. 
From a) (which is also used in the case where $X$ is replaced by a rational open subset),
we have
$B _x = A _x [\frac{1}{\varpi}]$.
Hence,  using the remark \ref{rem-B=cupA:In}, $B _x$ can be seen as a Tate ring
such that $A _x$ is a ring of definition
and $\varpi _x A _x$ is an ideal of definition, where $\varpi _x$ is the image of $\varpi$ 
via $A \to A _x$. 
Moreover, the valuation $v _x$ is continuous for this topology.
(Indeed, let 
$\gamma _x \in \Gamma _x$.
Since $v _{X,x} \colon B \to  \Gamma _x \cup \{0\}$ is continuous, 
then 
there exists an open
neighborhood $U$ of $0$ in $B$ such that 
$v _{X,x}(b) < \gamma$
for every $b \in U$.
Since $\varpi$ is topologically nilpotent,
for $n$ large enough,
$v _{X,x}(\varpi ^{n}) < \gamma_x$.
Hence, 
$v _{x}(\varpi _x ^{n}) < \gamma_x$.
This yields 
we have
$v _{x}(b _x) < \gamma _x$
for any $b _x \in \varpi _x ^{n} A _x$.)
Hence, it follows from \ref{dfncOX+}.\ref{dfncOX+-3} and \ref{lem-suppTatering}
that we have the equality
$ \cap _{n\geq 0} \varpi _x ^n \cO ^+ _{X,x} = 
\supp (v _x)$. 

c) 
Since $\supp (v _x)$ is equal to the maximal ideal of $B _x$,
then $K_x := B_x/\supp (v _x)$ is a field and we get the valuation 
$\overline{v} _x \colon K _x \to \Gamma _x \cup \{0\}$.
We denote by 
$V _x := \{s \in K _x \; ; \; \overline{v} _x (s) \leq 1\}$ its valuation ring. 
From part b), 
$V _x$ is $\overline{\varpi} _x$-adically separated, where
for any $b _x\in B_x$,  we denote by 
$\overline{b_x}$ the image of $b_x$ via the projection
$B_x \twoheadrightarrow K_x$.
We have
$A _x = \{ b _x \in B _x\; ; \; \overline{b _x} \in V _x\}$ 
(use \ref{dfncOX+}.\ref{dfncOX+-3}).
Moreover, since $K _x = V _{x} [ \frac{1}{\overline{\varpi} _x}]$,
then it follows from \ref{07878FK}.\ref{07878FK-2}, that $A _x$ is $\varpi_x$-valuative. 

d) Since the ring $A _x$ is $\varpi_x$-valuative, then we get the valuation 
$v _{(A _x, \varpi _x)}$ associated with $(A _x, \varpi _x A _x)$ (see \ref{valuation-v(A,I)}).
It follows from \ref{valuation-v(A,I)} (use also \ref{valuation2valuationring}.\ref{valuation2valuationring-b-ii}) 
that the valuations $v _{(A _x, \varpi _x)}$ and 
 $v _x $ are equivalent. 

\end{empt}

\begin{dfn}
\label{dfnTri}
We define the category of triples denoted by ${\bf Tri}$ as follows.
\begin{enumerate}[(a)]
\item A ``triple'' is a data $(X, \O _X ^+, \O _X)$
consisting of a topological space $X$ with a
topologically and locally ringed spaces of the form
$(X, \cO _X )$
(i.e. $(X, \cO _X )$ is both 
a topologically ringed spaces 
and 
a locally ringed spaces),
together with an injective morphism
$\iota 
\colon 
\cO _X ^+
\hookrightarrow
\cO _X$
of sheaves of rings which 
maps $\cO _X ^+$ onto an open subsheaf of $\cO _X$
and such that $(X,\cO _X ^+)$ is a locally ringed space.

\item A morphism of triples 
$(X, \O _X ^+, \O _X)
\to 
(Y, \O _Y ^+, \O _Y)$
is a morphism of topologically locally ringed spaces  
$$ (f, \phi) \colon 
(X, \O _X)
\to 
(Y, \O _Y)$$
(i.e. $(f,\phi)$ is both a morphism of topologically ringed spaces 
 and a morphism of locally ringed spaces)
which induces a morphism of locally ringed spaces of the form
$$ (f, \psi) \colon 
(X, \O _X ^+)
\to 
(Y, \O _Y ^+),$$
i.e. the homomorphism of sheaf of rings
$\phi \colon \cO _Y \to f _* \cO _X$ satisfies
$\phi ( \cO _Y ^+)
\subset
f _* \cO _X ^+$
and $(f,\psi) 
\colon 
(X, \O _X ^+)
\to 
(Y, \O _Y^+)$
is a morphism of locally ringed spaces, where
$\psi \colon \cO ^+ _Y \to f _* \cO ^+  _X$ 
is the restriction of $\phi$.
\end{enumerate}

\end{dfn}

\begin{rem}
In the definition of the category of triples of \cite[A.1.1]{FujiwaraKatoBookI},
$(X, \O _X ^+)$ is a topologically locally ringed space (hence their category seems bigger).
Since we would like to stick to Huber's adic spaces (see \ref{dfn-rho-AdTri}),
we avoid bothering with some topology on $\O _X ^+$.
\end{rem}

\begin{ntn}
We define here the following categories and functors.
\label{ntn-rho-AdTrietc}
Following 
\cite[2.1.(ii) and (iii)]{Huber-gen-rig-an-var}
 we get a fully faithful functor
\begin{equation}
\label{dfn-rho-AdTri}
\iota _\mathrm{Tri}
\colon {\bf Ad}
\to 
{\bf Tri},
\ 
\iota _\mathrm{ATri}
\colon 
{\bf AnAd}
\to 
{\bf Tri}.
\end{equation}
given by 
$(X, \cO _X, \{v _x \} _{x \in X})
\mapsto 
(X, \cO ^+ _X, \cO  _{X})$, where 
the sheaf $\cO ^+ _X$ is defined in \ref{dfncOX+}.
We denote by 
{\bf AdTri}
(resp. {\bf AnAdTri}) 
the essential image of the functor $\iota _\mathrm{Tri}$ (resp. 
$\iota _\mathrm{ATri}$).

\end{ntn}

\begin{empt}
We  construct canonically an inverse functor 
\begin{equation}
\label{eqcatAnAdTri}
\kappa _\mathrm{ATri}
\colon 
{\bf AnAdTri}
\cong
{\bf AnAd}
\end{equation}
of the equivalence of categories $\iota _\mathrm{ATri}$
as follows.

i) Let $(X, \cO ^+ _X, \cO  _{X}) 
\in {\bf AnAdTri}$.
Let $U$ be an open subset of $X$ containing $x$
which is an analytic affinoid adic space.
Let $\varpi \in \cO _X (U)$ be a unit which is topologically nilpotent.
Following \ref{anadsp-valuative}.c) and its notation,
$\cO ^{+} _{X,x}$ is $\varpi _x$-valuative, 
where $\varpi _x$ is the image of $\varpi$ 
via $\cO _X (U) \to \cO ^{+} _{X,x}$. 
We denote by 
$v _x := v _{(\cO ^{+} _{X,x}, \varpi _x \cO ^{+} _{X,x})}
\colon 
\cO _{X,x}
\to 
\Gamma _x \cup \{0\}$ the corresponding valuation (see \ref{valuation-v(A,I)}).
We set 
$\kappa _\mathrm{ATri} 
(X, \cO ^+ _X, \cO  _{X}) 
:= 
(X, \cO _X, \{v _x \} _{x \in X})
\in {\bf AnAd}$.
It follows from 
\ref{anadsp-valuative}.d)
that 
$\kappa _\mathrm{ATri}  \circ \iota _\mathrm{ATri} =id$.

ii) It remains to check that $\kappa _\mathrm{ATri}$
is functorial. 
Let $f \colon 
(X, \O _X ^+, \O _X)
\to 
(Y, \O _Y ^+, \O _Y)$
be a morphism of ${\bf AnAdTri}$.
We denote by 
$ (f, \phi) \colon 
(X, \O _X)
\to 
(Y, \O _Y)$ the underlying morphism of topologically locally ringed spaces  
which induces the morphism of locally ringed spaces 
$ (f, \psi) \colon 
(X, \O _X ^+)
\to 
(Y, \O _Y ^+)$.
Let $x \in X$ and $y := f (x)\in Y$.
We get the homomorphism of local rings
$\phi _x \colon 
\cO _{Y,y} 
\to \cO _{X,x}$
and
$\psi _x \colon 
\cO ^{+} _{Y,y} 
\to \cO ^{+} _{X,x}$.
Let $U _x$ (resp. $U _y$) be an open subset of $X$ (resp. $Y$) containing $x$
(resp. $y$)
which is an analytic affinoid adic space.
We may suppose that $f (U _x) \subset U _y$.
Let $\varpi _y $ be a unit of $\cO _Y (U _y)$ which is topologically nilpotent.
By continuity of $\cO _Y (U _y) \to \cO _X (U _x)$,
$\varpi _x:= \phi _x (\varpi _y)$ 
is a unit of $\cO _X (U _x)$ which is topologically nilpotent.
Hence, 
following \ref{anadsp-valuative}.b),
the maximal ideal of $\cO  _{X,x}$ (resp. $\cO  _{Y,y}$) is equal to 
$ \cap _{n\geq 0} \varpi _x ^n \cO ^+ _{X,x} = 
\supp (v _x)$
(resp. 
$ \cap _{n\geq 0} \varpi _y ^n \cO ^+ _{Y,y} = 
\supp (v _y)$). 
Let $K_x := \cO  _{X,x}/\supp (v _x)$ be the residue field, 
$\overline{v} _x \colon K _x \to \Gamma _x \cup \{0\}$
be the induced valuation 
and 
$V _x := \{s \in K _x \; ; \; \overline{v} _x (s) \leq 1\}$ be it valuation ring; 
and similarly by replacing $x$ or $X$ by $y$ or $Y$. 
We denote by 
$\overline{\psi _x}\colon V _y \to V _x$ the morphism induced by $\psi _x$
and 
by 
$\overline{\phi _x}\colon K _y \to K _x$ the homomorphism of residue fields induced by $\phi _x$.
In fact, $\overline{\phi _x} ( V _y ) \subset V _x$ and 
$\overline{\psi _x}$ is also the morphism induced by $\overline{\phi _x}$.
This yields the homomorphism of ordered groups
$\theta _x \colon 
K _y ^\times / V _y ^\times 
\cup \{0\}
\to
K _x ^\times / V _x ^\times 
\cup \{0\}$.
We get the commutative diagram
\begin{equation}
\label{eqcatAnAdTri-diag1}
\xymatrix{
{\cO  _{Y,y}} 
\ar[d] ^-{\phi _x}
\ar[r] ^-{\pi _y}
& 
{K _y} 
\ar[r] ^-{ \overline{v _y}}
\ar[d] ^-{\overline{\phi _x}}
&
{K _y ^\times / V _y ^\times 
\cup \{0\}}
\ar[d] ^-{\theta _x}
\\ 
{\cO  _{X,x}} 
\ar[r] ^-{\pi _x}
& 
{K _x} 
\ar[r] ^-{ \overline{v _x}}
&
{K _x ^\times / V _x ^\times 
\cup \{0\},}
}
\end{equation}
where 
$\pi _x$ and $\pi _y$ are the canonical projections.
Since $\psi _x$ is local, then so is 
$\overline{\psi _x}$.
By using 
\ref{lemm-VWdom-valequiv}.\ref{lemm-VWdom-valequiv-2}, 
this yields that 
$\overline{v _x}\circ \overline{\phi _x}$ and 
$\overline{v _y}$ are equivalent. 
Since by definition
$v _x= \overline{v _x} \circ \pi _x$ 
and
$v _y= \overline{v _y} \circ \pi _y$,
it follows from the commutative diagram \ref{eqcatAnAdTri-diag1}
that
$v _x \circ \phi _x$ and $v _y$ are equivalent.

\end{empt}

\subsection{Zariski-Riemann spaces}

\begin{dfn}
Let $\fX$ be a Noetherian formal $\fS$-scheme, 
and $\cI$ be an ideal sheaf of $\O _{\fX}$. 
Since $\fX$ is noetherian, 
then $\cI$ is an ``admissible ideal'' (see \cite[Definition 3.7.4]{FujiwaraKatoBookI})
if it satisfies the following conditions.
\begin{enumerate}[(a)]
\item Finiteness: $\cI$ is $\O _{\fX}$-coherent.
\item  Openness: $\cI$ contains a power of $(p)$.
\end{enumerate}
We denote by 
$\mathrm{AId} _{\fX}$ 
the set of all admissible ideals of $\fX$.
The blowing-up 
$\X' = \X _\cI \to \X$
with respect to the admissible ideal $\cI \in \mathrm{AId} _{\fX}$
 is the inductive limits of the family of morphisms
$X' _i = \mathrm{Proj} ( \oplus _{n\in \N} \cI _i ^n ) 
\to X _i$, 
where $\cI _i = \cI \otimes _{\O _{\X}} \O _{X _i}$.

\end{dfn}

\begin{empt}
[Topological visualization]
\label{empt-ZRT-pre}
Let $\X$ be a formal $\fS$-scheme locally of formal finite type.

\begin{enumerate}[(a)]
\item Following 
\cite[II.3.2.(a)]{FujiwaraKatoBookI},
taking projective limits
in the category of locally ringed spaces
we get the locally ringed spaces 
$$(\fX _K, \O _{\fX _K} ^{\mathrm{int}})
:=
\underset{\cI \in \mathrm{AId} _{\fX} }{\underleftarrow{\lim}} 
\X _{\cI}.$$  
The topological space $\fX _K$ is the ``Zariski-Riemann topological space''
associated to $\fX$.
This is also called the topological visualization of the rigid space $\fX ^{\mathrm{rid}}$ associated to $\fX$.
The sheaf 
$\O _{\fX _K} ^{\mathrm{int}}$ is its ``integral structure sheaf''.
The specialization map 
$\mathrm{sp} _{\fX}
\colon 
\fX _K
\to
\fX$
(or 
$\mathrm{sp} _{\fX _\cI}
\colon 
\fX _K
\to
\fX_\cI$)
is by definition the morphism of locally ringed spaces 
canonically induced by construction.
Since $\X$ is a $p$-adic formal scheme, 
then 
$p \O _{\fX _K} ^{\mathrm{int}}$ 
is an ideal of definition of 
$\O _{\fX _K} ^{\mathrm{int}}$
(see definition \cite[II.3.2.3]{FujiwaraKatoBookI}).
Hence, 
$p \O _{\fX _K} ^{\mathrm{int}}$ is an invertible ideal of $\O _{\fX _K} ^{\mathrm{int}}$
and then 
$\O _{\fX _K} ^{\mathrm{int}}$ is $p$-torsion free
(see \cite[II.3.2.5]{FujiwaraKatoBookI}).

\item Following \cite[II.3.2.9]{FujiwaraKatoBookI},
the sheaf of the local rings
$\O _{\fX _K}:= \O _{\fX _K} ^{\mathrm{int}} \otimes _{\cV} K$
is said to be the rigid structure sheaf.
\end{enumerate}

\end{empt}

\begin{empt}
[Some topological properties of $\fX _K$]
\label{coll-top-XK}
Let $\X$ be a formal $\fS$-scheme locally of formal finite type  (see definition \ref{dfn-fft}).
We collect below some basic topological properties of $\fX _K$.
\begin{enumerate}[(a)]
\item The topological space $\fX _K$ is coherent 
and sober (see \cite[II.3.1.2]{FujiwaraKatoBookI}).
\item 
The specialization map $\mathrm{sp} _{\fX}
\colon 
\fX _K
\to
\fX$
is quasi-compact, closed, surjective
(see \cite[II.3.1.2 and II.3.1.5]{FujiwaraKatoBookI}).

\item 
\label{qc-ZRqc} Let  $\fU$ be an open subscheme of $\fX$.
Following \cite[II.3.1.3.(2)]{FujiwaraKatoBookI},
since $\fU$ is quasi-compact,
then the canonical map ${\bf ZR} (\fU) \to {\bf ZR} (\fX) $
maps homeomorphically onto 
the quasi-compact open subset
$\sp ^{-1} (\fU)$ of ${\bf ZR} (\fX) $.

\item \label{openbasis-fXK}
For any quasi-compact open subset $U$ of $\fX _K$
there exist an admissible ideal $\cI \in \mathrm{AId} _{\fX}$
and a quasi-compact open subset 
$\fU _{\cI}$ of 
$\X _\cI$
such that 
$U = \sp _{\fX _{\cI}} ^{-1} ( \fU _{\cI})$
(see \cite[II.3.1.3.(1)]{FujiwaraKatoBookI}).
Hence, 
the collection of 
$\sp _{\fX _{\cI}} ^{-1} ( \fU _{\cI})$,
where $\cI$ runs through the admissible ideals of 
$\mathrm{AId} _{\fX}$
and $\fU _{\cI}$ runs through affine open formal subshemes of $\fX _{\cI}$,
forms an open basis of the topological space $\fX _K$.
Such an element of the form $ \fU _{\cI}$ is called an affinoid open subspace of $\fX ^\mathrm{rig}$.
Since $\fX$ is Noetherian, 
then following 
\cite[II.1.1.3]{FujiwaraKatoBookI}
we can restrict to affinoid open subspace of $\fX ^\mathrm{rig}$ of the form 
$\fU (\tfrac{f_1,\dots, f _n}{f _0})
:= \Spf A \{ T _1,\dots, T _n\} / ( f _0 T _1 -f _1,\dots,  f _0 T _n -f _n) ^{f _0 - \text{sat}}$,
where $\fU = \Spf A$ is an open affine of $\fX$, 
$A \{ T _1,\dots, T _n\}$ is the $p$-adic completion of the polynomial $A$-algebra with $n$ variables,
$(f _0, f _1, \dots, f _n)$ is an admissible ideal of $A$,
$( f _0 T _1 -f _1,\dots,  f _0 T _n -f _n) ^{f _0 - \text{sat}}$
is the $f_0$-saturation of the ideal of $\Spf A \{ T _1,\dots, T _n\}$
generated by 
$f _0 T _1 -f _1,\dots,  f _0 T _n -f _n$.

When $A$ is $p$-torsion free, then 
$( f _0 T _1 -f _1,\dots,  f _0 T _n -f _n) ^{f _0 - \text{sat}}
=
( f _0 T _1 -f _1,\dots,  f _0 T _n -f _n) ^{p - \text{sat}}$.
(Indeed, 
$A \{ T _1,\dots, T _n\} / ( f _0 T _1 -f _1,\dots,  f _0 T _n -f _n) ^{f _0 - \text{sat}}$
is $p$-torsion free (see \cite[II.1.1.3]{FujiwaraKatoBookI}), i.e. 
$(( f _0 T _1 -f _1,\dots,  f _0 T _n -f _n) ^{f _0 - \text{sat}} ) ^{p- \text{sat}}
=( f _0 T _1 -f _1,\dots,  f _0 T _n -f _n) ^{f _0 - \text{sat}}$.
Hence, 
$( f _0 T _1 -f _1,\dots,  f _0 T _n -f _n) ^{p - \text{sat}}
\subset
( f _0 T _1 -f _1,\dots,  f _0 T _n -f _n) ^{f _0 - \text{sat}}$.
Since $(f _0, f _1, \dots, f _n)$ is an admissible ideal of $A$, 
since 
$(\overline{f} _0)= (\overline{f} _0, \overline{f} _1, \dots, \overline{f} _n)$ 
where $\overline{f} _i$ means the image of $f _i$ via the projection 
$A \to A \{ T _1,\dots, T _n\} / ( f _0 T _1 -f _1,\dots,  f _0 T _n -f _n)$,
then there exists an integer $N \geq 0$ large enough such that
$p ^N \in (\overline{f} _0)$. This yields 
$( f _0 T _1 -f _1,\dots,  f _0 T _n -f _n) ^{f _0 - \text{sat}}
\subset
( f _0 T _1 -f _1,\dots,  f _0 T _n -f _n) ^{p - \text{sat}}$.)
In other words, we get
$$\fU (\tfrac{f_1,\dots, f _n}{f _0})
= \Spf ( \left (A \{ T _1,\dots, T _n\} / ( f _0 T _1 -f _1,\dots,  f _0 T _n -f _n) \right ) / \text{$p$-torsion}).$$

 \end{enumerate}
\end{empt}

\begin{empt}
[Local description and notation]
\label{ZR-x-ntn}
Let $\X$ be a formal $\fS$-scheme locally of formal finite type and $\V$-flat.
Let  $x \in \fX _K$. 

\begin{enumerate}[(a)]
\item 
The local ring 
$A _x := \O _{\fX _K, x} ^{\mathrm{int}}
\riso 
\underset{\cI \in \mathrm{AId} _{\fX} }{\underrightarrow{\lim}} 
\O _{\fX _\cI, \mathrm{sp} _{\fX _\cI} (x)} $ 
is $p$-valuative 
and $p$-adically henselian (see \cite[II.3.2.6]{FujiwaraKatoBookI}). 
The henselianity implies 
$I _x := p A _x$ is included in the Jacobson ideal of
$A _x$
(in particular the ideal $I _x$ is proper).
Set
$J _x = \cap _{n\in \N } I _x ^n$, 
$B _x := \O _{\fX _K, x} 
=
A _x \otimes _\V K
=
A _x [ \frac{1}{p}]
$.
Set  $K _x : =B _x / J _x$, $V _x := A _x /J _x$,
$\Gamma _x := K _x ^\times / V _x ^\times$.
Following \ref{07878FK},
since $A _x$ is a $I _x$-valuative local  ring then 
$B _x$ is a local ring whose  maximal ideal is equal to $J _x$,
$V _x$ is a $p$-adically separated valuation ring with field of fraction equal to $K _x$,
and we have the equality 
$A _x = \{ b \in B _x \; | \; b \mod J _x  \in V _x\}$.
Following \ref{valuation-v(A,I)}, we get the valuation
$v _{(A _x, I _x)} \colon B _x \to \Gamma _x \cup \{0\}$
such that
$A _x  = \{b \in B _x  \; | \; v _{(A _x,I _x)} (b) \leq 1 \} $,
$\supp v _{(A _x,I _x)} =  J _x$
and 
$ \{b \in B _x  \; | \; v _{(A _x,I _x)}  (b) < 1 \} = \fm _{A _x}$,
where $\fm _{A _x}$ is
the maximal ideal of $A _x$.
In particular $A _x \to V _x$ is a homomorphism of local rings.
Finally,
following \ref{valuation-v(A,I)-cont}, 
$(B _x , A _x)$ 
is a Tate affinoid ring such that 
$A _x$ is a ring of definition and $I _x$ is an ideal of definition of $B _x$.
Moreover, 
$v _{(A _x, I _x)} 
\in 
\Spa (A _x, I _x)$.

\item Since $V _x$ is $p$-adically separated and $p \in \fm _{V _x}\setminus \{0\}$, then
it has a unique height-one prime ideal ; explicitly 
$\sqrt{pV _x}$ 
(see \cite[0.6.7.3]{FujiwaraKatoBookI}).
We denote by $\fp _x : =\sqrt{pV _x}$ this ideal
and by $V _{[x]} := V _{x,\fp _x}$ 
the corresponding height-one valuation ring with fraction field 
$K _x$. 
This yields that the value group $\Gamma _{[x]}:=K _x ^\times / V _{[x]}  ^\times$
is isomorphic to a non-zero subgroup of the ordered 
multiplicative group $(\R _{>0}, \times)$ of real numbers greater than $0$.
We get the valuation of the form
$v _{V _{[x]} }
\colon 
K _x \to \R _{\geq 0}$,
uniquely determined by 
$v _{V _{[x]} } (p ) = \frac{1}{p}$
(indeed, this equality determines the choice of  the monomorphism
$\Gamma _{[x]} \hookrightarrow 
(\R _{>0}, \times)$ of ordered multiplicative groups).
By construction we have 
$V _{[x]} = \{z \in K _x \; | \; v _{V _{[x]} } (z) \leq 1 \} $.
We get the valuation
$ v _{[A _x, I _x]} \colon B _x \to\R _{\geq 0}$
induced by composing the valuation $v _{V _{[x]} }$ 
with the homomorphisms of rings
$B _x \to K _x$.

\end{enumerate}

\end{empt}

\begin{dfn}
[Rigid points]
\label{rigidpoints}
Let $\fX$ be a formal $\fS$-scheme locally of formal finite type and $\V$-flat.
We recall few things on rigid points (see \cite[II.3.3]{FujiwaraKatoBookI}).
\begin{enumerate}[(a)]
\item A rigid point of $\fX$
is a morphism of formal $\fS$-schemes of the form
$\alpha 
\colon 
\Spf (V) 
\to 
\fX$,
where 
$V$ is a $p$-adically complete valuation ring 
with $p \in \fm _V \setminus\{0\}$.

\item A rigid point of $(\fX _K, \O _{\fX _K} ^{\mathrm{int}})$
is a morphism of locally ringed spaces of the form
$\alpha 
\colon 
\Spf (V) 
\to 
(\fX _K, \O _{\fX _K} ^{\mathrm{int}})$,
where 
$V$ is a $p$-adically complete valuation ring 
with $p \in \fm _V \setminus\{0\}$.
Remark that since 
$\fX/\fS$ is adic
and since
$\Spf (V) /\fS$ is adic,
then
$\sp _{\X} \circ \alpha 
\colon 
\Spf (V) 
\to 
\fX$
is a rigid point of $\fX$.

\item \label{valuation-affinoid3}
Let 
$\beta  \colon \Spf (V) \to \X$ 
be a rigid point. 
Then there exist a unique rigid point
$\alpha 
\colon 
\Spf (V) 
\to 
(\fX _K, \O _{\fX _K} ^{\mathrm{int}})$
such that
$\beta = \sp _{\X} \circ \alpha$.
(Indeed, 
let 
$\fX '$ be an admissible blow-up
$\fX ' \to \fX$
of $\fX$.
Since $V$ is
a $p$-adically complete valuation ring, 
then following \cite[II.3.3.7]{FujiwaraKatoBookI} 
there exists a unique morphism 
$\Spf (V ) 
\to \fX'$
factorizing
$\beta$.)

Hence, the map $\alpha \mapsto \sp _{\X} \circ \alpha$
gives a bijection between rigid points of 
$\fX$ and rigid points of 
$(\fX _K, \O _{\fX _K} ^{\mathrm{int}})$.

\end{enumerate}

\end{dfn}

\begin{empt}
[Associated  rigid point of a point of $\fX _K$]
\label{assoc-rig-pt}
Let $\fX$ be a formal $\fS$-scheme locally of formal finite type and $\V$-flat.
Let  $x \in \fX _K$. We have the following facts (for more details, see \cite[II.3.3.6]{FujiwaraKatoBookI}). 
\begin{enumerate}[(i)]
\item 
\label{assoc-rig-pt-i}
With notation \ref{ZR-x-ntn}, 
let $\widehat{V} _x$ be the $p$-adic completion of 
$V _x$.
The composition of the canonical maps
$A _x \to V _x
\to 
\widehat{V} _x$
induces the inductive system of homomorphisms 
$\{ \O _{\X ' , \sp _{\fX '} (x)} 
\to 
\widehat{V} _x\}$,
where
$\fX '$ runs through all admissible blow-ups
$\fX ' \to \fX$
of $\fX$
and hence the projective system of adic morphisms 
$\{ \Spf (\widehat{V} _x) \to \fX'\}$.
Taking the projective limits in the category of local ringed spaces, we get the rigid point
$$\alpha _x 
\colon 
\Spf (\widehat{V} _x) 
\to 
(\fX _K, \O _{\fX _K} ^{\mathrm{int}}),$$
such that 
$\alpha _x ( \fm _{\widehat{V _x}}) 
=x$ 
and such that the induced map of stalks at $x$ is the canonical map 
$$A _x 
=\O _{\fX _K, x} ^{\mathrm{int}}
\to
\widehat{V _x}.$$
This morphism $\alpha _x$ is
the so called associated with $x$ rigid point of $(\fX _K, \O _{\fX _K} ^{\mathrm{int}})$.
We define the associated with $x$ rigid point of $\fX$ to be
$\beta _x := \sp _{\fX} \circ \alpha _x$.

\item 
\label{assoc-rig-pt-ii}
Conversely, 
let 
$\alpha 
\colon 
\Spf (V) 
\to 
(\fX _K, \O _{\fX _K} ^{\mathrm{int}})$
be a rigid point such that 
$\alpha (\fm _{V})
=
x$.
Then there
exists uniquely an injective homomorphism 
$j\colon \widehat{V _x} \hookrightarrow V$ such that 
$V$ dominates $\widehat{V _x}$
and $\alpha _x \circ \Spf (j) = \alpha$.

\item 
\label{assoc-rig-pt-iii}
We define on the set of rigid points 
$(\fX _K, \O _{\fX _K} ^{\mathrm{int}})$
an equivalence relation denoted by $\approx$ as follows. 
This is the equivalence relation generated by the relation $\sim$ defined as follows: 
for any rigid points 
$\alpha 
\colon 
\Spf (V) 
\to 
(\fX _K, \O _{\fX _K} ^{\mathrm{int}})$
and
$\beta 
\colon 
\Spf (W) 
\to 
(\fX _K, \O _{\fX _K} ^{\mathrm{int}})$,
we say that $\alpha \sim \beta$ if 
there exists an injective map $f \colon V \hookrightarrow W$
such that $W$ dominates $V$ and
$\alpha \circ \Spf f = \beta$.
By using (\ref{assoc-rig-pt-ii}), 
we can check that $\alpha \approx \beta$ if and only if 
there exists a rigid point 
$\alpha _0$ such that
$\alpha _0\sim \alpha$
and $\alpha _0\sim \beta$.
Moreover, 
we get a bijection between the points of $\fX _K$ and
the set of $\approx$-equivalence class of rigid points of $(\fX _K, \O _{\fX _K} ^{\mathrm{int}})$.
\end{enumerate}

\end{empt}

\begin{dfn}
\label{dfn-rig-Noeth}
We recall below the notion of ``coherent rigid formal spaces over $\fS$'' as in the book 
\cite{FujiwaraKatoBookI}. We will not really need this point of view  since we focus on adic spaces. 
But this will help the reader when he looks at 
the references appearing in \cite{FujiwaraKatoBookI} that we will use 
(for instance in the paragraph \ref{ntn-affinoid-ads} below).
\begin{enumerate}[(a)]
\item We denote by ${\bf AcFs ^*}/\cV$,
the category of coherent (i.e. quasi-compact and quasi-separated) 
formal $\Spf \cV$-schemes
(see \cite[Definitions II.2.1.1 and II.2.1.13]{FujiwaraKatoBookI}).

\item We denote by ${\bf CRf} _\V$ the category of coherent rigid (formal) spaces over $\fS$, i.e. 
the localisation of the category ${\bf AcFs ^*}/\cV$ by admissible blowing-ups 
(see \cite[Definition II.2.1.13 and Theorem II.2.1.14]{FujiwaraKatoBookI}).
We get the functor 
${} ^\mathrm{rig} 
\colon 
{\bf AcFs ^*}/\cV
\to 
{\bf CRf} _\V$, i.e.
if $\fX \in
{\bf AcFs ^*}/\cV$
then 
$\fX ^{\mathrm{rig}}$ is the corresponding rigid space. 
\end{enumerate}

\end{dfn}

\begin{empt}
[Analytic affinoid adic space associated to an affine formal scheme locally of formal finite type]
\label{ntn-affinoid-ads}
Let $\fX=\Spf A$ be a formal $\fS$-scheme locally of formal finite type and $\V$-flat.
We make $A _K$ a complete Tate ring with ring of definition $A$ and ideal of definition 
$pA$ (see \ref{fadic-dfnTate}).
Let $A _K ^0$ be the set of bounded power elements of $A _K$ for this topology induced by $A$.

\begin{enumerate}[(a)]
\item 
Since $\fX$ is Noetherian, then then $\fX ^\mathrm{rig}$ is a Stein affinoid of type $(N)$ 
(see definition \cite[II.2.5.2]{FujiwaraKatoBookI})
and we can apply Proposition  \cite[II.6.4.1]{FujiwaraKatoBookI} (to understand this proposition, see also 
Notation \cite[II.6.3.2]{FujiwaraKatoBookI}).
Hence, we get the last equality
$$A _K
:=
\Gamma ( \fX, \O _{\fX} ) \otimes _\V K 
= 
\Gamma ( \fX _K, \O _{\fX_K} ).$$ 

\item Moreover, since $A$ is $\V$-flat then we can apply
Proposition  \cite[II.6.2.6]{FujiwaraKatoBookI} and we get that
$\Gamma ( \fX _K, \O _{\fX _K} ^{\mathrm{int}})$  is canonically isomorphic to
the integral closure of $A$ in $A _K$. 

\item Since $\fX ^\mathrm{rig}$ is a Stein affinoid of type $(N)$, then 
from \cite[A.4.10]{FujiwaraKatoBookI},
we have
\begin{equation}
\label{A.4.10-FK}
\Gamma ( \fX _K, \O _{\fX _K} ^{\mathrm{int}})= A _K ^0.
\end{equation}

\item We get the affinoid ring $(A _K, A _K ^{0})$ in the sense of Huber (see \ref{dfn-HuberSpa}).
This yields the affinoid adic space
$
\Spa (A _K, A _K ^{0})
\in 
{\bf AnAd}$.

\end{enumerate}

\end{empt}

\begin{empt}
\label{ZR-adic-affinoid}
Let $\fX=\Spf A$ be a formal $\fS$-scheme locally of formal finite type and $\V$-flat.
Then the Zariski-Riemann space $\fX _K$ is homeomorphic to the topological space
underlying the analytic adic space
$\Spa (A _K, A _K ^{0})$.
Since we think it is important to understand this identification in our work,
we give some details below  concerning easy checks which was left to the reader in 
\cite[II.A.7]{FujiwaraKatoBookI}.
\begin{enumerate}[(a)]
\item 
Take 
$x \in \fX _K$. 
It follows from the equality  \ref{A.4.10-FK}
that the canonical map 
$A _K \to B _x$ sends an element of 
$A _K ^{0} $ to an element of $A _x$. 
Recall 
$(B _x, A _x)$ is an affinoid ring such that
$A _x$ is a ring of definition and $I _x = p A _x$ is an ideal of definition of $B_x$.
This yields the continuous morphism of affinoid rings
$(A _K, A _K ^0) \to (B _x, A _x)$.
Since $v _{(A _x, I _x)}$ is continuous (see \ref{valuation-v(A,I)-cont}), 
then by composing
$v _{(A _x, I _x)}$
with $(A _K, A _K ^0) \to (B _x, A _x)$ we get
$v _x \in \Spa (A _K, A _K ^0)$.
This yields the canonical map
$\phi \colon 
\fX _K
\to
\Spa (A _K, A _K ^{0})$, 
given by $x \mapsto v _x$.

\item 
\label{ZR-adic-affinoid-a}
We construct the  canonical inverse map
$\Spa (A _K, A _K ^{0})
\to 
\fX _K$ as follows.

Let $v \in \Spa (A _K, A _K ^{0})$, i.e. let
$v \colon A _K \to \Gamma \cup \{0\} $
be a continuous valuation such that $v (a) \leq 1$ for any $a \in A _K ^{0}$. 
The support of $v$ is by definition the prime ideal
$\supp (v) := v ^{-1} (\{0\})$
of $A _K$.
Let $K _v$ be the fraction field of 
$A _K / \supp (v)$.
We get the factorization 
$\overline{v}\colon K _v \to  \Gamma \cup \{0\} $ of $v$ which is also a valuation.
Let $V _v:= \{ x\in K _v \; ;\; \overline{v} (x) \leq 1\}$ be the valuation ring of $\overline{v}$.
Since $v$ is continuous,
then $V _v$ is $p$-adically separated.
(Indeed, 
let $x  _0 \in V _v \setminus\{0\}$. Then 
$\overline{v} (x_0) \not = 0$.
Since $\{ b\in A _K ~;~v (b) < \overline{v} (x_0)\}$ is open then 
$v (p ^n) <  \overline{v} (x_0)$ for some integer $n \geq 1$.
Since for any $x \in p^n V _v$, we have
$ \overline{v} (x) \leq \overline{v} (p ^n) $, 
this yields that 
$x _0 \not \in \cap _{n\geq 0} p^n V _v $.)
Following \ref{valuation-affinoid}.2),
the $p$-adic completion $\widehat{V} _v$ of 
$V _v$,
is  a 
$p$-adically complete valuation ring having the same value group as 
$V _v$.

Via the composite of the homomorphisms of rings
$A \to A ^0 _K \to V _v \to \widehat{V} _v$,
we get the morphism 
$\beta _v \colon \Spf (\widehat{V} _v) 
\to \X$ of formal schemes.
Following \ref{rigidpoints}.\ref{valuation-affinoid3}),
there exists a unique rigid point
$\alpha _v 
\colon 
\Spf (\widehat{V} _v) 
\to 
(\fX _K, \O _{\fX _K} ^{\mathrm{int}})$
factorizing
$\beta _v $. 
The image of the closed point is the desired element $x _v$ of $\fX _K$.
Hence, 
we have constructed the canonical map
$\psi \colon 
\Spa (A _K, A _K ^{0})
\to 
\fX _K$, given by $v \mapsto x _v$.

\item We check in this step that these maps $\phi$ and $\psi$ are inverse to each other.

i) Let $x \in \fX _K$. First we check the equality $x _{v _x}=x$. 
Let us denote by 
$\theta _x \colon A _K \to B _x \twoheadrightarrow K _x$
the composition morphism 
(hence, we have
$v _{V _x} \circ \theta _x = v _x \colon 
A _K \to \Gamma _x \cup \{0\}$).
Since 
$\theta _x ^{-1} (\{0\})
=
v _x ^{-1} (\{0\})$,
then 
there exists a unique homomorphism of fields $g _x \colon K _{v _x} \to K _x$
making commutative the left square below 
\begin{equation}
\label{ZR-adic-affinoid-diag1}
\xymatrix{
{A _K} 
\ar[r] ^-{\theta _x}
\ar[d] ^-{}
& 
{K _x} 
\ar[d] ^-{v _{V _x}}
\\ 
{K _{v _x}} 
\ar@{.>}[ur] ^-{g _x}
\ar[r] ^-{\overline{v_x}}
& 
{\Gamma _x\cup \{0\},} 
}
\xymatrix{
{A} 
\ar[r] ^-{}
\ar[d] ^-{}
& 
{V _x} 
\ar[d] ^-{v _{V _x}}
\\ 
{V _{v _x}} 
\ar@{.>}[ur] ^-{f _x}
\ar[r] ^-{\overline{v_x}}
& 
{\Gamma _x\cup \{0\}.} 
}
\end{equation}
This yields
$ g _x ^{-1} (V _x) = V _{v _x}$
and
$ g _x ^{-1} (V ^\times _x) = V ^\times _{v _x}$.
Hence, $g _x$ induces the homomorphism of local rings $f _x \colon V _{v _x}\to V _x$ 
making commutative the right square of \ref{ZR-adic-affinoid-diag1}.
Following
\ref{lemm-VWdom-valequiv}.\ref{lemm-VWdom-valequiv-3},
this implies that $\widehat{V _{x}}$ dominates $\widehat{V _{v _x}}$
(via the $p$-adic completion of $f _x$).
We get the commutative left square below
\begin{equation}
\label{ZR-adic-affinoid-diag2}
\xymatrix{
{ \Spf (\widehat{V} _{v _x}) } 
\ar[r] ^-{\beta _{v _x}}
& {\fX } 
\ar@{=}[d] ^-{}
\\ 
{ \Spf (\widehat{V} _{x}) } 
\ar[r] ^-{\beta _{x}}
\ar[u] ^-{}
& {\fX ,} 
}
\
\xymatrix{
{ \Spf (\widehat{V} _{v _x}) } 
\ar[r] ^-{\alpha _{v _x}}
& {(\fX _K, \O _{\fX _K} ^{\mathrm{int}})} 
\ar@{=}[d] ^-{}
\\ 
{ \Spf (\widehat{V} _{x}) } 
\ar[r] ^-{\alpha _{x}}
\ar[u] ^-{}
& 
{(\fX _K, \O _{\fX _K} ^{\mathrm{int}}),} 
}
\end{equation}
where $\beta _x$ is the associated with $x$ rigid point of $\fX$ (see \ref{assoc-rig-pt}.\ref{assoc-rig-pt-i}),
$\beta _{v _x}$ is the rigid point associated with the valuation $v _x$ as defined in (b)
and the vertical arrow is induced by 
$\widehat{V _{v _x}}\to \widehat{V _x}$.
Following \ref{rigidpoints}.\ref{valuation-affinoid3}),
this yields the commutativity of the right square of \ref{ZR-adic-affinoid-diag2}. 
Since 
$\widehat{V _{x}}$ dominates $\widehat{V _{v _x}}$,
then 
$x _{v _x}=x$.

ii) Let $v, w \in \Spa ( A_K, A _K ^0)$
such that $x _v = x _w$.
We  check  that the valuation $v$ and $w$ are equivalent as follows. 
Using \ref{assoc-rig-pt}.\ref{assoc-rig-pt-iii}
we can suppose 
$\alpha _v \sim \alpha _w$.
Hence, we suppose there exists 
an injective homomorphism of local rings
$f\colon \widehat{V} _{v}\to \widehat{V} _w$
making commutative the left diagram below :
\begin{equation}
\notag
\xymatrix{
{A} 
\ar@{=}[d] ^-{}
\ar[r] ^-{\beta _v}
& 
{\widehat{V} _v } 
\ar[d] ^-{f}
\\ 
{A} 
\ar[r] ^-{\beta _w}
& 
{\widehat{V} _w ,} 
}
\ 
\xymatrix{
{A _K} 
\ar@{=}[d] ^-{}
\ar[r] ^-{}
&
{\widehat{K} _v} 
\ar[d] ^-{g}
\ar[r] ^-{v _{\widehat{V} _v}}
& 
{\widehat{K} _v ^\times  / \widehat{V} _v ^\times\cup \{ 0\}} 
\ar[d] ^-{\overline{\phi}}
\\ 
{A _K} 
\ar[r] ^-{}
&
{\widehat{K} _w} 
\ar[r] ^-{v _{\widehat{V} _w}}
& 
{\widehat{K} _w ^\times  / \widehat{V}  _w ^\times\cup \{ 0\}.}
}
\end{equation}
Denoting by $\widehat{K} _v $ (resp. $\widehat{K} _w $)
the fraction field of 
$\widehat{V} _v $ (resp. $\widehat{V} _w $)
and by $g \colon \widehat{K} _v \to \widehat{K} _w$ the morphism induced by $f$, we get the commutativity of 
the left square of the right diagram. 
Using  \ref{lemm-VWdom-valequiv}.\ref{lemm-VWdom-valequiv-1},
we get the commutative right diagram.
Since $\overline{\phi}$ is injective, 
since the composition of the top arrows (resp. bottom arrows)
is $v$ (resp. $w$) of the right diagram,
then $v$ and $w$ are equivalent.

iii) Using i) and ii), we get that $x \mapsto v _x$ and $v \mapsto x _v$ are inverse to each other. 

\item The mapping $\phi$ and $\psi $ are homeomorphisms.
Indeed, let 
$f_0,\dots, f _n \in A _K$ such that  $(f _1,\dots, f_n) =A _K$.
We get the rational subset
$U ^{\ad}:=\Spa (A _K, A _K ^{0}) (\tfrac{f_1,\dots, f _n}{f _0})
\riso 
\Spa\left ( (A _K, A _K ^{0}) <\tfrac{f_1,\dots, f _n}{f _0}> \right ) $
(see \ref{15(ii)Huber}).
Recall rational open subspaces of $\Spa (A _K, A _K ^{0})$ form an open basis.
For any integer $N$ large enough, we get
$p ^N f _i \in A$ for any $i = 0,\dots, n$.
With notation 
\ref{coll-top-XK}.\ref{openbasis-fXK},
we have the affinoid open subspace of $\fX ^\mathrm{rig}$ of the form 
$\fX (\frac{p ^N f_1,\dots, p ^N f _n}{p ^N f _0})
:= \Spf A \{ T _1,\dots, T _n\} / ( p ^Nf _0 T _1 -p ^Nf _1,\dots,  p ^Nf _0 T _n -p ^Nf _n) ^{p - \text{sat}}$, which is an open formal subscheme of $\fX _{\cI}$,
where 
$\cI = (p ^N f_1,\dots, p ^N f _n)$.
Using the open basis of $\fX _K$ explained at
\ref{coll-top-XK}.\ref{openbasis-fXK},
it is sufficient to check the equality
\begin{equation}
\label{equ-openbases}
\phi ^{-1} (U ^{\ad}) 
=
\sp _{\fX _{\cI}}  ^{-1}( \fX (\frac{p ^N f_1,\dots, p ^N f _n}{p ^N f _0})).
\end{equation}
Let us check this equality : 
let $x \in \fX _K$. We get  $\phi (x)= v _x \in \Spa (A _K, A _K ^{0})$.
Let $\beta _x \colon \Spf (\widehat{V _x}) \to
\Spf (A) $ be the induced  rigid point associated to $x$. 
Then $x$ belong to the right term of \ref{equ-openbases} if and only if 
the morphism $\beta _x$ can be factorized through
$\fX (\frac{p ^N f_1,\dots, p ^N f _n}{p ^N f _0}) \to \fX$.
We compute that this property holds 
if and only if the ideal of $\widehat{V _x}$
generated by $ \beta _x ^{*} (p ^N f _0)$
contains the ideal generated by $ \beta _x ^{*} ( p ^N f _i)$ for any $i = 1,\dots, n$.
This is equivalent to saying that 
$v _x ( f _i) \leq v _x ( f _0)$ for any $i = 1, \dots, n$, i.e. 
$v _x \in U ^{\ad}$.
Hence, we are done. 

\item \label{ZR-adic-affinoid-c)} 
Hence, these maps $\phi$ and $\psi$ are homeomorphisms inverse to each other.
Moreover, under this identification, it is checked in \cite[A.4.7]{FujiwaraKatoBookI} that 
the sheaf $\cO _{\fX _K} ^{\mathrm{int}}$
(resp. $\cO _{\fX _K} $)
coincides with the presheaf $\O _{\Spa (A _K, A _K ^{0})} ^+$ (resp. $\cO _{\Spa (A _K, A _K ^{0})}$).
\end{enumerate}
\end{empt}

\begin{empt}
[Local notation for affinoid spaces]
\label{ntnvx&v[x]}
Let $\fX=\Spf A$ be a formal $\fS$-scheme locally of formal finite type and $\V$-flat.
Let $x \in \fX _K$. 
We keep notation \ref{ZR-x-ntn}.

\begin{enumerate}[(a)]
\item Recall (see \ref{ZR-adic-affinoid}.\ref{ZR-adic-affinoid-a}), we denote by 
$v _x \colon A _K \to \Gamma _x \cup \{0\}$
the valuation induced by composing the valuation $v _{(A _x, I _x)}$ 
with the homomorphism of rings
$A _K \to B _x$.

\item We denote by 
$v _{[x]} \colon A _K \to \R _{\geq 0}$
the valuation induced by composing the valuation $v _{[A _x, I _x]}$ 
with the homomorphism of rings
$A _K \to B _x$.

\end{enumerate}

\end{empt}

\begin{empt}
[Zariski-Riemann triple]
\label{ZRtriple}
Let $\X$ be a formal $\fS$-scheme locally of formal finite type.

\begin{enumerate}[(a)]
\item From \ref{ZR-adic-affinoid}.\ref{ZR-adic-affinoid-c)}, we can view the 
the sheaf  $\cO _{\fX _K} $
as a sheaf of complete topological local rings.
 We denote the corresponding triple (see definition \ref{dfnTri}) by 
$${\bf ZR} (\fX)
:=(\fX _K, \O _{\fX _K} ^{\mathrm{int}}, \O _{\fX _K}).$$
The letters ${\bf ZR}$ refer to the Zariski-Riemann space.

\item When $\fX=\Spf A$ is moreover affine, 
then with notation \ref{AffAdSp} and \ref{dfn-rho-AdTri}, it follows from \ref{ZR-adic-affinoid}.\ref{ZR-adic-affinoid-c)}
that we get the isomorphism
\begin{equation}
\label{Aff-Spa=ZR}
{\bf ZR} (\fX)
\riso 
\iota _\mathrm{ATri}
(\Spa (A _K, A _K ^{0})).
\end{equation}

\item This yields that in general (i.e. when $\fX$ is not necessary affine)
${\bf ZR} (\fX)
\in {\bf AnAdTri}$.

Then we set 
$$\fX ^{\mathrm{ad}} := \kappa _\mathrm{ATri} ({\bf ZR} (\fX)) 
= 
(\fX _K, \O _{\fX _K}, \{ v _x \} _{x \in \fX _K})\in {\bf AnAd},$$
where the functor $\kappa _\mathrm{ATri}$ 
is defined at 
\ref{eqcatAnAdTri}. 

\end{enumerate}
\end{empt}

\subsection{Tubes, overconvergent singularities}

\begin{empt}
[Tubes of a closed or open subset]
\label{dfntubes}
Let $\X$ be a formal $\fS$-scheme locally of formal finite type.

\begin{enumerate}[(a)]
\item Let $Z$ be a closed subscheme of $\fX$.
We define the tube of $Z$ in $\fX$ by setting 
$$]Z[ _{\fX} := (\sp _{\fX} ^{-1} (Z)  )^{\circ}$$
to be the interior of the closed subset 
$\sp _{\fX} ^{-1} (Z) $.

\item Let $U$ be an open subscheme of $X$.
We define the tube of $U$ in $\fX$ by setting 
$$]U[ _{\fX} := \overline{\sp _{\fX} ^{-1} (U) }$$
to be the closure of 
$\sp  _{\fX} ^{-1} (U) $.

\item Suppose $\fX=\Spf A$ is affine. 
For any $a \in A$, we denote by 
$\overline{a}$ the image of $a$ in $A/\pi A$.
Choose $f _1, \dots, f _r\in A$ such that 
the  ideal defined by $Z$ in $X$ is 
generated by $\overline{f}_1,\dots, \overline{f} _r$.
With notation \ref{ntnvx&v[x]} and \ref{ZRtriple}, 
following \cite[II.4.2.11]{FujiwaraKatoBookI},
we get
\begin{equation}
\label{]Z[formula1}
]Z[ _{\fX} ^{\ad} 
=
\{
x \in \fX ^{\mathrm{ad}}
\;
;
\;
v _{[x]} (f _i) 
< 1
\text{, for any $i=1,\dots, r$}
\}.
\end{equation}
Choose $g _1, \dots, g _s\in A$ such that 
$U = \cup _{j=1} ^{s} D ( \overline{g} _j )$.
Taking the complementary of 
\ref{]Z[formula1},
we get
\begin{equation}
\label{]Z[formula1bis}
]U[ _{\fX} ^{\ad} 
=
\{
x \in \fX ^{\mathrm{ad}}
\;
;
\;
\exists j \in \{ 1,\dots, s\}
\text{, such that }
v _{[x]} (g _j) 
\geq 1
\}.
\end{equation}

\end{enumerate}

\end{empt}

\begin{lem}
[Tubes of a locally closed subset]
\label{lem-cofinal-sys-gen-case-pre}
Let $\fP$ be a formal $\fS$-scheme locally of formal finite type.
Let $Y$ be a subscheme of $P$. 
Let $X$ (resp. $X'$) be a  closed subscheme of $P$ and $U$ (resp. $U'$) be an open subscheme of $P$
and $Y = U \cap X$ (resp. $Y = U '\cap X'$). 
We have the equality
$
]U [ ^{\ad} _\fP
\cap ]X [ ^{\ad} _\fP
=
]U '[ ^{\ad} _\fP
\cap ]X '[ ^{\ad} _\fP$. 

\end{lem}

\begin{proof}
Remark we can suppose all schemes reduced. 
1) Suppose $X= X'$. Since $] X [ ^{\ad} _\fP$ is open, 
since 
$\sp ^{-1} _{\fP} (U) \cap ] X [ ^{\ad} _\fP
\subset 
\sp ^{-1} _{\fP} (Y)$,
then 
$]U [ ^{\ad} _\fP
\cap ]X [ ^{\ad} _\fP
\subset
\overline{\sp ^{-1} _{\fP} (Y)}$. 
Hence,
$]U [ ^{\ad} _\fP
\cap ]X [ ^{\ad} _\fP
=
\overline{\sp ^{-1} _{\fP} (Y)}
\cap ]X [ ^{\ad} _\fP $. 
This yields the independence on $U$. 

2) Suppose $U= U'$. 
Replacing $X$ or $X'$ by  $X \cap X'$ if necessary, 
we reduce to the case where $X ' \subset X$. 
We have to check that the inclusion 
$]U [ ^{\ad} _\fP
\cap ]X '[ ^{\ad} _\fP
\subset
]U [ ^{\ad} _\fP
\cap ]X [ ^{\ad} _\fP$
is in fact an equality. 
Since this is local, we can suppose 
$\fP = \Spf(A)$ is an affine formal $\fS$-scheme locally of formal finite type.
For any $a \in A$, we denote by 
$\overline{a}$ the image of $a$ in $A/\pi A$.
Choose functions
$f _1,\dots, f _r \in A$, 
$g_1, g_2,\dots,g _s\in A$  such that
$\overline{f} _1,\dots, \overline{f} _r$ generates the ideal defined by the closed immersion
$X \hookrightarrow P$, and 
$U =\cup _{j  = 1} ^{s} D (\overline{g} _j)$.
We reduce to suppose there exists
$f \in A$ such 
$X ' = X \cap V (\overline{f})$.
Let $j\in \{ 1,\dots, s\}$.
Since $X \cap D (\overline{g _j}) 
\subset V (\overline{f})$, then 
 there exists an integer $r _j \geq 1$
large enough so that
$\overline{f} \overline{g _j} ^{r _j}
\in (\overline{f} _1,\dots, \overline{f} _r)$. 
Hence,
$(\overline{f} _1,\dots, \overline{f} _r)
=
(\overline{f} _1,\dots, \overline{f} _r, 
\overline{f} \overline{g _1} ^{r _1},
\dots, 
\overline{f} \overline{g _s} ^{r _s})$.
This yields
\begin{align}
\notag
]U [ ^{\ad} _\fP
\cap
]X[ ^{\ad} _\fP
&
=
\{ 
x \in \fP ^{\ad}
\;
|
\; 
\exists j
\
v _{[x]} (g _j) \geq 1,
\forall i 
\ 
v _{[x]} (f _i) < 1,
\forall j 
\ 
v _{[x]} (f g _j ^{r _j}) < 1
\}
\\
\notag
&
=
\{ 
x \in \fP ^{\ad}
\;
|
\; 
\exists j
\
v _{[x]} (g _j) \geq 1,
\forall i 
\ 
v _{[x]} (f _i) < 1,
\ 
v _{[x]} (f ) < 1
\}
\\
\notag
&
=
]U [ ^{\ad} _\fP
\cap
]X'[ ^{\ad} _\fP.
\end{align}
\end{proof}

\begin{dfn}
[Tubes of a locally closed subset]
\label{cofinal-sys-gen-case-pre}
Let $\fP$ be a formal $\fS$-scheme locally of formal finite type.
Let $X$ be a  closed subscheme of $P$ and $U$ be an open subscheme of $P$
and $Y = U \cap X$. 
We set 
$]Y [ ^{\ad} _\fP := 
]U [ ^{\ad} _\fP
\cap ]X [ ^{\ad} _\fP$. 
The tube $]Y [ ^{\ad} _\fP$ is independent of the choice of $U$ and $X$ such that 
$Y = U \cap X$.
\end{dfn}

\begin{lem}
\label{lem-vxVSv[x]}
Let $\fX=\Spf A$ be a formal $\fS$-scheme locally of formal finite type and $\V$-flat.
Let $x \in \fX ^{\ad}, f \in A _K$. 
With notation 
\ref{ntnvx&v[x]},
we have the properties
\begin{gather}
\label{vxVSv[x]}
v _x (f) \leq 1
\Rightarrow
v _{[x]} (f) \leq 1 
\
; 
\ 
v _{[x]} (f) > 1
\Rightarrow
v _x (f) > 1 ;
\\
\label{vxVSv[x]bis}
 1 \leq v _x (f) 
\Rightarrow
1\leq v _{[x]} (f) 
\
; 
\ 
v _{[x]} (f) < 1
\Rightarrow
v _x (f) < 1 .
\end{gather}
\end{lem}

\begin{proof}
We denote by $f _x\in K _x$ the image of $f$ via 
the canonical map $A _K \to K _x$.
By definition, 
$v _x (f) \leq 1$ if and only if 
$f _x\in V _x$. 
Hence, 
$f _x \in V _{x,\fp _x}$, 
i.e. 
$v _{[x]} (f) \leq 1$.
This yields \ref{vxVSv[x]}. 
Moreover, 
$ 1 \leq v _x (f) $ if and only if 
$1 \in f _xV _x$.
This yields
$1 \in f _x V _{x,\fp _x}$, 
i.e. 
$1\leq v _{[x]} (f) $.
Hence, we have checked \ref{vxVSv[x]bis}.
\end{proof}

\begin{lem}
\label{lem2-vxVSv[x]}
Let $\fX=\Spf A$ be a formal $\fS$-scheme locally of formal finite type and $\V$-flat.
Let $x \in \fX ^{\ad}, f \in A _K$. 
With notation 
\ref{ntnvx&v[x]},
for any integer $n \geq1$
we have the inclusions
\begin{gather}
\label{lem2-vxVSv[x]-incl}
\{x \in \fX ^{\ad} 
\; 
|
\; 
v _{[x]} (\pi ^{-1} f ^{n+1}) \geq 1\}
\subset
\{x \in \fX ^{\ad} 
\; 
|
\; 
v _{[x]} (\pi ^{-1} f ^{n}) >1\},
\\
\label{lem2-vxVSv[x]-incl2}
\{x \in \fX ^{\ad} 
\; 
|
\; 
v _{x} (\pi ^{-1} f ^{n+1}) \geq 1\}
\subset
\{x \in \fX ^{\ad} 
\; 
|
\; 
v _{x}  (\pi ^{-1} f ^{n}) >1\}.
\end{gather}

\end{lem}

\begin{proof}
The property 
$v _{[x]} (\pi ^{-1} f ^{n}) > 1$ 
(resp. $v _{[x]} (\pi ^{-1} f ^{n+1}) \geq 1$)
is equivalent to 
$v _{[x]} (f) > (\frac{1}{p}) ^{1/en}$
(resp. $v _{[x]} (f) \geq (\frac{1}{p}) ^{1/e(n+1)}$).
Since 
$(\frac{1}{p}) ^{1/e(n+1)} >  (\frac{1}{p}) ^{1/en}$, then we get
\ref{lem2-vxVSv[x]-incl}.
Following \ref{lem-vxVSv[x]},
we have the inclusions
\begin{gather}
\notag
\{x \in \fX ^{\ad} 
\; 
|
\; 
v _{x} (\pi ^{-1} f ^{n+1}) \geq 1\}
\subset
\{x \in \fX ^{\ad} 
\; 
|
\; 
v _{[x]} (\pi ^{-1} f ^{n+1}) \geq 1\},
\\
\notag 
\{x \in \fX ^{\ad} 
\; 
|
\; 
v _{[x]} (\pi ^{-1} f ^{n}) >1\}
\subset
\{x \in \fX ^{\ad} 
\; 
|
\; 
v _{x} (\pi ^{-1} f ^{n}) >1\}.
\end{gather}
Hence, we are one.
\end{proof}

\begin{empt}
\label{empt-tubZndfn}
Let $\fX$ be a formal $\fS$-scheme locally of formal finite type and $\V$-flat.
Let $Z$ be a closed subscheme of $X$.

\begin{enumerate}[(a)]
\item 
\label{empt-tubZndfn(a)}
Suppose $\fX=\Spf A$ is affine. 
For any $a \in A$, we denote by 
$\overline{a}$ the image of $a$ in $A/\pi A$.
Choose $f _1, \dots, f _r\in A$ such that 
the  ideal defined by $Z$ in $X$ is 
generated by $\overline{f}_1,\dots, \overline{f} _r$.
For any integer $n \geq 1$, we set
\begin{gather}
\label{tubZndfn}
[Z] _{\fX, n} ^{\ad}
:=
\{
x \in \fX ^{\mathrm{ad}}
\;
;
\;
v _{x} (\pi ^{-1}f _i ^n) 
\leq 1
\text{, for $i=1,\dots, r$}
\},
\\
\label{tubZndfnopen}
]Z[ _{\fX, n} ^{\ad}
:=
\{
x \in \fX ^{\mathrm{ad}}
\;
;
\;
v _{x} (\pi ^{-1}f _i ^n) 
<1
\text{, for $i=1,\dots, r$}
\}.
\end{gather}
Since 
$\fX ^{\mathrm{ad}}
\riso 
\Spa (A _K, A _K ^{0})$
(see \ref{Aff-Spa=ZR}),
since the underlying homeomorphism
is $x \mapsto v _x$ (see \ref{ZR-adic-affinoid}), 
then
we have the isomorphism 
$$[Z] _{\fX, n} ^{\ad}
\riso 
\Spa \left ( ( A _K , A _K ^0)
<\tfrac{f_1 ^{n},\dots, f _r ^{n}}{\pi}>\right ) ,$$
where the right side is a rational subset of $\Spa ( A _K, A _K ^0 )$
(see \ref{15(ii)Huber}).
In particular, 
$[Z] _{\fX, n} ^{\ad}$ is an affinoid open subspace of $\fX ^{\ad}$.

\item 
\label{empt-tubZndfn(b)}
Let us check that the definition of $[Z] _{\fX, n} ^{\ad}$ 
given at \ref{tubZndfn} does not depend on the choice. 
Let $g_1, g_2,\dots,g _s\in A$  such that
$\overline{g} _1,\dots, \overline{g} _s$ generate the ideal given by the closed immersion
$Z \hookrightarrow X$.
Then, there exists 
$(a _{ji})  _{ji} \in M _{s,r} (A)$
and 
$(a _{j})  _{j} \in M _{s,1} (A)$
such that 
$g _j = \sum _{i = 1} ^{r} a _{ji} f _i + \pi a _j$, 
for any $j = 1,\dots, s$.
Let $x \in \fX ^{\ad} $ such that
$v _{x} (\pi ^{-1}f _i ^n) \leq 1$.
We compute 
$v _{x} (\pi ^{-1}f _i ^n) \leq 1
\Leftrightarrow
 v _x (f _i) ^n
\leq 
v _{x} (\pi) < 1$
for any $i=1,\dots, r$.
Since 
$v _x (a _{ji}) \leq 1$ and
$v _x (a _{j}) \leq 1$, then we get 
$v _x ( g _j ) \leq \max \{  \{v _x ( f  _i)\; ; \; i =1,\dots, r\} \cup \{v _x ( \pi)\} \}$. 
This yields 
$v _x ( g _j ) ^n \leq \max \{  \{v _x ( f  _i) ^n \; ; \; i =1,\dots, r\} \cup \{v _x ( \pi) ^n \} \}
\leq v _x ( \pi) $.
Hence, by symmetry we get the equality
\begin{equation}
\label{ind-affine1}
\{
x \in \fX ^{\mathrm{ad}}
\;
;
\;
v _{x} (\pi ^{-1}f _i ^n) 
\leq 1
\text{, for $i=1,\dots, r$}
\}
=
\{
x \in \fX ^{\mathrm{ad}}
\;
;
\;
v _{x} (\pi ^{-1}g _j ^n) 
\leq 1
\text{, for $j=1,\dots, s$}
\},
\end{equation}
which means that the affinoid $[Z] _{\fX, n} ^{\ad}$ is independent of the choice of the $f _1, \dots, f _r$. 

\item
\label{empt-tubZndfn(c)}
From the independence of (b), 
the affinoids $[Z] _{\fX, n} ^{\ad}$ glue over an open affine covering of $\fX$.
 We still denote it by  $[Z] _{\fX, n} ^{\ad}$.
This is an open adic subspace of $\fX ^{\ad}$. 
 
 \item 
 \label{empt-tubZndfn(d)}
 Suppose $n \geq 2$. 
 Copying the computations of the part  (b), 
 since $v _x ( \pi) ^n < v _x ( \pi)$, then we can check that 
the set $]Z[ _{\fX, n} ^{\ad}$ 
is independent of the choice of the $f _1, \dots, f _r$. 
Hence, 
$]Z[ _{\fX, n} ^{\ad}$ glue over an open affine covering of $\fX$.
We still denote it by  $]Z[ _{\fX, n} ^{\ad}$.
Beware this is not necessary an adic space. 

\item Let $Y: = X \setminus Z$ be the open subscheme of $X$.
For any integer $n \geq 1$, 
we set 
$]Y[ _{\fX, n} ^{\ad}:= \X ^{\ad} \setminus [Z] _{\fX, n} ^{\ad}$.
For any $n \geq 2$, we set
$[Y] _{\fX, n} ^{\ad}:= \X ^{\ad} \setminus ]Z[ _{\fX, n} ^{\ad}$.
When $\fX = \Spf A$ is affine, 
$f _1, \dots, f _r\in A$ are such that 
the  ideal defined by $Z$ in $X$
is generated by $\overline{f}_1,\dots, \overline{f} _r$, 
then 
$[Y] _{\fX, n} ^{\ad}$ is the union of the $r$ affinoids
defined by the equation
$v _{x} (\pi ^{-1}f _i ^n) \geq 1$.
Hence, 
$[Y] _{\fX, n} ^{\ad}$ is endowed with a structure of adic space
so that is an open adic subspace of $\fX ^{\ad}$.

\end{enumerate}

\end{empt}

\begin{prop}
\label{prop-]Z[=cup[Z]n}
Let $\fX$ be a formal $\fS$-scheme locally of formal finite type and $\V$-flat.
Let $Z$ be a closed subscheme of $X$
and let 
$Y: = X \setminus Z$ be the corresponding open subscheme of $X$.
We get the equalities 
\begin{gather}
\label{prop-]Z[=cup[Z]n=1}
]Z[ ^{\mathrm{ad}}_{\fX} = \cup _{n\geq 1} [Z] _{\fX, n} ^{\mathrm{ad}}
,
\
]Z[ ^{\mathrm{ad}}_{\fX} = \cup _{n\geq 2} ]Z[ _{\fX, n} ^{\mathrm{ad}}
\\
\label{prop-]Z[=cup[Z]n=2}
]Y[ ^{\mathrm{ad}}_{\fX} = \cap _{n\geq 1} ]Y[ _{\fX, n} ^{\mathrm{ad}}
,
\ 
]Y[ ^{\mathrm{ad}}_{\fX} = \cap _{n\geq 2} [Y] _{\fX, n} ^{\mathrm{ad}}.
\end{gather}
\end{prop}

\begin{proof}
1) Let us check the equality of \ref{prop-]Z[=cup[Z]n=1}. 
Since this is local, we can suppose 
$\fX= \Spf A$ affine. 
For any $a \in A$, we denote by 
$\overline{a}$ the image of $a$ in $A/\pi A$.
Choose $f _1, \dots, f _r\in A$ such that 
the  ideal defined by $Z$ in $X$ is 
generated by $\overline{f}_1,\dots, \overline{f} _r$.

a) Let us check the inclusion 
$]Z[ ^{\mathrm{ad}}_{\fX} \subset \cup _{n\geq 1} ]Z[ _{\fX, n} ^{\mathrm{ad}}$.
Let $x \in ]Z[ ^{\mathrm{ad}} _{\fX} $. Let $\rho:= v _{[x]} ( \pi) <1$.
Then, for $n$ large enough,
$v _{[x]} (f _i) < \rho ^{1/n}$, i.e. 
$v _{[x]} (\pi ^{-1}f _i ^n)<1$.
Using \ref{vxVSv[x]bis}, 
this yields
$v _{x} (\pi ^{-1}f _i ^n)<1$,
i.e. 
$x \in ]Z[ ^{\mathrm{ad}} _{\fX, n} $.

b) Since
$]Z[ ^{\mathrm{ad}} _{\fX, n} 
\subset 
[Z] ^{\mathrm{ad}} _{\fX, n} $, 
then 
$ \cup _{n\geq 2} ]Z[ _{\fX, n} ^{\mathrm{ad}}
\subset
\cup _{n\geq 2} [Z] _{\fX, n} ^{\mathrm{ad}}$.

c) 
It remains to check 
$[Z] _{\fX, n} ^{\mathrm{ad}}
\subset ]Z[ ^{\mathrm{ad}} _{\fX}$.
Let  $x \in [Z] ^{\mathrm{ad}} _{\fX, n} $ for some $n$. 
Using \ref{vxVSv[x]}, 
we get $v _{[x]} (\pi ^{-1}f _i ^n) \leq 1$, 
which is equivalent to saying that
$v _{[x]} (f _i) \leq \rho ^{1/n} <1$.
Hence, $x \in  ]Z[ ^{\mathrm{ad}} _{\fX} $.

2) By taking the complementary, we get  \ref{prop-]Z[=cup[Z]n=2} from  \ref{prop-]Z[=cup[Z]n=1}.
\end{proof}

\begin{coro}
\label{coro-Vsupsetggeq1}
Let $\fX=\Spf A$ be a formal $\fS$-scheme locally of formal finite type and $\V$-flat.
Let $V$ be an open subset of $\fX ^{\ad}$
and $g \in A _K $ such that 
\begin{equation}
\label{Vsupsetggeq1}
V \supset 
\{x \in \fX ^{\ad} 
\; 
|
\; 
v _{[x]} ( g) \geq 1\}.
\end{equation}
Then, there exists some integer $m \geq 1$ such that
$$V \supset 
\{x \in \fX ^{\ad} 
\; 
|
\; 
v _{x} (\pi ^{-1} g ^{m}) \geq 1\}.$$
\end{coro}

\begin{proof}
Since $\fX ^{\ad} $ is quasi-compact (see \ref{coll-top-XK}.\ref{qc-ZRqc}), 
since 
$T := \fX ^{\ad} \setminus V$
is a closed subset, 
then $T$ is quasi-compact.
The inclusion 
\ref{Vsupsetggeq1}
is equivalent to the following one:
$T
 \subset 
\{x \in \fX ^{\ad} 
\; 
|
\; 
v _{[x]} ( g) <1\}$.
Let $Z$ be the closed subscheme of $\fX = \Spf A$ defined by the admissible ideal
$I = (\pi ,g ) \subset A$.
Hence, 
$]Z[ ^{\mathrm{ad}}_{\fX} 
=
\{x \in \fX ^{\ad} 
\; 
|
\; 
v _{[x]} ( g) <1\}$.
Using \ref{prop-]Z[=cup[Z]n}, 
we get the affinoid covering 
$T
 \subset 
 \cup _{n\geq 1} [Z] _{\fX, n} ^{\mathrm{ad}}$.
 Since $T$ is quasi-compact, 
 for $n $ large enough, 
 we get
 $T
 \subset 
  [Z] _{\fX, n} ^{\mathrm{ad}}$.
  Hence,
 $V \supset 
\{x \in \fX ^{\ad} 
\; 
|
\; 
v _{x} (\pi ^{-1} g ^{n}) >1\}$.
Following \ref{lem2-vxVSv[x]-incl2},
this implies 
$V \supset 
\{x \in \fX ^{\ad} 
\; 
|
\; 
v _{x} (\pi ^{-1} g ^{n+1}) \geq 1\}$.
\end{proof}

\begin{empt}
\label{cofinal-sys-gen-casebis}
Let $\fP$ be a formal $\fS$-scheme locally of formal finite type.
Let $X$ be a closed subscheme of $P$ and $U$ be an open subscheme of $P$. 
We put $Y := X \cap U$. 
We denote by $S ^+ (\bbN,\bbN)$
the set  of increasing sequence $\underline{m}\colon \N\to \N$ of integers 
such that $\underline{m}(n) \to \infty$ when $n \to \infty$.

\begin{enumerate}[(a)]
\item For integers $n \geq 1$ and $m \geq 2$, 
we have the adic spaces 
$V ^{\ad} _{n,m}
:=
[X] _{\fP, n} ^{\ad}
\cap
[U] _{\fP, m} ^{\ad}$
(see notation \ref{empt-tubZndfn}).
Now, for any $\underline{m}\in S ^+ (\bbN,\bbN)$, we set
\begin{equation}
\label{ntn-Vmad}
V ^{\ad} _{\underline{m}}
:=
\cup _{n\in \N} V ^{\ad} _{n,\underline{m}(n)}.
\end{equation}

\item Copying the proof of \cite[2.18]{Lazda-Pal-Book},
it follows from \ref{coro-Vsupsetggeq1} that,
as
$\underline{m}$ varies,
the
$V ^{\ad} _{\underline{m}}$
form a cofinal system of open neighbourhoods
of $]Y [ ^{\ad} _\fP$
inside 
$]X [ ^{\ad} _\fP$.

\item We set $j \colon 
]Y [ ^{\ad} _\fP
\to 
]X [ ^{\ad} _\fP$
and 
$j _{\underline{m}}
\colon
V ^{\ad} _{\underline{m}}
\to 
]X [ ^{\ad} _\fP$.
For any sheaf $\cF$ on $]X [ ^{\ad} _\fP$,
we set 
$ j ^\dag _Y  \cF 
:= 
j _* j ^{-1} (\cF)$.
Using  \cite[2.19]{Lazda-Pal-Book},
we get from part $(b)$ the canonical isomorphism
\begin{equation}
\label{2.19LP}
j ^\dag _Y  \cF 
\riso
\underset{\underline{m}\in S ^+ (\bbN,\bbN)}{\underrightarrow{\lim}} 
\,
j _{\underline{m}*} j _{\underline{m}} ^{-1} (\cF).
\end{equation}
\end{enumerate}

\end{empt}

\begin{ex}
\label{ntn-jm}
Let $\X $ be a formal $\fS$-scheme locally of formal finite type and $\V$-flat.
Let $Z$ be a divisor of $X$ and $\fY $ the open of $\fX $ complementary to the support of $Z$,
and $j\colon \fY \hookrightarrow \fX$ be the open immersion.

1) 
For any integer $m \geq 0$, we set
$r _m : = p ^{m+1}$
and 
$Y _m := [Y] _{\fX, r _m} ^{\ad}$.
When $\fX=\Spf A$ is affine and there exists
$g \in A$ such that $\fY = D (g)$, then  
$Y _m $ is an affinoid open subspace and more precisely
\begin{equation}
\label{ThA-B-affin-[Y]m}
Y _m 
=
\{
x \in \fX ^{\mathrm{ad}}
\;
;
\;
v _{x} (\pi ^{-1} g ^{p ^{m+1}}) \geq 1\}
\riso
\Spa \left ( ( A _K, A _K ^0) 
< \frac{\pi}{g ^{p ^{m+1}}}>\right ),
\end{equation}
where the right term is defined at 
\ref{ntn-Spa<>}.

2) 
Following \ref{prop-]Z[=cup[Z]n}, 
$]Y [ ^{\ad} _\fP= \cap _{m \in \N} 
Y _m$.
Following \ref{cofinal-sys-gen-casebis}, 
$(Y _m )_{m\in\N}$ forms a cofinal system of neighbourhoods of 
$]Y[  ^{\ad}_{\X} $ in $ \X ^{\ad}$.
We denote by $j _Y \colon ]Y[ ^{\ad} _{\X} \hookrightarrow \X _K$  the canonical open immersion 
induced by $j$.
We denote by $j _m \colon Y _m \hookrightarrow \X _K$ the canonical immersion. 
We set  $j ^\dag \G = j _{Y*} j ^{-1} _Y \G$ for any $\O _{\X _K}$-module $\G$.  
From 
\ref{2.19LP},
we get 
$j ^\dag \G \cong \underrightarrow{\lim} _m j _{m*} j ^{-1} _m \G$.

\end{ex}

\subsection{Affinity of the specialisation morphism, local cohomology,
Cech resolutions}

\begin{empt}
We keep notation \ref{ntn-jm} and we suppose $\fX$ affine.
We have the canonical isomorphism
$$
j ^\dag _Y  \cO _{\fX ^{\ad}}
\riso
\underrightarrow{\lim} _{m}
\,
j _{m*} 
\cO _{Y _m}.
$$
By a rational subset of 
$\X ^{\ad}$ we means an 
open subset of $\fX _K$ (equal to the underlying topological space of 
$\X ^{\ad}$)
which corresponds to 
a  rational subset of $\Spa (A _K, A _K ^{0})$
via the canonical homeomorphism
$\fX _K
\to
\Spa (A _K, A _K ^{0})$ (see \ref{ZR-adic-affinoid}).
Let  $m,r \geq 1$ be two integers, $U$ be a rational subset  of $\X ^{\ad}$.
Since $Y _m$ is a rational subset of $\X ^{\ad}$, then
$j _m ^{-1} (U) = U \cap Y _m$ is also a rational subset of $\X ^{\ad}$
and then is a rational subset of
$Y _m$ (see \cite[Lemma 1.5.(ii)]{Huber-gen-rig-an-var}).
Hence, 
following \cite[2.2]{Huber-gen-rig-an-var}, we have
$H ^r ( j _m ^{-1} (U)  , \cO _{Y _m})=0$.
This yields that 
$R ^r j _{m*} 
(\cO _{Y _m})=0$
(see \cite[0.12.2.1]{EGAIII1}).
This means that the canonical morphism
\begin{equation}
\label{Rjm*1}
j _{m*} 
(\cO _{Y _m})
\to
\R j _{m*} 
(\cO _{Y _m})
\end{equation}

is an isomorphism.
\end{empt}

The following Lemma will be useful to check the resolution
\ref{RspCheck-reso}.
\begin{lem}
\label{lem-Ru*ad=u*ad}
Let $\fX =\Spf A $ be an affine formal $\fS$-scheme locally of formal finite type and $\V$-flat.
Let $g \in A$, $\fY = D (g)$ be the corresponding standard open formal subscheme of $\fX$.
Let $u\colon \fX \hookrightarrow \fX '$ be an open immersion of 
 separated formal $\fS$-schemes locally of formal finite type and $\V$-flat.
Let 
$u ^{\ad}\colon \fX ^{\ad} \hookrightarrow \fX ^{ \prime \ad}$ be the induced morphism of adic spaces.
\begin{enumerate}[(a)]
\item The canonical morphism 
\begin{equation}
\label{Ru*ad=u*adpre}
u _{*} ^{\ad} (j ^\dag _Y  \cO _{\fX ^{\ad}})
\to
\R u _{*} ^{\ad} (j ^\dag _Y  \cO _{\fX ^{\ad}})
\end{equation}
is an isomorphism.

\item The module $u _{*} ^{\ad} (j ^\dag _Y  \cO _{\fX ^{\ad}})$ is acyclic for $\sp  _*$, i.e., 
the canonical morphism 
\begin{equation}
\label{Ru*ad=u*ad}
\sp _* u _{*} ^{\ad} (j ^\dag _Y  \cO _{\fX ^{\ad}})
\to
\R \sp _*  u _{*} ^{\ad} (j ^\dag _Y  \cO _{\fX ^{\ad}})
\end{equation}
is an isomorphism.
\end{enumerate}

\end{lem}

\begin{proof}
1) Let us check the first statement. 
This is local in $\fX ^{\prime \ad}$. 
This yields, since $\fX'$ is separated, that we can suppose both $\fX$ and $\fX '$ are affine (and $\fY$ is still a standard open formal subscheme  of $\fX$). 
Let  $r \geq 1$ be an integer, $U' $ be a rational subset  of $\X ^{\prime\ad }$.
Since $\X ^{\ad}$ and $Y _m$ are affinoid adic spaces (we use notation \ref{ntn-jm}),
then $U := (u ^{\ad}) {} ^{-1} (U')$ is a rational subset of $\X ^{\ad }$, 
and $j _m ^{-1} (U)$ is a rational subset of 
$Y _m$ (see \cite[Lemma 1.5.(ii)]{Huber-gen-rig-an-var}).
Following \cite[2.2]{Huber-gen-rig-an-var}, 
this implies
$H ^r ( (u ^{\ad} \circ j _m) ^{-1} (U ')  , \cO _{Y _m})=0$.
This yields that the canonical morphism
$u _{*} ^{\ad} \circ j _{m*} 
(\cO _{Y _m})
\to
\R ( u _{*} ^{\ad} \circ j _{m*} )
(\cO _{Y _m})$
is an isomorphism.
From the isomorphism \ref{Rjm*1}, 
we get
$\R ( u _{*} ^{\ad} \circ j _{m*} )
(\cO _{Y _m})
\riso 
\R  u _{*} ^{\ad} 
( j _{m*} (\cO _{Y _m}))$.
Hence the canonical morphism
$u _{*} ^{\ad} 
(j _{m*}  \cO _{Y _m})
\to
\R u _{*} ^{\ad} 
(j _{m*}  \cO _{Y _m})$
is an isomorphism, 
i.e. 
$R ^r u _{*} ^{\ad} 
(j _{m*}  \cO _{Y _m}) =0$ for any $r \geq 1$.
Since $u ^{\ad}$ is a coherent morphism of coherent topological spaces,
then inductive limits commutes with 
$R ^r u _{*} ^{\ad} $ (see \cite[VI.5.1]{sga4-2}, or also \cite[0.3.1.9]{FujiwaraKatoBookI}).
Hence, taking the inductive limit,
this yields 
$R ^r u _{*} ^{\ad} 
(j ^\dag _Y  \cO _{\fX ^{\ad}})=0$ for any $r \geq 1$, i.e. 
that 
the canonical morphism
\ref{Ru*ad=u*ad}
is an isomorphism.

2) Using the same arguments than in the first part (i.e. \cite[2.2]{Huber-gen-rig-an-var} and next taking the inductive limits), 
we check the canonical morphism 
\begin{equation}
\label{Ru*ad=u*ad2}
(\sp _* \circ u _{*} ^{\ad} )(j ^\dag _Y  \cO _{\fX ^{\ad}})
\to
\R (\sp _* \circ   u _{*} ^{\ad})  (j ^\dag _Y  \cO _{\fX ^{\ad}})
\end{equation}
is an isomorphism.
Using \ref{Ru*ad=u*adpre} and \ref{Ru*ad=u*ad2}, we get \ref{Ru*ad=u*ad}.
\end{proof}

\begin{empt}
Let $\fP$ be a formal $\fS$-scheme locally of formal finite type and $\V$-flat.
Let $X$ be a closed subscheme of $P$ and $Y$ be an open subscheme of $X$. 
\label{loc-coh&cech-reso}
\begin{enumerate}[(a)]
\item \label{loc-coh&cech-reso1}
If $Z = Y \setminus X$ and
$E$ is any sheaf on 
$]X [ ^{\ad} _\fP$ then we define
$\underline{\Gamma} ^{\dag} _Z E$ by the exact sequence 
\begin{equation}
\label{loc-coh&cech-reso1-ex-tri}
0\to\underline{\Gamma} ^{\dag} _Z E\to E \to j _Y ^\dag  E\to0 .
\end{equation}

Note that $j _{Y} ^\dag$ and
$\underline{\Gamma} ^{\dag} _Z $
are exact, and we have 
$j _{Y} ^\dag j _{Y'} ^\dag E 
\riso j _{Y \cap Y'} ^\dag$ 
and
$\underline{\Gamma} ^{\dag} _{Z} \underline{\Gamma} ^{\dag} _{Z'} E 
\riso 
\underline{\Gamma} ^{\dag} _{Z \cap Z'} E$, 
for any open subsets $Y$ and $Y'$ of $P$
and any closed subsets $Z$ and $Z'$ of $P$.

\item  \label{loc-coh&cech-reso2}
Let $\mathscr{Y}:= (Y _i) _{i=1,\dots, r}$ be a finite open covering of $Y$. 
Similarly to \cite[2.50]{Lazda-Pal-Book},
we check that for any sheaf $E$ on $]X [ ^{\ad} _\fP$ there is an exact sequence of sheaves 
\begin{equation}
\label{2.50Lazda-Pal-Book}
0 
\to
j _Y ^\dag  E
\to 
\prod _{i = 1} ^{r}
j _{Y _i} ^\dag  E
\to 
\prod _{1\leq i _0 < i _1\leq r}
j _{Y _{i _0}\cap Y _{i _1}} ^\dag  E
\to 
\cdots
\to 
j _{\cap _{i = 1} ^{r} Y _{i}} ^\dag  E
\to 
0
\end{equation}
on $]X [ ^{\ad} _\fP$.
Denoting by 
$\check{C} ^{\dag \bullet} 
(\fX, \mathscr{Y},E)$ 
the complex
\begin{equation}
\label{2.50Lazda-Pal-Book-ntn}
\dots \to 0 \to 
\prod _{i = 1} ^{r} j _{Y _i} ^\dag  E
\to 
\prod _{1\leq i _0 < i _1\leq r}
j _{Y _{i _0}\cap Y _{i _1}} ^\dag  E
\to 
\cdots
\to 
j _{\cap _{i = 1} ^{r} Y _{i}} ^\dag  E
\to 
0
\to \dots,
\end{equation}
whose $0th$ term is 
$\prod _{i = 1} ^{r} j _{Y _i} ^\dag  E$,
this means that 
$\check{C} ^{\dag \bullet} 
(\fX, \mathscr{Y},E)$
is a resolution of 
$j _Y ^\dag  E$.
\end{enumerate}

\end{empt}

\section{Sheaf of differential operators on schemes locally of formal finite type}

\subsection{Sheaf of differential operators}

\begin{empt}
Put $S:= S _i$.
Let $X$ be an $S$-scheme 
locally of formal finite type
and having  locally finite $p$-bases  over $S$ 
(hence $X/S$ is flat following Theorem \ref{f0formétale-fforméta0}).
Let $m \geq 0$ be an integer.
Following \cite{Caro-Vauclair},
the sheaf of differential operators of level $m$
on $X /S $ denoted by $\D ^{(m)} _{X /S }$  
is well defined and 
we have the usual local description 
(more precisely, see for instance \ref{ntn-loc-coor-TvsS}).

\end{empt}

\begin{empt}
\label{2.2Be1}
Let $ \X$ be a formal  $\fS$-scheme
locally of formal finite type
and having locally finite $p$-bases  over $\fS$.
Let $m \geq 0$ be an integer.
Following \cite{Caro-Vauclair},
the sheaf of differential operators 
$\D ^{(m)} _{\fX /\fS }$ 
on $\fX /\fS $ is well defined and 
we have the usual local description.
We get the sheaf 
$\D ^{\dag} _{\fX /\fS }
:=
\underset{m}{\underrightarrow{\lim}}
\,
\widehat{\D} ^{(m)} _{\fX /\fS }$.

\end{empt}

\subsection{PD-stratification of level $m$, left $\cD ^{(m)} _{X/T}$-modules, inverse images}
Let $T$ be an ${S _i}$-scheme of finite type (resp. a formal $\fS$-scheme of finite type)
for some integer $i \geq 0$.
Let $X$ be an $T$-scheme (resp. a formal $T$-scheme)
locally of formal finite type
and having  locally finite $p$-bases  over $T$.
For the reader, let us recall the definition of a PD-stratification of level $m$ relatively to $T$
and its link with a structure of left $\cD ^{(m)} _{X/T}$-module
(see \cite{Caro-Vauclair} for a logarithmic version).

\begin{ntn}
For any $n,n'\in \N$, the sheaf $\cP ^n _{X/T (m)} \otimes _{\O _X} \cP ^{n'} _{X/T (m)}$
has three structures of $\O _X$-algebras. The structure of  $\O _X$-algebra of 
$\cP ^n _{X/T (m)}\otimes _{\O _X} \cP ^{n'} _{X/T (m)}$ coming from the left structure of  $\cP ^n _{X/T (m)} $ is said to be the left one, 
 that coming from the tensor product is said to be the middle one, that coming from the right structure on
 $\cP ^{n'} _{X/T (m)}$ is said to be  the right one. 
We denote by $d ^{n,n'} _0$, $d ^{n,n'}_1$, $d ^{n,n'} _2$ the corresponding structural homomorphisms 
$\O _X \rightarrow \cP ^n _{X/T (m)} \otimes _{\O _X} \cP ^{n'} _{X/T (m)}$.
We denote by $\delta ^{n,n'} _{(m)}$ : $\cP ^{n + n'} _{X/T (m)} \rightarrow \cP ^n _{X/T (m)} \otimes _{\O _X} \cP ^{n'}
_{X/T (m)}$ the morphism defined in \cite[2.1.3]{Be1} (see also \cite[2.2.14]{Caro-Vauclair}). 
Finally, we denote by 
$q _0 ^{n,n'}$ et $q _1 ^{n,n'}$ : $\cP ^{n + n'} _{X/T (m)} \rightarrow \cP ^n _{X/T (m)} \otimes _{\O _X} \cP
^{n'} _{X/T (m)}$ the natural homomorphisms defined in \cite[2.3.1]{Be1}. 
\end{ntn}

\begin{dfn}
Let $\E$ be an $\O _X$-module. An $m$-PD-stratification (or a PD-stratification of level $m$)
relatively to $T$
is the data of a family of compatible 
$\cP ^{n} _{X/T, (m)}$-linear isomorphisms
$$\epsilon ^{\E} _n
\colon
\cP ^{n} _{X/T, (m)} \otimes _{\O _X} \E
\riso
\E \otimes _{\O _X} \cP ^{n} _{X/T, (m)}$$
satisfying the following conditions:
\begin{enumerate}[(a)]
\item $\epsilon ^{\E} _0 = \mathrm{Id} _\E$ ;
\item for any $n, n'$, the diagram
\begin{equation}
\notag
\xymatrix{
{\cP ^{n} _{X/T, (m)} \otimes _{\O _X} \cP ^{n'} _{X/T, (m)} \otimes _{\O _X} \E}
\ar[rr] ^-{\delta _{(m)} ^{n,n' *} (\epsilon ^{\E} _{n+n'})} _-{\sim}
\ar[rd] _-{q _{1 (m)} ^{n,n' *} (\epsilon ^{\E} _{n+n'})} ^-{\sim}
&&
{ \E \otimes _{\O _X} \cP ^{n} _{X/T, (m)} \otimes _{\O _X} \cP ^{n'} _{X/T, (m)}}
\\
&
{\cP ^{n} _{X/T, (m)}  \otimes _{\O _X} \E \otimes _{\O _X} \cP ^{n'} _{X/T, (m)}}
\ar[ur] _-{q _{0 (m)} ^{n,n' *} (\epsilon ^{\E} _{n+n'})} ^-{\sim}
&
}
\end{equation}
is commutative
\end{enumerate}

\end{dfn}

\begin{prop}
\label{strat-prop}
Let $\E$ be an $\O _X$-module.
The following datas are equivalent :
\begin{enumerate}[(a)]
\item A structure of left $\D ^{(m)} _{X/T}$-module on $\E$ extending its structure of
$\O _X$-module.
\item A family of compatible $\O _X$-linear homomorphisms
$\theta ^{\E} _n \colon \E \to \E \otimes _{\O _X} \cP ^{n} _{X/T, (m)} $
 such that $\theta ^{\E} _0= \mathrm{Id} _\E$ and for any integers $n,n'$ the diagram
 \begin{equation}
 \label{diag-theta}
 \xymatrix{
 {\E \otimes _{\O _X} \cP ^{n} _{X/T, (m)} }
\ar[r] ^-{\mathrm{Id}\otimes \delta _{(m)} ^{n,n'} }
 &
 { \E \otimes _{\O _X} \cP ^{n} _{X/T, (m)}   \otimes _{\O _X} \cP ^{n'} _{X/T, (m)}   }
 \\
 {\E}
 \ar[u] ^-{\theta ^{\E} _{n+n'}}
  \ar[r] ^-{\theta ^{\E} _{n'}}
 &
 {\E \otimes _{\O _X} \cP ^{n'} _{X/T, (m)}  }
  \ar[u] ^-{\theta ^{\E} _{n} \otimes \mathrm{Id}}
 }
 \end{equation}
is commutative.

\item An $m$-PD-stratification relatively to $T$ on $\E$.
\end{enumerate}

An $\O _X$-linear morphism
$\phi \colon \E \to \FF$ between two left $\D ^{(m)} _{X/T}$-modules is
$\D ^{(m)} _{X/T}$-linear if and only if it commutes with the homomorphisms
$\theta _n$ (resp. $\epsilon  _n$).
\end{prop}

\begin{proof}
The proof is identical to that of
\cite[2.3.2]{Be1}.
\end{proof}

\begin{empt}
[Local description, notation]
If $X \to T$ is endowed with the finite $p$-basis $(t  _\lambda) _{\lambda =1,\dots, n}$ of level $m$ then
$ \cP ^{n} _{X/T, (m)}$ is a free $\cO _{T}$-module with
the basis 
$\{ \underline{\tau}  ^{\{ \underline{k}\} _{(m)}}~; ~ |\underline{k} | \leq n\}$, 
with $\tau _\lambda = 1 \otimes t _{\lambda} - t _{\lambda} \otimes 1$.
We get the dual basis 
$\{ \underline{\partial} ^{< \underline{k}> _{(m)}}~; ~ |\underline{k} | \leq n\}$
of 
$\cD ^{(m)} _{X/T,n}$.
For any $x\in \E$ we have the Taylor development
\begin{equation}
\label{Taylordev}
\theta ^{\E} _n (x) =
\sum _{|\underline{k}|\leq n}
\underline{\partial} ^{< \underline{k}> _{(m)}} (x) \otimes \underline{\tau} ^{\{ \underline{k}\} _{(m)} }.
\end{equation}
\end{empt}

In order to define overconvergent isocrystals in our context (see \ref{Be1-421}), 
we will need the following definition and proposition.
\begin{dfn}
\label{dfn-algcomp-mod}
Let $\B$ be a commutative $\O _X$-algebra endowed with a  structure of left $\D ^{(m)} _{X/T}$-module.
We say that the structure of left $\D ^{(m)} _{X/T}$-module on $\B$ is compatible with its
structure of $\O _X$-algebra if the isomorphisms
$\epsilon _n ^{\B}$ are isomorphisms of $\cP ^{n} _{X/T, (m)} $-algebras.
This compatibility is equivalent to the following condition : for any $f, g \in \B$ and
$\underline{k}  \in \N ^{d}$,
$$\underline{\partial} ^{< \underline{k}> _{(m)}} (fg )
=
\sum _{\underline{i} \leq \underline{k}}
\left \{
\begin{smallmatrix}
 \underline{k}\\
 \underline{i}
\end{smallmatrix}
\right\}
\underline{\partial} ^{< \underline{i}> _{(m)}} (f)
\underline{\partial} ^{< \underline{k}- \underline{i}> _{(m)}} (g).
$$
\end{dfn}

\begin{prop}
\label{prop-algBotimes}
Let $\B$ be a commutative $\O _X$-algebra endowed with a compatible structure of left $\D ^{(m)} _{X/T}$-module.
Then there exists on the tensor product
$\B \otimes _{\O _X} \D ^{(m)} _{X/T}$ a unique ring structure satisfying the following conditions
\begin{enumerate}[(a)]
\item the canonical morphisms
$\B \to \B \otimes _{\O _X} \D ^{(m)} _{X/T}$
and
$\D ^{(m)} _{X/T} \to \B \otimes _{\O _X} \D ^{(m)} _{X/T}$
are homomorphisms of sheaf of rings,
\item if $X \to T$ is endowed with the finite $p$-basis $(t  _\lambda) _{\lambda =1,\dots, n}$ of level $m$, then, for any $b \in \B$ and $\underline{k}\in \N ^{n}$, we have
$(b \otimes 1) ( 1 \otimes P) = b \otimes P$
and
$$(1 \otimes \underline{\partial} ^{< \underline{k}> _{(m)}})
(b \otimes 1)
=
\sum _{\underline{i} \leq \underline{k}}
\left \{
\begin{smallmatrix}
 \underline{k}\\
 \underline{i}
\end{smallmatrix}
\right\}
\underline{\partial} ^{< \underline{i}> _{(m)}}(b)
\otimes
\underline{\partial} ^{< \underline{k}- \underline{i}> _{(m)}}.$$

\end{enumerate}

If $\B \to \B'$ is a morphism of $\O _X$-algebras with compatible structure of  left $\D ^{(m)} _{X/T}$-modules,
then the induced
morphism
$\B \otimes _{\O _X} \D ^{(m)} _{X/T}
\to
\B '\otimes _{\O _X} \D ^{(m)} _{X/T}$
is a homomorphism of rings.
\end{prop}

\begin{proof}
We copy \cite[2.3.5]{Be1}.
\end{proof}

\begin{ntn}
[Dimension and rank of a finite $p$-basis]
\label{ntn-deltaX/T}
Let $X$ be an $S$-scheme locally of formal finite type
and having locally finite $p$-bases over $S$
(resp. a formal  $\fS$-scheme locally of formal finite type
and having locally finite $p$-bases  over $\fS$).
We set either $T= S$  or in the respective case $T= \fS$.

\begin{enumerate}[(a)]
\item The $\cO _{X}$-module $\Omega _{X/T}$ is locally free of finite rank. 
We denote by	 
$\delta ^{T} _{X}\colon X \to \N$ the locally constant function given by 
$x \mapsto \mathrm{rank} _{\O _{X,x}}   \Omega _{X/T ,x}$.
Since $X _0$ is regular, recall that $X$ is sum of its irreducible components (see \cite[6.1.10]{EGAI}).
If $U$ is an irreducible component of $X$, then
$\delta ^{T} _{X} |U$ is a  constant function. 
If moreover $U/T$ has a finite $p$-basis, then 
$\delta ^{T} _{X} | U$ is the constant function equal to the rank of $\Omega _{U/T}$,
which is equal to the number of elements of the finite $p$-basis.
When there is no ambiguity with the base $T$, 
we simply write 
$\delta _{X}$ instead of $\delta ^{T} _{X}$.

\item We get the locally constant  function
$d _{X} \colon X \to \N$, given by 
$x \mapsto \dim _x X$.

\item Let $g \colon X' \to X$ be a morphism 
of $S$-schemes locally of formal finite type
and having locally finite $p$-bases over $S$
(resp. formal  $\fS$-schemes locally of formal finite type
and having locally finite $p$-bases  over $\fS$).
We set $\delta ^{T} _{X '/X}: = \delta ^{T}_{X '} - \delta ^{T}_{X} \circ g $ 
and  $d _{X '/X}: = d _{X '} - d _{X} \circ g $.
When there is no ambiguity with the base $T$, 
we simply write 
$\delta _{X '/X}$ instead of $\delta ^{T} _{X '/X}$.

\end{enumerate}

\end{ntn}

\begin{empt}
\label{EGA17.15.4dim-delta-formula}
With notation \ref{ntn-deltaX/T}, 
suppose $X$ is integral. 
Let $x$ be a point of $X$ such that 
$\dim X =  \dim \O _{X,x}$.
We denote by
$i _x \colon \Spec k(x) \hookrightarrow X$ 
 the closed immersion induced by $x$, and by 
 $\cI$ the ideal given by this closed immersion. 
Since $k$ is perfect, then 
$\Spec k(x) \to \Spec k$ is 
formally smooth (see \ref{ex-ext-fieldpbases}). 
Following \cite[0.20.5.14.1]{EGAIV1}, we get the exact sequence of $k(x)$-vector spaces.
\begin{equation}
\label{EGA17.15.4dim-delta-formula1}
0
\to 
\cI /\cI ^2
\to 
i _x ^* \Omega _{X/\Spec k}
\to 
\Omega _{\Spec k(x)/\Spec k}
\to 0
.
\end{equation}
Since $\Spec k(x) \to \Spec k$ is locally of formal finite type, then 
$\Omega _{\Spec k(x)/\Spec k}$ is a finite dimensional $k(x)$-vector space. 
Since $\O _{X,x}$ is regular, we get 
$\dim \O _{X,x} = \dim _{k(x)}\cI /\cI ^2$.
Hence, 
\begin{equation}
\dim X +  \delta _{\Spec k(x)/\Spec k}
=
\delta _{X/\Spec k}.
\end{equation}
It might happen that 
$\delta _{\Spec k(x)/\Spec k}\not = 0$.
For instance, when $X = \Spec k ((t))$, the Krull dimension of 
$X$ is $0$ but $\Omega _{X/S}$ is $\cO _{X}$-free of rank $1$.
The function $d _{X} $ is not the right one in our context and is different from 
$\delta _{X/\Spec k} $ which behaves as fine as the dimension in the case of 
smooth formal $\fS$-schemes.

 \end{empt}

\begin{empt}
[Inverse images]
\label{inv-image-def0}
Let $f\colon X' \to X$ be a morphism of $T$-schemes (resp. formal $T$-schemes)
locally of formal finite type
and having  locally finite $p$-bases  over $T$.
Let $\E$ be a left $\D ^{(m)} _{X/T}$-module. 
Following \ref{strat-prop}, this means that $\cE$ is endowed with a PD-stratification of level $m$. 
Taking the inverse image of this PD-stratification of level $m$, we get a PD-stratification of level $m$
on $f ^* (\cE)$, i.e. $f ^* (\cE)$ is canonically endowed with structure of 
left $\D ^{(m)} _{Y/T}$-module (for more details see \cite[2.1.1]{Be2}). 
By functoriality, viewing
$\D ^{(m)} _{X/T}$ as a $\D ^{(m)} _{X/T}$-bimodule, 
we get a structure of $(\D ^{(m)} _{X'/T}, f ^{-1} \D ^{(m)} _{X/T})$-bimodule on 
$f ^* (\D ^{(m)} _{X/T})$.
We denote this bimodule by 
$\D ^{(m)} _{X'\to X/T}$.

We can extend it when the bases vary as follows : suppose we have a morphism 
$f\colon X'/T' \to X/T$ such that 
$X/T$ (resp. $X'/T'$) is locally of formal finite type
and having  locally finite $p$-bases,
then $f ^* (\D ^{(m)} _{X/T})$ is in fact 
a $(\D ^{(m)} _{X'/T'}, f ^{-1} \D ^{(m)} _{X/T})$-bimodule 
and we denote it by 
$\D ^{(m)} _{X' \to X/T'\to T}$.
The canonical morphism
\begin{equation}
\label{inv-image-def0-iso1}
f ^{*} (\cE) \to \D ^{(m)} _{X' \to X/T'\to T} \otimes _{f ^{-1} \D ^{(m)} _{X/T}} f ^{-1} \cE
\end{equation}
is an isomorphism of left 
$\D ^{(m)} _{X'/T'}$-modules.
We define the functor 
$f ^{!}
\colon 
D ^- (\D ^{(m)} _{X/T})
\to 
D ^- (\D ^{(m)} _{X'/T'})$
by setting, 
for any $\cE \in D ^- (\D ^{(m)} _{X/T})$, 
\begin{equation}
\label{inv-image-def0f!}
f ^! (\cE) 
:= 
\D ^{(m)} _{X' \to X/T'\to T} \otimes ^{\bbL }_{f ^{-1} \D ^{(m)} _{X/T}} f ^{-1} \cE
[\delta ^{T'/T} _{X '/X}],
\end{equation}
where
$\delta ^{T'/T} _{X '/X}:= \delta ^{T'}_{X '} - \delta ^{T}_{X} \circ f $.
\end{empt}

\subsection{Overconvergent singularities}

\begin{prop}
\label{Be1-421}
Put $S:= S _i$ for some integer $i \geq 0$ (resp. $S := \fS$).
Let $X$ be an $S$-scheme (resp. a formal $S$-scheme)
locally of formal finite type
and having  locally finite $p$-bases  over $S$.
Let $m,r \in \N$ be two integers such that $p ^{m+1}$ divides $r$.
Fix $f \in \Gamma (X, \O _{X})$ and put 
$\B _X(f,r) := \O _{X} [T] / (f ^{r} T-p)$.
 \begin{enumerate}[(a)]
\item Their exists on $\B _X (f,r)$ a canonical structure of $\D ^{(m)} _{X /S }$-module compatible with its structure of $\O _X$-algebra
(see \ref{dfn-algcomp-mod}). 

\item If $g \in   \Gamma (\X, \O _{\X})$, and $f ' = gf$, the homomorphism 
$$ \rho _g \colon 
\B _X(f,r)
\to 
\B _X(f',r)$$
is $\D ^{(m)} _{X /S }$-linear.
\item If $r$ is divisible by $p ^{m'+1}$ with $m' \geq m$, then the structure of 
$\D ^{(m)} _{X /S }$-module of 
$\B _X(f,r)$ is equal to that induced by 
its structure of 
$\D ^{(m')} _{X /S }$-module.
\end{enumerate}

\end{prop}

\begin{proof}
Similarly to \cite[4.2.1]{Be1}, 
by functoriality we reduce to the case where 
$S = \Spec \Z _{(p)}$ and $X =  \Spec \Z _{(p)} [t]$, 
$f =t$. This latter case is already proved in \cite[4.2.1]{Be1}.
\end{proof}

\begin{prop}
\label{Be1-422}
With the hypotheses \ref{Be1-421}, let $\I \subset \O _X$ be an $m$-PD-nilpotent quasi-coherent ideal,
$f,g \in \Gamma (X, \O _X)$,
$h \in \Gamma (X, \I )$,
and $f' = gf +h$. 
There exists the  canonical $\D _{X  / S } ^{(m)}$-linear homomorphism of $\O _X$-algebras
$$\eta _{g,h}
\colon 
\B _X (f,r) 
\to 
\B _X (f',r),$$
satisfying the following properties:
\begin{enumerate}[(a)]
\item 
If $g ' \in \Gamma (X, \O _X)$, $h ' \in \Gamma (X, \I \O _X)$,
and 
$f'' = g' f ' +h'$,
$g'' = g' g$,
$h'' = g' h +h'$, then
$\eta _{g'',h''}=
\eta _{g',h'}
\circ
\eta _{g,h}$.

\item $\eta _{g,0} = \rho _g$, $\eta _{1,0} = \mathrm{Id}$.

\item If $f $ is not a divisor of $0$ in $\O _X / \I \O _X$,
$\eta _{g,h}$ only depend on $f$ and $f'$.
\item If $r$ is divisible by $p ^{m' +1}$, with $m' \geq m$, 
$\eta _{g,h}$ is independ on $m \leq m'$.

\end{enumerate}
 
\end{prop}

\begin{proof}
This is checked similarly to \cite[4.2.2]{Be1}.
For the reader, we will only  recall below the construction of $\eta _{g,h}$.

1) Suppose $g=1$.
Let $u\colon Z  \hookrightarrow  X $ be the closed immersion defined by 
$\I $.
Put $S _0= \Spec \Z _{(p)}$ and $X _0=  \Spec \Z _{(p)} [t]$.
Let $\overline{f}$ and $\overline{f'}$ be the image of $f$ and $f'$ via the morphism 
$\Gamma ( X, \cO _X) 
\to 
\Gamma ( Z, \cO _Z) $
induced by $u$.
Since $\overline{f}=\overline{f'}$,  then $f$ an $f'$ induce the same morphism
$\overline{f}\colon Z  \to X _0$. 
Since $\I$ is an $m$-PD-nilpotent ideal, then by using the universal property of the $m$-PD-enveloppe,
for any integer $n$ large enough, we get a unique factorization
$\theta \colon X \to \Delta ^{n} _{X_0/S _0, (m)}$ making commutative the following diagram
\begin{equation}
\xymatrix{
{Z } 
\ar@{^{(}->}[rr] ^-{}
\ar[d] ^-{\overline{f}}
&
& 
{X } 
\ar[d] ^-{f '\times f}
\ar@{.>}[ld] _-{\theta}
\\ 
{X _0} 
\ar[r] ^-{}
& 
{\Delta ^{n} _{X_0/S _0, (m)}}
\ar[r] ^-{}
&
{X _0 \times _{S _0} X _0} 
}
\end{equation}
Let $\epsilon _n \colon 
\mathcal{P} ^{n} _{X _0/S _0, (m)}
\otimes _{\O _{X _0}}
\B _{X _0} (t,r) 
\riso 
\B _{X _0} (t,r) 
\otimes _{\O _{X _0}}
\mathcal{P} ^{n} _{X _0/S  _0, (m)} $
be the isomorphism given by the $\D ^{(m)} _{X _0/S _0}$-module structure of 
$\B _{X _0} (t,r) $.
Taking the inverse image by $\theta$ we get the isomorphism
$\epsilon _h \colon 
\B _{X } (f,r) 
\riso 
\B _{X} (f',r) $.

 2) In general, $\eta _{g,h}:= \epsilon _h \circ \rho _g$. 
\end{proof}

\begin{ntn}
\label{ntnODdagZ}
Let $ \X$ be a formal  $\fS$-scheme
locally of formal finite type
and having locally finite $p$-bases  over $\fS$.
Let $Z$ be a divisor of $X$.
Let $\U $ be an open set of $\X $, $f \in \Gamma (\U, \O _{\X})$ such that 
the closed immersion $Z \cap U \hookrightarrow U$ is given by $\overline{f} \in \Gamma (\U, \O _{X})$ the image of 
$f$ via $\Gamma (\U, \O _{\X}) \to \Gamma (\U, \O _{X})$.
Following \ref{Be1-422}, 
 $\B _{U _i} (f,r) $ (resp.  $\B _{\U} (f,r) $) only depends on $Z$. 
 Hence, glueing  $\B _{U _i} (f,r) $ (resp.  $\B _{\U} (f,r) $) we get the $\O _{X _i}$-algebra
 (resp. $\O _{\X}$-algebra)
 $\B _{X _i} (Z,r) $ (resp.  $\B _{\X} (Z,r) $).
 Put
$\B _{X _i}^{(m)}   (Z)   := \B _{X _i} (Z,p ^{m+1})  $,
$\B ^{(m)}  _{\X} (Z)  := \B _{\X} (Z,p ^{m+1})  $,
and 
$\widehat{\B} ^{(m)} _{\X} (Z) :=\underleftarrow{\lim} _i \B  ^{(m)} _{X _i} (Z)  $,
the $p$-adic completion of 
$\B ^{(m)}  _{\X} (Z) $.
Finally, we set
$$\O _{\X} (\hdag Z) := \underleftarrow{\lim} _m \widehat{\B}  ^{(m)} _{\X} (Z)  ,
\D ^\dag _{\X /\fS } (\hdag Z) := 
\underleftarrow{\lim} _m \widehat{\B}  ^{(m)} _{\X} (Z)  \widehat{\otimes} _{\O _\X} 
\widehat{\D} ^{(m)} _{\X /\fS } .$$

If $Z \subset T$ are two divisors of $X$, we get from \ref{Be1-422} the canonical morphisms
$\B _{X _i}^{(m)}   (Z)   \to \B _{X _i}^{(m)}   (T)  $,
 $\B ^{(m)}  _{\X} (Z) \to \B ^{(m)}  _{\X} (T)$,
 and 
$\O _{\X} (\hdag Z) 
\to 
\O _{\X} (\hdag T)$. 
\end{ntn}

\begin{thm}
\label{Be1-4310&12}
Let $ \X$ be a formal  $\fS$-scheme
locally of formal finite type
and having locally finite $p$-bases  over $\fS$.
Let $Z$ be a divisor of $X$.
Let $\Y $ be the open subset of $\X $ complementary to the support of $Z$,
and $j \colon \Y  \hookrightarrow \X $ be the open immersion.
\begin{enumerate}[(a)]
\item The homomorphisms 
 $\O _{\X} (\hdag Z)  _\Q \to j _*  \O _{\Y} (\hdag Z)  _\Q $
 and 
 $\D ^\dag _{\X /\fS } (\hdag Z) _\Q
 \to 
j _*  \D ^\dag _{\Y /\fS ,\Q} $
are faithfully flat. 
\item For any coherent
$\D ^\dag _{\X /\fS } (\hdag Z) _\Q$-module $\E$, 
the canonical homomorphism
$$
j _*  \D ^\dag _{\Y /\fS ,\Q} 
\otimes _{\D ^\dag _{\X /\fS } (\hdag Z) _\Q}
\E
\to 
j _* j ^* \E$$
is an isomorphism.
\item 
\label{Be1-4310&12-item3}
A coherent
$\D ^\dag _{\X /\fS } (\hdag Z) _\Q$-module $\E$ is null if and only if 
$j ^* \E$ is null.
\end{enumerate}
\end{thm}

\begin{proof}
We can follow the proof of \cite[4.3.10 and 4.3.12]{Be1}.
\end{proof}

\begin{prop}
\label{432-442Be1}
We keep notation \ref{ntn-jm}.
\begin{enumerate}[(a)]
\item There exist canonical isomorphisms of $\O _{\X}$-algebras
\begin{gather}
\widehat{\B} ^{(m)} _{\X} (Z) _{\Q}
\riso
\sp _* j _{m *} j ^{*} _m \O _{\X _K}
\\
\O _{\X} (\hdag Z) _\Q \riso \sp _* j ^\dag  \O _{\X _K}.
\end{gather}

\item For any affine open formal subscheme $\U \subset \X$, 
$\Gamma ( \U, \widehat{\B} ^{(m)} _{\X} (Z))$,
and 
$\Gamma ( \U, \widehat{\B} ^{(m)} _{\X} (Z) _\Q)$
are noetherian. The extensions 
$\O _{\X,\Q}
\to 
\widehat{\B} ^{(m)} _{\X} (Z) _\Q$
and 
$\widehat{\B} ^{(m)} _{\X} (Z) _\Q
\to 
\widehat{\B} ^{(m+1)} _{\X} (Z) _\Q $
are flat. 

The sheaves $\widehat{\B} ^{(m)} _{\X} (Z) $, 
$\widehat{\B} ^{(m)} _{\X} (Z) _\Q$,
and
$\O _{\X} (\hdag Z) _\Q$ 
are coherent. Moreover, coherent modules over these sheaves 
satisfy theorems $A$ and $B$.

\end{enumerate}

\end{prop}

\begin{proof}
We can copy the proof of 
\cite[4.3.2]{Be1}.
\end{proof}

\subsection{PD-costratification of level $m$ and right $\cD ^{(m)} _{X/T}$-modules}
Let $T$ be an ${S _i}$-scheme of finite type (resp. a formal $\fS$-scheme of finite type)
for some integer $i \geq 0$.
Let $X$ be an $T$-scheme (resp. a formal $T$-scheme)
locally of formal finite type
and having  locally finite $p$-bases  over $T$.
Similarly to \cite[1.1]{Be2}, we define the notion of 
$m$-PD-{\it costratifications relatively to $T$} on an 
$\O _X$-module $\M$.
\begin{dfn}
\label{defincostrat}
  Let  $\M$ be an  $\O _X$-module. An $m$-PD-{\it costratification on $\M$ relatively to $T$}
 is the data of a compatible family of  $\cP ^n _{X/T (m)}$-linear isomorphisms
  $$ \epsilon _n \ :\
  \mathcal{H} om _{\O _X} ( p ^{n} _{0*} \cP ^n _{X/T (m)},\,\M)
  \riso
  \mathcal{H} om _{\O _X} ( p ^{n} _{1*} \cP ^n _{X/T (m)},\,\M),$$
  satisfying the following conditions:
  \begin{enumerate}[(a)]
    \item $\epsilon _0 = Id _{\M}$ ;
    \item For any $n$, $n'$, the diagram 
    \begin{equation}
      \label{definstratdiag1}
      \xymatrix@C=-2,5cm {
 {  \mathcal{H} om _{\O _X}
  ( d ^{n,n'} _{0*}( \cP ^n _{X/T (m)} \otimes _{\O _X} \cP ^{n'} _{X/T (m)}),\,\M)}
  \ar[rr]^{\delta ^{n,n'\flat} _{(m)}(\epsilon _{n +n'})}
  _{\widetilde{{}\hspace{0,3cm}  }}
   \ar[rd]^{q _0 ^{n,n'\flat}(\epsilon _{n +n'})}
   _{\widetilde{{}\hspace{0,3cm}  }}
   &&
 {  \mathcal{H} om _{\O _X}
  ( d ^{n,n'} _{2*}( \cP ^n _{X/T (m)} \otimes _{\O _X} \cP ^{n'} _{X/T (m)}),\,\M)}
\\
&
 {  \mathcal{H} om _{\O _X}
  ( d ^{n,n'} _{1*} (\cP ^n _{X/T (m)} \otimes _{\O _X} \cP ^{n'} _{X/T (m)}),\,\M)}
  \ar[ur] ^{q _1 ^{n,n'\flat}(\epsilon _{n +n'})}
  _{\widetilde{{}\hspace{0,3cm}  }}
  }
    \end{equation}
is commutative. 
  \end{enumerate}
 This latter condition is equivalent to the following one: for any $n\in \N$, 
 the diagram
    \begin{equation}
      \label{definstratdiag1bis}
      \xymatrix@C=-1cm {
 {  \mathcal{H} om _{\O _X}
  ( p ^n _{0*}( \cP ^n _{X/T (m)} (2) ),\,\M)}
  \ar[rr]^{p ^{n\flat} _{02}(\epsilon _{n})}
  _{\widetilde{{}\hspace{0,3cm}  }}
   \ar[rd]^{p ^{n\flat} _{01}(\epsilon _{n})}
   _{\widetilde{{}\hspace{0,3cm}  }}
   &&
 {  \mathcal{H} om _{\O _X}
  ( p ^{n} _{2*}( \cP ^n _{X/T (m)} (2) ),\,\M)}
\\
&
 {  \mathcal{H} om _{\O _X}
  ( p ^{n} _{1*} (\cP ^n _{X/T (m)} (2) ),\,\M)}
  \ar[ur] ^{p ^{n\flat} _{12}(\epsilon _{n})}
  _{\widetilde{{}\hspace{0,3cm}  }}
  }
    \end{equation}
is commutative.  
\end{dfn}

\begin{prop}
\label{prop-rightD-costrat}
Let   $\M$ be an $\O _X$-module. 
The following data are equivalent :

\begin{enumerate}[(a)]
\item  A structure of right  $\D _{X/T} ^{(m)}$-module on $\M$
extending its structure de $\O _X$-module ;

\item An $m$-PD-costratification $( \epsilon _n ^{\M})$ relatively to $T$ on $\M$.

\end{enumerate}

An $\O _X$-linear homomorphism between two 
right $\D _{X/T} ^{(m)}$-modules is
$\D _{X/T} ^{(m)}$-linear if and only if it commutes with the  isomorphisms
$\epsilon _n ^{\M}$.

\end{prop}

\begin{proof}
We can copy word by word the proof of 
\cite[1.1.4]{Be2}.
\end{proof}

\begin{lem}
Let $Y$ be an $S$-scheme
locally of formal finite type
and having  locally finite $p$-bases  over $S$.
There exists a canonical structure of right $\D ^{(m)} _{Y/S} $-module on $\omega _{Y /S}$.
It is characterized by the following local formula:
suppose that $Y$ is endowed with a finite $p$-basis $(b  _{\lambda}) _{\lambda =1,\dots, n}$.
Then, for any differential operator $P \in \D ^{(m)} _{Y /S}$ and $a \in \O _{Y }$ we have
\begin{equation}
\label{coro-rightD-mod-omegaformal}
(a \ d\, b_1 \wedge \cdots \wedge d\, b_n) \cdot P := 
{} ^t 
P (a) \ d\, b_1 \wedge \cdots \wedge d\, b_n.
\end{equation}

\end{lem}

\begin{proof}
By canonicity, this is local in $Y$. 
Hence, we can suppose $Y$ 
is endowed with a finite $p$-basis $(b  _{\lambda}) _{\lambda =1,\dots, n}$.
Following 
\ref{lifting-pbasis}, 
there exists a (unique up to isomorphisms) 
formal $\fS$-scheme locally of formal finite type $\Y$ having finite $p$-basis
and such that $\fY \times _{\Spf \V} S  \riso Y$.
We conclude using \cite[3.6.3]{Caro-Vauclair}.
\end{proof}

\begin{empt}
\label{lem-rightD-petale}
Let $f \colon X \to Y$ be a $p$-étale morphism of 
$T$-schemes (resp. a formal $T$-schemes) locally of formal finite type
and having  locally finite $p$-bases  over $T$.
Let $\cM $ be right $\D ^{(m)} _{Y/T}$-module. 
Let 
$\epsilon _n \ :\
  \mathcal{H} om _{\O _Y} ( p ^{n} _{0*} \cP ^n _{Y/T,\,(m)},\,\M)
  \riso
  \mathcal{H} om _{\O _Y} ( p ^{n} _{1*} \cP ^n _{Y/T,\,(m)},\,\M)$ 
  be the corresponding $m$-PD-costratification relatively to $T$. 
For $i = 0,1$, 
we have the isomorphism of $\O _X$-algebras
$f ^* p ^{n} _{i*} \cP ^n _{Y/T,\,(m)}
\riso 
p ^{n} _{i*} \cP ^n _{X/T (m)}$.
Hence, by applying the functor $f ^*$ to $\epsilon _n$ 
we get the isomorphism
$f ^* (\epsilon _n)
\colon 
\mathcal{H} om _{\O _Y} ( p ^{n} _{0*} \cP ^n _{Y/T,\,(m)},\, f ^*\M)
 \riso
\mathcal{H} om _{\O _Y} ( p ^{n} _{1*} \cP ^n _{Y/T,\,(m)},\, f^* \M)$.
We check that 
$f ^* (\epsilon _n)$ is an $m$-PD-costratification of $f^* \M$, i.e. 
$f^* \M$ is endowed with a canonical structure of 
right $\D ^{(m)} _{X/T}$-module.
Moreover, the canonical morphism
\begin{equation}
\notag
f^* \M
\to 
f ^{-1} \cM
\otimes _{f ^{-1}\D ^{(m)} _{Y/T}}
\D ^{(m)} _{X/T}
\end{equation}
is an isomorphism of right $\D ^{(m)} _{X/T}$-modules.
Moreover, the canonical isomophism
\begin{equation}
\label{lem-rightD-petale-omega}
f ^*  (\omega _{Y/T})
\riso
\omega _{X/T}
\end{equation}
is an isomorphism of right $\D ^{(m)} _{X/T}$-modules.
\end{empt}

\begin{empt}
[Inverse images]
\label{inv-image-def0-rightGen}
Let $f\colon X'/T' \to X/T$ be a morphism such that 
$X/T$ (resp. $X'/T'$) is locally of formal finite type
and has  locally finite $p$-bases.
Viewing
$\D ^{(m)} _{X/T}$ as a
$\D ^{(m)} _{X/T}$-bimodule, 
we get a structure of left
$\D ^{(m)} _{X/T}$-bimodule on 
$\D _{X  / T  } ^{(m)}\otimes _{\O _{X }} \omega ^{-1} _{X / T}$ (see \cite[3.6.4]{Caro-Vauclair}).
By functoriality, we get a structure of 
left $(f ^{-1} \D ^{(m)} _{X/T},\D ^{(m)} _{X'/T'})$-bimodule on 
$f ^* _l \left ( \D _{X  / T  } ^{(m)}\otimes _{\O _{X }} \omega ^{-1} _{X / T} \right)$,
where ``$l$'' means that we choose the left structure of 
left $\D ^{(m)} _{X/T}$-module to get a structure of 
left $\D ^{(m)} _{X'/T'}$-module
on 
$f ^* _l \left ( \D _{X  / T  } ^{(m)}\otimes _{\O _{X }} \omega ^{-1} _{X / T} \right)$.
We get the
$(f ^{-1} \D ^{(m)} _{X/T},\D ^{(m)} _{X'/T'})$-bimodule
$\cD ^{(m)}_{X \leftarrow X'/T\leftarrow T'}
:=
\omega _{X ^{\prime } /T'} \otimes _{\O _{X ' }}f ^* _l \left ( \D _{X  / T  } ^{(m)}\otimes _{\O _{X }} \omega ^{-1} _{X / T} \right)$.
We define the functor 
$f ^{!}
\colon 
D ^- ({} ^{\mathrm{r}} \D ^{(m)} _{X/T})
\to 
D ^- ({} ^{\mathrm{r}}\D ^{(m)} _{X'/T'})$
by setting, 
for any $\cM \in D ^- ({} ^{\mathrm{r}} \D ^{(m)} _{X/T})$, 
\begin{equation}
\label{inv-image-def0-rightGenf!}
f ^! (\cM) 
:= 
f ^{-1} \cM
\otimes ^{\bbL }_{f ^{-1} \D ^{(m)} _{X/T}} 
\cD ^{(m)}_{X \leftarrow X'/T\leftarrow T'}
[\delta ^{T'/T} _{X '/X}],
\end{equation}
where
$\delta ^{T'/T} _{X '/X}:= \delta ^{T'}_{X '} - \delta ^{T}_{X} \circ f $.

With notation \ref{inv-image-def0f!} and \ref{inv-image-def0-rightGenf!},
for any $\cE \in D ^{-} ( \widetilde{\D} _X )$
we have the canonical isomorphisms
\begin{gather}
\notag
f ^{!} ( \omega _{X/T} \otimes _{\O _X} \E)
=
f ^{-1} \left ( \omega _{X/T} \otimes _{\O _X} \E \right )
\otimes ^{\bbL }_{f ^{-1} \D ^{(m)} _{X/T}} 
\left (\omega _{X ^{\prime } /T'} \otimes _{\O _{X ' }}f ^* _l \left ( \D _{X  / T  } ^{(m)}\otimes _{\O _{X }} \omega ^{-1} _{X / T} \right) \right)
[\delta ^{T'/T} _{X '/X}]
\\
\label{fund-isom2bisprepre}
\riso 
\left (\omega _{X ^{\prime } /T'} \otimes _{\O _{X ' }}f ^*  \D _{X  / T  } ^{(m)}  \right)
\otimes ^{\bbL }_{f ^{-1} \D ^{(m)} _{X/T}} 
f ^{-1} \left ( \E \right )
[\delta ^{T'/T} _{X '/X}]
\riso 
\omega _{X'/T'} \otimes _{\O _{X'}} f ^{!} (\E ) .
\end{gather}

\end{empt}

\begin{empt} 
\label{inv-image-def0-right}
Let $f\colon X' \to X$ be a finite morphism of $T$-schemes (resp. a formal $T$-schemes)
locally of formal finite type
and having  locally finite $p$-bases  over $T$.
Let $\cM$ be a right $\D ^{(m)} _{X/T}$-module. 
Following \ref{prop-rightD-costrat}, this means that $\cM$ is endowed with a PD-costratification of level $m$. 
Since $f$ is finite, then following \cite[1.1.1]{Be2} (which is some kind of survey of \cite[III.6]{HaRD}),
the functor $f ^{\flat}$ is defined by setting
 $$f ^{\flat} (\cM)
:= 
\overline{f} ^* \bbR\cH om _{\cO _{X }} ( f _* \cO _{X'}, \cM ').$$
Similarly to \ref{inv-image-def0} (we just have to replace functors of the form $f ^*$
by functors of the forms $f ^\flat$), 
by applying the functors of the form $f ^\flat$ to the PD-costratification of level $m$ of $\cM$, 
we get a structure of PD-costratification of level $m$ on 
$f ^{\flat} (\cM)$,
 i.e. $f ^{\flat} (\cM)$ is canonically endowed with a structure of 
right $\D ^{(m)} _{X'/T}$-module (we copy word by word  \cite[2.1.1]{Be2}). 
By functoriality, viewing
$\D ^{(m)} _{X/T}$ as a $\D ^{(m)} _{X/T}$-bimodule, 
we get a structure of $(f ^{-1} \D ^{(m)} _{X/T}, \D ^{(m)} _{X'/T})$-bimodule on 
$f ^\flat (\D ^{(m)} _{X/T})$.

We can extend it when the bases vary as follows : let  
$f\colon X'/T' \to X/T$ be a morphism such that 
$X/T$ (resp. $X'/T'$) is locally of formal finite type
and having  locally finite $p$-bases,
$X' \to X$ is a finite morphism.
Then $f ^\flat (\D ^{(m)} _{X/T})$ is in fact 
a $(f ^{-1} \D ^{(m)} _{X/T}, \D ^{(m)} _{X'/T'})$-bimodule.
\end{empt}

\subsection{Extraordinary inverse image, direct image :algebraic case}

\label{subsec-/Svs/T}

Let $i \geq 0$ be some integer, 
$S$ be a $\Spec \V / \pi ^{i+1} \V$-scheme of finite type. 
Let $T$ be an $S$-scheme locally of formal finite type 
and having locally finite $p$-bases over $S$.
Let $h \colon X  \to Y $ be a morphism of 
$T$-schemes locally of formal finite type over $S$
and having locally finite $p$-bases over $T$.
We denote by $g \colon Y  \to T$ and $f \colon X  \to T$ the structural morphisms.

\begin{empt}
\label{empt-DeltaY/T2Y}
With notation \cite[2.2.2 and 2.2.4]{Caro-Vauclair},
using the universal properties of $m$-PD-envelops , 
we get the commutative diagram
\begin{equation}
\label{DeltaY/T2Y}
\xymatrix{
{\Delta ^{n} _{Y  / T , (m)}(2)} 
\ar@{.>}[d] ^-{}
\ar@<1ex>[r]
\ar[r]
\ar@<-1ex>[r]
& 
{\Delta ^{n} _{Y  / T , (m)}} 
\ar@{.>}[d] ^-{}
\ar@<0.5ex>[r]
\ar@<-0.5ex>[r]
&
{id _{Y  }}
\ar@{=}[d] ^-{}
\\ 
{\Delta ^{n} _{Y /S, (m)}(2)} 
\ar@<1ex>[r]
\ar[r]
\ar@<-1ex>[r]
& 
{\Delta ^{n} _{Y /S , (m)}} 
\ar@<0.5ex>[r]
\ar@<-0.5ex>[r]
&
{id _{Y  }.}
}
\end{equation}
This yields that
we get the homomorphisms of rings
$\mathcal{P} ^n _{Y /S , (m)}
\to 
\mathcal{P} ^n _{Y  / T , (m)}$. 
By duality, this yields the homomorphism of $\O_Y$-modules
$\D ^{(m)} _{Y  / T , n}
\to 
\D ^{(m)} _{Y /S , n}$.
Using the commutativity of the diagram
\ref{DeltaY/T2Y}, by definition of their ring structures,
we can check the induced homomorphism 
of $\cO _Y$-modules
$\D ^{(m)} _{Y  / T }
\to 
\D ^{(m)} _{Y/S}$
is in fact a morphism of rings.

Let $\cB  _T$ be an $\cO _T$-algebra endowed with a compatible structure of 
left $\D ^{(m)} _{T /S}$-module. 
Set 
$\widetilde{\D} ^{(m)} _{T / S}
:= 
\cB  _T \otimes _{\cO _T} \D ^{(m)} _{T / S }$,
and for any $n \in \N$, 
$\widetilde{\D} ^{(m)} _{T / S ,n}
:= 
\cB  _T \otimes _{\cO _T} \D ^{(m)} _{T / S ,n}$,
$\widetilde{\cP}^{n} _{T/S,(m)}:=
\cB  _T \otimes _{\cO _T} \cP ^{n} _{T/S,(m)}$.

Let  $\cB  _Y$ be a $g ^{*} (\cB _T)$-algebra which is 
endowed with a compatible structure of 
left $\D ^{(m)} _{Y / S }$-module. 
Since 
$\D ^{(m)} _{Y  / T }
\to 
\D ^{(m)} _{Y/S}$
is in fact a morphism of rings, then  
$\cB  _Y$ is also an $g ^{*} (\cB _T)$-algebra which is endowed with a compatible structure of 
left $\D ^{(m)} _{Y /T}$-module. 
Set 
$\widetilde{\D} ^{(m)} _{Y / S}
:= 
\cB  _Y \otimes _{\cO _Y} \D ^{(m)} _{Y / S }$,
and for any $n \in \N$, 
$\widetilde{\D} ^{(m)} _{Y / S ,n}
:= 
\cB  _Y \otimes _{\cO _Y} \D ^{(m)} _{Y / S ,n}$,
$\widetilde{\cP}^{n} _{Y/S,(m)}
:=
\cB  _Y \otimes _{\cO _Y} \cP ^{n} _{Y/S,(m)}$.
Set $\widetilde{\D} ^{(m)} _{Y / T }
:= 
\cB  _Y \otimes _{\cO _Y} \D ^{(m)} _{Y / T }$,
and for any $n \in \N$, 
$\widetilde{\D} ^{(m)} _{Y / T ,n}
:= 
\cB  _Y \otimes _{\cO _Y} \D ^{(m)} _{Y / T ,n}$,
$\widetilde{\cP}^{n} _{Y/T,(m)}
:=
\cB  _Y \otimes _{\cO _Y} \cP  ^{n} _{Y/T,(m)}$.

Similarly, let 
$\cB  _X$ be a $f ^{*} (\cB _T)$-algebra which is endowed with a compatible structure of 
left $\D ^{(m)} _{X / S }$-module.
Set 
$\widetilde{\D} ^{(m)} _{X / S}
:= 
\cB  _X \otimes _{\cO _X} \D ^{(m)} _{X / S }$,
and for any $n \in \N$, 
$\widetilde{\D} ^{(m)} _{X / S ,n}
:= 
\cB  _X \otimes _{\cO _X} \D ^{(m)} _{X / S ,n}$,
$\widetilde{\cP}^{n} _{X/S,(m)}
:=
\cB  _Y \otimes _{\cO _X} \cP  ^{n} _{X/S,(m)}$.
Set $\widetilde{\D} ^{(m)} _{X / T }
:= 
\cB  _X \otimes _{\cO _X} \D ^{(m)} _{X / T }$,
and for any $n \in \N$, 
$\widetilde{\D} ^{(m)} _{X / T ,n}
:= 
\cB  _X \otimes _{\cO _X} \D ^{(m)} _{X / T ,n}$,
$\widetilde{\cP}^{n} _{X/T,(m)}
:=
\cB  _Y \otimes _{\cO _X} \cP  ^{n} _{X/T,(m)}$.

We denote by $\mathrm{oub} _{Y/T/S}$ the 
forgetful functor (via the canonical morphism 
$\widetilde{\cD} ^{(m)} _{Y  / T } \to \widetilde{\cD} ^{(m)} _{Y/S  }$)
from the category of 
left (resp. right) $ \widetilde{\cD} ^{(m)} _{Y /S }$-modules
to that of 
left (resp. right) $ \widetilde{\cD} ^{(m)} _{Y  / T }$-modules ; 
and similarly by replacing $Y$ by $X$.

Let us denote by 
$\widetilde{h}  ^* _{/S}:= \cB _{X} \otimes _{h ^{-1} \cB _Y} h ^{-1}(-)$ 
the functor  from the category of 
left $ \widetilde{\cD} ^{(m)} _{Y/S}$-modules
to that of left  $\widetilde{\cD} ^{(m)} _{X/S} $-modules
and 
by 
$\widetilde{h}  ^* _{/T}:= \cB _{X} \otimes _{h ^{-1} \cB _Y} h ^{-1}(-)$ 
the functor  from the category of 
left $ \widetilde{\cD} ^{(m)} _{Y/T}$-modules
to that of left  $\widetilde{\cD} ^{(m)} _{X/T} $-modules
From the commutative diagram \ref{DeltaY/T2Y} (still valid with some tildes), 
we get the commutation 
\begin{equation}
\label{pre-u!=u!/T(m)-0}
\mathrm{oub} _{X/T/S}  \circ h ^* _{/T}
\riso 
h ^* _{/S} \circ \mathrm{oub} _{Y/T/S}.
\end{equation}

By functoriality, 
we get the $(\widetilde{\cD} ^{(m)} _{X / S }, h ^{-1}\widetilde{\cD} ^{(m)} _{Y/S})$-bimodule
$\widetilde{\cD} ^{(m)} _{X \to Y  /S} : = 
h ^* _{/S} \widetilde{\cD} ^{(m)} _{Y/S} $
and the $(\widetilde{\cD} ^{(m)} _{X / T}, h ^{-1}\widetilde{\cD} ^{(m)} _{Y/T})$-bimodule
$\widetilde{\cD} ^{(m)} _{X \to Y  /T} : = 
h ^* _{/T} \widetilde{\cD} ^{(m)} _{Y/T} $
\end{empt}

\begin{lemm}
We have the isomorphism of 
$(\widetilde{\cD} ^{(m)} _{X/ T }, h ^{-1}\widetilde{\cD} ^{(m)} _{Y/S})$-bimodules 
\begin{equation}
\label{pre-u!=u!/T(m)}
\widetilde{\cD} ^{(m)} _{X \to Y  / T }
\otimes _{h ^{-1}\widetilde{\cD} ^{(m)} _{Y  / T}} h ^{-1}\widetilde{\cD} ^{(m)} _{Y /S}
\riso 
\widetilde{\cD} ^{(m)} _{X \to Y/S}.
\end{equation}

\end{lemm}

\begin{proof}
By functoriality, 
this is a consequence of 
\ref{pre-u!=u!/T(m)-0}.
\end{proof}

\begin{ntn}
[Local description of $\widetilde{\cD} ^{(m)} _{Y/T}$]
\label{ntn-loc-coor-TvsS}
Suppose 
$Y/T$ has the finite $p$-basis 
$t _1,\cdots, t _d$.
We set $\tau _i := 1 \otimes t _i - t _i \otimes 1 \in \cO _{Y \times _T Y}$ for any 
$i= 1,\cdots, d$. 
For any $\underline{i} = (i _1, \cdots, i _d) \in \bbN ^d$, 
let
$\underline{\tau}  {} ^{\{\underline{i}\} _{(m)}}:= 
\tau _1 ^{\{i _1\} _{(m)}}
\cdots
\tau _d ^{\{i _d\} _{(m)}}\in \widetilde{\cP}^{n} _{Y/T,(m)}$.
We get the basis of the free $\cB _Y$-module
$\widetilde{\cP}^{n} _{Y/T,(m)}$ given by 
$\underline{\tau} ^{\{ \underline{i}\} _{(m)}}$, with $| \underline{i}| \leq n$.
By taking the dual basis and taking the inverse limits,  we get a basis on 
the free (for the left or right structure) $\cB _Y$-module $\widetilde{\cD} ^{(m)} _{Y/T}$
(for its left structure this is by definition but this is also the case its right structure) 
which is  denoted by
$\{ 
\underline{\partial}  ^{<\underline{i}> _{(m)}}
\ | \ 
\underline{i} \in \N ^d
\} $.
Hence, a section of 
$\widetilde{\cD} ^{(m)} _{Y/T}$ can uniquely be written as a {\it finite} sum of the form
$\sum  _{\underline{i} \in \N ^d} 
a _{\underline{i}} 
\underline{\partial}  ^{<\underline{i}> _{(m)}}$
(resp. 
$\sum  _{\underline{i} \in \N ^d} 
\underline{\partial}  ^{<\underline{i}> _{(m)}}
 a _{\underline{i}}$)
with 
$a _{\underline{i}} \in \cB _Y$.

\end{ntn}

\begin{ntn}
[Local description of $\widetilde{\cD} ^{(m)} _{T/S}$]
\label{ntn-loc-coor-TvsSbis}
Suppose $T/S$ has the finite $p$-basis 
$\tilde{t} _1,\cdots, \tilde{t} _r$.
We set $\widetilde{\tau} _j := 1 \otimes \widetilde{t} _j - \widetilde{t} _j \otimes 1 \in \cO _{T \times _S T}$ for any 
$j= 1,\cdots, r$. 
For any $\underline{j} = (j _1, \cdots, j _r) \in \bbN ^r$, 
let 
$\underline{\widetilde{\tau}}  {} ^{\{\underline{j}\} _{(m)}}:= 
\widetilde{\tau}  _1 ^{\{j _1\} _{(m)}}
\cdots
\widetilde{\tau}  _r ^{\{j _r\} _{(m)}}\in \widetilde{\cP}^{n} _{T/S,(m)}$.
We get the basis of the free $\cB _T$-module
$\widetilde{\cP}^{n} _{T/S,(m)}$ given by 
$\underline{\widetilde{\tau}}  {} ^{\{ \underline{j}\} _{(m)}}$, with $| \underline{j}| \leq n$.
By taking the dual basis and taking the inverse limits,  we get a basis on 
the free  (for the left or right structure) $\cB _T$-module $\widetilde{\cD} ^{(m)} _{T/S}$ 
which is  denoted by
$\{ 
\underline{\widetilde{\partial}}  {}  ^{<\underline{j}> _{(m)}}
\ | \ 
\underline{j} \in \N ^r
\} $.
Hence, a section of 
$\widetilde{\cD} ^{(m)} _{T/S}$ can uniquely be written as a {\it finite} sum of the form
$\sum  _{\underline{j} \in \N ^r} 
a _{\underline{j}} 
\underline{\widetilde{\partial}}  {}  ^{<\underline{j}> _{(m)}}$
(resp. 
$\sum  _{\underline{j} \in \N ^r} 
\underline{\widetilde{\partial}}  {}  ^{<\underline{j}> _{(m)}}
 a _{\underline{j}}$)
with 
$a _{\underline{j}} \in \cB _T$.

\end{ntn}

\begin{ntn}
[Local description of $\widetilde{\cD} ^{(m)} _{Y/S}$]
\label{ntn-loc-coor-TvsS2}
Suppose $T/S$ has the finite $p$-basis 
$\tilde{t} _1,\cdots, \tilde{t} _r$.
and suppose moreover that
$Y/T$ has the finite $p$-basis 
$t _1,\cdots, t _d$.
By abuse of notation, 
we denote by 
$\tilde{t} _1,\cdots, \tilde{t} _r$ the element of $\Gamma (Y, \cO _Y)$
induced by 
$\tilde{t} _1,\cdots, \tilde{t} _r$
via $g$.
We get the finite $p$-basis
$\tilde{t} _1,\cdots, \tilde{t} _r, t _1,\cdots, t _d$
of $Y/S$.
We set $\tau _i := 1 \otimes t _i - t _i \otimes 1 \in \cO _{Y \times _T Y}$ for any 
$i= 1,\cdots, d$ ; 
$\widetilde{\tau} _j := 1 \otimes \widetilde{t} _j - \widetilde{t} _j \otimes 1 \in \cO _{Y \times _T Y}$ for any 
$j= 1,\cdots, r$.  
For any $\underline{i} = (i _1, \cdots, i _d) \in \bbN ^d$, 
let
$\underline{\tau}  {} ^{\{\underline{i}\} _{(m)}}:= 
\tau _1 ^{\{i _1\} _{(m)}}
\cdots
\tau _d ^{\{i _d\} _{(m)}}$ ; 
for any $\underline{j} = (j _1, \cdots, j _r) \in \bbN ^r$, 
let 
$\underline{\widetilde{\tau}}  {} ^{\{\underline{j}\} _{(m)}}:= 
\widetilde{\tau}  _1 ^{\{j _1\} _{(m)}}
\cdots
\widetilde{\tau}  _r ^{\{j _r\} _{(m)}}\in \widetilde{\cP}^{n} _{Y/S,(m)}$.
We get the basis of the free $\cB _Y$-module
$\widetilde{\cP}^{n} _{Y/S,(m)}$ given by 
$\underline{\tau} ^{\{ \underline{i}\} _{(m)}}\underline{\widetilde{\tau}}  {} ^{\{ \underline{j}\} _{(m)}}$, 
with $| \underline{i}| + | \underline{j}| \leq n$.
We denote by 
$\{ \underline{\partial}  ^{<\underline{i}> _{(m)}}
\underline{\widetilde{\partial}}  {} ^{<\underline{j}> _{(m)}} \text{, with $| \underline{i}| + | \underline{j}| \leq n$}\} $ 
the corresponding dual basis of $\widetilde{\cD} ^{(m)} _{Y/S,n}$.
By taking the inductive limits (i.e. this is simply a union),
this yields the basis
$\{ \underline{\partial}  ^{<\underline{i}> _{(m)}}
\underline{\widetilde{\partial}}  {} ^{<\underline{j}> _{(m)}} \text{, with $ \underline{i} \in \bbN ^d$ and $ \underline{j}\in \bbN ^r$}\} $
of the free $\cB _Y$-module $\widetilde{\cD} ^{(m)} _{Y/S}$. 
In other words, 
a section of the sheaf
$\widetilde{\cD} ^{(m)} _{Y/S}$ can uniquely be written as a {\it finite} sum of the form
$\sum  _{\underline{i} \in \N ^r, \underline{j} \in \N ^r} 
a _{\underline{i}, \underline{j}} 
\underline{\partial}  ^{<\underline{i}> _{(m)}}
\underline{\widetilde{\partial}}  {} ^{<\underline{j}> _{(m)}}$
(resp. 
$\sum  _{\underline{i} \in \N ^r, \underline{j} \in \N ^r} 
\underline{\partial}  ^{<\underline{i}> _{(m)}}
\underline{\widetilde{\partial}}  {} ^{<\underline{j}> _{(m)}}
 a _{\underline{i}, \underline{j}}$)
with 
$a _{\underline{i}, \underline{j}} \in \cB _Y$.

We hope this is not too confusing that 
$\underline{\partial}  {} ^{<\underline{i}> _{(m)}} $
(resp. $\underline{\widetilde{\partial}}  {} ^{<\underline{j}> _{(m)}} $)
is either a global section of 
$\cD ^{(m)} _{Y/S}$
or of $\cD ^{(m)} _{Y/T}$
(resp. of $\cD ^{(m)} _{T/S}$).

\end{ntn}

\begin{empt}
\label{empt-DeltaY/T2Y-cor}
Suppose $T/S$ has the finite $p$-basis 
$\tilde{t} _1,\cdots, \tilde{t} _r$
and 
$Y/T$ has the finite $p$-basis 
$t _1,\cdots, t _d$.
We keep notations 
\ref{ntn-loc-coor-TvsS}, \ref{ntn-loc-coor-TvsSbis}
and \ref{ntn-loc-coor-TvsS2}.

\begin{enumerate}[(a)]
\item 
\label{empt-DeltaY/T2Y-cor-item} 
Following \ref{empt-DeltaY/T2Y}, 
we have the homomorphisms of rings
$\mathcal{P} ^n _{Y /S , (m)}
\to 
\mathcal{P} ^n _{Y  / T , (m)}$.
We compute this map sends 
$\underline{\tau}  {} ^{\{\underline{i}\} _{(m)}}$ to 
$\underline{\tau}  {} ^{\{\underline{i}\} _{(m)}}$, 
which justifies why we took the same notation. 
Hence, the morphism
$\widetilde{\cD} ^{(m)} _{Y  / T }
\to 
\widetilde{\cD} ^{(m)} _{Y/S}$
corresponds to the inclusion given by
$$\sum  _{\underline{i} \in \N ^d}
a _{\underline{i}} 
\underline{\partial}  ^{<\underline{i}> _{(m)}}
\mapsto
\sum  _{\underline{i} \in \N ^d}
a _{\underline{i}} 
\underline{\partial}  ^{<\underline{i}> _{(m)}}
,$$
where 
$a _{\underline{i}}$ are global sections of $\cB _Y$.
Since $\mathcal{P} ^n _{Y /S , (m)}
\to 
\mathcal{P} ^n _{Y  / T , (m)}$
is a homomorphism of $\cB _Y$-algebras for the right structure 
(and also for the left one, but this is useless here), 
then the action of 
$\widetilde{\cD} ^{(m)} _{Y  / T }$ on 
$\cB _Y$ and of
$\widetilde{\cD} ^{(m)} _{Y/S}$
on 
$\cB _Y$
are compatible with the canonical inclusion
$\widetilde{\cD} ^{(m)} _{Y  / T }
\to 
\widetilde{\cD} ^{(m)} _{Y/S}$.
Hence, we get the homomorphism
$\widetilde{\cD} ^{(m)} _{Y  / T }
\to 
\widetilde{\cD} ^{(m)} _{Y/S}$
is also given by the formula
$$
\sum  _{\underline{i} \in \N ^d}
\underline{\partial}  ^{<\underline{i}> _{(m)}}
a _{\underline{i}} 
\mapsto
\sum  _{\underline{i} \in \N ^d}
\underline{\partial}  ^{<\underline{i}> _{(m)}}a _{\underline{i}},$$
where 
$a _{\underline{i}}$ are global sections of $\cB _Y$.

\item 
\label{empt-DeltaY/T2Y-cor-item2} 
Using the universal property of $m$-PD-envelopes, 
we get the homomorphisms of rings
$g ^* \mathcal{P} ^n _{T/S , (m)}
\to 
\mathcal{P} ^n _{Y / S , (m)}$.
We compute this map sends 
$1 \otimes \underline{\widetilde{\tau}}  {} ^{\{\underline{j}\} _{(m)}}$ 
to 
$\underline{\widetilde{\tau}}  {} ^{\{\underline{j}\} _{(m)}}$, 
which justifies a bit why we took the same notation. 
This yields that 
the homomorphism 
$\widetilde{\cD} ^{(m)} _{Y/S} 
\to 
g ^* \widetilde{\cD} ^{(m)} _{T/S}$
is given by 
$$\sum  _{\underline{i} \in \N ^r, \underline{j} \in \N ^r} 
\alpha _{\underline{i}, \underline{j}} 
\underline{\partial}  ^{<\underline{i}> _{(m)}}
\underline{\widetilde{\partial}}  {} ^{<\underline{j}> _{(m)}}   
\mapsto
\sum  _{\underline{j} \in \N ^r} 
\alpha _{\underline{0}, \underline{j}} \otimes \underline{\widetilde{\partial}}  {} ^{<\underline{j}> _{(m)}}  ,
$$
where 
$\alpha _{\underline{i}, \underline{j}} \in \cB _Y$.

\item 
\label{empt-DeltaY/T2Y-cor-item3} 
The left $\widetilde{\cD} ^{(m)} _{Y  / T }$-module (resp. right $\widetilde{\cD} ^{(m)} _{Y  / T }$-module)
$\widetilde{\cD} ^{(m)} _{Y /S}$ canonically splits as follows:
\begin{equation}
\label{filt-DYm}
 \widetilde{\cD} ^{(m)} _{Y  /S}
=
\oplus _{\underline{j} \in \N ^r}
 \widetilde{\cD} ^{(m)} _{Y  / T } \underline{\widetilde{\partial}}  {} ^{<\underline{j}> _{(m)}}, 
 \
\widetilde{\cD} ^{(m)} _{Y /S }
=
 \oplus _{\underline{j} \in \N ^r}
 \underline{\widetilde{\partial}}  {} ^{<\underline{j}> _{(m)}} 
 \widetilde{\cD} ^{(m)} _{Y  / T },
\end{equation}
where 
$ \widetilde{\cD} ^{(m)} _{Y  / T } \underline{\widetilde{\partial}}  {} ^{<\underline{j}> _{(m)}}$
(resp. 
$ \underline{\widetilde{\partial}}  {} ^{<\underline{j}> _{(m)}} \widetilde{\cD} ^{(m)} _{Y  / T }$) 
is the left (resp. right) 
free  $\widetilde{\cD} ^{(m)} _{Y  / T }$-submodule of 
$\widetilde{\cD} ^{(m)} _{Y /S }$
generated by 
$\underline{\widetilde{\partial}}  {} ^{<\underline{j}> _{(m)}}$.
We get the exhausted filtration of $\widetilde{\cD} ^{(m)} _{Y/S  } $ by 
left $\widetilde{\cD} ^{(m)} _{Y  / T }$-submodules 
(resp. right $\widetilde{\cD} ^{(m)} _{Y  / T }$-submodules)
$F ^{\mathrm{l}} _n \widetilde{\cD} ^{(m)} _{Y /S } := 
 \oplus _{|\underline{j} | \leq n}
 \widetilde{\cD} ^{(m)} _{Y  / T } 
 \underline{\widetilde{\partial}}  {} ^{<\underline{j}> _{(m)}} $
(resp. $F ^{\mathrm{r}} _n \widetilde{\cD} ^{(m)} _{Y /S } := 
 \oplus _{|\underline{j} | \leq n}
 \underline{\widetilde{\partial}}  {} ^{<\underline{j}> _{(m)}} 
 \widetilde{\cD} ^{(m)} _{Y  / T }$).
\end{enumerate}

\end{empt}

\begin{empt}
\label{coro-flat-DYTS}
It follows from \ref{empt-DeltaY/T2Y-cor}.\ref{empt-DeltaY/T2Y-cor-item}
that $\widetilde{\cD} ^{(m)} _{Y  / T }
\to 
\widetilde{\cD} ^{(m)} _{Y/S}$
is a monomorphism, 
from \ref{empt-DeltaY/T2Y-cor}.\ref{empt-DeltaY/T2Y-cor-item2}
that 
$\widetilde{\cD} ^{(m)} _{Y/S} 
\to 
g ^* \widetilde{\cD} ^{(m)} _{T/S}$
is an epimorphism. 
From \ref{empt-DeltaY/T2Y-cor}.\ref{empt-DeltaY/T2Y-cor-item3}, 
we check that 
$\widetilde{\cD} ^{(m)} _{Y/S }$
is a left (resp. right) flat $\widetilde{\cD} ^{(m)} _{Y  / T } $-module. 
This yields from \ref{pre-u!=u!/T(m)} the isomorphism 
\begin{equation}
\label{pre-u!=u!/T(m)bis}
\widetilde{\cD} ^{(m)} _{X \to Y  / T }
\otimes ^{\bbL}_{h ^{-1}\widetilde{\cD} ^{(m)} _{Y  / T}} h ^{-1}\widetilde{\cD} ^{(m)} _{Y /S}
\riso 
\widetilde{\cD} ^{(m)} _{X \to Y  /S}.
\end{equation}
\end{empt}

\begin{empt}
For any 
$\E \in  D ^{-} (\overset{^\mathrm{l}}{} \widetilde{\cD} ^{(m)} _{Y  / T })$, 
we will write
$ h ^{(m)!} _{/T} (\E) := 
\widetilde{\cD} ^{(m)} _{X \to Y  / T }
\otimes ^{\L} _{h ^{-1}\widetilde{\cD} ^{(m)} _{Y  / T }} 
h ^{-1} \E [\delta _{X/Y}] $, 
and 
for any 
$\E \in  D ^{-} 
(\overset{^\mathrm{l}}{} \widetilde{\cD} ^{(m)} _{Y/S  })$,
we will write
$ h ^{(m)!} _{/S}(\E) := 
\widetilde{\cD} ^{(m)} _{X \to Y  }
\otimes ^{\L} _{h ^{-1}\widetilde{\cD} ^{(m)} _{Y/S  }} 
h ^{-1} \E [\delta _{X/Y}]$.
We denote by 
$\mathrm{oub} _{Y/T/S}
\colon 
D ^{-} (\overset{^\mathrm{l}}{} \widetilde{\cD} ^{(m)} _{Y /S })
\to 
D ^{-} (\overset{^\mathrm{l}}{} \widetilde{\cD} ^{(m)} _{Y  / T })$
the canonical forgetful functor (and similarly by replacing $Y$ with $X $).

\end{empt}

\begin{prop}
\label{prop-u!=u!/T(m)}
For any $\E \in  D ^{-} 
(\overset{^\mathrm{l}}{} \widetilde{\cD} ^{(m)} _{Y/S })$,
we get the isomorphism
\begin{equation}
\label{u!=u!/T(m)}
\mathrm{oub} _{X/T/S}  \circ h ^{(m)!} _{/S} (\E) 
\riso 
h ^{(m)!} _{/T} \circ  
\mathrm{oub} _{Y/T/S} (\E).
\end{equation}
\end{prop}

\begin{proof}
By associativity of the tensor product, 
this is a consequence of \ref{pre-u!=u!/T(m)bis}.
\end{proof}

\begin{empt}
\label{partialR2Lcomp-empt}
We have the commutative diagram of left $\widetilde{\cD} ^{(m)} _{X/S}$-modules
\begin{equation}
\label{partialR2Lcomp-diag1}
\xymatrix{
{\widetilde{\cD} ^{(m)} _{X/S}} 
\ar[r] ^-{}
\ar[d] ^-{}
& 
{h ^* \widetilde{\cD} ^{(m)} _{Y /S } } 
\ar[d] ^-{}
\ar[dl] ^-{\psi}
\\ 
{f ^* \widetilde{\cD} ^{(m)} _{T/S}} 
\ar[r] ^-{\sim}
& 
{ h ^* g^* \widetilde{\cD} ^{(m)} _{T/S} ,} 
}
\end{equation}
where $\psi$ is the composition $\psi \colon 
h ^* \widetilde{\cD} ^{(m)} _{Y /S } 
\to 
h ^* g ^* \widetilde{\cD} ^{(m)} _{T/S}
\riso 
f ^*  \widetilde{\cD} ^{(m)} _{T/S}$.
Indeed, we check the commutativity of the square as follows: 
since both composition morphisms are epimorphisms of left $\widetilde{\cD} ^{(m)} _{X/S}$-modules,
we reduce to check that the images of $1$ via both paths $\widetilde{\cD} ^{(m)} _{X/S}\to  h ^* g ^* \widetilde{\cD} ^{(m)} _{T/S} $
are the same. 
We compute $1\mapsto 1 \otimes 1 \otimes 1$ via both paths.

Suppose $T/S$ has the finite $p$-basis 
$\tilde{t} _1,\cdots, \tilde{t} _r$, 
$Y/T$ has the finite $p$-basis 
$t _1,\cdots, t _d$,
$X/T$ has the finite $p$-basis 
$t ' _1,\cdots, t '_{d'}$.
By abuse of notation, 
we denote by 
$\tilde{t} _1,\cdots, \tilde{t} _r$ the element of $\Gamma (Y, \cO _Y)$
induced by 
$\tilde{t} _1,\cdots, \tilde{t} _r$
via $g$.
We keep notation \ref{ntn-loc-coor-TvsS2}: 
we get the basis
$\{ \underline{\partial}  ^{<\underline{i}> _{(m)}}
\underline{\widetilde{\partial}}  {} ^{<\underline{j}> _{(m)}} \text{, with $ \underline{i} \in \bbN ^d$ and $ \underline{j}\in \bbN ^r$}\} $
of the free $\cB _Y$-module $\widetilde{\cD} ^{(m)} _{Y/S}$. 

We denote by 
$\tilde{t} '_1,\cdots, \tilde{t} '_r$ the elements of $\Gamma (X, \cO _X)$
induced by 
$\tilde{t} _1,\cdots, \tilde{t} _r$
via $f$ (we add some prime to avoid any confusion).
Similarly to notation \ref{ntn-loc-coor-TvsS2},
we get the basis
$\{ \underline{\partial} ' {}  ^{<\underline{i}'> _{(m)}}
\underline{\widetilde{\partial}} '  {} ^{<\underline{j}> _{(m)}} \text{, with $ \underline{i} '\in \bbN ^{d'}$ and $ \underline{j}\in \bbN ^r$}\} $
of the free $\cB _X$-module $\widetilde{\cD} ^{(m)} _{X/S}$.

Let $n \in \N$. 
Fix $\underline{l}\in \bbN ^r$ such that 
$| \underline{l}| = n$.

i) The morphism of left $\widetilde{\cD} ^{(m)} _{X}$-modules
$\widetilde{\cD} ^{(m)} _{X/S }  \to h ^*( \widetilde{\cD} ^{(m)} _{Y /S } )$
factorizes through
$\widetilde{\cD} ^{(m)} _{X/S,n}  \to h ^*( \widetilde{\cD} ^{(m)} _{Y/S  ,n} )$.
This yields
$\underline{\widetilde{\partial}} ' {} ^{<\underline{l}> _{(m)}}
\cdot( 1 \otimes1) \in h ^*( \widetilde{\cD} ^{(m)} _{Y /S,n} )$. 
Hence, we can write uniquely 
\begin{equation}
\label{partialR2Lcomp-form1}
\underline{\widetilde{\partial}} '  {} ^{<\underline{l}> _{(m)}}
\cdot( 1 \otimes1)
=
\sum _{\underline{i} \in \N ^d, \underline{j} \in \N ^r,
|\underline{i}|+|\underline{j}|\leq n} 
a _{\underline{i}, \underline{j}} 
\otimes  
\underline{\partial}  ^{<\underline{i}> _{(m)}}
\underline{\widetilde{\partial}}  {} ^{<\underline{j}> _{(m)}},
\end{equation}
where the sum is finite and where
$a _{\underline{i}, \underline{j}}  \in \cB _X$.

ii) By using  \ref{empt-DeltaY/T2Y-cor}.\ref{empt-DeltaY/T2Y-cor-item2}
we compute 
\begin{equation}
\label{partialR2Lcomp-mapsto1}
\psi 
\left ( 
 \sum  _{\underline{i} \in \N ^d, \underline{j} \in \N ^r} 
a _{\underline{i}, \underline{j}} 
\otimes
\underline{\partial}  ^{<\underline{i}> _{(m)}}
\underline{\widetilde{\partial}}  {} ^{<\underline{j}> _{(m)}} 
\right ) 
=
\sum  _{\underline{j} \in \N ^r} 
a _{\underline{0}, \underline{j}} 
 \otimes \underline{\widetilde{\partial}}  {} ^{<\underline{j}> _{(m)}}   .
\end{equation}

iii) Since $\psi$ is $\widetilde{\cD} ^{(m)} _{X/S }$-linear, 
$\psi (\underline{\widetilde{\partial}} ' {} ^{<\underline{l}> _{(m)}}
\cdot( 1 \otimes1) )
=
\underline{\widetilde{\partial}}' {} ^{<\underline{l}> _{(m)}}
\cdot
 \psi( 1 \otimes1)
 =
 \underline{\widetilde{\partial}} ' {} ^{<\underline{l}> _{(m)}} \cdot ( 1 \otimes1)$.
By using \ref{empt-DeltaY/T2Y-cor}.\ref{empt-DeltaY/T2Y-cor-item2} (applied to $X/S$ instead of $Y/S$),
we get  
$
 \underline{\widetilde{\partial}} ' {} ^{<\underline{l}> _{(m)}} \cdot ( 1 \otimes1)
 =
1 \otimes  \underline{\widetilde{\partial}}  {} ^{<\underline{l}> _{(m)}}$.
Hence
\begin{equation}
\label{partialR2Lcomp-form2}
\psi (\underline{\widetilde{\partial}} ' {} ^{<\underline{l}> _{(m)}}
\cdot( 1 \otimes1) )
= 
1 \otimes  \underline{\widetilde{\partial}}  {} ^{<\underline{l}> _{(m)}}.
\end{equation}

iv) It follows from  \ref{partialR2Lcomp-form1}, \ref{partialR2Lcomp-mapsto1} and \ref{partialR2Lcomp-form2} 
that we have the formula
$\sum _{\underline{j} \in \N ^r} 
a _{\underline{0}, \underline{j}} 
 \otimes \underline{\widetilde{\partial}}  {} ^{<\underline{j}> _{(m)}} 
 =1 \otimes  \underline{\widetilde{\partial}}  {} ^{<\underline{l}> _{(m)}}$.  
This yields $a _{\underline{0}, \underline{l}}  = 1$ 
and 
$a _{\underline{0}, \underline{j}} = 0$ if $\underline{j}\not = \underline{l}$.
Hence, 
we have in 
$h ^*( \widetilde{\cD} ^{(m)} _{Y/S  ,n} )$ the equality :
\begin{equation}
\label{partialR2Lcomp}
\underline{\widetilde{\partial}} ' {} ^{<\underline{l}> _{(m)}}
\cdot( 1 \otimes1)
=
1 \otimes  \underline{\widetilde{\partial}}  {} ^{<\underline{l}> _{(m)}}
+
\sum _{\underline{i} \in \N ^d, \underline{j} \in \N ^r,|\underline{i}|+|\underline{j}|\leq n,
|\underline{i}|\not =0} 
a _{\underline{i}, \underline{j}} 
\otimes
\underline{\partial}  ^{<\underline{i}> _{(m)}}
\underline{\widetilde{\partial}}  {} ^{<\underline{j}> _{(m)}}.
\end{equation}
Hence, we have in 
$h ^*( \widetilde{\cD} ^{(m)} _{Y/S  ,n} )$ the congruence :
\begin{equation}
\label{partialR2Lcompcor}
\underline{\widetilde{\partial}}  {} ^{<\underline{l}> _{(m)}}
\cdot( 1 \otimes1)
\equiv
1 \otimes  \underline{\widetilde{\partial}}  {} ^{<\underline{l}> _{(m)}}
\mod h ^* ( F ^{\mathrm{l}} _{n-1} \widetilde{\cD} ^{(m)} _{Y /S }),
\end{equation}
where $(F ^{\mathrm{l}} _{n} \widetilde{\cD} ^{(m)} _{Y /S } )_n$
is the filtration defined at \ref{empt-DeltaY/T2Y-cor}.\ref{empt-DeltaY/T2Y-cor-item3}.
\end{empt}

\begin{lem}
\label{lem-u+=u+/T}
The canonical morphism of 
$(\widetilde{\cD} ^{(m)} _{X/S}, h ^{-1}\widetilde{\cD} ^{(m)} _{Y  / T })$-bimodules 
\begin{equation}
\label{lem-u+=u+/Tiso}
\widetilde{\cD} ^{(m)} _{X/S} \otimes _{\widetilde{\cD} ^{(m)} _{X/ T}}
\widetilde{\cD} ^{(m)} _{X \to Y  / T }
\to 
\widetilde{\cD} ^{(m)} _{X \to Y /S }
\end{equation}
is an isomorphism.
\end{lem}

\begin{proof}
The canonical homomorphism \ref{lem-u+=u+/Tiso}
is constructed as follows. 
By applying the functor $h ^*$
to the homomorphism 
$\widetilde{\cD} ^{(m)} _{Y /T}
\to 
\widetilde{\cD} ^{(m)} _{Y /S}$, 
we get the homomorphism of 
$(\widetilde{\cD} ^{(m)} _{X/T}, h ^{-1}\widetilde{\cD} ^{(m)} _{Y  / T })$-bimodules 
$\widetilde{\cD} ^{(m)} _{X \to Y  / T }
=
h ^* \widetilde{\cD} ^{(m)} _{Y  /T}
\to 
h ^* \widetilde{\cD} ^{(m)} _{Y  /S}
=
\widetilde{\cD} ^{(m)} _{X \to Y  /S}$.
This yields the homomorphism  of $(\widetilde{\cD} ^{(m)} _{X/S}, h ^{-1}\widetilde{\cD} ^{(m)} _{Y  / T })$-bimodules 
$$\phi 
\colon 
\widetilde{\cD} ^{(m)} _{X/S} \otimes _{\widetilde{\cD} ^{(m)} _{X/ T}}
\widetilde{\cD} ^{(m)} _{X \to Y  / T }
\to 
\widetilde{\cD} ^{(m)} _{X \to Y  /S}
.$$
We have to check that this is an isomorphism.
Since this is local, we can suppose
$T/S$ has the finite $p$-basis 
$\tilde{t} _1,\cdots, \tilde{t} _r$, 
$Y/T$ has the finite $p$-basis 
$t _1,\cdots, t _d$,
$X/T$ has the finite $p$-basis 
$t ' _1,\cdots, t '_{d}$. 
We follow notation \ref{partialR2Lcomp-empt}.

Let 
$P \in \widetilde{\cD} ^{(m)} _{X/S}\otimes _{\widetilde{\cD} ^{(m)} _{X/ T}}\widetilde{\cD} ^{(m)} _{X \to Y  / T }$.
By using \ref{ntn-loc-coor-TvsS}
and by using \ref{empt-DeltaY/T2Y-cor}.\ref{empt-DeltaY/T2Y-cor-item3} for $X/S$ instead of $Y/S$, 
we can  uniquely write (use \ref{filt-DYm}) $P$ of the form 
$$P=\sum  _{\underline{i} \in \N ^r, \underline{j} \in \N ^r} 
\underline{\widetilde{\partial}} ' {} ^{<\underline{j}> _{(m)}}  \otimes ( a _{\underline{i}, \underline{j}} \otimes 
\underline{\partial}  ^{<\underline{i}> _{(m)}})$$
where the sum is finite and 
$a _{\underline{i}, \underline{j}} \in \cB _X$.
We get 
$$\phi (P) = 
\sum  _{\underline{i} \in \N ^r, \underline{j} \in \N ^r} 
\underline{\widetilde{\partial}} ' {} ^{<\underline{j}> _{(m)}}
a _{\underline{i}, \underline{j}}  \cdot( 1 \otimes   
\underline{\partial}  ^{<\underline{i}> _{(m)}})
=
\sum  _{\underline{i} \in \N ^r, \underline{j} \in \N ^r} 
\underline{\widetilde{\partial}} ' {} ^{<\underline{j}> _{(m)}} 
a _{\underline{i}, \underline{j}}  \cdot( 1 \otimes1)  \cdot \underline{\partial}  ^{<\underline{i}> _{(m)}}.$$
Let $n:= \max \{k\in \bbN ~| ~\exists \underline{j} \in \N ^r, 
\exists \underline{i}\in \N ^d,
\text{such that 
$| \underline{j}|=k $
and 
$a _{\underline{i}, \underline{j}} 
\not =0 $}\}$.
Let $\underline{l}\in \bbN ^r$ be such that 
$| \underline{l}| = n$.
For any integer $s$, we denote by 
$\cD ^{(m)} _{X,T/S,s}$ the free $\cB _X$-submodule (for both structure) of 
$\cD ^{(m)} _{X/S}$ whose basis is given by 
$\underline{\widetilde{\partial}} ' {} ^{<\underline{j}> _{(m)}}$
for any 
$\underline{j} \in \bbN ^r$ 
such that 
$|\underline{j}| \leq s$.
We remark that 
$\underline{\widetilde{\partial}} ' {} ^{<\underline{l}> _{(m)}}
a _{\underline{i},\underline{l}} 
- 
a _{\underline{i},\underline{l}} 
\underline{\widetilde{\partial}} ' {} ^{<\underline{l}> _{(m)}}
\in 
\cD ^{(m)} _{X,T/S,n-1}$.
Hence, 
by using \ref{partialR2Lcompcor}, 
we compute
$$\underline{\widetilde{\partial}} ' {} ^{<\underline{l}> _{(m)}}a _{\underline{i},\underline{l}}  \cdot( 1 \otimes1)
\equiv
a _{\underline{i},\underline{l}} \underline{\widetilde{\partial}} ' {} ^{<\underline{l}> _{(m)}}  \cdot( 1 \otimes1)
\equiv
a _{\underline{i},\underline{l}}  \otimes \underline{\widetilde{\partial}}  {} ^{<\underline{l}> _{(m)}} 
\mod h ^ * ( F ^{\mathrm{l}} _{n-1} \widetilde{\cD} ^{(m)} _{Y/S  }).$$
Since the action of $\underline{\partial}  ^{<\underline{i}> _{(m)}}$
via the right $h ^{-1}\widetilde{\cD} ^{(m)} _{Y/S }$-module structure
of $h ^{*}\widetilde{\cD} ^{(m)} _{Y /S }$ preserves 
$h ^* ( F ^{\mathrm{l}} _{n-1} \widetilde{\cD} ^{(m)} _{Y  /S})$ (because 
$\underline{\partial}  ^{<\underline{i}> _{(m)}}$ and 
$\underline{\widetilde{\partial}}  {} ^{<\underline{j}> _{(m)}}$ commute),
we get
$$\underline{\widetilde{\partial}} ' {} ^{<\underline{l}> _{(m)}}a _{\underline{i},\underline{l}}  
\cdot( 1 \otimes1)  \cdot \underline{\partial}  ^{<\underline{i}> _{(m)}}
\equiv
a _{\underline{i},\underline{l}}  
\otimes 
\underline{\widetilde{\partial}}  {} ^{<\underline{l}> _{(m)}}
\underline{\partial}  ^{<\underline{i}> _{(m)}} 
\mod h ^* ( F ^{\mathrm{l}} _{n-1} \widetilde{\cD} ^{(m)} _{Y /S }).$$
Since 
$h ^{*}\widetilde{\cD} ^{(m)} _{Y /S }$ is a free $\cB _X$-module
with the basis 
$\{
\underline{\partial}  ^{<\underline{i}> _{(m)}} 
\underline{\widetilde{\partial}}  {} ^{<\underline{j}> _{(m)}}
~
|
~
\underline{i} \in \bbN ^d,
~\underline{j} \in \bbN ^r
 \}$
then 
from this latter congruence, we check easily by induction in $n$ the injectivity and the surjectivity of $\phi$.
\end{proof}

\begin{prop}
For any $\E \in  D ^{-} 
(\overset{^\mathrm{l}}{} \widetilde{\cD} ^{(m)} _{Y  / T})$,
we get the isomorphism 
of 
$ D ^{-} 
(\overset{^\mathrm{l}}{} \widetilde{\cD} ^{(m)} _{X/ S})$
\begin{equation}
\label{u!=u!/T(m)2}
\widetilde{\cD} ^{(m)} _{X/ S} \otimes _{\widetilde{\cD} ^{(m)} _{X/ T}} h ^{(m)!} _{/T} (\E)
\riso 
h ^{(m)!} _{/S} (\widetilde{\cD} ^{(m)} _{Y / S } \otimes _{\widetilde{\cD} ^{(m)} _{Y  / T}} \E) .
\end{equation}

\end{prop}

\begin{proof}
By associativity of the tensor product, we get 
\begin{gather}
\notag
\widetilde{\cD} ^{(m)} _{X/ S} \otimes _{\widetilde{\cD} ^{(m)} _{X/ T}} h ^{(m)!} _{/T} (\E)
=
\widetilde{\cD} ^{(m)} _{X/ S} \otimes _{\widetilde{\cD} ^{(m)} _{X/ T}} 
\left ( 
\widetilde{\cD} ^{(m)} _{X \to Y  / T }
\otimes ^{\L} _{h ^{-1}\widetilde{\cD} ^{(m)} _{Y  / T }} 
h ^{-1} \E
\right ) [\delta _{X/Y}]
\\
\notag
\underset{\ref{lem-u+=u+/Tiso}}{\riso}
\widetilde{\cD} ^{(m)} _{X \to Y  / S}
\otimes ^{\L} _{h ^{-1}\widetilde{\cD} ^{(m)} _{Y  / T }} 
h ^{-1} \E [\delta _{X/Y}]
\\
\notag
\riso
\widetilde{\cD} ^{(m)} _{X \to Y  / S}
\otimes ^{\L} _{h ^{-1}\widetilde{\cD} ^{(m)} _{Y  / S }} 
h ^{-1} \left( \widetilde{\cD} ^{(m)} _{Y  / S }
\otimes ^{\L} _{\widetilde{\cD} ^{(m)} _{Y  / T }} 
\E \right) [\delta _{X/Y}]
=
h ^{(m)!} _{/S}(\widetilde{\cD} ^{(m)} _{Y  / S} \otimes _{\widetilde{\cD} ^{(m)} _{Y  / T}} \E) .
\end{gather}
\end{proof}

\begin{prop}
\label{prop-u+(m)=u+(m)/T}
Let $\M \in  D ^{-} _{\mathrm{qc}}
(\overset{^\mathrm{r}}{} \widetilde{\cD} ^{(m)} _{X  / S})$. 
Then the canonical morphism of 
$D ^{-} _{\mathrm{qc}}
(\overset{^\mathrm{r}}{} h ^{-1}\widetilde{\cD} ^{(m)} _{Y / T })$
\begin{equation}
\M \otimes ^\L _{\widetilde{\cD} ^{(m)} _{X/ T}}\widetilde{\cD} ^{(m)} _{X \to Y  / T }
\to 
\M \otimes ^\L _{\widetilde{\cD} ^{(m)} _{X/S}}\widetilde{\cD} ^{(m)} _{X \to Y/S }
\end{equation}
is an isomorphism.
\end{prop}

\begin{proof}
Since this is local, we can suppose $X$ affine.
Using the way-out left version of  \cite[I.7.1.(iv)]{HaRD},
since the functors
$\M \mapsto \M \otimes ^\L _{\widetilde{\cD} ^{(m)} _{X/ T}}\widetilde{\cD} ^{(m)} _{X \to Y  / T }$
and 
$\M \mapsto \M \otimes ^\L _{\widetilde{\cD} ^{(m)} _{X/S}}\widetilde{\cD} ^{(m)} _{X \to Y  /S}$ are way-out left, 
we reduce to check the isomorphism
when $\M$ is a free right $\widetilde{\cD} ^{(m)} _{X/S  }$-module.
Hence, we come down to the case where
$\M = \widetilde{\cD} ^{(m)} _{X/S  }$. 
In that case, $\M$ is a flat right $\widetilde{\cD} ^{(m)} _{X /S }$-module and a 
flat right $\widetilde{\cD} ^{(m)} _{X  / T }$-module (see \ref{empt-DeltaY/T2Y-cor}.\ref{empt-DeltaY/T2Y-cor-item3}).
Hence, we conclude using \ref{lem-u+=u+/T}.\end{proof}

\begin{empt}
We define the functor 
$h ^{(m)}  _{/S ~+}
\colon 
D ^{-} _{\mathrm{qc}}
(\overset{^\mathrm{r}}{} \widetilde{\cD} ^{(m)} _{X /S})
\to 
D ^{-} _{\mathrm{qc}}
(\overset{^\mathrm{r}}{} \widetilde{\cD} ^{(m)} _{Y/S})$
by setting 
$$h ^{(m)}  _{/S ~+} (\M) 
:=
\R h _* \left ( \M \otimes ^\L _{\widetilde{\cD} ^{(m)} _{X/S}}\widetilde{\cD} ^{(m)} _{X \to Y  /S}
\right) 
$$
for 
 $\M \in  
D ^{-} _{\mathrm{qc}}
(\overset{^\mathrm{r}}{} \widetilde{\cD} ^{(m)} _{X/S})$.
We define the functor 
$h ^{(m)}  _{/T ~+}
\colon 
D ^{-} _{\mathrm{qc}}
(\overset{^\mathrm{r}}{} \widetilde{\cD} ^{(m)} _{X /T })
\to 
D ^{-} _{\mathrm{qc}}
(\overset{^\mathrm{r}}{} \widetilde{\cD} ^{(m)} _{Y  /T })$
by setting 
$h ^{(m)}  _{/T ~+} (\M) 
:=
\R h _* \left ( 
\M \otimes ^\L _{\widetilde{\cD} ^{(m)} _{X/ T}}\widetilde{\cD} ^{(m)} _{X \to Y  / T }
\right) $
for  $\M \in  D ^{-} _{\mathrm{qc}}
(\overset{^\mathrm{r}}{} \widetilde{\cD} ^{(m)} _{X /T })$.
Following \ref{prop-u+(m)=u+(m)/T}, we have 
for any $\M \in  
D ^{-} _{\mathrm{qc}}
(\overset{^\mathrm{r}}{} \widetilde{\cD} ^{(m)} _{X/S})$ 
the isomorphism
\begin{equation}
\label{isooub-u+(m)=u+(m)/T}
\mathrm{oub} _{Y/T/S}  \circ h ^{(m)}  _{/S ~+} (\M)
\riso 
h ^{(m)}  _{/T \,+} \circ 
\mathrm{oub} _{X/T/S}   (\M) .
\end{equation}
\end{empt}

\begin{prop}
For $\M \in  D ^{-} _{\mathrm{qc}}
(\overset{^\mathrm{r}}{} \widetilde{\cD} ^{(m)} _{Y  /T })$, we have the canonical isomorphism
\begin{equation}
\label{extT-u+(m)=u+(m)/T}
h ^{(m)}  _{/S \,+} (\M \otimes _{\widetilde{\cD} ^{(m)} _{X / T}} \widetilde{\cD} ^{(m)} _{X /S} )
\riso 
h ^{(m)}  _{/T \,+} (\M )  \otimes _{\widetilde{\cD} ^{(m)} _{Y  / T}} \widetilde{\cD} ^{(m)} _{Y /S } .
\end{equation}
\end{prop}

\begin{proof}
Using the projection isomorphism, we get
$$h ^{(m)}  _{/T \,+} (\M)   \otimes _{\widetilde{\cD} ^{(m)} _{Y  / T}} \widetilde{\cD} ^{(m)} _{Y /S } 
\riso 
\R h _* \left (  \left ( \M \otimes ^\L _{\widetilde{\cD} ^{(m)} _{X / T}}
\widetilde{\cD} ^{(m)} _{X \to Y  / T }
\right )
\otimes _{h ^{-1}\widetilde{\cD} ^{(m)} _{Y  / T}} h ^{-1}\widetilde{\cD} ^{(m)} _{Y  /S}
\right).$$
We conclude by using \ref{pre-u!=u!/T(m)}.
\end{proof}

\subsection{Spencer resolutions, level $0$ case}
We keep notation \ref{subsec-/Svs/T}.
We suppose $Y= T$ and $g = id$ (and then $f =h$).

\begin{empt}
\label{ntn-Spencer-pre}

 Following 
 \cite[2.2.10]{Caro-Vauclair},
 we have
 $\Omega ^1 _{X  /T  }
 =
 \ker ( \mathcal{P} ^1 _{X / T , (m)}
 \to 
\mathcal{P} ^0 _{X / T , (m)}
=
\cO _X)$. 
We set
$\widetilde{\Omega} ^1 _{X  /T  }
:= 
\cB  _X \otimes _{\cO _X}
\Omega ^1 _{X  /T  }$, and
$\widetilde{\mathcal{T}} _{X  /T  }
:= 
\mathcal{H} om _{\cB  _X} (\widetilde{\Omega} ^1 _{X  /T  }, 
\cB  _X)$. When $\cB  _X = \cO _X$, we remove the tilde symbole. 
From 
$\Omega ^1 _{X  /T  }
 \hookrightarrow 
\mathcal{P} ^1 _{X / T , (m)}$, we obtain by duality 
the canonical epimorphism 
$ \widetilde{\D} ^{(m)} _{X / T ,1}
 \twoheadrightarrow
 \widetilde{\mathcal{T}} _{X  /T  }$ whose kernel is
$\widetilde{\D} ^{(m)} _{X / T ,0}$.
Hence, 
$\gr _1 \widetilde{\D} ^{(m)} _{X / T }
\riso 
\widetilde{\mathcal{T}} _{X  /T  }$.
\end{empt}

\begin{empt}
\label{ntn-Spencer}
For any sections $v _1, v _2$ of 
$\widetilde{\mathcal{T}} _{X  /T  }$, we write
$[v _1,v _2]$ the section of  $\widetilde{\mathcal{T}} _{X  /T  }$
which corresponds to the section $v _1 v _2 -v _2 v _1$ of 
$ \widetilde{\D} ^{(m)} _{X / T ,1}$
modulo $\cO _X$.

Let $\E=(\E _n) _{n\in \N}$ 
be a filtered left 
$\widetilde{\D} ^{(0)} _{X / T }$-module, i.e a filtration so that 
$\widetilde{\D} ^{(0)} _{X / T ,n'} \cdot \E _n \subset \E _{n +n'}$.
Similarly to \cite[1.6]{kashiwarathesis}, 
we define the morphism of left 
$\widetilde{\D} ^{(0)} _{X / T }$-modules
\begin{equation}
\delta
\colon 
\widetilde{\D} ^{(0)} _{X / T } 
\otimes _{\cB  _X}
\wedge ^{i} \widetilde{\mathcal{T}} _{X  /T  }
\otimes _{\cB  _X}
\E _{j-1}
\to 
\widetilde{\D} ^{(0)} _{X / T } 
\otimes _{\cB  _X}
\wedge ^{i-1} \widetilde{\mathcal{T}} _{X  /T  }
\otimes _{\cB  _X}
\E _{j}
\end{equation}
by 
\begin{align}
\notag
\delta ( P \otimes ( v _1 \wedge \cdots \wedge v _i) \otimes u
= &
\sum _{a=1} ^i 
(-1) ^{a-1} P v _a 
 \otimes ( v _1 \wedge \cdots \wedge \widehat{v _a} \wedge \cdots \wedge v _i) \otimes u
 \\
\notag
& - \sum _{a=1} ^i
(-1) ^{a-1} P 
 \otimes ( v _1 \wedge \cdots \wedge \widehat{v _a} \wedge \cdots \wedge v _i) \otimes v _a  u
\\
\notag
&
+ \sum _{1 \leq a <b\leq i} 
(-1) ^{a-1} P 
 \otimes ( [v _a,v _b] \wedge v _1 \wedge \cdots \wedge \widehat{v _a}  \wedge \cdots \wedge \widehat{v _b} \wedge \cdots \wedge v _i) \otimes u.
\end{align}
We compute easily that that we get the following complex of left $\widetilde{\D} ^{(0)} _{X / T }$-modules
\begin{equation}
\label{Spencer-sequence}
0 
\to 
\widetilde{\D} ^{(0)} _{X / T } 
\otimes _{\cB  _X}
\wedge ^{d} \widetilde{\mathcal{T}} _{X  /T  }
\otimes _{\cB  _X}
\E _{n-d}
\cdots
\underset{\delta}{\longrightarrow}
\widetilde{\D} ^{(0)} _{X / T } 
\otimes _{\cB  _X}
\wedge \widetilde{\mathcal{T}} _{X  /T  }
\otimes _{\cB  _X}
\E _{n-1}
\underset{\delta}{\longrightarrow}
\widetilde{\D} ^{(0)} _{X / T } 
\otimes _{\cB  _X}
\E _{n}
\to 
\E
\to 0,
\end{equation}
where $d$ is the cardinal of thebasis of level $m$.
We call \ref{Spencer-sequence}
the first Spencer sequence of degree $n$ of $\E$
and denote it by 
$ \mathrm{Sp} _n (\E)$.

\end{empt}

\begin{thm}
Let $\E=(\E _n) _{n\in \N}$ be a good filtered left $\widetilde{\D} ^{(0)} _{X / T }$-module, i.e. 
$\gr \E $ is a $\gr \widetilde{\D} ^{(0)} _{X / T }$-module coherent.
With the notation \ref{ntn-Spencer}, 
$ \mathrm{Sp} _n (\E)$ is exact 
for $n$ large enough. 

\end{thm}

\begin{proof}
Similarly to \cite[5.1.1]{these_montagnon}, we check 
$\gr \D ^{(0)} _{X / T }
\riso \mathbb{S} (\mathcal{T} _{X  /T  }) $. 
Hence, 
$\gr \widetilde{\D} ^{(0)} _{X / T }
\riso \mathbb{S} (\widetilde{\mathcal{T}} _{X  /T  }) $. 
Hence, 
we can copy word by word the proof of \cite[1.6.1]{kashiwarathesis}.
\end{proof}

\begin{empt}
In particular, taking the trivial filtration of 
$\cB  _X$, we get the exact sequence of left $\widetilde{\D} ^{(0)} _{X / T }$-modules
\begin{equation}
\label{Spencer-sequence-B}
0 
\to 
\widetilde{\D} ^{(0)} _{X / T } 
\otimes _{\cB  _X}
\wedge ^{d} \widetilde{\mathcal{T}} _{X  /T  }
\cdots
\underset{\delta}{\longrightarrow}
\widetilde{\D} ^{(0)} _{X / T } 
\otimes _{\cB  _X}
\wedge \widetilde{\mathcal{T}} _{X  /T  }
\underset{\delta}{\longrightarrow}
\widetilde{\D} ^{(0)} _{X / T } 
\to 
\cB  _X
\to 0.
\end{equation}

\end{empt}

\begin{empt}
Following \ref{lem-u+=u+/Tiso},
since 
$\widetilde{\cD} ^{(m)} _{X \to T  / T } 
=
\cB _X$,
we have the canonical isomorphism of 
$(\smash{\widetilde{\D}} _{X /S } ^{(m)} ,f ^{-1} \cB _T )$-bimodules 
\begin{equation}
\label{DXY2DXT}
\smash{\widetilde{\D}} _{X /S } ^{(m)} 
\otimes _{\smash{\widetilde{\D}} _{X /T } ^{(m)} }
\cB _X
\riso
\smash{\widetilde{\D}} _{X \to T /S } ^{(m)} .
\end{equation}
By applying the exact functor 
$\smash{\widetilde{\D}} _{X /S } ^{(0)} 
\otimes _{\smash{\widetilde{\D}} _{X /T } ^{(0)} } -$
to the exact sequence \ref{Spencer-sequence-B},
by using the isomorphism \ref{DXY2DXT}, 
we get the exact sequence of left 
$\widetilde{\D} ^{(0)} _{X  /S }  $-modules:
\begin{equation}
\label{Spencer-sequence-B-prehat2}
0 
\to 
\widetilde{\D} ^{(0)} _{X  /S }  
\otimes _{\widetilde{\B} ^{(0)}_X }
\wedge ^{d} \widetilde{\mathcal{T}} _{X  /T  }
\cdots
\underset{\delta}{\longrightarrow}
\widetilde{\D} ^{(0)} _{X  /S }  
\otimes _{\widetilde{\B} ^{(0)}_X }
\widetilde{\mathcal{T}} _{X  /T  }
\underset{\delta}{\longrightarrow}
\widetilde{\D} ^{(0)} _{X  /S }  
\to 
\smash{\widetilde{\D}} _{X \to T /S } ^{(0)} 
\to 0.
\end{equation}

\end{empt}

\subsection{Projection formula}
Let $T$ be an ${S _i}$-scheme of finite type 
for some integer $i \geq 0$.
Let $u \colon Y  \to X $ be a morphism of 
$T$-schemes locally of formal finite type
and having locally finite $p$-bases over $T$.
Recall that following the remark \ref{rem-Dnn'Ynoeth}.\ref{rem-Dnn'Ynoeth-i)} and Theorem
\ref{f0formétale-fforméta0},
$X$ and $Y$ are noetherian,
$X/T$ and $Y/T$ are flat.
Hence, $u$ is quasi-separated and quasi-compact.
 Let $\B _X$ be an $\O _X$-algebra endowed with a compatible structure of 
left $\D ^{(m)} _{X/ T }$-module. 
Put 
$\widetilde{\D} ^{(m)} _{X/ T}
:= 
\B _X \otimes _{\O _X} \D ^{(m)} _{X/ T}$,
$\B _Y := u ^* (\B _X)$,
 $\widetilde{\D} ^{(m)} _{Y/ T}
:= 
\B _Y \otimes _{\O _Y} \D ^{(m)} _{Y/ T}$,
$\D _{Y\to X/T } ^{(m)}:=
u ^* \D _{X/T } ^{(m)}$,
$\smash{\widetilde{\D}} _{Y\to X/T} ^{(m)}:=
\B ^{(m)} _{Y}   \otimes _{\O _{Y}} 
\D _{Y\to X/T } ^{(m)}$.

\begin{empt}
Following \cite[3.6.5]{Tohoku}, since
$X$ is noetherian of finite Krull dimension $d _X$, then for $i> d _X$, 
for every sheaf $\E$ of abelian groups
we have 
$H ^i (X, \E)=0$. 
Then, following \cite[12.2.1]{EGAIII1},
we get 
that 
$R ^i u _* (\E) = 0$ for $i> d _X$
and every sheaf $\E$ of abelian groups. 
In particular, by definition (see \cite[12.1.1]{EGAIII1}),
the functor $u _*$ has finite (bounded by $d _X$) cohomological dimension on $\mathrm{Mod} (u ^{-1}\O _X)$, the category 
of $u ^{-1}\O _X$-modules, or on $\mathrm{Mod} (u ^{-1}\widetilde{\D} ^{(m)} _{X / T })$.

Let $P$ be the subset of objects of $\mathrm{Mod} (u ^{-1}\widetilde{\D} ^{(m)} _{X / T })$ which are $u _*$-acyclic. 
Remark that $P$ contains injective $u ^{-1}\widetilde{\D} ^{(m)} _{X / T }$-modules. 
Using the cohomological dimension finiteness of $u _*$, if
\begin{equation}
\G ^0 \to \G ^1 \to \cdots 
\to \G ^{d _X}
\to 
\E
\to 0
\end{equation}
is an exact sequence of $\mathrm{Mod} (u ^{-1}\O _X)$,
and $\G ^0,\dots, \G ^{d _X} \in P$, 
then $\E \in P$.
Using \cite[Lemma I.4.6.2]{HaRD},
this implies that for any complex
$\E \in K (u ^{-1}\widetilde{\D} ^{(m)} _{X / T })$ 
(resp.  $\E \in K ^{-} (u ^{-1}\widetilde{\D} ^{(m)} _{X / T })$,
resp. $\E \in K ^{+} (u ^{-1}\widetilde{\D} ^{(m)} _{X / T })$,
resp. $\E \in K ^{\mathrm{b}} (u ^{-1}\widetilde{\D} ^{(m)} _{X / T })$) 
there exists a quasi-isomorphism
$\E \riso \I$ where $\I \in  K (u ^{-1}\widetilde{\D} ^{(m)} _{X / T })$ 
(resp.  $\I \in K ^{-} (u ^{-1}\widetilde{\D} ^{(m)} _{X / T })$,
resp. $\I \in K ^{+} (u ^{-1}\widetilde{\D} ^{(m)} _{X / T })$,
resp. $\I \in K ^{\mathrm{b}} (u ^{-1}\widetilde{\D} ^{(m)} _{X / T })$) 
is a complex whose modules belong to $P$.
We get the functor
$\R u _* \colon D (u ^{-1}\widetilde{\D} ^{(m)} _{X / T }) \to D (\widetilde{\D} ^{(m)} _{X / T }) $
(resp. 
$\R u _* \colon D ^- ( u ^{-1}\widetilde{\D} ^{(m)} _{X / T }) \to D ^-( \widetilde{\D} ^{(m)} _{X / T }) $,
resp. 
$\R u _* \colon D ^+ (u ^{-1}\widetilde{\D} ^{(m)} _{X / T }) \to D ^+(\widetilde{\D} ^{(m)} _{X / T }) $,
resp. 
$\R u _* \colon D ^\mathrm{b} (u ^{-1}\widetilde{\D} ^{(m)} _{X / T }) \to D ^\mathrm{b}(\widetilde{\D} ^{(m)} _{X / T }) $)
which is computed by taking a resolution with objects in $P$.

Moreover, following \cite[II.2.1]{HaRD}
$\R u_*$ takes
$D ^{?} _{\mathrm{qc}} ( \O _X) $
into 
$D ^{?} _{\mathrm{qc}} ( \O _Y) $
with $? \in \{\emptyset, +,- ,\mathrm{b}\}$.
\end{empt}

\begin{prop}
\label{Prop1.2.21-caro-surc}
Suppose one of the following conditions: 
\begin{enumerate}[(a)]
\item 
Let $\FF 
\in 
D _{\mathrm{qc},\mathrm{tdf}}
( \overset{^\mathrm{r}}{} \widetilde{\D} ^{(m)} _{X / T })$, 
and 
$\G 
\in 
D 
( \overset{^\mathrm{l}}{} u ^{-1} \widetilde{\D} ^{(m)} _{X / T })$.
\item
Let $\FF 
\in 
D ^- _{\mathrm{qc}}
( \overset{^\mathrm{r}}{} \widetilde{\D} ^{(m)} _{X / T })$, 
and 
$\G 
\in 
D ^-
( \overset{^\mathrm{l}}{} u ^{-1} \widetilde{\D} ^{(m)} _{X / T })$.
\end{enumerate}
Then we have the following isomorphism
\begin{gather}
\label{Prop1.2.21-caro-surc-iso2}
\FF \otimes ^{\L} _{\widetilde{\D} ^{(m)} _{X / T }} \R u _* ( \G)
\riso 
\R u _* \left (u ^{-1}  \FF \otimes ^{\L} _{u ^{-1} \widetilde{\D} ^{(m)} _{X / T }} \G \right).
\end{gather}
Inverting $\mathrm{r}$ and $\mathrm{l}$ in the hypotheses, 
we get the isomorphism 
\begin{gather}
\label{Prop1.2.21-caro-surc-iso2bis}
\R u _* ( \G) \otimes ^{\L} _{\widetilde{\D} ^{(m)} _{X / T }} \FF
\riso 
\R u _* \left (\G  \otimes ^{\L} _{u ^{-1} \widetilde{\D} ^{(m)} _{X / T }} u ^{-1}  \FF\right).
\end{gather}

\end{prop}

\begin{proof}
Taking a left resolution of $\FF$ by flat $\widetilde{\D} ^{(m)} _{X / T }$-modules,
and a right resolution of $\G$ by $u ^{-1} \widetilde{\D} ^{(m)} _{X / T }$-modules which are
$u _*$-acyclic, we construct the 
morphism \ref{Prop1.2.21-caro-surc-iso2}.
To check that this is an isomorphism,
using 
\cite[I.7.1 (ii), (iii) and (iv)]{HaRD}
and 
\cite[VI.5.1]{sga4-2},
we reduce to the case where 
$\FF
=
\widetilde{\D} ^{(m)} _{X / T }$, which is obvious.
\end{proof}

\begin{cor}
\label{cor-Prop1.2.21-caro-surc}
Let $*,** \in \{ \mathrm{l},\mathrm{r}\}$ such that both are not equal to $\mathrm{r}$.
Suppose one of the following conditions: 
\begin{enumerate}[(a)]
\item 
Let $\FF 
\in 
D _{\mathrm{qc},\mathrm{tdf}}
( \overset{^\mathrm{*}}{} \widetilde{\D} ^{(m)} _{X / T })$, 
and 
$\G 
\in 
D 
( \overset{^\mathrm{**}}{} u ^{-1} \widetilde{\D} ^{(m)} _{X / T })$.
\item
Let $\FF 
\in 
D ^- _{\mathrm{qc}}
( \overset{^\mathrm{*}}{} \widetilde{\D} ^{(m)} _{X / T })$, 
and 
$\G 
\in 
D ^-
( \overset{^\mathrm{**}}{} u ^{-1} \widetilde{\D} ^{(m)} _{X / T })$.
\end{enumerate}
Then we have the following isomorphism
\begin{gather}
\label{Prop1.2.21-caro-surc-iso1}
\FF \otimes ^{\L} _{\B _X} \R u _* ( \G)
\riso 
\R u _* \left (u ^{-1}  \FF \otimes ^{\L} _{u ^{-1} \B _X} \G \right).
\end{gather}

\end{cor}

\begin{proof}
For instance, if $** = \mathrm{l}$, we get 
\begin{gather}
\notag
\FF \otimes ^{\L} _{\B _X} \R u _* ( \G)
\riso 
(\FF \otimes ^{\L} _{\B _X} 
\widetilde{\D} ^{(m)} _{X / T })
 \otimes ^{\L} _{\widetilde{\D} ^{(m)} _{X / T }} 
\R u _* ( \G)
\\
\notag
\underset{\ref{Prop1.2.21-caro-surc-iso2}}{\riso}
\R u _* \left (u ^{-1} (\FF \otimes ^{\L} _{\B _X} 
\widetilde{\D} ^{(m)} _{X / T }) \otimes ^{\L} _{u ^{-1} \widetilde{\D} ^{(m)} _{X / T }} \G \right)
\riso 
\R u _* \left (u ^{-1}  \FF \otimes ^{\L} _{u ^{-1} \B _X} \G \right).
\end{gather}
\end{proof}

\begin{ntn}
\label{ntn-utilde*+}
For $\E
\in D ^-
( \overset{^\mathrm{l}}{} \widetilde{\D} ^{(m)} _{X / T })$,
we set 
$\L \widetilde{u}  ^* (\E) 
:= 
\widetilde{\D} ^{(m)} _{Y\to X / T } 
\otimes _{u ^{-1}\widetilde{\D} ^{(m)} _{X / T } }
^{\bbL}
(\E) $.
For $\M \in D ^- 
( \overset{^\mathrm{r}}{} \widetilde{\D} ^{(m)} _{Y / T })$,
we set 
$\widetilde{u} _+ ^{(m)} (\cM )
:= 
 \R u _* \left ( \cM 
\otimes ^{\L} _{\widetilde{\D} ^{(m)} _{Y / T }}
\widetilde{\D} ^{(m)} _{Y  \to X / T }
\right ).$

\end{ntn}

\begin{lemm}
\label{u^*otimes}
For $\E$ and $\FF$ two objects of 
$D ^-
( \overset{^\mathrm{l}}{} \widetilde{\D} ^{(m)} _{X / T })$,
with notation \ref{ntn-utilde*+}, 
we have the isomorphism of $D ^-
( \overset{^\mathrm{l}}{} \widetilde{\D} ^{(m)} _{Y / T })$
\begin{equation}
\label{u^*otimes-iso}
\L \widetilde{u}  ^* (\E) 
\otimes ^\L _{\B _Y}
\L  \widetilde{u}  ^* (\FF) 
\riso 
\L \widetilde{u}  ^* (\E 
\otimes ^\L _{\B _X}
\FF) .
\end{equation}
\end{lemm}

\begin{proof}
Left to the reader.\end{proof}

\begin{prop}
\label{prop-u+otimes-u!}
For $\M \in D ^- 
( \overset{^\mathrm{r}}{} \widetilde{\D} ^{(m)} _{Y / T })$ and
$\E \in D ^- _{\mathrm{qc}}
( \overset{^\mathrm{l}}{} \widetilde{\D} ^{(m)} _{X / T })$,
with notation \ref{ntn-utilde*+}
we have the canonical isomorphism
\begin{equation}
\widetilde{u} _+ ^{(m)}
\left (\M 
\otimes ^{\L} _{\B _Y}
\L \widetilde{u}  ^* (\E) 
\right )
\riso
\widetilde{u} _+ ^{(m)}  (\M )
\otimes ^{\L} _{\B _X}
\E.
\end{equation}

\end{prop}

\begin{proof}
This is proved similary to \cite[3.3.6]{caro-6operations}.
\end{proof}

\section{Frobenius descent}
Let 
$m,s\geq 0$
be two integers, $T$ be an ${S_i}$-scheme of finite type.
We suppose $\O _T$ is endowed with a quasi-coherent $m$-PD-ideal
$(\mathfrak{a},\mathfrak{b},\alpha)$ such that $p \in \mathfrak{a}$.
Since $T$ is noetherian and $\mathfrak{a}$ is an nilideal, then 
$\mathfrak{a}$ is nilpotent.
We set 
$\mathfrak{b} _1:= \mathfrak{b}+p \O _T$,
$T _0 := V ( \mathfrak{a})$.

\subsection{Complements on the $m$-PD-enveloppe of $\Delta _{X/T,(m)} (\nu)$}
Let $\nu \geq 1$ be an integer. 

\begin{ntn}
If $X$ is a $T$-scheme,
we will denote by 
$X ^{\nu +1}= X ^{\nu +1}/T$, 
$\I _\nu$
the ideal of the diagonal
$X \hookrightarrow X ^{\nu +1}$
and
$(\cP _{X, (m)} (\nu), \overline{\I} _{\nu}, \widetilde{\I} _\nu)$
the $m$-PD-enveloppe
of 
$\I _\nu$
When $\nu =1$, we 
simply write 
$(\cP _{X, (m)} , \overline{\I} , \widetilde{\I})$.

\end{ntn}

\begin{prop}
[Local description of the $m$-PD-enveloppe]
\label{prop-localdescr(nu)pre}
Let $g \colon X  \to \A ^d _{T}$ be a $p$-étale morphism.
Let  $t _1,\dots, t _d$ be the element of $\Gamma (X, \O _X)$ defining $g$.
Set 
$\tau _i := 1 \otimes t _i - t _i \otimes 1\in \I _1$
For any $j= 0,\dots, \nu$, let $p _j \colon X ^{ \nu+1}\to X$ be the index $j$ projection. 
For any $1\leq i \leq d ,\, 1 \leq j \leq \nu$, set 
$ \tau  _{i,j}
= 
p ^{*} _{j} ( t   _i)
-
p ^{*} _{j-1} ( t _i)
=
1 \otimes \cdots \otimes \tau _i \otimes \cdots \otimes 1$.
We have  the following  isomorphism of $\O _{X}$-$m$-PD-algebras
   \begin{align}
   \notag
   \O _{X} <T _{ij} ,\, 1\leq i \leq d, \, 1 \leq j \leq \nu > _{(m)}
   &\riso 
(\cP _{X, (m)} (\nu), \overline{\I} _{\nu}, \widetilde{\I} _\nu)\\
   \label{loc-desc-P}
T _{ij} &\mapsto
\tau  _{i,j},
   \end{align}
   where the structure of $\O _X$-module of 
$\cP _{X, (m)} (\nu)$
is given by its left structure. 
\end{prop}

\begin{proof}
We are in the situation of the proposition
\cite[1.6.6]{Caro-Vauclair}
where $u=\Delta$ and $f$ is the left projection
$p _0\colon X \times _{S} X ^{\nu} \to X$.
Hence, we get an isomorphism of the form 
\ref{loc-desc-P} where $\tau _{ij}$ is replaced by 
$p ^{*} _{j} ( t   _i)
-
p ^{*} _0 ( t _i)$.
Since
$p ^{*} _{j} ( t   _i)
-
p ^{*} _0 ( t _i)=
\tau _{i 0}
+
\tau _{i 1}
+ \dots
+ \tau _{i j}$, we are done. 
\end{proof}

\begin{cor}
\label{mPDenv-logpet}
Let $f \colon X \to Y$ be a $p$-étale morphism of $p$-smooth $T$-schemes. 
Then the canonical homomorphism
$f ^* \mathcal{P}_{Y/T,(m)} (\nu)
\to 
 \mathcal{P} _{X/T,(m)}  (\nu)$
is an isomorphism.
\end{cor}

\begin{proof}
Since this is local then we can suppose there exists a $p$-étale morphism of the form
$Y \to \A ^{d} _T$.
Hence,
this follows from 
\ref{prop-localdescr(nu)pre}.\end{proof}
\subsection{Level rise by Frobenius}
\label{subsecElevFrob}

Let $X$ be a $T$-scheme locally of formal finite type having locally finite $p$-basis.
Following \ref{rem-Dnn'Ynoeth}.\ref{rem-Dnn'Ynoeth-i)} and \ref{f0formétale-fforméta0}, $X/T$ is flat and $X$ is noetherian. 
Let $X _0 := X \times _{T} T _0$ 
and
$X _0  ^{(s)}$ be the base change of $X _0$ by the $s$-th power of the absolute Frobenius of $T _0$.
We denote by 
$F ^s _{X _0/T _0}\colon X _0 \to X _0 ^{(s)}$ the relative Frobenius.
We suppose there exists $X'$ a $T$-scheme locally of formal finite type having locally finite $p$-basis
such that $X ' \times _T T _0 \riso X _0 ^{(s)}$ (recall following  \ref{lifting-pbasis} such a lifting exists when $X/T$ has a finite $p$-basis).
Beware that, even locally, this is not clear that there exists  a lifting $F _T\colon T \to T$ of 
the absolute Frobenius 
$F ^{s} _{T _0}
\colon 
T _0 \to T _0$.

We suppose there exists a morphism 
$F \colon X \to X'$ which is a lifting of 
$F ^s _{X _0/T _0}$ 
(e.g. since $X'/T$ is formally smooth, such a lifting exists when $X$ is affine) 
i.e. $F$ makes commutative in the category of $T$-schemes the diagram 
\begin{equation}
\label{Fliftingdiag}
\xymatrix{
{X}
\ar[r] ^-{F} 
& 
{X'} 
\\ 
{X _0} 
\ar[r] ^-{F ^s _{X _0/T _0}}
\ar@{^{(}->}[u] ^-{}
& 
{X _0 ^{(s)}.} 
\ar@{^{(}->}[u] ^-{}
}
\end{equation}

\begin{ntn}
Let $\nu \geq 1$ be an integer. 
We will denote by 
$X ^{\nu +1}= X ^{\nu +1}/T$, 
$X ^{\prime \nu +1}= X ^{\prime \nu +1}/T$, 
$F _{\nu}
\colon 
X ^{\nu +1}
\to 
X ^{\prime \nu +1}$
the morphism 
$F \times \dots \times F$ induced by $F$, 
$\I _\nu$ (resp. $\I ' _{\nu}$) 
the ideal of the diagonal
$X \hookrightarrow X ^{\nu +1}$
(resp. 
$X '\hookrightarrow X ^{\prime \nu +1}$), 
and
$(\cP _{X, (m+s)} (\nu), \overline{\I} _{\nu}, \widetilde{\I} _\nu)$
(resp. 
$(\cP _{X', (m)} (\nu), \overline{\I} ' _{\nu}, \widetilde{\I} ' _\nu)$)
the $(m+s)$-enveloppe (resp. $m$-enveloppe)
of 
$\I _\nu$ (resp. $\I ' _{\nu}$).
When $\nu =1$, we simply write 
$(\cP _{X, (m+s)} , \overline{\I} , \widetilde{\I})$
(resp. $(\cP _{X', (m)} , \overline{\I} ' , \widetilde{\I}')$).

\end{ntn}

\begin{ntn}
[Local coordinates]
\label{loc-desc-m+s-m}
Suppose we have a relatively perfect morphism of the form $g _0\colon X _0 \to \A ^d _{T _0}$. 
Since 
$ (\A ^d _{T _0} )^{(s)}=\A ^d _{T _0}$, then we get by definition the right  cartesian squares of the diagram
\begin{equation}
\label{relperf-Frob}
\xymatrix@ R= 0,4cm{
{X _0} 
\ar[r] ^-{F ^s _{X _0/T _0}}
\ar[d] ^-{g _0}
\ar@{}[dr] ^-{} |\square
& 
{X _0 ^{(s)}} 
\ar[r] ^-{}
\ar[d] ^-{g _0 ^{(s)}}
\ar@{}[dr] ^-{} |\square
& 
{X _0}
\ar[d] ^-{g _0}
\\ 
{ \A ^d _{T _0}} 
\ar[r] _-{F ^s _{ \A ^d _{T _0}/T _0}}
& 
{\A ^d _{T _0} } 
\ar[r] ^-{}
\ar[d] ^-{}
\ar@{}[dr] ^-{} |\square
& 
{ \A ^d _{T _0}}
\ar[d] ^-{}
\\
& 
{T _0} 
\ar[r] ^-{F _{T _0}}
& 
 { T _0.}
}
\end{equation}
Since $g _0\colon X _0 \to \A ^d _{T _0}$ is relatively perfect, 
then the left square is also cartesian.

Choose coordinates $t _1,\dots, t _d$ inducing a lifting 
$g \colon X  \to \A ^d _{T}$ of $g _0$.
It follows from \ref{rem-dfn-pbasispadic}.\ref{rem-dfn-pbasispadic2)}.
that
$g$ is formally étale, i.e. $g$ is relatively perfect.

Choose coordinates $t  '_1,\dots, t '_d$ inducing a lifting 
$g '\colon X ' \to \A ^d _{T}$ of $g _0 ^{(s)}$.
We check similarly that $g'$ is flat and relatively perfect.

Set 
$\tau _i := 1 \otimes t _i - t _i \otimes 1\in \I _1$
and 
$\tau ' _i := 1 \otimes t '_i - t '_i \otimes 1 \in \I ' _{1}$. 
For any $j= 0,\dots, \nu$, let $p _j \colon X ^{ \nu+1}\to X$ be the index $j$ projection. 
For any $j= 0,\dots, \nu$, let $p ' _j \colon X ^{\prime \nu+1}\to X ^{\prime}$ be the index $j$ projection. 
For any $j= 1,\dots, \nu$, let 
$q _j = ( p _j, p _{j-1}) 
\colon 
X ^{ \nu+1}\to X ^2$,
and 
$q ' _j = ( p '_{j}, p ' _{j-1}) 
\colon 
X ^{\prime \nu+1}\to X ^{\prime 2}$.
For any $1\leq i \leq d ,\, 1 \leq j \leq \nu$, set 
$ \tau  _{i,j}
=
q ^{*} _j ( \tau _i)
= 
p ^{*} _{j} ( t   _i)
-
p ^{*} _{j-1} ( t _i)
=
1 \otimes \cdots \otimes \tau _i \otimes \cdots \otimes 1
$,
and 
$ \tau ' _{i,j}= 
q ^{\prime*} _{j} ( \tau ' _i)
=
p ^{\prime*} _{j} ( t ' _i)
-
p ^{\prime*} _{j-1} ( t ' _i)
=
1 \otimes \cdots \otimes \tau ' _i \otimes \cdots \otimes 1
$.

Since the above left square of \ref{relperf-Frob} is cartesian,
we get the commutative diagram
\begin{equation}
\label{relperf-Frob-X(1)pre}
\xymatrix@ R= 0,4cm
@ C=2cm{
{X _0} 
\ar@{^{(}->}[r] _-{}
\ar[d] ^-{g _0}
\ar@{}[dr] ^-{} |\square
& 
{X _0 \times _{X _0 ^{(s)}} X _0} 
\ar[d] ^-{g _0 \times _{g _0 ^{(s)}} g _0}
\ar@<0.5ex>[r]
\ar@<-0.5ex>[r]
\ar@{}[dr] ^-{} |\square
&
{X _0}
\ar[d] ^-{g _0} 
\\ 
{ \A ^d _{T _0}} 
\ar@{^{(}->}[r] _-{}
& 
{\A ^d _{T _0} \times _{ \A  _{T _0} ^{d (s)}} \A ^d _{T _0}} 
\ar@<0.5ex>[r]
\ar@<-0.5ex>[r]
&
{ \A ^d _{T _0}} 
}
\end{equation}
whose squares are cartesian.
Similarly, we get the cartesian square
\begin{equation}
\label{relperf-Frob-X(1)}
\xymatrix@ R= 0,4cm{
{X _0} 
\ar@{^{(}->}[r] _-{}
\ar[d] ^-{g _0}
\ar@{}[dr] ^-{} |\square
& 
{X _0 (\nu) } 
\ar[d] ^-{g _0 (\nu) }
\\ 
{ \A ^d _{T _0}} 
\ar@{^{(}->}[r] _-{}
& 
{\A ^d _{T _0} (\nu) } 
}
\end{equation}
where
$X _0 (\nu) 
:= 
X _0 \times _{X _0 ^{(s)}} 
\dots
\times _{X _0 ^{(s)}} X _0$ 
is the fibered product of $\nu +1$-copies of $X _0$ above $X _0 ^{(s)}$
and 
$\A ^d _{T _0} (\nu) := 
\A ^d _{T _0} 
\times _{ \A  _{T _0} ^{d (s)}}
\dots
\times _{ \A  _{T _0} ^{d (s)}}
 \A ^d _{T _0}$
 and
 $g _0 (\nu):= 
 g _0 \times _{g _0 ^{(s)}} 
\dots
\times _{g _0 ^{(s)}} g _0$ is the morphism induced by 
 $g _0$ and $g _0 ^{(s)}$.

\end{ntn}

\begin{prop}
[Local description of the $m$-PD-enveloppe]
\label{prop-localdescr(nu)}
Suppose we are in the local situation of \ref{loc-desc-m+s-m}. 
\begin{enumerate}[(a)]
\item We have  the following $\O _{X}$-$(m+s)$-PD isomorphism
   \begin{align}
   \notag
   \O _{X} <T _{ij} ,\, 1\leq i \leq d, \, 1 \leq j \leq \nu > _{(m+s)}
   &\riso 
(\cP _{X, (m+s)} (\nu), \overline{\I} _{\nu}, \widetilde{\I} _\nu)\\
   \label{loc-desc-P1}
T _{ij} &\mapsto
\tau  _{i,j},
   \end{align}
   where the structure of $\O _X$-module of 
$\cP _{X, (m+s)} (\nu)$
is given by the left structure. 
\item We have  the following $\O _{X}$-$m$-PD isomorphism
   \begin{align}
   \notag
   \O _{X} <T ' _{ij} ,\, 1\leq i \leq d, \, 1 \leq j \leq \nu > _{(m)}
   &\riso 
(\cP _{X', (m)} (\nu), \overline{\I} ' _{\nu}, \widetilde{\I} ' _\nu)\\
   \label{loc-desc-P2}
T '_{ij} &\mapsto
\tau  ' _{i,j},
   \end{align}
   where the structure of $\O _X$-module of 
$\cP _{X', (m)} (\nu)$
is given by the left structure. 
\end{enumerate}
\end{prop}

\begin{proof}
This follows from \ref{prop-localdescr(nu)pre}.
\end{proof}

\begin{empt}
\label{localdescItildenu}
Suppose we are in the local situation of \ref{loc-desc-m+s-m}

\begin{enumerate}[(a)]
\item Following \cite[1.5.1.(i)]{Be1} and \ref{prop-localdescr(nu)}, 
$\smash{\overline{\I} '}  ^{\{ n \} _{(m)}} _{\nu}$ is a free $\O _{X'}$-module and has the basis
$\prod _{i,j} \smash{\tau ' _{i,j}} ^{ \{ n _{i,j}\} _{(m)}}$ 
where $\sum _{i,j} n _{i,j} \geq n$. 
In particular, 
$\smash{\overline{\I} '}   _{\nu}$ is a free $\O _{X'}$-module and has the basis
$\prod _{i,j} \smash{\tau ' _{i,j}} ^{ \{ n _{i,j}\} _{(m)}}$ 
such that $\sum _{i,j} n _{i,j} \geq 1$. 
Moreover, 
$\widetilde{\I} ' _{\nu}$ is the ideal of 
$\cP _{X', (m)} (\nu)$
generated by 
$p\tau ' _{i,j}$ 
and $(\tau  _{i,j} ^{\prime p ^m}) ^{[q]}$ 
for $i= 1,\dots, d$,
$j = 1,\dots, \nu$,
and 
$q \geq 1$.

\item We have a similar description of 
$\overline{\I}  _{\nu}$,
$\smash{\overline{\I}  }^{\{ n \} _{(m+s)}} _{\nu}$,
and
$\widetilde{\I}  _{\nu}$ by removing some prime and replacing $m$ by $m+s$. 

\end{enumerate}

\end{empt}

\begin{lem}
\label{Idnu-I-link}
Suppose we are in the local situation of 
 \ref{loc-desc-m+s-m}.
 We denote by 
$g _\nu 
\colon 
X ^{\nu +1}
\to 
\A ^{d(\nu +1)} _{T}$
(resp. 
$g ' _\nu 
\colon 
X ^{\prime \nu +1}
\to 
\A ^{d(\nu +1)} _{T}$)
the morphism
$g\times g \times \dots \times g$ (resp. $g'\times g '\times \dots \times g'$)
induced by $g$ (resp. $g'$).
Let $\I _{d,\nu}$ be the ideal of the diagonal
$\A ^d _{T} \hookrightarrow \A ^{d(\nu +1)} _{T}$.

\begin{enumerate}[(a)]
\item 
The canonical homomorphisms
$$g _\nu ^{ *}  \cP _{\A ^d _{T}, (m+s)} (\nu)
\to 
\cP _{(m+s),\alpha} ( \I _{d,\nu} \O _{X ^{\nu +1}})
\to 
\cP _{X, (m+s)} (\nu)$$
are isomorphisms.

\item The canonical homomorphisms
$$g _\nu ^{\prime *} \cP _{\A ^d _{T}, (m)} (\nu)
\to
\cP _{(m),\alpha} ( \I _{d,\nu} \O _{X ^{\prime \nu +1}})
\to 
\cP _{X', (m)} (\nu)$$
are isomorphisms.
\end{enumerate}

\end{lem}

\begin{proof}
Since $g _{\nu}$ is flat, then following \cite[1.4.6]{Be1} 
the canonical homomorphism
$$g _\nu ^* \cP _{\A ^d _{T}, (m+s)} (\nu)
\to 
\cP _{(m+s),\alpha} ( \I _{d,\nu} \O _{X ^{\nu +1}})$$
is an isomorphism.
Since $g _\nu$ is relatively perfect and then $p$-étale, then following \ref{mPDenv-logpet} the homomorphism
$g _\nu ^* \cP _{\A ^d _{T}, (m+s)} (\nu)
\to 
\cP _{X, (m+s)} (\nu)$ is an isomorphism.
This yields that the homomorphism
$\cP _{(m+s),\alpha} ( \I _{d,\nu} \O _{X ^{\nu +1}})
\to 
\cP _{X, (m+s)} (\nu)$ is an isomorphism.

Similarly, we get the second part of the lemma.\end{proof}

\begin{empt}
\label{ntnIf}
Since the left top square of \ref{relperf-Frob} is commutative, 
similarly to \cite[2.2.2.1]{Be2}, we compute 
$F ^* _1 ( \tau ' _i ) = \tau ^{p ^s} _i + \zeta _i$, with $\zeta _i \in \mathfrak{a} \I$.
Let $\I _f$ be a finitely generated ideal of $\O _{X ^{2}}$ which contains 
$\I _{d,1} \O _{X ^{2}}$,
which is contained in $\I$ and such that 
$\zeta _i \in \mathfrak{a} \I _f $.
Similarly to \cite[2.2.2.2]{Be2}, this yields
\begin{equation}
\label{2222Be2}
F ^* _1 ( \tau  _i ^{\prime p ^m} ) = \tau ^{p ^{m+s}} _i + \sigma _i, 
\end{equation}
with $\sigma _i \in \mathfrak{b} \I _f ^{p ^m}$.
Remark that to prove some local theorems, 
it might be possible to reduce to the case where
$\I _f=\I _{d,1} \O _{X ^{2}}$ (see \ref{rem-Frob-desc-simple} below).
To make things as general as possible we do not make such assumption.
\end{empt}

\begin{rem}
\label{rem-Frob-desc-simple}
As in \ref{loc-desc-m+s-m}, 
suppose we have a relatively perfect morphism of the form $g _0\colon X _0 \to \A ^d _{T _0}$. 
Choose coordinates $t _1,\dots, t _d$ (resp. $t  '_1,\dots, t '_d$)
inducing a lifting 
$g \colon X  \to \A ^d _{T}$ of $g _0$
(resp. $g '\colon X ' \to \A ^d _{T}$ of $g _0 ^{(s)}$).

Suppose fixed the coordinates
$x _1,\dots, x _d$ of 
$\mathbb{A} ^d _{T}$.
We denote by 
$F \colon \mathbb{A} ^d _{T} \to \mathbb{A} ^d _{T}$
the canonical $T$-morphism given by 
$x _i \to x _i ^{p^s}$, for any $i =1,\dots, d$.
Since 
$X'
\to \mathbb{A} ^d _{T}
$
is relatively perfect,
using the universal property of a formally étale morphism,
there exists a unique morphism
$F \colon X \to X'$ which is a lifting of 
$F ^{s} _{X_0/T _0}$ and making commutative the diagram
\begin{equation}
\label{F!}
\xymatrix@R=0,3cm{
{X} 
\ar@{.>}[d] ^-{F}
\ar[r] ^-{t _1,\dots,t _d}
& 
{ \mathbb{A} ^d _{T}} 
\ar[d] ^-{F}
\ar[rd] ^-{}
\\ 
{X'} 
\ar[r] ^-{t '_1,\dots,t '_d}
& 
{\mathbb{A} ^d _{T}} 
\ar[r] ^-{}
& 
{T .} 
}
\end{equation}
Hence, we get the  relation
$F ^* ( t ' _i) =  t _i ^{p ^s}$.
This yields
$F _1 ^* ( \tau '_i ) = \tau _i ^{p ^s} + \sum _{j = 1} ^{p ^s -1} 
\left (
\begin{smallmatrix}
p ^s  \\
j
\end{smallmatrix}
\right)
t _i ^{p^s -j} \tau _i ^j$.
Hence, 
we can choose in this case 
$\I _f=\I _{d,1} \O _{X ^{2}}$ (see notation \ref{ntnIf}).

Since the square of \ref{F!} is cartesian modulo $\fa$, since its horizontal arrows are formally étale,
then the square of \ref{F!} is cartesian.
Moreover, similarly to \ref{relperf-Frob-X(1)}
we get the cartesian square
\begin{equation}
\label{relperf-Frob-X(1)bis}
\xymatrix@ R= 0,4cm{
{X } 
\ar@{^{(}->}[r] _-{}
\ar[d] ^-{g }
\ar@{}[dr] ^-{} |\square
& 
{X  (\nu) } 
\ar[d] ^-{g  (\nu) }
\\ 
{ \A ^d _{T }} 
\ar@{^{(}->}[r] _-{}
& 
{\A ^d _{T } (\nu) } 
}
\end{equation}
where
$X  (\nu) 
:= 
X  \times _{X  ^{\prime}} 
\dots
\times _{X  ^{\prime}} X $ 
the fibered product of $\nu +1$-copies of $X$ above $X'$ 
and 
$\A ^d _{T } (\nu) := 
\A ^d _{T } 
\times _{ \A  _{T } ^{d (s)}}
\dots
\times _{ \A  _{T } ^{d(s)}}
 \A ^d _{T }$
 and
 $g  (\nu)$ is the morphism induced by 
 $g $ and $g  ^{\prime}$.
\end{rem}

\begin{lem}
\label{F-locfreeft}
The morphism $F\colon X \to X'$ is finite and 
is localement free 
of finite type, i.e.
$F _* \O _X$ is a locally free $\O _{X'}$-module
of finite type. 
\end{lem}

\begin{proof}
Since this is local, we come down to the local context of \ref{loc-desc-m+s-m}.
Since $F ^s _{ {\A} ^d _{T _0}/T _0}$ is locally free of rank $d s$,
using the cartesian left square of \ref{relperf-Frob}, 
then so is $F ^s _{X _0/T _0}$.
Since $X/T$ is flat and $X$ is noetherian, then using 
\cite[11.3.10]{EGAIV3}, we get that $F$ is flat.
Since $\mathfrak{a}$ is nilpotent and $X$ is noetherian, this yields that 
$F$ is also finite and then $F$ is free of rank $d s$. 
\end{proof}

\begin{rem}
In another context, we have a similar to  \ref{F-locfreeft} lemma  in \ref{F-locfreeft2} .
\end{rem}

\begin{lem}
\label{222Be2}
\begin{enumerate}[(a)]
\item There exists a unique PD-morphism
$$\Phi ^* _\nu
\colon 
F ^{-1} _\nu \cP _{X', (m)} (\nu)
\to 
\cP _{X, (m+s)} (\nu)$$
sending 
$F ^{-1} _\nu \widetilde{\I} ' _{\nu}$
to 
$ \widetilde{\I}  _{\nu} + \mathfrak{b} \overline{\I}  _{\nu}$.
This yields the morphism
$\Phi _\nu 
\colon 
\Delta _{X,(m+s)} (\nu)
\to 
\Delta _{X',(m)} (\nu)$.

\item For any $n \in \N$, we have the inclusion
$$\Phi ^* _\nu ( F ^{-1} _\nu \overline{\I}  _{\nu} ^{\prime \{ n\} _{(m)}})
\subset 
\overline{\I}  _{\nu} ^{\{ n\} _{(m+s)}}.$$
\end{enumerate}
\end{lem}

\begin{proof}
Since this is local, we come down to the local context of \ref{loc-desc-m+s-m}.
Recall (see \cite[1.4.1]{Be1}) that 
$\cP _{(m),\alpha} ( \I _{d,\nu} \O _{X ^{\prime \nu +1}})$
is the PD-enveloppe of 
$( \I _{d,\nu} \O _{X ^{\prime \nu +1}}) ^{(p ^m)} + (\fb +p \O _T) \O _{X ^{\prime \nu +1}}$
and 
$\cP _{(m+s),\alpha} ( \I _{d,\nu} \O _{X ^{\nu +1}})$
is the PD-enveloppe of the ideal
$( \I _{d,\nu} \O _{X ^{\nu +1}}) ^{(p ^{m+s})} + (\fb +p \O _T) \O _{X ^{\nu +1}}$.
Moreover, the ideal 
$\I _{d,\nu} \O _{X ^{\nu +1}}$ 
(resp. 
$\I _{d,\nu} \O _{X ^{\prime \nu +1}}$)
is generated by the sections
$\tau  _{i,j}$
(resp. $\tau  _{i,j} ^{\prime })$
for 
$i= 1,\dots, d,
\
j = 1,\dots, \nu$.
Using the formula \ref{2222Be2}, this yields that the image of 
$( \I _{d,\nu} \O _{X ^{\prime \nu +1}}) ^{(p ^m)}$ 
via 
$F ^{-1} _\nu  \O _{X ^{\prime \nu +1}}
\to 
\O _{X ^{\nu +1}}$
is included in
$( \I _{d,\nu} \O _{X ^{\nu +1}}) ^{(p ^{m+s})} + (\fb +p \O _T) \O _{X ^{\nu +1}}$.
Using Lemma \ref{Idnu-I-link}
and using  the universal property of the PD-enveloppe, 
we get uniquely a PD-morphism
$\Phi ^* _\nu
\colon 
F ^{-1} _\nu \cP _{X', (m)} (\nu)
\to 
\cP _{X, (m+s)} (\nu)$.
To check the inclusion, using the local description 
\ref{localdescItildenu}, we reduce to the case $\nu = 1$.
Then, using \ref{2222Be2}, we conclude 
(for more details, we can copy the proof of \cite[2.2.2]{Be2}).\end{proof}

\subsection{Frobenius descent for left $\cD$-modules}
\label{subsecFrobDesc}
We keep notation and hypotheses of \ref{subsecElevFrob}.

\begin{lem}
Let $\nu \geq 0$. The commutative square
\begin{equation}
\label{Be2-232}
\xymatrix{
{\Delta _{X,(m+s)} (\nu)} 
\ar[d] ^-{\Phi _\nu }
\ar[r] ^-{}
& 
{X ^{\nu +1}} 
\ar[d] ^-{F _\nu }
\\ 
{\Delta _{X',(m)} (\nu)} 
\ar[r] ^-{}
& 
{X ^{\prime \nu +1}}
}
\end{equation}

is cartesian.
\end{lem}

\begin{proof}
We follow the ingredients of the proof of \cite[2.3.2]{Be2} : 

1) Since this is local we come down to the local context of \ref{loc-desc-m+s-m}. 
Set 
$\fI ' _{\nu}:= \I _{d,\nu} \O _{X ^{\prime \nu +1}}$,
and 
$\fI  _{\nu}
:=
\I _{d,\nu} \O _{X ^{ \nu +1}}$ (beware it depends on the choice of the finite $p$-basis).
Since $F _\nu$ is flat, then the canonical homomorphism
$$\O _{X ^{\nu +1}}
\otimes _{\O _{X ^{\prime \nu +1}}}
\cP _{(m),\alpha} ( \fI ' _{\nu})
\to 
\cP _{(m),\alpha} ( \fI ' _{\nu}  \O _{X ^{ \nu +1}})$$
is an isomorphism.

2) By construction (see \cite[1.4.1]{Be1}),
$\cP _{(m),\alpha} ( \fI ' _{\nu}  \O _{X ^{ \nu +1}})$
is the PD-enveloppe of
$\fI  _{\nu}   ^{\prime (p^m)} \O _{X ^{ \nu +1}}+ \mathfrak{b} _1 \O _{X ^{\nu +1}}$, with compatibility with respect to the PD-structure of 
$\mathfrak{b} _1$.
On the other hand, by construction, 
$\cP _{(m+s),\alpha} (\fI  _{\nu})$
is the PD-enveloppe of
$\fI  _{\nu} ^{(p ^{m+s})} + \mathfrak{b} _1 \O _{X ^{\nu +1}}$, with compatibility with respect to the PD-structure of 
$\mathfrak{b} _1$.

3) The ideal $\fI  ' _{\nu}$ is generated by the sections 
$ \tau ' _{i,j}= q ^{\prime*} _j ( \tau  ' _i),
\
i= 1,\dots, d,
\
j = 1,\dots, \nu.$
Hence, 
$ \fI ' _{\nu}  \O _{X ^{ \nu +1}}$ is the ideal generated by
$F ^* _\nu( \tau ^{\prime p ^m } _{i,j})
= q ^{*} _j (F ^* _1 (  \tau  ^{\prime p ^m }  _i))$,
for 
$i= 1,\dots, d,
\
j = 1,\dots, \nu.$
Using the formula \ref{2222Be2}, 
we get 
$q ^{*} _j (F ^* _1 (  \tau  ^{\prime p ^m }  _i)) 
=
q ^{*} _j (  \tau  ^{p ^{m+s} }  _i) 
 + q ^{*} _j (\sigma _i)$, 
 with 
 $q ^{*} _j (\sigma _i) \in \mathfrak{b} _1 \O _{X ^{\nu +1}}$.
This implies 
$ \fI  ^{\prime (p^m)} _{\nu}  \O _{X ^{ \nu +1}}  + \mathfrak{b} _1 \O _{X ^{\nu +1}}
=
\fI  _{\nu} ^{(p ^{m+s})} + \mathfrak{b} _1 \O _{X ^{\nu +1}}$.

4) Using 2) and 3) we get 
$\cP _{(m),\alpha} ( \fI ' _{\nu}  \O _{X ^{ \nu +1}})
=
\cP _{(m+s),\alpha} (\fI  _{\nu})$.
Since 
$\cP _{(m+s),\alpha} (\fI  _{\nu})
\riso\cP _{X, (m+s)} (\nu)$
and 
$\cP _{(m),\alpha} ( \fI ' _{\nu}  \O _{X ^{ \nu +1}})
\riso\cP _{X', (m)} (\nu)$ (see Lemma \ref{Idnu-I-link}),
then we conclude by using 1).
\end{proof}

\begin{lem}
\label{233Be2}
Let $\nu \geq 0$ be an integer; 
$X (\nu) := X \times _{X'}\cdots _{X'} X$ the fibered product of $\nu +1$-copies of $X$ above $X'$, 
$\K _\nu$ be the ideal defined by the diagonal immersion $X \hookrightarrow X (\nu)$.
Then 
$\fb _1 \O _{X (\nu)} \cap \K _\nu$
is a sub-PD-ideal of 
$\fb _1 \O _{X (\nu)}$,
and endows 
$\K _\nu$ with a canonical nilpotent $(m+s)$-PD-structure, 
compatible with $\alpha$.
\end{lem}

\begin{proof}
1) We follow the ingredients of the proof of \cite[2.3.2]{Be2} : 
Since $\O _X = \O _{X (\nu)}/\K _\nu$ is $\O _T$-flat then following 
\cite[Corollary of Proposition 7 of 1 of I.2.6]{bourbaki}
we get
$\fb _1 \O _{X (\nu)} \cap \K _\nu 
=
\fb _1 \K _\nu $.
Since $X (\nu)/T$ is flat, then we get 
a canonical PD-structure on $\fb _1  \O _{X (\nu)}$
extending that of $\fb _1$.
This yields that  $\fb _1 \K _\nu $ is a sub-PD-ideal of $\fb _1  \O _{X (\nu)}$.
Hence, it remains to check 
$\K _\nu ^{(p ^{m+s})}
\subset 
\fb _1 \O _{X (\nu)} \cap \K _\nu $
and that the corresponding $(m+s)$-PD-structure is nilpotent.

2) Since the lemma is local we come down to the local context of \ref{loc-desc-m+s-m}. 
Following \ref{relperf-Frob-X(1)} and with its notation, the left square of the commutative diagram 
\begin{equation}
\label{relperf-Frob-bis}
\xymatrix@ R= 0,4cm{
{X _0} 
\ar@{^{(}->}[r] _-{}
\ar[d] ^-{g _0}
\ar@{}[dr] ^-{} |\square
& 
{X _0 (\nu)} 
\ar[d] ^-{g _0 (\nu)}
\ar@{^{(}->}[r] _-{}
&
\ar@{^{(}->}[r] _-{}
{X (\nu)}
&
{X ^{\nu+1} }
\ar[d] ^-{g _\nu } 
\\ 
{ \A ^d _{T _0}} 
\ar@{^{(}->}[r] _-{}
& 
{\A ^d _{T _0} (\nu)} 
\ar@{^{(}->}[rr] _-{}
&&
{ \A  _{T}  ^{d (\nu+1)}} 
}
\end{equation}
is cartesian.
Let 
$\overline{\K} _{d,\nu}$ be the ideal of  the diagonal immersion $\A ^{d} _{T _0} \hookrightarrow \A ^{d} _{T _0} (\nu)$,
and
$\overline{\K} _\nu$
be the ideal of  the diagonal immersion $X _{0} \hookrightarrow X _{0} (\nu)$. 
Then we get
$\overline{\K} _\nu
=
( \K _\nu + \fa \O _{X (\nu)}) \O _{X _0 (\nu)}$,
and 
$\overline{\K} _{d,\nu}
=
\left (
\I _{d,\nu}+ \fa \O _{\A  _{T} ^{d (\nu+1)}} 
\right ) \O _{\A ^{d} _{T _0} (\nu)}$,
where 
according to notation \ref{loc-desc-m+s-m}
the sheaf $\I _{d,\nu}$ is the ideal of the diagonal
$\A ^d _{T} \hookrightarrow \A ^{d(\nu +1)} _{T}$.
Since the left square of \ref{relperf-Frob-bis}
is cartesian, this yields
$\overline{\K} _\nu
= 
\I _{d,\nu}  \O _{X _0(\nu)}+ \fa \O _{X _0 (\nu)}$
and then 
$$\K _\nu + \fa \O _{X (\nu)}
= 
\I _{d,\nu}  \O _{X (\nu)}+ \fa \O _{X (\nu)} .$$
Let $r _j \colon X (\nu) \to X ^{2}$ be the index $j-1$ and $j$ projection. 
The ideal $\I _{d,\nu} \O _{X (\nu)}$ 
is generated by the sections 
$
 r ^{*} _j ( \tau   _i),
\
i= 1,\dots, d,
\
j = 1,\dots, \nu.$
Hence, $\I ^{(p ^{m+s})}_{d,\nu}\O _{X (\nu)}$ is generated by 
$r ^{*} _j (\tau ^{p ^{m+s}} _i) ,
\
i= 1,\dots, d,
\
j = 1,\dots, \nu.$
Using \ref{2222Be2},
we get
$0= r ^{*} _j (F ^* _1 ( \tau  _i ^{\prime p ^m} )) = 
 r ^{*} _j (\tau ^{p ^{m+s}} _i) +  r ^{*} _j (\sigma _i)$,
with $\sigma _i \in \mathfrak{b} \I$.
Hence, 
$r ^{*} _j (\tau ^{p ^{m+s}} _i) 
\in 
\fb _1 \O _{X (\nu)} \cap \K _\nu $.
This yields 
$\I ^{(p ^{m+s})}_{d,\nu}\O _{X (\nu)}
\subset 
\fb _1 \O _{X (\nu)} \cap \K _\nu $.
Since
$\K _\nu 
\subset
\I _{d,\nu}  \O _{X (\nu)}+ \fa \O _{X (\nu)}$, 
this implies
$\K _\nu  ^{(p ^{m+s})}
\subset
\fb _1 \O _{X (\nu)} \cap \K _\nu $.
As for \cite[2.3.2]{Be2}, we check 
the $(m+s)$-PD-nilpotence.
\end{proof}

\begin{lem}
\label{234Be2}
Let $\Phi ^* _\nu
\colon 
 F ^{-1} _\nu \cP _{X', (m)} (\nu)
\to 
\cP _{X, (m+s)} (\nu)$
be the factorization of \ref{222Be2}. 
For any $n \in \N$, there exists an integer $n' \geq n$ 
(only depending on $T$, $n$, $m$, $s$, $\nu$ and the cardinality $d$ of 
the finite $p$-basis $X/T$)
such that
\begin{equation}
\notag
\overline{\I} _\nu ^{\{n' \} _{(m+s)}}
\subset
\Phi ^* _\nu (\overline{\I} _\nu ^{\prime \{n \} _{(m)}})
\cP _{X, (m+s)} (\nu).
\end{equation}

\end{lem}

\begin{proof}
Since the lemma is local we come down to the local context of \ref{loc-desc-m+s-m}. 
Using 
$\cP _{(m+s),\alpha} ( \I _{d,\nu} \O _{X ^{\nu +1}})
\riso
\cP _{X, (m+s)} (\nu)$
and 
$\cP _{(m),\alpha} ( \I _{d,\nu} \O _{X ^{\prime \nu +1}})
\riso\cP _{X', (m)} (\nu)$
and 
the relation 
\ref{2222Be2},
we can follow the proof of 
\cite[2.3.4]{Be2} where we replace $\I$ by 
the finitely generated ideal $\I _f $
(more precisely the only change is the following :  since $p$ is nilpotent and $\I _f $ is a finitely generated ideal then, 
still denoting by $\sigma _i$ its image in
$ \mathfrak{b} \I _f ^{p ^m} \cP _{X, (m+s)}$, we get 
$\sigma _i ^{[N]}=0$ for $N$ large enough).
\end{proof}

\begin{lem}
\label{235Be2}
Let $n \in \N$ be an integer, 
$\Delta : = \Delta _{X,(m+s),\alpha}$,
$\Delta ': = \Delta _{X',(m),\alpha}$,
$\Delta ^{\prime n} : = \Delta _{X',(m),\alpha} ^n $
and 
$\Delta ^n: = \Delta ^{\prime n}  \times _{\Delta'} \Delta$.
\begin{enumerate}[(a)]

\item The $(m+s)$-PD-structure of the ideal $\overline{\I}$ induces a nilpotent  $(m+s)$-PD-structure compatible with $\alpha$ on the ideal of the immersion
$X \hookrightarrow \Delta ^n$
such that 
$\Delta ^n \to \Delta$ is an $(m+s)$-PD-morphism.

\item The ideal $\mathcal{J}$ (resp. $\mathcal{J} _n$) of the diagonal immersion
$X \hookrightarrow \Delta  \times _{\Delta ^{\prime} } \Delta $
(resp. $X \hookrightarrow \Delta ^n \times _{\Delta ^{\prime n} } \Delta ^{n} $)
is canonically endowed with a 
$(m+s)$-PD-structure 
(resp. nilpotent $(m+s)$-PD-structure )
compatible with $\alpha$, such that
both projections
$\Delta  \times _{\Delta ^{\prime} } \Delta \to \Delta $
(resp. $\Delta ^n \times _{\Delta ^{\prime n} } \Delta ^{n} \to \Delta ^n$)
are $(m+s)$-PD-morphisms.
\end{enumerate}

\end{lem}

\begin{proof}
Following \cite[1.3.4]{Be1}, the first assertion (except the nilpotence) is equivalent to the property 
that $\widetilde{\I}\cap (\overline{\I}  ^{\prime \{n \} _{(m)}}\cP _{X, (m+s)})$ is a sub-PD-ideal of 
$\widetilde{\I}$. We can copy the proof of \cite[2.3.5]{Be2}. Similarly for the second result.
\end{proof}

\begin{theo}
\label{thm236Be2}
The functor $F ^*$ is an equivalence between the category of 
left (resp. quasi-coherent) $\D ^{(m)} _{X'/T}$-modules 
and that of 
left (resp. quasi-coherent) $\D ^{(m+s)} _{X/T}$-modules.
\end{theo}

\begin{proof}
Using Lemmas \ref{222Be2}, \ref{233Be2}, \ref{234Be2}, \ref{235Be2}
we can copy the proof of Theorem 
\cite[2.3.6]{Be2}.
\end{proof}

\begin{cor}
[Homological dimension]
\label{finitecohdimDdagpre}
Suppose $T$ is affine and regular, 
$f\colon X\to T$ is affine, locally of formal finite type, has finite $p$-basis.
Suppose  the fibers of 
$f\colon X\to T$ are of dimension $d$. Let $r:= \sup _{t \in f(X) } \O _{T,t}$.
Then for any integer $m \in \N$, the ring $D _{X/T} ^{(m)} := \Gamma (X, \D _{X/T} ^{(m)})$ has homological dimension
equal to $2d +r$. 
\end{cor}

\begin{proof}
Using the Frobenius descent Theorem \ref{thm236Be2},
we reduce to the case $m=0$. Then, this is standard (see \cite[4.4.3]{Be2}).
\end{proof}

\begin{cor}
\label{finitecohdimDdag}
Let $\X$ be an affine formal $\fS$-scheme
locally of formal finite type
and having locally finite $p$-bases  over $\fS$.
Then, we have the following properties. 
\begin{enumerate}[(a)]
\item For any integer $m \in \N$, the ring $\widehat{D} _{\fX} ^{(m)} := \Gamma (\fX, \widehat{\D} _{\X} ^{(m)})$ has homological dimension
equal to $2d +1$. 

\item The ring $D ^\dag _{\fX,\Q}  := \Gamma (\fX, \D ^\dag _{\fX,\Q}  )$ has homological dimension
equal to $d''$ with $d \leq d '' \leq 2d +1$. 
\end{enumerate}

\end{cor}

\begin{proof}
We can copy the proof of \cite[4.4.7]{Be2}.
\end{proof}

\begin{cor}
\label{thm236Be2-cordag}
Let $\X$ be a formal $\fS$-scheme
locally of formal finite type
and having locally finite $p$-bases  over $\fS$.
Let $X _0$ be its special fiber and 
$X _0  ^{(s)}$ be the base change of $X _0$ by the $s$-th power of the absolute Frobenius of $S _0$.
Suppose there exists 
$F \colon \X \to \X ^{\prime}$ a morphism of 
formal $\fS$-schemes
locally of formal finite type
and having locally finite $p$-bases  over $\fS$
which is a lifting of the relative Frobenius
$F ^s _{X _0/S _0}\colon X _0 \to X _0 ^{(s)}$.
Then 
$F ^* $ induces an equivalence between the category of 
left $\D ^\dag _{\fX ^{\prime},\Q} $-modules 
and that of 
left $\D ^\dag _{\fX,\Q} $-modules.

\end{cor}

\subsection{Frobenius descent for right $\D$-modules}
We keep notation and hypotheses of \ref{subsecElevFrob}.

\begin{empt}
Following \ref{F-locfreeft}, 
$F _*\cO _X$ is an $\cO _{X '}$-module of finite type. 
In fact, 
since $F _*$ is the identity, we get a structure of 
$\cO _{X '}$-module on $\cO _{X}$ via $F$.
Since $F$ is supposed to be fixed, 
we simply write $\cO _X$ instead of $F _*\cO _X$.
For any $\O _{X'}$-module $\cM'$, this yields the isomorphism
$$F ^{\flat} \cM '
=
\bbR\cH om _{\cO _{X '}} ( \cO _X, \cM ')
\riso 
\cH om _{\cO _{X '}} ( \cO _X, \cM ').$$
\end{empt}

\begin{prop}
\label{241Be2}
Let $\cM '$ be a right $ \D ^{(m)} _{X'/T}$-module (resp. 
a $ \D ^{(m)} _{X'/T}$-bimodule etc.).
Then  $F ^{\flat} \cM '$ is canonically endowed with a structure of 
right $ \D ^{(m+s)} _{X/T}$-module (resp. 
a $(\D ^{(m)} _{X'/T} , \D ^{(m+s)} _{X/T})$-bimodule etc.)
\end{prop}

\begin{proof}
Let us check the non respective case. 
Following \ref{prop-rightD-costrat}, $\cM'$ has a structural 
 $m$-PD-costratifcation. By applying $F ^{\flat}$ and using  \ref{222Be2}, 
 we get a canonical structure of  $(m+s)$-PD-costratifcation on
 $F ^{\flat} \cM '$, i.e.  $F ^{\flat} \cM '$ is endowed with a structure 
 of right $ \D ^{(m+s)} _{X/T}$-module. 
 By functoriality, we get the respective case from the non respective case.
\end{proof}

\begin{thm}
\label{thm236Be2-cor-cohrightpre}
The functor $F ^{\flat} $ is an equivalence between the category of 
right (resp. quasi-coherent) $\D ^{(m)} _{X'/T}$-modules 
and that of 
left (resp. quasi-coherent) $\D ^{(m+s)} _{X/T}$-modules.
\end{thm}

\begin{proof}
By exchanging ``stratifications'' by 
``costratification'', by exchanging the functors of form $f ^*$
by $f ^\flat$, 
this is just a matter of copying the proof of Berthelot of \ref{thm236Be2}.
\end{proof}

\begin{cor}
\label{thm236Be2-cordag-cohrightpre}
Let $\X$ be a formal $\fS$-scheme locally of formal finite type
and having locally finite $p$-bases over $\fS$.
Let $X _0$ be its special fiber and 
$X _0  ^{(s)}$ be the base change of $X _0$ by the $s$-th power of the absolute Frobenius of $S _0$.
Suppose there exists 
$F \colon \X \to \X ^{\prime}$ a morphism of 
formal $\fS$-schemes
locally of formal finite type
and having locally finite $p$-bases  over $\fS$
which is a lifting of the relative Frobenius
$F ^s _{X _0/S _0}\colon X _0 \to X _0 ^{(s)}$.
The functor
$F ^{\flat} $ induces an equivalence between the category of 
right $\D ^\dag _{\fX ^{\prime},\Q} $-modules 
and that of 
right $\D ^\dag _{\fX,\Q} $-modules.

\end{cor}

\subsection{Quasi-inverse functor}
We keep notation and hypotheses of \ref{subsecElevFrob}.

\begin{prop}
\label{252Be2}
There exists an isomorphism of 
$\D ^{(m+s)} _{X/T}$-bimodules of the form
\begin{equation}
\label{252Be2-iso}
\D ^{(m+s)} _{X/T}
\riso
F ^* F ^{\flat} \D ^{(m)} _{X'/T} .
\end{equation}
\end{prop}

\begin{proof}
Using \ref{222Be2} and 
\ref{234Be2}, 
we can copy the proof of \cite[2.5.2]{Be2}.
\end{proof}

\begin{coro}
\label{253Be2}
\begin{enumerate}[(a)]
\item The
$\D ^{(m+s)} _{X/T}$-modules 
$F ^*  \D ^{(m)} _{X'/T} $
and
$F ^{\flat} \D ^{(m)} _{X'/T} $
are locally projective of finite type

\item A left (resp. right) $\D ^{(m)} _{X'/T} $-module $\cE'$ (resp. $\cM'$)
is coherent if and only if 
$F ^* ( \cE ') $
(resp. $F ^{\flat} \cM '$) is $\D ^{(m+s)} _{X/T}$-coherent.

\end{enumerate}
\end{coro}

\begin{proof}
We can copy the proof of \cite[2.5.3]{Be2}.
\end{proof}

\begin{cor}
\label{thm236Be2-cordag-cohright}
Let $\X$ be a formal $\fS$-scheme locally of formal finite type
and having locally finite $p$-bases over $\fS$.
Let $X _0$ be its special fiber and 
$X _0  ^{(s)}$ be the base change of $X _0$ by the $s$-th power of the absolute Frobenius of $S _0$.
Suppose there exists 
$F \colon \X \to \X ^{\prime}$ a morphism of 
formal $\fS$-schemes
locally of formal finite type
and having locally finite $p$-bases  over $\fS$
which is a lifting of the relative Frobenius
$F ^s _{X _0/S _0}\colon X _0 \to X _0 ^{(s)}$.
\begin{enumerate}[(a)]
\item The functor
$F ^* $ induces an equivalence between the category of 
(coherent) left $\D ^\dag _{\fX ^{\prime},\Q} $-modules 
and that of 
(coherent) left $\D ^\dag _{\fX,\Q} $-modules.
\item The functor
$F ^{\flat} $ induces an equivalence between the category of 
(coherent) right $\D ^\dag _{\fX ^{\prime},\Q} $-modules 
and that of 
(coherent) right $\D ^\dag _{\fX,\Q} $-modules.

\end{enumerate}

\end{cor}

\begin{proof}
The first (resp. second) statement 
is a consequence of \ref{thm236Be2-cordag} (resp. \ref{thm236Be2-cordag-cohrightpre})
and of \ref{253Be2}.
\end{proof}

\begin{coro}
\label{2567Be2}
Let $\cE'$ be a left $\D ^{(m)} _{X'/T} $-module.
Let $\cM'$ be a right $\D ^{(m)} _{X'/T} $-module.
\begin{enumerate}[(a)]
\item We have the functorial isomorphisms 
\begin{equation}
\label{256Be2-iso1Pre}
F ^{\flat} \D ^{(m)} _{X'/T}
\otimes ^{\bbL}_{\D ^{(m+s)} _{X/T}}
F ^* \cE '
\riso 
F ^{\flat} \D ^{(m)} _{X'/T}
\otimes _{\D ^{(m+s)} _{X/T}}
F ^* \cE '
\riso 
\cE'.
\end{equation}
\item 
We have the functorial isomorphisms 
\begin{equation}
\label{256Be2-iso2Pre}
F ^{\flat} \cM '
\otimes ^{\bbL} _{\D ^{(m+s)} _{X/T}}
F ^* \D ^{(m)} _{X'/T}
\riso 
F ^{\flat} \cM '
\otimes _{\D ^{(m+s)} _{X/T}}
F ^* \D ^{(m)} _{X'/T}
\riso 
\cM'.
\end{equation}

\item Denoting by $f \colon X \to T$ the structural morphism, 
we have the functorial isomorphism in $D (f ^{-1}\O _T)$
\begin{equation}
\label{256Be2-iso3Pre}
F ^{\flat} \cM '
\otimes ^{\bbL} _{\D ^{(m+s)} _{X/T}}
F ^* \cE'
\riso 
 \cM '
\otimes ^{\bbL} _{\D ^{(m)} _{X'/T}}
 \cE'.
\end{equation}
\end{enumerate}
\end{coro}

\begin{proof}
We can copy the proof of \cite[2.5.6-7]{Be2}.
\end{proof}

\begin{coro}
\label{F+invF*-form}
Let $\X$ be a formal $\fS$-scheme locally of formal finite type
and having locally finite $p$-bases  over $\fS$.
Let $X _0$ be its special fiber and 
$X _0  ^{(s)}$ be the base change of $X _0$ by the $s$-th power of the absolute Frobenius of $S _0$.
Suppose there exists 
$F \colon \X \to \X ^{\prime}$ a morphism of 
formal $\fS$-schemes
locally of formal finite type
and having locally finite $p$-bases  over $\fS$
which is a lifting of the relative Frobenius
$F ^s _{X _0/S _0}\colon X _0 \to X _0 ^{(s)}$.
\begin{enumerate}[(a)]
\item The functor
$F _+ := F ^\flat \D ^\dag _{\fX ^{\prime},\Q}  \otimes _{\D ^\dag _{\fX ^{\prime},\Q} } -$ induces an equivalence between the category of 
(coherent) left $\D ^\dag _{\fX ^{\prime},\Q} $-modules 
and that of 
(coherent) left $\D ^\dag _{\fX,\Q} $-modules, which is a quasi-inverse equivalence given by  $F ^*$ 
(see \ref{thm236Be2-cordag-cohright}).
\item The functor
$F_+ : = -  \otimes _{\D ^\dag _{\fX ^{\prime},\Q} }  F ^* \D ^\dag _{\fX ^{\prime},\Q}$ induces an equivalence between the category of 
(coherent) right $\D ^\dag _{\fX ^{\prime},\Q} $-modules 
and that of 
(coherent) right $\D ^\dag _{\fX,\Q} $-modules
which is a quasi-inverse equivalence given by  $F ^\flat$ (see \ref{thm236Be2-cordag-cohright}).

\end{enumerate}

\end{coro}

\subsection{Exchanging left and right $\cD$-modules, commutation with Frobenius}

\begin{lem}
\label{lem-f*Fflatcomm}
Let $Y$ be a smooth $T$-scheme of finite type. 
Let $f \colon X \to Y$ be a relatively perfect morphism locally of formal finite type.
We suppose there exists $X'$ a $T$-scheme locally of  formal finite type  having locally finite $p$-bases
(resp. $Y'$ a smooth $T$-scheme of finite type)
such that 
$X ' \times _T T _0 \riso X _0 ^{(s)}$
(resp. $Y ' \times _T T _0 \riso Y _0 ^{(s)}$).
We suppose there exists a lifting
$F _X \colon X \to X'$ 
 of 
$F ^s _{X _0/T _0}$,
and a lifting 
$F _Y \colon Y \to Y'$ 
 of 
$F ^s _{Y _0/T _0}$.
We suppose there exists 
$f '\colon X '\to Y'$ a (relatively perfect) morphism 
which is a lifting of $f _0 ^{(s)}$ 
and is such that 
$f '  \circ F _X= F _Y \circ f$.
Let $\cM '$ be right $\D ^{(m)} _{Y'/T}$-module. 
Then we have the isomorphism of right $\D ^{(m+s)} _{Y/T}$-modules
of the form
\begin{equation}
\label{lem-f*Fflatcomm-iso}
F _X ^{\flat} f ^{\prime *} (\cM ')
\riso 
f ^{*} F _Y ^{\flat}  (\cM '),
\end{equation}
where the structure of 
right $\D ^{(m+s)} _{Y/T}$-modules comes from \ref{lem-rightD-petale} and \ref{241Be2}.

\end{lem}

\begin{proof}
1) We check that the square 
\begin{equation}
\label{lem-f*Fflatcomm-iso-sq1}
\xymatrix{
{X} 
\ar[r] ^-{f}
\ar[d] ^-{F _X}
& 
{Y} 
\ar[d] ^-{F _Y}
\\ 
{X'} 
\ar[r] ^-{f'}
& 
{Y'} }
\end{equation}
is cartesian. Indeed, since $f _0$ is relatively perfect, 
then by definition of the notion of relative perfectness, 
the square \ref{lem-f*Fflatcomm-iso-sq1} is cartesian modulo $\pi$.
Since $f$ and $f'$ are formally étale, 
then we get the cartesianity of \ref{lem-f*Fflatcomm-iso-sq1}.

2) We construct the isomorphism \ref{lem-f*Fflatcomm-iso} as follows. 
Since the diagram \ref{lem-f*Fflatcomm-iso-sq1} is cartesian, 
the functors 
$f ^{\prime *}$ and 
$f ^{*}$ are equal on the category of $\O _{Y}$-modules
(viewing an $\O _{Y}$-module as an $\O _{Y'}$-module via $F _Y$). 
Hence, we get the isomorphism 
$F _X ^{\flat} f ^{\prime *} (\cM ')
=
\cH om _{\cO _{X '}} ( \cO _X, f ^{\prime *} (\cM '))
\riso 
\cH om _{\cO _{X '}} ( f ^{\prime *} (\cO _Y), f ^{\prime *} (\cM '))
\riso
f ^{\prime *} \cH om _{\cO _{Y '}} (  \cO _Y, \cM ')
\riso
f ^{*} F _Y ^{\flat}  (\cM ')$.

3) It remains to check that the isomorphism \ref{lem-f*Fflatcomm-iso}
is horizontal, i.e. commutes with 
$(m+s)$-PD-costratifcations. This is easy and left to the reader.
\end{proof}

\begin{rem}
\label{remlem-f*Fflatcomm}
With notation \ref{lem-f*Fflatcomm}, it follows from
\ref{f0formétale-fforméta0},
that $f$ and $f'$ are flat.
The isomorphism \ref{lem-f*Fflatcomm-iso}
is equal to that of 
\cite[III.6.3]{HaRD}.
\end{rem}

\begin{empt}
\label{ntn-Frob-desc-simple}
We keep notation \ref{subsecElevFrob}.
Suppose we have a relatively perfect morphism of the form $g _0\colon X _0 \to \A ^d _{T _0}$. 
Choose coordinates $t _1,\dots, t _d$ (resp. $t  '_1,\dots, t '_d$)
inducing a lifting 
$g \colon X  \to \A ^d _{T}$ of $g _0$
(resp. $g '\colon X ' \to \A ^d _{T}$ of $g _0 ^{(s)}$).
Following remark \ref{rem-Frob-desc-simple}, 
there exists a unique morphism
$F \colon X \to X'$ which is a lifting of 
$F ^{s} _{X_0/T _0}$ and making commutative the diagram
\begin{equation}
\label{F!-appl1}
\xymatrix@R=0,3cm{
{X} 
\ar@{.>}[d] _-{F} 
\ar[r] ^-{t _1,\dots,t _d} _-{g}
& 
{ \mathbb{A} ^d _{T}} 
\ar[d] ^-{F} 
\ar[rd] ^-{}
\\ 
{X'} 
\ar[r] ^-{t '_1,\dots,t '_d} _-{g'}
& 
{\mathbb{A} ^d _{T}} 
\ar[r] ^-{}
& 
{T .} 
}
\end{equation}

\end{empt}

\begin{lem}
\label{242Be2}
With notation and hypothesis of \ref{ntn-Frob-desc-simple}, 
there exists a canonical isomorphism of right $ \D  _{X/T}$-modules
\begin{equation}
\label{242Be2iso}
\mu _X \colon 
F ^\flat (\omega _{X'/T})
\riso
\omega _{X/T}.
\end{equation}

\end{lem}

\begin{proof}
Following
\cite[2.4.2]{Be2}, we have the canonical isomorphism
$F ^\flat (\omega _{\mathbb{A} ^d _{T}/T})
\riso
\omega _{\mathbb{A} ^d _{T}/T}$
of right $ \D  _{\mathbb{A} ^d _{T}/T}$-modules.
By applying $g ^* $ this yields the isomorphism
$g ^* F ^\flat (\omega _{\mathbb{A} ^d _{T}/T})
\riso
g ^* \omega _{\mathbb{A} ^d _{T}/T}$
of right $ \D  _{X/T}$-modules.
Since the square of \ref{F!-appl1} is cartesian, 
following 
\ref{lem-f*Fflatcomm}, 
we get the canonical isomorphism
$g ^* F ^\flat (\omega _{\mathbb{A} ^d _{T}/T})
\riso 
F ^\flat g ^{\prime *}  (\omega _{\mathbb{A} ^d _{T}/T})$.
We get the canonical isomorphism
$F ^\flat g ^{\prime *}  (\omega _{\mathbb{A} ^d _{T}/T})
\riso
g ^* \omega _{\mathbb{A} ^d _{T}/T}$.
Following \ref{lem-rightD-petale-omega}, we have 
$g ^{\prime *}  (\omega _{\mathbb{A} ^d _{T}/T})
\riso 
\omega _{X'/T}$
adn
$g ^* \omega _{\mathbb{A} ^d _{T}/T}
\riso 
\omega _{X/T}$.
Hence we get the isomorphism
$\mu _X \colon 
F ^\flat (\omega _{X'/T})
\riso
\omega _{X/T}$.
\end{proof}

\begin{prop}
\label{243Be2}
We keep notation and hypothesis of \ref{ntn-Frob-desc-simple}.
For any  left $\D ^{(m)} _{X'/T}$-module $\cE '$,
we have the canonical isomorphism of right $\D ^{(m+s)} _{X/T}$-modules of the form
\begin{equation}
\label{243Be2iso}
\omega _{X/T} \otimes _{\O _X}
F  ^{ *} (\cE ')
\riso 
F ^{\flat}  (\omega _{X'/T} \otimes _{\O _{X'}} \cM ').
\end{equation}

\end{prop}

\begin{proof}
By using \ref{242Be2},
we can copy the proof of \cite[2.4.3]{Be2}.
\end{proof}

Similarly to \cite[2.4.4--5]{Be2},
we get the following corollaries.
\begin{coro}
\label{244Be2}
We keep notation and hypothesis of \ref{ntn-Frob-desc-simple}.
For any  right $\D ^{(m)} _{X'/T}$-module $\cM '$,
we have the canonical isomorphism of right $\D ^{(m+s)} _{X/T}$-modules of the form
\begin{equation}
\label{244Be2iso}
F  ^{ *} (\cM' \otimes _{\O _X} \omega _{X/T} ^{-1})
\riso 
F ^{\flat}  (\cM ') \otimes _{\O _{X'}} \omega _{X'/T} ^{-1}.
\end{equation}

\end{coro}

\subsection{Glueing isomorphisms and Frobenius}

\begin{prop}
\label{215Be2}
Let $f,f ' \colon X \to Y$ be two morphisms of $p$-smooth $T$-schemes inducing the same restriction
$f _0 \colon X _0 \to Y$. 
Suppose the $m$-PD-ideal $\fa$ is $m$-PD-nilpotent. 

\begin{enumerate}[(a)]
\item Let $\E $ be a left $\D ^{(m)} _{Y/T}$-module. Then, we have a canonical isomorphism of left $\D ^{(m)} _{X/T}$-modules of the form
\begin{equation}
\label{215Be2-tauff'}
\tau _{f,f'}
\colon 
f ^* ( \E) 
\riso 
f ^{\prime *} ( \E)
\end{equation}
such that $\tau _{f,f}=id$,
and, for any third morphism $f '' \colon X \to Y$ 
inducing the same restriction
$f _0 \colon X _0 \to Y$,
we have 
$\tau _{f,f''} = \tau _{f,f'} \circ \tau _{f',f''}$.

\item Suppose that $f$ is finite. 
Let $\cM $ be right $\D ^{(m)} _{Y/T}$-module. 
Then, we have a canonical isomorphism of right $\D ^{(m)} _{X/T}$-modules of the form
\begin{equation}
\label{215Be2-tauff'right}
\sigma _{f,f'}
\colon 
f ^\flat ( \E) 
\riso 
f ^{\prime \flat} ( \E)
\end{equation}
such that $\sigma _{f,f}=id$,
and, for any third morphism $f '' \colon X \to Y$ 
inducing the same restriction
$f _0 \colon X _0 \to Y$,
we have 
$\sigma _{f,f''} = \sigma _{f,f'} \circ \sigma _{f',f''}$.

\end{enumerate}

\end{prop}

\begin{proof}
For the first assertion, we can copy word by word the proof of \cite[2.1.5]{Be2}. 
For the second one, we can copy this proof above by replacing the functor $f \mapsto f ^*$
by $f \mapsto f ^\flat$ (and by replacing the use of $m$-PD-stratification by that of $m$-PD-costratification).
\end{proof}

\begin{prop}
\label{225Be2}
We keep notation and hypotheses of \ref{subsecElevFrob}.
Suppose 
there exists a second morphism 
$F ' \colon X \to X'$ (e.g. since $X'/T$ is formally smooth, such a lifting exists when $X$ is affine) making commutative 
the diagram \ref{Fliftingdiag}.

\begin{enumerate}[(a)]
\item Let $\E '$ be a left $\D ^{(m)} _{X'/T}$-module. 
Then the glueing isomorphism
$\tau _{F,F'}
\colon 
F ^* ( \E') 
\riso 
F^{\prime *} ( \E')$
defined in \ref{215Be2-tauff'} is 
$\D ^{(m+s)} _{X/T}$-linear.

\item Let $\cM '$ be a right $\D ^{(m)} _{X'/T}$-module. 
Then the glueing isomorphism
$\sigma _{F,F'}
\colon 
F ^\flat ( \cM') 
\riso 
F^{\prime \flat} ( \cM ')$
defined in \ref{215Be2-tauff'right} is 
$\D ^{(m+s)} _{X/T}$-linear.
\end{enumerate}
\end{prop}

\begin{proof}
We can copy the proof of \cite[2.2.5]{Be2}.
\end{proof}

\section{Derived categories of inductive systems of arithmetic $\D$-modules}
\label{ntn-tildeD(Z)}
Let $\fP$ be a formal  $\fS$-scheme locally of formal finite type
and having locally finite $p$-bases  over $\fS$.
Let 
$T$ be a divisor of $P$.
Divisors of $P$ will be supposed to be reduced divisors (in our context, this is not really less general). 
Remark that since $P$ is regular (see \ref{regularity/formalsm}), then Weil divisors correspond to Cartiel divisors. Hence, in our context, 
a divisor is determined by its irreducible components. 
To reduce the amount of notation, we set 
$\smash{\widehat{\D}} _{\fP /\fS } ^{(m)} (T):=
\widehat{\B} ^{(m)} _{\fP } ( T)  \smash{\widehat{\otimes}} _{\O _{\fP}} \smash{\widehat{\D}} _{\fP /\fS } ^{(m)}$, 
where $\widehat{\B} ^{(m)} _{\fP} ( T) $ is the sheaf constructed in
\ref{ntnODdagZ}
and
$\smash{\D} _{\fP /\fS } ^{(m)}$ is the sheaf of differential operators of level $m$ over $\fP /\fS $
(see \ref{2.2Be1}). 
We fix  $\lambda _0\colon \N \to \N$ an increasing map such that 
$\lambda _{0} (m) \geq m$ for any $m \in \N$. 
We set 
$\widetilde{\B} ^{(m)} _{\fP} ( T):= \widehat{\B} ^{(\lambda _0 (m))} _{\fP} ( T)$ 
et
$\smash{\widetilde{\D}} _{\fP /\fS } ^{(m)} (T):=
\widetilde{\B} ^{(m)} _{\fP} ( T)  \smash{\widehat{\otimes}} _{\O _{\fP}} \smash{\widehat{\D}} _{\fP /\fS } ^{(m)}$.
Finally, we set 
$\smash{\D} _{P  _i /S  _i} ^{(m)} (T):= \V / \pi ^{i+1} \otimes _{\V} \smash{\widehat{\D}} _{\fP /\fS  } ^{(m)} (T) 
=
\B ^{(m)} _{P _i} ( T)  \otimes _{\O _{P _i}} \smash{\D} _{P  _i/S  _i} ^{(m)}$
and
$\smash{\widetilde{\D}} _{P  _i/S  _i} ^{(m)} (T):=\widetilde{\B} ^{(m)} _{P _i} ( T)  \otimes _{\O _{P _i}} \smash{\D} _{P  _i/S  _i} ^{(m)}$.
When $T$ is empty, we get

\subsection{Localisation of derived categories of inductive systems of arithmetic $\D$-modules}

\begin{empt}
[Berthelot's localized categories of the form $\smash{\underrightarrow{LD}} _{\Q}$]
\label{loc-LM}
We recall below some constructions of Berthelot of \cite[4.2.1 and 4.2.2]{Beintro2} 
which are still valid in our context of formal schemes locally of formal finite type having locally finite $p$-bases
and also by adding singularities along a divisor. 
We have the inductive system of rings 
$\smash{\widetilde{\D}} _{\fP /\fS } ^{(\bullet)}(T) 
: =
(\smash{\widetilde{\D}} _{\fP /\fS } ^{(m)}(T) )_{m\in \N}$. 
We get the derived categories
$D ^{\sharp} ( \smash{\widetilde{\D}} _{\fP /\fS } ^{(\bullet)}(T))$,  
where $\sharp \in \{\emptyset, +,-, \mathrm{b}\}$.
The objects of $D ^{\sharp}( \smash{\widetilde{\D}} _{\fP /\fS } ^{(\bullet)}(T))$
are denoted by
 $\E ^{(\bullet)}= (\E ^{(m)} , \alpha ^{(m',m)})$, 
where $m,m'$ run over non negative integers such that  $m' \geq m$,
where $\E ^{(m)} $ is a complex of $\smash{\widetilde{\D}} _{\fP /\fS } ^{(m)}(T)$-modules
and $\alpha ^{(m',m)} \colon \E ^{(m)}\to \E ^{(m')}$ are $\smash{\widetilde{\D}} _{\fP /\fS } ^{(m)}(T)$-linear morphisms.

\begin{enumerate}[(a)]
\item Let $M$ bet the filtrant set (endowed with the canonical order) 
of increasing maps $\chi \colon \N \to \N$. 
For any map $\chi \in M$, we set
$\chi ^{*} (\E ^{(\bullet)}) := (\E ^{(m)} , p ^{\chi (m') -\chi (m)}\alpha ^{(m',m)})$.
We obtain the functor
$\chi ^{*} \colon D ( \smash{\widetilde{\D}} _{\fP /\fS } ^{(\bullet)}(T))\to D ( \smash{\widetilde{\D}} _{\fP /\fS } ^{(\bullet)}(T))$ 
as follows:
if $f ^{(\bullet)} \colon \E ^{(\bullet)} \to \FF ^{(\bullet)}$ is a morphism of $D ( \smash{\widetilde{\D}} _{\fP /\fS } ^{(\bullet)}(T))$, 
then the morphism of level $m$ of  $\chi ^{*} f ^{(\bullet)} 
\colon 
\chi ^{*} (\E ^{(\bullet)}) 
\to 
\chi ^{*} (\FF ^{(\bullet)})$ 
is $f ^{(m)}$.
If $\chi _1, \chi _2 \in M$, we compute
$\chi _1 ^* \circ \chi _2 ^* = (\chi _1 +\chi _2)^*$, and in particular $\chi _1 ^*$ and $\chi _2 ^*$ commute.
Moreover, if  $\chi _1 \leq \chi _2$, then we get the morphism 
$\chi _1 ^* ( \E ^{(\bullet)})\to \chi _2 ^* ( \E ^{(\bullet)})$ defined at the level $m$ by  
$p ^{\chi _2 (m) -\chi _1(m)}\colon \E ^{(m)} \to \E ^{(m)}$. 
A morphism $f ^{(\bullet)} \colon \E ^{(\bullet)} \to \FF ^{(\bullet)}$ of $D ( \smash{\widetilde{\D}} _{\fP /\fS } ^{(\bullet)}(T))$
is an {``ind-isogeny''} if there exist  $\chi \in M$ 
and a morphism
$g ^{(\bullet)} \colon \FF ^{(\bullet)} \to \chi ^{*} \E ^{(\bullet)}$ of $D ( \smash{\widetilde{\D}} _{\fP /\fS } ^{(\bullet)}(T))$
such that 
$g ^{(\bullet)}\circ f ^{(\bullet)}$ and $\chi ^{*} (f ^{(\bullet)}) \circ g ^{(\bullet)}$ 
are the canonical morphisms described above (in the case $\chi _1 =0$ and $\chi _2=\chi$).
The subset of ind-isogenies is a multiplicative system (this follows from Proposition
\cite[I.4.2]{HaRD} and the analogue of Lemma \cite[1.1.2]{caro-stab-sys-ind-surcoh} still valid without the hypothesis that $k$ is perfect). 
The localisation of $D ^{\sharp} ( \smash{\widetilde{\D}} _{\fP /\fS } ^{(\bullet)}(T))$
with respect to ind-isogenies is denoted by
$\smash{\underrightarrow{D}} ^{\sharp} _{\Q} ( \smash{\widetilde{\D}} _{\fP /\fS } ^{(\bullet)}(T))$.

\item Let $L$ be the filtrant set of increasing maps 
$\lambda \colon \N \to \N$ such that $\lambda (m ) \geq m$.
For any $\lambda \in L$, we put
$\lambda ^{*} (\E ^{(\bullet)}) := (\E ^{(\lambda(m))} , \alpha ^{(\lambda(m'),\lambda(m))})_{m'\geq m}$.
When $\lambda _1, \lambda _2 \in L$, we compute
$\lambda _1 ^* \circ \lambda _2 ^* = (\lambda _1 \circ\lambda _2)^*$.
When $\lambda _1 \leq \lambda _2$, we have the canonical morphism 
$\lambda _1 ^* (\E ^{(\bullet)}) \to 
\lambda _2 ^* (\E ^{(\bullet)})$ defined at the level $m$ by the morphism
$\alpha ^{(\lambda _2(m),\lambda _1(m))} \colon 
\E ^{(\lambda _1(m))} \to \E ^{(\lambda _2(m))}$. 
Similarly to \cite[4.2.2]{Beintro2}, 
we denote by $\Lambda ^{\sharp}$ the set of morphisms
$f ^{(\bullet)} \colon \E ^{(\bullet)} \to \FF ^{(\bullet)}$
of 
$\smash{\underrightarrow{D}} ^{\sharp} _{\Q}  ( \smash{\widetilde{\D}} _{\fP /\fS } ^{(\bullet)}(T))$
such that there exist $\lambda \in L$ and a morphism
$g ^{(\bullet)} \colon \FF ^{(\bullet)} \to \lambda ^{*} \E ^{(\bullet)}$ of $\smash{\underrightarrow{D}} _{\Q} ( \smash{\widetilde{\D}} _{\fP /\fS } ^{(\bullet)}(T))$
such that 
the morphisms 
$g ^{(\bullet)}\circ f ^{(\bullet)}$ 
and $\lambda ^{*} (f ^{(\bullet)}) \circ g ^{(\bullet)}$ 
of $\smash{\underrightarrow{D}} ^{\sharp} _{\Q} ( \smash{\widetilde{\D}} _{\fP /\fS } ^{(\bullet)}(T))$
are the canonical morphisms (i.e. we take $\lambda _1 =id$ and $\lambda _2=\lambda$).
The morphisms belonging to $\Lambda $ are called {``lim-isomorphisms''}.
We check that  $\Lambda ^{\sharp} $ is a multiplicative system
(again, use  \cite[I.4.2]{HaRD} and the analogue of Lemma \cite[1.1.2]{caro-stab-sys-ind-surcoh}). 
By localizing 
$\smash{\underrightarrow{D}}  ^{\sharp}_{\Q}
( \smash{\widetilde{\D}} _{\fP /\fS } ^{(\bullet)}(T))$
with respect to  lim-isomorphisms 
we get a category denoted by
$\smash{\underrightarrow{LD}} ^{\sharp} _{\Q}
( \smash{\widetilde{\D}} _{\fP /\fS } ^{(\bullet)}(T))$.

\item Let  $\chi _1\leq \chi _2 $ in $M$ and $\lambda _1 \leq \lambda _2$ in $L$.
We get by composition the canonical morphism
$\lambda  _1^{*} \chi  _1^{*} \to \lambda _2^{*} \chi _2^{*}$.
By considering  $\chi _1 \circ \lambda _1$ as an element of $M$, we get the equality
$\lambda _1 ^* \chi _1 ^* = (\chi _1 \circ \lambda _1) ^* \lambda _1 ^*$.
Let  $S ^{\sharp}$ be the set of morphisms 
$f ^{(\bullet)} \colon \E ^{(\bullet)} \to \FF ^{(\bullet)}$
of 
$D ^{\sharp} ( \smash{\widetilde{\D}} _{\fP /\fS } ^{(\bullet)}(T))$
such that there exist $\chi \in M$, 
$\lambda \in L$ and a morphism
$g ^{(\bullet)} \colon \FF ^{(\bullet)} \to \lambda ^{*} \chi ^{*}\E ^{(\bullet)}$ of $D ( \smash{\widetilde{\D}} _{\fP /\fS } ^{(\bullet)}(T))$
such that 
$g ^{(\bullet)}\circ f ^{(\bullet)}$ and $\lambda ^{*} \chi ^{*}(f ^{(\bullet)}) \circ g ^{(\bullet)}$ 
are the canonical morphisms. 
The elements of   $S ^{\sharp}$ are called {`` lim-ind-isogenies''}.
We check as usual that $S ^{\sharp}$ is a multiplicative system.

\end{enumerate}

\end{empt}

\begin{empt}
\label{S=LDQ}
Similarly to \cite[1.1.5]{caro-stab-sys-ind-surcoh},
we check the canonical equivalence of categories
$S ^{\sharp -1} D ^{\sharp}
( \smash{\widetilde{\D}} _{\fP /\fS } ^{(\bullet)}(T))
\cong
\smash{\underrightarrow{LD}} ^{\sharp} _{\Q} 
( \smash{\widetilde{\D}} _{\fP /\fS } ^{(\bullet)}(T))$,
which is the identity over the objects.

\end{empt}

\begin{empt}
\label{HomLDQ}
Similarly to \cite[1.1.6]{caro-stab-sys-ind-surcoh},
for any  $\E ^{(\bullet)}, \FF ^{(\bullet)} \in \underrightarrow{LD} ^{\sharp} _{\Q} (\smash{\widetilde{\D}} _{\fP /\fS } ^{(\bullet)} (T))$,
we have the equality
\begin{equation}
\label{4.2.2Beintro}
\mathrm{Hom} _{\underrightarrow{LD} ^{\sharp} _{\Q} (\smash{\widetilde{\D}} _{\fP /\fS } ^{(\bullet)} (T))}
(\E ^{(\bullet)}, \FF ^{(\bullet)} )
=
\underset{\lambda \in L}{\underrightarrow{\lim}}~
\underset{\chi \in M}{\underrightarrow{\lim}}~
\mathrm{Hom} _{D ^{\sharp} (\smash{\widetilde{\D}} _{\fP /\fS } ^{(\bullet)} (T))}
(\E ^{(\bullet)}, \lambda ^{*} \chi ^{*}\FF ^{(\bullet)} ).
\end{equation}

\end{empt}

\begin{empt}
\label{defi-M-L}
We denote by $M (\smash{\widetilde{\D}} _{\fP /\fS } ^{(\bullet)} (T))$ the abelian category of
$\smash{\widetilde{\D}} _{\fP /\fS } ^{(\bullet)} (T)$-modules.
The $\smash{\widetilde{\D}} _{\fP /\fS } ^{(\bullet)} (T)$-modules
are denoted by 
$\E ^{(\bullet)}= (\E ^{(m)} , \alpha ^{(m',m)})$, 
where $m,m'$ run through non negative integers  $m' \geq m$,
where $\E ^{(m)} $ is a $\smash{\widetilde{\D}} _{\fP /\fS } ^{(m)}(T)$-module
and $\alpha ^{(m',m)} \colon \E ^{(m)}\to \E ^{(m')}$ are
$\smash{\widetilde{\D}} _{\fP /\fS } ^{(m)}(T)$-linear morphisms.
For any  $\chi \in M$, we denote similarly to \ref{loc-LM} the object
$\chi ^{*} (\E ^{(\bullet)}) := (\E ^{(m)} , p ^{\chi (m') -\chi (m)}\alpha ^{(m',m)})$.
In fact, we get the functor
$\chi ^{*}\colon M (\smash{\widetilde{\D}} _{\fP /\fS } ^{(\bullet)} (T)) \to M (\smash{\widetilde{\D}} _{\fP /\fS } ^{(\bullet)} (T))$. 
Moreover, similarly to \ref{loc-LM},
for any $\lambda \in L$, we set 
$\lambda ^{*} (\E ^{(\bullet)}) := (\E ^{(\lambda(m))} , \alpha ^{(\lambda(m'),\lambda(m))})$.

Similarly to 
\ref{loc-LM},
we can define the saturated multiplicative system of ``lim-ind-isogenies'' of 
$M( \smash{\widetilde{\D}} _{\fP /\fS } ^{(\bullet)}(T))$
and we get the corresponding localized category
$\smash{\underrightarrow{LM}} _{\Q}
( \smash{\widetilde{\D}} _{\fP /\fS } ^{(\bullet)}(T))$.
\end{empt}

\begin{empt}
The results of \cite[1.2.1]{caro-stab-sys-ind-surcoh} are still valid in our context:
we check the canonical equivalence of categories
$S ^{-1}M (\smash{\widetilde{\D}} _{\fP /\fS } ^{(\bullet)} (T))
\cong 
\underrightarrow{LM} _{\Q} (\smash{\widetilde{\D}} _{\fP /\fS } ^{(\bullet)} (T))$.
Moreover, for any $\E ^{(\bullet)}, \FF ^{(\bullet)} \in \underrightarrow{LM} _{\Q} (\smash{\widetilde{\D}} _{\fP /\fS } ^{(\bullet)} (T))$ 
we have 
\begin{equation}
\label{4.2.2BeintroLMQ}
\mathrm{Hom} _{\underrightarrow{LM} _{\Q} (\smash{\widetilde{\D}} _{\fP /\fS } ^{(\bullet)} (T))}
(\E ^{(\bullet)}, \FF ^{(\bullet)} )
=
\underset{\lambda \in L}{\underrightarrow{\lim}}\;
\underset{\chi \in M}{\underrightarrow{\lim}}\;
\mathrm{Hom} _{M (\smash{\widetilde{\D}} _{\fP /\fS } ^{(\bullet)} (T))}
(\E ^{(\bullet)}, \lambda ^{*} \chi ^{*}\FF ^{(\bullet)} ).
\end{equation}

The category $\underrightarrow{LM} _{\Q} (\smash{\widetilde{\D}} _{\fP /\fS } ^{(\bullet)} (T))$ is abelian 
and the multiplicative system of  lim-ind-isogenies of
$M (\smash{\widetilde{\D}} _{\fP /\fS } ^{(\bullet)} (T))$ is saturated
(we can copy the proof of \cite[1.2.4]{caro-stab-sys-ind-surcoh}).

\end{empt}

\begin{empt}
We denote by  $M ( \smash{\D} ^\dag _{\fP } (\hdag T) _{\Q} )$ 
the abelian category of 
$\smash{\D} ^\dag _{\fP } (\hdag T) _{\Q}$-modules.
By tensoring by $\Q$ and next by applying the inductive limit on the level, we get the functor
$\underrightarrow{\lim} 
\colon
M
(\smash{\widetilde{\D}} _{\fP /\fS } ^{(\bullet)}(T))
\to
M ( \smash{\D} ^\dag _{\fP } (\hdag T) _{\Q} )$.
Since this functor sends a lim-ind-isomorphism to an  isomorphism, it factorizes canonically through the functor 
\begin{equation}
\label{M-eq-lim}
\underrightarrow{\lim} 
\colon
\smash{\underrightarrow{LM}}  _{\Q}
(\smash{\widetilde{\D}} _{\fP /\fS } ^{(\bullet)}(T))
\to
M ( \smash{\D} ^\dag _{\fP } (\hdag T) _{\Q} ).
\end{equation}
Similarly, we get 
\begin{equation}
\label{D-eq-lim}
\underrightarrow{\lim} 
\colon
\smash{\underrightarrow{LD}}  ^\mathrm{b} _{\Q}
(\smash{\widetilde{\D}} _{\fP /\fS } ^{(\bullet)}(T))
\to
D ^\mathrm{b} ( \smash{\D} ^\dag _{\fP } (\hdag T) _{\Q} ).
\end{equation}
\end{empt}

\begin{prop}
\label{eqcatLD=DSM}
The canonical functor 
$D ^{\mathrm{b}} (\smash{\widetilde{\D}} _{\fP /\fS } ^{(\bullet)} (T))
\to
D ^{\mathrm{b}}
(\underrightarrow{LM} _{\Q} (\smash{\widetilde{\D}} _{\fP /\fS } ^{(\bullet)} (T)))$
of triangulated categories
induced by the functor of abelian categories
$M(\smash{\widetilde{\D}} _{\fP /\fS } ^{(\bullet)} (T))
\to
\underrightarrow{LM} _{\Q} (\smash{\widetilde{\D}} _{\fP /\fS } ^{(\bullet)} (T))$
factorizes canonically through the equivalence of triangulated categories
\begin{equation}
\label{eqcatLD=DSM-fonct}
\underrightarrow{LD} ^{\mathrm{b}} _{\Q} (\smash{\widetilde{\D}} _{\fP /\fS } ^{(\bullet)} (T))
\cong 
D ^{\mathrm{b}}
(\underrightarrow{LM} _{\Q} (\smash{\widetilde{\D}} _{\fP /\fS } ^{(\bullet)} (T))).
\end{equation}

\end{prop}

\begin{proof}
We can copy the proof of \cite[1.2.11]{caro-stab-sys-ind-surcoh}.
\end{proof}

\begin{empt}
\label{empt-diag-Hn-comp}
The equivalence \ref{eqcatLD=DSM-fonct} commutes with cohomological functors, 
i.e. we have for any $n\in \N$ the commutative diagram
\begin{equation}
\label{diag-Hn-comp}
\xymatrix{
{D ^{\mathrm{b}} (\smash{\widetilde{\D}} _{\fP /\fS } ^{(\bullet)} (T))} 
\ar[r] ^-{}
\ar[d] ^-{\H ^n}
&
{\underrightarrow{LD} ^{\mathrm{b}} _{\Q} (\smash{\widetilde{\D}} _{\fP /\fS } ^{(\bullet)} (T))} 
\ar[r] ^-{\cong}
\ar@{.>}[d] ^-{\H ^n}
& 
{D ^{\mathrm{b}} (\underrightarrow{LM}  _{\Q} (\smash{\widetilde{\D}} _{\fP /\fS } ^{(\bullet)} (T))) } 
\ar[d] ^-{\H ^n}
\\ 
{M(\smash{\widetilde{\D}} _{\fP /\fS } ^{(\bullet)} (T))} 
\ar[r] ^-{}
& 
{\underrightarrow{LM}  _{\Q} (\smash{\widetilde{\D}} _{\fP /\fS } ^{(\bullet)} (T))} 
\ar@{=}[r] ^-{}
& 
{\underrightarrow{LM}  _{\Q} (\smash{\widetilde{\D}} _{\fP /\fS } ^{(\bullet)} (T))} 
}
\end{equation}
where the middle vertical arrow is the one making 
commutative by definition the left square (see \cite[1.2.6]{caro-stab-sys-ind-surcoh}).
\end{empt}

\subsection{Coherence}

Similarly to \cite[2.2.1]{caro-stab-sys-ind-surcoh}, we have the following definition. 
\begin{dfn}
[Coherence up to lim-ind-isogeny]
\label{coh-loc-pres-fini-lim-ind-iso}
Let $\E ^{(\bullet)}$ be a
$\smash{\widetilde{\D}} _{\fP /\fS } ^{(\bullet)} (T)$-module.
The module 
$\E ^{(\bullet)}$
is said to be a
$\smash{\widetilde{\D}} _{\fP /\fS } ^{(\bullet)} (T)$-module 
of finite type up to  lim-ind-isogeny
if there exists an open covering 
$(\fP _i ) _{i\in I}$ of  $\fP$ 
such that, for any $i \in I$, 
there exists an exact sequence of 
$\underrightarrow{LM} _{\Q}  (\smash{\widetilde{\D}} _{\fP /\fS } ^{(\bullet)} (T))$ of the form: 
$\left ( \smash{\widetilde{\D}} _{\fP  _i} ^{(\bullet)} (T  \cap P _i) \right) ^{r _i}
\to 
\E ^{(\bullet)} | \fP _i
\to 0$,
where $r _i\in \N$.
Similarly, we get the notion of
$\smash{\widetilde{\D}} _{\fP /\fS } ^{(\bullet)} (T)$-module 
locally of finite presentation up to  lim-ind-isogeny 
(resp. coherence up lim-ind-isogeny). 
\end{dfn}

\begin{ntn}
\label{nota-(L)Mcoh}
We denote by 
 $\underrightarrow{LM} _{\Q, \mathrm{coh}} (\smash{\widetilde{\D}} _{\fP /\fS } ^{(\bullet)} (T))$
the full subcategory of 
$\underrightarrow{LM} _{\Q} (\smash{\widetilde{\D}} _{\fP /\fS } ^{(\bullet)} (T))$
consisting of 
coherent $\smash{\widetilde{\D}} _{\fP /\fS } ^{(\bullet)} (T)$-modules up to lim-ind-isogeny.
\end{ntn}

\begin{prop}
\label{LQ-coh-stab}
The full subcategory 
$\underrightarrow{LM} _{\Q, \mathrm{coh}} ( \smash{\widetilde{\D}} _{\fP /\fS } ^{(\bullet)} (T))$
of 
$\underrightarrow{LM} _{\Q}  (\smash{\widetilde{\D}} _{\fP /\fS } ^{(\bullet)} (T))$
is stable by isomorphisms, 
kernels, cokernels, extensions.
\end{prop}

\begin{proof}
We can copy the proof of \cite[2.2.8]{caro-stab-sys-ind-surcoh}.
\end{proof}

\begin{ntn}
\label{ntn-cohDLM}
For any $\sharp  \in \{ 0,+,-, \mathrm{b}, \emptyset\}$, 
we denote by  
$D  ^\sharp _{\mathrm{coh}}
(\underrightarrow{LM} _{\Q} (\smash{\widetilde{\D}} _{\fP /\fS } ^{(\bullet)} (T)))$
the full subcategory of 
$D ^\sharp
(\underrightarrow{LM} _{\Q} (\smash{\widetilde{\D}} _{\fP /\fS } ^{(\bullet)} (T)))$
consisting of complexes $\E ^{(\bullet)}$ such that, for any $n \in \Z$, 
$\mathcal{H} ^{n} (\E ^{(\bullet)}) \in 
\underrightarrow{LM} _{\Q, \mathrm{coh}} (\smash{\widetilde{\D}} _{\fP /\fS } ^{(\bullet)} (T))$
(see notation  \ref{nota-(L)Mcoh}). 
These objects are called coherent complexes of 
$D ^\sharp
(\underrightarrow{LM} _{\Q} (\smash{\widetilde{\D}} _{\fP /\fS } ^{(\bullet)} (T)))$.
\end{ntn}

\begin{empt}
\label{cohDLMislocal}
By definition, the property that an object of 
$\underrightarrow{LM} _{\Q}  (\smash{\widetilde{\D}} _{\fP /\fS } ^{(\bullet)} (T))$
is an object of 
$\underrightarrow{LM} _{\Q, \mathrm{coh}} ( \smash{\widetilde{\D}} _{\fP /\fS } ^{(\bullet)} (T))$
is local in $\fP$. 
This yields that the notion of coherence of \ref{ntn-cohDLM} is local in  $\fP$, i.e. 
the fact that a complex $\E ^{(\bullet)}$ of 
$D ^\sharp
(\underrightarrow{LM} _{\Q} (\smash{\widetilde{\D}} _{\fP /\fS } ^{(\bullet)} (T)))$
is coherent is local.
\end{empt}

\begin{dfn}
[Coherence in the sense of Berthelot]
\label{defi-LDQ0coh}
Let  $\sharp  \in \{\emptyset, +,-, \mathrm{b}\}$.
Let
$\E ^{(\bullet)} \in \smash{\underrightarrow{LD}} ^{\sharp} _{\Q}  ( \smash{\widetilde{\D}} _{\fP /\fS } ^{(\bullet)}(T))$.
The complex
$\E ^{(\bullet)}$
is said to be coherent if 
there exist
$\lambda \in L$
and
$\FF ^{(\bullet)}\in \smash{\underrightarrow{LD}} _{\Q} ^{\sharp} (\lambda ^{*} \smash{\widetilde{\D}} _{\fP /\fS } ^{(\bullet)}(T))$ 
together with an  isomorphism in 
$ \smash{\underrightarrow{LD}} _{\Q} ^{\sharp}  ( \smash{\widetilde{\D}} _{\fP /\fS } ^{(\bullet)}(T))$ of the form
$\E ^{(\bullet)} \riso \FF ^{(\bullet) }$,
such that $\FF ^{(\bullet)}$ satisfies the following conditions:
\begin{enumerate}[(a)]
\item For any $m \in \N$, $\FF ^{(m)} \in D  _{\mathrm{coh}} ^{\sharp} (\smash{\widetilde{\D}} _{\fP /\fS } ^{(\lambda (m))} (T))$ ; 
\item For any  $0\leq m \leq m'$,  
the canonical morphism
\begin{equation}
\label{Beintro-4.2.3M}
\smash{\widetilde{\D}} _{\fP /\fS } ^{(\lambda (m'))} (T) \otimes ^\L _{\smash{\widetilde{\D}} _{\fP /\fS } ^{(\lambda (m))} (T)}
\FF ^{(m)} \to \FF ^{(m')} 
\end{equation}
is an isomorphism.
\end{enumerate}
\end{dfn}

\begin{ntn}
\label{nota-LDQ0coh}
Let $\sharp \in \{\emptyset, +,-, \mathrm{b}\}$.
We denote by  $\underrightarrow{LD}  ^{\sharp} _{\Q, \mathrm{coh}} (\smash{\widetilde{\D}} _{\fP /\fS } ^{(\bullet)} (T))$
the strictly full subcategory of 
$\underrightarrow{LD}  ^{\sharp} _{\Q} (\smash{\widetilde{\D}} _{\fP /\fS } ^{(\bullet)} (T))$
consisting of coherent complexes. 
\end{ntn}

\begin{prop}
\label{eqcat-limcoh}
\begin{enumerate}[(a)]
\item The functor  \ref{M-eq-lim} induces the equivalence of categories
\begin{equation}
\label{M-eq-coh-lim}
\underrightarrow{\lim} 
\colon
\smash{\underrightarrow{LM}}  _{\Q, \mathrm{coh}}
(\smash{\widetilde{\D}} _{\fP /\fS } ^{(\bullet)}(T))
\cong
\mathrm{Coh} ( \smash{\D} ^\dag _{\fP } (\hdag T) _{\Q} ),
\end{equation}
where $\mathrm{Coh} ( \smash{\D} ^\dag _{\fP } (\hdag T) _{\Q} )$ 
is the category of coherent (left)
$\smash{\D} ^\dag _{\fP } (\hdag T) _{\Q}$-modules.

\item The functor 
\ref{D-eq-lim} induces the equivalence of triangulated categories
\begin{equation}
\label{eqcatcoh}
\underrightarrow{\lim} 
\colon 
D ^{\mathrm{b}} _{\mathrm{coh}} (\underrightarrow{LM} _{\Q} (\smash{\widetilde{\D}} _{\fP /\fS } ^{(\bullet)} (T)))
\cong
D ^{\mathrm{b}} _{\mathrm{coh}}( \smash{\D} ^\dag _{\fP } (\hdag T) _{\Q} ).
\end{equation}

\item The equivalence of triangulated categories 
$\underrightarrow{LD} ^{\mathrm{b}} _{\Q} (\smash{\widetilde{\D}} _{\fP /\fS } ^{(\bullet)} (T))
\cong 
D ^{\mathrm{b}}
(\underrightarrow{LM} _{\Q} (\smash{\widetilde{\D}} _{\fP /\fS } ^{(\bullet)} (T)))$
of 
\ref{eqcatLD=DSM-fonct}
induces the equivalence of triangulated categories
\begin{equation}
\label{eqcatLD=DSM-fonct-coh}
\underrightarrow{LD} ^{\mathrm{b}} _{\Q, \mathrm{coh}} (\smash{\widetilde{\D}} _{\fP /\fS } ^{(\bullet)} (T))
\cong
D ^{\mathrm{b}} _{\mathrm{coh}}
(\underrightarrow{LM} _{\Q} (\smash{\widetilde{\D}} _{\fP /\fS } ^{(\bullet)} (T))).
\end{equation}
\end{enumerate}
\end{prop}

\begin{proof}
We can copy the proof of Theorems \cite[2.4.4, 2.5.7]{caro-stab-sys-ind-surcoh}.
\end{proof}

\begin{empt}
\label{thick-subcat}
\begin{enumerate}[(a)]
\item 
\label{thick-subcat1}
Using \ref{LQ-coh-stab}, 
we get that 
$D ^{\mathrm{b}} _{\mathrm{coh}}
(\underrightarrow{LM} _{\Q} (\smash{\widetilde{\D}} _{\fP /\fS } ^{(\bullet)} (T)))$ 
is a thick triangulated subcategory (some authors say saturated or épaisse) of 
$D ^{\mathrm{b}} 
(\underrightarrow{LM} _{\Q} (\smash{\widetilde{\D}} _{\fP /\fS } ^{(\bullet)} (T)))$,
i.e. is a strict triangulated subcategory closed under direct summands. 
Hence, using \ref{eqcatLD=DSM-fonct} and 
\ref{eqcatLD=DSM-fonct-coh}, we get that 
$\underrightarrow{LD} ^{\mathrm{b}} _{\Q, \mathrm{coh}} (\smash{\widetilde{\D}} _{\fP /\fS } ^{(\bullet)} (T))$
is a thick triangulated subcategory of 
$\underrightarrow{LD} ^{\mathrm{b}} _{\Q} (\smash{\widetilde{\D}} _{\fP /\fS } ^{(\bullet)} (T))$. 
\item 
\label{thick-subcat2}
Using the same arguments, 
it follows from \ref{cohDLMislocal} the following local property : the fact that 
a complex of 
$\underrightarrow{LD} ^{\mathrm{b}} _{\Q} (\smash{\widetilde{\D}} _{\fP /\fS } ^{(\bullet)} (T))$
is a coherent complex (i.e. 
a complex of 
$\underrightarrow{LD} ^{\mathrm{b}} _{\Q,\mathrm{coh}} (\smash{\widetilde{\D}} _{\fP /\fS } ^{(\bullet)} (T))$)
is local in $\fP$. 

\end{enumerate}
\end{empt}

\begin{empt}
For any $n \in \N$, 
the cohomological functor 
$\H ^n \colon 
{\underrightarrow{LD} ^{\mathrm{b}} _{\Q} (\smash{\widetilde{\D}} _{\fP /\fS } ^{(\bullet)} (T))}
\to 
{\underrightarrow{LM}  _{\Q} (\smash{\widetilde{\D}} _{\fP /\fS } ^{(\bullet)} (T))} $
of \ref{diag-Hn-comp}
induces 
$\H ^n \colon 
{\underrightarrow{LD} ^{\mathrm{b}} _{\Q,\mathrm{coh}} (\smash{\widetilde{\D}} _{\fP /\fS } ^{(\bullet)} (T))}
\to 
{\underrightarrow{LM}  _{\Q,\mathrm{coh}} (\smash{\widetilde{\D}} _{\fP /\fS } ^{(\bullet)} (T))} $
and we have the commutative diagram (up to canonical isomorphism)
\begin{equation}
\label{diag-Hn-comp-coh}
\xymatrix{
{D ^{\mathrm{b}} (\underrightarrow{LM}  _{\Q,\mathrm{coh}} (\smash{\widetilde{\D}} _{\fP /\fS } ^{(\bullet)} (T))) } 
\ar[d] ^-{\H ^n}
\ar[r] ^-{}
& 
{D ^{\mathrm{b}} _{\mathrm{coh}}(\underrightarrow{LM}  _{\Q} (\smash{\widetilde{\D}} _{\fP /\fS } ^{(\bullet)} (T))) } 
\ar[d] ^-{\H ^n}
& 
{\underrightarrow{LD} ^{\mathrm{b}} _{\Q,\mathrm{coh}} (\smash{\widetilde{\D}} _{\fP /\fS } ^{(\bullet)} (T))} 
\ar[d] ^-{\H ^n}
\ar[l] ^-{\cong}
\\ 
{\underrightarrow{LM}  _{\Q,\mathrm{coh}} (\smash{\widetilde{\D}} _{\fP /\fS } ^{(\bullet)} (T))}
\ar@{=}[r] ^-{}
& 
{\underrightarrow{LM}  _{\Q,\mathrm{coh}} (\smash{\widetilde{\D}} _{\fP /\fS } ^{(\bullet)} (T))}
\ar@{=}[r] ^-{}
& 
{\underrightarrow{LM}  _{\Q,\mathrm{coh}} (\smash{\widetilde{\D}} _{\fP /\fS } ^{(\bullet)} (T)).}
}
\end{equation}
Indeed, the commutativity of the left square is obvious 
and that of the right one is almost tautological (see the commutative diagram \ref{diag-Hn-comp}).
\end{empt}

\subsection{Indcoherence}
We denote by 
$\D ^{(m)}$ either 
$\widetilde{\D} ^{(m)}_{\fP /\fS } (T)$ 
or 
$\widetilde{\D} ^{(m)}_{\fP /\fS } (T) _\Q$.
We denote by 
$\D $ either 
$\widetilde{\D} ^{(m)}_{\fP /\fS } (T)$ 
or 
$\widetilde{\D} ^{(m)}_{\fP /\fS } (T) _\Q$
or 
$\D ^\dag _{\fP /\fS } (\hdag T) _\Q$.
We put 
$D^{(m)} := \Gamma (\fP, \D ^{(m)})$,
$D:= \Gamma (\fP, \D)$.
By copying word by word their proofs,
we get an analogue of the section 
\cite[1.4]{caro-6operations}
in our context. 
For the reader, let us collect them below.

\begin{empt}
\label{IndCoh}
We denote by $\Mod ( \D)$ (resp. $\coh ( \D)$)
the abelian category of left $\D$-modules (resp. coherent left $\D$-modules).
We denote by $\iota \colon \coh ( \D) \to \Mod ( \D)$
the canonical fully faithful functor. 
Since 
$\Mod ( \D)$ admits small filtrant inductive limits,
from \cite[6.3.2]{Kashiwara-schapira-book}
we get a functor denoted by 
 $J \iota  \colon \ind ( \coh ( \D) )\to \Mod ( \D)$
 such that $J \iota$ commutes with small filtrant inductive limits and the composition 
$ \coh ( \D)\to 
  \ind ( \coh ( \D) )\to \Mod ( \D)$
  is isomorphic to $\iota$. 
 The functor $J \iota$ is fully faithful. 
 We denote by $\indcoh  ( \D) $ the essential image of 
 $J \iota$.
 By definition, the category 
 $\indcoh  ( \D) $
 is the subcategory of 
  $ \Mod ( \D) $
  consisting of objects which are filtrant inductive limits of objects of 
  $\coh ( \D)$.
 Since $\fP $ is noetherian, 
 the category $\coh ( \D)$ is essentially small. 
 From  
 \cite[8.6.5.(vi)]{Kashiwara-schapira-book}, 
 this yields that 
 $\indcoh  ( \D) $ is a Grothendieck category.

We set $D ^\mathrm{b} 
_{\mathrm{indcoh}}
( \D)
:=D ^\mathrm{b} _{ \indcoh  ( \D)} ( \Mod (\D))$.

Replacing $\D$ by $D$, we define the categories 
$\Mod  ( D)$,$\coh  ( D)$, $\indcoh  ( D) $.
\end{empt}

\begin{lem}
We keep the notation of \ref{IndCoh}.

\begin{enumerate}[(a)]
\item We have the equalities 
$\coh  ( \D)
=
 \Mod  ( \D) ^{\mathrm{fp}}
= \indcoh  ( \D) ^{\mathrm{fp}}$.

\item Suppose $\fP$ affine. 
We have the equalities
$\indcoh  ( D) 
=
\Mod  ( D)$,
$\coh  ( D)
= 
\Mod  ( D) ^{\mathrm{fp}}$.
\end{enumerate}

\end{lem}

\begin{lem}
\label{eqcatindcoh}
We keep the notation of \ref{IndCoh}.
We suppose $\fP$ affine. 
\begin{enumerate}[(a)]
\item The functors $\D 
\otimes _{D} -$ and $\Gamma (\fP, -)$ 
induce quasi-inverse
equivalences of categories between
$ \indcoh  ( \D)$
and 
$ \Mod ( D)$
(resp. $ \coh  ( \D)$
and 
$ \coh  ( D)$).
Moreover
$ \coh  ( D)$
(resp. $ \coh  ( D ^{(m)})$)
is equal to category of 
finitely presented $ D$-module 
(resp. 
the category of 
finitely generated $ D ^{(m)}$-module).

\item For any 
$\E \in \indcoh  ( D)$, 
$q \geq 1$, 
$H ^q ( \fP, \E)  = 0$.
\end{enumerate}

\end{lem}

\begin{prop}
\label{eq-cat-indcoh-m}
We keep the notation of \ref{IndCoh}.
We suppose $\fP$ affine. 
The canonical functor
\begin{equation}
\label{eq-cat-indcoh-m-funct}
D ^\mathrm{b} ( \indcoh  ( \D)) 
\to 
D ^\mathrm{b} 
_{\mathrm{indcoh}}
( \D)
\end{equation}
is an equivalence of categories.
\end{prop}

\begin{lem}
\label{Serresubcat-cohindcoh}
We suppose $\fP$ affine. The category 
$\coh  ( \D ^{(m)} )$ 
is a Serre subcategory of
$\indcoh  ( \D ^{(m)})$. 
\end{lem}

\begin{proof}
This is a consequence of \ref{eqcatindcoh}
and of the fact that 
$D ^{(m)}$ is noetherian (see \cite[3.3-3.4]{Be1}).
\end{proof}

\begin{rem}
It seems false that the category 
$\coh  ( \D ^\dag _{\fP /\fS } (\hdag T) _\Q)$
is a Serre subcategory of
$\indcoh  ( \D ^\dag _{\fP /\fS } (\hdag T) _\Q)$.

\end{rem}

\begin{prop}
\label{eq-cat-coh-m}
We keep the notation of \ref{IndCoh}.
We suppose $\fP$ affine. 
The canonical functor
$$
D ^\mathrm{b} ( \coh  ( \D^{(m)})  )
\to 
D ^\mathrm{b} _{\mathrm{coh}}
( \D ^{(m)})$$ 
is an equivalence of categories.
\end{prop}

\begin{cor}
\label{cor-eq-cat-coh-m}
We keep the notation of \ref{IndCoh}.
We suppose $\fP$ affine. 
The canonical functors
\begin{gather}
\label{cor-eq-cat-coh-m-1}
D ^\mathrm{b} ( \coh  ( \D ^\dag _{\fP /\fS } (\hdag T) _\Q ))
\to 
D ^\mathrm{b} _{\mathrm{coh}}
(\D ^\dag _{\fP /\fS } (\hdag T) _\Q),
\\
\label{cor-eq-cat-coh-m-2}
D ^{\mathrm{b}}  (\underrightarrow{LM} _{\Q,\mathrm{coh}} (\smash{\widetilde{\D}} _{\fP /\fS } ^{(\bullet)} (T)))
\to 
\underrightarrow{LD} ^{\mathrm{b}} _{\Q, \mathrm{coh}} (\smash{\widetilde{\D}} _{\fP /\fS } ^{(\bullet)} (T))
\end{gather}
are essentially surjective.
\end{cor}

\begin{prop}
Let $\fU:= \fP \setminus T$ be the open formal $\fS $-scheme.
Let $\E \in  \indcoh  ( \D ^\dag _{\fP /\fS } (\hdag T) _\Q)$.
If $\E | \fU \in  \coh  ( \D ^\dag _{\fU/\fS,\Q} )$
then 
$\E \in  \coh  ( \D ^\dag _{\fP /\fS } (\hdag T) _\Q)$.
\end{prop}

\section{Localization functor outside a divisor}

We keep the notation of chapter \ref{ntn-tildeD(Z)}.

\subsection{Tensor products, quasi-coherence, forgetful functor,
localization functor outside a divisor}

\begin{empt}
For any $\E,\FF \in D  ^-
(\overset{^\mathrm{l}}{} \smash{\widetilde{\D}} _{\fP /\fS } ^{(m)} (T ))$
and $\M \in D  ^-
({}^r \smash{\widetilde{\D}} _{\fP /\fS } ^{(m)} (T ))$,
we set:

\begin{gather} \notag
\M _i := \M \otimes ^\L _{\smash{\widetilde{\D}} _{\fP /\fS } ^{(m)} (T )} \smash{\widetilde{\D}} _{P  _i/ S  _i} ^{(m)} (T ),\
\E _i := \smash{\widetilde{\D}} _{P  _i/ S  _i} ^{(m)} (T ) \otimes ^\L _{\smash{\widetilde{\D}} _{\fP /\fS } ^{(m)} ( T )} \E,\\
\notag
\M \smash{\widehat{\otimes}} ^\L _{\smash{\widetilde{\B}} _{\fP} ^{(m)} (T )} \E :=
\R \underset{\underset{i}{\longleftarrow}}{\lim}\, ( \M _i \otimes ^\L  _{\widetilde{\B} _{P _i} ^{(m)} (T )} \E _i)
,\,
\E \smash{\widehat{\otimes}} ^\L _{\smash{\widetilde{\B}} _{\fP} ^{(m)} (T )} \FF :=
\R \underset{\underset{i}{\longleftarrow}}{\lim}\, ( \E _i \otimes ^\L  _{\widetilde{\B} _{P _i} ^{(m)} (T )} \FF _i),
\\
\M \smash{\widehat{\otimes}} ^\L _{\smash{\widetilde{\D}} _{\fP /\fS } ^{(m)} (T )} \E :=
\R \underset{\underset{i}{\longleftarrow}}{\lim}\, ( \M _i \otimes ^\L  _{\widetilde{\D} _{P  _i/S  _i} ^{(m)} (T )} \E _i).
\end{gather}
\end{empt}

\begin{empt}
For any
$\E ^{(\bullet)} \in 
D ^{-}
(\overset{^\mathrm{l}}{} \smash{\widetilde{\D}} _{\fP /\fS } ^{(\bullet)} (T ))$,
$\M ^{(\bullet)} \in 
D ^{-}
( {} ^r \smash{\widetilde{\D}} _{\fP /\fS } ^{(\bullet)} (T )  )$,
we set
\begin{gather}
\label{predef-otimes-coh0}
  \M ^{(\bullet)}
\smash{\widehat{\otimes}} ^\L _{\smash{\widetilde{\D}} ^{(\bullet)} _{\fP /\fS }( T) }\E ^{(\bullet)}
:=
(\M ^{(m)}  \smash{\widehat{\otimes}} ^\L _{\widetilde{\D} ^{(m)} _{\fP /\fS } ( T) } \E ^{(m)}) _{m\in \N}.
\end{gather}
For $? = r$ or $? = l$, 
we define the following tensor product bifunctor
\begin{align}
\label{predef-otimes-coh1}
-
 \smash{\widehat{\otimes}}
^\L _{\widetilde{\B} ^{(\bullet)}  _{\fP} ( T) }
 -
 \colon
D ^-
(\overset{^\mathrm{?}}{} \smash{\widetilde{\D}} _{\fP /\fS } ^{(\bullet)} (T ))
\times 
D ^-
(\overset{^\mathrm{l}}{} \smash{\widetilde{\D}} _{\fP /\fS } ^{(\bullet)} (T ))
&
\to 
D ^-
(\overset{^\mathrm{?}}{} \smash{\widetilde{\D}} _{\fP /\fS } ^{(\bullet)} (T )),
\end{align}
by setting, for any
$\E ^{(\bullet)}\in 
D ^{-}
(\overset{^\mathrm{?}}{} \smash{\widetilde{\D}} _{\fP /\fS } ^{(\bullet)} (T ))$,
$\FF ^{(\bullet)} \in 
D ^{-}
(\overset{^\mathrm{l}}{} \smash{\widetilde{\D}} _{\fP /\fS } ^{(\bullet)} (T ))$,
\begin{equation}
\notag
\E ^{(\bullet)}
\smash{\widehat{\otimes}}
^\L _{\widetilde{\B} ^{(\bullet)}  _{\fP} ( T) }\FF ^{(\bullet)}
:=
(\E ^{(m)}  
\smash{\widehat{\otimes}}
^\L _{\widetilde{\B} ^{(m)}  _{\fP} ( T) }
\FF ^{(m)}) _{m\in \N}.
\end{equation}

When $T$ is empty, $\widetilde{\B} ^{(\bullet)}  _{\fP} ( T) $
will simply be denoted by 
$\O _{\fP} ^{(\bullet)}$, i.e. 
$\O _{\fP} ^{(\bullet)}$ is the subring of 
$\widehat{\D} ^{(\bullet)} _{\fP/\fS}$ whose transition morphisms
are the identity of $\O _{\fP}$.

\end{empt}

\begin{ntn}
[Quasi-coherence and partial forgetful functor of the divisor]
Let $D \subset T$ be a second divisor.

\begin{itemize}
\item Let  $\E ^{(m)} \in 
D ^{\mathrm{b}}
(\overset{^\mathrm{l}}{} \smash{\widetilde{\D}} _{\fP /\fS } ^{(m)} (T ))$.
Since  $\smash{\widetilde{\D}} _{\fP /\fS } ^{(m)}(T)$ (resp. $\smash{\widetilde{\B}} _{\fP} ^{(m)}(T)$)
has not $p$-torsion,
using the Theorem \cite[3.2.2]{Beintro2} 
we get that 
$\E ^{(m)} $ is quasi-coherent in the sense of Berthelot 
as object of
$D ^{\mathrm{b}}
(\overset{^\mathrm{l}}{} \smash{\D} _{\fP} ^{(m)})$ (see his definition \cite[3.2.1]{Beintro2}) 
if and only if 
$\E ^{(m)} _0 \in D ^{\mathrm{b}} _{\mathrm{qc}}
(\O  _{P})$
and 
the canonical morphism 
$\E ^{(m)} \to 
\smash{\widetilde{\D}} _{\fP /\fS } ^{(m)} (T ) 
\smash{\widehat{\otimes}} ^\L _{\smash{\widetilde{\D}} _{\fP /\fS } ^{(m)} (T )} \E ^{(m)} $
(resp. $\E ^{(m)} \to 
\smash{\widetilde{\B}} _{\fP} ^{(m)} (T ) 
\smash{\widehat{\otimes}} ^\L _{\smash{\widetilde{\B}} _{\fP} ^{(m)} (T )} \E ^{(m)} $)
is an isomorphism. 
In particular, this does not depend on the divisor $T$.
We denote by 
$D ^{\mathrm{b}} _{\mathrm{qc}}
(\overset{^\mathrm{l}}{} \smash{\widetilde{\D}} _{\fP /\fS } ^{(m)} (T ))$, 
the full subcategory of
$D ^{\mathrm{b}}
(\overset{^\mathrm{l}}{} \smash{\widetilde{\D}} _{\fP /\fS } ^{(m)} (T ))$
of quasi-coherent complexes.
We get the
{\it partial forgetful functor of the divisor} 
$$\mathrm{oub} _{D,T}\colon 
D ^{\mathrm{b}} _{\mathrm{qc}}
(\overset{^\mathrm{l}}{} \smash{\widetilde{\D}} _{\fP /\fS } ^{(m)} (T))
\to
D ^{\mathrm{b}} _{\mathrm{qc}}
(\overset{^\mathrm{l}}{} \smash{\widetilde{\D}} _{\fP /\fS } ^{(m)} (D ))$$
which is induced by the canonical forgetful functor
$\mathrm{oub} _{D,T}\colon 
D ^{\mathrm{b}} 
(\overset{^\mathrm{l}}{} \smash{\widetilde{\D}} _{\fP /\fS } ^{(m)} (T))
\to
D ^{\mathrm{b}} 
(\overset{^\mathrm{l}}{} \smash{\widetilde{\D}} _{\fP /\fS } ^{(m)} (D ))$.

\item Similarly, we denote by
$D ^{\mathrm{b}} _{\mathrm{qc}}
(\overset{^\mathrm{l}}{} \smash{\widetilde{\D}} _{\fP /\fS } ^{(\bullet)} (T ))$
the full subcategory of 
$D ^{\mathrm{b}}
(\overset{^\mathrm{l}}{} \smash{\widetilde{\D}} _{\fP /\fS } ^{(\bullet)} (T ))$
of complexes $\E ^{(\bullet)} $ such that, for any $m\in \Z$,  
$\E ^{(m)} _0 \in D ^{\mathrm{b}} _{\mathrm{qc}}
(\O  _{P})$ 
and the canonical morphism 
$\E ^{(\bullet)}  \to 
\smash{\widetilde{\D}} _{\fP /\fS } ^{(\bullet)} (T )
\smash{\widehat{\otimes}} ^\L _{\smash{\widetilde{\D}} _{\fP /\fS } ^{(\bullet)} (T )} \E ^{(\bullet)} $
is an isomorphism of $D ^{\mathrm{b}}
(\overset{^\mathrm{l}}{} \smash{\widetilde{\D}} _{\fP /\fS } ^{(\bullet)} (T ))$. 
We get the
{\it partial forgetful functor of the divisor} 
$$\mathrm{oub} _{D,T}\colon 
D ^{\mathrm{b}} _{\mathrm{qc}}
(\overset{^\mathrm{l}}{} \smash{\widetilde{\D}} _{\fP /\fS } ^{(\bullet)} (T))
\to
D ^{\mathrm{b}}  _{\mathrm{qc}}
(\overset{^\mathrm{l}}{} \smash{\widetilde{\D}} _{\fP /\fS } ^{(\bullet)} (D )).$$

\item We denote by 
$\smash{\underrightarrow{LD}}  ^{\mathrm{b}} _{\Q,\mathrm{qc}}
( \smash{\widetilde{\D}} _{\fP /\fS } ^{(\bullet)} (T ))$ 
the strictly full subcategory of 
$\smash{\underrightarrow{LD}}  ^{\mathrm{b}} _{\Q}
( \smash{\widetilde{\D}} _{\fP /\fS } ^{(\bullet)} (T ))$
of complexes which are isomorphic in
$\smash{\underrightarrow{LD}}  ^{\mathrm{b}} _{\Q}
( \smash{\widetilde{\D}} _{\fP /\fS } ^{(\bullet)} (T ))$ 
to a complex belonging to 
$D ^{\mathrm{b}}  _{\mathrm{qc}}
(\overset{^\mathrm{l}}{} \smash{\widetilde{\D}} _{\fP /\fS } ^{(\bullet)} (T ))$.
Since the functor $\mathrm{oub} _{D,T}$
sends a lim-ind-isogeny to a lim-ind-isogeny, 
we obtain the factorization of the form : 
\begin{equation}
\label{def-oubDT}
\mathrm{oub} _{D,T}
\colon 
\smash{\underrightarrow{LD}}  ^{\mathrm{b}} _{\Q,\mathrm{qc}}
( \smash{\widetilde{\D}} _{\fP /\fS } ^{(\bullet)} (T ))
\to
\smash{\underrightarrow{LD}}  ^{\mathrm{b}} _{\Q,\mathrm{qc}}
( \smash{\widetilde{\D}} _{\fP /\fS } ^{(\bullet)} (D )).
\end{equation}

\item We still denote by 
$\mathrm{oub} _{D, T}\colon 
D ^{\mathrm{b}} (\D ^\dag _{\fP /\fS } (\hdag T) _\Q) \to 
D ^{\mathrm{b}} (\D ^\dag _{\fP /\fS } (\hdag D) _\Q)$
the partial forgetful functor of the divisor.

\end{itemize}
\end{ntn}

\begin{rem}
\begin{enumerate}[(a)]
\item A morphism 
$\E ^{(\bullet)}
\to \cF ^{(\bullet)}$ 
of 
$D ^{\mathrm{b}}
(\overset{^\mathrm{l}}{} \smash{\widetilde{\D}} _{\fP /\fS } ^{(\bullet)} (T ))$
is an isomorphism if and only if 
the induced morphism
$\E ^{(m)}
\to \cF ^{(m)}$ is an isomorphism of 
$D ^{\mathrm{b}}
(\overset{^\mathrm{l}}{} \smash{\widetilde{\D}} _{\fP /\fS } ^{(m)} (T ))$ for every $m\in \Z$.

\item Let 
$\E ^{(\bullet)}
\in 
D ^{\mathrm{b}}
(\overset{^\mathrm{l}}{} \smash{\widetilde{\D}} _{\fP /\fS } ^{(\bullet)} (T ))$.
Using the first remark, 
we check the property $\E ^{(\bullet)}
\in 
D ^{\mathrm{b}} _{\mathrm{qc}}
(\overset{^\mathrm{l}}{} \smash{\widetilde{\D}} _{\fP /\fS } ^{(\bullet)} (T ))$
is equivalent to
the property that,  
for any $m \in \Z$, 
$\E ^{(m)}
\in 
D ^{\mathrm{b}} _{\mathrm{qc}}
(\overset{^\mathrm{l}}{} \smash{\widetilde{\D}} _{\fP /\fS } ^{(m)} (T ))$.
Hence, the above definition of 
$\smash{\underrightarrow{LD}}  ^{\mathrm{b}} _{\Q,\mathrm{qc}}
( \smash{\widetilde{\D}} _{\fP /\fS } ^{(\bullet)} (T ))$ 
corresponds to that of Berthelot's one formulated in 
 \cite[4.2.3]{Beintro2} without singularities along a divisor. 
\end{enumerate}

\end{rem}

\begin{lemm}
\label{lemm-def-otimes-coh1}
The bifunctor \ref{predef-otimes-coh1} induces
\begin{align}
\label{def-otimes-coh1}
-
 \smash{\widehat{\otimes}}
^\L _{\widetilde{\B} ^{(\bullet)}  _{\fP} ( T) }
 -
 \colon
 \smash{\underrightarrow{LD}}  ^- _{\Q}
(\overset{^\mathrm{?}}{} \smash{\widetilde{\D}} _{\fP /\fS } ^{(\bullet)} (T ))
\times 
\smash{\underrightarrow{LD}}  ^- _{\Q}
(\overset{^\mathrm{l}}{} \smash{\widetilde{\D}} _{\fP /\fS } ^{(\bullet)} (T ))
&
\to 
\smash{\underrightarrow{LD}}  ^- _{\Q}
(\overset{^\mathrm{?}}{} \smash{\widetilde{\D}} _{\fP /\fS } ^{(\bullet)} (T )).
\end{align}
\end{lemm}

\begin{proof}
This is similar to \cite[2.1.5]{caro-6operations}.
\end{proof}

\begin{empt}
Let $D \subset T$ be a second divisor. 
For any
$\E ^{(\bullet)} \in 
D ^{-}
(\overset{^\mathrm{l}}{} \smash{\widetilde{\D}} _{\fP /\fS } ^{(\bullet)} (D))$, 
similarly to \cite[1.1.8]{caro_courbe-nouveau} we get 
the commutative diagram in $D ^{-}
(\overset{^\mathrm{l}}{} \smash{\widetilde{\D}} _{\fP /\fS } ^{(\bullet)} (T))$:
\begin{equation}
\label{ODdivcohe}
  \xymatrix @ R=0,4cm {
{(\widetilde{\B} ^{(m)} _{\fP} ( T)  \smash{\widehat{\otimes}} ^\L
_{\widetilde{\B} ^{(m)}  _{\fP} ( D) } \E ^{(m)}) _{m\in \N}}
\ar@{=}[r] ^-{\mathrm{def}} \ar[d] _\sim
&
{\widetilde{\B} ^{(\bullet)} _{\fP} ( T)  \smash{\widehat{\otimes}} ^\L
_{\widetilde{\B} ^{(\bullet)}  _{\fP} ( D) } \E ^{(\bullet)}}
\ar@{.>}[d] _\sim
\\
{(\widetilde{\D} ^{(m)} _{\fP } ( T)  \smash{\widehat{\otimes}} ^\L
_{\widetilde{\D} ^{(m)} _{\fP } ( D) } \E ^{(m)}) _{m\in \N}}
\ar@{=}[r] ^-{\mathrm{def}}
&
{\smash{\widetilde{\D}} ^{(\bullet)} _{\fP /\fS }( T)  \smash{\widehat{\otimes}} ^\L
_{\smash{\widetilde{\D}} ^{(\bullet)} _{\fP /\fS }( D) }\E ^{(\bullet)}=: (\hdag T, D) (\E ^{(\bullet)})}.
}
\end{equation}
As for Lemma \ref{lemm-def-otimes-coh1}, 
we get the the localization outside $T$ functor :
\begin{gather}
\label{hdag-def}
 (\hdag T ,\,D) 
 :=\smash{\widetilde{\D}} ^{(\bullet)} _{\fP /\fS }( T)  \smash{\widehat{\otimes}} ^\L
_{\smash{\widetilde{\D}} ^{(\bullet)} _{\fP /\fS }( D) }-
\colon
\smash{\underrightarrow{LD}} ^{-} _{\Q} ( \smash{\widetilde{\D}} _{\fP /\fS } ^{(\bullet)}(D))
\to
\smash{\underrightarrow{LD}} ^{-} _{\Q} ( \smash{\widetilde{\D}} _{\fP /\fS } ^{(\bullet)}(T)).
\end{gather}

\end{empt}

\subsection{Preservation of the quasi-coherence}
\label{section3.2}
Let $m' \geq m \geq 0$ be two integers,
$D ' \subset D \subset T$ be three (reduced) divisors of  $P$. 
We have the canonical morphisms
$\widetilde{\B} _{P _i} ^{(m)} (D ') \to \widetilde{\B} _{P _i} ^{(m)} (D )\to \widetilde{\B} _{P _i} ^{(m')} (T)$.
Similarly to the notation of \cite{Beintro2},
we denote by $D _{\Q,\mathrm{qc}} ^{-} (\smash{\widetilde{\B}} _{\fP} ^{(m)} (D))$ 
(resp. 
$D _{\Q,\mathrm{qc}} ^{-} (\widetilde{\B} ^{(m')} _{\fP} ( D)  \smash{\widehat{\otimes}} _{\O _{\fP}} \smash{\widehat{\D}} _{\fP /\fS } ^{(m)})$)
the localization of the category 
$D ^{-} _{\mathrm{qc}} (\smash{\widetilde{\B}} _{\fP} ^{(m)} (D))$ 
(resp.  $D _{\mathrm{qc}} ^{-} (\widetilde{\B} ^{(m')} _{\fP} ( D)  \smash{\widehat{\otimes}} _{\O _{\fP}} \smash{\widehat{\D}} _{\fP /\fS } ^{(m)})$)
by isogenies.

\begin{lemm}
\label{lem1-hdagDT}
\begin{enumerate}[(a)]
\item The kernel of the canonical epimorphism
$\smash{\widetilde{\B}} _{\fP} ^{(m)} (D ) \widehat{\otimes} _{\O _{\fP} }
\smash{\widetilde{\B}} _{\fP} ^{(m')} (T )
\to 
\smash{\widetilde{\B}} _{\fP} ^{(m')} (T )$
is a quasi-coherent $\O _{P}$-module. 

\item The canonical morphism
$\smash{\widetilde{\B}} _{\fP} ^{(m)} (D ) 
\widehat{\otimes} ^{\L}  _{\O _{\fP} }
\smash{\widetilde{\B}} _{\fP} ^{(m')} (T )
\to 
\smash{\widetilde{\B}} _{\fP} ^{(m)} (D ) \widehat{\otimes} _{\O _{\fP} }
\smash{\widetilde{\B}} _{\fP} ^{(m')} (T )$
is an isomorphism.

\end{enumerate}
\end{lemm}

\begin{proof}
We can copy word by word the proof of  \cite[3.2.1]{caro-stab-sys-ind-surcoh}.\end{proof}

\begin{empt}
Let us clarify some terminology.

\begin{enumerate}[(a)]
\item A morphism of rings 
$f \colon \AA \to \B$ is a $p ^n$-isogeny if there exists 
a morphisms of rings $g \colon \B \to \AA$ such that $f \circ g = p ^n id$
and
$g \circ f = p ^n id$.

\item A morphism $f \colon \AA \to \B$ of 
$D ^- ( \smash{\widetilde{\B}} _{\fP} ^{(m')} (T ))$
is a $p ^n$-isogeny if there exists a morphisms  $g \colon \B \to \AA$ of $D ^- ( \smash{\widetilde{\B}} _{\fP} ^{(m')} (T ))$
such that $f \circ g = p ^n id$
and
$g \circ f = p ^n id$.
\end{enumerate}

\end{empt}

\begin{prop}
\label{lem3-hdagDT}
The canonical homomorphisms of 
$D ^- ( \smash{\widetilde{\B}} _{\fP} ^{(m')} (T ))$
or respectively of rings
\begin{equation}
\label{lem3-hdagDT-iso}
\smash{\widetilde{\B}} _{\fP} ^{(m')} (T ) \to
\smash{\widetilde{\B}} _{\fP} ^{(m)} (D ) \widehat{\otimes} ^{\L}_{\smash{\widetilde{\B}} _{\fP} ^{(m)} (D' )}
\smash{\widetilde{\B}} _{\fP} ^{(m')} (T )
\to 
\smash{\widetilde{\B}} _{\fP} ^{(m)} (D ) \widehat{\otimes} _{\smash{\widetilde{\B}} _{\fP} ^{(m)} (D' )}
\smash{\widetilde{\B}} _{\fP} ^{(m')} (T )
\to 
\smash{\widetilde{\B}} _{\fP} ^{(m')} (T )
\end{equation}
are $p$-isogenies.
\end{prop}

\begin{proof}
We can copy word by word the proof of \cite[3.2.2]{caro-stab-sys-ind-surcoh}.\end{proof}
\begin{coro}
\label{rema-dim-coh-finie}
\begin{enumerate}[(a)]
\item The functors of the form
$\B _{P _i} ^{(m')} (T) \otimes ^{\L}_{\O _{P _i} }- $ have cohomological dimension 
$1$. 
The functor $\smash{\widetilde{\B}} _{\fP} ^{(m')} (T) 
\widehat{\otimes} ^{\L} _{\O _{\fP} }-$
is way-out over
$D  ^{-} (\O _{\fP} )$
with bounded amplitude independent of $m'$ and $m$.

\item 
The functor 
$\smash{\widetilde{\B}} _{\fP} ^{(m')} (T) \widehat{\otimes} ^{\L} _{\smash{\widetilde{\B}} _{\fP} ^{(m)} (D)}-
\colon 
D _{\Q,\mathrm{qc}} ^{\mathrm{b}} (\smash{\widetilde{\B}} _{\fP} ^{(m)} (D))
\to 
D _{\Q,\mathrm{qc}} ^{\mathrm{b}} (\smash{\widetilde{\B}} _{\fP} ^{(m')} (T))$ 
is way-out
with bounded amplitude independent of $m'$ and $m$.
We have the factorization 
$\smash{\widetilde{\B}} _{\fP} ^{(m+\bullet)} (T) \widehat{\otimes} ^{\L} _{\smash{\widetilde{\B}} _{\fP} ^{(m)} (D)}-
\colon 
D _{\Q,\mathrm{qc}} ^{\mathrm{b}} (\smash{\widetilde{\B}} _{\fP} ^{(m)} (D))
\to
\smash{\underrightarrow{LD}} _{\Q,\mathrm{qc}} ^{\mathrm{b}} (\smash{\widetilde{\B}} _{\fP} ^{(m+\bullet)} (T ))$.

\item 
The functor
$$(\widetilde{\B} ^{(m')} _{\fP} ( T)  \smash{\widehat{\otimes}} _{\O _{\fP}} \smash{\widehat{\D}} _{\fP /\fS } ^{(m)})
\widehat{\otimes} ^{\L} _{(\widetilde{\B} ^{(m)} _{\fP} ( D)  \smash{\widehat{\otimes}} _{\O _{\fP}} \smash{\widehat{\D}} _{\fP /\fS } ^{(m)})}
-
\colon 
D _{\Q,\mathrm{qc}} ^{\mathrm{b}} (\widetilde{\B} ^{(m)} _{\fP} ( D)  \smash{\widehat{\otimes}} _{\O _{\fP}} \smash{\widehat{\D}} _{\fP /\fS } ^{(m)})
\to 
D _{\Q,\mathrm{qc}} ^{\mathrm{b}} (\widetilde{\B} ^{(m')} _{\fP} (T)  \smash{\widehat{\otimes}} _{\O _{\fP}} \smash{\widehat{\D}} _{\fP /\fS } ^{(m)})$$
is way-out
with bounded amplitude independent of $m'$ and $m$.
\end{enumerate}

\end{coro}

\begin{proof}
We can copy word by word the proof of \cite[3.2.3]{caro-stab-sys-ind-surcoh}.\end{proof}

\begin{empt}
\label{hdagT-nota}
With Corollary  \ref{rema-dim-coh-finie} which implies the stability of the boundedness of the cohomology,
we check the factorization of the functor of \ref{hdag-def} as follows:
\begin{gather}
\label{hdag-def-qc}
 (\hdag T ,\,D) 
 :=\smash{\widetilde{\D}} ^{(\bullet)} _{\fP /\fS }( T)  \smash{\widehat{\otimes}} ^\L
_{\smash{\widetilde{\D}} ^{(\bullet)} _{\fP /\fS }( D) }-
\colon
\smash{\underrightarrow{LD}} ^{\mathrm{b}} _{\Q,\mathrm{qc}} ( \smash{\widetilde{\D}} _{\fP /\fS } ^{(\bullet)}(D))
\to
\smash{\underrightarrow{LD}} ^{\mathrm{b}} _{\Q,\mathrm{qc}} ( \smash{\widetilde{\D}} _{\fP /\fS } ^{(\bullet)}(T)).
\end{gather}
We also write 
$ \E ^{(\bullet)} (\hdag D ,\,T) :=(\hdag T ,\,D) (\E ^{(\bullet)})$. 
This functor  $(\hdag T ,\,D)$ is 
 {\it the localization outside $T$ functor}. 
When  $D=\emptyset $, we omit writing it. 
We write in the same way the associated functor for coherent complexes:
\begin{equation}
\label{hdag-def-coh}
(\hdag T, D) := 
\D ^\dag _{\fP /\fS } (\hdag T) _\Q \otimes _{ \D ^\dag _{\fP /\fS } (\hdag D) _\Q} - \colon 
D ^\mathrm{b} _\mathrm{coh} ( \D ^\dag _{\fP /\fS } (\hdag D) _\Q)
\to 
D ^\mathrm{b} _\mathrm{coh} (\D ^\dag _{\fP /\fS } (\hdag T) _\Q).
\end{equation}
The  functor \ref{hdag-def-coh} is exact, which justifies the absence of the symbol 
$\L$. 
\end{empt}

\begin{prop}
\label{oub-pl-fid}
Let  $\E ^{(\bullet)}\in 
\smash{\underrightarrow{LD}} ^{\mathrm{b}} _{\Q,\mathrm{qc}}
 ( \smash{\widetilde{\D}} _{\fP /\fS } ^{(\bullet)}(T))$.
\begin{enumerate}[(a)]
\item The functorial in $\E ^{(\bullet)}$ canonical morphism :
\begin{equation}
\label{oub-pl-fid-iso1}
(\hdag T ,\,D) \circ \mathrm{oub} _{D,T} (\E ^{(\bullet)})
\to 
\E ^{(\bullet)}
\end{equation}
is an isomorphism of
$\smash{\underrightarrow{LD}} ^{\mathrm{b}} _{\Q,\mathrm{qc}} ( \smash{\widetilde{\D}} _{\fP /\fS } ^{(\bullet)}(T))$.
\item The functorial in $\E ^{(\bullet)}$ canonical morphism :
\begin{equation}
\label{oub-pl-fid-iso2}
\mathrm{oub}_{D, T}  (\E ^{ (\bullet)} )
\to
\mathrm{oub}_{D, T} \circ (\hdag T, D) \circ \mathrm{oub}_{D, T}(\E ^{ (\bullet)} ) 
\end{equation}
is an isomorphism
of $\underrightarrow{LD} ^{\mathrm{b}}  _{\Q, \mathrm{qc}}  ( \smash{\widetilde{\D}} _{\fP /\fS } ^{(\bullet)}(D))$.
\item 
\label{oub-pl-fid-iso3}
The functor
$\mathrm{oub} _{D,T}\colon 
\underrightarrow{LD} ^{\mathrm{b}}  _{\Q, \mathrm{qc}} ( \smash{\widetilde{\D}} _{\fP /\fS } ^{(\bullet)}(T))
\to 
\underrightarrow{LD} ^{\mathrm{b}}  _{\Q, \mathrm{qc}} ( \smash{\widetilde{\D}} _{\fP /\fS } ^{(\bullet)}(D))$
is fully faithful.
\end{enumerate}

\end{prop}

\begin{proof}
We can copy word by word the proof of \cite[3.2.6]{caro-stab-sys-ind-surcoh}.\end{proof}

\begin{coro}
\label{gen-oub-pl-fid}
Let $\E ^{(\bullet)}\in \smash{\underrightarrow{LD}} ^{\mathrm{b}} _{\Q,\mathrm{qc}} ( \smash{\widetilde{\D}} _{\fP /\fS } ^{(\bullet)}(D))$.
The functorial in $\E ^{(\bullet)}$ canonical morphism 
\begin{equation}
\label{gen-oub-pl-fid-iso1}
(\hdag T ,\,D ') \circ \mathrm{oub} _{D',D} (\E ^{(\bullet)})
\to 
(\hdag T ,~D) (\E ^{(\bullet)})
\end{equation}
is an isomorphism of
$\smash{\underrightarrow{LD}} ^{\mathrm{b}} _{\Q,\mathrm{qc}} ( \smash{\widetilde{\D}} _{\fP /\fS } ^{(\bullet)}(T))$.
\end{coro}

\begin{proof}
We can copy word by word the proof of \cite[3.2.7]{caro-stab-sys-ind-surcoh}.\end{proof}

\begin{ntn}
\label{nota-hag-sansrisque}
Let $D \subset T \subset T' $ be some divisors of $P$.
Following \ref{gen-oub-pl-fid},
by forgetting to write some forgetful functors, 
the functors $(\hdag T',~D) $ 
and
$(\hdag T', ~T) $ are canonically isomorphic over
$\smash{\underrightarrow{LD}} ^{\mathrm{b}} _{\Q ,\mathrm{qc}}
(\smash{\widetilde{\D}} _{\fP } ^{(\bullet)}(T))$.
Hence, we can simply write
$(\hdag T') $ in both case.

\end{ntn}

\begin{ntn}
\label{nota-tdf}
We denote by 
$D  ^{\mathrm{b}} _{\mathrm{tdf}}
( \smash{\widetilde{\D}} _{\fP /\fS } ^{(\bullet)} (T ))$
the full subcategory of 
$D  ^{\mathrm{b}} 
(\smash{\widetilde{\D}} _{\fP /\fS } ^{(\bullet)} (T ))$
consisting of complexes of finite Tor-dimension. 
We denote by 
$\smash{\underrightarrow{LD}}  ^{\mathrm{b}} _{\Q, \mathrm{qc}, \mathrm{tdf}}
(\smash{\widetilde{\D}} _{\fP /\fS } ^{(\bullet)} (T ))$
the strictly full subcategory of 
$\smash{\underrightarrow{LD}}  ^{\mathrm{b}} _{\Q, \mathrm{qc}}
(\smash{\widetilde{\D}} _{\fP /\fS } ^{(\bullet)} (T ))$
consisting of objects isomorphic in 
$\smash{\underrightarrow{LD}}  ^{\mathrm{b}} _{\Q, \mathrm{qc}}
(\smash{\widetilde{\D}} _{\fP /\fS } ^{(\bullet)} (T ))$
to an object of $D  ^{\mathrm{b}} _{\mathrm{tdf}}
( \smash{\widetilde{\D}} _{\fP /\fS } ^{(\bullet)} (T ))$.

\end{ntn}

\begin{coro}
\begin{enumerate}[(a)]
\label{def-otimes-coh1&2qc}
\item 
The bifunctor \ref{def-otimes-coh1} factorizes throught the bifunctor 
\begin{align}
\label{def-otimes-coh1qc}
-
 \smash{\widehat{\otimes}}
^\L _{\widetilde{\B} ^{(\bullet)}  _{\fP} ( T) }
 -
\colon
\smash{\underrightarrow{LD}}  ^{\mathrm{b}} _{\Q, \mathrm{qc}}
(\overset{^\mathrm{?}}{} \smash{\widetilde{\D}} _{\fP /\fS } ^{(\bullet)} (T ))
\times 
\smash{\underrightarrow{LD}}  ^{\mathrm{b}} _{\Q, \mathrm{qc}}
(\overset{^\mathrm{l}}{} \smash{\widetilde{\D}} _{\fP /\fS } ^{(\bullet)} (T ))
&
\to 
\smash{\underrightarrow{LD}}  ^{\mathrm{b}}  _{\Q, \mathrm{qc}}
(\overset{^\mathrm{?}}{} \smash{\widetilde{\D}} _{\fP /\fS } ^{(\bullet)} (T )).
\end{align}

\item 
\label{def-otimes-coh2qc}
With notation  \ref{nota-tdf}, 
we have the equality 
$\smash{\underrightarrow{LD}}  ^{\mathrm{b}} _{\Q, \mathrm{qc}}
(\smash{\widetilde{\D}} _{\fP /\fS } ^{(\bullet)} (T ))
=\smash{\underrightarrow{LD}}  ^{\mathrm{b}} _{\Q, \mathrm{qc}, \mathrm{tdf}}
(\smash{\widetilde{\D}} _{\fP /\fS } ^{(\bullet)} (T ))$.

\end{enumerate}
\end{coro}

\begin{proof}
We can copy word by word the proof of \cite[3.2.9]{caro-stab-sys-ind-surcoh} (for the second statement, the careful reader might notice in fact
we need the slightly more precise argument that the cohomological dimension of our rings can be bounded independently of the level $m$).\end{proof}

\begin{rem}
We have $\smash{\widetilde{\D}} _{\fP /\fS } ^{(m)} (T ) \in 
D  ^{\mathrm{b}} _{\mathrm{tdf}}
( \smash{\widehat{\D}} _{\fP /\fS } ^{(0)} )$, with flat amplitude bounded independently of the level $m$.
Indeed, following 
\ref{finitecohdimDdagpre}, 
$ \D _{P /S } ^{(0)} $ has finite homological dimension. 
Hence, 
$\D _{P /S } ^{(m)} \in
D  ^{\mathrm{b}} _{\mathrm{qc}}
( \D _{P /S } ^{(0)} )
= 
D  ^{\mathrm{b}} _{\mathrm{qc}, \mathrm{tdf}}
( \D _{P /S } ^{(0)} )$, with flat amplitude bounded by 
the cohomological dimension of $ \D _{P /S } ^{(0)}$
(see \cite[I.5.9]{sga6}).
Since 
$\smash{\widehat{\D}} _{\fP /\fS } ^{(m)} \in 
D  ^{\mathrm{b}} _{\mathrm{qc}}
( \smash{\widehat{\D}} _{\fP /\fS } ^{(0)} )$,
then from \cite[3.2.3]{Beintro2} (still valid in our context), $\smash{\widehat{\D}} _{\fP /\fS } ^{(m)} \in 
D  ^{\mathrm{b}} _{\mathrm{tdf}}
( \smash{\widehat{\D}} _{\fP /\fS } ^{(0)} )$, with flat amplitude bounded by 
the cohomological dimension of $ \smash{\widehat{\D}} _{\fP /\fS } ^{(0)}$
(in fact the proof of  \cite[3.2.3]{Beintro2} shows more precisely the preservation of flat amplitude).
Then, using 
\ref{rema-dim-coh-finie},
$\smash{\widetilde{\D}} _{\fP /\fS } ^{(m)} (T ) \in 
D  ^{\mathrm{b}} _{\mathrm{tdf}}
( \smash{\widehat{\D}} _{\fP /\fS } ^{(0)} )$, with flat amplitude bounded independently of the level $m$.
Hence, 
$\smash{\widetilde{\D}} _{\fP /\fS } ^{(\bullet)} (T ) \in 
D  ^{\mathrm{b}} _{\mathrm{tdf}}
( \smash{\widehat{\D}} _{\fP /\fS } ^{(\bullet)} )$.

\end{rem}

\begin{coro}
Let
$\M ^{(\bullet)}
\in \smash{\underrightarrow{LD}}  ^{\mathrm{b}} _{\Q, \mathrm{qc}}
(\overset{^\mathrm{?}}{} \smash{\widetilde{\D}} _{\fP /\fS } ^{(\bullet)}(D))$,
and
$\E ^{(\bullet)} \in \smash{\underrightarrow{LD}}  ^{\mathrm{b}} _{\Q, \mathrm{qc}}
(\smash{\widetilde{\D}} _{\fP /\fS } ^{(\bullet)}(D))$.
We have the canonical isomorphism in 
$\smash{\underrightarrow{LD}} ^{\mathrm{b}} _{\Q, \mathrm{qc}} 
(\overset{^\mathrm{?}}{}  \smash{\widetilde{\D}} _{\fP /\fS } ^{(\bullet)}(T))$
of the form 
\begin{equation}
\label{hdagTDotimes}
 (\hdag T ,\,D) ( \M ^{(\bullet)} )
 \smash{\widehat{\otimes}} ^\L
_{\widetilde{\B} ^{(\bullet)}  _{\fP} ( T) } 
  (\hdag T ,\,D) ( \E ^{(\bullet)})
  \riso 
   (\hdag T ,\,D)  
   \left (
   \M ^{(\bullet)}
   \smash{\widehat{\otimes}} ^\L
_{\widetilde{\B} ^{(\bullet)}  _{\fP} ( D) } 
   \E ^{(\bullet)}
      \right ) .
\end{equation}
\end{coro}

\begin{proof}
Using the bounded quasi-coherence of our objects, 
this is straightforward from the associativity of the tensor products
(use the equivalence of categories of 
\cite[3.2.3]{Beintro2} to reduce to the case of usual tensor products of complexes).
\end{proof}

\begin{coro}
Let
$\M ^{(\bullet)}
\in \smash{\underrightarrow{LD}}  ^{\mathrm{b}} _{\Q, \mathrm{qc}}
(\overset{^\mathrm{?}}{} \smash{\widetilde{\D}} _{\fP /\fS } ^{(\bullet)}(T))$,
and
$\E ^{(\bullet)} \in \smash{\underrightarrow{LD}}  ^{\mathrm{b}} _{\Q, \mathrm{qc}}
(\smash{\widetilde{\D}} _{\fP /\fS } ^{(\bullet)}(T))$.
We have the isomorphism
\begin{equation}
\label{oubTDDotimes}
\mathrm{oub} _{D,T} ( \M ^{(\bullet)} )
 \smash{\widehat{\otimes}} ^\L
_{\widetilde{\B} ^{(\bullet)}  _{\fP} ( D) } 
\mathrm{oub} _{D,T} ( \E ^{(\bullet)})
  \riso 
\mathrm{oub} _{D,T}
   \left (
   \M ^{(\bullet)}
   \smash{\widehat{\otimes}} ^\L
_{\widetilde{\B} ^{(\bullet)}  _{\fP} ( T) } 
   \E ^{(\bullet)}
      \right ) .
\end{equation}

\end{coro}

\begin{proof}
Using \ref{oub-pl-fid-iso1}, we get 
$   \M ^{(\bullet)}
   \smash{\widehat{\otimes}} ^\L
_{\widetilde{\B} ^{(\bullet)}  _{\fP} ( T) } 
   \E ^{(\bullet)}
\riso
 \M ^{(\bullet)} 
 \smash{\widehat{\otimes}} ^\L
_{\widetilde{\B} ^{(\bullet)}  _{\fP} ( T) } 
\left (
\widetilde{\B} ^{(\bullet)}  _{\fP} ( T) 
 \smash{\widehat{\otimes}} ^\L
_{\widetilde{\B} ^{(\bullet)}  _{\fP} ( D) } 
 (\mathrm{oub} _{D,T} ( \E ^{(\bullet)}))
 \right ) $.
We conclude by associativity of the tensor product.\end{proof}

\subsection{Composition of localisation functors}

\begin{lemm}
\label{lem1-hdagT1T2}
Let $T, ~T'$ be two divisors of  $P$ 
whose irreducible components are distinct, 
$\U''$ the open set of  $\fP$ complementary to $T \cup T'$.

\begin{enumerate}[(a)]
\item For any $i \in \N$, 
the canonical morphism
$\smash{\widetilde{\B}} _{P _i} ^{(m)} (T ) \otimes ^{\L}_{\O _{P _i} }
\smash{\widetilde{\B}} _{P _i} ^{(m)} (T' ) \to 
\smash{\widetilde{\B}} _{P _i} ^{(m)} (T ) \otimes _{\O _{P _i} }
\smash{\widetilde{\B}} _{P _i} ^{(m)} (T' )$
is an isomorphism.

\item The canonical morphism
$\smash{\widetilde{\B}} _{\fP} ^{(m)} (T ) \widehat{\otimes} ^{\L} _{\O _{\fP} }
\smash{\widetilde{\B}} _{\fP} ^{(m)} (T ')
\to
\smash{\widetilde{\B}} _{\fP} ^{(m)} (T ) \widehat{\otimes} _{\O _{\fP} }
\smash{\widetilde{\B}} _{\fP} ^{(m)} (T ')$
is an isomorphism
and 
the $\O _{\fP} $-algebra
$\smash{\widetilde{\B}} _{\fP} ^{(m)} (T ) \widehat{\otimes} _{\O _{\fP} }
\smash{\widetilde{\B}} _{\fP} ^{(m)} (T ')$
has no $p$-torsion.

\item The canonical morphism of  $\O _{\fP}$-algebras
$\smash{\widetilde{\B}} _{\fP} ^{(m)} (T ) \widehat{\otimes} _{\O _{\fP} }
\smash{\widetilde{\B}} _{\fP} ^{(m)} (T ')
\to 
j _* \O _{\U ''} $,
where $j\colon \U '' \hookrightarrow \fP$ is the inclusion,
is a monomorphism.

\item Let  $\chi,~\lambda \colon \N \to \N$ defined respectively by setting for any integer $m\in \N$ 
$\chi (m) := p ^{p-1}$ and $\lambda (m) := m +1$.
We have two canonical monomorphisms 
$\alpha ^{(\bullet)}
\colon
\smash{\widetilde{\B}} _{\fP} ^{(\bullet)} (T ) \widehat{\otimes} _{\O _{\fP} }
\smash{\widetilde{\B}} _{\fP} ^{(\bullet)} (T ')
\to 
\smash{\widetilde{\B}} _{\fP} ^{(\bullet)} (T \cup T')$
and
$\beta ^{(\bullet)}\colon 
\smash{\widetilde{\B}} _{\fP} ^{(\bullet)} (T \cup T')
\to 
\lambda ^{*} \chi ^{*} (
\smash{\widetilde{\B}} _{\fP} ^{(\bullet)} (T ) \widehat{\otimes} _{\O _{\fP} }
\smash{\widetilde{\B}} _{\fP} ^{(\bullet)} (T '))$
such that 
$\lambda ^{*} \chi ^{*} (\alpha  ^{(\bullet)})
\circ \beta  ^{(\bullet)}$ 
and
$\beta  ^{(\bullet)} \circ \alpha  ^{(\bullet)}$ 
are the canonical morphisms.
\end{enumerate}

\end{lemm}

\begin{proof}
We can copy word by word the proof of \cite[3.2.10]{caro-stab-sys-ind-surcoh}.\end{proof}

\begin{prop}
\label{hdagT'T=cup}
Let  $T', T $ be two divisors of $P$.
For any $\E ^{(\bullet)}
\in 
\underrightarrow{LD} ^{\mathrm{b}} _{\Q,\mathrm{qc}} (\smash{\widetilde{\D}} _{\fP /\fS } ^{(\bullet)})$, 
we have the isomorphism
$(\hdag T ') \circ (\hdag T) (\E ^{(\bullet)})
\to
(T '\cup T) (\E ^{(\bullet)})$
functorial in  $T,~T',~\E ^{(\bullet)}$.
\end{prop}

\begin{proof}
Using \ref{lem1-hdagT1T2}, 
we can copy word by word the proof of \cite[3.2.11]{caro-stab-sys-ind-surcoh}.\end{proof}

\subsection{A coherence stability criterion by localisation outside a divisor}

\begin{thm}
\label{limTouD}
Let  $T' \supset T$ be a divisor,
$\E ^{(\bullet)}
\in 
\underrightarrow{LD} ^{\mathrm{b}} _{\Q, \mathrm{coh}} (\smash{\widetilde{\D}} _{\fP } ^{(\bullet)} (T))$
and
$\E := 
\underrightarrow{\lim}
\E ^{(\bullet)}
\in D ^{\mathrm{b}} _{\mathrm{coh}} (\smash{\D} ^\dag _{\fP } (\hdag T) _{\Q})$.
We suppose that the morphism 
$\E \to (\hdag T',T) (\E)$
is an isomorphism of 
$D ^{\mathrm{b}}  (\smash{\D} ^\dag _{\fP } (\hdag T) _{\Q})$.
Then, the canonical morphism
$\E ^{(\bullet)} \to 
(\hdag T',T) (\E ^{(\bullet)})$
is an isomorphism of 
$\underrightarrow{LD} ^{\mathrm{b}} _{\Q, \mathrm{coh}} (\smash{\widetilde{\D}} _{\fP } ^{(\bullet)}(T))$.
\end{thm}

\begin{proof}
We can copy \cite[3.5.1]{caro-stab-sys-ind-surcoh}.
\end{proof}

\begin{cor}
\label{coro1limTouD}
Let  $T' \supset T$ be a divisor,
$\E ^{\prime (\bullet)}
\in 
\underrightarrow{LD} ^{\mathrm{b}} _{\Q, \mathrm{coh}} (\smash{\widetilde{\D}} _{\fP } ^{(\bullet)} (T'))$
and
$\E ':= 
\underrightarrow{\lim}
\E ^{\prime (\bullet)}
\in D ^{\mathrm{b}} _{\mathrm{coh}} (\smash{\D} ^\dag _{\fP } (\hdag T ') _{\Q})$.
If
$\E '\in D ^{\mathrm{b}}  _{\mathrm{coh}} (\smash{\D} ^\dag _{\fP } (\hdag T) _{\Q})$,
then
$\E ^{\prime (\bullet)}
\in 
\underrightarrow{LD} ^{\mathrm{b}} _{\Q, \mathrm{coh}} (\smash{\widetilde{\D}} _{\fP } ^{(\bullet)}(T))$.
\end{cor}

\begin{coro}
\label{coro2limTouD}
Let  $T' \supset T$ be a divisor,
$\E \in D ^{\mathrm{b}}  _{\mathrm{coh}} (\smash{\D} ^\dag _{\fP } (\hdag T) _{\Q}) \cap 
D ^{\mathrm{b}}  _{\mathrm{coh}} (\smash{\D} ^\dag _{\fP } (\hdag T') _{\Q})$.
Let 
$\E ^{(\bullet)}
\in 
\underrightarrow{LD} ^{\mathrm{b}} _{\Q, \mathrm{coh}} (\smash{\widetilde{\D}} _{\fP } ^{(\bullet)} (T))$
and
$\E ^{\prime (\bullet)}
\in 
\underrightarrow{LD} ^{\mathrm{b}} _{\Q, \mathrm{coh}} (\smash{\widetilde{\D}} _{\fP } ^{(\bullet)} (T'))$
such that we have the $\smash{\D} ^\dag _{\fP } (\hdag T) _{\Q}$-linear isomorphisms of the form
$\underrightarrow{\lim}
\E ^{(\bullet)}
\riso 
\E$
and
$\underrightarrow{\lim}
\E ^{\prime (\bullet)}
\riso 
\E$.
Then, we have the isomorphism
$\E ^{(\bullet)} \riso \E ^{\prime (\bullet)}$
of
$\underrightarrow{LD} ^{\mathrm{b}} _{\Q, \mathrm{coh}} (\smash{\widetilde{\D}} _{\fP } ^{(\bullet)} (T))$.
\end{coro}

\begin{proof}
This is a straightforward consequence of \ref{coro1limTouD} and 
of the full faithfulness of the functor
$\underrightarrow{\lim}$ 
on 
$\underrightarrow{LD} ^{\mathrm{b}} _{\Q, \mathrm{coh}} (\smash{\widetilde{\D}} _{\fP } ^{(\bullet)} (T))$.
\end{proof}

\begin{prop}
\label{stab-coh-oub-DT}
Let  $T \subset D \subset T' $ be some divisors of  $P$.
\begin{enumerate}[(a)]
\item Let 
$\E ^{(\bullet)}
\in 
\underrightarrow{LD} ^{\mathrm{b}} _{\Q, \mathrm{coh}} (\smash{\widetilde{\D}} _{\fP } ^{(\bullet)}(T))
\cap \underrightarrow{LD} ^{\mathrm{b}} _{\Q, \mathrm{coh}} (\smash{\widetilde{\D}} _{\fP } ^{(\bullet)} (T'))$.
Then 
$\E ^{(\bullet)}
\in 
 \underrightarrow{LD} ^{\mathrm{b}} _{\Q, \mathrm{coh}} (\smash{\widetilde{\D}} _{\fP } ^{(\bullet)} (D))$.
\item Let 
$\E \in 
D ^{\mathrm{b}} _{\mathrm{coh}}( \smash{\D} ^\dag _{\fP } (\hdag T) _{\Q} )
\cap 
D ^{\mathrm{b}} _{\mathrm{coh}}( \smash{\D} ^\dag _{\fP } (\hdag T') _{\Q} )$.
Then 
$\E 
\in 
D ^{\mathrm{b}} _{\mathrm{coh}}( \smash{\D} ^\dag _{\fP } (\hdag D) _{\Q} )$.
\end{enumerate}
\end{prop}

\begin{proof}
Using \ref{oub-pl-fid-iso1}, we check that the canonical morphism 
$(\hdag D,\,T) \circ 
\mathrm{oub} _{T ,D} 
( \mathrm{oub} _{D ,T'} (\E ^{(\bullet)}))
\to 
\mathrm{oub} _{D ,T'} (\E ^{(\bullet)})$
of $ \underrightarrow{LD} ^{\mathrm{b}} _{\Q, \mathrm{qc}} (\smash{\widetilde{\D}} _{\fP } ^{(\bullet)} (D))$
is an isomorphism.
Hence, we get the first assertion. 
Using \ref{coro1limTouD}, this yields the second one.
\end{proof}

\begin{ntn}
\label{fct-qcoh2coh}
Let $\fP$ and $\mathfrak{Q}$ be two formal  $\fS$-schemes locally of formal finite type
and having locally finite $p$-bases  over $\fS$.
Let $T$ be a divisor of $P$,
$U$ be a divisor of $Q$,
and 
$\phi ^{(\bullet)}
\colon 
\smash{\underrightarrow{LD}} ^{\mathrm{b}} _{\Q,\mathrm{qc}}
( \smash{\widetilde{\D}} _{\fP /\fS } ^{(\bullet)}(T))
\to 
\smash{\underrightarrow{LD}} ^{\mathrm{b}} _{\Q,\mathrm{qc}}
( \smash{\widetilde{\D}} _{\mathfrak{Q} /\fS } ^{(\bullet)}(U))$
be a functor.
We denote by 
$\mathrm{Coh} _{T} ( \phi ^{(\bullet)})
\colon 
D ^\mathrm{b} _\mathrm{coh} ( \smash{\D} ^\dag _{\fP /\fS } (\hdag T) _{\Q} ) \to
D ^\mathrm{b}  ( \smash{\D} ^\dag _{\mathfrak{Q} /\fS } (\hdag U) _{\Q} )$
the functor 
$\mathrm{Coh} _{T} ( \phi ^{(\bullet)}) := \underrightarrow{\lim} \circ \phi ^{(\bullet)}\circ (\underrightarrow{\lim} _T ) ^{-1}$,
where $(\underrightarrow{\lim} _T ) ^{-1}$ is a quasi-inverse functor of 
the equivalence of categories
\begin{equation}
\label{eq-catLDBer-LD-D}
\underrightarrow{\lim} 
\colon 
\underrightarrow{LD} ^{\mathrm{b}}  _{\Q, \mathrm{coh}} (\smash{\widetilde{\D}} _{\fP } ^{(\bullet)} (T))
\overset{\ref{eqcatLD=DSM-fonct-coh}}{\cong} 
D ^{\mathrm{b}} _{\mathrm{coh}}
(\underrightarrow{LM} _{\Q} (\smash{\widetilde{\D}} _{\fP } ^{(\bullet)} (T)))
\overset{\ref{eqcatcoh}}{\cong}
D ^{\mathrm{b}} _{\mathrm{coh}}( \smash{\D} ^\dag _{\fP } (\hdag T) _{\Q} ).
\end{equation}
\end{ntn}

\begin{rem}
\label{rema-fct-qcoh2coh}
Let  $T \subset T'$ be a second divisor. 
Let  $\E \in 
D ^\mathrm{b} _\mathrm{coh} ( \smash{\D} ^\dag _{\fP } (\hdag T') _{\Q} ) 
\cap D ^\mathrm{b} _\mathrm{coh} ( \smash{\D} ^\dag _{\fP } (\hdag T) _{\Q} )$.
Using \ref{coro2limTouD},  
the corresponding objects 
of $\smash{\underrightarrow{LD}} ^{\mathrm{b}} _{\Q ,\mathrm{coh}}
(\smash{\widetilde{\D}} _{\fP } ^{(\bullet)}(T))$
and
$\smash{\underrightarrow{LD}} ^{\mathrm{b}} _{\Q ,\mathrm{coh}}
(\smash{\widetilde{\D}} _{\fP } ^{(\bullet)}(T'))$
(via the equivalence of categories \ref{eq-catLDBer-LD-D})
are isomorphic. 
With notation \ref{fct-qcoh2coh}, the functors 
$\mathrm{Coh} _{T} ( \phi ^{(\bullet)}) $
and
$\mathrm{Coh} _{T'} ( \phi ^{(\bullet)}) $
are then isomorphic over 
$D ^\mathrm{b} _\mathrm{coh} ( \smash{\D} ^\dag _{\fP } (\hdag T') _{\Q} ) 
\cap D ^\mathrm{b} _\mathrm{coh} ( \smash{\D} ^\dag _{\fP } (\hdag T) _{\Q} )$.
\end{rem}

\begin{rem}
\label{rem-hag-sansrisque}

\begin{itemize}

\item For any divisors $D \subset T$, 
we have the isomorphism of functors 
$\mathrm{Coh} _{D} ((\hdag T' ,~D) ) \riso (\hdag T' ,~D) $
(see notation \ref{hdagT-nota}) 
Hence, both notation of \ref{hdagT-nota} are compatible.

\item

Let $T$ and $D \subset D' $ be some divisors of $P$.
We obtain the functor
$(\hdag T ) := \mathrm{Coh} _{D} ((\hdag T) )
\colon D ^\mathrm{b} _\mathrm{coh} ( \smash{\D} ^\dag _{\fP } (\hdag D) _{\Q} ) 
\to D ^\mathrm{b} _\mathrm{coh} ( \smash{\D} ^\dag _{\fP } (\hdag T \cup D) _{\Q} ) $
(see notation \ref{nota-hag-sansrisque}).
With the remark \ref{rema-fct-qcoh2coh}, 
since the functors 
$\mathrm{Coh} _{D} ((\hdag T) )$ and
$\mathrm{Coh} _{D'} ((\hdag T) )$
are isomorphic over
 $D ^\mathrm{b} _\mathrm{coh} ( \smash{\D} ^\dag _{\fP } (\hdag D) _{\Q} ) 
\cap D ^\mathrm{b} _\mathrm{coh} ( \smash{\D} ^\dag _{\fP } (\hdag D') _{\Q} )$,
then it is not necessary to clarify $D$.

\end{itemize}

\end{rem}

\section{Extraordinary inverse image, direct image, duality, base change}

\subsection{Definitions of the functors}
\label{subsection3.1}
Let  $f \colon \fP ^{\prime } \to \fP $ be a morphism of 
 formal  $\fS$-schemes locally of formal finite type
and having locally finite $p$-bases  over $\fS$,
$T$ and $T'$ be some divisors of respectively $P$ and $P'$ such that 
$f ( P '\setminus T' ) \subset P \setminus T$.

We define in this section the extraordinary inverse image and 
direct image by $f$ with overconvergent singularities along $T$ and $T'$,
and the dual functor.

We fix  $\lambda _0\colon \N \to \N$ an increasing map such that 
$\lambda _{0} (m) \geq m$ for any $m \in \N$. 
We set 
$\widetilde{\B} ^{(m)} _{\fP} ( T):= \widehat{\B} ^{(\lambda _0 (m))} _{\fP} ( T)$ 
et
$\smash{\widetilde{\D}} _{\fP /\fS } ^{(m)} (T):=
\widetilde{\B} ^{(m)} _{\fP} ( T)  \smash{\widehat{\otimes}} _{\O _{\fP}} \smash{\widehat{\D}} _{\fP /\fS } ^{(m)}$.
Finally, we set 
$\smash{\D} _{P  _i /S  _i} ^{(m)} (T):= \V / \pi ^{i+1} \otimes _{\V} \smash{\widehat{\D}} _{\fP /\fS  } ^{(m)} (T) 
=
\B ^{(m)} _{P _i} ( T)  \otimes _{\O _{P _i}} \smash{\D} _{P  _i/S  _i} ^{(m)}$
and
$\smash{\widetilde{\D}} _{P  _i/S  _i} ^{(m)} (T):=\widetilde{\B} ^{(m)} _{P _i} ( T)  \otimes _{\O _{P _i}} \smash{\D} _{P  _i/S  _i} ^{(m)}$.
We use similar notation by adding some primes,
e.g. 
$\widetilde{\B} ^{(m)} _{\fP'} ( T'):= \widehat{\B} ^{(\lambda _0 (m))} _{\fP'} ( T')$ .

\begin{ntn}
\begin{enumerate}[(a)]
\item Since $f ^{-1} (T) \subset T'$, 
we get the canonical morphism 
$f ^{-1} \widetilde{\B} _{P _i} ^{(m)} ( T ) \rightarrow \widetilde{\B} _{P' _i} ^{(m)} ( T ')$.
Hence, the sheaf
$\widetilde{\B} _{P' _i} ^{(m)} ( T ') \otimes _{\O _{P ^\prime _i}} f _i ^* \D _{P  _i/ S  _i} ^{(m)}
\riso 
\widetilde{\B} _{P' _i} ^{(m)} ( T ') \otimes _{f ^{-1} \widetilde{\B} _{P _i} ^{(m)} ( T ) } f ^{-1} \widetilde{\D} _{P  _i/ S  _i} ^{(m)} (T)$
is endowed with a canonical structure 
of ($ \widetilde{\D} _{P  ^{\prime } _i/ S  _i} ^{(m)} ( T ')$, $ f  ^{-1} \widetilde{\D} _{P  _i/ S  _i} ^{(m)}(T)$)-bimodule.
We denote this bimodule by $\widetilde{\D} ^{(m)} _{P ^{\prime } _i \rightarrow P  _i/ S  _i} ( T' , T)$.

\item By $p$-adic completion, we get the following
$(\smash{\widetilde{\D}} _{\fP ^{\prime }/\fS }  ^{(m)}(T') , f ^{-1} \smash{\widetilde{\D}} _{\fP / \fS } ^{(m)} (T))$-bimodule :
$\smash{\widetilde{\D}} _{\fP ^{\prime }\rightarrow \fP / \fS } ^{(m)} ( T' , T):=
\underset{\underset{i}{\longleftarrow}}{\lim}\, \widetilde{\D}_{P ^{\prime } _i \rightarrow P  _i/ S  _i} ^{(m)}( T' , T)$.

\item We get a $(\D ^{\dag } _{\fP ^{\prime }/ \fS } (\hdag T' )_{\Q}, f ^{-1} \D ^{\dag }_{\fP ^{ }/ \fS } (\hdag T) _{\Q})$-bimodule
by setting
$\D ^{\dag} _{\fP ^{\prime }\rightarrow \fP ^{ }/ \fS } (\hdag T' ,T) _{\Q}:=\underset{\underset{m}{\longrightarrow}}{\lim}\,
\smash{\widetilde{\D}} _{\fP ^{\prime }\rightarrow \fP ^{ }/ \fS } ^{(m)} ( T' ,T)_{\Q}$.

\end{enumerate}

\end{ntn}

\begin{empt}
[Extraordinary inverse image]
\label{ntn-Lf!+*}
\begin{enumerate}[(a)]
\item  The extraordinary inverse image by $f$ with overconvergent singularities along $T$ and $T'$ 
is a functor of the form
$f  ^{!(\bullet)} _{T',T} \colon
\smash{\underrightarrow{LD}} ^{\mathrm{b}} _{\Q,\mathrm{qc}} ( \smash{\widetilde{\D}} _{\fP /\fS } ^{(\bullet)}(T))
\to
\smash{\underrightarrow{LD}} ^{\mathrm{b}} _{\Q,\mathrm{qc}} ( \smash{\widetilde{\D}} _{\fP ^{\prime }/\fS }  ^{(\bullet)}(T'))$
which is defined 
for any $\E ^{(\bullet)} \in \smash{\underrightarrow{LD}}  ^\mathrm{b} _{\Q, \mathrm{qc}}
( \smash{\widetilde{\D}} _{\fP /\fS } ^{(\bullet)} (T ))$ by setting:
\begin{equation}
\notag
f _{T',T} ^{ !(\bullet)} ( \E ^{(\bullet)}) :=
\smash{\widetilde{\D}} ^{(\bullet)} _{\fP ^{\prime } \rightarrow \fP /\fS } (T',T)
\smash{\widehat{\otimes}} ^\L _{f ^{-1} \smash{\widetilde{\D}} ^{(\bullet)} _{\fP /\fS } (T)}
f ^{-1} \E ^{(\bullet)} [  \delta _{\fP ' /\fP}],
\end{equation}
where the tensor product is defined similarly to \ref{predef-otimes-coh0}.
\item The extraordinary inverse image by $f$ with overconvergent singularities along $T$ and $T'$ 
is also a functor of the form
$f  ^{ !} _{T',T} \colon
D ^\mathrm{b} _\mathrm{coh} ( \smash{\D} ^\dag _{\fP /\fS } (\hdag T) _{\Q} )
\to 
 D ^\mathrm{b} ( \smash{\D} ^\dag _{\fP ^{\prime }/\fS }  (\hdag T') _{\Q} )$
which is defined 
for any $\E \in D ^{\mathrm{b}} _{\mathrm{coh}} ( \D ^{\dag} _{\fP /\fS } (\hdag T ) _{\Q})$ by setting:
\begin{equation}
\label{def-image-inv-extr}
f  ^{ !} _{T' , T} (\E ):=\D ^{\dag} _{\fP ^{\prime }\rightarrow \fP  }  ( \hdag T' , T ) _{\Q}
\otimes ^{\L} _{ f ^{-1} \D ^{\dag} _{\fP /\fS } (\hdag T ) _{\Q}} f ^{-1} \E [\delta _{\fP '/\fP} ].
\end{equation}

\item Mostly when $f$ is flat, we can also consider the functors
$\L f _{T',T} ^{*(\bullet)} 
:= 
f _{T',T} ^{ !(\bullet)} 
[ -  \delta _{\fP ' /\fP}]$,
and
$\L f  ^{ *} _{T',T}  : = f  ^{ !} _{T',T} [ -  \delta _{\fP ' /\fP}]$.
Beware that our notation might be misleading since 
$\L f _{T',T} ^{*(\bullet)} $ is not  necessarily a left derived functor of some functor (except for coherent complexes).
When $f$ is flat, these  functors are t-exact over coherent complexes, and we denote them respectively
$f _{T',T} ^{*(\bullet)} $ and $f  ^{ *} _{T',T} $.

\item When $T '= f ^{-1} (T)$, we simply write respectively
$  f _{T} ^{ !(\bullet)} $, $f   _{T} ^!$, and $f _T ^*$. If moreover $T$ is empty, 
we write $  f^{ !(\bullet)} $, 
$f ^!$, and $f ^*$.

\end{enumerate}

\end{empt}

\begin{ntn}

\begin{enumerate}[(a)]
\item We define a ($f ^{-1} \widetilde{\D} _{P  _i/ S  _i} ^{(m)}(T)$, $ \widetilde{\D} _{P ^{\prime } _i} ^{(m)} ( T ')$)-bimodule by setting
$$\widetilde{\D} ^{(m)} _{P  _i \leftarrow P ^{\prime   }_i/ S  _i} ( T , T'):=
\widetilde{\B} _{P' _i} ^{(m)} ( T ') \otimes _{\O _{P ^\prime _i}}
\left (
\omega _{P ^{\prime } _i/S  _i} \otimes _{\O _{P ' _i}}f ^* _l \left ( \D _{P  _i/ S  _i} ^{(m)}(T) \otimes _{\O _{P _i}} \omega ^{-1} _{P  _i/ S  _i} \right)\right),$$
where the symbol $l$ means that we choose the left structure of left $\D _{P  _i/ S  _i} ^{(m)}(T) $-module.

\item This yields by completion the 
$(f ^{-1} \smash{\widetilde{\D}} _{\fP / \fS } ^{(m)} (T),~\smash{\widetilde{\D}} _{\fP ^{\prime }/ \fS } ^{(m)}(T') )$-bimodule :
$$\smash{\widetilde{\D}} _{\fP  \leftarrow \fP ^{\prime }/ \fS } ^{(m)} ( T, T'):=
\underset{\underset{i}{\longleftarrow}}{\lim}\, \widetilde{\D}_{P  _i \leftarrow P ^{\prime } _i/ S  _i} ^{(m)}( T , T').$$

\item We get the
($f ^{-1} \D ^{\dag }_{\fP  / \fS } (\hdag T) _{\Q}$, $\D ^{\dag } _{\fP ^{\prime }/ \fS } (\hdag T' )_{\Q}$)-bimodule
$\D ^{\dag} _{\fP  \leftarrow \fP ^{\prime }/ \fS } (\hdag T ,T ')_{\Q}:=
\underset{\underset{m}{\longrightarrow}}{\lim}
\smash{\widetilde{\D}} _{\fP  \leftarrow \fP ^{\prime }/ \fS } ^{(m)} ( T, T') _\Q$.

\end{enumerate}

\end{ntn}

\begin{empt}
\begin{enumerate}[(a)]
\item The
direct image  by $f$ with overconvergent singularities along $T$ and $T'$ 
is a functor of the form
$f  ^{ (\bullet)} _{T,T',+} \colon
\smash{\underrightarrow{LD}} ^{\mathrm{b}} _{\Q,\mathrm{qc}} ( \smash{\widetilde{\D}} _{\fP ^{\prime }/\fS }  ^{(\bullet)}(T'))
\to 
\smash{\underrightarrow{LD}} ^{\mathrm{b}} _{\Q,\mathrm{qc}} ( \smash{\widetilde{\D}} _{\fP /\fS } ^{(\bullet)}(T))$
defined by setting,
for any $\E ^{\prime (\bullet)} \in \smash{\underrightarrow{LD}}  ^\mathrm{b} _{\Q, \mathrm{qc}}
( \smash{\widetilde{\D}} _{\fP ^{\prime }/\fS }  ^{(\bullet)} (T '))$:
\begin{gather}\notag
f _{T,T',+} ^{ (\bullet)} ( \E ^{\prime (\bullet)} ):= 
\R f _* (
\smash{\widetilde{\D}} ^{(\bullet)} _{\fP  \leftarrow \fP ^{\prime }/\fS } (T,T')
\smash{\widehat{\otimes}} ^\L _{\smash{\widetilde{\D}} ^{(\bullet)} _{\fP ^{\prime }/\fS }  (T')}
\E ^{\prime (\bullet)}).
\end{gather}

\item The
direct image by $f$ with overconvergent singularities along $T$ and $T'$ is a functor of the form
$f _{T,T',+} \colon
D ^\mathrm{b} _\mathrm{coh} ( \smash{\D} ^\dag _{\fP ^{\prime }/\fS }  (\hdag T') _{\Q} )
\to 
 D ^\mathrm{b} ( \smash{\D} ^\dag _{\fP /\fS } (\hdag T) _{\Q} )$,
defined by setting,
for any $\E '\in D ^{\mathrm{b}} _{\mathrm{coh}} ( \D ^{\dag} _{\fP ^{\prime }/\fS }  (\hdag T' ) _{\Q})$ :
\begin{equation}
\label{ftt'+}
f   _{T , T ', +}( \E '):=
 \R f_* (\D ^{\dag} _{\fP  \leftarrow \fP ^{\prime }/\fS }  ( \hdag T , T ') _{\Q}
\otimes ^{\L} _{ \D ^{\dag} _{\fP ^{\prime }/\fS }  (\hdag T ') _{\Q}} \E ').
\end{equation}

\item When $T '= f ^{-1} (T)$, we simply write respectively
$f  ^{ (\bullet)} _{T,+} $ and $f   _{T ,+}$. If moreover $T$ is empty, 
we write $f  ^{ (\bullet)} _{+} $ and $f   _{+}$. 
\end{enumerate}

\end{empt}

\begin{empt}
\label{coh-Qcoh}
With notation \ref{fct-qcoh2coh}, 
we have the isomorphism of functors
$\mathrm{Coh} _{T'} (f  ^{ (\bullet)}_{T , T ', +}) \riso f  ^{} _{T , T ', +}$ 
and
$\mathrm{Coh} _{T} ( f _{T',T} ^{ !(\bullet)}) \riso  f _{T',T} ^{ !}$
(this is checked similarly to \cite[4.3.2.2 and 4.3.7.1]{Beintro2}).

\end{empt}

\begin{ntn}
[Duality]
\label{ntn-dualfunctor}
\begin{enumerate}[(a)]
\item Let $\E \in D ^{\mathrm{b}} _{\mathrm{coh}}( \smash{\D} ^\dag _{\fP /\fS \Q} )$.
The $\smash{\D} ^\dag _{\fP /\fS ,\Q}$-linear dual of $\E$ is defined by setting
$$\DD  (\E):= \R \mathcal{H} om _{ \smash{\D} ^\dag _{\fP /\fS ,\Q} }
(\E,
 \smash{\D} ^\dag _{\fP /\fS ,\Q} 
\otimes  _{\O _{\fP}} \omega _{\fP ^{ }/\fS } ^{-1})[\delta _{P}].$$
Following \ref{finitecohdimDdag}, 
we get 
$D ^{\mathrm{b}} _{\mathrm{coh}}( \smash{\D} ^\dag _{\fP /\fS ,\Q} )
=
D ^{\mathrm{b}} _{\mathrm{parf}}( \smash{\D} ^\dag _{\fP /\fS ,\Q} )$, 
where the right category is that of perfect bounded complexes of 
$ \smash{\D} ^\dag _{\fP /\fS ,\Q} $-modules. 
This yields 
$ \DD  (\E)
 \in 
 D ^{\mathrm{b}} _{\mathrm{coh}}( \smash{\D} ^\dag _{\fP /\fS ,\Q} )$.
 Hence, by biduality, we get the equivalence of categories
$ \DD 
\colon 
 D ^{\mathrm{b}} _{\mathrm{coh}}( \smash{\D} ^\dag _{\fP /\fS ,\Q} )
 \cong
  D ^{\mathrm{b}} _{\mathrm{coh}}( \smash{\D} ^\dag _{\fP /\fS ,\Q} )$.
  
\item   We denote by
$\DD ^{(\bullet)} 
\colon 
\smash{\underrightarrow{LD}} ^{\mathrm{b}} _{\Q,\mathrm{coh}}
( \smash{\widetilde{\D}} _{\fP /\fS } ^{(\bullet)}(T))
\to 
\smash{\underrightarrow{LD}} ^{\mathrm{b}} _{\Q,\mathrm{coh}}
( \smash{\widetilde{\D}} _{\fP /\fS } ^{(\bullet)}(T))$
the equivalence of categories such that
$\mathrm{Coh} _{T} (\DD ^{(\bullet)} )
\riso 
\DD $.

\end{enumerate}

\end{ntn}

\begin{dfn}
[Base change]
\label{chg-base}
Let $\alpha \colon \V \to \W$ be a morphism of local algebras
 such that 
 $\cV$ and $\W$ are  complete discrete valued ring  of mixed characteristic $(0,p)$ with perfect residue fields.
We set $\fS: =\Spf \cV$ and $\fT: =\Spf \cW$.
Let $r\geq 0$ be an integer, 
let $\X$ be a 
formal  $\bbD ^r _{\fS}$-scheme of finite type, having locally finite $p$-bases  over $\fS$,
$\E ^{(\bullet)} \in \smash{\underrightarrow{LD}} ^{\mathrm{b}} _{\Q,\mathrm{qc}} ( \smash{\widehat{\D}} _{\X/\fS} ^{(\bullet)})$,
$\fY := \X \times _{\bbD ^r _{\fS}} \bbD ^r _{\fT}$, and $\varpi  \colon \fY \to \fX$ 
be the projection.
Following \ref{ex-ft-Dn}, 
$\varpi$ is flat and 
$\fY$ is a formal  $\bbD ^r _{\fT}$-scheme of finite type, having locally finite $p$-bases  over $\fT$.

The ``base change of $\E ^{(\bullet)} $ induced by $\alpha$'' is 
the object $\varpi  ^{*(\bullet)} (\E ^{(\bullet)} )$ of 
$\smash{\underrightarrow{LD}} ^{\mathrm{b}} _{\Q,\mathrm{qc}} ( \smash{\widehat{\D}} _{\fY/\fT} ^{(\bullet)})$ (see \cite[2.2.2]{Beintro2}).
The object 
$\varpi  ^{*(\bullet)} (\E ^{(\bullet)} )$ can simply be denoted by 
$$ \cO _{\bbD ^r _{\fT}} \smash{\widehat{\otimes}}^\L
_{\cO _{\bbD ^r _{\fS }}}  \E^{(\bullet) }.$$

Following \ref{subsect-comm-bc},
base changes 
 commute with push forwards
 base change commutes with quasi-projective extraordinary pullbacks, local cohomological functors,
duals functors (for coherent complexes), and tensor products is 
straightforward. 

\end{dfn}

\subsection{Commutation of pullbacks with localization functors outside of a divisor}

We keep notation \ref{subsection3.1}.

\begin{lemm}
\label{lem-f!B}
Suppose $T':= f ^{-1} (T)$.
We have the canonical isomorphism 
$$\O _{P ' _i} \otimes  ^\L
_{f ^{-1} \O _{P _i}} f ^{-1} \B ^{(m)}  _{P _i} ( T) 
\riso 
\B ^{(m)}  _{P' _i} ( T') .$$
We have also the canonical isomorphism 
$f ^{!(\bullet)} (\widetilde{\B} ^{(\bullet)}  _{\fP } ( T) )
\riso 
\widetilde{\B} ^{(\bullet)}  _{\fP '} ( T')[\delta _{P'/P}]$
in
$\smash{\underrightarrow{LD}} ^\mathrm{b} _{\Q, \mathrm{qc}}
(\overset{^\mathrm{l}}{} \smash{\widetilde{\D}} _{\fP ^{\prime}/\fS } ^{(\bullet)}( T'))$.

\end{lemm}

\begin{proof}
This is checked similarly to  \cite[5.2.1]{caro-stab-sys-ind-surcoh}.
\end{proof}

\begin{empt}
\label{oub-div-opcoh}
\begin{enumerate}[(a)]
\item  \label{oub-div-opcoha)} 
Let  $\E ^{\prime (\bullet)} \in \underrightarrow{LD} ^{\mathrm{b}} _{\Q ,\mathrm{qc}}
( \smash{\widetilde{\D}} _{\fP ^{\prime  }} ^{(\bullet)}(T'))$.
Similarly to \cite[1.1.9 ]{caro_courbe-nouveau}, we check that we have the canonical isomorphism
$oub _{T} \circ f ^{(\bullet)} _{T,T',+} (\E ^{\prime (\bullet)} )\riso f ^{(\bullet)} _+  \circ oub _{T'} (\E ^{\prime (\bullet)})$.
Hence, it is harmless to write by abuse of notation
$f ^{(\bullet)} _{+}$ instead of  $f ^{(\bullet)} _{T, T',+}$.

Using the remark 
\ref{rema-fct-qcoh2coh}
this yields that the functors 
$\mathrm{Coh} _{T'} (f ^{(\bullet)} _{T,T',+}) $
and
$\mathrm{Coh}  (f ^{(\bullet)} _{+}) $ 
are isomorphic over
$D ^\mathrm{b} _\mathrm{coh} ( \smash{\D} ^\dag _{\fP ^{\prime  },\Q} ) 
\cap D ^\mathrm{b} _\mathrm{coh} ( \smash{\D} ^\dag _{\fP ^{\prime  }} (\hdag T') _{\Q} )$.
Since 
we have the canonical isomorphisms of functors 
$\mathrm{Coh} _{T'} (f  ^{ (\bullet)}_{T , T ', +}) \riso f  _{T , T ', +}$ 
and
$\mathrm{Coh} (f  ^{ (\bullet)}_{+}) \riso f _+$ (see \ref{coh-Qcoh}), 
then it is harmless to write 
$f _+$ instead of  $f   _{T, T',+}$
and we get the functor
$f _+
\colon 
D ^\mathrm{b} _\mathrm{coh} ( \smash{\D} ^\dag _{\fP ^{\prime  },\Q} ) 
\cap D ^\mathrm{b} _\mathrm{coh} ( \smash{\D} ^\dag _{\fP ^{\prime  }} (\hdag T') _{\Q} )
\to 
D ^\mathrm{b} _\mathrm{coh} ( \smash{\D} ^\dag _{\fP ,\Q} ) 
\cap D ^\mathrm{b} _\mathrm{coh} ( \smash{\D} ^\dag _{\fP}  (\hdag T) _{\Q} ) )$.

\item 
\label{oub-div-opcohb)}
Let $D$ and $D'$ be some divisors of respectively $P$ and $P'$ such that 
$f ( P '\setminus D' ) \subset P \setminus D$, $D \subset T$, and
$D' \subset T'$.
Let  $\E ^{(\bullet)} \in \underrightarrow{LD} ^{\mathrm{b}} _{\Q ,\mathrm{qc}}
( \smash{\widetilde{\D}} _{\fP /\fS } ^{(\bullet)} (D))$.
Similarly to \cite[1.1.9]{caro_courbe-nouveau}, we check easily the isomorphism
$(\hdag T',D') \circ f ^{!(\bullet)} _{D',D} (\E^{(\bullet)} )
\riso 
f ^{!(\bullet)} _{T',T} \circ (\hdag T, D) (\E^{(\bullet)}) $.

\end{enumerate}
\end{empt}

\begin{empt}
Let
$\FF ^{(\bullet)}
,\E ^{(\bullet)} \in \smash{\underrightarrow{LD}}  ^{\mathrm{b}} _{\Q ,\mathrm{qc}}
(\smash{\widetilde{\D}} _{\fP /\fS } ^{(\bullet)}(T))$.
We easily check  (see \cite[2.1.9.1]{caro-stab-prod-tens})  the following isomorphism of
$\smash{\underrightarrow{LD}} ^{\mathrm{b}} _{\Q ,\mathrm{qc}} ( \smash{\widetilde{\D}} _{\fP ^{\prime }/\fS } ^{(\bullet)}(T'))$
\begin{equation}
\label{f!T'Totimes}
 f ^{!(\bullet)} _{T',T} ( \FF ^{(\bullet)} )
 \smash{\widehat{\otimes}} ^\L
_{\widetilde{\B} ^{(\bullet)}  _{\fP '} ( T ') } 
 f ^{!(\bullet)} _{T',T} ( \E ^{(\bullet)})
  \riso 
f ^{!(\bullet)} _{T',T}
   \left (
   \FF ^{(\bullet)}
   \smash{\widehat{\otimes}} ^\L
_{\widetilde{\B} ^{(\bullet)}  _{\fP} ( T) } 
   \E ^{(\bullet)}
\right ) [\delta _{P'/P}].
\end{equation}

\end{empt}

\begin{prop}
\label{f!commoub}
Suppose  $T ' = f ^{-1} (T)$.
\begin{enumerate}[(a)]
\item 
\label{f!commoub-item1}
Let  $\E ^{(\bullet)} \in \underrightarrow{LD} ^{\mathrm{b}} _{\Q ,\mathrm{qc}}
( \smash{\widehat{\D}} _{\fP /\fS } ^{(\bullet)})$.
We have the canonical isomorphism 
$$f ^{ !(\bullet)} \circ  oub _T \circ (\hdag T) (\E^{(\bullet)}) 
\riso 
oub _{T'} \circ (\hdag T ') \circ f ^{ !(\bullet)}   (\E^{(\bullet)}),$$ 
which we can simply write 
$f ^{ !(\bullet)}  \circ (\hdag T) (\E^{(\bullet)}) 
\riso 
 (\hdag T ') \circ f ^{ !(\bullet)}   (\E^{(\bullet)}) $. 
\item Let  $\E ^{(\bullet)} \in \underrightarrow{LD} ^{\mathrm{b}} _{\Q ,\mathrm{qc}}
( \smash{\widetilde{\D}} _{\fP /\fS } ^{(\bullet)} (T))$.
We have the canonical isomorphism
$$oub _{T'} \circ f ^{ !(\bullet)} _{T} (\E^{(\bullet)} )\riso f ^{!(\bullet)} \circ oub _T (\E^{(\bullet)}).$$
Hence, it is harmless to write by abuse of notation
$f ^{!(\bullet)} $ instead of $f ^{!(\bullet)} _T$. 

\end{enumerate}

\end{prop}

\begin{proof}
Using \ref{f!T'Totimes},\ref{lem-f!B}, for any $\E ^{(\bullet)} \in \underrightarrow{LD} ^{\mathrm{b}} _{\Q ,\mathrm{qc}}
( \smash{\widehat{\D}} _{\fP /\fS } ^{(\bullet)})$, 
we get the isomorphism
\begin{equation}
\notag
f ^{ !(\bullet)} \circ  oub _T \circ (\hdag T) (\E^{(\bullet)}) 
=
f ^{!(\bullet)}
   \left (
   \widetilde{\B} ^{(\bullet)}  _{\fP} ( T) 
   \smash{\widehat{\otimes}} ^\L
_{\O ^{(\bullet)}  _{\fP} } 
   \E ^{(\bullet)}
\right ) 
\riso
   \widetilde{\B} ^{(\bullet)}  _{\fP'} ( T') 
   \smash{\widehat{\otimes}} ^\L
_{\O ^{(\bullet)}  _{\fP'} } 
f ^{!(\bullet)} 
(  \E ^{(\bullet)})
=
oub _{T'} \circ (\hdag T ') \circ f ^{ !(\bullet)}   (\E^{(\bullet)}) .
\end{equation}
By using 
\ref{oub-pl-fid-iso1} and \ref{oub-div-opcoh}.\ref{oub-div-opcohb)},
we check the second part from the first one.
\end{proof}

\begin{rem}
\label{rem-f!TornoT}
With notation \ref{f!commoub}, 
using the remark 
\ref{rema-fct-qcoh2coh}
we check that the functors 
$\mathrm{Coh} _{T} (f ^{!(\bullet)} _T) $
and
$\mathrm{Coh}  (f ^{!(\bullet)} ) $ 
are isomorphic over
$D ^\mathrm{b} _\mathrm{coh} ( \smash{\D} ^\dag _{\fP,\Q} ) 
\cap D ^\mathrm{b} _\mathrm{coh} ( \smash{\D} ^\dag _{\fP } (\hdag T) _{\Q} )$.
Since 
we have the canonical isomorphisms of functors 
$\mathrm{Coh} _{T} (f ^{!(\bullet)} _T)  \riso f  _{T} ^!$ 
and
$\mathrm{Coh} _{T} (f ^{!(\bullet)} )  \riso f ^!$ (\ref{coh-Qcoh}), 
then it is harmless to write 
$f ^!$ instead of  $f   _{T} ^!$. 

\end{rem}

\subsection{Spencer resolutions, finite Tor-dimension}

\begin{ntn}
\label{ntn-fXY/T}
We keep notation \ref{subsection3.1}.
When 
$f$ has locally $p$-bases,
for all $m  \leq m'$,
we set 
$\smash{\D} _{\fP ' /\fP   } ^{(m,m')} (T'):=
\widehat{\B} ^{(m')} _{\fP '} ( T')  \otimes _{\cO _{\fP '}} \smash{\D} _{\fP ' /\fP   } ^{(m)}$,
$\widetilde{\B} ^{(m)} _{\fP '} ( T'):= \widehat{\B} ^{(n_m)} _{\fP '} ( T')$,
$\smash{\widetilde{\D}} _{\fP ' /\fP   } ^{(m)} (T'):=
\smash{\widehat{\D}} _{\fP ' /\fP   } ^{(m,n_m)} (T')$.
\end{ntn}

\begin{lem}
With notation \ref{ntn-fXY/T}, 
we suppose 
$f$ has locally $p$-bases. 
Let $m' \geq m\geq 0$ be two integers. 

\begin{enumerate}[a)]
\item We have the canonical isomorphism
\begin{equation}
\label{Btildem2m'}
\widetilde{\B} ^{(m')}_{\fP '} (T') _\Q
\riso 
\widetilde{\D} ^{(m')} _{\fP '  /\fP   } (T') _\Q
\otimes _{\widetilde{\D} ^{(m)} _{\fP '  /\fP   } (T') _\Q}
\widetilde{\B} ^{(m)}_{\fP '} (T') _\Q.
\end{equation}

\item We have the canonical isomorphism
\begin{equation}
\label{DXY2DXThatcor}
\smash{\widetilde{\D}} _{\fP ' \to \fP   /\fS } ^{(m')} (T') _\Q
\riso 
\widetilde{\D} ^{(m')} _{\fP '  /\fS } (T') _\Q
\otimes _{\widetilde{\D} ^{(m)} _{\fP '  /\fS } (T') _\Q}
\smash{\widetilde{\D}} _{\fP ' \to \fP   /\fS } ^{(m)} (T') _\Q.
\end{equation}

\end{enumerate}

\end{lem}

\begin{proof}
a) For all $m''\geq m ' \geq m$, since the morphism
$$\widehat{\B} ^{(m'')}_{\fP '} (T')
\to 
\smash{\widehat{\D}} _{\fP ' /\fP   } ^{(m',m'')} ( T')
\otimes _{\D _{\fP ' /\fP   } ^{(m',m'')} ( T')}
\widehat{\B} ^{(m'')}_{\fP '} (T')$$ 
is an isomorphism (see \cite[4.4.9]{Be1}),
since 
$\D _{\fP ' /\fP   } ^{(m,m'')} ( T') _\Q
=
\D _{\fP ' /\fP   } ^{(m',m'')} ( T') _\Q$,
then we get 
$$\widehat{\B} ^{(m'')}_{\fP '} (T') _\Q
\riso
\smash{\widehat{\D}} _{\fP ' /\fP   } ^{(m',m'')} ( T') _\Q
\otimes _{\smash{\widehat{\D}} _{\fP ' /\fP   } ^{(m,m'' )} ( T') _\Q}
\widehat{\B} ^{(m'')}_{\fP '} (T')_\Q.$$

Since 
$\widehat{\B} ^{(m'')} _{\fP '} (T')
\to 
\smash{\widehat{\D}} _{\fP ' /\fP   } ^{(m,m'')} ( T')
\otimes _{\smash{\widehat{\D}} _{\fP ' /\fP   } ^{(m,m')} ( T')}
\widehat{\B} ^{(m')}_{\fP '} (T')$ is an isomorphism (use the arguments of the proof of \cite[4.4.8]{Be1}),
then 
$\widehat{\B} ^{(m'')} _{\fP '} (T') _\Q
\to 
\smash{\widehat{\D}} _{\fP ' /\fP   } ^{(m,m'')} ( T') _\Q
\otimes _{\smash{\widehat{\D}} _{\fP ' /\fP   } ^{(m,m')} ( T')_\Q}
\widehat{\B} ^{(m')}_{\fP '} (T')_\Q$.
Hence, we are done.

b) From \ref{DXY2DXT}, we get by projective limit the isomorphism
\begin{equation}
\label{DXY2DXThat}
\smash{\widetilde{\D}} _{\fP ' /\fS } ^{(m)} (T')
\otimes _{\smash{\widetilde{\D}} _{\fP ' /\fP   } ^{(m)} (T')}
\widetilde{\B} ^{(m)} _{\fP '} ( T')
\riso 
\smash{\widetilde{\D}} _{\fP ' \to \fP   /\fS } ^{(m)} (T').
\end{equation}
Hence, using \ref{Btildem2m'}, we get the isomorphism \ref{DXY2DXThatcor}.
\end{proof}

\begin{empt}
\label{Spencer-XY}
With notation \ref{ntn-fXY/T}, 
we suppose  $f$ has locally $p$-bases. 
Taking the inverse limit of the exact sequences of the form 
\ref{Spencer-sequence-B-prehat2}, 
 we get the exact sequence
\begin{equation}
\label{Spencer-sequence-B-hat}
0 
\to 
\widetilde{\D} ^{(0)} _{\fP '  /\fS } (T') 
\otimes _{\cO _{\fP '}}
\wedge ^{d} \mathcal{T} _{\fP '  /\fP    }
\cdots
\underset{\delta}{\longrightarrow}
\widetilde{\D} ^{(0)} _{\fP '  /\fS } (T') 
\otimes _{\cO _{\fP '}}
\mathcal{T} _{\fP '  /\fP    }
\underset{\delta}{\longrightarrow}
\widetilde{\D} ^{(0)} _{\fP '  /\fS } (T') 
\to 
\smash{\widetilde{\D}} _{\fP ' \to \fP   /\fS } ^{(0)} (T')
\to 0.
\end{equation}

Since 
$\widetilde{\D} ^{(0)} _{\fP '  /\fS } (T') 
\to 
\widetilde{\D} ^{(m)} _{\fP '  /\fS } (T') _\Q$
is flat, from \ref{DXY2DXThatcor} we get by extension 
the exact sequence
\begin{equation}
\label{Spencer-sequence-B-hat(m)}
0 
\to 
\widetilde{\D} ^{(m)} _{\fP '  /\fS } (T') _\Q
\otimes _{\cO _{\fP '}}
\wedge ^{d} \mathcal{T} _{\fP '  /\fP    }
\cdots
\underset{\delta}{\longrightarrow}
\widetilde{\D} ^{(m)} _{\fP '  /\fS } (T') _\Q
\otimes _{\cO _{\fP '}}
\mathcal{T} _{\fP '  /\fP    }
\underset{\delta}{\longrightarrow}
\widetilde{\D} ^{(m)} _{\fP '  /\fS } (T') _\Q
\to 
\smash{\widetilde{\D}} _{\fP ' \to \fP   /\fS } ^{(m)} (T') _\Q
\to 0.
\end{equation}
We denote by $\smash{\widetilde{\mathrm{Sp}}} ^{(m)} _{\fP '  /\fP   } (T') $ the complex
$\widetilde{\D} ^{(m)} _{\fP '  /\fS } (T') 
\otimes _{\cO _{\fP '}}
\wedge ^{d} \mathcal{T} _{\fP '  /\fP    }
\cdots
\underset{\delta}{\longrightarrow}
\widetilde{\D} ^{(m)} _{\fP '  /\fS } (T') 
\otimes _{\cO _{\fP '}}
\mathcal{T} _{\fP '  /\fP    }
\underset{\delta}{\longrightarrow}
\widetilde{\D} ^{(m)} _{\fP '  /\fS } (T') $.

Via the equivalence of categories
$D _{\Q,\mathrm{coh}} ^{\mathrm{b}}
(\widetilde{\D} ^{(m)} _{\fP '  /\fS } (T') )
\cong
D _{\mathrm{coh}} ^{\mathrm{b}}
(\widetilde{\D} ^{(m)} _{\fP '  /\fS } (T') _\Q)$, we get the isomorphism in 
$D _{\Q,\mathrm{coh}} ^{\mathrm{b}}
(\widetilde{\D} ^{(m)} _{\fP '  /\fS } (T') )$:
\begin{equation}
\label{Spencer-XY-iso(m)}
\smash{\widetilde{\mathrm{Sp}}} ^{(m)} _{\fP '  /\fP   } (T') 
\riso 
\smash{\widetilde{\D}} _{\fP ' \to \fP   /\fS } ^{(m)} (T') .
\end{equation}
We get in 
$\underrightarrow{LD} ^{\mathrm{b}} _{\Q,\mathrm{coh}} (\widetilde{\D} ^{(\bullet)} _{\fP '  /\fS } (T') )$
the isomorphism
\begin{equation}
\label{Spencer-XY-isobullet}
\smash{\widetilde{\mathrm{Sp}}} ^{(\bullet)} _{\fP '  /\fP   } (T') 
\riso 
\smash{\widetilde{\D}} _{\fP ' \to \fP   /\fS } ^{(\bullet)} (T').
\end{equation}

\end{empt}

\begin{lem}
\label{DX2Yflatclosedimm}
With notation \ref{ntn-fXY/T}, 
we suppose $f$ is a closed immersion. 
\begin{enumerate}[(a)]
\item The left 
$\smash{\widetilde{\D}} _{P ' _i /S _i } ^{(m)} (T') $-module
$\smash{\widetilde{\D}} _{P ' _i \to P _i /S _i } ^{(m)} (T') $
is flat.
\item The left 
$\smash{\widetilde{\D}} _{\fP ' /\fS } ^{(m)} (T') $-module
$\smash{\widetilde{\D}} _{\fP ' \to \fP   /\fS } ^{(m)} (T') $
is flat.
\end{enumerate}

\end{lem}

\begin{proof}
Since 
$\smash{\widetilde{\D}} _{\fP ' \to \fP   /\fS } ^{(m)} (T') 
\riso 
\underleftarrow{\lim} _i 
\smash{\widetilde{\D}} _{P ' _i \to P _i /S _i } ^{(m)} (T') $, 
using \cite[3.2.4]{Be1}, we reduce to check the first assertion. 
Since this is local, we can suppose $\fP   /\fS $
has some $p$-basis  $t _1,\dots, t _d$ such that
the image of 
$t _1,\dots, t _{d'}$, via
$f ^{-1} \cO _{\fP   } \to \cO _{\fP ' }$,
is a $p$-basis of 
$\fP ' /\fS $.
In that case, 
$\smash{\widetilde{\D}} _{P ' _i \to P _i /S _i } ^{(m)} (T') $
is a free 
$\smash{\widetilde{\D}} _{P ' _i /P _i } ^{(m)} (T') $-module.
\end{proof}

\begin{prop}
Let $\alpha \colon \fP ' \to \bbD ^r _{\fS}$ and
$\beta \colon \fP   \to \bbD ^s _{\fS}$ be two objects of 
$\scr{C} _{\fS}$ (see notation \ref{dfn-CfS}).
Let $f \colon \alpha \to \beta$ be a morphism of $\scr{C} _{\fS}$.
\begin{enumerate}[(a)]
\item We have
$\smash{\widetilde{\D}} _{\fP ' \to \fP   /\fS } ^{(m)} (T') 
\in D _{\Q,\mathrm{tdf}} ^{\mathrm{b}}
(\widetilde{\D} ^{(m)} _{\fP '  /\fS } (T') )
$.
\item We have
$\smash{\widetilde{\D}} _{\fP ' \to \fP   /\fS } ^{(\bullet)} (T') 
\in 
\underrightarrow{LD} ^{\mathrm{b}} _{\Q,\mathrm{tdf}} (\widetilde{\D} ^{(\bullet)} _{\fP '  /\fS } (T') )$.

\end{enumerate}
\end{prop}

\begin{proof}
Following \ref{Spencer-XY},
we have the canonical isomorphism
$$\smash{\widetilde{\mathrm{Sp}}} ^{(m)} _{\fP ' \times _{\scr{C} _{\fS}}\fP    /\fP   } (T') 
\riso 
\smash{\widetilde{\D}} _{\fP ' \times _{\scr{C} _{\fS}} \fP   \to \fP   /\fS } ^{(m)} (T') $$
in 
$D _{\Q,\mathrm{tdf}} ^{\mathrm{b}}
(\widetilde{\D} ^{(m)} _{\fP '  \times _{\scr{C} _{\fS}} \fP   /\fS } (T') )$.
Let $u \colon 
\fP ' 
\to 
\fP '  \times _{\scr{C} _{\fS}} \fP   $ be the graph of $f$.
Using \ref{DX2Yflatclosedimm}, 
since 
$\widetilde{f} ^* \widetilde{\D} ^{(m)} _{\fP '  \times _{\scr{C} _{\fS}} \fP   /\fS } (T') )
\riso
\widetilde{\D} ^{(m)} _{\fP ' \to \fP   /\fS } (T') )$, then
$\L \widetilde{f} ^* \smash{\widetilde{\mathrm{Sp}}} ^{(m)} _{\fP ' \times _{\scr{C} _{\fS}} \fP    /\fP   } (T') $
is bounded complex with Tor amplitude in $[0,\delta _{P'/S}]$.
\end{proof}

\subsection{Projection formula : commutation of pushforwards with localization functors outside of a divisor}

Let  $f \colon \fP ^{\prime } \to \fP $ be a morphism of formal schemes 
locally of formal finite type
and having  locally finite $p$-bases over $\fS $,
$T$ and $T'$ be some divisors of respectively $P$ and $P'$ such that 
$f ( P '\setminus T' ) \subset P \setminus T$.
We finish this subsection by giving some applications of the projection formula.
\begin{prop}
\label{surcoh2.1.4}
Let $\E ^{ (\bullet)}
\in  
\smash{\underrightarrow{LD}} ^\mathrm{b} _{\Q, \mathrm{qc}}
(\overset{^\mathrm{l}}{} \smash{\widetilde{\D}} _{\fP ^{ }/\fS } ^{(\bullet)} (T))$,
and 
$\E ^{\prime (\bullet)}
\in  \smash{\underrightarrow{LD}} ^\mathrm{b} _{\Q, \mathrm{qc}}
(\overset{^\mathrm{l}}{} \smash{\widetilde{\D}} _{\fP ^{\prime }/\fS } ^{(\bullet)}(T'))$. 
We have the following isomorphism of $\smash{\underrightarrow{LD}} ^\mathrm{b} _{\Q, \mathrm{qc}}
(\overset{^\mathrm{l}}{} \smash{\widetilde{\D}} _{\fP ^{ }/\fS } ^{(\bullet)} (T))$ 
\begin{equation}
\label{surcoh2.1.4-iso}
f _{T,T',+} ^{(\bullet)} ( \E ^{\prime (\bullet)} )
\smash{\widehat{\otimes}}^\L 
_{\widetilde{\B} ^{(\bullet)}  _{\fP} ( T) } 
\E ^{ (\bullet)} [\delta _{P'/P}]
\riso 
f _{T,T',+} ^{(\bullet)} 
\left ( 
\E ^{\prime (\bullet)} 
\smash{\widehat{\otimes}}^\L 
_{\widetilde{\B} ^{(\bullet)}  _{\fP'} ( T') } 
f _{T',T} ^{!(\bullet)} (\E ^{ (\bullet)} )
\right ) .
\end{equation}

\end{prop}

\begin{proof}
Similarly to \cite[3.3.7]{caro-6operations}, we check this is a consequence of \ref{prop-u+otimes-u!}.
\end{proof}

\begin{cor}
\label{surcoh2.1.4-cor1}
Let $\E ^{ (\bullet)}
\in  \smash{\underrightarrow{LD}} ^\mathrm{b} _{\Q, \mathrm{qc}}
(\overset{^\mathrm{l}}{} \smash{\widetilde{\D}} _{\fP ^{ }/\fS } ^{(\bullet)} (T))$.
We have the isomorphism 
\begin{equation}
\label{surcoh2.1.4-isocor1}
f _{T,T',+} ^{(\bullet)} \left (\widetilde{\B} ^{(\bullet)}  _{\fP'} ( T') \right )
\smash{\widehat{\otimes}}^\L 
_{\widetilde{\B} ^{(\bullet)}  _{\fP} ( T) } 
\E ^{ (\bullet)} [\delta _{P'/P}]
\riso 
f _{T,T',+} ^{(\bullet)} \circ  f _{T',T} ^{!(\bullet)} (\E ^{ (\bullet)} ) .
\end{equation}
\end{cor}

\begin{proof}
We apply \ref{surcoh2.1.4}  
to the case where
$\E ^{\prime (\bullet)}
=
\widetilde{\B} ^{(\bullet)}  _{\fP'} ( T') $.
\end{proof}

\begin{cor}
\label{surcoh2.1.4-cor}
Suppose  $T ' = f ^{-1} (T)$.
Let $\E ^{\prime (\bullet)}
\in  \smash{\underrightarrow{LD}} ^\mathrm{b} _{\Q, \mathrm{qc}}
(\overset{^\mathrm{l}}{} \smash{\widetilde{\D}} _{\fP ^{\prime }/\fS } ^{(\bullet)})$.
We have the isomorphism of
$\smash{\underrightarrow{LD}} ^\mathrm{b} _{\Q, \mathrm{qc}}
(\overset{^\mathrm{l}}{} \smash{\widetilde{\D}} _{\fP ^{ }/\fS } ^{(\bullet)})$: 
$$f _{T,T'+} ^{(\bullet)}  \circ (\hdag T') ( \E ^{\prime (\bullet)} )
\riso 
(\hdag T ) \circ f _{+} ^{(\bullet)} ( \E ^{\prime (\bullet)} ).$$

\end{cor}

\begin{proof}
Use \ref{surcoh2.1.4}
and \ref{lem-f!B}, we get the isomorphism
\begin{equation}
\label{surcoh2.1.4-cor-iso}
f _{+} ^{(\bullet)} ( \E ^{\prime (\bullet)} )
\smash{\widehat{\otimes}}^\L _{\O ^{(\bullet)}  _{\fP} }  
\widetilde{\B} ^{(\bullet)}  _{\fP} ( T) 
\riso 
f _{+} ^{(\bullet)} ( \E ^{\prime (\bullet)} 
\smash{\widehat{\otimes}}^\L 
_{\O ^{(\bullet)}  _{\fP'} }  
\widetilde{\B} ^{(\bullet)}  _{\fP'} ( T') ).
\end{equation}
We conclude using \ref{oub-div-opcoh}.\ref{oub-div-opcoha)}.
\end{proof}

\begin{rem}
Using \ref{oub-div-opcoh}, the isomorphism 
of \ref{surcoh2.1.4-cor} could be written
$f _{+} ^{(\bullet)}  \circ (\hdag T') ( \E ^{\prime (\bullet)} )
\riso 
(\hdag T ) \circ f _{+} ^{(\bullet)} ( \E ^{\prime (\bullet)} ).$

\end{rem}

\subsection{On the stability of the coherence}
Let  $f \colon \fP ^{\prime } \to \fP $ be a morphism of 
formal  $\fS$-schemes locally of formal finite type
and having locally finite $p$-bases  over $\fS$,
$T$ and $T'$ be some divisors of respectively $P$ and $P'$ such that 
$f ( P '\setminus T' ) \subset P \setminus T$.

\begin{lem}
\label{stab-coh-f^!pre}
Suppose $f _i\colon P _i ^{\prime} \to P  _i$ 
has locally finite $p$-bases.
For any $\E \in D ^{-} _{\mathrm{coh}} (\widetilde{\D} ^{(m)} _{P _i /S _i }(T))$, 
we have 
$f _{i,T',T}^{(m)!} (\E )\in D ^{-} _{\mathrm{coh}} (\widetilde{\D} ^{(m)} _{P _i ^{\prime}/S _i }(T'))$. 
\end{lem}

\begin{proof}
Since this is local in $P _i ^{\prime}$, using locally free resolution, 
we reduce to the case $\E= \widetilde{\D} ^{(m)} _{P _i /S _i }(T)$. 
Then we compute in local coordinates 
that 
the canonical morphism
$\D ^{(m)} _{P _i ^{\prime}/S _i } \to f _i ^{*} \D ^{(m)} _{P _i /S _i }$ is surjective
whose kernel has the usual description in local 
coordinates.\end{proof}
\begin{prop}
\label{stab-coh-f^!}
Suppose $f$ is is flat and has locally finite $p$-bases.
\begin{enumerate}[(a)]
\item For  $\E \in D ^{\mathrm{b}} _{\mathrm{coh}} (\widetilde{\D} ^{(m)} _{\fP /\fS }(T))$, 
we have 
$f _{T',T}^{(m)!} (\E )\in D ^{\mathrm{b}} _{\mathrm{coh}} (\widetilde{\D} ^{(m)} _{\fP ^{\prime}/\fS }(T'))$. 

\item For  $\E \in D ^{\mathrm{b}} _{\mathrm{coh}} (\widetilde{\D} ^{(m)} _{\fP /\fS }(T) _\Q)$, 
we have
\begin{equation}
\notag
\widetilde{\D} ^{(m+1)} _{\fP ^{\prime}/\fS }(T') _\Q
\otimes ^{\L} _{\widetilde{\D} ^{(m)} _{\fP ^{\prime}/\fS }(T')_\Q}
 f _{T',T}^{(m)!} (\E )
\riso 
f _{T',T}^{(m+1)!} (
\widetilde{\D} ^{(m+1)} _{\fP /\fS }(T)_\Q
\otimes ^{\L} _{\widetilde{\D} ^{(m)} _{\fP /\fS }(T)_\Q}
\E ).
\end{equation}

\item The functor $ f _{T',T}^{!(\bullet)}$
sends 
$\underrightarrow{LD} ^{\mathrm{b}}  _{\Q, \mathrm{coh}} (\smash{\widetilde{\D}} _{\fP /\fS } ^{(\bullet)} (T))$
to 
$\underrightarrow{LD} ^{\mathrm{b}}  _{\Q, \mathrm{coh}} (\smash{\widetilde{\D}} _{\fP ^{\prime }/\fS } ^{(\bullet)} (T'))$.
\item 
For $\E \in D ^{\mathrm{b}} _{\mathrm{coh}} ( \D ^{\dag} _{\fP  } (\hdag T ) _{\Q})$,
we have
$f  ^{ !} _{T' , T} (\E ) \in
D ^{\mathrm{b}} _{\mathrm{coh}} ( \D ^{\dag} _{\fP ^{\prime }} (\hdag T' ) _{\Q})$.
\end{enumerate}

\end{prop}

\begin{proof}
The first part is a consequence of \ref{stab-coh-f^!pre}. 
We check the second part similarly to  \cite[3.4.6]{Beintro2},
i.e. this is an easy consequence of the Spencer resolution 
\ref{Spencer-XY-iso(m)}.
The third and forth parts are a consequence of the previous ones.\end{proof}

\begin{lem}
\label{stab-coh-f_+}
Suppose $f$ is proper, and 
$T ' = f ^{-1}(T)$.
\begin{enumerate}[(a)]
\item The functor $f ^{(m)} _{i,T+}$ sends 
$D ^{-} _{\mathrm{coh}} (\widetilde{\D} ^{(m)} _{P _i ^{\prime}/S _i }(T'))$
to 
$D ^{-} _{\mathrm{coh}} (\widetilde{\D} ^{(m)} _{P _i /S _i }(T))$. 
\item For $\E ' \in D ^{-} _{\mathrm{coh}} (\widetilde{\D} ^{(m)} _{P _i ^{\prime}/S _i }(T'))$,
we have the canonical isomorphism
\begin{equation}
\notag
\widetilde{\D} ^{(m+1)} _{P _i /S _i }(T')
\otimes ^{\L} _{\widetilde{\D} ^{(m)} _{P _i /S _i }(T')}
f ^{(m)} _{i,T+}
(\E ' )
\riso 
f ^{(m+1)} _{i,T+}
\left (\widetilde{\D} ^{(m+1)} _{P _i ^{\prime}/S _i }(T')
\otimes ^{\L} _{\widetilde{\D} ^{(m)} _{P _i ^{\prime}/S _i }(T')}
\E ' \right ).
\end{equation}

\end{enumerate}
\end{lem}

\begin{proof}
This is proved similary to \cite[3.4.3]{caro-6operations}.
\end{proof}

\begin{prop}
Suppose $f$ is proper, and 
$T ' = f ^{-1}(T)$.
\begin{enumerate}[(a)]
\item For  $\E '\in D ^{\mathrm{b}} _{\mathrm{coh}} (\widetilde{\D} ^{(m)} _{\fP ^{\prime}/\fS }(T'))$, 
we have 
$f _{T,+}^{(m)} (\E' )\in D ^{\mathrm{b}} _{\mathrm{coh}} (\widetilde{\D} ^{(m)} _{\fP /\fS }(T))$. 

\item For  $\E' \in D ^{\mathrm{b}} _{\mathrm{coh}} (\widetilde{\D} ^{(m)} _{\fP ^{\prime}/\fS }(T'))$, 
we have
\begin{equation}
\notag
\widetilde{\D} ^{(m+1)} _{\fP /\fS }(T)
\otimes ^{\L} _{\widetilde{\D} ^{(m)} _{\fP /\fS }(T)}
  f _{T,+}^{(m)} (\E )
\riso 
  f _{T,+}^{(m+1)} (
\widetilde{\D} ^{(m+1)} _{\fP ^{\prime}/\fS }(T')
\otimes ^{\L} _{\widetilde{\D} ^{(m)} _{\fP ^{\prime}/\fS }(T')}
\E ).
\end{equation}

\item The functor $ f _{T,+} ^{(\bullet)}$
sends 
$\underrightarrow{LD} ^{\mathrm{b}}  _{\Q, \mathrm{coh}} (\smash{\widetilde{\D}} _{\fP ^{\prime }/\fS } ^{(\bullet)} (T'))$
to 
$\underrightarrow{LD} ^{\mathrm{b}}  _{\Q, \mathrm{coh}} (\smash{\widetilde{\D}} _{\fP /\fS } ^{(\bullet)} (T))$.

\item For $\E '\in D ^{\mathrm{b}} _{\mathrm{coh}} ( \D ^{\dag} _{\fP ^{\prime }} (\hdag T' ) _{\Q})$, 
we have
$f  _{T, +}( \E ')\in D ^{\mathrm{b}} _{\mathrm{coh}} ( \D ^{\dag} _{\fP  } (\hdag T ) _{\Q})$.
\end{enumerate}

\end{prop}

\begin{proof}
This is a consequence of \ref{stab-coh-f_+}.
\end{proof}

\subsection{Base change and their commutation with cohomological operations}
\label{subsect-comm-bc}

Let $\alpha \colon \V \to \W$ be a morphism of local algebras
 such that 
 $\cV$ and $\W$ are  complete discrete valued ring  of mixed characteristic $(0,p)$ with perfect residue fields.
We set $\fS: =\Spf \cV$ and $\fT: =\Spf \cW$.

Let $f \colon \X ^{\prime } \to \X $
 be a morphism of formal $\fS$-scheme locally of formal finite type and having locally finite $p$-bases over $\fS$.
Let $\fY$ be a formal $\fT$-scheme locally of formal finite type and having locally finite $p$-bases over $\fT$
such that there exists a relatively perfect morphism of formal $\fT$-schemes of the form
$\vartheta \colon \fY \to \fX \times _{\fS} \fT$ (this is the product in the category of $p$-adic formal schemes,
i.e. $\fX \times _{\fS} \fT := \underset{i}{\underrightarrow{\lim}} X  _i \times _{S _i} T _i$). Let 
$\varpi \colon \fY \to \fX$ be the composition of $\vartheta$ with the projection 
$\mathrm{pr}\colon \fX \times _{\fS} \fT \to \fX$.
Let 
$\fY ^{ \prime } : = \X ^{ \prime } \times _{\fX }\fY  $, 
and $\varpi '\colon  \fY ^{ \prime } \to \X ^{ \prime }$,
$g  \colon  \fY ^{ \prime } \to \fY $
be the projections.
We suppose $\varpi$ (and hence $\varpi '$) is flat. 

\begin{ex}
\label{ex-ft-Dn}
Set 
$\cV [[\underline{t}]]:=\cV [[t_1, \dots, t _n]]$
and
$\cW [[\underline{t}]]:=\cW [[t_1, \dots, t _n]]$.
The canonical morphism
$\cV [[\underline{t}]] \to \cW [[\underline{t}]]$ is flat. 
(Indeed, let $\pi$ be a uniformiser of $\cV$ and $k := \cV /\pi \cV$ be its residue field. 
Modulo $\pi$, we get the morphism
$k [[\underline{t}]] \to (\cW/\pi \cW)  [[\underline{t}]]$.
Since $(\cW/\pi \cW)  [[\underline{t}]]$ is noetherian, $t$-adically complete and without $t$-torsion, 
then
$k [[\underline{t}]] \to \cW/\pi \cW  [[\underline{t}]]$ is flat
(use \cite[3.5,Theorem 1]{bourbaki}).
Hence, since 
$\cW [[\underline{t}]]$ is noetherian, $\pi$-adically complete and without $\pi$-torsion, 
this yields that
$\cV [[\underline{t}]] \to \cW [[\underline{t}]]$ is flat
(again, use \cite[3.5,Theorem 1]{bourbaki}).)

The canonical morphism 
\begin{equation}
\label{ex-ft-Dn-relperf1}
\Spf \cW [[\underline{t}]] \to \Spf (\cV [[\underline{t}]]) \times _\fS \fT
\end{equation}
is relatively perfect.
(Indeed, 
since relative perfect morphisms are stable under base change, 
since the canonical morphism
$\Spf (\cV [[\underline{t}]])  \to \widehat{\bbA} ^n _{\fS}$
is relatively perfect, then so is 
$\Spf (\cV [[\underline{t}]]) \times _{\fS} \fT 
\to 
\widehat{\bbA} ^n _{\fT}$.
Since $\Spf \cW [[\underline{t}]] \to \widehat{\bbA} ^n _{\fT}$
is also relatively perfect, then so is the morphism
$\Spf \cW [[\underline{t}]] \to \Spf (\cV [[\underline{t}]]) \times _{\fS} \fT.$)

Let $\fX$ be a formal  $\Spf \cV [[\underline{t}]]$-scheme of finite type, having locally finite $p$-bases  over $\Spf \V$.
By applying the functor
$\fX \times _{\Spf \cV [[\underline{t}]]} -$ to \ref{ex-ft-Dn-relperf1}, 
we get the relatively perfect morphism
$\fY : = \fX \times _{\Spf \cV [[\underline{t}]]} \Spf \cW [[\underline{t}]]
\to 
\fX\times _{\fS} \fT$.
Since 
$\fX$ has locally finite $p$-bases  over $\fS$,
then 
$\fX\times _{\fS} \fT$ 
has locally finite $p$-bases  over $\fT$.
Hence, so is $\fY/\fT$.
Let 
$\varpi \colon \fY \to \fX$
be the canonical projection
$\fX \times _{\Spf \cV [[\underline{t}]]} \Spf \cW [[\underline{t}]]
\to 
\fX$.
Then $\varpi$ is flat.
(Indeed, this is a consequence of \ref{lem-flat-ff0} :
$\cV [[\underline{t}]] \to \cW [[\underline{t}]]$ is a flat morphism of noetherian rings,
$\fX$ is of finite type over $\Spf \cV [[\underline{t}]]$,
$\fX$ and 
$\fY$ have no $p$-torsion following 
\ref{f0formétale-fforméta0-cor}).
\end{ex}

\begin{empt}
\label{rel-pdim-bc}
Since $\vartheta$ is relatively perfect, 
then with notation \ref{ntn-deltaX/T} we have $ \delta ^{\fT} _{\fY}
= 
\vartheta ^{-1} \delta ^{\fT} _{\fX \times _{\fS} \fT}$.
We check easily the formula 
$\delta ^{\fS} _{\fX} \circ \mathrm{pr}
= 
\delta ^{\fT} _{\fX \times _{\fS} \fT}$.
Hence 
$ \delta ^{\fS} _{\fX} \circ \varpi
=
\delta ^{\fT} _{\fY}$.
Similarly, 
$   \delta ^{\fS} _{\fX '} \circ \varpi ^{\prime}
=
\delta ^{\fT} _{\fY '}$.
This yields 
\begin{equation}
\label{rel-pdim-bc-=}
\delta ^{\fS} _{\fX '/\fX}
\circ \varpi ^{\prime}
=
\delta ^{\fT} _{\fY '/\fY}.
\end{equation}

\end{empt}

\begin{empt}
\label{dft-bc}
For any integer $i $, the canonical morphism 
$\smash{\D} _{ Y _i/T _i} ^{(m)} 
\to 
\varpi _i ^*
\smash{\D} _{X _i/S _i} ^{(m)} 
=
\smash{\D} _{Y _i \to X _i/T _i \to S _i} ^{(m)} 
$
is an isomorphism. 
Let  $\E _i$ be a  left $\smash{\D} _{ X _i/S _i} ^{(m)} $-module
and 
$\M _i $ be a right 
$\smash{\D} _{ X _i/S _i} ^{(m)} $-module.
The canonical homomorphisms
\begin{gather}
\notag
\varpi _i ^* (\E _i) 
=
\smash{\O} _{ Y _i} 
\otimes  _{\varpi ^{-1}\smash{\O} _{ X _i}}
\varpi ^{-1} \E _i
\to 
\smash{\D} _{ Y _i/T _i} ^{(m)} 
\otimes  _{\varpi ^{-1}\smash{\D} _{ X _i/S _i} ^{(m)} }
\varpi ^{-1} \E _i
\leftarrow
\smash{\D} _{ Y _i/T _i} ^{(m)} 
\otimes ^\L _{\varpi ^{-1}\smash{\D} _{ X _i/S _i} ^{(m)} }
\varpi ^{-1} \E _i
=:
\varpi _i ^{* (m)} (\E _i),
\\
\notag
\varpi _i ^* (\M _i) 
=
\varpi ^{-1} \M _i
\otimes  _{\varpi ^{-1}\smash{\O} _{ X _i}}
\smash{\O} _{ Y _i} 
\to 
\varpi ^{-1} \M _i
\otimes   _{\varpi ^{-1}\smash{\D} _{ X _i/S _i} ^{(m)} }
\smash{\D} _{ Y _i/T _i} ^{(m)} 
\leftarrow
\varpi ^{-1} \M _i
\otimes ^\L  _{\varpi ^{-1}\smash{\D} _{ X _i/S _i} ^{(m)} }
\smash{\D} _{ Y _i/T _i} ^{(m)} 
=:
\varpi _i ^{* (m)} (\cM _i) 
\end{gather}
are isomorphisms.
By computing in local coordinates, 
we can check that the canonical isomorphism of $\O _{Y _i}$-modules
\begin{equation}
\label{formula-bc-otimes}
\varpi ^{*(m)} _i ( \M _i \otimes _{\O _{X _i}} \E _i)
\riso 
\varpi ^{*(m)} _i ( \M _i ) \otimes _{\O _{Y _i}} \varpi ^{*(m)} _i (  \E _i)
\end{equation}
is $\smash{\D} _{ Y _i/T _i} ^{(m)} $-linear.
We have similar results by adding some primes in the notations.

We get the functor
$\varpi ^{* (\bullet)}
\colon 
\smash{\underrightarrow{LD}} ^\mathrm{b} _{\Q, \mathrm{qc}}
(\overset{^\mathrm{l}}{} \smash{\widehat{\D}} _{\X /\fS } ^{(\bullet)})
\to 
\smash{\underrightarrow{LD}} ^\mathrm{b} _{\Q, \mathrm{qc}}
(\overset{^\mathrm{l}}{} \smash{\widehat{\D}} _{\Y/\fT } ^{(\bullet)})$, 
given by 
$\varpi ^{* (\bullet)} 
(\E ^{(\bullet)})
=
(\varpi ^{* (m)} 
(\E ^{(m)})) _m$,
for any 
$\E ^{(\bullet)} 
\in 
\smash{\underrightarrow{LD}} ^\mathrm{b} _{\Q, \mathrm{qc}}
(\overset{^\mathrm{l}}{} \smash{\widehat{\D}} _{\X /\fS } ^{(\bullet)})$.
This is the base change functor. 
This functor preserves the coherence, i.e. this induces the functor 
$\varpi ^{* (\bullet)}
\colon 
\smash{\underrightarrow{LD}} ^\mathrm{b} _{\Q, \mathrm{coh}}
(\overset{^\mathrm{l}}{} \smash{\widehat{\D}} _{\X /\fS } ^{(\bullet)})
\to 
\smash{\underrightarrow{LD}} ^\mathrm{b} _{\Q, \mathrm{coh}}
(\overset{^\mathrm{l}}{} \smash{\widehat{\D}} _{\Y/\fT } ^{(\bullet)})$.
Via the equivalence of categories 
\ref{eq-catLDBer-LD-D}, 
this yields the functor
$\varpi ^{*}
\colon 
D ^{\mathrm{b}} _{\mathrm{coh}}( \smash{\D} ^\dag _{\fX /\fS, \bbQ} )
\to 
D ^{\mathrm{b}} _{\mathrm{coh}}( \smash{\D} ^\dag _{\fY /\fT ,\bbQ} )
$
which is canonically isomorphic to 
$\E \mapsto \smash{\D} ^\dag _{\fY /\fS ,\Q} 
\otimes _{\varpi ^{-1}\smash{\D} ^\dag _{\fX /\fS ,\Q} }
\varpi ^{-1} \E $.

\end{empt}

\begin{prop}
\label{prop-iso-chgtbasef!}
Let 
$\E ^{ (\bullet)}
\in  \smash{\underrightarrow{LD}} ^\mathrm{b} _{\Q, \mathrm{qc}}
(\overset{^\mathrm{l}}{} \smash{\widehat{\D}} _{\X /\fS } ^{(\bullet)})$. 
There exists a canonical isomorphism in 
$\smash{\underrightarrow{LD}} ^\mathrm{b} _{\Q, \mathrm{qc}}
(\overset{^\mathrm{l}}{} \smash{\widehat{\D}} _{\fY ^{\prime } /\fT } ^{(\bullet)})$ of the form:
\begin{equation}
\label{iso-chgtbase1}
\varpi  ^{\prime * (\bullet) } \circ f ^{! (\bullet)}_{/\fS } (\E ^{(\bullet)})
\riso
g  ^{! (\bullet)}_{/\fT}  \circ 
\varpi ^{* (\bullet)} 
 (\E ^{ (\bullet)}). 
\end{equation}
\end{prop}

\begin{proof}
We reduce to check that, for any integer $i \in \N$, we have in $D^\mathrm{b} _{\mathrm{qc}}
(\smash{\D} _{ Y '_i/T _i} ^{(m)})$ the canonical isomorphism 
\begin{equation}
\label{iso-chgtbase1-bis}
\varpi _i ^{\prime * (m)}  \circ f ^{! (m)} _{i /S _i}( \E _i^{(m)})
\riso 
g ^{! (m)} _{i /T _i}  \circ \varpi _i ^{ * (m) } (\E _i^{(m)}).
\end{equation}
Since
$f ^{! (m)} _{i /S _i} \riso
\bbL f ^{*} [  \delta ^{\fS} _{\fX '/\fX}]$
and 
$g ^{! (m)} _{i /S _i} \riso
\bbL g ^{*} [\delta ^{\fT} _{\fY '/\fY}]$, 
via \ref{rel-pdim-bc-=}, we conclude by transitivity of the inverse image.\end{proof}

\begin{prop}
\label{bc-com-tensprod}
Tensor products (see \ref{def-otimes-coh1qc}) 
commutes with base change, i.e. we have the canonical isomorphism in 
$\smash{\underrightarrow{LD}}  ^{\mathrm{b}}  _{\Q, \mathrm{qc}}
(\overset{^\mathrm{?}}{} \smash{\widetilde{\D}} _{\fX /\fS } ^{(\bullet)})$
for any $\cM ^{ (\bullet)}
\in \smash{\underrightarrow{LD}}  ^{\mathrm{b}} _{\Q, \mathrm{qc}}
(\overset{^\mathrm{?}}{} \smash{\widetilde{\D}} _{\fX /\fS } ^{(\bullet)} )
$
and 
$\cE ^{ (\bullet)} 
\in 
\smash{\underrightarrow{LD}}  ^{\mathrm{b}} _{\Q, \mathrm{qc}}
(\overset{^\mathrm{l}}{} \smash{\widetilde{\D}} _{\fX /\fS } ^{(\bullet)} )$: 
\begin{equation}
\label{bc-com-tensprod-iso}
\varpi  ^{* (\bullet) }
(
\cM ^{ (\bullet)}
 \smash{\widehat{\otimes}}
^\L _{\cO ^{(\bullet)}  _{\fX}}
\cE ^{ (\bullet)} 
)
\riso 
\varpi  ^{* (\bullet) }
(\cM ^{ (\bullet)})
 \smash{\widehat{\otimes}}
^\L _{\cO ^{(\bullet)}  _{\fY}}
\varpi  ^{* (\bullet) } (\cE ^{ (\bullet)} )
.
\end{equation}

\end{prop}

\begin{proof}
We reduce to check that, for any integer $i \in \N$, we have in $D^\mathrm{b} _{\mathrm{qc}}
(\smash{\D} _{ Y _i/T _i} ^{(m)})$ the canonical isomorphism 
\begin{equation}
\label{bc-com-tensprod-iso1}
\varpi _i  ^{* (m) }
(
\cM _i ^{ (m)}
\otimes 
^\L _{\cO   _{X _i}}
\cE _i ^{ (m)} 
)
\riso 
\varpi   _i ^{* (m) }
(\cM _i ^{ (m)})
\otimes 
^\L _{\cO  _{Y _i}}
\varpi  _i ^{* (m) } (\cE _i ^{ (m)} )
.
\end{equation}
which is obvious.
\end{proof}

\begin{theo}
\label{theo-iso-chgtbase1}
Let 
$\E ^{\prime (\bullet)}
\in  \smash{\underrightarrow{LD}} ^\mathrm{b} _{\Q, \mathrm{qc}}
(\overset{^\mathrm{l}}{} \smash{\widehat{\D}} _{\X ^{\prime }/\fS } ^{(\bullet)})$. 
There exists a canonical isomorphism in 
$\smash{\underrightarrow{LD}} ^\mathrm{b} _{\Q, \mathrm{qc}}
(\overset{^\mathrm{l}}{} \smash{\widehat{\D}} _{\fY ^{ } /\fT } ^{(\bullet)})$ of the form:
\begin{equation}
\label{iso-chgtbase1bis}
\varpi ^{* (\bullet)} \circ f ^{(\bullet)}_{/\fS  +} (\E ^{\prime (\bullet)})
\riso
g  ^{(\bullet)}_{/\fT  +}  \circ \varpi  ^{\prime * (\bullet) } (\E ^{\prime (\bullet)}). 
\end{equation}
\end{theo}

\begin{proof}
We reduce to check that, for any integer $i \in \N$, we have in $D^\mathrm{b} _{\mathrm{qc}}
(\smash{\D} _{ Y _i/T _i} ^{(m)})$ the canonical isomorphism 
\begin{equation}
\label{pre-adj-morph-gen-Rham3pre}
\varpi _i ^{* (m)}  \circ f ^{ (m)} _{i /S _i+}( \E _i^{\prime (m)})
\riso 
g ^{ (m)} _{i /T _i+}  \circ \varpi _i ^{\prime  * (m) } (\E _i^{\prime (m)}).
\end{equation}

a) By computing in local coordinates, 
we can check the canonical isomorphisms of $\O _{Y _i}$-modules
(resp. of $\O _{Y '_i}$-modules)
$\varpi ^* _i (\omega _{X _i /S _i})
\riso 
\omega _{Y _i /T _i}$
(resp. 
$\varpi ^{\prime *} _i (\omega _{X '_i /S _i})
\riso 
\omega _{Y' _i /T _i}$)
is $\smash{\D} _{ Y _i/T _i} ^{(m)} $-linear
(resp. $\smash{\D} _{ Y ' _i/T _i} ^{(m)} $-linear).
Hence, we get the isomorphism of right $\smash{\D} _{ Y '_i/T _i} ^{(m)} $-modules:
\begin{gather}
\notag
\varpi ^{\prime *} _i 
(\smash{\D} _{ X _i \leftarrow X ^{ \prime } _i/S _i} ^{(m)})
=
\varpi ^{\prime *} _i 
(\omega _{X '_i /S _i}
 \otimes _{\O _{ X ^{\prime }_i}}  
   f ^{ *} _{i \mathrm{r}} ( \smash{\D} _{ X _i/S _i} ^{(m)}  
 \otimes _{\O _{ X  _i}} \omega _{X ^{ }_i/S _i} ^{-1}) )
 \\
 \notag 
 \riso
 \omega _{Y ^{ \prime }_i /T _i} \otimes _{\O _{ Y ^{\prime }_i}}   
 \varpi ^{\prime *} _i 
   f ^{ *} _{i \mathrm{r}} ( \smash{\D} _{ X _i/S _i} ^{(m)}  
 \otimes _{\O _{ X  _i}} \omega _{X ^{ }_i/S _i} ^{-1}) 
 \riso
 \omega _{Y ^{ \prime }_i /T _i} \otimes _{\O _{ Y ^{\prime }_i}}   
g ^{ *} _{i \mathrm{r}}  \varpi ^{*} _{i, \mathrm{r}} 
    ( \smash{\D} _{ X _i/S _i} ^{(m)}  
 \otimes _{\O _{ X  _i}} \omega _{X ^{ }_i/S _i} ^{-1})  
 \\ 
 \label{5114}
\riso 
\omega _{Y ^{ \prime }_i /T _i} \otimes _{\O _{ Y ^{\prime }_i}}  
   g ^{ *} _{i \mathrm{r}} ( \smash{\D} _{ Y _i/T _i} ^{(m)}  
 \otimes _{\O _{ Y _i}} \omega _{Y _i/T _i} ^{-1}) )
 =
 \smash{\D} _{Y _i \leftarrow Y ^{ \prime } _i/T _i} ^{(m)}.
\end{gather}
In fact, 
since 
$\smash{\D} _{ X _i \leftarrow X ^{ \prime } _i/S _i} ^{(m)}$ 
is a 
$(f ^{-1}\smash{\D} _{ X _i/S _i} ^{(m)},
\smash{\D} _{ X ^{ \prime } _i/S _i} ^{(m)})$-bimodule, 
by functoriality, the homomorphisms of \ref{5114}
are homomorphisms of 
$(  ( \varpi \circ g) ^{-1}\smash{\D} _{ X _i/S _i} ^{(m)},
\smash{\D} _{ Y ^{ \prime } _i/S _i} ^{(m)})$-bimodules.
Let $\mathcal{P} _i^{\prime (m)}$ be a left resolution of 
$\E _i^{\prime (m)}$
by flat left $ \smash{\D} _{ X ^{ \prime } _i/S _i} ^{(m)}$-module.
Hence, we get the morphism of
$D (  ( \varpi \circ g) ^{-1} \smash{\D} _{ X _i/S _i} ^{(m)})$:
\begin{gather}
\label{iso1-theo-iso-chgtbase1pre}
\varpi _i ^{ \prime -1}
 (\smash{\D} _{ X _i \leftarrow X ^{ \prime } _i/S _i} ^{(m)}
\otimes _{ \smash{\D} _{ X ^{\prime } _i/S _i} ^{(m)}} \mathcal{P} _i^{\prime (m)})
\to
 \varpi _i ^{ \prime *}
 (\smash{\D} _{ X _i \leftarrow X ^{ \prime } _i/S _i} ^{(m)})
 \otimes _{ \smash{\D} _{ Y ^{\prime } _i/T _i} ^{(m)}} \varpi _i ^{ \prime *} (\mathcal{P} _i^{\prime (m)})
 \underset{\ref{5114}}{\riso}
     \smash{\D} _{Y _i \leftarrow Y ^{ \prime } _i/T _i} ^{(m)}
 \otimes _{ \smash{\D} _{ Y ^{\prime } _i/T _i} ^{(m)}} \varpi _i ^{ \prime *} (\mathcal{P} _i^{\prime (m)}).
\end{gather}

b) We have the adjunction morphism 
$\mathrm{adj} \colon \varpi ^{-1} \R f _* \to \R g _* \varpi ^{\prime -1}$ of functors
$D ( f ^{-1} \smash{\D} _{ X _i/S _i} ^{(m)} ) 
\to 
D ( \varpi ^{-1} \smash{\D} _{ X _i/S _i} ^{(m)})$.
Hence, we get the morphism of
$D ( \varpi ^{-1} \smash{\D} _{ X _i/S _i} ^{(m)})$:
\begin{gather}
\notag
 \varpi _i ^{-1} \circ \R f _{*}
 (\smash{\D} _{ X _i \leftarrow X ^{ \prime } _i/S _i} ^{(m)}
 \otimes  _{ \smash{\D} _{ X ^{\prime } _i/S _i} ^{(m)}} \mathcal{P} _i^{\prime (m)})
\overset{\mathrm{adj}}{\longrightarrow}
\R g _{*}  \circ \varpi _i ^{ \prime -1}
 (\smash{\D} _{ X _i \leftarrow X ^{ \prime } _i/S _i} ^{(m)}
\otimes _{ \smash{\D} _{ X ^{\prime } _i/S _i} ^{(m)}} \mathcal{P} _i^{\prime (m)})
\to 
\\
%
%
\label{iso1-theo-iso-chgtbase1}
 \underset{\ref{iso1-theo-iso-chgtbase1pre}}{\longrightarrow}
\R g _{*} (    \smash{\D} _{Y _i \leftarrow Y ^{ \prime } _i/T _i} ^{(m)}
 \otimes _{ \smash{\D} _{ Y ^{\prime } _i/T _i} ^{(m)}} \varpi _i ^{ \prime *} (\mathcal{P} _i^{\prime (m)})).
\end{gather}
This yields the morphism
of 
$D (\smash{\D} _{Y _i/T _i} ^{(m)})$
\begin{equation}
\label{iso1-theo-iso-chgtbase2}
\varpi _i ^{ *} \circ f ^{ (m)} _{i /S _i +}( \E _i^{\prime (m)})
=
 \varpi _i ^{ *} \circ \R f _{*}
 (\smash{\D} _{ X _i \leftarrow X ^{ \prime } _i/S _i} ^{(m)}
 \otimes ^{\L}_{ \smash{\D} _{ X ^{\prime } _i/S _i} ^{(m)}} \E _i^{\prime (m)})
\to 
 \R g_{*} (    \smash{\D} _{Y _i \leftarrow Y ^{ \prime } _i/T _i} ^{(m)}
 \otimes ^{\L}_{ \smash{\D} _{ Y ^{\prime } _i/T _i} ^{(m)}} \varpi _i ^{ \prime *} (\E _i^{\prime (m)}))
=
g ^{ (m)} _{i /T _i+}  \circ \varpi _i ^{\prime *} (\E _i^{\prime (m)}).
\end{equation}
It remains to check that this morphism is an isomorphism.
Since the functors
$\varpi _i ^{ *} \circ f ^{ (m)} _{i /S _i +}$
and
$g ^{ (m)} _{i /T _i+}  \circ \varpi _i ^{\prime *}$ 
are way out left, using (the way out left version of)
Proposition \cite[I.7.1.(iv)]{HaRD},
we reduce to the case where
$\E _i^{\prime (m)}$ is of the form 
$\smash{\D} _{ X ^{\prime } _i/S _i} ^{(m)}  
\otimes _{\O _{X ^{\prime} _i}}
\FF ' _i$,
where 
$\FF ' _i$ is a quasi-coherent 
$\O _{X ' _i}$-module.
The morphism \ref{iso1-theo-iso-chgtbase1} is canonically isomorphism to the composite of the top arrow of the following diagram:
\small
\begin{equation}
\xymatrix @C=0,3cm {
{ \varpi _i ^{-1} \R f _{*}
 (\smash{\D} _{ X _i \leftarrow X ^{ \prime } _i/S _i} ^{(m)}
 \otimes ^{\L}_{ \smash{\D} _{ X ^{\prime } _i/S _i} ^{(m)}} \E _i^{\prime (m)})} 
\ar[r] ^-{\mathrm{adj}}
\ar[d] ^-{\sim}
& 
{\R g _{*} \varpi _i ^{ \prime -1}
(\smash{\D} _{ X _i \leftarrow X ^{ \prime } _i/S _i} ^{(m)}
 \otimes ^{\L}_{ \smash{\D} _{ X ^{\prime } _i/S _i} ^{(m)}} \E _i^{\prime (m)})} 
\ar[r] ^-{}
\ar[d] ^-{\sim}
& 
{\R g _{*} (    \smash{\D} _{Y _i \leftarrow Y ^{ \prime } _i/T _i} ^{(m)}
 \otimes ^{\L}_{ \smash{\D} _{ Y ^{\prime } _i/T _i} ^{(m)}} \varpi _i ^{ \prime *} (\E _i^{\prime (m)}))
} 
\ar[d] ^-{\sim}
\\ 
{ \varpi _i ^{-1} \circ \R f _{*}
 (\smash{\D} _{ X _i \leftarrow X ^{ \prime } _i/S _i} ^{(m)}
\otimes _{\O _{X ^{\prime} _i}}
\FF ' _i)
} 
\ar[r] ^-{\mathrm{adj}}
\ar[d] ^-{}
& 
{\R g _{*}  \circ \varpi _i ^{ \prime -1}
 (\smash{\D} _{ X _i \leftarrow X ^{ \prime } _i/S _i} ^{(m)}
\otimes _{\O _{X ^{\prime} _i}}
\FF ' _i)
 } 
\ar[r] ^-{}
\ar[d] ^-{}
& 
{\R g _{*} (    \smash{\D} _{Y _i \leftarrow Y ^{ \prime } _i/T _i} ^{(m)}
\otimes _{\O _{U ^{\prime} _i}}
\varpi _i ^{ \prime *} (\FF ' _i))} 
\ar@{=}[d] ^-{}
\\ 
{ \varpi _i ^{*} \circ \R f _{*}
 (\smash{\D} _{ X _i \leftarrow X ^{ \prime } _i/S _i} ^{(m)}
\otimes _{\O _{X ^{\prime} _i}}
\FF ' _i)
} 
\ar[r] ^-{\mathrm{adj}} _-{\sim}
& 
{\R g _{*}  \circ \varpi _i ^{ \prime *}
 (\smash{\D} _{ X _i \leftarrow X ^{ \prime } _i/S _i} ^{(m)}
\otimes _{\O _{X ^{\prime} _i}}
\FF ' _i)
 } 
\ar[r] ^-{\sim}
& 
{\R g _{*} (    \smash{\D} _{Y _i \leftarrow Y ^{ \prime } _i/T _i} ^{(m)}
\otimes _{\O _{U ^{\prime} _i}}
\varpi _i ^{ \prime *} (\FF ' _i)),} 
}
\end{equation}
\normalsize
where the adjunction isomorphism of the bottom line is the one in the categories of $\O$-modules.
This yields the commutative diagram:
\begin{equation}
\xymatrix{
{\varpi _i ^{ *} \circ f ^{ (m)} _{i /S _i +}( \E _i^{\prime (m)})} 
\ar[rr] ^-{}
\ar[d] ^-{\sim}
&& 
{ g ^{ (m)} _{i /T _i+}  \circ \varpi _i ^{\prime *} (\E _i^{\prime (m)})} 
\ar[d] ^-{\sim}
\\ 
{\varpi _i ^{*} \circ \R f _{*}
 (\smash{\D} _{ X _i \leftarrow X ^{ \prime } _i/S _i} ^{(m)}
\otimes _{\O _{X ^{\prime} _i}}
\FF ' _i)} 
\ar[r] ^-{\mathrm{adj}} _-{\sim}
& 
{\R g_{*}  \circ \varpi _i ^{ \prime *}
 (\smash{\D} _{ X _i \leftarrow X ^{ \prime } _i/S _i} ^{(m)}
\otimes _{\O _{X ^{\prime} _i}}
\FF ' _i)
 } 
\ar[r] ^-{\sim}
&
{\R g_{*} (    \smash{\D} _{Y _i \leftarrow Y ^{ \prime } _i/T _i} ^{(m)}
\otimes _{\O _{U ^{\prime} _i}}
\varpi _i ^{ \prime *} (\FF ' _i)).} 
}
\end{equation}
Recall (e.g. see Lemma 30.5.2 of the stack project), since $\varpi _i$ is flat, 
then 
we get the isomorphism
$ \varpi _i ^{ *} \circ \R f _{*} \riso  \R g_{*}  \circ \varpi _i ^{ \prime *}$, 
where 
$\R f _{*}  
\colon 
D ( \O _{X ^{\prime} _i}) 
\to 
D ( \O _{X  _i}) $, 
$\varpi _i ^{ *} 
\colon 
D (\O _{X  _i}) 
\to 
D ( \O _{Y _i} )$,
$\varpi _i ^{ \prime *}
\colon 
D (\O _{X ^{\prime} _i}) 
\to 
D ( \O _{U ^{\prime} _i} )$,
$\R f _{*}  '
\colon 
D ( \O _{U ^{\prime} _i}) 
\to 
D (\O _{Y _i}) $.
\end{proof}

\begin{prop}
\label{com-dual-bc}
Let $\E \in D ^{\mathrm{b}} _{\mathrm{coh}}( \smash{\D} ^\dag _{\fX /\fS ,\Q} )$.
We have the canonical isomorphism
\begin{equation}
\label{com-dual-bc-iso}
\varpi ^{*}(\DD _{\fX /\fS }  (\E))
\riso 
\DD _{\fY /\fT}
(\varpi ^{*} (\E)).
\end{equation}

\end{prop}

\begin{proof}
Since $D ^{\mathrm{b}} _{\mathrm{coh}}( \smash{\D} ^\dag _{\fX /\fS ,\Q} )= D ^{\mathrm{b}} _{\mathrm{parf}}( \smash{\D} ^\dag _{\fX /\fS ,\Q} )$, 
then we have the last canonical isomorphism
\begin{gather}
\varpi ^{*}(\DD _{\fX /\fS }  (\E))
=
\smash{\D} ^\dag _{\fY /\fS ,\Q} 
\otimes _{\varpi ^{-1}\smash{\D} ^\dag _{\fX /\fS ,\Q} }
\varpi ^{-1} 
\R \mathcal{H} om _{ \smash{\D} ^\dag _{\fX /\fS ,\Q} }
(\E,
 \smash{\D} ^\dag _{\fX /\fS ,\Q} 
\otimes  _{\O _{\fX}} \omega _{\fX /\fS } ^{-1}) [\delta _{\fX} ^{\fS}]
\\
\riso 
\smash{\D} ^\dag _{\fY /\fS ,\Q} 
\otimes _{\varpi ^{-1}\smash{\D} ^\dag _{\fX /\fS ,\Q} }
\R \mathcal{H} om _{ \varpi ^{-1}\smash{\D} ^\dag _{\fX /\fS ,\Q} }
(\varpi ^{-1}\E,
\varpi ^{-1}( \smash{\D} ^\dag _{\fX /\fS ,\Q} 
\otimes  _{\O _{\fX}} \omega _{\fX ^{ }/\fS } ^{-1})) [\delta _{\fX} ^{\fS}]
\\
\underset{\cite[2.1.12]{caro-frobdualrel}}{\riso} 
\R \mathcal{H} om _{ \varpi ^{-1}\smash{\D} ^\dag _{\fX /\fS ,\Q} }
(\varpi ^{-1}\E,
\varpi ^{*}( \smash{\D} ^\dag _{\fX /\fS ,\Q} 
\otimes  _{\O _{\fX}} \omega _{\fX ^{ }/\fS } ^{-1})) [\delta _{\fX} ^{\fS}].
\end{gather}
Using \ref{formula-bc-otimes}, 
since 
$\varpi ^{*(m)}( \smash{\widehat{\D}} _{\X/\fS } ^{(m)} )
\riso
\smash{\widehat{\D}} _{\Y/\fT } ^{(m)}$, then 
we get the isomorphism of 
left $\smash{\widehat{\D}} _{\Y/\fT } ^{(\bullet)}$-bimodules
$$\varpi ^{*(\bullet)}( \smash{\widehat{\D}} _{\X/\fS } ^{(\bullet)}
\otimes  _{\O _{\fX}} \omega _{\fX ^{ }/\fS } ^{-1})
\riso 
\smash{\widehat{\D}} _{\Y/\fT } ^{(\bullet)}
\otimes  _{\O _{\fY}} \omega _{\fY /\fT } ^{-1}.$$
Via the equivalence of categories 
\ref{limeqcat}, 
this is translated by the canonical isomorphism
of left $\smash{\D} ^\dag _{\fY /\fT ,\Q} $-bimodules 
$$\varpi ^{*}( \smash{\D} ^\dag _{\fX /\fS ,\Q} 
\otimes  _{\O _{\fX}} \omega _{\fX ^{ }/\fS } ^{-1})
\riso 
\smash{\D} ^\dag _{\fY /\fT ,\Q} 
\otimes  _{\O _{\fY}} \omega _{\fY /\fT } ^{-1}).$$
Using \ref{rel-pdim-bc}, 
this yields the first isomorphism
\begin{gather}
\R \mathcal{H} om _{ \varpi ^{-1}\smash{\D} ^\dag _{\fX /\fS ,\Q} }
(\varpi ^{-1}\E,
\varpi ^{*}( \smash{\D} ^\dag _{\fX /\fS ,\Q} 
\otimes  _{\O _{\fX}} \omega _{\fX ^{ }/\fS } ^{-1})) [\delta _{\fX} ^{\fS}]
\\
\riso 
\R \mathcal{H} om _{ \varpi ^{-1}\smash{\D} ^\dag _{\fX /\fS ,\Q} }
(\varpi ^{-1}\E,
\smash{\D} ^\dag _{\fY /\fT ,\Q} 
\otimes  _{\O _{\fY}} \omega _{\fY ^{ }/\fT } ^{-1})) [\delta _{\fY} ^{\fT}]
\\
\riso 
\R \mathcal{H} om _{\smash{\D} ^\dag _{\fY /\fT ,\Q} }
(\varpi ^{*}\E,
\smash{\D} ^\dag _{\fY /\fT ,\Q} 
\otimes  _{\O _{\fY}} \omega _{\fY /\fT } ^{-1})) [\delta _{\fY} ^{\fT}]
=
\DD _{\fY /\fT}
(\varpi ^{*} (\E)).
\end{gather}
\end{proof}

\section{Closed immersions : pushforwards and extraordinary pullbacks}

\subsection{The fundamental isomorphism for schemes}
\label{ntnclosedimmer-1}
Put $S:= S _i$.
Let $u \colon Z \hookrightarrow X$ be a closed immersion of 
$S$-schemes 
locally of formal finite type
and having  locally finite $p$-bases  over $S$.
Let $\I$ be the ideal defining $u$. 
The subsection \cite[4.1]{caro-6operations} is still valid without new argument in our context. 
For the reader, let us collect its results and local description below.

\begin{empt}
[Some notation with local coordinates]
\label{locdesc-climm}
Suppose
$X$ is affine and there exist
$t _{r +1},\dots , t _{d}  \in \Gamma (X,\I)$
generating 
$I:=\Gamma (X,\I)$,
$t _{1},\dots , t _{r}\in \Gamma ( X,\O _{X})$
such that
$t _{1},\dots ,t  _{d}$ form a finite $p$-basis of $X$ over $S$,
$\overline{t} _{1},\dots ,\overline{t} _{r}$ 
form a finite $p$-basis of $Z$ over $S$,
and 
$\overline{t} _{r +1},\dots ,\overline{t} _{d}$ is a basis of $\I /\I ^2$,
where $\overline{t} _{1},\dots , \overline{t} _{r} \in\Gamma ( Z ,\O _{Z})$
(resp. $\overline{t} _{r +1},\dots ,\overline{t} _{d}\in\Gamma ( X ,\cI /\cI ^2)$)
are the images of 
$t _{1},\dots , t _{r}$
(resp. $t _{r+1},\dots , t _{d}$)
via 
$\Gamma ( X,\O _{X})
\to
\Gamma ( Z ,\O _{Z})$
(resp. 
$\Gamma ( X,\cI)
\to 
\Gamma ( X ,\cI /\cI ^2)$).

We denote by 
$\tau _i := 1 \otimes t _i -t _i \otimes 1$, 
$\overline{\tau} _j := 1 \otimes \overline{t} _j -\overline{t} _j \otimes 1$, 
for any $i= 1,\dots, d$, $j= 1,\dots, r$.
The sheaf of $\O _X$-algebras
$\cP ^n  _{X/S,(m)}$ is 
a free $\O _X$-module with the basis 
$\{ \underline{\tau} ^{\{\underline{k}\} _{(m)}}\ | \ \underline{k}\in \N ^d \text{ such that } | \underline{k}| \leq n\} $,
and
$\cP ^n  _{Z/S,(m)}$ is 
a free $\O _Z$-module with the basis 
$\{ \underline{\overline{\tau}} ^{\{\underline{l}\} _{(m)}}\ | \ \underline{l}\in \N ^r \text{ such that } | \underline{l}| \leq n\} $.
We denote by 
$\{ 
\underline{\partial} ^{<\underline{k}> _{(m)}}
\ | \ \underline{k}\in \N ^d,\ | \underline{k} | \leq n
\}$
the corresponding dual basis  of 
$\D ^{(m)} _{X/S,n}$ 
and
by 
$\{ \underline{\partial} ^{<\underline{l}> _{(m)}}\  
| \ \underline{l}\in \N ^r,\ | \underline{l} | \leq n\} $
the corresponding dual basis of $\D ^{(m)} _{Z/S,n}$ (we hope the similar notation is not too confusing).
The sheaf
$\D ^{(m)} _{X/S}$ is 
a free $\O _X$-module with the basis 
$\{ 
\underline{\partial} ^{<\underline{k}> _{(m)}}
\ | \ \underline{k}\in \N ^d
\}$,
and
$\D ^{(m)} _{Z/S}$ is 
a free $\O _Z$-module with the basis 
$\{ \underline{\partial} ^{<\underline{l}> _{(m)}}\ | \ \underline{l}\in \N ^r\} $.

a) We compute the canonical homomorphism
$u ^* \cP ^n  _{X/S,(m)}
\to \cP ^n  _{Z/S,(m)}$
sends 
$\underline{\tau} ^{\{(\underline{l}, \underline{h})\} _{(m)}}$
where 
$\underline{l} \in \N ^r$
and 
$\underline{h} \in \N ^{d-r}$
to
$\underline{\overline{\tau}} ^{\{\underline{l}\} _{(m)}}$
if $\underline{h} = (0,\dots,0)$
and to $0$ otherwise.

b) We denote by 
$\theta \colon 
\cD ^{(m)} _{Z/S}
\to 
\cD ^{(m)} _{Z \to X/S}$
the canonical homomorphism of left $\cD ^{(m)} _{Z/S}$-modules
(which is built by duality from the 
canonical homomorphisms $u ^* \cP ^n  _{X/S,(m)}
\to \cP ^n  _{Z/S,(m)}$).
For any $P \in D ^{(m)} _{X/S}$, we denote by 
$\overline{P}$ its image via
the canonical morphism of left $D ^{(m)} _{X/S}$-modules
$D ^{(m)} _{X/S}
\to 
D ^{(m)} _{X/S} / I D ^{(m)} _{X/S}
= 
D ^{(m)} _{Z \to X/S}$.
We set 
$\underline{\xi} ^{<\underline{k}> _{(m)}}:= 
\overline{\underline{\partial} ^{<\underline{k}> _{(m)}}}$. 
By duality from a), we compute
$\theta (\underline{\partial} ^{<\underline{l}> _{(m)}})
=
\underline{\xi} ^{<(\underline{l}, \underline{0})> _{(m)}}$,
for any 
$\underline{l}\in \N ^r$.

\end{empt}

\begin{empt}
\label{locdesc-climm2}
Suppose we are in the local situation of \ref{locdesc-climm}. 
We denote by 
$\cD ^{(m)} _{X,Z,\underline{t}/S}$ the subring of 
$\cD ^{(m)} _{X/S}$ which is a 
 free $\O _X$-module with the basis 
$\{ \underline{\partial} ^{<(\underline{l}, \underline{0})> _{(m)}}\ | \ \underline{l}\in \N ^r\} $, 
where 
$\underline{0}:=(0,\dots, 0) \in 
\N ^{d-r}$.
If there is no ambiguity concerning the local coordinates (resp. and $S$),
we might simply denote $\cD ^{(m)} _{X,Z,\underline{t}/S}$ by 
$\cD ^{(m)} _{X,Z/S}$
(resp. $\cD ^{(m)} _{X,Z}$).
\begin{enumerate}[(a)]
\item 
We have the following factorization
\begin{equation}
\label{diag-DZ2X-t-comm}
\xymatrix{
{ \cD ^{(m)} _{X,Z,\underline{t}/S}/\cI\cD ^{(m)} _{X,Z,\underline{t}/S}} 
\ar@{^{(}->}[r] ^-{}
& 
{\cD ^{(m)} _{X/S} / \cI \cD ^{(m)} _{X/S}} 
\\ 
{u _* \cD ^{(m)} _{Z/S}} 
\ar[u] ^-{\sim} _-{\theta}
\ar[r] ^-{\theta}
& 
{u _* \cD ^{(m)} _{Z \to X/S},} 
\ar@{=}[u] ^-{}
}
\end{equation}
where both horizontal morphisms  are canonical.
Both rings 
$u _* \cD ^{(m)} _{Z/S}$ and 
$\cD ^{(m)} _{X,Z,\underline{t}/S}$ are 
$\O _X$-rings (i.e. they are rings endowed with a structural 
homomorphism of rings 
$\O _X \to \cD ^{(m)} _{X,Z,\underline{t}/S}$
and 
$\O _X \to u _* \cD ^{(m)} _{Z/S}$).
The vertical arrow
$u _* \cD ^{(m)} _{Z/S}
\to
 \cD ^{(m)} _{X,Z,\underline{t}/S}/\cI\cD ^{(m)} _{X,Z,\underline{t}/S}$
is an isomorphism of $\cO _X$-rings.
Moreover, the sheaf
$\cD ^{(m)} _{X,Z,\underline{t}/S}/\cI\cD ^{(m)} _{X,Z,\underline{t}/S} $
is a $(u _*\cD ^{(m)} _{Z/S} , \D ^{(m)} _{X,Z,\underline{t}/S})$-subbimodule
of $\cD ^{(m)} _{X/S} / \cI \cD ^{(m)} _{X/S}$
and 
the vertical arrow
$u _* \cD ^{(m)} _{Z/S}
\to
 \cD ^{(m)} _{X,Z,\underline{t}/S}/\cI\cD ^{(m)} _{X,Z,\underline{t}/S}$
is also 
an isomorphism of 
left $u _* \cD ^{(m)} _{Z/S}$-modules.

\item \label{g}
We get the composite 
\begin{equation}
\label{morp-rho}
\rho\colon \cD ^{(m)} _{X,Z,\underline{t}/S}
\to
\cD ^{(m)} _{X,Z,\underline{t}/S}
/\I \cD ^{(m)} _{X,Z,\underline{t}/S}
\underset{\theta}{\liso}  
u _* \cD ^{(m)} _{Z/S}
\end{equation}
is an homomorphism of $\O _X$-rings.

\item Since 
$\cD ^{(m)} _{X/S} $
is a free left $\cD ^{(m)} _{X,Z,\underline{t}/S}$-module
with the basis
$\{ \underline{\partial} ^{<(\underline{0},\underline{h})> _{(m)}}\ | \ \underline{h}\in \N ^{d-r}\} $,
where 
$\underline{0}:=(0,\dots, 0) \in 
\N ^{r}$,
then from the commutativity of \ref{diag-DZ2X-t-comm},
we get that 
$\cD ^{(m)} _{Z \to X/S}$ is a free 
left 
$\cD ^{(m)} _{Z/S}$-module 
with the basis
$\{ \underline{\xi} ^{<(\underline{0},\underline{h})> _{(m)}}\ | \ \underline{h}\in \N ^{d-r}\} $,
where 
$\underline{0}:=(0,\dots, 0) \in 
\N ^{r}$.

\item We have the transposition automorphism
${} ^t \colon 
D ^{(m)} _{X/S}
\to 
D ^{(m)} _{X/S}$
given by 
$P 
= 
\sum _{\underline{k} \in \N ^d}
a _{\underline{k}} \underline{\partial} ^{<\underline{k}> _{(m)}}
\mapsto
{} ^t P:= 
\sum _{\underline{k} \in \N ^d}
(-1) ^{| \underline{k}|}\underline{\partial} ^{<\underline{k}> _{(m)}}
a _{\underline{k}} $.
Beware that this transposition depends on the choice of the local coordinates $t _1,\dots, t _d$.
This transposition automorphism induces 
${} ^t \colon 
D ^{(m)} _{X,Z,\underline{t}/S}
\to 
D ^{(m)} _{X,Z,\underline{t}/S}$
such that
${} ^t ( ID ^{(m)} _{X,Z,\underline{t}/S})
=ID ^{(m)} _{X,Z,\underline{t}/S}$.
This yields the automorphism
${} ^t \colon 
D ^{(m)} _{X,Z,\underline{t}/S}/I D ^{(m)} _{X,Z,\underline{t}/S}
\to 
D ^{(m)} _{X,Z,\underline{t}/S}
/
ID ^{(m)} _{X,Z,\underline{t}/S}$.
On the other hand, 
via the local coordinates $\overline{t} _{1},\dots ,\overline{t} _{r}$ 
of $Z$ over $S$,
we get
the transposition automorphism
${} ^t \colon 
D ^{(m)} _{Z/S}
\to 
D ^{(m)} _{Z/S}$
given by 
$Q 
= 
\sum _{\underline{k} \in \N ^r}
b _{\underline{k}} \underline{\partial} ^{<\underline{k}> _{(m)}}
\mapsto
{} ^t Q:= 
\sum _{\underline{k} \in \N ^r}
(-1) ^{| \underline{k}|}\underline{\partial} ^{<\underline{k}> _{(m)}}
b _{\underline{k}} $.
We compute the following diagram 
\begin{equation}
\label{transp-closedimmersion}
\xymatrix{
{D ^{(m)} _{X,Z,\underline{t}/S}/I D ^{(m)} _{X,Z,\underline{t}/S}} 
\ar[r] _-{\sim} ^-{{}^t}
& 
{D ^{(m)} _{X,Z,\underline{t}/S}/I D ^{(m)} _{X,Z,\underline{t}/S}} 
\\ 
{D ^{(m)} _{Z/S}} 
\ar[u] ^-{\sim} _-{\theta}
\ar[r] _-{\sim} ^-{{}^t}
& 
{D ^{(m)} _{Z/S}} 
\ar[u] ^-{\sim} _-{\theta}
}
\end{equation}
is commutative.
\end{enumerate}

\end{empt}

\begin{empt}
\label{rightD-moduleflat}

We denote $\overline{u}\colon (Z, \O _Z) \to (X, u _* \O _Z)$ the morphism of ringed spaces induced by $u$.
We remark that $\overline{u}$ is flat and 
that 
$\overline{u} ^* = u ^{-1} \colon 
D ^{+}  ( u _* \O _Z )
\to 
D ^{+}  (  \O _Z )$.
Recall that for any $\M \in D ^{+}  ( \O _X )$, by definition
$u ^\flat ( \M) := u ^{-1} \R \mathcal{H}om _{\O _X} ( u _* \O _{Z}, \M)$ (see \cite[III.6]{HaRD}).

If $\M$ is a right $\D _X ^{(m)}$-module, we denote by 
$u ^{\flat 0} (\M)
:=
u ^{-1} \mathcal{H}om _{\O _X} ( u _* \O _{Z}, \M)$.
To simplify notation, we will write 
$u ^{\flat 0} (\M)
:=
\mathcal{H}om _{\O _X} ( \O _{Z}, \M)$.
We have a canonical structure of 
right $\D _Z ^{(m)}$-module on
$u ^{\flat 0} (\M)
:=
u ^{-1} \mathcal{H}om _{\O _X} ( u _* \O _{Z}, \M)$
which is constructed by using the $m$-PD-costratification associated to $\M$.
In fact, using the canonical homomorphism
$\cD ^{(m)} _{Z/S}
\to 
 \cD ^{(m)} _{Z \to X/S}
 =
u ^{-1}( \cD ^{(m)} _{X/S}
 /\cI \cD ^{(m)} _{X/S})$,
there is another canonical way to give a structure of 
right $\D _Z ^{(m)}$-module on
$u ^{\flat 0} (\M)$.
Indeed, suppose $X$ affine. 
Let $x \in \Gamma (Z,  u ^{\flat 0} (\M))$ and 
$Q \in D ^{(m)} _{Z}$. 
For any 
$Q _X\in D ^{(m)} _{X}$ 
such that  $\theta (Q ) = \overline{Q _X}$,
we define $x  \cdot Q$ so that we get the equality
\begin{equation}
\label{rem-locadesc-uflat1-formula}
\mathrm{ev} _1
(x  \cdot Q)
: = 
\mathrm{ev} _1 (x)  \cdot Q _X,
\end{equation}
where $\mathrm{ev} _1 \colon \Gamma ( Z, u ^{\flat 0} (\M) )
\hookrightarrow \Gamma (X, \M)$ is 
the evaluation at $1$ homomorphism (which is injective).
Since $I$ annihilates $\mathrm{ev} _1 (x) $,
we remark that this is well defined. 
Both canonical structures of
right $\D _Z ^{(m)}$-module on
$u ^{\flat 0} (\M)$ are identical.

Since $\D _X ^{(m)}$ is a flat $\O _X$-module, then an injective 
right $\D _X ^{(m)}$-module is an injective $\O_X$-module. 
Hence, taking an injective resolution of a complex 
of $D ^{+} ( {} ^r \D _X ^{(m)})$, 
we check the functor $u ^\flat$ sends
$D ^{+}  ( {} ^r \D _X ^{(m)})$
to
$D ^{+} ({} ^r \D _Z ^{(m)} )$, i.e. 
it induces
\begin{equation}
\label{dfn-uflat}
u ^\flat
\colon 
D ^{+}  ( {} ^r \D _X ^{(m)})
\to
D ^{+} ({} ^r \D _Z ^{(m)} ).
\end{equation}
When the level $m$ is ambiguous, 
we denote it more specifically by
$u ^{\flat (m)}$.

Since $X$ is locally noetherian, then $u ^\flat$ preserves the quasi-coherence and
sends 
$D ^{+} _{\mathrm{qc}} ( {} ^r \D _X ^{(m)})$
to 
$D ^{+} _{\mathrm{qc}} ({} ^r \D _Z ^{(m)} )$.

\end{empt}

\begin{empt}
[Local description of the right $\D _X ^{(m)}$-module structure of $u ^{\flat 0} (\M) $]
\label{locadesc-uflat2} 
Suppose we are in the local situation of \ref{locdesc-climm}. 
Let $\M$ be a right $\D _X ^{(m)}$-module.
We have the isomorphism  
\begin{gather}
\label{locadesc-uflat2-iso2pre}
\rho _* u ^{\flat 0} (\M)
\riso 
u ^{-1}
\mathcal{H}om _{\D ^{(m)} _{X,Z,\underline{t}/S} } ( \D ^{(m)} _{X,Z,\underline{t}/S}/\cI \D ^{(m)} _{X,Z,\underline{t}/S}, \M)
\end{gather}
of right $u _*\D ^{(m)} _{Z}$-modules, where
$\rho$ is defined in \ref{morp-rho}.
If there is no ambiguity, 
we can avoid writing $u ^{-1}$, $u _*$ and $\rho _*$.
\end{empt}

\begin{empt}
Suppose we are in the local situation of \ref{locdesc-climm}. 
Let $\M$ be a right $\D _X ^{(m)}$-module.
For simplicity, 
we remove $\rho _*$ in the notation
and 
we  view 
$u ^{\flat} (\M)$
as an object of 
$D ^{\mathrm{b}} ({} ^r\D ^{(m)} _{X,Z,\underline{t}/S})$
(via the isomorphism \ref{morp-rho}).
By derivating 
\ref{locadesc-uflat2-iso2pre},
we get the isomorphism of $D ^{\mathrm{b}} ({} ^r\D ^{(m)} _{X,Z,\underline{t}/S})$ of the form
\begin{equation}
\label{Rlocadesc-uflat2-iso2pre}
u ^{\flat} (\M) \riso \R \mathcal{H}om _{\D ^{(m)} _{X,Z,\underline{t}/S} } (  \D ^{(m)} _{X,Z,\underline{t}/S}/\cI \D ^{(m)} _{X,Z,\underline{t}/S} , \M).
\end{equation}

Let $s := d -r$, and $f _1= t _{r+1},\dots, f _s := t _{d}$. 
Let $K _{\bullet} (\underline{f})$ be the Koszul complex
of $\underline{f}=(f_1,\dots,f_s)$. Let $e _1, \dots, e _s$ be the canonical basis of $\O _X ^s$. Recall 
$K _i (\underline{f}) = \wedge ^i (\O _X ^s)$ and 
$d _{i, \underline{f}} \colon K _i (\underline{f}) \to K _{i-1} (\underline{f})$ 
(or simply $d _i$) is the $\O _X$-linear map defined by 
$$d _i ( e _{n _1} \wedge \cdots \wedge e _{n _i})
= 
\sum _{j=1} ^{i}
(-1) ^{j-1}
f _{n _j} e _{n _1} \wedge \cdots \wedge \widehat{e} _{n _j} \wedge \cdots \wedge e _{n _i}.$$
This yields the  isomorphism of right $\D ^{(m)} _{X,Z,\underline{t}/S}$-modules
\begin{equation}
\label{isos-fund-isom1pre}
\phi ^s _{\underline{t}} 
=
\mathcal{H} ^s(\phi _{\underline{t}})
\colon 
R ^s u ^{\flat 0} (\M)
\riso 
\mathcal{H} ^s \mathcal{H}om _{\D ^{(m)} _{X,Z,\underline{t}/S} } ( \D ^{(m)} _{X,Z,\underline{t}/S} \otimes _{\O _X} K _{\bullet} (\underline{f}) , \M)
\riso  \M / \I \M.
\end{equation}

\end{empt}

\begin{ntn}
\label{dfn-u!leftclimm}
If $\cE$ is a left $\D _X ^{(m)}$-module, we denote by 
$u ^{*} (\cE)
:=
 \O _{Z}
 \otimes _{ u ^{-1}\O _X}
 u ^{-1}\cE$.
Using $m$-PD-stratifications, 
we get a structure of left $\D _Z ^{(m)}$-module 
on $u ^{*} (\cE)$.
Via the homomorphisms of left $\D _X ^{(m)}$-modules of the form $\D  ^{(m)} _X\to \cE$, 
we check by functoriality that 
the canonical homomorphism 
\begin{equation}
\label{dfn-u!leftclimm-iso}
\O _{Z}
 \otimes _{ u ^{-1}\O _X}
 u ^{-1}\cE
 \to \D  ^{(m)} _{Z \to X}
 \otimes  _{ u ^{-1}\D  ^{(m)} _X}
 u ^{-1}\cE
\end{equation}
 is an isomorphism of 
 left $\D _Z ^{(m)}$-modules.
By deriving, we get the functor
$\L u ^* \colon 
D ^{+}  ( {} ^l \D _X ^{(m)}) 
\to
D ^{+} ({} ^l \D _Z ^{(m)} )$
defined by setting 
\begin{equation}
\label{ntnu*form-dfn1pre}
\L u ^{*} (\cE)
:=
\D  ^{(m)} _{Z \to X}
 \otimes ^\L _{ u ^{-1}\D  ^{(m)} _X}
 u ^{-1}\cE.
\end{equation}
Finally, we set
$u ^! (\cE) := \L u ^{*} (\E ) [\delta _{Z/X}]$.

Suppose now we are in the local situation of \ref{locdesc-climm}. 
Let $Q \in D ^{(m)} _{Z}$. 
Choose $Q _X\in D ^{(m)} _{X,Z,\underline{t}/S}$ 
such that  $\overline{Q _X}= \theta (Q)$.
From \ref{dfn-u!leftclimm-iso}, we check the formula
\begin{equation}
\label{dfn-u!leftclimm-pre2}
Q (u ^* (x))
=
u ^* (Q _X \cdot x)).
\end{equation}
Let $\E \in D ^{+}  ( {} ^l  \D  ^{(m)} _X) $.
We have the canonical isomorphisms of $D ^{+}  ( {} ^l  \D  ^{(m)} _Z) $:
$$
 ( \D  ^{(m)} _{X, Z ,\underline{t}} \otimes _{\O _{X}} K _{\bullet} (\underline{f}) ) 
\otimes _{ u ^{-1}\D  ^{(m)} _{X, Z ,\underline{t}} } u ^{-1}\E
\to 
\D  ^{(m)} _{ Z } \otimes ^\L _{u ^{-1} \D  ^{(m)} _{X, Z ,\underline{t}}} u ^{-1}\E
 \riso
 \L u ^* (\E). $$
\end{ntn}

\begin{prop}
\label{fund-isom}
Let $\E$ be a left $\D _X ^{(m)}$-module (resp. a $\D _X ^{(m)}$-bimodule). 
We have the canonical isomorphism of right $\D _Z ^{(m)}$-modules (resp. 
of right $( \D _{Z} ^{(m)}, u ^{-1}\D _{X} ^{(m)})$-bimodules):
\begin{equation}
\label{fund-isom1}
R ^{-\delta _{Z/X}} u ^{\flat 0} ( \omega _{X} \otimes _{\O _X} \E)
\riso 
\omega _{Z} \otimes _{\O _Z} u ^{*} (\E ).
\end{equation}
\end{prop}

\begin{proof}
We can copy the proof of 
\cite[4.1.9]{caro-6operations}.
\end{proof}

\begin{coro}
\label{fund-isom-thmpre}
Let $\E \in D ( {} ^l \D _X ^{(m)})$ (resp. $\E \in D  ( {} ^l \D _X ^{(m)}, {} ^r \D _X ^{(m)})$). 
We have the canonical isomorphism of 
$D  ({} ^r \D _Z ^{(m)})$
(resp. $D  ({} ^r \D _Z ^{(m)},{} ^r u ^{-1}\D _{X} ^{(m)})$)
\begin{equation}
\label{fund-isom2pre}
\omega _{Z} \otimes _{\O _Z} u ^{!} (\E )
\riso
u ^{\flat} ( \omega _{X} \otimes _{\O _X} \E).
\end{equation}
\end{coro}

\begin{proof}
We can copy the proof of 
\cite[4.1.10]{caro-6operations}.
\end{proof}

\begin{empt}
With notation \ref{fund-isom-thmpre}, 
similarly to \ref{fund-isom2bisprepre}, 
the functor
$u ^! 
\colon 
D  ({} ^r \widetilde{\D} _X )
\to 
D  ({} ^r \widetilde{\D} _Z )$
(resp. 
$u ^! 
\colon 
D  ({} ^r \widetilde{\D} _X ,{} ^r \widetilde{\D} _{X} )
\to 
D  ({} ^r \widetilde{\D} _Z ,{} ^r u ^{-1}\widetilde{\D} _{X} )$)
satisfies the isomorphism
\begin{equation}
\label{fund-isom2bispre}
\omega _{Z} \otimes _{\O _Z} u ^{!} (\E ) 
\riso
u ^{!} ( \omega _{X} \otimes _{\O _X} \E).
\end{equation}
Hence, with \ref{fund-isom2pre}, 
we get the isomorphism
\begin{equation}
\label{fflat=f!-immer}
u ^{\flat} \riso u ^!
\end{equation}
of functors
$D  ({} ^r \widetilde{\D} _X )
\to 
D  ({} ^r \widetilde{\D} _Z )$
(resp. 
$D  ({} ^r \widetilde{\D} _X ,{} ^r \widetilde{\D} _{X} )
\to 
D  ({} ^r \widetilde{\D} _Z ,{} ^r u ^{-1}\widetilde{\D} _{X} )$).
\end{empt}

\begin{coro}
\begin{enumerate}[(a)]
\item We have the canonical isomorphism of right $( \D _{Z} ^{(m)}, u ^{-1}\D _{X} ^{(m)}) $-bimodules of the form
\begin{equation}
\label{fund-isom2-corpre}
\omega _{Z} \otimes _{\O _Z} \D _{Z\to X} ^{(m)} 
\riso
u ^{\flat} _{\mathrm{l}} ( \omega _{X} \otimes _{\O _X}  \D _X ^{(m)}) [-\delta _{Z/X}],
\end{equation}
where ``${\mathrm{l}}$'' means that in we have chosen the left structure of 
right $\D _X ^{(m)}$-module
of the right $\D _X ^{(m)} $-bimodule
$\omega _{\X} \otimes _{\O _X} \D _X ^{(m)}$.

 \item We have the canonical isomorphism of  $( u ^{-1}\D _{X} ^{(m)}, \D _{Z} ^{(m)}) $-bimodules of the form
\begin{equation}
\label{fund-isom2-corpre-bis}
 \D _{X\leftarrow Z} ^{(m)} 
\riso
u ^{\flat}  (  \D _X ^{(m)}) [-\delta _{Z/X}].
\end{equation}
\end{enumerate}
\end{coro}

\begin{proof}
By apply Theorem \ref{fund-isom-thmpre} in the case
$\cE= \D _X ^{(m)}$, we get the isomorphism \ref{fund-isom2-corpre}.
By apply Theorem \ref{fund-isom-thmpre} in the case
$\cE= \D _X ^{(m)}$ and
by using the transposition isomorphism
$\omega _{X} \otimes _{\O _X} ( \D _X ^{(m)}  \otimes _{\O _X}  \omega _{X}  ^{-1})
\riso 
\D _X ^{(m)}$
we get the isomorphism \ref{fund-isom2-corpre-bis}.
\end{proof}

\subsection{Adjunction, relative duality isomorphism for schemes}
We keep notation \ref{ntnclosedimmer-1}.
The subsection \cite[4.2]{caro-6operations} is still valid without new argument in our context. 
For the reader, let us collect its results and local description below.
\begin{ntn}
We get the functors
$u ^{(m)} _+\colon 
D  ( {} ^*  \D _Z ^{(m)}) 
\to
D  ({} ^*  \D _X^{(m)} )$ 
by setting
for any 
$\cE \in D   ( {} ^l  \D _Z ^{(m)}) $
and
$\cN \in D   ( {} ^r  \D _Z ^{(m)}) $
by setting
\begin{gather}
\label{dfnu+(m)}
u ^{(m)} _+ (\cE) := u _* 
\left ( 
\D  _{X \leftarrow Z} ^{(m)}
\otimes _{\D  _Z ^{(m)}} 
\cE 
\right ),
u ^{(m)} _+ (\cN) := u _* 
\left ( \cN \otimes _{\D  _Z ^{(m)}} \D  _{Z \to X} ^{(m)}
\right ).
\end{gather}

Moreover, we get the dual functors
$\DD ^{(m)} \colon 
D  ( {} ^*  \D _X ^{(m)}) 
\to
D  ({} ^*  \D _X^{(m)} )$ by setting
for any 
$\cE \in D   ( {} ^l  \D _X ^{(m)}) $
and
$\cM \in D   ( {} ^r  \D _X ^{(m)}) $,
\begin{gather}
\label{dfnuD(m)}
\DD ^{(m)}  (\cE) := 
\R \mathcal{H}om _{\D _X ^{(m)}}
(\cE, \D _X ^{(m)} \otimes _{\O _X} \omega _{X} ^{-1}) [\delta _{X}],
\
\
\DD ^{(m)}  (\cM) := 
\R \mathcal{H}om _{\D _X ^{(m)}}
(\cM,\omega _{X} \otimes _{\O _X}  \D _X ^{(m)}) [\delta _{X}],
\end{gather}
which are respectively computed by taking an injective resolution of
$\D _X ^{(m)} \otimes _{\O _X} \omega _{X} ^{-1}$ 
and $\omega _{X} \otimes _{\O _X}  \D _X ^{(m)}$.
These functors preserve the coherence. 
We can remove $(m)$ in the notation if there is no ambiguity with the level.

These functors are compatible with the quasi-inverse functors
$-\otimes _{\O _X} \omega _{X} ^{-1}$ and
$\omega _{X}  \otimes _{\O _X} -$ exchanging left and right 
$\D ^{(m)} _X$-modules structures.
More precisely, we have the canonical isomorphism
\begin{equation}
\label{u+D-left2right}
\omega _{X}  \otimes _{\O _X} u  ^{(m)} _+ (\cE)
\riso 
u ^{(m)} _+ (\omega _{Z}  \otimes _{\O _Z}  \cE),
\end{equation}
which is constructed as follows :
\begin{gather}
\notag
\omega _{X}  \otimes _{\O _X} u _* 
\left ( 
\D ^{(m)}  _{X \leftarrow Z}  
\otimes _{\D ^{(m)}  _Z } 
\cE
\right )
\riso
u _* 
\left ( 
( u ^{-1}\omega _{X}  \otimes _{u ^{-1} \O _X}  \D ^{(m)}  _{X \leftarrow Z}  )
\otimes _{\D ^{(m)}  _Z } 
\cE
\right )
\\
\notag
\riso
u _* 
\left ( 
(\omega _{Z}  \otimes _{\O _Z}   \cE) 
\otimes _{\D ^{(m)}  _Z } 
( u ^{-1}\omega _{X}  \otimes _{u ^{-1} \O _X}  \D ^{(m)}  _{X \leftarrow Z}  
 \otimes _{\O _Z}  \omega _{Z}  ^{-1})
\right )
\riso
u _* 
\left ( 
(\omega _{Z}  \otimes _{\O _Z}   \cE) 
\otimes _{\D ^{(m)}  _Z } 
 \D ^{(m)}  _{Z \to X}  
\right ).
\end{gather}
More easily, we can check the canonical isomorphism
$\omega _{X}  \otimes _{\O _X} \DD ^{(m)} (\cE)
\riso 
\DD^{(m)} (\omega _{X}  \otimes _{\O _X}  \cE)$.

\end{ntn}

\begin{prop}
\label{u+uflat-lem}
Let $\M$ be a 
right $\D _X ^{(m)}$-module, 
$\cN$ be a 
right  $\D _Z ^{(m)}$-module.
\begin{enumerate}[(a)]
\item We have the canonical adjunction morphisms
$\mathrm{adj} 
\colon 
u _+ u ^{\flat 0} (\M)
\to 
\M$
and 
$\mathrm{adj} 
\colon 
\cN
\to 
 u ^{\flat 0} u _+  (\cN)$.
Moreover, the compositions
$ u ^{\flat 0} (\M)
 \overset{\mathrm{adj}}{\longrightarrow} 
u ^{\flat 0} u _+ u ^{\flat 0} (\M)
\overset{\mathrm{adj}}{\longrightarrow} 
 u ^{\flat 0} (\M)$
and 
$u _+  (\cN)
\overset{\mathrm{adj}}{\longrightarrow} 
u _+ u ^{\flat 0} u _+  (\cN)
\overset{\mathrm{adj}}{\longrightarrow} 
u _+  (\cN)$
are the identity.

\item Using the above  adjunction morphisms, we construct maps 
$$\mathcal{H}om _{\D _X ^{(m)}}
(u _+ (\cN), \M)
\to 
u _* \mathcal{H}om _{\D _Z ^{(m)} }
(\cN, u ^{\flat 0} (\M)),
\
u _* \mathcal{H}om _{\D _Z ^{(m)} }
(\cN, u ^{\flat 0} (\M))
\to 
\mathcal{H}om _{\D _X ^{(m)}}
(u _+ (\cN), \M),$$
which are  inverse of each other.
\item The functor $u ^\flat$ transforms 
$K$-injective complexes into
$K$-injective complexes.
\end{enumerate}

\end{prop}

\begin{proof}
We can copy the proof of \cite[4.2.2]{caro-6operations}.
\end{proof}

\begin{coro}
\label{Radj-u+flat}
Let 
$\M \in D 
( {} ^r \D _X ^{(m)})$,
$\cN\in D 
({} ^r \D _Z ^{(m)} )$. 
Let 
$\cE \in D 
( {} ^l \D _X ^{(m)})$,
$\cN\in D 
({} ^l \D _Z ^{(m)} )$. 
We have the isomorphisms
\begin{gather}
\notag
\R \mathcal{H}om _{\D _X ^{(m)}}
(u _+ (\cN), \M)
\riso 
u _* \R \mathcal{H}om _{\D _Z ^{(m)} }
(\cN, u ^{\flat} (\M));
\\
\notag
\R \mathcal{H}om _{\D _X ^{(m)}}
(u _+ (\cE), \cF)
\riso 
u _* \R \mathcal{H}om _{\D _Z ^{(m)} }
(\cE, u ^{!} (\cF)).
\end{gather}

\end{coro}

\begin{proof}
Taking an injective resolution of 
$\M$, the first isomorphism is a consequence of \ref{u+uflat-lem}.2--3.
This yields the second one by using \ref{fund-isom-thmpre}
and \ref{u+D-left2right}.
\end{proof}

\begin{coro}
\label{rel-dual-isom-imm}
Let 
$\cN\in 
D ^{\mathrm{b}} _{\mathrm{coh}}({} ^* \D _Z ^{(m)} )$.
We have the isomorphism of $D ^{\mathrm{b}} _{\mathrm{coh}}({} ^* \D _X ^{(m)} )$:
\begin{equation}
\label{rel-dual-isom-imm1}
\DD ^{(m)} \circ u ^{(m)} _+ (\cN)
\riso 
 u ^{(m)} _+ \circ \DD ^{(m)} (\cN).
\end{equation}
\end{coro}

\begin{proof}
By using 
\ref{u+D-left2right}; we reduce to the case $* = r$.
In this case, the isomorphism \ref{rel-dual-isom-imm1} 
is the composition of the following  isomorphisms :
\begin{gather}
\notag
\R \mathcal{H}om _{\D _X ^{(m)}}
(u _+ (\cN),\omega _{X} \otimes _{\O _X}  \D _X ^{(m)}) [\delta _{X}]
\overset{\ref{Radj-u+flat}}{\riso} 
u _* \R \mathcal{H}om _{\D _Z ^{(m)} }
(\cN, u ^{\flat} ( \omega _{X} \otimes _{\O _X}  \D _X ^{(m)}) )
[\delta _{X}]
\overset{\ref{fund-isom2-corpre}}{\riso} 
\\
\notag
u _* \R \mathcal{H}om _{\D _Z ^{(m)} }
(\cN, \omega _{Z} \otimes _{\O _Z} \D _{Z\to X} ^{(m)} ) 
[\delta _{Z}]
\underset{\cite[2.1.17]{caro_comparaison}}{\riso} 
u _* 
\left ( 
\R \mathcal{H}om _{\D _Z ^{(m)} }
(\cN,  \omega _{Z} \otimes _{\O _Z} \D _{Z} ^{(m)} 
[\delta _{Z}]
)
\otimes _{\D _{Z} ^{(m)} }
 \D _{Z\to X} ^{(m)} 
 \right ).
\end{gather}
\end{proof}

\subsection{The fundamental isomorphism for formal schemes}
\label{ntn-form-sche-closimm}
The subsection \cite[4.3]{caro-6operations} is still valid without new argument in our context. 
For the reader, let us collect its results and local descriptions below.
Let $u \colon \ZZ \hookrightarrow \X$ be a closed immersion of 
 formal  $\fS$-schemes 
locally of formal finite type
and having locally finite $p$-bases  over $\fS$.
Let $\I$ be the ideal defining $u$. 
The level $m\in \N$ is fixed.
In this subsection, by 
the letter $\widetilde{\D} $ we mean
$\widehat{\D}  ^{(m)}$ or respectively
$\D ^{\dag } \otimes _\bbZ \bbQ$.
For instance, 
$\widetilde{\D} _{\X/\fS} $ is 
$\widehat{\D} ^{(m)} _{\X/\fS} $ 
(resp. 
$\D ^{\dag} _{\X/\fS,\Q} $).

\begin{empt}
[Local description]
\label{locdesc-climm-form}
Suppose
$\fX$ is affine and there exist
$t _{r +1},\dots , t _{d}  \in \Gamma (\fX ,\I)$
generating 
$I:=\Gamma (\fX \,,\I)$,
$t _{1},\dots , t _{r}\in \Gamma (  \fX,\O _{ \fX})$
such that
$t _{1},\dots ,t  _{d}$ form a finite $p$-basis of $ \fX$ over $S$,
$\overline{t} _{1},\dots ,\overline{t} _{r}$ 
form a finite $p$-basis of $ \fZ$ over $\fS$,
and 
$\overline{t} _{r +1},\dots ,\overline{t} _{d}$ is a basis of $\I /\I ^2$,
where $\overline{t} _{1},\dots , \overline{t} _{r} \in\Gamma (  \fZ ,\O _{ \fZ})$
(resp. $\overline{t} _{r +1},\dots ,\overline{t} _{d}\in\Gamma (  \fX ,\cI /\cI ^2)$)
are the images of 
$t _{1},\dots , t _{r}$
(resp. $t _{r+1},\dots , t _{d}$)
via 
$\Gamma (  \fX,\O _{ \fX})
\to
\Gamma (  \fZ ,\O _{ \fZ})$
(resp. 
$\Gamma (  \fX,\cI)
\to 
\Gamma (  \fX ,\cI /\cI ^2)$).

We denote by 
$\tau _i := 1 \otimes t _i -t _i \otimes 1$, 
$\overline{\tau} _j := 1 \otimes \overline{t} _j -\overline{t} _j \otimes 1$, 
for any $i= 1,\dots, d$, $j= 1,\dots, r$.
The sheaf of $\O _ \fX$-algebras
$\cP ^n  _{ \fX/\fS,(m)}$ is 
a free $\O _ \fX$-module with the basis 
$\{ \underline{\tau} ^{\{\underline{k}\} _{(m)}}\ | \ \underline{k}\in \N ^d \text{ such that } | \underline{k}| \leq n\} $,
and
$\cP ^n  _{ \fZ/\fS,(m)}$ is 
a free $\O _ \fZ$-module with the basis 
$\{ \underline{\overline{\tau}} ^{\{\underline{l}\} _{(m)}}\ | \ \underline{l}\in \N ^r \text{ such that } | \underline{l}| \leq n\} $.
We denote by 
$\{ 
\underline{\partial} ^{<\underline{k}> _{(m)}}
\ | \ \underline{k}\in \N ^d,\ | \underline{k} | \leq n
\}$
the corresponding dual basis  of 
$\D ^{(m)} _{ \fX/\fS,n}$ 
and
by 
$\{ \underline{\partial} ^{<\underline{l}> _{(m)}}\  
| \ \underline{l}\in \N ^r,\ | \underline{l} | \leq n\} $
the corresponding dual basis of $\D ^{(m)} _{ \fZ/ \fS,n}$ (if there is no possible confusion).
The sheaf
$\D ^{(m)} _{ \fX/ \fS}$ is 
a free $\O _ \fX$-module with the basis 
$\{ 
\underline{\partial} ^{<\underline{k}> _{(m)}}
\ | \ \underline{k}\in \N ^d
\}$,
and
$\D ^{(m)} _{ \fZ/ \fS}$ is 
a free $\O _ \fZ$-module with the basis 
$\{ \underline{\partial} ^{<\underline{l}> _{(m)}}\ | \ \underline{l}\in \N ^r\} $.

a) We compute the canonical homomorphism
$u ^* \cP ^n  _{ \fX/ \fS,(m)}
\to \cP ^n  _{ \fZ/ \fS,(m)}$
sends 
$\underline{\tau} ^{\{(\underline{l}, \underline{h})\} _{(m)}}$
where 
$\underline{l} \in \N ^r$
and 
$\underline{h} \in \N ^{d-r}$
to
$\underline{\overline{\tau}} ^{\{\underline{l}\} _{(m)}}$
if $\underline{h} = (0,\dots,0)$
and to $0$ otherwise.

b) We denote by 
$\theta \colon 
\cD ^{(m)} _{ \fZ/ \fS}
\to 
\cD ^{(m)} _{ \fZ \to  \fX/ \fS}$
the canonical homomorphism of left $\cD ^{(m)} _{ \fZ/ \fS}$-modules
(which is built by duality from the 
canonical homomorphisms $u ^* \cP ^n  _{ \fX/ \fS,(m)}
\to \cP ^n  _{ \fZ/ \fS,(m)}$).
For any $P \in D ^{(m)} _{ \fX/ \fS}$, we denote by 
$\overline{P}$ its image via
the canonical morphism of left $D ^{(m)} _{ \fX/ \fS}$-modules
$D ^{(m)} _{ \fX/ \fS}
\to 
D ^{(m)} _{ \fX/ \fS} / I D ^{(m)} _{ \fX/ \fS}
= 
D ^{(m)} _{ \fZ \to  \fX/ \fS}$.
We set 
$\underline{\xi} ^{<\underline{k}> _{(m)}}:= 
\overline{\underline{\partial} ^{<\underline{k}> _{(m)}}}$. 
By duality from a), we compute
$\theta (\underline{\partial} ^{<\underline{l}> _{(m)}})
=
\underline{\xi} ^{<(\underline{l}, \underline{0})> _{(m)}}$,
for any 
$\underline{l}\in \N ^r$.

\end{empt}

\begin{empt}
\label{locdesc-climm2-form}
Suppose we are in the local situation of \ref{locdesc-climm-form}. 
We denote by 
$\cD ^{(m)} _{ \fX, \fZ,\underline{t}/ \fS}$ the subring of 
$\cD ^{(m)} _{ \fX/ \fS}$ which is a 
 free $\O _ \fX$-module with the basis 
$\{ \underline{\partial} ^{<(\underline{l}, \underline{0})> _{(m)}}\ | \ \underline{l}\in \N ^r\} $, 
where 
$\underline{0}:=(0,\dots, 0) \in 
\N ^{d-r}$.
If there is no ambiguity concerning the finite $p$-basis (resp. and $\fS$),
we might simply denote $\cD ^{(m)} _{ \fX, \fZ,\underline{t}/ \fS}$ by 
$\cD ^{(m)} _{ \fX, \fZ/ \fS}$
(resp. $\cD ^{(m)} _{ \fX, \fZ}$).
The properties of \ref{locdesc-climm2} are still valid in the context of formal schemes, 
we have only to replace
respectively 
$X$, $Z$, $S$ by $\fX$, $\fZ$, $\fS$.

\end{empt}

\begin{empt}
[Local description of $u ^{\flat}$]
\label{locadesc-uflat-form}
Suppose we are in the local situation of \ref{locdesc-climm-form}. 
Let $\M$ be a right $\D _{\X/\fS} ^{(m)}$-module.
Let $x \in \Gamma (\fZ,  u ^{\flat 0} (\M))$ and 
$Q \in D ^{(m)} _{\fZ}$. 
For any 
$Q _\fX\in D ^{(m)} _{\fX,\fZ,\underline{t}/\fS}$ 
such that  $\theta (Q ) = \overline{Q _\fX}$,
we compute
\begin{equation}
\label{locadesc-uflat1-formula-form}
\mathrm{ev} _1
(x  \cdot Q)
= 
\mathrm{ev} _1 (x)  \cdot Q _\fX.
\end{equation}

\end{empt}

\begin{empt}
\label{rem-rightD-moduleflat-form}
Let $\M$ be a right $\widetilde{\D} _\X $-module. 
There is a canonical way to endow 
$u ^{\flat 0} (\M)$
with a structure of 
right $\widetilde{\D} _\fZ$-module.
Indeed, suppose $\fX$ affine. 
Let $x \in \Gamma (\fZ,  u ^{\flat 0} (\M))$ and 
$Q \in \widetilde{D} _{\fZ}$. 
For any 
$Q _\fX\in \widetilde{D} _{\fX}$ 
such that  $\theta (Q ) = \overline{Q _\fX}$,
we define $x  \cdot Q$ so that we get the equality
\begin{equation}
\label{rem-locadesc-uflat1-formula-form}
\mathrm{ev} _1
(x  \cdot Q)
: = 
\mathrm{ev} _1 (x)  \cdot Q _\fX,
\end{equation}
where $\mathrm{ev} _1 \colon \Gamma ( \fZ, u ^{\flat 0} (\M) )
\hookrightarrow \Gamma (\fX, \M)$ is 
the evaluation at $1$ homomorphism (which is injective).
Since $I$ annihilates $\mathrm{ev} _1 (x) $,
we remark that this is well defined.

\end{empt}

\begin{empt}
\label{dfntheta-rhotilde}
Suppose we are in the local situation of \ref{locdesc-climm-form}.
We keep notation \ref{locdesc-climm-form}
and \ref{locdesc-climm2-form}.

\begin{enumerate}\item A section of the sheaf
$\widehat{\D} ^{(m)} _{\X/\fS}$ 
can uniquely be written in the form 
$\sum 
_{\underline{k}\in \N ^d}
a _{\underline{k}}
\underline{\partial} ^{<\underline{k}> _{(m)}}$
such that 
$a _{\underline{k}} \in \O _\X$
converges to $0$ when 
$| \underline{k}|\to \infty$.
A section of the sheaf
$\widehat{\D} ^{(m)} _{\fZ/\fS}$ 
can uniquely be written in the form 
$\sum 
_{\underline{l}\in \N ^r}
b _{\underline{l}}
\underline{\partial} ^{<\underline{l}> _{(m)}}$
such that 
$b _{\underline{l}} \in \O _\fZ$
converges to $0$ when 
$| \underline{l}|\to \infty$.
Let 
$\widehat{\D} ^{(m)} _{\X,\fZ,\underline{t}}$ be
the $p$-adic completion of
$\D ^{(m)} _{\X,\fZ,\underline{t}}$.
Then $\widehat{\cD} ^{(m)} _{\X,\fZ,\underline{t}}$
is a subring of 
$\widehat{\cD} ^{(m)} _{\X/\fS}$ whose elements can uniquely be written in the form
$\sum 
_{\underline{l}\in \N ^r}
a _{\underline{l}}
\underline{\partial} ^{<(\underline{l}, \underline{0})> _{(m)}}$
(recall $\underline{0}:=(0,\dots, 0) \in 
\N ^{d-r}$)
where 
$a _{\underline{l}} \in \O _\fX$
converges to $0$ when 
$| \underline{l}|\to \infty$.
Taking the $p$-adic completion of the diagram \ref{diag-DZ2X-t-comm} (still valid for formal schemes),
we get the canonical diagram
\begin{equation}
\label{diag-DZ2X-t-comm-form}
\xymatrix{
{\widehat{\cD} ^{(m)} _{\fX,\fZ,\underline{t}/\fS}/\cI\widehat{\cD} ^{(m)} _{\fX,\fZ,\underline{t}/\fS}} 
\ar@{^{(}->}[r] ^-{}
& 
{\widehat{\cD} ^{(m)} _{\fX/\fS} / \cI \widehat{\cD} ^{(m)} _{\fX/\fS}} 
\\ 
{u _* \widehat{\cD} ^{(m)} _{\fZ/\fS}} 
\ar[u] ^-{\sim} _-{\widehat{\theta}}
\ar[r] ^-{\widehat{\theta}}
& 
{u _* \widehat{\cD} ^{(m)} _{\fZ \to \fX/\fS}} 
\ar@{=}[u] 
}
\end{equation}
where 
$\widehat{\theta}
\colon 
u _* \widehat{\cD} ^{(m)} _{\ZZ/\fS}
\riso 
\widehat{\cD} ^{(m)} _{\X,\fZ,\underline{t}}
/\cI\widehat{\cD} ^{(m)} _{\X,\fZ,\underline{t}}
$
is an isomorphism of $\V$-algebras.

\item We set 
$\D ^{\dag} _{\X,\fZ,\underline{t},\Q}
:=
\underrightarrow{\lim}
\widehat{\D} ^{(m)} _{\X,\fZ,\underline{t},\bbQ}$.
We get a similar diagram than \ref{diag-DZ2X-t-comm-form}
by replacing $\widehat{\cD} ^{(m)}$ with $\D ^{\dag}$
and by adding some $\bbQ$.

\item The isomorphism of $\cV$-algebras 
$u _* \widetilde{\cD} _{ \fZ/ \fS}
\riso 
\widetilde{\cD} _{ \fX, \fZ,\underline{t}/ \fS}
/\I \widetilde{\cD} _{ \fX, \fZ,\underline{t}/ \fS}$
induced by $\theta$ will be denoted by 
$\widetilde{\theta}$.
This yields by composition the  homomorphism of $\O _ \fX$-rings :
\begin{equation}
\label{morp-rhof-tilde}
\widetilde{\rho} 
\colon
\widetilde{\cD}  _{ \fX, \fZ,\underline{t}/ \fS}
\to
\widetilde{\cD} _{ \fX, \fZ,\underline{t}/ \fS}
/\I \widetilde{\cD} _{ \fX, \fZ,\underline{t}/ \fS}
\underset{\widetilde{\theta}}{\liso}  
u _* \widetilde{\cD} _{ \fZ/ \fS}.
\end{equation}

\end{enumerate}

\end{empt}

\begin{empt}
We have the canonical isomorphism of $\widetilde{\cD} _{\fX,\fZ,\underline{t}/\fS} $-modules
\begin{equation}
\label{locadesc-uflat2-iso2pre-form-tilde}
\widetilde{\rho} _* u ^{\flat 0} (\M)
\riso 
u ^{-1}
\mathcal{H}om _{\widetilde{\cD} _{\fX, \fZ ,\underline{t}} } ( \widetilde{\cD} _{\fX, \fZ ,\underline{t}} /\cI \widetilde{\cD} _{\fX, \fZ ,\underline{t}}, \M)
\end{equation}
If there is no ambiguity, 
we can avoid writing $u ^{-1}$,  $u _*$,

\end{empt}

\begin{empt}
\label{hat-uflat-desc}
Suppose we are in the local situation of \ref{locdesc-climm-form}. Let $\M$ be a right $\widetilde{\D} _{\X/\fS} $-module.
\begin{enumerate}[(a)]

\item Let $s := d -r$, and $f _1= t _{r+1},\dots, f _s := t _{d}$. 
Let $K _{\bullet} (\underline{f})$ be the Koszul complex
of $\underline{f}=(f_1,\dots,f_s)$. 
We have the isomorphism 
 of $D ^{\mathrm{b}} ( \widetilde{\D} _{\fX, \fZ ,\underline{t}})$: 
\begin{equation}
\label{locadesc-uflat-iso-form}
\phi _{\underline{t}}
\colon 
u ^{\flat} (\M)
\riso 
\mathcal{H}om _{ \widetilde{\D} _{\fX, \fZ ,\underline{t}} } (  \widetilde{\D} _{\fX, \fZ ,\underline{t}} \otimes _{\O _{\X}} K _{\bullet} (\underline{f}) , \M).
\end{equation}
This yields the  isomorphisms of right $ \widetilde{\D} _{\fX, \fZ ,\underline{t}}$-modules
\begin{equation}
\label{isos-fund-isom1pre-form}
\phi ^s _{\underline{t}} 
=
\mathcal{H} ^s(\phi _{\underline{f}})
\colon 
R ^s u ^{\flat 0} (\M)
\riso 
\mathcal{H} ^s \mathcal{H}om _{ \widetilde{\D} _{\fX, \fZ ,\underline{t}} } (  \widetilde{\D} _{\fX, \fZ ,\underline{t}} 
\otimes _{\O _{\X}} K _{\bullet} (\underline{f}) , \M)
\riso  \M / \I \M.
\end{equation}

\end{enumerate}

\end{empt}

\begin{ntn}
\label{ntnu*form}
If $\cE$ is a left $\D _\X ^{(m)}$-module, we set
$u ^{*} (\cE)
:=
 \O _{\fZ}
 \otimes _{ u ^{-1}\O _\X}
 u ^{-1}\cE$.
Using $m$-PD-stratifications, 
we get a structure of left $\D _\fZ ^{(m)}$-module 
on $u ^{*} (\cE)$.
This yields the functor
$\L u ^* \colon 
D ^{+}  ( {} ^l \D _\X ^{(m)}) 
\to
D ^{+} ({} ^l \D _\ZZ ^{(m)} )$
(resp. 
$\L u ^* \colon 
D  ( {} ^l \D _\X ^{(m)}) 
\to
D  ({} ^l \D _\ZZ ^{(m)} )$).
Similarly, we get the functor
$\L u ^* \colon 
D ^{+}  ( {} ^l \widetilde{\D} _\X ) 
\to
D ^{+} ({} ^l \widetilde{\D} _\ZZ )$
(resp. 
$\L u ^* \colon 
D  ( {} ^l \widetilde{\D} _\X ) 
\to
D  ({} ^l \widetilde{\D} _\ZZ )$)
defined by setting 
\begin{equation}
\label{ntnu*form-dfn1}
\L u ^{*} (\cE)
:=
\widetilde{\D} _{\fZ \to \fX}
 \otimes ^\L _{ u ^{-1}\widetilde{\D} _\X}
 u ^{-1}\cE.
\end{equation}
Finally, we set
$u ^! (\cE) := \L u ^{*} (\E ) [\delta _{Z/X}]$.

Suppose we are in the local situation of \ref{locdesc-climm-form}. 
Let $\E \in D ( {} ^l  \widetilde{\D} _\X) $.
The canonical homorphism 
$$
 \widetilde{\D} _{ \fZ } \otimes ^\L _{u ^{-1} \widetilde{\D} _{\fX, \fZ ,\underline{t}}} u ^{-1}\E
 \to
 \L u ^* (\E) \riso ( \widetilde{\D} _{\fX, \fZ ,\underline{t}} \otimes _{\O _{\X}} K _{\bullet} (\underline{f}) ) 
\otimes _{ u ^{-1}\widetilde{\D} _{\fX, \fZ ,\underline{t}} } u ^{-1}\E.$$
 is an isomorphism of $D   ( {} ^l  \widetilde{\D} _\fZ) $.

\end{ntn}

\begin{prop}
\label{fund-isom-form}
Let $\E$ be a left $ \widetilde{\D} _{\X}$-module (resp. a $ \widetilde{\D} _{\X}$-bimodule). 
Set $n := -\delta _{Z/X} \in \N$.
We have the canonical isomorphism of right $ \widetilde{\D} _\fZ $-modules (resp. 
of right $(  \widetilde{\D} _{\fZ}, u ^{-1} \widetilde{\D} _{\X})$-bimodules):
\begin{equation}
\label{fund-isom1-form}
R ^{n} u ^{\flat 0} ( \omega _{\X} \otimes _{\O _{\X}} \E)
\riso 
\omega _{\fZ} \otimes _{\O _{\fZ}} u ^{*} (\E ).
\end{equation}
\end{prop}

\begin{proof}
Using \ref{hat-uflat-desc} and \ref{ntnu*form}, we proceed as \ref{fund-isom}.
\end{proof}

\begin{cor}
\label{fund-isom-thm-form}
Let $* \in \{\mathrm{l}, \mathrm{r}\}$ and 
let $\E \in D  ( {} ^l \widetilde{\D} _\fX )$ (resp. $\E \in D  ( {} ^l \widetilde{\D} _\fX , {} ^* \widetilde{\D} _\fX )$). 
With notation \ref{ntnu*form}, 
we have the canonical isomorphism of 
$D  ({} ^r \widetilde{\D} _\fZ )$
(resp. $D  ({} ^r \widetilde{\D} _\fZ ,{} ^* u ^{-1}\widetilde{\D} _{\fX} )$)
of the form
\begin{equation}
\label{fund-isom2}
\omega _{\ZZ} \otimes _{\O _\ZZ} u ^{!} (\E ) 
\riso
u ^{\flat} ( \omega _{\X} \otimes _{\O _\X} \E).
\end{equation}
\end{cor}

\begin{proof}
Using
\cite[I.7.4]{HaRD},
this is a consequence of \ref{fund-isom-form}.
\end{proof}

\begin{empt}
With notation \ref{fund-isom-thm-form}, 
the functor
$u ^! 
\colon 
D  ({} ^r \widetilde{\D} _\fX )
\to 
D  ({} ^r \widetilde{\D} _\fZ )$
(resp. 
$u ^! 
\colon 
D ({} ^r \widetilde{\D} _\fX ,{} ^* \widetilde{\D} _{\fX} )
\to 
D  ({} ^r \widetilde{\D} _\fZ ,{} ^* u ^{-1}\widetilde{\D} _{\fX} )$)
satisfies the isomorphism
\begin{equation}
\label{fund-isom2bis}
\omega _{\ZZ} \otimes _{\O _\ZZ} u ^{!} (\E ) 
\riso
u ^{!} ( \omega _{\X} \otimes _{\O _\X} \E).
\end{equation}
Hence, with \ref{fund-isom2}, 
we get the isomorphism
\begin{equation}
\label{fflat=f!-immer-f}
u ^{\flat} \riso u ^!
\end{equation}
of functors
$D ({} ^r \widetilde{\D} _\fX )
\to 
D  ({} ^r \widetilde{\D} _\fZ )$
(resp. 
$D  ({} ^r \widetilde{\D} _\fX ,{} ^* \widetilde{\D} _{\fX} )
\to 
D  ({} ^r \widetilde{\D} _\fZ ,{} ^* u ^{-1}\widetilde{\D} _{\fX} )$).

\end{empt}

\begin{coro}
\begin{enumerate}[(a)]
\item We have the canonical isomorphism of right $( \widetilde{\D} _{\fZ} , u ^{-1}\widetilde{\D} _{\fX} ) $-bimodules of the form
\begin{equation}
\label{fund-isom2-corpref}
\omega _{\fZ} \otimes _{\O _\fZ} \widetilde{\D} _{\fZ\to \fX}  
\riso
u _{\mathrm{l}} ^{\flat} ( \omega _{\fX} \otimes _{\O _\fX}  \widetilde{\D} _\fX ) [-\delta _{Z/X}],
\end{equation}
where ``${\mathrm{l}}$'' means that in we have chosen the left structure of 
right $ \widetilde{\D} _\fX $-module
of the right $\widetilde{\D} _\fX $-bimodule
$\omega _{\fX} \otimes _{\O _\fX}  \widetilde{\D} _\fX $.

\item We have the canonical isomorphism of $( u ^{-1}\widetilde{\D} _{\fX} , \widetilde{\D} _{\fZ}) $-bimodules of the form
\begin{equation}
\label{fund-isom2-corpref-bis}
 \widetilde{\D} _{\fX \leftarrow \fZ}  
\riso
u ^{\flat} ( \widetilde{\D} _\fX )[-\delta _{Z/X}].
\end{equation}

\end{enumerate}

\end{coro}

\begin{proof}
By applying Theorem \ref{fund-isom-thm-form} in the case
$\cE= \widetilde{\D} _\fX $, we get the isomorphism \ref{fund-isom2-corpref}.
By applying Theorem \ref{fund-isom-thm-form} in the case
$\cE= \widetilde{\D} _\fX  \otimes _{\O _\fX}  \omega _{\fX} ^{-1} $, 
and by using the transposition isomorphism
$\omega _{\fX} \otimes _{\O _\fX} ( \widetilde{\D} _\fX  \otimes _{\O _\fX}  \omega _{\fX} ^{-1} )
\riso 
 \widetilde{\D} _\fX $,
we get the isomorphism \ref{fund-isom2-corpref-bis}.
\end{proof}

\subsection{Adjunction, relative duality isomorphism}
We keep notation \ref{ntn-form-sche-closimm}.
The (end of the) subsection \cite[4.3]{caro-6operations} is still valid without new argument in our context. 
For the reader, let us collect its results and local description below.

\begin{ntn}
We get the functor
$u  _+\colon 
D  ( {} ^*  \widetilde{\D} _\fZ ) 
\to
D  ({} ^*  \widetilde{\D} _\fX )$ by setting
for any 
$\cE \in D   ( {} ^l  \widetilde{\D} _\fZ ) $
and
$\cN \in D   ( {} ^r  \widetilde{\D} _\fZ ) $,
\begin{gather}
\label{dfnu+(m)f}
u  _+ (\cN) := u _* 
\left ( \cN \otimes _{\widetilde{\D}  _\fZ } \widetilde{\D}  _{\fZ \to \fX} 
\right ),
u  _+ (\cE) := u _* 
\left ( 
\widetilde{\D}  _{\fX \leftarrow \fZ}  
\otimes _{\widetilde{\D}  _\fZ } 
\cE
\right ).
\end{gather}
Moreover, we get the functor
$\DD  \colon 
D  ( {} ^*  \widetilde{\D} _\fX ) 
\to
D  ({} ^*  \widetilde{\D} _\fX )$ by setting
for any 
$\cM \in D   ( {} ^r  \widetilde{\D} _\fX ) $,
$\cE \in D   ( {} ^l  \widetilde{\D} _\fX ) $
\begin{gather}
\label{dfnuD(m)f}
\DD   (\cM) := 
\R \mathcal{H}om _{\widetilde{\D} _\fX }
(\cM,\omega _{\fX} \otimes _{\O _\fX}  \widetilde{\D} _\fX ) [\delta _{X}],
\DD   (\cE) := 
\R \mathcal{H}om _{\D _\fX }
(\cE,  \widetilde{\D} _\fX\otimes _{\O _\fX} \omega _{\fX} ^{-1} ) [\delta _{X}],
\end{gather}
which are computed respectively by taking an injective resolution of 
$\omega _{\fX} \otimes _{\O _\fX}  \widetilde{\D} _\fX $
and
$\widetilde{\D} _\fX\otimes _{\O _\fX} \omega _{\fX} ^{-1}$.
These functors preserves the coherence and are compatible with the quasi-inverse functors 
$-\otimes _{\O _\fX} \omega _{\fX} ^{-1}$ and
$\omega _{\fX}  \otimes _{\O _\fX} -$ exchanging left and right 
$\widetilde{\D} _\fX$-modules structure.
More precisely, we have the canonical isomorphisms
\begin{equation}
\label{u+D-left2rightf}
\omega _{\fX}  \otimes _{\O _\fX} u  _+ (\cE)
\riso 
u  _+ (\omega _{\fZ}  \otimes _{\O _\fZ}  \cE),
\omega _{\fX}  \otimes _{\O _\fX} \DD (\cE)
\riso 
\DD (\omega _{\fX}  \otimes _{\O _\fX}  \cE)
\end{equation}
whose first one is constructed as \ref{u+D-left2right}.

\end{ntn}

\begin{prop}
\label{u+uflat-lem-form}
Let $\M$ be a 
right $\widetilde{\D} _\fX $-module, 
$\cN$ be a 
right  $\widetilde{\D} _\fZ $-module.
We keep notations \ref{ntnu*form-dfn1} and \ref{dfnu+(m)f}.
\begin{enumerate}[(a)]
\item We have the canonical adjunction morphisms
$\mathrm{adj} 
\colon 
u _+ u ^{\flat 0} (\M)
\to 
\M$
and 
$\mathrm{adj} 
\colon 
\cN
\to 
 u ^{\flat 0} u _+  (\cN)$.
Moreover, the compositions
$ u ^{\flat 0} (\M)
 \overset{\mathrm{adj}}{\longrightarrow} 
u ^{\flat 0} u _+ u ^{\flat 0} (\M)
\overset{\mathrm{adj}}{\longrightarrow} 
 u ^{\flat 0} (\M)$
and 
$u _+  (\cN)
\overset{\mathrm{adj}}{\longrightarrow} 
u _+ u ^{\flat 0} u _+  (\cN)
\overset{\mathrm{adj}}{\longrightarrow} 
u _+  (\cN)$
are the identity.

\item Using the above  adjunction morphisms, we construct maps 
$$\mathcal{H}om _{\widetilde{\D} _\fX }
(u _+ (\cN), \M)
\to 
u _* \mathcal{H}om _{\widetilde{\D} _\fZ  }
(\cN, u ^{\flat 0} (\M)),
\
u _* \mathcal{H}om _{\widetilde{\D} _\fZ  }
(\cN, u ^{\flat 0} (\M))
\to 
\mathcal{H}om _{\widetilde{\D} _{\fX} }
(u _+ (\cN), \M),$$
which are  inverse of each other.
\item If $\M$ is an injective right $\widetilde{\D} _{\fX} $-module, 
then 
$u ^{\flat 0} (\M)$ is an injective right $\widetilde{\D} _\fZ $-module.
\end{enumerate}

\end{prop}

\begin{proof}
We can copy the proof of \ref{u+uflat-lem}.
\end{proof}

\begin{coro}
\label{Radj-u+flatf}
Let 
$\M \in D  
( {} ^r \widetilde{\D} _{\fX} )$,
$\cN\in D 
({} ^r \widetilde{\D} _\fZ  )$. 
Let 
$\cE \in D 
( {} ^l \widetilde{\D} _{\fX} )$,
$\cF \in D 
({} ^l \widetilde{\D} _\fZ  )$. 
We have the isomorphisms
\begin{gather}
\R \mathcal{H}om _{\widetilde{\D} _{\fX} }
(u _+ (\cN), \M)
\riso 
u _* \R \mathcal{H}om _{\widetilde{\D} _\fZ  }
(\cN, u ^{\flat} (\M));
\\
\R \mathcal{H}om _{\widetilde{\D} _{\fX} }
(u _+ (\cE), \cF)
\riso 
u _* \R \mathcal{H}om _{\widetilde{\D} _\fZ  }
(\cE, u ^{!} (\cF)).
\end{gather}

\end{coro}

\begin{proof}
Taking a K-injective resolution of 
$\M$ (see 13.33.5 of the stack project), 
the first isomorphism  is a consequence of \ref{u+uflat-lem-form}.2--3.
Using \ref{fund-isom2} and \ref{u+D-left2rightf}, we deduce  the second isomorphism from the first one.
\end{proof}

\begin{coro}
\label{rel-dual-isom-immf}
Let 
$\cN\in D ^{\mathrm{b}} _{\mathrm{coh}}  ({} ^* \widetilde{\D} _\fZ  )$
with $*=r$ or $*=l$.
We have the isomorphism of $D ^{\mathrm{b}} _{\mathrm{coh}}  ({} ^* \widetilde{\D} _\fX  )$:
\begin{equation}
\label{rel-dual-isom-immf1}
\DD  \circ u _+ (\cN)
\riso 
 u _+ \circ \DD  (\cN).
\end{equation}
\end{coro}

\begin{proof}
Using \ref{fund-isom2-corpref} and \ref{Radj-u+flatf}, 
we can copy the proof of \ref{rel-dual-isom-imm}.\end{proof}

\begin{prop}
\label{u+uflat-lem-formal}
Let 
$\cN$ be a right coherent $\widetilde{\D} _{\ZZ} $-module without $p$-torsion.
The canonical homomorphism of $\widetilde{\D} _{\ZZ} $-modules 
$\mathrm{adj} 
\colon 
\cN
\to
 u ^{\flat 0} u _+  (\cN)$
 is an isomorphism.
\end{prop}

\begin{proof}
We proceed similarly to \cite[2.3.1]{surcoh-hol}.
\end{proof}

\subsection{Glueing isomorphisms, base change isomorphisms for pushforwards by a closed immersion}

\begin{prop}
\label{prop-glueiniso-coh}
Let $f,f',f''\colon \X  \to \Y $ be three morphisms of 
formal $\fS$-schemes locally of formal finite type
and having locally finite $p$-bases  over $\fS$
such that $f _0=  f' _0=f'' _0$. 
Let $g,g'\colon \Y  \to \ZZ $ be two morphisms of formal schemes 
locally of formal finite type
and having locally finite $p$-bases over $\fS $
such that $g _0=  g' _0$. 
Let $T _Z$ be a divisor of $Z$ such that $T _Y:= g _0 ^{-1} (T)$ is a divisor of $Y$
and 
$T _X := f _0 ^{-1} (T _Y)$ is a divisor of $X$.

\begin{enumerate}[(a)]
\item 
We have the canonical isomorphism of functors 
$\smash{\underrightarrow{LD}} ^{\mathrm{b}} _{\Q,\mathrm{qc}}
 ( \smash{\widetilde{\D}} _{\Y /\fS } ^{(\bullet)}(T _Y))
 \to 
 \smash{\underrightarrow{LD}} ^{\mathrm{b}} _{\Q,\mathrm{qc}}
 ( \smash{\widetilde{\D}} _{\X /\fS } ^{(\bullet)}(T _X))$
 of the form
\begin{equation}
\notag
\tau _{f,f'} ^{(\bullet)}
\colon 
f  _{T _Y} ^{\prime !(\bullet)} \riso f _{T _Y} ^{ !(\bullet)}.
\end{equation}
These isomorphisms satisfy the following formulas
$\tau ^{(\bullet)} _{f,f}=\mathrm{Id}$,
$\tau ^{(\bullet)} _{f,f''}= \tau ^{(\bullet)} _{f,f'} \circ \tau ^{(\bullet)} _{f',f''}$,
$\tau ^{(\bullet)} _{f,f'} \circ g  _{T _Z} ^{!(\bullet)} = \tau ^{(\bullet)} _{g \circ f ,g \circ f'} $ 
and
$f _{T _Y} ^{ !(\bullet)} \circ \tau ^{(\bullet)} _{g,g'}= \tau ^{(\bullet)} _{g \circ f ,g '\circ f} $.

\item There exists a canonical glueing isomorphism of functors 
$D ^{\mathrm{b}} _{\mathrm{coh}} ( \smash{\D} ^{\dag} _{\Y} (\hdag T _Y) _{\Q} )
\to
D ^{\mathrm{b}}  ( \smash{\D} ^{\dag} _{\X } (\hdag T _X ) _{\Q}  )$
of the form
\begin{equation}
\label{prop-glueiniso-coh1}
\tau _{f,f'}\colon 
f  _{T _Y} ^{\prime!} \riso f _{T _Y} ^{!},
\end{equation}
such that $\tau _{f,f}=\mathrm{Id}$,
$\tau _{f,f''}= \tau _{f,f'} \circ \tau _{f',f''}$,
$\tau _{f,f'} \circ g  _{T _Z} ^{!} = \tau _{g \circ f ,g \circ f'} $ 
and
$f _{T _Y} ^{!} \circ \tau _{g,g'}= \tau _{g \circ f ,g '\circ f} $.

\item The diagram 
 of functors
$\smash{\underrightarrow{LD}} ^{\mathrm{b}} _{\Q,\mathrm{coh}}
 ( \smash{\widetilde{\D}} _{\Y /\fS } ^{(\bullet)}(T _Y))
 \to 
 D ^{\mathrm{b}}  ( \smash{\D} ^{\dag} _{\X } (\hdag T _X) _{\Q}  )$
 \begin{equation}
 \notag
 \xymatrix {
 {\underrightarrow{\lim} \circ f  _{T _Y} ^{\prime !(\bullet)} } 
 \ar[r] ^-{\sim} _-{\underrightarrow{\lim} \circ \tau _{f,f'} ^{(\bullet)}}
 \ar[d] ^-{\sim}
 & 
 {\underrightarrow{\lim} \circ f _{T _Y} ^{ !(\bullet)}} 
 \ar[d] ^-{\sim} 
 \\ 
 {f  _{T _Y} ^{\prime!} \circ \underrightarrow{\lim}} 
\ar[r] ^-{\tau _{f,f'} \circ \underrightarrow{\lim} } & 
 {f  _{T _Y} ^{!}\circ \underrightarrow{\lim}} 
 }
 \end{equation}
is commutative up to canonical isomorphism.
\end{enumerate}

\end{prop}

\begin{proof}
1) By copying \cite[2.1.3 and 2.1.10]{caro-construction} (still valid in our context), 
we check the first statement from \ref{215Be2}.
 
2) and 3) Let $\FF\in D ^{\mathrm{b}} _{\mathrm{coh}} ( \smash{\D} ^{\dag} _{\Y} (\hdag T _Y) _{\Q} )$.
Taking inductive limits of the completion of the glueing isomorphisms  \ref{215Be2}, we get the isomorphism 
$\tau _{f,f '}
\colon  \smash{\D} ^{\dag} _{\X \overset{f'}{\rightarrow}\Y}  (\hdag T _Y) _{ \Q}  
\riso 
 \smash{\D} ^{\dag} _{\X \overset{f}{\rightarrow}\Y} (\hdag T _Y) _{ \Q}   $.
It follows from \ref{215Be2} that these isomorphisms satisfies the desired properties.  
Finally, we still denote by 
$\tau _{f,f '}$
the composition
$f  ^{\prime !} _{T _Y} \FF = 
\smash{\D} ^{\dag} _{\X \overset{f'}{\rightarrow}\Y}(\hdag T _Y) _{ \Q}      \otimes ^{\L} _{f _0 ^{-1}\smash{\D} ^{\dag} _{\Y} (\hdag T _Y) _{ \Q}  }
f _0 ^{-1} \FF [\delta _{X/Y}]
\underset{\tau _{f,f'} \otimes ^{\L} id}{\riso}
\smash{\D} ^{\dag} _{\X \overset{f}{\rightarrow}\Y} (\hdag T _Y) _{ \Q}     \otimes ^{\L} _{f _0 ^{-1}\smash{\D} ^{\dag} _{\Y} (\hdag T _Y) _{ \Q} }
f _0 ^{-1} \FF [\delta _{X/Y}]
=
f _{T _Y} ^{!} \FF$.
They also satisfy the desired properties.
\end{proof}

\begin{empt}
\label{2.1.5Be2-empt}
We keep notation \ref{prop-glueiniso-coh}.

a) For any  $\smash{\D} ^{\dag} _{\Y} (\hdag T _Y) _{\Q}$-module $\G$, 
we set 
$$ f ^{\dag *} _{T _Y} (\G)  := 
\smash{\D} ^{\dag} _{\X \overset{f}{\rightarrow}\Y} (\hdag T _Y) _{ \Q}     \otimes  _{f _0 ^{-1}\smash{\D} ^{\dag} _{\Y} (\hdag T _Y) _{ \Q} }
f _0 ^{-1} \G.$$
Similarly to \ref{prop-glueiniso-coh},
we construct isomorphisms
$\tau _{f,f'}\colon 
f ^{\prime \dag *} _{T _Y}  (\G) \riso 
f ^{\dag *} _{T _Y}(  \G)$
functorial in $\G$ and 
such that such that $\tau _{f,f}=\mathrm{Id}$,
$\tau _{f,f''}= \tau _{f,f'} \circ \tau _{f',f''}$.
We have the isomorphism 
of functors
$D ^{\mathrm{b}} _{\mathrm{coh}} ( \smash{\D} ^{\dag} _{\Y} (\hdag T _Y) _{\Q} )
\to 
D ^{\mathrm{b}}  ( \smash{\D} ^{\dag} _{\X } (\hdag T _X ) _{\Q}  )$
of the form
$ f _{T _Y} ^! \riso \L 
f ^{\dag *} _{T _Y} [\delta _{X/Y}]$.

b) Suppose $f$ is finite.
Then using \cite[3.2.4]{Be1},
we check that the canonical morphism
$$\widetilde{\B} ^{(m)} _{\X} ( T _X)
\otimes _{f ^{-1} \widetilde{\B} ^{(m)} _{\Y} ( T _Y) }f _0 ^{-1} \smash{\widetilde{\D}} _{\Y / \fS } ^{(m)} (T _Y)
\to
\smash{\widetilde{\D}} _{\X \overset{f}{\rightarrow}\Y} ^{(m)} ( T _Y)$$
is an isomorphism.
Hence, so is the canonical morphism
$$
\O _{\X} (\hdag T _X) _{\Q} \otimes _{f _0 ^{-1} \O _{\Y} (\hdag T _Y) _{\Q} } f _0 ^{-1}\smash{\D} ^{\dag} _{\Y}  (\hdag T _Y) _{ \Q}  
\to
\smash{\D} ^{\dag} _{\X \overset{f}{\rightarrow}\Y}  (\hdag T _Y) _{ \Q} . $$
Tensoring by $\Q$ and taking the inductive limit over the level, 
this yields the canonical morphism
$$f _{T _Y} ^* (\G) := \O _{\X} (\hdag T _X) _{\Q} \otimes _{f _0 ^{-1} \O _{\Y} (\hdag T _Y) _{\Q} } f _0 ^{-1}\G 
\to 
f ^{\dag *} _{T _Y} (\G)  $$
is an isomorphism.
Hence, if 
$\FF\in D ^{\mathrm{b}} _{\mathrm{coh}} ( \smash{\D} ^{\dag} _{\Y} (\hdag T _Y) _{\Q} )$ 
has a resolution $\cP$ by $\smash{\D} ^{\dag} _{\Y} (\hdag T _Y) _{\Q}$-modules  which are 
$\O _{\Y} (\hdag T _Y) _{\Q}$-flat , then we get the isomorphism
$f ^{\dag*} _{T _Y} (\FF ) 
\riso 
\L f ^{ *} _{T _Y} (\FF) $.
\end{empt}

\begin{rem}
\label{rem-flat-resol-tau}
Let $\FF\in D ^{\mathrm{b}} _{\mathrm{coh}} ( \smash{\D} ^{\dag} _{\Y} (\hdag T _Y) _{\Q} )$.
\begin{enumerate}[(a)]
\item Suppose $\FF$ has a resolution $\cP$ by flat coherent $\smash{\D} ^{\dag} _{\Y} (\hdag T) _{\Q}$-modules.
Via $ f _{T _Y} ^! (\FF) \riso 
f ^{\dag *} _{T _Y} ( \cP) [\delta _{X/Y}]$
and 
$f _{T _Y} ^{\prime !} (\FF) \riso 
f ^{\prime \dag *} _{T _Y} ( \cP) [\delta _{X/Y}]$
(see \ref{2.1.5Be2-empt}), 
the isomorphism 
$\tau _{f,f'}\colon 
f  _{T _Y} ^{\prime !} \FF \riso f _{T _Y} ^{!} \FF$ is the same (up to the shift $[\delta _{X/Y}]$)
than 
that 
$\tau _{f,f'}\colon 
f ^{\prime \dag *} _{T _Y} ( \cP)  \riso f ^{\dag *} _{T _Y} ( \cP)$,
which is 
computed term by term.

\item  Suppose $\FF$ has a resolution $\cP$ by coherent $\smash{\D} ^{\dag} _{\Y} (\hdag T) _{\Q}$-modules which are 
$\O _{\Y} (\hdag T) _{\Q}$-flat and suppose 
$f$ and $g$ are finite morphisms.
Via $ f _{T _Y} ^! (\FF) \riso 
f ^{ *} _{T _Y} ( \cP) [\delta _{X/Y}]$
and 
$f _{T _Y} ^{\prime !} (\FF) \riso 
f ^{\prime *} _{T _Y} ( \cP) [\delta _{X/Y}]$
(see \ref{2.1.5Be2-empt}), 
the isomorphism 
$\tau _{f,f'}\colon 
f  _{T _Y} ^{\prime !} \FF \riso f _{T _Y} ^{!} \FF$ is the same (up to the shift $[\delta _{X/Y}]$)
than 
that 
$\tau _{f,f'}\colon 
f ^{\prime *} _{T _Y} ( \cP)  \riso f ^{ *} _{T _Y} ( \cP)$,
which is 
computed term by term.

\end{enumerate}
\end{rem}

\begin{prop}
\label{comp-comp-adj-immf}
  Consider the following diagram in the category of 
   formal $\fS$-schemes locally of formal finite type
and having locally finite $p$-bases  over $\fS$:
\begin{equation}\label{deuxcarresadj}
  \xymatrix  @R=0,3cm {
  {\fP ^{\prime \prime }} \ar[r] ^g
  &
  {\fP ^{\prime }} \ar[r] ^f
  &
  {\fP }
  \\
  {\X ^{\prime \prime }} \ar[u]  ^{u''}\ar[r]^b
  &
  {\X ^{\prime }} \ar[u] ^-{u'} \ar[r]^a
  &
  {\X  ,} \ar[u] ^u
  }
\end{equation}
where $f$, $g$, $a$ and $b$ are  flat and have locally finite $p$-bases,
where $u$, $u'$ and $u''$ are some closed immersions. 
We suppose that the diagram \ref{deuxcarresadj} 
is commutative modulo $\pi$. 

\begin{enumerate}[(i)]
 \item We have the canonical adjunction morphism
\begin{equation}
\label{comp-comp-adj-immf-morp1}
u ^\prime _+ \circ a ^!\rightarrow f ^{!}\circ u _+
\end{equation}
of functors
  $D ^{\mathrm{b}} _{\mathrm{coh}} (\smash{\D} ^\dag _{\X ,\Q })
  \to 
  D ^{\mathrm{b}} _{\mathrm{coh}} (\smash{\D} ^\dag _{\fP '  ,\bbQ} )$.
If the right square of \ref{deuxcarresadj} 
is cartesian modulo $\pi$ then 
\ref{comp-comp-adj-immf-morp1}  is an isomorphism.

\item Denoting by 
$\phi \colon u ^\prime _+ \circ a ^! \rightarrow f ^{!}\circ u _+$,
(resp. 
$\phi '\ : \ u ^{''} _+ \circ b ^! \rightarrow g ^! \circ u ^\prime _+$,
resp. $\phi ''\ : \ u ^{''} _+ \circ (a \circ b ) ^! \rightarrow (f \circ g) ^!\circ u  _+$)
the morphism of adjunction of the right square 
 \ref{deuxcarresadj} (resp. the left square, resp. the outline of \ref{deuxcarresadj}),
then the following diagram
$$\xymatrix  @R=0,3cm {
{ u ''  _+ \circ (a\circ b) ^!}
\ar[r] _\sim
\ar[d] ^-{\phi ''}
&
{ u ''  _+ \circ b ^!\circ  a ^! }
\ar[d]^{( g ^! \circ \phi)\circ (\phi ' \circ a ^!)}
\\
{(f\circ g) ^! \circ  u _+  }
\ar[r] _-{\sim}
&
{ g ^!  \circ f ^!\circ  u _+,}
}
$$
is commutative. 
By abuse of notation, 
we get the transitivity equality
$\phi ''=( g ^! \circ \phi)\circ (\phi ' \circ a ^!)  $.

\item Let $a '$ : $ \X ^{\prime } \rightarrow \X $ (resp. $f'$ : $\fP ^{\prime } \rightarrow \fP $)
be a morphism whose reduction $X ^{\prime } \rightarrow \X $ (resp. $P ^{\prime } \rightarrow \fP $)
is equal to that of $a$ (resp. $f$). Then the following diagram
$$\xymatrix  @R=0,3cm {
{ u ^\prime_+ a ^{ !}}
\ar[r] ^-{\phi}
&
{f ^{!}\circ u _+}
\\
{ u' _+ a ^{\prime !} }
\ar[r]^{\psi}
\ar[u] ^-{ u ^\prime_+ (\tau _{a,a'})} _-{\sim}
&
{f ^{\prime !}\circ u  _+,}
\ar[u] ^-{\tau _{f,f'} u _+} _-{\sim}
}$$
where $\psi$ means the morphism of adjunction of the right square of  \ref{deuxcarresadj} whose 
$a$ and  $f$ have been replaced respectively by $a'$ and  $f'$,
is commutative.
\end{enumerate}

\end{prop}

\begin{proof}
We build \ref{comp-comp-adj-immf-morp1} using the adjoint pairs
$(u _+, u ^!)$ and $(u '_+, u ^{\prime !})$
(see \ref{Radj-u+flatf}).
If the right square of \ref{deuxcarresadj} 
is cartesian modulo $\pi$ then using Berthelot-Kashiwara's theorem \ref{exact-Berthelot-Kashiwara-full} (whose proof does not use 
\ref{comp-comp-adj-immf}), 
we check 
\ref{comp-comp-adj-immf-morp1}  is an isomorphism.
We proceed similarly to \cite[2.2.2]{caro-construction} to check the other properties.
\end{proof}

\subsection{Berthelot-Kashiwara theorem}
Let $u \colon \ZZ \hookrightarrow \X$ be a closed immersion of 
 formal $\fS$-schemes locally of formal finite type
and having locally finite $p$-bases  over $\fS$.
Let $\I$ be the ideal defining $u$. 
Let  $\Y$ be the open formal subscheme of $\X$ whose underlying topological space is complementary to that of $\ZZ$.
Let  $(\B ^{(m)}) _{m\in \N}$ be an inductive system of coherent separated complete (for the $p$-adic topology) commutative $\O _{\X}$-algebras. 
We suppose $\B ^{(m)}$ is endowed with a compatible structure of 
left $\D ^{(m)} _{\X/\fS}$-module (see the definition \ref{dfn-algcomp-mod}) such that the homomorphism of 
$\O _\X$-algebras
$\B ^{(m)}\to \B ^{(m+1)}$ is an monomorphism of 
$\D ^{(m)} _{\X/\fS}$-modules.  
We set 
$\widetilde{\D} ^{(m)} _{\X/\fS}: = 
\B ^{(m)} \widehat{\otimes} _{\O _\X}
\widehat{\D} ^{(m)} _{\X/\fS}$.
We suppose that the family 
$(u ^*\B ^{(m)}) _{m\in \N}$ satisfies the same properties over $\ZZ$. 
We set 
$\widetilde{\D} ^{(m)} _{\ZZ/\fS}: = 
(u ^*\B ^{(m)})\widehat{\otimes} _{\O _\ZZ}
\widehat{\D} ^{(m)} _{\ZZ/\fS}$.
We set $\smash{\widetilde{\D}} ^{\dag} _{\X/\fS,\Q}:=\underset{\underset{m}{\longrightarrow}}{\lim}\,
\smash{\widetilde{\D}} ^{(m)} _{\X/\fS,\Q}$
and
$\smash{\widetilde{\D}} ^{\dag} _{\ZZ/\fS,\Q}:=\underset{\underset{m}{\longrightarrow}}{\lim}\,
\smash{\widetilde{\D}} ^{(m)}  _{\ZZ/\fS,\Q}$.

\begin{thm}
[Berthelot]
\label{Berthelot-Kashiwara}
Let  $\E $ be a coherent $\widetilde{\D} ^{(m)} _{\X/\fS,\Q}$-module with support in $\ZZ$ (i.e. such that $\E | \Y = 0$).
Then, there exists a large enough integer $m' \geq m$, 
a coherent $\widetilde{\D} ^{(m')} _{\ZZ/\fS,\Q}$-module $\FF$,
and an isomorphism of $\widetilde{\D} ^{(m')} _{\X/\fS,\Q}$-modules of the form
\begin{equation}
\notag
u _{+} ^{(m')} (\FF) \riso 
\widetilde{\D} ^{(m')} _{\X/\fS,\Q} \otimes _{\widetilde{\D} ^{(m)} _{\X/\fS,\Q}} \E.
\end{equation}

\end{thm}

\begin{proof}
We can copy the proof of \cite[A.6]{caro-stab-sys-ind-surcoh}: 
since the theorem is local, 
using \ref{cor-closed-immer-local} we can suppose that 
 $\X$ is affine and there exists integers
 $n \geq r$ and a cartesian diagram of formal $\fS$-schemes
of the form:
$$\xymatrix{
{\X} 
\ar[r] ^-{}
\ar@{}[rd] ^-{}|\square
& {\widehat{\A} _\fS  ^{d} } 
\\ 
{\ZZ} 
\ar@{^{(}->}[u] ^-{u}
\ar[r] ^-{}
& 
{\widehat{\A} _\fS ^{r},} 
\ar@{^{(}->}[u] ^-{}
}
$$
where the horizontal arrows are relatively perfect, 
the right vertical arrow is given by 
the identification
$\widehat{\A} _\fS ^{r}
=
V (x _{r+1},\dots, x _d)$
if $x _1,\dots, x _d$ are the coordinates of $\widehat{\A} _\fS  ^{d} /\fS$.
\end{proof}

\begin{theo}
[Berthelot-Kashiwara]
\label{exact-Berthelot-Kashiwara-full}
Let $u \colon \ZZ  \to \X $ be a closed immersion of 
formal  $\fS$-schemes locally of formal finite type
and having locally finite $p$-bases  over $\fS$.
Let $D$ be a divisor of $X$ such that $Z\cap D $ is a divisor of $Z$.

The functors $u ^{!}$ and $u _{+}$ induce quasi-inverse equivalences between the category of 
coherent $\D ^{\dag} _{\fX /\fS } (\hdag D) _{\Q}$-modules with support in $Z$ 
and that of coherent 
$\D ^{\dag} _{\fZ } (\hdag D \cap Z) _{\Q}$-modules.
These functors $u ^{!}$ and  $u  _{+}$ are exact over these categories. 

\end{theo}

\begin{proof}
We can copy word by word the proof of \cite[A.8]{caro-stab-sys-ind-surcoh}.
\end{proof}

\begin{rem}
\label{rem-exact-Berthelot-Kashiwara-full}
With notation \ref{exact-Berthelot-Kashiwara-full}, 
by copying the proof of \cite[A.8]{caro-stab-sys-ind-surcoh},
we check more precisely that the adjunction morphism of 
$u _+ u ^! (\E) \to \E$ (see \ref{u+uflat-lem-form}) is an isomorphism for any 
coherent $\D ^{\dag} _{\fX /\fS,\Q}$-module $\E$ with support in $Z$. 
\end{rem}

\begin{thm}
[Inductive system version of Berthelot-Kashiwara's theorem]
\label{u!u+=id}
We keep notation \ref{exact-Berthelot-Kashiwara-full}.
Set $\fY: =\X \setminus Z$. 
Let 
$\FF ^{(\bullet)} 
\in \smash{\underrightarrow{LD}}  ^\mathrm{b} _{\Q, \mathrm{coh}}
(\overset{^\mathrm{l}}{} \smash{\widehat{\D}} _{\ZZ /\fS } ^{(\bullet)} (D \cap Z))$,
$\E ^{(\bullet)}  
\in \smash{\underrightarrow{LD}}  ^\mathrm{b} _{\Q, \mathrm{coh}}
(\overset{^\mathrm{l}}{} \smash{\widehat{\D}} _{\X /\fS } ^{(\bullet)} (D ))$
such that
$\E ^{(\bullet)} |\fY \riso 0$
in $\smash{\underrightarrow{LD}}  ^\mathrm{b} _{\Q, \mathrm{coh}}
(\overset{^\mathrm{l}}{} \smash{\widehat{\D}} _{\X /\fS  } ^{(\bullet)} (D ))$. 

\begin{enumerate}[(a)]
\item We have the canonical isomorphism in 
$\smash{\underrightarrow{LD}}  ^\mathrm{b} _{\Q, \mathrm{coh}}
(\overset{^\mathrm{l}}{} \smash{\widehat{\D}} _{\ZZ /\fS  } ^{(\bullet)} (D \cap Z))$ of the form:
\begin{equation}
\label{u!u+=id-iso}
u ^{ !(\bullet)} \circ u _{+} ^{ (\bullet)} (\FF ^{(\bullet)})
\riso \FF ^{(\bullet)}.
\end{equation}

\item We have $u ^{ !(\bullet)} (\E ^{(\bullet)} )\in \smash{\underrightarrow{LD}}  ^\mathrm{b} _{\Q, \mathrm{coh}}
(\overset{^\mathrm{l}}{} \smash{\widehat{\D}} _{\ZZ /\fS } ^{(\bullet)} (D \cap Z))$ and 
we have the canonical isomorphism :
\begin{equation}
\label{u!u+=id-isobis}
u _{+} ^{ (\bullet)} \circ  u ^{ !(\bullet)}  (\E ^{(\bullet)})
\riso \E ^{(\bullet)}.
\end{equation}

\item 
\label{u!u+=id-eq-cat}
The functors $u _{+} ^{ (\bullet)} $ and $u ^{ !(\bullet)}  $ induce t-exact quasi-inverse equivalences of categories 
between 
\begin{enumerate}[(a)]
\item $\smash{\underrightarrow{LD}}  ^\mathrm{b} _{\Q, \mathrm{coh}}
(\overset{^\mathrm{l}}{} \smash{\widehat{\D}} _{\ZZ /\fS } ^{(\bullet)} (D \cap Z))$
(resp. 
$\smash{\underrightarrow{LD}}  ^0 _{\Q, \mathrm{coh}}
(\overset{^\mathrm{l}}{} \smash{\widehat{\D}} _{\ZZ /\fS } ^{(\bullet)} (D \cap Z))$)

\item and the subcategory of 
$\smash{\underrightarrow{LD}}  ^\mathrm{b} _{\Q, \mathrm{coh}}
(\overset{^\mathrm{l}}{} \smash{\widehat{\D}} _{\X /\fS } ^{(\bullet)} (D ))$
(resp. $\smash{\underrightarrow{LD}}  ^0 _{\Q, \mathrm{coh}}
(\overset{^\mathrm{l}}{} \smash{\widehat{\D}} _{\X /\fS } ^{(\bullet)} (D ))$) 
of complexes 
$\E ^{(\bullet)} $ so that 
$\E ^{(\bullet)} |\fY \riso 0$.

\end{enumerate}
\end{enumerate}

\end{thm}

\begin{proof}
Using Theorem \ref{exact-Berthelot-Kashiwara-full}, 
we can copy the proof of  \cite[5.3.7]{caro-stab-sys-ind-surcoh}.
\end{proof}

\begin{cor}
\label{u+-com-otimes}
We keep notation \ref{exact-Berthelot-Kashiwara-full}.
Let 
$\cF ^{(\bullet)} ,
\cG ^{(\bullet)} 
\in \smash{\underrightarrow{LD}}  ^\mathrm{b} _{\Q, \mathrm{coh}}
(\overset{^\mathrm{l}}{} \smash{\widehat{\D}} _{\fZ /\fS } ^{(\bullet)} (D \cap Z))$.
We have the canonical isomorphism in 
$\smash{\underrightarrow{LD}}  ^\mathrm{b} _{\Q, \mathrm{coh}}
(\overset{^\mathrm{l}}{} \smash{\widehat{\D}} _{\ZZ /\fS  } ^{(\bullet)} (D \cap Z))$ of the form:
\begin{equation}
\label{u!u+=id-iso-coro}
u _{D,+} ^{(\bullet)} ( \cG ^{ (\bullet)} )
\smash{\widehat{\otimes}}^\L 
_{\widetilde{\B} ^{(\bullet)}  _{\fX} ( D) } 
u _{D,+} ^{(\bullet)} ( \cF ^{(\bullet)} ) [\delta _{Z/X}]
\riso 
u _{D,+} ^{(\bullet)} 
\left ( 
\cG ^{ (\bullet)} 
\smash{\widehat{\otimes}}^\L 
_{\widetilde{\B} ^{(\bullet)}  _{\fZ} (D \cap Z) } 
\cF ^{ (\bullet)} 
\right ) .
\end{equation}

\end{cor}

\begin{proof}
We apply the projection isomorphism
\ref{surcoh2.1.4} and we use the isomorphism \ref{u!u+=id-iso}.
\end{proof}

\subsection{Coherent arithmetic $\D$-modules over a realizable scheme having locally finite $p$-bases}
\label{ntnPPalpha}

Let $\fP$ be a formal  $\fS$-scheme locally of formal finite type
and having locally finite $p$-bases  over $\fS$.
Let $u _0\colon X  \to P $ be a closed immersion of schemes locally of formal finite type
and having locally finite $p$-bases over $S $.

Let $(\fP  _{\alpha}) _{\alpha \in \Lambda}$ be an open covering of  $\fP $.
We set $\fP  _{\alpha \beta}:= \fP  _\alpha \cap \fP  _\beta$,
$\fP  _{\alpha \beta \gamma}:= \fP  _\alpha \cap \fP  _\beta \cap \fP  _\gamma$,
$X  _\alpha := X  \cap P  _\alpha$,
$X _{\alpha \beta } := X  _\alpha \cap X  _\beta$ and
$X _{\alpha \beta \gamma } := X  _\alpha \cap X  _\beta \cap X  _\gamma $.

We suppose the covering $(\fP  _{\alpha}) _{\alpha \in \Lambda}$ satisfies the following lifting properties 
(such coverings exist following : see example \ref{section-BK-var-ex}).
For any 3uple $(\alpha, \, \beta,\, \gamma)\in \Lambda ^3$, we suppose there exists
$\X  _\alpha$ (resp. $\X  _{\alpha \beta}$, $\X  _{\alpha \beta \gamma}$)
some lifting of  $X  _\alpha$ (resp. $X  _{\alpha \beta}$, $X  _{\alpha \beta \gamma}$) 
which is a formal $\fS $-scheme
locally of formal finite type
and having finite $p$-bases,
$p _1 ^{\alpha \beta}$ :
$\X   _{\alpha \beta} \rightarrow \X  _{\alpha}$
(resp. $p _2 ^{\alpha \beta}$ :
$\X   _{\alpha \beta} \rightarrow \X  _{\beta}$)
some flat lifting 
of 
$X   _{\alpha \beta} \rightarrow X  _{\alpha}$
(resp. $X   _{\alpha \beta} \rightarrow X  _{\beta}$).
Similarly, for any $(\alpha,\,\beta,\,\gamma )\in \Lambda ^3$, fix some lifting 
$p _{12} ^{\alpha \beta \gamma}$ : $\X   _{\alpha \beta \gamma} \rightarrow \X   _{\alpha \beta} $,
$p _{23} ^{\alpha \beta \gamma}$ : $\X   _{\alpha \beta \gamma} \rightarrow \X   _{\beta \gamma} $,
$p _{13} ^{\alpha \beta \gamma}$ : $\X   _{\alpha \beta \gamma} \rightarrow \X   _{\alpha \gamma} $,
$p _1 ^{\alpha \beta \gamma}$ : $\X   _{\alpha \beta \gamma} \rightarrow \X   _{\alpha} $,
$p _2 ^{\alpha \beta \gamma}$ : $\X   _{\alpha \beta \gamma} \rightarrow \X   _{\beta} $,
$p _3 ^{\alpha \beta \gamma}$ : $\X   _{\alpha \beta \gamma} \rightarrow \X   _{\gamma} $,
$u _{\alpha}$ : $\X  _{\alpha } \hookrightarrow \fP  _{\alpha }$,
$u _{\alpha \beta}$ : $\X  _{\alpha \beta} \hookrightarrow \fP  _{\alpha \beta}$
and
$u _{\alpha \beta \gamma}$ : $\X  _{\alpha \beta \gamma } \hookrightarrow \fP  _{\alpha \beta \gamma}$.

\begin{ex}
\label{section-BK-var-ex}
Using \ref{cor-closed-immer-local}, 
we can choose such covering $(\fP  _{\alpha}) _{\alpha \in \Lambda}$
so that for every $\alpha\in \Lambda$, $X  _\alpha$ is affine and has finite $p$-basis.
Since $P$ is separated (recall our convention at the beginning of the paper), for any $\alpha,\beta ,\gamma \in \Lambda$,
$X _{\alpha \beta }$ and  $X _{\alpha \beta \gamma }$ are also affine and have finite $p$-bases.
Hence,  following \ref{lifting-pbasis}.\ref{lifting-pbasis-p2} such liftings exists.
\end{ex}

\begin{dfn}\label{defindonnederecol}
For any $\alpha \in \Lambda$, let $\E _\alpha$ be a coherent
$\D ^{\dag} _{\X  _{\alpha}\Q} $-module.
A \textit{glueing data} on $(\E _{\alpha})_{\alpha \in \Lambda}$
is the data for any $\alpha,\,\beta \in \Lambda$ of a
$\D ^{\dag} _{\X  _{\alpha \beta}\Q} $-linear
isomorphism
$$ \theta _{  \alpha \beta} \ : \  p _2  ^{\alpha \beta !} (\E _{\beta}) \riso p  _1 ^{\alpha \beta !} (\E _{\alpha}),$$
satisfying the cocycle condition:
$\theta _{13} ^{\alpha \beta \gamma }=
\theta _{12} ^{\alpha \beta \gamma }
\circ
\theta _{23} ^{\alpha \beta \gamma }$,
where $\theta _{12} ^{\alpha \beta \gamma }$, $\theta _{23} ^{\alpha \beta \gamma }$
and $\theta _{13} ^{\alpha \beta \gamma }$ are the isomorphisms making commutative the following diagram
\begin{equation}
  \label{diag1-defindonnederecol}
\xymatrix  @R=0,3cm {
{  p _{12} ^{\alpha \beta \gamma !} p  _2 ^{\alpha \beta !}  (\E _\beta )}
\ar[r] ^-{\tau} _-{\sim}
\ar[d] ^-{p _{12} ^{\alpha \beta \gamma !} (\theta _{\alpha \beta})} _-{\sim}
&
{p _2 ^{\alpha \beta \gamma!}  (\E _\beta )}
\ar@{.>}[d] ^-{\theta _{12} ^{\alpha \beta \gamma }}
\\
{ p _{12} ^{\alpha \beta \gamma !}  p  _1 ^{\alpha \beta !}  (\E _\alpha)}
\ar[r]^{\tau} _-{\sim}
&
{p _1 ^{\alpha \beta \gamma!}(\E _\alpha),}
}
\xymatrix  @R=0,3cm {
{  p _{23} ^{\alpha \beta \gamma !} p  _2 ^{\beta \gamma!}  (\E _\gamma )}
\ar[r] ^-{\tau} _-{\sim}
\ar[d] ^-{p _{23} ^{\alpha \beta \gamma !} (\theta _{ \beta \gamma})} _-{\sim}
&
{p _3 ^{\alpha \beta \gamma!}  (\E _\gamma )}
\ar@{.>}[d] ^-{\theta _{23} ^{\alpha \beta \gamma }}
\\
{ p _{23} ^{\alpha \beta \gamma !}  p  _1 ^{ \beta \gamma !}  (\E _\beta)}
\ar[r]^{\tau} _-{\sim}
&
{p _2 ^{\alpha \beta \gamma!}(\E _\beta),}
}
\xymatrix  @R=0,3cm {
{  p _{13} ^{\alpha \beta \gamma !} p  _2 ^{\alpha \gamma !}  (\E _\gamma )}
\ar[r] ^-{\tau} _-{\sim}
\ar[d] ^-{p _{13} ^{\alpha \beta \gamma !} (\theta _{\alpha \gamma})} _-{\sim}
&
{p _3 ^{\alpha \beta \gamma!}  (\E _\gamma )}
\ar@{.>}[d]^{\theta _{13} ^{\alpha \beta \gamma }}
\\
{ p _{13} ^{\alpha \beta \gamma !}  p  _1 ^{\alpha \gamma !}  (\E _\alpha)}
\ar[r]^{\tau} _-{\sim}
&
{p _1 ^{\alpha \beta \gamma!}(\E _\alpha),}
}
\end{equation}
where $\tau$ are the glueing isomorphisms defined in \ref{prop-glueiniso-coh1}.

\end{dfn}
\begin{dfn}
We define the category $\mathrm{Coh} ((\X   _\alpha )_{\alpha \in \Lambda}/K)$ as follows: 

\begin{itemize}
\item [-] an object is a family $(\E _\alpha) _{\alpha \in \Lambda}$
of coherent  $\D ^{\dag} _{\X  _{\alpha}\Q} $-modules
together with a glueing data $ (\theta _{\alpha\beta}) _{\alpha ,\beta \in \Lambda}$,

\item [-] a morphism
$((\E _{\alpha})_{\alpha \in \Lambda},\, (\theta _{\alpha\beta}) _{\alpha ,\beta \in \Lambda})
\rightarrow
((\E ' _{\alpha})_{\alpha \in \Lambda},\, (\theta '_{\alpha\beta}) _{\alpha ,\beta \in \Lambda})$
is a familly of morphisms $f _\alpha$ : $\E _\alpha \rightarrow \E '_\alpha$
of coherent  $\D ^{\dag} _{\X  _{\alpha}\Q} $-modules
commuting with glueing data, i.e., such that the following diagrams are commutative : 
\begin{equation}
  \label{diag2-defindonnederecol}
\xymatrix  @R=0,3cm {
{ p _2  ^{\alpha \beta !} (\E _{\beta}) }
\ar[d] _-{p _2  ^{\alpha \beta !} (f _{\beta}) }
\ar[r] ^-{\theta _{\alpha\beta}} _-{\sim}
&
{  p  _1 ^{\alpha \beta !} (\E _{\alpha}) }
\ar[d] ^-{p  _1 ^{\alpha \beta !} (f _{\alpha})}
\\
{p _2  ^{\alpha \beta !} (\E '_{\beta})  }
\ar[r]^{\theta '_{\alpha\beta}} _-{\sim}
&
{ p  _1 ^{\alpha \beta !} (\E '_{\alpha})  .}
}
\end{equation}

\end{itemize}
\end{dfn}

\begin{rem}
\label{rem-ind-ariYXP}
We can consider the category $\mathrm{Coh} ((\X   _\alpha )_{\alpha \in \Lambda}/K)$
as the category of arithmetic $\D$-modules over $(X,\fP)/\V$ or over $X/\V$ (we can check that, up to canonical equivalence of categories,
this is independent of 
the choice of the closed immersion $X \hookrightarrow \fP$ and of the liftings $\X _\alpha$ etc.).
\end{rem}

\begin{thm}
\label{prop1}
We denote by 
$\mathrm{Coh} (X, \fP /K)$ the category 
of 
coherent $\D ^{\dag} _{\fP\Q} $-modules with support in $X$. 
We have the following properties. 
\begin{enumerate}[(a)]
\item 
There exists
a canonical functor 
\begin{equation}
\label{const-u0!}
u _0 ^! \colon 
\mathrm{Coh} (X, \fP /K)
\to 
\mathrm{Coh} ((\X   _\alpha )_{\alpha \in \Lambda}/K)
\end{equation}
extending the usual functor $u _0 ^!$ when $X$ has a lifting formal $\fS$-scheme locally of formal finite type
and having locally finite $p$-bases.

\item There exists a canonical functor  
\begin{equation}
\label{const-u0+}
u _{0+}
\colon 
\mathrm{Coh} ((\X   _\alpha )_{\alpha \in \Lambda}/K)
\rightarrow
\mathrm{Coh} (X,\, \fP  /K)
\end{equation}
extending the usual functor $u _{0+}$ when $X$ can lift to a formal $\fS$-scheme locally of formal finite type
and having locally finite $p$-bases.

\item The functors $u ^!  _0$ and $ u _{0+}$ constructed in respectively \ref{const-u0!} and \ref{const-u0+} 
are quasi-inverse equivalences of categories between 
$\mathrm{Coh} ((\X   _\alpha )_{\alpha \in \Lambda}/K)$
and
$\mathrm{Coh} (X,\, \fP  /K)$.

\end{enumerate}

\end{thm}

\begin{proof}
We can copy the proof of \cite[5.3.7]{caro-6operations}.
\end{proof}

\section{Convergent isocrystals and arithmetic $\D$-modules}

\subsection{Convergent isocrystals on formal schemes having locally finite $p$-bases}
Let $\fY$ be a formal  $\fS$-scheme locally of formal finite type
and having locally finite $p$-bases  over $\fS$.
The results of \cite[6.2]{caro-6operations} are still true in our context. 
We collect them  below  concerning convergent isocrystals.

\begin{ntn}
\label{ntnMICdag2fs}
Let $\mathrm{MIC} ^{\dag \dag} (\Y/K)$ 
be the full subcategory of
$\D ^\dag _{\Y, \Q} $-modules
consisting of 
$\D ^\dag _{\Y, \Q} $-modules
which are 
$\O _{\Y,\Q} $-coherent.
\end{ntn}

\begin{thm}
[Berthelot]
\label{thm-eqcat-cvisoc}
\begin{enumerate}[(a)]
\item \label{thm-eqcat-cvisoc1}
Let $\E \in \mathrm{MIC} ^{\dag \dag} (\Y/K)$.
Then $\E$ 
is $\D ^\dag _{\Y, \Q} $-coherent
and 
$\O _{\Y,\Q} $-locally projective of finite type.

\item 
\label{thm-eqcat-cvisoc2}
Let $\E$ be a coherent 
$\D ^\dag _{\Y, \Q} $-module
which is $\O _{\Y,\Q} $-locally projective of finite type. We have the following properties.
\begin{enumerate}[(a)]
\item For any $m \in \N$, there exists 
a (coherent) $\widehat{\D} ^{(m)} _{\Y} $-module
$\overset{\circ}{\E}$, coherent over $\O _\Y$ together with 
an isomorphism of $\widehat{\D} ^{(m)} _{\Y, \Q} $-modules
$\overset{\circ}{\E} _\Q \riso \E$.
\item 
\label{thm-eqcat-cvisoc-arrows}
The module $\E$ is $\D  _{\Y, \Q} $-coherent 
 and for any $m \in \N$
the canonical homomorphisms 
\begin{gather}
\notag
\E 
\to 
\widehat{\D} ^{(m)} _{\Y, \Q}  \otimes _{\D  _{\Y, \Q} }
\E
,
\
\
\E
\to 
\D ^\dag _{\Y, \Q} 
 \otimes _{\widehat{\D} ^{(m)} _{\Y, \Q} }
\E
\end{gather}
are isomorphisms.

\end{enumerate}

\end{enumerate}

\end{thm}

\begin{proof}
We copy \cite[4.1.4]{Be1} and \cite[3.1.2 and 3.1.4]{Be0}.
\end{proof}

\begin{empt}
\label{empt-projff}
Let $\E \in \mathrm{MIC} ^{\dag\dag} (\fY/K) $.
Since  $\E $ is a locally projective $\O _{\Y,\Q}$-module of finite type (see \ref{thm-eqcat-cvisoc}),
then we have the following property: 
 $\E =0$ if and only if there exists an open dense subset $\U$ of $\fY$
such that $\E | \U =0$.
\end{empt}

\begin{lem}
Let $\E ^{(m)}$ be a coherent $\widehat{\D} ^{(m)} _{\Y/\fS, \Q}  $-module. 
For any $m' \geq m$, we set
$\E ^{(m')}:=
\widehat{\D} ^{(m')} _{\Y/\fS, \Q}   
\otimes _{\widehat{\D} ^{(m)} _{\Y/\fS, \Q}  }
\E ^{(m)}$,
and 
$\E := 
\D ^\dag _{\Y/\fS, \Q} \otimes _{\widehat{\D} ^{(m)} _{\Y/\fS, \Q}  }
\E ^{(m)}$.

If $\E$ is $\O _{\Y,\Q}$-coherent, then for $m'$ large enough
the canonical homomorphism 
$\E ^{(m')}
\to 
\E$
is an isomorphism.

\end{lem}

\begin{proof}
This is a consequence of Proposition \cite[3.6.2]{Be1} and of \ref{thm-eqcat-cvisoc}.2.
\end{proof}

\begin{lem}
\label{letterBerthelot-Lem3}
Let $\E$ be a coherent $\D ^\dag _{\Y/\fS, \Q} $-module which is 
$\O _{\Y,\Q}$-coherent, and 
$\overset{\circ}{\E}$ be a 
coherent $\widehat{\D} ^{(m)} _{\Y/\fS}$-module without $p$-torsion together with a 
$\widehat{\D} ^{(m)} _{\Y/\fS,\Q}$-linear isomorphism of the form
$\E \riso \overset{\circ}{\E} _\Q$.
Then $\overset{\circ}{\E}$ is $\O _{\Y}$-coherent, and this is a 
locally topologically nilpotent $\widehat{\D} ^{(m)} _{\Y/\fS}$-module.
\end{lem}

\begin{proof}
We can copy the proof of \cite[6.2.7]{caro-6operations}.
\end{proof}

\begin{prop}
\label{lem-projff}
Let $\E \in \mathrm{MIC} ^{\dag \dag} (\fY/K) $.
\begin{enumerate}[(a)]
\item 
\label{lem-projff-it1}
If $\fY$ is affine, then $\Gamma (\fY, \E)$ is a projectif 
$\Gamma (\fY, \O _{\fY,\Q} )$-module of finite type.
\item 
\label{lem-projff-it2}
The object $\E$
is a locally projective $\O _{\fY,\Q}$-module of finite type. 
\item 
\label{lem-projff-it3}
We have $\E =0$ if and only if there exists an open dense subset $\U $ of $\fY $
such that $\E | \U =0$.
\end{enumerate}
\end{prop}

\begin{proof}
We can copy the proof of \cite[6.2.9]{caro-6operations}.
\end{proof}

\begin{ntn}
\label{cvisoc-f*}
\begin{enumerate}[(a)]
\item 
\label{cvisoc-f*1}
Similarly to \ref{defi-M-L},
we denote by $M (\O _{\Y} ^{(\bullet)})$ the category of
$\O _{\Y} ^{(\bullet)}$-modules.
We get a canonical functor 
$\mathrm{cst}
\colon 
M (\O _{\Y} )
\to 
M (\O _{\Y} ^{(\bullet)})$ defined by 
$\FF \mapsto \FF ^{(\bullet)}$ so that 
$\FF ^{(m)}
\to 
\FF ^{(m+1)}$
is the identity of 
$\FF$.
Since this functor is exact, 
this yields the t-exact functor
$\mathrm{cst}
\colon 
D (\O _{\Y} )
\to 
D (\O _{\Y} ^{(\bullet)})$.
Similarly to \ref{defi-M-L},
we define the notion of ind-isogenies  (resp. of lim-ind-isogenies) of
$M (\O _{\Y} ^{(\bullet)})$.
Similarly to \ref{nota-(L)Mcoh},
we define the category
$\underrightarrow{LM} _{\Q, \mathrm{coh}} ( \O ^{(\bullet)} _{\Y})$.
We remark that 
$\underrightarrow{LM} _{\Q, \mathrm{coh}} ( \O ^{(\bullet)} _{\Y})$
is the  subcategory of
$\underrightarrow{LM} _{\Q} ( \O ^{(\bullet)} _{\Y})$
consisting of objects which are locally isomorphic to an object of the form 
$\mathrm{cst} (\G)$ where $\G$ is a coherent $\O _{\Y}$-module
(use analogous versions of \cite[2.1.7 and 2.2.2]{caro-stab-sys-ind-surcoh}).

\item 
\label{cvisoc-f*2}
 Following notation \ref{ntnMICdag2fs}, 
we denote by 
$\mathrm{MIC} ^{\dag \dag}(\Y/\V)$ the category of $\D ^\dag _{\Y/\fS, \Q} $-modules which are also
$\O _{\Y,\Q}$-coherent. Recall these objects are necessarily 
$\D ^\dag _{\Y, \Q} $-coherent,
and 
$\O _{\Y,\Q} $-locally projective of finite type.
We denote by 
$\mathrm{MIC}  ^{(\bullet)}(\Y/\V)$
the full subcategory of 
$\underrightarrow{LM} _{\Q, \mathrm{coh}} ( \widehat{\D} ^{(\bullet)} _{\Y/\fS})$
consisting of objects 
$\E ^{(\bullet)}$ such that 
$\underrightarrow{\lim} \E ^{(\bullet)}$
are 
$\O _{\Y,\Q}$-coherent.

\end{enumerate}

\end{ntn}

\begin{rem}
\label{remMICY/K}
Let $\E \in \mathrm{MIC} ^{\dag \dag}(\Y/\V)$.
Let 
$\widetilde{\D} := \D ^\dag _{\Y, \Q}$ or 
$\widetilde{\D} := \widehat{\D} ^{(m)} _{\Y, \Q}$. 
Let 
$\D := \D _{\Y, \Q}$ or 
$\D := \widehat{\D} ^{(0)} _{\Y, \Q}$.
By using the isomorphisms of  \ref{thm-eqcat-cvisoc}.\ref{thm-eqcat-cvisoc-arrows}, 
we check that both morphisms 
$\E 
\to 
\widetilde{\D} \otimes _{\D  }
\E
\to 
\E$
are isomorphisms.
This yields that the first morphism 
is in fact
$\widetilde{\D}$-linear. 
Hence, if $\FF$ is a $\widetilde{\D}$-module, 
then any 
$\D $-linear morphism
$\E \to \FF$ is necessarily 
$\widetilde{\D}$-linear. 
\end{rem}

\begin{lem}
\label{remMICY/K2}
Let 
$\FF ^{(m)}$ be a coherent $\widehat{\D} ^{(m)} _{\Y} $-module
et $f \colon \FF ^{(m)} \to \FF ^{(m)}$ be a $\V$-linear morphism
such that $f _\Q \colon \FF ^{(m)} _\Q\to \FF ^{(m)} _\Q$ is equal to $p ^N id$ for some $N \in \N$.
Then, for $N' \in \N$ large enough, we have
$p ^{N'} f = p ^{N'+N} id$.
\end{lem}

\begin{proof}
Since $\fY$ is quasi-compact and
$\FF ^{(m)}$ is a coherent $\widehat{\D} ^{(m)} _{\Y} $-module,
then
the $p$-torsion part of  
$\FF ^{(m)}$ is killed by some power of $p$.
Hence, we are done.
\end{proof}

\begin{prop}
\label{MICbullet-lem}
Let $\E \in \mathrm{MIC} ^{\dag \dag}(\Y/\V)$.
Let 
$\FF ^{(0)}$ be a $\widehat{\D} ^{(0)} _{\Y} $-module, coherent over $\O _\Y$ together with 
an isomorphism of $\widehat{\D} ^{(0)} _{\Y, \Q} $-modules of the form
$\FF ^{(0)} _\Q \riso \E$.
For any $m\in \N$, let $\G ^{(m)} $ be the quotient of 
$\widehat{\D} ^{(m)} _{\Y/\fS} \otimes _{\widehat{\D} ^{(0)} _{\Y/\fS}} \FF ^{(0)}$ by its $p$-torsion part. 
The following conditions are satisfied. 
\begin{enumerate}[(a)]

\item 
The module $\G ^{(m) }$ 
is $\O _{\Y}$-coherent.

\item The first (resp. second) canonical morphism
$$\FF ^{(0)} \to 
\widehat{\D} ^{(m)} _{\Y} \otimes _{\widehat{\D} ^{(0)} _{\Y} } \FF ^{(0)}
\to \G ^{(m)}$$
is an isogeny in the category of
$\widehat{\D} ^{(0)} _{\Y} $-modules
(resp. of coherent $\widehat{\D} ^{(m)} _{\Y} $-modules).

\item $\widehat{\D} ^{(\bullet)} _{\Y/\fS} \otimes _{\widehat{\D} ^{(0)} _{\Y/\fS}} \FF ^{(0)}
\in 
\mathrm{MIC}  ^{(\bullet)}(\Y/\V)$
and 
$ \underrightarrow{\lim} 
\,
(\widehat{\D} ^{(\bullet)} _{\Y/\fS} \otimes _{\widehat{\D} ^{(0)} _{\Y/\fS}} \FF ^{(0)})
\riso 
\E$.

\end{enumerate}

\end{prop}

\begin{proof}
Thanks to 
\ref{letterBerthelot-Lem3},
\ref{remMICY/K2}
and \ref{remMICY/K},  
we can copy word by word 
the proof of \cite[6.2.14]{caro-6operations}.
\end{proof}

\begin{cor}
\label{MICbullet-prop}
Let $\E ^{(\bullet)}
\in
\underrightarrow{LM} _{\Q} ( \widehat{\D} ^{(\bullet)} _{\Y/\fS})$.
The object 
$\E ^{(\bullet)}$ belongs to 
$\mathrm{MIC}  ^{(\bullet)}(\Y/\V)$ if and only if the following condition is satisfied: 
There exists  
a $\widehat{\D} ^{(0)} _{\Y/\fS}$-module $\FF ^{(0)}$,
coherent over $\O _\Y$ 
such that 
$\widehat{\D} ^{(\bullet)} _{\Y/\fS} \otimes _{\widehat{\D} ^{(0)} _{\Y/\fS}} \FF ^{(0)}$
is isomorphic in 
$\underrightarrow{LM} _{\Q} ( \widehat{\D} ^{(\bullet)} _{\Y/\fS})$ to 
$\E ^{(\bullet)}$
and such that 
the canonical morphism 
$\mathrm{cst} (\FF ^{(0)})
\to 
\widehat{\D} ^{(\bullet)} _{\Y/\fS} \otimes _{\widehat{\D} ^{(0)} _{\Y/\fS}} \FF ^{(0)}$
is an ind-isogeny in  
$M ( \O ^{(\bullet)} _{\Y})$.
Moreover, when 
$\E ^{(\bullet)}
\in
\mathrm{MIC}  ^{(\bullet)}(\Y/\V)$,
we can choose such $\FF ^{(0)}$ without $p$-torsion.
\end{cor}

\begin{proof}
We can copy word by word 
the proof of \cite[6.2.15]{caro-6operations}.
\end{proof}

\begin{empt}
Let $f \colon \Y ' \to \Y $ be a morphism of 
formal  $\fS$-schemes locally of formal finite type
and having locally finite $p$-bases  over $\fS$.
Let 
$\E ^{(\bullet)}
\in 
M( \widehat{\D} ^{(\bullet)} _{\Y/\fS})$.
We set
$f ^{*(m)} _{\mathrm{alg}}(\E ^{(m)})
:=
\smash{\widehat{\D}} ^{(m)} _{\Y ^{\prime } \rightarrow \Y /\fS }
\otimes   _{f ^{-1} \smash{\widehat{\D}} ^{(m)} _{\Y /\fS } }
f ^{-1} \E ^{(m)}$.
We denote by 
$f ^{*(\bullet)} _{\mathrm{alg}}(\E ^{(\bullet)}):=
\smash{\widehat{\D}} ^{(\bullet)} _{\Y ^{\prime } \rightarrow \Y /\fS }
\otimes   _{f ^{-1} \smash{\widehat{\D}} ^{(\bullet)} _{\Y /\fS } }
f ^{-1} \E ^{(\bullet)}$
the object of 
$M( \widehat{\D} ^{(\bullet)} _{\Y'/\fS})$
whose transition morphisms are
$f ^{*(m)} _{\mathrm{alg}} ( \E ^{(m)})
\to 
f ^{*(m+1)} _{\mathrm{alg}} ( \E ^{(m+1)})$.
By left deriving the functor
$f ^{*(\bullet)} _{\mathrm{alg}}$, 
this yields the functor
$\L f ^{*(\bullet)} _{\mathrm{alg}}
\colon 
D ^- ( \widehat{\D} ^{(\bullet)} _{\Y/\fS})
\to 
D ^- ( \widehat{\D} ^{(\bullet)} _{\Y'/\fS})$, 
defined by setting 
$\L f ^{*(\bullet)} _{\mathrm{alg}}
(\cF ^{(\bullet)}):=
\smash{\widehat{\D}} ^{(\bullet)} _{\Y ^{\prime } \rightarrow \Y /\fS }
\otimes   ^\L _{f ^{-1} \smash{\widehat{\D}} ^{(\bullet)} _{\Y /\fS } }
f ^{-1} \cF ^{(\bullet)}$
for any
$\cF ^{(\bullet)} \in 
D ^- ( \widehat{\D} ^{(\bullet)} _{\Y/\fS})$.
Since it preserves lim-ind-isogenies, 
this induces the functor
$\L f ^{*(\bullet)} _{\mathrm{alg}}
\colon 
\underrightarrow{LD} _{\Q}  ^{-} (\widehat{\D} ^{(\bullet)} _{\Y/\fS} )
\to 
\underrightarrow{LD} _{\Q}  ^{-} (\widehat{\D} ^{(\bullet)} _{\Y'/\fS} )$.

Following notation \ref{ntn-Lf!+*},
we set 
$\L f ^{*(\bullet)} ( \cF ^{(\bullet)}) :=
\smash{\widehat{\D}} ^{(\bullet)} _{\Y ^{\prime } \rightarrow \Y /\fS }
\smash{\widehat{\otimes}} ^\L _{f ^{-1} \smash{\widehat{\D}} ^{(\bullet)} _{\Y /\fS } }
f ^{-1} \cF ^{(\bullet)}$,
for any
$\cF ^{(\bullet)} \in 
\underrightarrow{LD} _{\Q,\mathrm{qc}}  ^{\mathrm{b}} (\widehat{\D} ^{(\bullet)} _{\Y/\fS} )$.
Beware the notation is slightly misleading since $\L f ^{*(\bullet)}$
is not necessarily the left derived functor of a functor.
We get the morphism
$\L f ^{*(\bullet)} _{\mathrm{alg}} (\cF ^{(\bullet)})
\to 
\L f ^{*(\bullet)} ( \cF ^{(\bullet)})$

\end{empt}

\begin{lem}
\label{f*-OcohDm}
Let $f \colon \Y ' \to \Y $ be a morphism of 
formal  $\fS$-schemes locally of formal finite type
and having locally finite $p$-bases  over $\fS$.
We have the following properties.
\begin{enumerate}[(a)]
\item Let 
$\FF ^{(\bullet)} \in 
\underrightarrow{LD} _{\Q,\mathrm{qc}}  ^{\mathrm{b}} (\widehat{\D} ^{(\bullet)} _{\Y/\fS} )$.
The canonical morphism
$$\O _{\Y ^{\prime } }  ^{(\bullet)}
\smash{\widehat{\otimes}} ^\L _{f ^{-1}\O   ^{(\bullet)} _{\Y /\fS }}
f ^{-1}\FF ^{(\bullet)}
\to
\smash{\widehat{\D}} ^{(\bullet)} _{\Y ^{\prime } \rightarrow \Y /\fS }
\smash{\widehat{\otimes}} ^\L _{f ^{-1} \smash{\widehat{\D}} ^{(\bullet)} _{\Y /\fS } }
f ^{-1} \FF ^{(\bullet)}$$
is an isomorphism.

\item Let 
$\FF ^{(\bullet)} \in 
\underrightarrow{LD} _{\Q,\mathrm{coh}}  ^{\mathrm{b}} (\widehat{\D} ^{(\bullet)} _{\Y/\fS} )$.
The canonical morphism
$$
\L f ^{*(\bullet)} _{\mathrm{alg}} (\FF ^{(\bullet)}):=
\smash{\widehat{\D}} ^{(\bullet)} _{\Y ^{\prime } \rightarrow \Y /\fS }
\otimes ^\L _{f ^{-1} \smash{\widehat{\D}} ^{(\bullet)} _{\Y /\fS } }
f ^{-1} \FF ^{(\bullet)}
\to
\smash{\widehat{\D}} ^{(\bullet)} _{\Y ^{\prime } \rightarrow \Y /\fS }
\smash{\widehat{\otimes}} ^\L _{f ^{-1} \smash{\widehat{\D}} ^{(\bullet)} _{\Y /\fS } }
f ^{-1} \FF ^{(\bullet)}
=:\L f ^{*(\bullet)} (\FF ^{(\bullet)})$$
is an isomorphism.

\item 
\label{f*-OcohDm3}
Let 
$\G ^{(\bullet)} \in 
\underrightarrow{LD} _{\Q,\mathrm{coh}}  ^{\mathrm{b}} ( \O ^{(\bullet)} _{\Y})$.
Then, 
the canonical morphism
$$\O  ^{(\bullet)} _{\Y ^{\prime } }
\otimes  ^\L _{f ^{-1} \O   ^{(\bullet)} _{\Y /\fS } }
f ^{-1}\G ^{(\bullet)}
\to 
\O  ^{(\bullet)} _{\Y ^{\prime } }
\smash{\widehat{\otimes}} ^\L _{f ^{-1}\O  ^{(\bullet)} _{\Y /\fS }}
f ^{-1}\G ^{(\bullet)}$$
is an isomorphism of 
$\underrightarrow{LD} _{\Q,\mathrm{coh}}  ^{\mathrm{b}} ( \O ^{(\bullet)} _{\Y'})$.

\end{enumerate}
\end{lem}

\begin{proof}
This is left to the reader and easy (hint : to check 1) use \cite[2.3.5.2]{Be1}, and
the proof of 2) and 3) is identical to that of  \cite[3.4.2.2]{Beintro2}).
\end{proof}

\begin{prop}
\label{ntn-f*}
Let $f \colon \Y ' \to \Y $ be a morphism of 
formal  $\fS$-schemes locally of formal finite type
and having locally finite $p$-bases  over $\fS$.
\begin{enumerate}[(a)]
\item 
\label{ntn-f*1}
Let $\E \in \mathrm{MIC} ^{\dag \dag}(\Y/\V)$. Then the canonical last morphism
$$
\O _{\Y ^{\prime},\Q}
\otimes  _{ f ^{-1} \O _{\Y /\fS,\Q}} f ^{-1} \E
\liso 
\O _{\Y ^{\prime},\Q}
\otimes ^{\L} _{ f ^{-1} \O _{\Y /\fS,\Q}} f ^{-1} \E
\to
\D ^{\dag} _{\Y ^{\prime }\rightarrow \Y ,\Q}
\otimes ^{\L} _{ f ^{-1} \D ^{\dag} _{\Y /\fS,\Q}} f ^{-1} \E$$
is an isomorphism.
Hence, we can set $f ^* (\E):= \D ^{\dag} _{\Y ^{\prime }\rightarrow \Y ,\Q}
\otimes  _{ f ^{-1} \D ^{\dag} _{\Y /\fS,\Q}} f ^{-1} \E$ without ambiguity. 
We have also $f ^* (\E) \in \mathrm{MIC} ^{\dag \dag}(\Y'/\V)$.

\item \label{ntn-f*2}
Let 
$\FF $ be a $\widehat{\D} ^{(m)} _{\Y} $-module, coherent over $\O _\Y$.
Then the morphisms 
\begin{equation}
\notag
\O _{\Y ^{\prime}}
\otimes  _{ f ^{-1} \O _{\Y /\fS}} f ^{-1} \FF
\to 
\O _{\Y ^{\prime}}
\widehat{\otimes}  _{ f ^{-1} \O _{\Y /\fS}} f ^{-1} \FF
\to 
\smash{\widehat{\D}} ^{(m)} _{\Y ^{\prime } \rightarrow \Y /\fS }
\smash{\widehat{\otimes}}  _{f ^{-1} \smash{\widehat{\D}} ^{(m)} _{\Y /\fS } }
f ^{-1} \FF 
\leftarrow
\smash{\widehat{\D}} ^{(m)} _{\Y ^{\prime } \rightarrow \Y /\fS }
\otimes  _{f ^{-1} \smash{\widehat{\D}} ^{(m)} _{\Y /\fS } }
f ^{-1} \FF 
\end{equation}
are isomorphisms.
Hence, we can set $f ^* (\FF):= 
\smash{\widehat{\D}} ^{(m)} _{\Y ^{\prime } \rightarrow \Y /\fS }
\otimes  _{f ^{-1} \smash{\widehat{\D}} ^{(m)} _{\Y /\fS } }
f ^{-1} \FF $ 
without ambiguity.
Moreover, 
$f ^* (\FF)$ is  a $\widehat{\D} ^{(m)} _{\Y'} $-module, coherent over $\O _{\Y'}$.
\end{enumerate}

\end{prop}

\begin{proof}
We can copy word by word 
the proof of \cite[6.2.18]{caro-6operations}.
\end{proof}

\begin{prop}
\label{stab-MIC-f*}
Let $f \colon \Y ' \to \Y $ be a morphism of formal  $\fS$-schemes locally of formal finite type
and having locally finite $p$-bases  over $\fS$. 
Let $\FF ^{(0)}$ 
be a $\widehat{\D} ^{(0)} _{\Y/\fS}$-module,
coherent over $\O _\Y$ 
and such that 
the canonical morphism 
$\mathrm{cst} (\FF ^{(0)})
\to 
\widehat{\D} ^{(\bullet)} _{\Y/\fS} \otimes _{\widehat{\D} ^{(0)} _{\Y/\fS}} \FF ^{(0)}=:\FF ^{(\bullet)}$
is an ind-isogeny in  
$M ( \O ^{(\bullet)} _{\Y})$.
For any $m\in\N$,
let $\G ^{(m)} $ be the quotient of 
$\widehat{\D} ^{(m)} _{\Y/\fS} \otimes _{\widehat{\D} ^{(0)} _{\Y/\fS}} \FF ^{(0)}$ by its $p$-torsion part. 

\begin{enumerate}[(a)]

\item 
The canonical morphism 
$\mathrm{cst} (f ^* (\FF ^{(0)}))
\to 
\widehat{\D} ^{(\bullet)} _{\Y'/\fS} \otimes _{\widehat{\D} ^{(0)} _{\Y'/\fS}} f ^* (\FF ^{(0)})$
is an ind-isogeny of  
$M ( \O ^{(\bullet)} _{\Y'})$.

\item 
The canonical morphisms
$f ^{*(\bullet)} _{\mathrm{alg}} (\FF ^{(\bullet)})
\to
f ^{*(\bullet)} _{\mathrm{alg}} (\G ^{(\bullet)})$,
and 
$\widehat{\D} ^{(\bullet)} _{\Y'/\fS} \otimes _{\widehat{\D} ^{(0)} _{\Y'/\fS}} f ^* (\FF ^{(0)})
\to 
f ^{*(\bullet)} _{\mathrm{alg}} (\G ^{(\bullet)})$
are ind-isogenies of 
$M (\widehat{\D} ^{(\bullet)} _{\Y'/\fS} )$.
\item The canonical morphisms
$\L f ^{*(\bullet)} _{\mathrm{alg}} (\FF ^{(\bullet)})
\to 
\L f ^{*(\bullet)} ( \FF ^{(\bullet)})$
and 
$\L f ^{*(\bullet)} _{\mathrm{alg}} (\FF ^{(\bullet)}) 
\to 
 f ^{*(\bullet)} _{\mathrm{alg}} (\FF ^{(\bullet)})$ 
are isomorphisms of
$\underrightarrow{LD} ^\mathrm{b} _{\Q} ( \widehat{\D} ^{(\bullet)} _{\Y'/\fS})$.
\end{enumerate}
\end{prop}

\begin{proof}

We can copy word by word 
the proof of \cite[6.2.19]{caro-6operations}.
\end{proof}

\begin{coro}
\label{corostab-MIC-f*}
Let $f \colon \Y ' \to \Y $ be a morphism of formal  $\fS$-schemes locally of formal finite type
and having locally finite $p$-bases  over $\fS$. 
Let $\E ^{(\bullet)}
\in
\mathrm{MIC}  ^{(\bullet)}(\Y/\V)$, and
$\E := \underrightarrow{\lim}\, \E ^{(\bullet)}\in \mathrm{MIC} ^{\dag \dag}(\Y/\V)$.

\begin{enumerate}[(a)]

\item $\L f ^{*(\bullet)} ( \E ^{(\bullet)})
\in \mathrm{MIC}  ^{(\bullet)}(\Y'/\V)$ (i.e. is isomorphic to such an object)
and 
$ \underrightarrow{\lim} \L f ^{*(\bullet)} ( \E ^{(\bullet)})
\riso 
f ^* (\E)$.

\item 
Choose
a $\widehat{\D} ^{(0)} _{\Y/\fS}$-module $\FF ^{(0)}$,
coherent over $\O _\Y$ 
such that 
$\widehat{\D} ^{(\bullet)} _{\Y/\fS} \otimes _{\widehat{\D} ^{(0)} _{\Y/\fS}} \FF ^{(0)}$
is isomorphic in 
$\underrightarrow{LM} _{\Q} ( \widehat{\D} ^{(\bullet)} _{\Y/\fS})$ to 
$\E ^{(\bullet)}$
and such that 
the canonical morphism 
$\mathrm{cst} (\FF ^{(0)})
\to 
\widehat{\D} ^{(\bullet)} _{\Y/\fS} \otimes _{\widehat{\D} ^{(0)} _{\Y/\fS}} \FF ^{(0)}$
is an ind-isogeny in  
$M ( \O ^{(\bullet)} _{\Y})$.
Then 
$\L f ^{*(\bullet)} ( \E ^{(\bullet)})
\riso 
\widehat{\D} ^{(\bullet)} _{\Y'/\fS} \otimes _{\widehat{\D} ^{(0)} _{\Y'/\fS}} 
f ^* (\FF ^{(0)})$.

\end{enumerate}
\end{coro}

\subsection{Duality, inverse images on formal schemes having locally finite $p$-bases}
Let $\fX$ be a formal  $\fS$-scheme locally of formal finite type
and having locally finite $p$-bases  over $\fS$.
Let $\E \in \mathrm{MIC} ^{\dag \dag} (\X /K) $.
We have the equalities 
  $D ^\mathrm{b} _{\mathrm{coh}}(\O _{\X ,\Q}) =
  D _{\mathrm{parf}} (\O _{\X ,\Q})$,
  $D ^\mathrm{b} _{\mathrm{coh}}(\D _{\X /\fS, \Q }) =
  D _{\mathrm{parf}} (\D _{\X /\fS, \Q })$,
  and 
    $D ^\mathrm{b} _{\mathrm{coh}}(\D ^\dag  _{\X /\fS, \Q }) =
  D _{\mathrm{parf}} (\D ^\dag _{\X /\fS, \Q })$ (see \ref{finitecohdimDdagpre} and \ref{finitecohdimDdag}).
We get $\E \in D _{\mathrm{parf}} (\O _{\X ,\Q})$,
$\E \in D _{\mathrm{parf}} (\D  _{\X /\fS, \Q })$
and 
$\E \in D _{\mathrm{parf}} (\D ^\dag _{\X /\fS, \Q })$.

\begin{ntn}
For any $\FF \in D  ^\mathrm{b} _{\mathrm{coh}}
( \D  _{\X /\fS, \Q } ) $, we set
$\DD ^{\mathrm{alg}}  (\FF)=
\R \mathcal{H}om  _{\D  _{\X /\fS, \Q } }
( \FF , \,\D  _{\X /\fS, \Q } \otimes _{\O _{\X}}\omega _{\X /\fS  } ^{-1})[\delta _X]$ 
and
$\FF ^{\vee} =
\R \mathcal{H}om  _{\O _{\X ,\Q} }
( \FF , \,\O_{\X , \Q} ( \hdag Z ))$.
For any $\G \in D  ^\mathrm{b} _{\mathrm{coh}} ( \D ^\dag  _{\X /\fS, \Q } ) $, we set
$\DD   (\G)=
\R \mathcal{H}om  _{\D ^\dag _{\X /\fS, \Q } }
( \G , \,\D  ^\dag _{\X /\fS, \Q } \otimes _{\O _{\X}}\omega _{\X /\fS  } ^{-1})[\delta _X]$. 

\end{ntn}

\begin{prop}\label{Doevee=De}
There exists a canonical isomorphism 
$$ \theta \ : \ \DD ^{\mathrm{alg}}  ( \O _{\X ,\Q} ) \otimes ^\L _{\O _{\X ,\Q}} \E ^\vee
\riso \DD ^{\mathrm{alg}}  (\E).$$
\end{prop}
\begin{proof}
It is sufficient to copy  \cite[2.2.1]{caro_comparaison}. \end{proof}

\begin{lem}
\begin{enumerate}[(i)]
\item  $\O _{\X , \Q } \in D _{\mathrm{parf}}( \D  _{\X /\fS, \Q })$.

\item We have the canonical isomorphism:
\begin{equation}\label{dualB=B}
  \DD ^{\mathrm{alg}}  ( \O _{\X ,\Q})
\riso
\O _{\X ,\Q}.
\end{equation}

\end{enumerate}
\end{lem}
\begin{proof}
It is sufficient to copy  \cite[5.20]{caro_log-iso-hol}.
\end{proof}

\begin{rem}
\label{evee=De0}
From \ref{dualB=B} and \ref{Doevee=De},
we get the canonical isomorphism
$\E ^\vee \riso \DD ^{\mathrm{alg}}  (\E)$.
\end{rem}

\begin{empt}\label{rhoisom}
Consider the following morphism:
$$\rho ^\dag \ :\ \DD ^{\mathrm{alg}}  (\E) \rightarrow
\D _{\X ,\Q } ^{\dag} \otimes _{\D _{\X ,\Q}} \DD ^{\mathrm{alg}}  (\E)
\riso
\DD   ( \D _{\X ,\Q } ^{\dag} \otimes _{\D _{\X ,\Q}} \E)
\rightarrow \DD   (\E).$$
Since $\E$ is locally projective of finite type over $\O _{\X , \Q} $,
then the morphism
$ \mathcal{H}om _ {\O _{\X , \Q } } ( \E , \,\O _{\X , \Q }  )
\rightarrow
\R \mathcal{H}om _ {\O _{\X , \Q } } ( \E , \,\O _{\X , \Q }  )= \E ^\vee$ is 
an isomorphism.
This yields $\E ^\vee \in \mathrm{MIC} ^{\dag \dag} (\X /K) $.
Since $\E ^\vee \riso \DD ^{\mathrm{alg}}  (\E)$ (see \ref{evee=De0}), 
via \ref{thm-eqcat-cvisoc}.\ref{thm-eqcat-cvisoc-arrows}
we check that 
$\rho ^\dag $ is an isomorphism.
\end{empt}

\begin{empt}\label{constrhodag}
Let  $\theta ^\dag$ :
$\DD ^{\mathrm{alg}}   ( \O _{\X ,\Q} ) \otimes _{\O _{\X ,\Q}} \E ^\vee
\riso \DD ^{\mathrm{alg}}  (\E)$
be the isomorphism making commutative the following diagram:
$$\xymatrix {
{\DD ^{\mathrm{alg}}  ( \O _{\X ,\Q} ) \otimes _{\O _{\X ,\Q}} \E ^\vee}
\ar[r] ^(0.68){\theta} _(0.68){\sim}
\ar[d] ^{\rho ^\dag \otimes id} _{\sim}
&
{\DD ^{\mathrm{alg}}  (\E)}
\ar[d] ^{\rho ^\dag } _{\sim}
\\
{\DD   ( \O _{\X ,\Q} ) \otimes _{\O _{\X ,\Q}} \E ^\vee }
\ar@{.>}[r] ^(0.68){\theta ^\dag} _(0.68){\sim}
&
{\DD  (\E)}
}$$
\end{empt}

\begin{empt}
  \label{dualisoscvdag} 
  From  \ref{dualB=B} and \ref{rhoisom}, we get the isomorphism
  $\DD   ( \O _{\X ,\Q} ) \riso \O _{\X ,\Q}$.
Hence, the isomorphism $\theta ^\dag$ induces the following one
  $\E ^\vee \riso \DD  (\E)$. 
  Hence, we get the functor
  $ \DD \colon \mathrm{MIC} ^{\dag \dag} (\X/K) 
\to 
\mathrm{MIC} ^{\dag \dag} (\X/K) $.
\end{empt}

\begin{empt}
\label{sp*f*com-empt}
Let $u \colon \X' \to \X$ be a morphism of 
a formal  $\fS$-schemes locally of formal finite type
and having locally finite $p$-bases  over $\fS$.
Let $\E \in \mathrm{MIC} ^{\dag \dag} (\X/K) $ (see Notation \ref{ntnMICdag2fs}).
We have the functor 
$$u ^! [-\delta _{X'/X}] \colon \mathrm{MIC} ^{\dag \dag} (\X/K) 
\to 
\mathrm{MIC} ^{\dag \dag} (\X'/K) $$ 
which is compatible with $u ^*$, 
i.e. 
there exist a canonical isomorphism respectively of 
$\mathrm{MIC} ^{\dag \dag} (\X'/K) $
 of the form
\begin{equation}
\label{sp*f*com}
u  ^* (\E) \riso u ^! (\E) [-\delta _{X'/X}].
\end{equation}
Moreover, these isomorphisms are transitive with respect to the composition of morphisms (see \cite[2.4.1]{caro-construction}).
\end{empt}

\subsection{Direct image by the specialization morphism 
of the constant coefficient
when the boundary is not a divisor}
\label{directimagespovcv}
Let $\X$ be a (separated) formal $\fS$-scheme locally of formal finite type (see definition \ref{dfn-fft}).

\begin{empt}
[Cech complexes of the constant coefficient associated with divisors]
Let $\mathscr{T}:= (T _i) _{i\in I}$ be a finite set of divisor of $X$.
We can suppose $I = \{ 0,\dots, r\}$.
For each $h \in I$,
for any $i _0, \dots, i _h \in I$, put 
$T _{i _0,\dots, i _{h}}
:= 
T _{i _0}\cup \dots \cup T _{i _h}$.
For each $h \in I$, 
set 
\begin{equation}
\label{check dag h}
\check{C} ^{\dag h}  (\fX,\mathscr{T}, \O _{\fX,\bbQ})
:=
\prod _{ i _0<i _1 < \dots < i _h} 
\O _{\fX} 
(\hdag T _{i _0 \dots i _h})
_{\bbQ}.
\end{equation}
Let 
$\alpha \in  \check{C} ^{\dag h} (\fX,\mathscr{T}, \O _{\fX,\bbQ})$.
For any $h \in I$,
for any $i _0< \dots< i _h \in I$,
we denote by 
$\alpha_{ i _0, \dots,  i _h} $ the coefficient of $\alpha$
in $\O _{\fX} 
(\hdag T _{i _0 \dots i _h})
_{\bbQ}$.

We define the coboundary map 
$d \colon 
\check{C} ^{\dag h} (\fX,\mathscr{T}, \O _{\fX,\bbQ})
\to 
\check{C} ^{\dag h+1} (\fX,\mathscr{T}, \O _{\fX,\bbQ})$
by setting 
$$( d \alpha ) _{i _0, \dots, i _{h+1}}
:=
\sum _{j = 0} ^{h+1}
(-1) ^j 
\alpha _{i _0, \dots, \widehat{i} _j ,\dots, i _{h+1}}.$$

This yields the complex
$$
\cdots 0
\to 
\check{C} ^{\dag 0} (\fX,\mathscr{T}, \O _{\fX,\bbQ})
\to 
\check{C} ^{\dag 1} (\fX,\mathscr{T}, \O _{\fX,\bbQ})
\to 
\dots
\to 
\check{C} ^{\dag h} (\fX,\mathscr{T}, \O _{\fX,\bbQ})
\to 
0 \cdots 
$$
that we will denote by 
$\check{C} ^{\dag \bullet}(\fX,\mathscr{T}, \O _{\fX,\bbQ})$.

Let $Y _i := X \setminus T _i$ the open subscheme of $X$ et $Y := \cup _{i\in I} Y _i$.
We get the finite open covering  $\mathscr{Y}:= (Y _i) _{i=1,\dots, r}$ of $Y$. 
Since $\sp _* ( j _{Y _{i}}  ^{\dag} 
\O _{\X  ^{\ad}} )
\riso
\O _{\X   }( \hdag T _{i}) _\bbQ $, 
then 
\begin{equation}
\label{sp*Check=Check}
\sp _* \check{C} ^{\dag \bullet} (\fX,\mathscr{Y},\O _{\X  ^{\ad}} )
\riso 
\check{C} ^{\dag \bullet}(\fX, \mathscr{T}, \O _{\fX,\bbQ}),
\end{equation}
where 
$\check{C} ^{\dag \bullet} (\fX, \mathscr{Y},\O _{\X  ^{\ad}} )$
is defined in \ref{2.50Lazda-Pal-Book-ntn}.
\end{empt}

\begin{empt}
Let $\sp \colon \X  ^{\ad} \to \X$ be the specialization morphism.
Let $Y$ be an open subscheme of $X$.
Let $\mathscr{X}:= (\X _i) _{i\in I}$ be a finite affine covering of $\X$. 
For any $i\in I$, let $\mathscr{Y} _{i} :=
(Y _{i\,j _i}) _{j_i\in J _i}$ be a finite covering of $Y _i:= Y \cap \X _i$ such that 
there exists $f _{i\, j _i}\in \Gamma (\X _i, \O _{\X})$ satisfying
$Y _{i\,j _i} = D ( f _{i\,j _i})\cap X _i$.
We get the divisor 
$T _{i\,j _i} := V (f _{i\,j _i})$ 
of $X _i$ such that 
$Y _{i\,j _i} = X _i \setminus T _{i\,j _i}$.

Let $\underline{i} =(i_0,\dots, i _h) \in I ^{1+h}$. 
We set $\X _{\underline{i}}:= \X _{i _0}\cap \dots \cap \X _{i _h}$,
$Y _{\underline{i}}:= Y \cap \X _{\underline{i}}$,
$u _{\underline{i}} 
\colon 
\X _{\underline{i}}
\to 
\X$,
$u _{\underline{i}}  ^{\ad}
\colon 
\X _{\underline{i}} ^{\ad}
\to 
\X ^{\ad}$,
and
$J _{\underline{i}}:= J _{i _0}\times \dots \times J _{i _h}$.
For any $\underline{j} = (j _{i _0},\dots, j _{i _h})\in J _{\underline{i}}$, we set
$Y _{\underline{i}\,\underline{j}}:= Y _{i _0\, j _{i _0}} \cap \dots \cap Y _{i _h\, j _{i _h}}$,
$f _{\underline{i}\,\underline{j}}:= f _{i _0\, j _{i _0}} | _{\X _{\underline{i}}}  \cdots f _{i _h\, j _{i _h}} | _{\X _{\underline{i}}}$.
Denoting by 
$T _{\underline{i}\,\underline{j}} 
:= 
V (f _{\underline{i}\,\underline{j}})$ 
the divisor of $X _{\underline{i}}$, 
we have 
$Y _{\underline{i}\,\underline{j}} = X _{\underline{i}} \setminus T _{\underline{i}\,\underline{j}}$.

We get the covering 
$\mathscr{Y} _{\underline{i}}:= (Y _{\underline{i}\,\underline{j}}) _{\underline{j} \in J _{\underline{i}}}$ of 
$Y _{\underline{i}}$.
For any $\underline{\underline{j}} = (\underline{j} _0,\dots, \underline{j} _l)  \in (J _{\underline{i}} ) ^{1+l}$, 
we set 
$Y _{\underline{i},\underline{\underline{j}}}:= Y _{\underline{i}\,\underline{j} _0}\cap \dots \cap Y _{\underline{i}\,\underline{j} _{l}}$,
$f _{\underline{i},\underline{\underline{j}}}: = f _{\underline{i}\,\underline{j} _0}\cdots f _{\underline{i}\,\underline{j} _{l}}$,
and 
$v _{\underline{i},\underline{\underline{j}}} 
\colon 
] Y _{\underline{i},\underline{\underline{j}}} [ _{\X _{\underline{i}}} 
\to 
\X _{\underline{i}} ^{\ad}$.
With these notations, we get the functors
$j _{Y _{\underline{i},\underline{\underline{j}}}}  ^{\dag}
:=
v _{\underline{i},\underline{\underline{j}}*} 
v _{\underline{i},\underline{\underline{j}}} ^{-1} $ 
from the category of sheaves on 
$\X _{\underline{i}} ^{\ad}$.
Denoting by 
$T _{\underline{i}\,\underline{\underline{j}}} 
:= 
V (f _{\underline{i}\,\underline{\underline{j}}})$ 
the divisor of $X _{\underline{i}}$, 
we have 
$Y _{\underline{i}\,\underline{\underline{j}}} = X _{\underline{i}} \setminus T _{\underline{i}\,\underline{\underline{j}}}$.
We get
$\mathscr{T} _{\underline{i}}
:= (T _{\underline{i}\,\underline{\underline{j}}}  ) _{\underline{\underline{j}} \in (J _{\underline{i}} ) ^{1+l}}$ 
be a finite set of divisor of $X _{\underline{i}}$.

Let $E$ be an abelian sheaf on $\X ^{\ad}$.
As in \cite[4.1]{Be0}, we get the Cech bicomplexes 
$\check{C} ^{\dag \bullet \bullet} 
(\mathscr{X}, (\mathscr{Y} _{i} ) _{i\in I},E)$ 
associated with the coverings 
$\mathscr{X}, \mathscr{Y} _{\underline{i}}$ of $E$ by setting 
$$\check{C} ^{\dag hl} 
(\mathscr{X}, (\mathscr{Y} _{i} )_{i\in I},E)
:=
\prod _{\underline{i} \in I ^{1+h}} u _{\underline{i} *}  ^{\ad}
\check{C} ^{\dag l} 
( \fX _{\underline{i}}, \mathscr{Y} _{\underline{i}}, u _{\underline{i}} ^{\ad *} (E) )
=
\prod _{\underline{i} \in I ^{1+h}} u _{\underline{i} *}  ^{\ad}
\left (
\prod _{\underline{\underline{j}} \in J _{\underline{i}} ^{1+l}}
j _{Y _{\underline{i},\underline{\underline{j}}}}  ^{\dag} 
u _{\underline{i}} ^{\ad *} (E) 
\right),
$$
where $\check{C} ^{\dag l} ( \fX _{\underline{i}}, \mathscr{Y} _{\underline{i}} ,u _{\underline{i}} ^{\ad *} (E) )$ 
is 
defined in \ref{2.50Lazda-Pal-Book-ntn}.
We denote by 
$\check{C} ^{\dag \bullet} 
(\mathscr{X}, (\mathscr{Y} _{i} ) _{i\in I},E)$ 
the total complex of 
$\check{C} ^{\dag \bullet \bullet} 
(\mathscr{X}, (\mathscr{Y} _{i} ) _{i\in I},E)$.

Let us give the two extreme examples. 
On one hand, when $J _i$ has only one element for any $i \in I $, then 
$\check{C} ^{\dag \bullet \bullet} 
(\mathscr{X}, (\mathscr{Y} _{i} ) _{i\in I},E)$
is equal to the usual topological Check complex
$\check{C} ^{\dag \bullet} 
(\mathscr{X},E)$
given by 
$\check{C} ^{\dag h} 
(\mathscr{X},E)
:=
\prod _{\underline{i} \in I ^{1+h}} u _{\underline{i} *}  ^{\ad}
u _{\underline{i}} ^{\ad *} (E)$. 
On the other hand, when $I=\{i\}$ has only one element, 
the complex
$\check{C} ^{\dag \bullet} 
(\mathscr{X}, (\mathscr{Y} _{i} ) _{i\in I},E)$ 
is equal to the complex 
$\check{C} ^{\dag \bullet} 
(\fX, \mathscr{Y} _{i} ,E)$.

\begin{enumerate}[(a)]
\item We check similarly to \cite[4.1.3]{Be0} that 
$\check{C} ^{\dag \bullet} (\mathscr{X}, \mathscr{Y} _{\underline{i}},E)$ is a resolution of $j ^\dag _Y E$ 
(use \ref{2.50Lazda-Pal-Book}).

\item Since 
$Y _{\underline{i},\underline{\underline{j}}} = D (f_{\underline{i},\underline{\underline{j}}})$,
following \ref{Ru*ad=u*ad},
$u _{\underline{i} *}  ^{\ad}
\left (j _{Y _{\underline{i},\underline{\underline{j}}}}  ^{\dag} 
\O _{\X _{\underline{i}} ^{\ad}} 
\right)$ is acyclic for the functor $\sp _*$.
Hence, 
$$\check{C} ^{\dag hl} 
(\mathscr{X}, (\mathscr{Y} _{i} )_{i\in I},  \O _{\X ^{\ad}})
=
\prod _{\underline{i} \in I ^{1+h}} 
u _{\underline{i} *}  ^{\ad}
\left (
\prod _{\underline{\underline{j}} \in J _{\underline{i}} ^{1+l}}
\left (j _{Y _{\underline{i},\underline{\underline{j}}}}  ^{\dag} 
\O _{\X _{\underline{i}} ^{\ad}} 
\right)
\right) $$ 
is acyclic for the functor $\sp _*$.

\item Hence, we get in $D ^\mathrm{b} (\O _{\X,\Q})$ the isomorphism 
\begin{equation}
\label{RspCheck-reso}
\R \sp _* (j ^\dag _Y \O _{\X ^{\ad}}) 
\riso 
\sp _* \check{C} ^{\dag \bullet} 
(\mathscr{X}, (\mathscr{Y} _{i} )_{i\in I},\O _{\X ^{\ad}}).
\end{equation}

\end{enumerate}

Since 
$\sp _* u _{\underline{i} *}  ^{\ad}
\riso 
 u _{\underline{i} *}   \sp _*$,
then we get 
$$\sp _* u _{\underline{i} *}  ^{\ad}
\check{C} ^{\dag l} 
( \fX _{\underline{i}}, \mathscr{Y} _{\underline{i}}, \O _{\X _{\underline{i}} ^{\ad}}  )
\riso
u _{\underline{i} *}  \sp _* 
\check{C} ^{\dag l} 
( \fX _{\underline{i}}, \mathscr{Y} _{\underline{i}}, \O _{\X _{\underline{i}} ^{\ad}}  )
\underset{\ref{sp*Check=Check}}{\riso} 
u _{\underline{i} *} 
\check{C} ^{\dag l}  (\fX _{\underline{i}},\mathscr{T} _{\underline{i}}, \O _{\fX _{\underline{i}},\bbQ})
=
u _{\underline{i} *} 
\prod _{\underline{\underline{j}} \in J _{\underline{i}} ^{1+l}}
\O _{\X _{\underline{i}}  }( \hdag T _{\underline{i},\underline{\underline{j}}}) _\bbQ 
.$$
Let 
$\check{C} ^{\dag \bullet \bullet} 
(\mathscr{X}, (\mathscr{T} _{i} ) _{i\in I},\O _{\X})$
be the bicomplex defined similarly
and 
$\check{C} ^{\dag \bullet } 
(\mathscr{X}, (\mathscr{T} _{i} ) _{i\in I},\O _{\X})$
be its total complex. 
We can check the isomorphism
$$\sp _* \check{C} ^{\dag \bullet \bullet} 
(\mathscr{X}, (\mathscr{Y} _{i} )_{i\in I},\O _{\X ^{\ad}})
\riso
\check{C} ^{\dag \bullet \bullet} 
(\mathscr{X}, (\mathscr{T} _{i} ) _{i\in I},\O _{\X}).$$ 
Hence, 
$$\R \sp _* (j ^\dag _Y \O _{\X ^{\ad}}) 
\riso 
\check{C} ^{\dag \bullet } 
(\mathscr{X}, (\mathscr{T} _{i} ) _{i\in I},\O _{\X}).$$

\end{empt}

\subsection{Local cohomology with support in a closed subscheme having locally finite $p$-bases 
of the constant coefficient}
\label{loc-coh-const-subsec}
Let $\fP$ be a formal $\fS$-scheme locally of formal finite type
and having locally finite $p$-bases  over $\fS$.

\begin{empt}
\label{com-invMIC-ntn}
Let $u _0\colon X  \hookrightarrow P $ be a purely of codimension $r$ closed immersion of schemes having locally finite $p$-bases over $\Spec k$.
Choose $(\fP  _{\alpha}) _{\alpha \in \Lambda}$ an open affine covering of  $\fP $
and let us use the corresponding notation of \ref{ntnPPalpha} 

Similarly to the construction of
$u _0 ^! \colon 
\mathrm{Coh} (X, \fP/K)
\to 
\mathrm{Coh} ((\X   _\alpha )_{\alpha \in \Lambda}/K)$
of   \ref{const-u0!}, 
we can define the functor 
$u _0  ^*
  \colon
\mathrm{MIC} ^{\dag \dag} (\fP/K) 
\to 
\mathrm{MIC} ^{\dag \dag} ( (\X   _\alpha )_{\alpha \in \Lambda}/K)$
    as follows. 
    Let $\E \in  \mathrm{MIC} ^{\dag \dag} (\fP/K) $, i.e. 
    a coherent $\D ^{\dag} _{\fP ,\Q}$ which is also 
    $\O _{\fP  ,\Q}$-coherent.     
We set 
$\E _\alpha : =
u _{\alpha } ^{*} ( \E | \fP _\alpha) := 
\mathcal{H} ^{-r}u _{\alpha } ^{!} ( \E | \fP _\alpha)
\riso
u _{\alpha } ^{!} ( \E | \fP _\alpha) [-r]$.
Then 
$\E _\alpha$ is a coherent
$\D ^{\dag} _{\X  _{\alpha}\Q} $-module, 
which is also 
$\O  _{\X _{\alpha} ,\Q}$-coherent.
Via the isomorphisms of the form $\tau $ 
(\ref{prop-glueiniso-coh1}),
we obtain the glueing 
$\D ^{\dag} _{\X  _{\alpha \beta},\Q} $-linear
isomorphism
$ \theta _{  \alpha \beta} \ : \  p _2  ^{\alpha \beta !} (\E _{\beta}) \riso p  _1 ^{\alpha \beta !} (\E _{\alpha}),$
satisfying the cocycle condition:
$\theta _{13} ^{\alpha \beta \gamma }=
\theta _{12} ^{\alpha \beta \gamma }
\circ
\theta _{23} ^{\alpha \beta \gamma }$.

\end{empt}

\begin{prop}
[Berthelot]
\label{NCDgencoh}
Suppose there exists 
a finite $p$-basis 
$t _1, \dots, t _d $ of $\fP$ over $\fS$.
Let $T$ be the divisor of $P$ defined by setting
$T := V ( \overline{t} _1\cdots \overline{t} _r)$ with $r \leq d$,
where
$\overline{t} _1, \dots, \overline{t} _r$ are the images of $t _1,\dots, t _r$ in $\Gamma (P ,\O _{P})$.
We have  the exact sequence
\begin{equation}
(\D ^\dag _{\fP /\fS ,\Q}) ^{d}
\overset{\psi}{\longrightarrow} 
\D ^\dag _{\fP /\fS ,\Q}
\overset{\phi}{\longrightarrow}
\O _{\fP} (\hdag T ) _\Q
\to 0,
\end{equation}
where 
$\phi (P)= P \cdot (1/t _1\cdots t _r)$, and 
$\psi$ is defined by
\begin{equation}
\psi ( P _1,\dots, P _d) = \sum _{i=1} ^{r} P _i \partial _i t _i + \sum _{i=r+1} ^d P _i \partial _i.
\end{equation}
\end{prop}

\begin{proof}
This is checked similarly to \cite[4.3.2]{Be0}.
\end{proof}

\begin{dfn}
\label{st-nice-div}
Let $T$ be a divisor of $P$.
We say that $T$ is a ``strictly nice'' divisor of $P/S$ if 
for any $x \in T$, there exists an open subset $U$ of $P$ such that 
$U$ admits 
a finite $p$-basis 
$\overline{t} _1, \dots, \overline{t} _d$ of $P$ over $S$
satisfying 
$T\cap U = V ( \overline{t} _1\cdots \overline{t} _r)$ with $r \leq d$.

Remark that we can choose 
a finite $p$-basis 
$t _1, \dots, t _d $ of $\fP$ over $\fS$
such that 
$\overline{t} _1, \dots, \overline{t} _d$ are the images of $t _1,\dots, t _d$ in $\Gamma (P ,\O _{P})$
(see \ref{rem-dfn-pbasispadicbis}).
Hence, following \ref{NCDgencoh}, 
if $T$ is a strictly nice divisor of $P/S$, 
then 
$\O _{\fP} (\hdag T ) _\Q$ is 
$\D ^\dag _{\fP /\fS ,\Q}$-coherent.
\end{dfn}

\begin{ntn}
\label{ntn-GammaZO-rig}
Before defining local cohomology in the context of quasi-coherent complexes
(see \ref{dfn-Gamma-adm},
we will need to 
focus on the case of a 
$k$-scheme locally of formal finite type
and having locally finite $p$-bases  over $k$ 
for the constant coefficient as follows. 
We will see via \ref{cor-induction-div-coh2} that both local cohomology are canonically compatible, 
which justifies using the same notation. 

Let $u \colon X  \hookrightarrow P $ be a purely of codimension $r$ closed immersion of schemes having locally finite $p$-bases over $\Spec k$.
Let $j _X \colon  P   \setminus X\to P$ be the open immersion. 
We set 
$(\hdag X) ( \O _{\fP,\Q}) 
:=\R \sp _* j _X ^\dag (\O _{\fP _K})$
and
$\R \underline{\Gamma} ^\dag _X \O _{\fP,\Q} 
:=\R \sp _* \underline{\Gamma} ^\dag _X ( \O _{\fP _K})$.
By definition, 
$\R \underline{\Gamma} ^\dag _X \O _{\fP,\Q} $ is the 
local cohomology with support in $X$ of $\O _{\fP,\Q}$.
The exact sequence 
$0 \to \underline{\Gamma} ^\dag _X (\O _{\fP _K}) 
\to 
\O _{\fP _K}
\to 
j _X ^\dag (  \O _{\fP _K})
\to 
0$
induces the exact triangle 
\begin{equation}
\label{extriangleloc}
\R \underline{\Gamma} ^\dag _X \O _{\fP,\Q} 
\to 
\O _{\fP,\Q}
\to 
(\hdag X) ( \O _{\fP,\Q}) 
\to 
\R \underline{\Gamma} ^\dag _X \O _{\fP,\Q}  [1].
\end{equation}
For any integer $i\in \Z$, we set 
$\mathcal{H} ^{\dag i } _X (\O _{\fP,\Q}  ):= 
\mathcal{H} ^{i} \R \underline{\Gamma} ^\dag _X \O _{\fP,\Q} $.
\end{ntn}

\begin{rem}
\label{rem-why-2dfn-loc-coh}
Beware that in our work first we do need to use the left arrow of the exact triangle \ref{extriangleloc}
before being able to define local cohomology in the wider context of quasi-coherent complexes
(more precisely, see the proof of \ref{coh-cst-div} which is the main ingredient to define
the local cohomology in a wider context).

\end{rem}

\begin{prop}
[Berthelot]
\label{coh-smoothsubsch}
We keep notation \ref{ntn-GammaZO-rig}.
\begin{enumerate}[(a)]
\item $(\hdag X) ( \O _{\fP,\Q}) ,
\R \underline{\Gamma} ^\dag _X \O _{\fP,\Q}  
\in 
D ^{\mathrm{b}} _{\mathrm{coh}}
(\D ^\dag _{\fP, \Q})$, and
$\mathcal{H} ^{\dag i } _X (\O _{\fP,\Q}  )=0$
for any $i\not = r$.

\item Let $x \in P$. 
Following \ref{cor-closed-immer-local}, 
choose  an open affine formal subscheme $\U $ of $\fP $ containing $x$
such that 
there exist coordinates
$t _1, \dots, t _d \in \Gamma (\U ,\O _{\U})$ 
such that $X\cap U = V ( \overline{t} _1, \dots, \overline{t} _r)$ where
$r \leq d$ and
$\overline{t} _1, \dots, \overline{t} _r$ are the image of $t _1,\dots, t _r$ in $\Gamma (U ,\O _{U})$.
We have the exact sequence
\begin{equation}
\label{exseq-HrZ}
(\D ^\dag _{\U /\fS ,\Q}) ^{d}
\overset{\psi}{\longrightarrow} 
\D ^\dag _{\U /\fS ,\Q}
\overset{\phi}{\longrightarrow}
\mathcal{H} ^{\dag r} _{X \cap U} (\O _{\U,\Q}  )
\to 0,
\end{equation}
where 
$\phi (P)= P \cdot (1/t _1\cdots t _r)$, and 
$\psi$ is defined by
\begin{equation}
\psi ( P _1,\dots, P _d) = \sum _{i=1} ^{r} P _i t _i + \sum _{i=r+1} ^d P _i \partial _i.
\end{equation}
\end{enumerate}
\end{prop}

\begin{proof}
Similar to \cite[4.3.4]{Be0}.
\end{proof}

\begin{empt}
\label{empt-resolRsP*}
With the notation \ref{coh-smoothsubsch}, 
suppose, $\U  = \fP $.
For $i = 1,\dots, r$, put 
$X _i:  = V (\overline{t} _i)$, and 
$X _{i _0,\dots, i _{k}}
:= 
X _{i _0}\cup \dots \cup X _{i _k}$
(i.e. $V ( \overline{t} _{i _0} \cdots \overline{t} _{i _k})= X _{i _0,\dots, i _{k}}$).
Then $(\hdag X) ( \O _{\fP,\Q})$
is represented by the complex 
\begin{equation}
\label{coh-smoothsubsch-es1}
\prod _{i=1} ^{d}
\O _{\fP} (\hdag X _i) _{\Q}
\to 
\prod _{i _0 < i _1}
\O _{\fP} (\hdag X _{i _0 i _1}) _{\Q}
\to
\dots
\to 
\O _{\fP} (\hdag X _{1 \dots r}) _{\Q}
\to 0,
\end{equation}
whose first term is at degree $0$.
This yields that
$\R \underline{\Gamma} ^\dag _X \O _{\fP,\Q}$ is represented by the complex 
\begin{equation}
\label{coh-smoothsubsch-es2}
\O _{\fP,\Q}
\to 
\prod _{i=1} ^{d}
\O _{\fP} (\hdag X _i) _{\Q}
\to 
\prod _{i _0 < i _1}
\O _{\fP} (\hdag X _{i _0 i _1}) _{\Q}
\to
\dots
\to 
\O _{\fP} (\hdag X _{1 \dots r}) _{\Q}
\to 0,
\end{equation}
whose first term is at degree $0$.
Using \ref{NCDgencoh}, 
this is how Berthelot checked in \cite[4.3.4]{Be0} that 
$\R \underline{\Gamma} ^\dag _X \O _{\fP,\Q}  \in 
D ^{\mathrm{b}} _{\mathrm{coh}}
(\D ^\dag _{\fP, \Q})$.
\end{empt}

\begin{coro}
\label{coro-trace-upre}
Let $u \colon \fX  \hookrightarrow \fP $ be a purely of codimension $r$ closed immersion of formal schemes having locally finite $p$-bases over $\fS$.
\begin{enumerate}[(a)]
\item We have $u ^! (\hdag X) ( \O _{\fP,\Q})=0$, i.e. 
by applying the functor $u ^!$ to the canonical morphism
$\R \underline{\Gamma} ^\dag _X \O _{\fP,\Q}  
\to 
\O _{\fP,\Q}$, we get an isomorphism.

\item We have the canonical isomorphism
$u ^! (\O _{\fP,\Q})
\riso 
(\O _{\X,\Q}) [\delta _{X/P}]$.
We have the canonical isomorphism 
$\R \underline{\Gamma} ^\dag _X \O _{\fP,\Q}  
\riso 
u _+ u ^! (\O _{\fP,\Q})$ 
making commutative the canonical diagram
\begin{equation}
\label{coro-trace-upre-BK}
\xymatrix @ R=0,3cm {
{\R \underline{\Gamma} ^\dag _X \O _{\fP,\Q}  } 
\ar[r] ^-{\sim}
\ar[d] ^-{}
& 
{u _+ u ^! (\O _{\fP,\Q})} 
\ar[d] ^-{\mathrm{adj}} _-{\ref{Radj-u+flatf}}
\\ 
{\O _{\fP,\Q}} 
\ar@{=}[r] ^-{}
& 
{\O _{\fP,\Q}} 
}
\end{equation}

\end{enumerate}

\end{coro}

\begin{proof}
We can copy word by word 
the proof of \cite[8.1.10]{caro-6operations}.
\end{proof}

\begin{coro}
\label{coro-sp+jdagO}
Let $u _0 \colon X  \hookrightarrow P $ be a purely of codimension $r$ closed immersion of schemes having locally finite $p$-bases over $\Spec k$.
We have the isomorphism of $\mathrm{Coh} (X, \fP/K) $
of the form
$$
u _{0+} u _{0} ^*  (  \O  _{\fP,\Q})
\riso 
\mathcal{H} ^{\dag ,r} _X \O _{\fP,\Q} .$$

\end{coro}

\begin{proof}
We can copy word by word 
the proof of \cite[8.1.11]{caro-6operations}, we can check 
the  isomorphism 
\cite[8.1.11.4]{caro-6operations} is still valid, i.e. 
we have the canonical isomorphism 
$$ u ^! _0 (\mathcal{H} ^{\dag ,r} _X \O _{\fP,\Q} )
\riso 
u _0 ^*( \O _{\fP,\Q})$$
of 
$\mathrm{MIC} ^{\dag \dag} ( (\X   _\alpha )_{\alpha \in \Lambda}/K)$.
Then  we get
$$u _{0+} u _{0} ^*  (  \O  _{\fP,\Q})
\riso 
 u _{0+}  u ^! _0 (\mathcal{H} ^{\dag ,r} _X \O _{\fP,\Q} )
\underset{\ref{prop1}}{\riso} 
\mathcal{H} ^{\dag ,r} _X \O _{\fP,\Q}.$$ 
\end{proof}

\begin{prop}
\label{gluedautoduality}
Let $u _0 \colon X  \hookrightarrow P $ be a purely of codimension $r$ closed immersion of schemes having locally finite $p$-bases over $\Spec k$.
We have the isomorphism of $\mathrm{Coh} (X, \fP/K) $
of the form
$$
\DD u _{0} ^*  (  \O  _{\fP,\Q})
\riso 
u _{0} ^*  ( \O  _{\fP,\Q}) .$$
\end{prop}

\begin{proof}
The isomorphisms \ref{dualB=B} glue.
\end{proof}

\begin{prop}
\label{propspetdualsansfrob}
We have the functorial canonical isomorphism 
$\DD ( \R \underline{\Gamma} ^\dag _X \O _{\fP,\Q} [r] ) \riso \R \underline{\Gamma} ^\dag _X \O _{\fP,\Q} [r]$.
\end{prop}

\begin{proof}
This is a consequence of 
\ref{coro-sp+jdagO}
and of 
\ref{gluedautoduality}.
\end{proof}

\subsection{Convergent isocrystals on schemes having locally finite $p$-bases
and relative duality isomorphism}
\label{ntnsp+}
Let $\fP $ be 
a formal  $\fS$-scheme locally of formal finite type
and having locally finite $p$-bases  over $\fS$.
Let $u _0\colon X  \to P $ be a closed immersion of schemes having locally finite $p$-bases  over $S $.
Choose $(\fP  _{\alpha}) _{\alpha \in \Lambda}$ an open affine covering of  $\fP $.
We use the corresponding notation of \ref{ntnPPalpha}.

\begin{ntn}
\label{ntnMICalpha}
We denote by 
$\mathrm{MIC} ^{\dag \dag} ( (\X   _\alpha )_{\alpha \in \Lambda}/K)$
the full subcategory of 
$\mathrm{Coh} ((\X   _\alpha )_{\alpha \in \Lambda}/K)$ whose objects
$((\E _{\alpha})_{\alpha \in \Lambda},\, (\theta _{\alpha\beta}) _{\alpha ,\beta \in \Lambda})$
are such that, 
for all $\alpha \in \Lambda$,
$\E _{\alpha}$ 
is 
$\O _{\X _\alpha,\Q}$-coherent.
\end{ntn}

\begin{ntn}
\label{ntnMICdag2fs2}
We denote by $\mathrm{MIC} ^{\dag \dag} (X, \fP/K) $ the full subcategory of 
$\mathrm{Coh} (X, \fP/K)$ whose objects $\E$ satisfy the following condition:
for any affine open formal subscheme  
 $\fP ^{\prime }$ of $\fP $, for any morphism of formal schemes
  $v$ : $ \fX ^{\prime } \hookrightarrow \fP ^{\prime } $ which reduces modulo $\pi$ to the closed imbedding 
  $X ^{ } \cap P ^{\prime } \hookrightarrow P ^{\prime }$,
the sheaf $v ^! (\E |_{\fP^{\prime }}) $ is $\O _{\Y',\,\Q}$-coherent.
Finally, according to notation \ref{ntnMICdag2fs}, 
when $X= P$, we remove $X$ in the notation.

\end{ntn}

\begin{empt}
The functors $u ^!  _0$ and $ u _{0+}$ constructed in respectively \ref{const-u0!} and \ref{const-u0+} 
induce quasi-inverse equivalence of categories between 
$\mathrm{MIC} ^{\dag \dag} (X, \fP/K)$
and
$\mathrm{MIC} ^{\dag \dag} ( (\X   _\alpha )_{\alpha \in \Lambda}/K)$, i.e., 
we have the commutative diagram
\begin{equation}
\label{eqcat-u0+!}
\xymatrix{
{\mathrm{MIC} ^{\dag \dag} (X, \fP/K) } 
\ar@{^{(}->}[r] ^-{}
\ar@{.>}@<4ex>[d] ^-{\cong} _-{u _{0} ^!} 
& 
{\mathrm{Coh} (X, \fP /K)} 
\ar@{.>}@<4ex>[d] ^-{\cong} _-{u _{0} ^!} 
\\ 
{\mathrm{MIC} ^{\dag \dag} ( (\X   _\alpha )_{\alpha \in \Lambda}/K)}
\ar@{^{(}->}[r] ^-{}
\ar@{.>}@<4ex>[u] ^-{\cong} _-{u _{0+}} 
& {\mathrm{Coh} ((\X   _\alpha )_{\alpha \in \Lambda}/K).} 
\ar@{.>}@<4ex>[u] ^-{\cong} _-{u _{0+}} 
}
\end{equation}

\end{empt}

\begin{empt}
Let $f \colon \X ^{\prime } \to \X $ be an open immersion of 
formal  $\fS$-schemes locally of formal finite type
and having locally finite $p$-bases  over $\fS$.
Let $\E \in D ^{\mathrm{b}} _{\mathrm{coh}} ( \D ^\dag _{\X} (\hdag Z) _{\Q})$.
Similarly to \cite[3.2.8]{caro-construction}, we define the following isomorphism
\begin{gather}\notag
\xi \ :\   f ^!  \DD   (\E)
\riso
\R \mathcal{H} om _{\D ^\dag _{\X^{\prime },\Q}}
(f ^{!} ( \E),
f ^! _{\mathrm{r}}   (\D ^\dag _{\X,\Q}  \otimes _{\O _{\X}} \omega _{\X, \Q} ^{-1}))[\delta _X]
\\
\label{defDf!=f!D1bis}
\riso
\R \mathcal{H} om _{\D ^\dag _{\X^{\prime },\Q}}
(f ^{!}  ( \E),
(\D ^\dag _{\X^{\prime },\Q} \otimes  _{\O _{\X'}} \omega _{\X ^{\prime }/\fS } ^{-1} ) _\mathrm{t})[\delta _X]
\underset{\beta}{\riso}
\DD   f ^!  (\E),
\end{gather}
where $\beta$ is the transposition isomorphism exchanging both structures of left 
$\D ^\dag _{\X^{\prime },\Q} $-modules of 
$\D ^\dag _{\X^{\prime },\Q} \otimes  _{\O _{\X'}} \omega _{\X ^{\prime }/\fS } ^{-1} $.
\end{empt}

\begin{empt}
With notation \ref{ntnPPalpha},
let 
$((\E _{\alpha}) _{\alpha \in \Lambda}, (\theta _{\alpha \beta }) _{\alpha, \beta \in \Lambda})\in 
\mathrm{MIC} ^{\dag \dag}  ( (\X   _\alpha )_{\alpha \in \Lambda}/K)$.
Via the isomorphisms 
\ref{defDf!=f!D1bis}, the inverse of the isomorphism
$\DD (\theta _{\alpha \beta})$ is canonically isomorphic to 
$ \theta ^* _{  \alpha \beta} \ : \  p _2  ^{\alpha \beta !} (\DD (\E _{\beta})) \riso p  _1 ^{\alpha \beta !} (\DD (\E _{\alpha}))$.
These isomorphisms satisfy the cocycle condition 
(for more details, see  
 \cite[4.3.1]{caro-construction}).
Hence, we get the dual functor 
$$\DD \colon \mathrm{MIC} ^{\dag \dag}  ( (\X   _\alpha )_{\alpha \in \Lambda}/K)
\to \mathrm{MIC} ^{\dag \dag}  ( (\X   _\alpha )_{\alpha \in \Lambda}/K)$$
defined by 
$\DD ((\E _{\alpha}) _{\alpha \in \Lambda}, (\theta _{\alpha \beta }) _{\alpha, \beta \in \Lambda})
:=
((\DD (\E _{\alpha})) _{\alpha \in \Lambda}, (\theta ^*_{\alpha \beta }) _{\alpha, \beta \in \Lambda})$.

\end{empt}

\begin{empt}
\label{propspetdualsansfrob725}
With notation \ref{ntnPPalpha},
let 
$((\E _{\alpha}) _{\alpha \in \Lambda}, (\theta _{\alpha \beta }) _{\alpha, \beta \in \Lambda})\in 
\mathrm{MIC} ^{\dag \dag}  ( (\X   _\alpha )_{\alpha \in \Lambda}/K)$.
From the relative duality isomorphism (see 
\ref{rel-dual-isom-immf}),
we have the isomorphism
$u _{\alpha +} \circ \DD (\E _\alpha) 
\riso 
\DD \circ u _{\alpha +}  (\E _\alpha) $
These isomorphisms satisfy the cocycle condition 
(for more details, see  
 \cite[4.3.1]{caro-construction}), i.e.
we get the commutation isomorphism :
$$u _{0+} \circ \DD ((\E _{\alpha}) _{\alpha \in \Lambda}, (\theta _{\alpha \beta }) _{\alpha, \beta \in \Lambda}) 
\riso \DD \circ u _{0+} ((\E _{\alpha}) _{\alpha \in \Lambda}, (\theta _{\alpha \beta }) _{\alpha, \beta \in \Lambda})).$$
\end{empt}

\subsection{Convergent isocrystals, stability}
Convergent isocrystals in our context was defined in \ref{ntnsp+}.
We introduce here its inductive system avatar (see \ref{ntnMICdag2fs3}),
and we study its stability. This subsection can be avoided in a first reading. 

\begin{ntn}
\label{ntnMICdag2fs3}
Let  $\fP $ be a 
formal  $\fS$-scheme
of formal finite type
and having locally finite $p$-bases  over $\fS$.
Let $X$ be a closed subscheme of $P$
having locally finite $p$-bases  over $\Spec k$.

\begin{enumerate}[(a)]
\item We denote by $\mathrm{MIC} ^{(\bullet)} (X, \fP/K) $ the full subcategory of 
$\smash{\underrightarrow{LM}}  _{\Q, \mathrm{coh}}
(\smash{\widehat{\D}} _{\fP /\fS } ^{(\bullet)})$
consisting of objects 
$\E ^{(\bullet)}$ with support in $X$ 
and such that 
$\underrightarrow{\lim} (\E ^{(\bullet)})
\in \mathrm{MIC} ^{\dag \dag} (X, \fP/K) $
where
$\underrightarrow{\lim} 
\colon
\smash{\underrightarrow{LM}}  _{\Q, \mathrm{coh}}
(\smash{\widehat{\D}} _{\fP /\fS } ^{(\bullet)})
\cong
\mathrm{Coh} ( \smash{\D} ^\dag _{\fP,\Q} )$
is the equivalence of categories
of \ref{M-eq-coh-lim}, 
and where 
$\mathrm{MIC} ^{\dag \dag} (X, \fP/K) $
is defined in \ref{ntnMICdag2fs2}.
When $X=P$, we remove $X$ in the notation so that in this case 
we retrieve Notation \ref{cvisoc-f*}.\ref{cvisoc-f*2}.

\item \label{dfn-Gamma-pbases}
Let $\R \underline{\Gamma} ^\dag _{X} ( \cO _{\fP} ^{(\bullet)})
\in \smash{\underrightarrow{LD}} ^{\mathrm{b}} _{\Q,\mathrm{coh}} 
(\overset{^\mathrm{l}}{} \smash{\widehat{\D}} _{\fP /\fS  } ^{(\bullet)} )$ 
such that 
$\R \underline{\Gamma} ^\dag _{X} ( \cO _{\fP} ^{(\bullet)}) \riso 
\R \underline{\Gamma} ^\dag _X \O _{\fP,\Q} $, where this latter complex is defined 
at \ref{ntn-GammaZO-rig} (and is coherent thanks to \ref{coh-smoothsubsch}).
Then we can define the functor 
$\R\underline{\Gamma} ^\dag _{X}
\colon 
\smash{\underrightarrow{LD}} ^{\mathrm{b}} _{\Q,\mathrm{qc}} 
(\overset{^\mathrm{l}}{} \smash{\widehat{\D}} _{\fP /\fS  } ^{(\bullet)} )
\to 
\smash{\underrightarrow{LD}} ^{\mathrm{b}} _{\Q,\mathrm{qc}} 
(\overset{^\mathrm{l}}{} \smash{\widehat{\D}} _{\fP /\fS  } ^{(\bullet)} )
$ by setting for any 
$\E  ^{(\bullet)} 
\in 
\smash{\underrightarrow{LD}} ^{\mathrm{b}} _{\Q,\mathrm{qc}} 
(\overset{^\mathrm{l}}{} \smash{\widehat{\D}} _{\fP /\fS  } ^{(\bullet)} )$
\begin{equation}
\notag
\R \underline{\Gamma} ^\dag _{X} ( \E  ^{(\bullet)})
:=
\R \underline{\Gamma} ^\dag _{X} ( \cO _{\fP} ^{(\bullet)})
 \smash{\widehat{\otimes}}
^\L _{\cO ^{(\bullet)}  _{\fP}}
\E  ^{(\bullet)}.
\end{equation}
We will extend this  functor later in 
\ref{3.2.1caro-2006-surcoh-surcv} for any any subscheme $Y$ of $P$ in the case of overcoherent complexes,
but only in this subsection we consider this particular case. 
\end{enumerate}

\end{ntn}

\begin{lem}
\label{pre-loc-tri-B-t1Tbis}
Let $u \colon \fX  \to \fP $ be a closed immersion of 
 formal  $\fS$-schemes of formal finite type
and having locally finite $p$-bases  over $\fS$.
We suppose $(X \subset \fX)$ weak admissible.
For any 
$\E ^{(\bullet)} \in \smash{\underrightarrow{LD}} ^{\mathrm{b}} _{\Q,\mathrm{qc}}(\overset{^\mathrm{l}}{} \smash{\widehat{\D}} _{\fP /\fS  } ^{(\bullet)} )$,
we have the isomorphism
\begin{equation}
\label{pre-loc-tri-B-t1Tbis-iso}
\R \underline{\Gamma} ^\dag _{X} (\E ^{(\bullet)}) 
\riso 
u _{+} ^{ (\bullet)} \circ  u ^{ !(\bullet)} (\E ^{(\bullet)}),
\end{equation}
where by abuse of notation we denote $u (X)$ by $X$.
\end{lem}

\begin{proof}
Using \ref{surcoh2.1.4-cor1},
we reduce to the case where
$\E ^{(\bullet) }= \O _{\fP}^{(\bullet)} $.
Then the Lemma follows from 
\ref{coro-trace-upre}.
\end{proof}

\begin{prop}
\label{st-isoc-boxtimes}
Let $f\colon \fP \to \bbD ^r _{\fS}$ and
$g\colon \fQ \to \bbD ^s _{\fS}$ be two objects of 
$\scr{C} _{\fS}$ (see notation \ref{dfn-CfS}).
We suppose that $\fP $ and $\fQ$ 
have locally finite $p$-bases  over $\fS$.
Let $X$ (resp. $Y$) be a 
closed subscheme of $P$ (resp. $Q$)
and having locally finite $p$-bases  over $\Spec k$.
They induce the objects 
$X \to \bbD ^r _{S _0}$
and 
$Y \to \bbD ^s _{S _0}$
of $\fC _{S _0}$.
Let 
$\E ^{(\bullet)}$ be an object of 
$\mathrm{MIC} ^{(\bullet)} (X, \fP/K) $,
and  
$\cF ^{(\bullet)}$ be an object of 
$\mathrm{MIC} ^{(\bullet)} (Y, \fQ/K) $.
Then 
$\E ^{(\bullet)} \smash{\widehat{\boxtimes}} ^\L _{\O _{\fS }} \cF ^{(\bullet)}
\in 
\mathrm{MIC} ^{(\bullet)} (X\times _{\fC _{S _0}} Y, \fP\times _{\scr{C} _{\fS}} \fQ/K) $ (see notation later \ref{dfnboxtimes}).
\end{prop}

\begin{proof}
Following Lemma \ref{exact-boxtimes} (notice that this Lemma does not use this subsection), we already know
$\E ^{(\bullet)} \smash{\widehat{\boxtimes}} ^\L _{\O _{\fS }} \cF ^{(\bullet)}
\in 
\smash{\underrightarrow{LM}}  _{\Q, \mathrm{coh}}
( \smash{\widehat{\D}} _{\fP \times _{\scr{C} _{\fS}} \fQ/\fS }  ^{(\bullet)} )$.
Since the proposition is local, 
using 
\ref{prop-boxtimes-v+}, 
we reduce to the case where $X=P$ and $Y=Q$. 
Then this is obvious.
\end{proof}

\begin{prop}
\label{cor-com-sp+-f*}
Let  $f \colon \fP ^{\prime } \to \fP $
be a morphism having locally finite $p$-bases of formal $\fS$-schemes
of formal finite type and having locally finite $p$-bases. 
Let $X$ (resp. $X'$) be a closed subscheme of $P$ (resp. $P'$)
having locally finite $p$-bases  over $\Spec k$.
We suppose  $f (X ') \subset X$.
Let 
$\E ^{(\bullet)}$ and $\cF ^{(\bullet)}$ be two objects of 
$\mathrm{MIC} ^{(\bullet)} (X, \fP/K) $.
\begin{enumerate}[(a)]
\item 
\label{cor-com-sp+-f*-1}
$\R \underline{\Gamma} ^\dag _{X'} f ^{!(\bullet)} \E ^{(\bullet)} [-\delta _{X'/X}]
\in 
\mathrm{MIC} ^{(\bullet)} (X', \fP'/K) $.
\item $\DD ^{(\bullet)} (\E ^{(\bullet)})
\in
\mathrm{MIC} ^{(\bullet)} (X, \fP/K) $.
\item 
We have 
$\E ^{(\bullet)}
\smash{\widehat{\otimes}}^\L  _{\O ^{(\bullet)}  _{\fP} }  \FF ^{(\bullet)} [-\delta _{X/P}]
\in \mathrm{MIC} ^{(\bullet)} (X, \fP/K) $.
\end{enumerate}
\end{prop}

\begin{proof}
The fact that
$\R \underline{\Gamma} ^\dag _{X'} f ^{!(\bullet)} \E ^{(\bullet)} [- \delta _{X'/X}]
\in 
\mathrm{MIC} ^{(\bullet)} (X', \fP'/K) $
is local in $\fP '$.
Hence, we can suppose there exists a closed immersion of 
 formal  $\fS$-schemes of formal finite type
and having locally finite $p$-bases  over $\fS$ of the form
$\mathfrak{u} \colon \X \hookrightarrow \fP$ 
(resp. $\mathfrak{u} ' \colon \X ' \hookrightarrow \fP' $, 
resp. $\mathfrak{a} \colon \X' \to \X$)
which reduces modulo $\pi$ to $u _0$ (resp. $u ' _0$, resp. $a$).
Following \ref{pre-loc-tri-B-t1Tbis-iso}, 
$\R \underline{\Gamma} ^{\dag} _{X '}  f  ^{!(\bullet)} (   \E ^{(\bullet)})
\riso 
\mathfrak{u} _{+} ^{\prime (\bullet)} \circ  \mathfrak{u} ^{\prime !(\bullet)} \circ f  ^{!(\bullet)}(\E ^{(\bullet)})
\riso 
\mathfrak{u} _{+} ^{\prime (\bullet)} \circ  \mathfrak{a}   ^{!(\bullet)} \circ \mathfrak{u} ^{!(\bullet)} (\E ^{(\bullet)})$.
Since
$\mathfrak{u} ^{!(\bullet)} (\E ^{(\bullet)}) 
\in 
\mathrm{MIC} ^{(\bullet)} (\X/K) $,
then
$\L \mathfrak{a}   ^{*(\bullet)} \circ \mathfrak{u} ^{!(\bullet)} (\E ^{(\bullet)})
\in 
\mathrm{MIC} ^{(\bullet)} (\X'/K) $ (see \ref{corostab-MIC-f*}).
Since 
$\L \mathfrak{a}   ^{*(\bullet)} 
=
\mathfrak{a}   ^{!(\bullet)}[- \delta _{X'/X}] $, we get the first statement. 

The second statement is a consequence of 
\ref{dualisoscvdag}.
The last one is a consequence of \ref{st-isoc-boxtimes} and of the first statement.
\end{proof}

\begin{prop}
\label{cor-com-sp+-f*2}
With notation \ref{cor-com-sp+-f*}, 
we have the isomorphism of 
$\mathrm{MIC} ^{(\bullet)} (X', \fP'/K) $ 
of the form
\begin{equation}
\label{cor-com-sp+-f*2iso}
\DD ^{(\bullet)}
\left (
\R \underline{\Gamma} ^\dag _{X'} f ^{!(\bullet)} \E ^{(\bullet)} [-\delta _{X'/X}]
\right )
\riso 
\R \underline{\Gamma} ^\dag _{X'} f ^{!(\bullet)} (\DD ^{(\bullet)} \E ^{(\bullet)} ) [-\delta _{X'/X}].
\end{equation}

\end{prop}

\begin{proof}
Following \ref{cor-com-sp+-f*}, 
the objects appearing in \ref{cor-com-sp+-f*2iso}
belong to $\mathrm{MIC} ^{(\bullet)} (X', \fP'/K) $.
Hence, it is sufficient to check the isomorphism
\ref{cor-com-sp+-f*2iso}
in $\mathrm{MIC} ^{\dag \dag} (X', \fP'/K)$ 
(i.e.after applying the functor
$\underrightarrow{\lim}$ which is an equivalence of categories). 
We denote by $a \colon X' \to X$ the induced morphism. 
We get the commutative diagram
\begin{equation}
\label{com-sp+-f*-squ}
\xymatrix  
{
{\fP ^{\prime }}
\ar@{=}[r] ^-{}
&
{\fP ^{\prime }}
\ar[r] ^-{f}
\ar@{}[rd] ^-{} |\square
&
{\fP}
\\
{X ^{\prime }}
\ar@{^{(}->}[u] ^-{u'}
\ar[r] ^-{b}
&
{X ^{\prime \prime}}
\ar@{^{(}->}[u] ^{u'' }
\ar[r] ^-{c}
&
{X,}
\ar@{^{(}->}[u] ^-{u}
}
\end{equation}
where $X '': =f ^{-1} (X)$, vertical arrows are the canonical closed immersions. 
Hence to prove the isomorphism \ref{cor-com-sp+-f*2iso} we reduce to the following two cases.

1) We suppose $X '=f ^{-1} (X)$. In that case, notice that the functor
$\R \underline{\Gamma} ^\dag _{X'} $ is useless in the isomorphism
\ref{cor-com-sp+-f*2iso}.
Let $(\fP  _{\alpha}) _{\alpha \in \Lambda}$ be an open covering of  $\fP $
satisfying the condition of \ref{ntnPPalpha}.
We fix some liftings as in \ref{ntnPPalpha} and we use the same notation.
Moreover, we denote by
$\fP ' _\alpha := f ^{-1} (\fP _\alpha)$,
$\fX ' _\alpha := \fP ' _\alpha   \times _{\fP  _\alpha  }\fX _\alpha$,
$a _\alpha \colon \fX ' _\alpha \to \fX  _\alpha$ the projection, and similarly for other notations.
Let $((\E _{\alpha})_{\alpha \in \Lambda},\, (\theta _{\alpha\beta}) _{\alpha ,\beta \in \Lambda})$
be an object of 
$\mathrm{MIC} ^{\dag \dag} ( (\X   _\alpha )_{\alpha \in \Lambda}/K)$
(see notation \ref{ntnMICalpha})
We get canonically an object  of 
$\mathrm{MIC} ^{\dag \dag} ( (\X  ' _\alpha )_{\alpha \in \Lambda}/K)$
 of the form
$ (a _{\alpha} ^* (\E _{\alpha})_{\alpha \in \Lambda},\, (\theta ' _{\alpha\beta}) _{\alpha ,\beta \in \Lambda})$.
This yields the functor 
$a ^* 
\colon 
 \mathrm{MIC} ^{\dag \dag} ( (\X   _\alpha )_{\alpha \in \Lambda}/K)
 \to 
  \mathrm{MIC} ^{\dag \dag} ( (\X  ' _\alpha )_{\alpha \in \Lambda}/K)$.
    Consider the following diagram :
\begin{equation}
\label{com-sp+-f*diag1}
\xymatrix{
{ \mathrm{MIC} ^{\dag \dag} ( (\X   _\alpha )_{\alpha \in \Lambda}/K)} 
\ar[r] ^-{a ^*}
\ar[d] ^-{\bbD} _-{\ref{dualisoscvdag}}
&
{ \mathrm{MIC} ^{\dag \dag} ( (\X  ' _\alpha )_{\alpha \in \Lambda}/K)} 
\ar[d] ^-{\bbD} _-{\ref{dualisoscvdag}}
\\
{ \mathrm{MIC} ^{\dag \dag} ( (\X   _\alpha )_{\alpha \in \Lambda}/K)} 
\ar[r] ^-{a ^*}
\ar@<4ex>[d] ^-{\ref{eqcat-u0+!}} _-{u _{0+}}
& 
{ \mathrm{MIC} ^{\dag \dag} ( (\X  ' _\alpha )_{\alpha \in \Lambda}/K)} 
\ar@<4ex>[d] ^-{\ref{eqcat-u0+!}} _-{u ' _{0+}}
\\
{\mathrm{MIC} ^{\dag \dag} (X, \fP/K) }
\ar@<4ex>[u] ^-{\ref{eqcat-u0+!}} _-{u _{0} ^!}
\ar[r] ^-{f ^*}
&
{\mathrm{MIC} ^{\dag \dag} (X', \fP'/K) .}
\ar@<4ex>[u] ^-{\ref{eqcat-u0+!}} _-{u _{0} ^{\prime !}}
}
\end{equation}
We have the canonical isomorphism
$( a _{\alpha} ^* (\E _{\alpha} ) )  ^\vee
\riso 
a _{\alpha} ^* (\E _{\alpha } ^\vee)$.
Via the canonical isomorphisms  of
\ref{evee=De0}, this yields the isomorphisms
$\bbD ( a _{\alpha} ^* (\E _{\alpha} ) )  
\riso 
a _{\alpha} ^* (\bbD (\E _{\alpha }))$
which commute with the glueing data. Hence, the top square is commutative. 
By transitivity of the inverse image with respect to the composition, 
the bottom square involving $u _0 ^!$ and $u _0 ^{\prime !}$ is commutative up to canonical isomorphism. 
Since $u _{0+} $ and $u _0 ^!$ (resp. $u _{0+} '$ and $u _0 ^{\prime !}$)
are canonically quasi-inverse equivalences of categories,
this yields the bottom square involving $u _{0+} $ and $u _{0+} '$ is commutative up to canonical isomorphism. 

Using the commutativity of $\bbD$ with $u _{0+}$ (see \ref{propspetdualsansfrob725}),
using the commutativity of the diagram \ref{com-sp+-f*diag1}, 
this yields the isomorphism
\ref{cor-com-sp+-f*2iso}.

2) Now suppose $f = id$ and $a$ is a closed immersion.   
Then, 
we can fix some liftings (separately) for both $u$ and $u'$ (for the later case, add some primes in notation)
and we use notation \ref{ntnPPalpha} as follows.
By using \ref{cor-closed-immer-local}, 
we can choose such covering $(\fP  _{\alpha}) _{\alpha \in \Lambda}$
so that for every $\alpha\in \Lambda$, $X  _\alpha$ is affine and has finite $p$-basis.
Since $P$ is separated, for any $\alpha,\beta ,\gamma \in \Lambda$,
$X _{\alpha \beta }$ and  $X _{\alpha \beta \gamma }$ are also affine and have finite $p$-bases.
Hence,  following \ref{lifting-pbasis}.\ref{lifting-pbasis-p2} such liftings exists.
Moreover,  choose some lifting morphisms
$a _\alpha \colon \fX ' _\alpha \to \fX  _\alpha$, 
and similarly for other notations.
Let $(\fP  _{\alpha}) _{\alpha \in \Lambda}$ be an open covering of  $\fP $
satisfying the condition of \ref{ntnPPalpha} for both $X$ an $X'$.
Consider the following diagram.
\begin{equation}
\label{com-sp+-f*diag2}
\xymatrix{
{ \mathrm{MIC} ^{\dag \dag} ( (\X   _\alpha )_{\alpha \in \Lambda}/K)} 
\ar[r] ^-{a ^*}
\ar[d] ^-{\bbD} _-{\ref{dualisoscvdag}}
&
{ \mathrm{MIC} ^{\dag \dag} ( (\X  ' _\alpha )_{\alpha \in \Lambda}/K)} 
\ar[d] ^-{\bbD} _-{\ref{dualisoscvdag}}
\\
{ \mathrm{MIC} ^{\dag \dag} ( (\X   _\alpha )_{\alpha \in \Lambda}/K)} 
\ar[r] ^-{a ^*}
\ar@<4ex>[d] ^-{\ref{eqcat-u0+!}} _-{u _{0+}}
& 
{ \mathrm{MIC} ^{\dag \dag} ( (\X  ' _\alpha )_{\alpha \in \Lambda}/K)} 
\ar@<4ex>[d] ^-{\ref{eqcat-u0+!}} _-{u ' _{0+}}
\\
{\mathrm{MIC} ^{\dag \dag} (X, \fP/K) }
\ar@<4ex>[u] ^-{\ref{eqcat-u0+!}} _-{u _{0} ^!}
\ar[r] ^-{\R \underline{\Gamma} ^\dag _{X'} [-\delta _{X'/X}]}
&
{\mathrm{MIC} ^{\dag \dag} (X', \fP'/K) .}
\ar@<4ex>[u] ^-{\ref{eqcat-u0+!}} _-{u _{0} ^{\prime !}}
}
\end{equation}
The commutativity up to a canonical isomorphism
of the top square of \ref{com-sp+-f*diag2} is checked as for 
\ref{com-sp+-f*diag1}. It remains to look at the bottom square. 
 Let 
 $\E \in \mathrm{MIC} ^{\dag \dag} (X, \fP/K) $. 
 The canonical morphism
 $$ u ^{\prime !} _{\alpha}
 \left (
 \R \underline{\Gamma} ^\dag _{X' }( \E) | \fP _{\alpha} 
 \right)  
 [-\delta _{X'/X}]
 \to
u ^{\prime !} _{\alpha}
 \left (
  \E| \fP _{\alpha} 
 \right)
 [-\delta _{X'/X}] 
 $$
 is an isomorphism.
Moreover, 
$u ^{\prime !} _{\alpha}
 \left (
  \E| \fP _{\alpha}
 \right)
  [-\delta _{X'/X}] 
\riso
a ^{ !} _{\alpha} u ^{ !} _{\alpha}
 \left (
 \E| \fP _{\alpha}  
 \right)
 [-\delta _{X'/X}]
 \riso 
a ^{ *} _{\alpha}  \left ( u ^{ !} _{\alpha}
 ( \E| \fP _{\alpha})
 \right)$.
 By composition, this yields the isomorphism
 \begin{equation}
  u ^{\prime !} _{\alpha}
 \left (
 \R \underline{\Gamma} ^\dag _{X' }( \E) | \fP _{\alpha} 
 \right)  
 [-\delta _{X'/X}]
\riso 
a ^{ *} _{\alpha}  \left ( u ^{ !} _{\alpha}
 ( \E| \fP _{\alpha})
 \right).
 \end{equation}
These isomorphisms glue, hence we get the commutativity up to canonical isomorphism of the bottom square.
\end{proof}

\section{Exterior tensor products}

\subsection{On the exactness of the exterior tensor product}

\begin{lemm}
\label{lemm-flat-k2R[[t]]}
Let $R= \cV$ or $R = \cV / \pi ^{i+1} \cV$. 
Set $R  [[ \underline{v}]]:=R [[ v_1, \dots, v _{s}]]$. 
Let $M$ be an $R  [[ \underline{v}]]$-module. 
The following conditions are equivalent.
\begin{enumerate}[(a)]
\item The $R  [[ \underline{v}]]$-module $M$ is flat. 

\item The $k  [[ \underline{v}]]$-module $M \otimes _{R [[ \underline{v}]]} k  [[ \underline{v}]]$ is flat
and 
$\mathrm{Tor} _1 ^{R  [[ \underline{v}]]} ( k  [[ \underline{v}]] , M) = 0$.

\item The $k  [[ \underline{v}]]$-module $M \otimes _{R} k$ is flat
and 
$\mathrm{Tor} _1 ^{R} ( k , M) = 0$.
\end{enumerate}

\end{lemm}

\begin{proof}
Using the Krull intersection theorem, the noetherianity of $R$ and $R  [[ \underline{v}]]$, the separated completeness for the $p$-adic topology
of $R  [[ \underline{v}]]$, 
we can check that for any ideal $\mathfrak{a}$ of $R$, 
$\mathfrak{a} \otimes _R R  [[ \underline{v}]]$ is an $R  [[ \underline{v}]]$-module of finite type and is Hausdorff for the $p$-adic topology. 
Hence, thanks to \cite[Theorem 1 of III.5.2]{bourbaki}, 
we get the equivalence between $(a)$ and $(b)$.
Moreover, since 
$R \to k$ is finite then the canonical morphism 
$R [[ \underline{v}]] \otimes _{R} k
\to k  [[ \underline{v}]]$
is an isomorphism.
This yields 
$M \otimes ^{\bbL} _{R [[ \underline{v}]]} k  [[ \underline{v}]]
\riso 
M \otimes ^{\bbL} _{R} k$.
Hence, we get the equivalence between $(b)$ and $(c)$.
\end{proof}

\begin{lemm}
\label{N-flat-v-torsion}
Let $R= \cV$ or $R = \cV / \pi ^{i+1} \cV$. 
Let $N$ be a 
$R [[ \underline{v}]]$-module, where 
$R [[ \underline{v}]]:=R [[ v_1, \dots, v _{s}]]$. 
The following condition are equivalent. 
\begin{enumerate}[(a)]
\item The $R [[ \underline{v}]]$-module $N$ is flat.

\item The $k [[ \underline{v}]]$-module $N \otimes _{R} k$  has no $v _j$-torsion for any $j = 1,\dots, s$
and we have $\mathrm{Tor} _1 ^{R} ( k , N) = 0$.
\end{enumerate}

\end{lemm}

\begin{proof}
Let 
$\mathrm{ev} _0 \colon k [[ \underline{v}]] \to k$ be the homomorphism
of $k$-algebras defined by setting 
$P \to P (0)$. 
We have the exact sequence 
\begin{equation}
\label{N-flat-v-torsion-exseq}
0 
\to 
\left ( k [[ \underline{v}]] \right) ^{s}
\overset{(v _1, \dots, v _s)}{\longrightarrow}
k [[ \underline{v}]]
\overset{\mathrm{ev} _0}{\longrightarrow}
k
\to 0.
\end{equation}

1) Suppose $N$ is flat.
Set $N _0 : = N \otimes _{R} k$. 
Then
$\mathrm{Tor} _1 ^{R} ( k, N) = 0$
and 
$N _0$ is a flat $k [[ \underline{v}]]$-module (see \ref{lemm-flat-k2R[[t]]}). 
By using \ref{N-flat-v-torsion-exseq}, we can check 
$N _0 \otimes _{k [[ \underline{v}]]} ^\bbL k $ 
is isomorphic to the complex
$N _0 ^s \overset{(v _1, \dots, v _s)}{\longrightarrow} N _0$ such that 
$N _0$ is at the $0$th place.
Hence, 
$N$ has no $v _j$-torsion.

2) Conversely, suppose $N _0 = N \otimes _{R} k$ has no $v _j$-torsion. 
Let $O$ be a 
$k [[ \underline{v}]]$-module.
Since $k$ is a field, we have 
$N _0 \otimes _{k [[ \underline{v}]]} ^\bbL O
\riso 
\left (N _0 \otimes _{k [[ \underline{v}]]} ^\bbL 
k \right ) \otimes _{k} 
O$.
Since $N _0$ has no $v _j$-torsion,
then
$N _0 \otimes _{k [[ \underline{v}]]} ^\bbL k 
\riso 
N _0 \otimes _{k [[ \underline{v}]]} k $. 
This yields, 
$N _0 \otimes _{k [[ \underline{v}]]} ^\bbL O
\riso  
N _0 \otimes _{k [[ \underline{v}]]}  O$. 
Hence, $N _0$ is flat. 
We conclude by using \ref{lemm-flat-k2R[[t]]} that 
if moreover $\mathrm{Tor} _1 ^{R} ( k , N) = 0$
then $N$ is flat.
\end{proof}

\begin{rem}
\label{rem-loc-hom-fft}
Let $R$ be a local commutative ring.

\begin{enumerate}[(a)]
\item A morphism of $R$-algebras of the form 
$\phi \colon R [[ t _1,\dots, t _r]] \to R [[ u _1, \dots, u _s]]$
is necessarily a local homomorphism of complete local rings. 
Indeed, 
let $\psi \colon R [[ t _1,\dots, t _r]] \twoheadrightarrow R/\fm _R$
be the composition of $\phi $ with the morphism of $R$-algebras
$R [[ u _1, \dots, u _s]] \twoheadrightarrow R/\fm _R $ given by 
$u _i \mapsto 0$.
Then $\Ker \psi$ is the maximal ideal of 
$ R [[ t _1,\dots, t _r]]$, i.e. 
$\ker \psi
=
\fm _R +(t _1,\dots, t _r)$.
Hence, we are done.

\item More precisely, the data of a morphism of the form
$\phi \colon R [[ t _1,\dots, t _r]] \to R [[ u _1, \dots, u _s]]$
is equivalent to that of the data of $r$ elements of the maximal ideal of 
$R [[ u _1, \dots, u _s]]$ (indeed, 
$R [[ u _1, \dots, u _s]]$ is a complete local ring).
\end{enumerate}

\end{rem}

\begin{prop}
\label{prop-flat-u2uv}
Let $R= \cV$ or $R = \cV / \pi ^{i+1} \cV$. 
We set 
$R [[ \underline{u}]] := R [[ u_1, \dots, u _{r}]]$,
$R [[ \underline{v}]] := R [[ v_1, \dots, v _{s}]]$
and
$R [[ \underline{u},  \underline{v}]] 
:=
R [[ u_1, \dots, u _{r}, v _1, \dots, v _{s}]]$. 
Let $M$ be an $R [[ \underline{u}]]$-module such that $\mathrm{Tor} _1 ^{R} ( k , M) = 0$.
Then $M \otimes _{R [[ \underline{u}]]}R [[ \underline{u},  \underline{v}]] $ is a flat $R [[ \underline{v}]]$-module
(via the underlying structure given by $R [[ \underline{v}]] \to R [[\underline{u}, \underline{v}]]$).
\end{prop}

\begin{proof}
1) First we reduce to the case where $R=k$ as follows.
Since 
$k [[ \underline{u}]]
\to 
k [[ \underline{u},  \underline{v}]]$
and 
$R [[ \underline{u}]]
\to 
R [[ \underline{u},  \underline{v}]]$
are flat, then we get the isomorphisms
\begin{gather}
\notag
( M \otimes _{R [[ \underline{u}]]}R [[ \underline{u},  \underline{v}]] ) 
\otimes ^{\bbL} _{R }k 
\riso 
 M \otimes ^{\bbL}  _{R [[ \underline{u}]]}k  [[ \underline{u},  \underline{v}]] 
\\
 \riso 
( M \otimes ^{\bbL }_{R [[ \underline{u}]]} k [[ \underline{u}]])
\otimes _{k [[ \underline{u}]]}
k [[ \underline{u},  \underline{v}]] 
\riso 
( M \otimes ^{\bbL }_{R} k)
\otimes _{k [[ \underline{u}]]}
k [[ \underline{u},  \underline{v}]] .
\end{gather}
Since 
$\mathrm{Tor} _1 ^{R} ( k , M) = 0$, 
then 
$ M \otimes ^{\bbL }_{R} k
\riso 
 M \otimes _{R} k$. 
This yields
 $\mathrm{Tor} _1 ^{R} ( k , M \otimes _{R [[ \underline{u}]]}R [[ \underline{u},  \underline{v}]] ) = 0$.
Hence, 
by using \ref{lemm-flat-k2R[[t]]} we reduce to the case $R=k$.

2) Set $N: = M \otimes _{k [[ \underline{u}]]}k [[ \underline{u},  \underline{v}]] $.
Let  $y \in N$ such that 
$ v _j y = 0$ for some $j$. 
Following \ref{N-flat-v-torsion}, we have to check that $y = 0$.
We can write $y$ as a finite sum of the form 
$y = \sum _i x _i \otimes P _{i} $, where
$x _{i} \in M$ 
and 
$P _{i} \in k [[ \underline{u},  \underline{v}]] $.
Let $M' $ be the $k [[ \underline{u}]]$-submodule
of $M$ generated by the finite family $(x _{i}) _i$.
Set 
$N' : = 
M ' 
\otimes _{k [[ \underline{u}]]}
k [[ \underline{u},  \underline{v}]] $
and let 
$y ' \in N'$ be the element of $N'$ which can be written of the form 
$\sum _i x _i \otimes P _{i} $, i.e. 
the image of $y'$  via the injection
$N ' \hookrightarrow N$
(indeed 
$k [[ \underline{u}]]
\to 
k [[ \underline{u},  \underline{v}]]$ is flat)
is $y$. 
Since $N ' \hookrightarrow N$ is injective, then
$ v _j y' = 0$ in $N'$ for some $j$.
Hence, we reduce to the case where $M' = M$, i.e. to the case where $M$ is finitely generated.
Set $A : =k [[ \underline{u}]]$. 
Since $M$ is an $A$-module of finite type, then
$N= M \otimes _{A} A [[ \underline{v}]]$ is complete for the $(v _1, \dots, v _s)$-adic topology
and an element of $M \otimes _{A} A [[ \underline{v}]]$ can be written uniquely 
of the form
$\sum _{\underline{i} \in \bbN ^s} x _{\underline{i}} \underline{v} ^{\underline{i}}$ 
with $x _{\underline{i}} \in M$ (for instance, 
use \cite[3.2.3]{Be1}).
Hence, $N$ has no $v _j$-torsion.\end{proof}

\begin{empt}
\label{ntn-CT}
Fix some integer $j\geq 0$ and 
set
$T = S _j$. 

\begin{enumerate}[(a)]
\item We denote by $ \scr{F}  _T$ the full subcategory of the category
of $T$-schemes consisting in objects $X/T$ of formal finite type, i.e. such that there exists a finite type morphism of the form 
$X \to \bbD ^r _T$ for some integer $r$.

\item We denote by $\scr{C} _{T}$ the category whose objects are
finite type morphisms of the form 
$X \to \bbD ^r _T$ for some integer $r$.
A morphism $f \to g$ of $\scr{C} _{T}$, where 
$f\colon X \to \bbD ^r _T$ and
$g\colon Y \to \bbD ^s _T$ are objects of $\scr{C} _{T}$, consists in 
two morphisms
$\alpha \colon 
X \to Y$ and 
$\beta \colon \bbD ^r _T \to \bbD ^s _T$ of $ \scr{F}  _T$ making commutative the diagram
\begin{equation}
\xymatrix{
{X} 
\ar[r] ^-{\alpha}
\ar[d] ^-{f}
& 
{Y} 
\ar[d] ^-{g}
\\ 
{\bbD ^r _T} 
\ar[r] ^-{\beta}
& 
{\bbD ^s _T.} 
}
\end{equation}
We notice that $\beta$ is associated to 
a local homomorphism of local rings (see \ref{rem-loc-hom-fft}).
We denote by 
$(\alpha, \beta) $ such a morphism.

\item 
\label{ntn-CT-item3}
We get a functor 
$\mathscr{S}
\colon 
\scr{C} _{T} \to  \scr{F}  _T$ defined by setting 
$\mathscr{S} 
(
X \to \bbD ^r _T
)
= X$.

\end{enumerate}

\end{empt}

\begin{empt}
\label{priX,Y}
Fix some integer $j\geq 0$ and 
set
$T = S _j$. 
Let 
$f _1\colon X _1 \to \bbD ^{r _1} _T$ and 
$f _2 \colon X _2 \to \bbD ^{r _2} _T$ be two finite type morphisms.
We denote by 
$f \times _{\scr{C} _{T}} g$ the projection 
\begin{equation}
\label{ftimesgdfn}
f _1 \times _{\scr{C} _{T}} f _2
\colon 
\left ( X _2 \times _{\bbD ^{r _1} _T} \bbD ^{r _1+r _2} _T \right ) 
\times _{\bbD ^{r _1+r _2} _T}
\left ( \bbD ^{r _1+r _2} _T \times _{\bbD ^{r _2} _T} X _2 \right ) 
\to 
\bbD ^{r _1+r _2} _T,
\end{equation}
where $pr _{1} ^{r _1,r _2}\colon 
\bbD ^{r _1+r _2} _T \to \bbD ^{r _1} _T$ is the projection on the first $r _1$-coordinates
and 
$pr _{2} ^{r _1,r _2}\colon  \bbD ^{r _1 +r _2} _T \to \bbD ^{r _2} _T$ is the projection on the last $r _2$-coordinates.
In other words, they are the respective compositions
$pr _{1} ^{r _1,r _2}
\colon 
\bbD ^{r _1 +r _2} _T \to \bbD ^{r _1} _T \times _{T} \bbD ^{r _2} _T \to \bbD ^{r _1} _T$
and 
$pr _{2} ^{r _1 ,r _2}
\colon 
\bbD ^{r _1 +r _2} _T \to \bbD ^{r _1} _T \times _{T} \bbD ^{r _2} _T \to \bbD ^{r _2} _T$.
We put 
\begin{equation}
\label{ntnXCTY}
X  _1 \times _{\scr{C} _{T}} X _2 : = 
\left ( X _1 \times _{\bbD ^{r _1} _T} \bbD ^{r _1 +r _2} _T \right ) 
\times _{\bbD ^{r _1 +r _2} _T}
\left ( \bbD ^{r _1 +r _2} _T \times _{\bbD ^{r _2} _T} X _2 \right )
\riso 
 X _1 \times _{\bbD ^{r _1} _T} 
 \bbD ^{r _1 +r _2} _T \times _{\bbD ^{r _2} _T} X _2.
\end{equation}
We denote by 
$pr _1 
\colon 
X  _1 \times _{\scr{C} _{T}} X _2
\to X _1$
and 
by 
$pr _2 
\colon 
X  _1 \times _{\scr{C} _{T}} X _2
\to X _2$
the canonical projection. 
We get the morphisms
$(pr _1 , pr _{1} ^{r _1 ,r _2}) 
\colon 
f _1 \times _{\scr{C} _{T}} f _2 
\to f _1$
and 
$(pr _2 , pr _{2} ^{r _1 ,r _2}) 
\colon 
f _1 \times _{\scr{C} _{T}} f _2 
\to f _2$
 of $\scr{C} _{T}$.

We can check that $f _1\times f _2$ is the product of $f _1$ and $f _2$ in the category
$\scr{C} _{T}$ as follows, which justifies the notation.
Let 
$h \colon Z \to \bbD _T ^u$ 
be an object of $\scr{C} _{T}$, 
$(\alpha _1, \beta _1) \colon h \to f _1$ 
and 
$(\alpha _2, \beta _2) \colon h \to f _2$ 
be two morphisms of $\scr{C} _{T}$.
Using the remark \ref{rem-loc-hom-fft}, 
we can check that 
there exists a unique morphism
$\beta \colon \bbD _T ^u \to \bbD _T ^{r _1 +r _2}$
whose composition with 
the canonical map with 
$\bbD _T ^{r _1 +r _2} \to \bbD _T ^{r _1}$
(resp. $\bbD _T ^{r+s} \to \bbD _T ^{ r _2}$)
is $\beta _1$ (resp. $\beta _2$).
Via $\beta \circ h$, we can view $Z$ as a
$\bbD _T ^{r _1 +r _2}$-scheme. 
We get a unique morphism of $\bbD _T ^{r _1 +r _2}$-schemes
$\widetilde{\alpha} _1 \colon Z \to X _1 \times _{\bbD ^r _T} \bbD ^{r _1 +r _2} _T $
(resp. $\widetilde{\alpha} _2 \colon Z \to X _2 \times _{\bbD ^{r _1} _T} \bbD ^{r _1 +r _2} _T $)
whose composition with the projection 
$X  _1\times _{\bbD ^{r _1} _T} \bbD ^{r _1 +r _2} _T \to X _1$
(resp. 
$X _2 \times _{\bbD ^{r _1} _T} \bbD ^{r _1 +r _2} _T \to X _2$)
is 
$\alpha _1$ (resp. $\alpha _2$).
This yields the morphism of $\bbD _T ^{r _1 +r _2}$-schemes
$\alpha 
: =
( \widetilde{\alpha} _1, \widetilde{\alpha} _2)
\colon
Z
\to  
\left ( X _1 \times _{\bbD ^{r _1} _T} \bbD ^{r _1 +r _2} _T \right ) 
\times _{\bbD ^{r _1 +r _2} _T}
\left ( \bbD ^{r _1 +r _2} _T \times _{\bbD ^{r _2} _T} X _2 \right ) $
making commutative the following diagram
\begin{equation}
\notag
\xymatrix{
{Z} 
\ar@/^{0,4cm}/[rr] ^-{\alpha  _i}
\ar@{.>}[r] _-{\alpha}
\ar[d] ^-{h}
& 
{X  _1 \times _{\scr{C} _{T}} X _2}
\ar[r] _-{pr _i} 
\ar[d] ^-{f _1 \times f _2}
& 
{X _i}
\ar[d] ^-{f _i} 
\\ 
{\bbD ^{r _1 +r _2} _T}
\ar@{.>}[r] ^-{\beta} 
\ar@/_{0,4cm}/[rr] _-{\beta  _i}
& 
{ \bbD ^{r _1 +r _2} _T} 
\ar[r] ^-{pr _{i} ^{r _1 ,r _2}}
& 
{ \bbD ^{r _i} _T} 
}
\end{equation}
for any $i= 1,2$.
Moreover, the morphism 
$(\alpha, \beta)$ is the only one making commutative the above diagram

Hence, we are done.

\end{empt}

\begin{rem}
Let $f\colon X \to \bbD ^r _T$ and 
$g\colon X \to \bbD ^s _T$ be two finite type morphisms.
Then this is not clear that we can define a morphism
$(f,g) \colon X \to \bbD ^{r+s} _T$ whose composition with 
$\bbD ^{r+s} _T \to \bbD ^{r} _T$
(resp. $\bbD ^{r+s} _T \to \bbD ^{s} _T$)
is $f$ (resp. $g$).
In order to define products, 
this is why we have introduced the category $\scr{C} _{T}$. 
\end{rem}

\begin{empt}
\label{graphCT}
Fix some integer $j\geq 0$ and 
set
$T = S _j$. 
Let $(\alpha, \beta)  \colon f \to g$ be a morphism of $\scr{C} _{T}$, with
$f\colon X \to \bbD ^r _T$ and
$g\colon Y \to \bbD ^s _T$.
Using the universal property of the product in $\scr{C} _{T}$, 
there exists a unique morphism $(\gamma _{(\alpha, \beta)}, \gamma _\beta)$ making commutative the following diagrams:
\begin{equation}
\notag
\xymatrix{
{X} 
\ar@/^0,5cm/[rr] ^-{id}
\ar[d] ^-{f}
\ar@{.>}[r] _-{\gamma _{(\alpha, \beta)}}
&
{X \times _{\scr{C} _{T}} Y}
\ar[r] _-{pr _1}
\ar[d] ^-{}
& 
{X} 
\ar[d] ^-{f}
\\ 
{\bbD ^r _T} 
\ar@/_0,5cm/[rr] _-{id}
\ar@{.>}[r] ^-{\gamma _{\beta}}
& 
{\bbD ^{r+s} _T} 
\ar[r] ^-{pr _1 ^{r,s}}
&
{\bbD ^r  _T,} 
}
\xymatrix{
{X} 
\ar@/^0,5cm/[rr] ^-{\alpha}
\ar[d] ^-{f}
\ar@{.>}[r] _-{\gamma _{(\alpha, \beta)}}
&
{X \times _{\scr{C} _{T}} Y}
\ar[r] _-{pr _2}
\ar[d] ^-{}
& 
{Y} 
\ar[d] ^-{g}
\\ 
{\bbD ^r _T} 
\ar@/_0,5cm/[rr] _-{\beta}
\ar@{.>}[r] ^-{\gamma _{\beta}}
& 
{\bbD ^{r+s} _T} 
\ar[r] ^-{pr _2 ^{r,s}}
&
{\bbD ^s _T.} 
}
\end{equation}
Remark that the morphism $\gamma _\beta$ only depends on $\beta$ which justifies the notation. 
Since the composition of 
$\gamma _{(\alpha, \beta)}$ 
with 
$X \times _{\scr{C} _{T}} Y
\to 
X \times _T Y$ is an immersion, 
then so is 
$\gamma _{(\alpha, \beta)}$.
\end{empt}

\begin{lem}
\label{Xr,s2s-flat}
Fix some integer $j\geq 0$ and put $S = S _j$. 
Let 
$f\colon X \to \bbD ^r _S$ be a finite type morphism such that 
$X/S$ has locally finite $p$-bases.
Let us denote by 
$X \times _{\bbD ^r _S} \bbD ^{r+s} _S $
the base change of $X$ via the first projection
$pr _{1} ^{r ,s}\colon 
\bbD ^{r +s} _T \to \bbD ^{r} _T$.
Then the canonical morphism 
$X \times _{\bbD ^r _S} \bbD ^{r+s} _S 
\to \bbD ^{s} _S$,
which is the composition of the projection
$X \times _{\bbD ^r _S} \bbD ^{r+s} _S 
\to 
\bbD ^{r+s} _S $
with $pr _{2} ^{r,s}$, 
is flat. 
\end{lem}

\begin{proof}
Since $X/S$ is flat (see \ref{f0formétale-fforméta0})
and 
$X \times _{\bbD ^r _S} \bbD ^{r+s} _S \to X$ is flat, 
then so is 
$X \times _{\bbD ^r _S} \bbD ^{r+s} _S/S$.
Since $X$ and $\bbD ^{r+s} _S $ are noetherian,
since $\bbD ^{r+s} _S \times _{S} S _0 \riso \bbD ^{r+s} _{S _0}$
(because $S _0 \to S$ is finite), then by using the 
flatness criterium by fiber (see \cite[Theorem 11.3.10]{EGAIV3}),
we reduce to the case where $S= S _0$. 

We can suppose $X = \Spec A$. 
We set 
$k [[ \underline{u}]] := k [[ u_1, \dots, u _{r}]]$,
$k [[ \underline{v}]] := k [[ v_1, \dots, v _{s}]]$
and
$k [[ \underline{u},  \underline{v}]] 
:=
k [[ u_1, \dots, u _{r}, v _1, \dots, v _{s}]]$. 
We have to check that the homomorphism
$k [[ \underline{v}]] \to A \otimes _{k [[ \underline{u}]]}k [[ \underline{u},  \underline{v}]] $ is flat,
which follows from
Proposition \ref{prop-flat-u2uv}.
\end{proof}

\begin{prop}
\label{corXtimeCY2XtimesY}
We keep notation \ref{graphCT}.
\begin{enumerate}[(a)]
\item The canonical  morphism 
\begin{equation}
\label{XtimeCY2XtimesY}
X \times _{\scr{C} _{T}} Y
\to 
X \times _{T} Y
\end{equation}
is relatively perfect.

\item 
\label{corXtimeCY2XtimesY-item2}
Suppose 
$X/T$ and $Y/T$ have locally finite $p$-bases.
Then $X \times _{\scr{C} _{T}} Y/T$ have locally finite $p$-bases and 
the projections
\begin{gather}
\notag
pr _1 
\colon 
X \times _{\scr{C} _{T}} Y
\to X,
\\
\label{XtimeCY2XtimesY-corprN2}
pr _2 
\colon 
X \times _{\scr{C} _{T}} Y
\to Y
\end{gather}
are flat and have locally finite $p$-bases. 

\item 
\label{corXtimeCY2XtimesY-item3}
Suppose 
$X/T$ and $Y/T$ have locally finite $p$-bases.
Set 
$Z:= X \times _{\scr{C} _{T}} Y$. Let $\cE$ be a quasi-coherent $\cO _{X}$-module and
$\cF$ be a quasi-coherent $\cO _{Y}$-module.
If 
$Tor _1 ^{\cO _T} ( \cO _{S _0}, \cE) = 0$ 
and
$Tor _1 ^{\cO _T} ( \cO _{S _0}, \cF) = 0$,
then the canonical morphism
$$pr _1 ^* (\cE)
\otimes ^{\bbL} _{\cO _Z}
pr _2 ^* (\cF)
\riso 
pr _1 ^* (\cE)
\otimes  _{\cO _Z}
pr _2 ^* (\cF)$$
is an isomorphism.
\end{enumerate}
\end{prop}

\begin{proof}
1) Since $ \bbD ^{r} _T \to \bbA ^{r} _T$
and $ \bbD ^{s} _T \to \bbA ^{s} _T$ are relatively perfect, then so is
$\bbD ^{r} _T \times _{T} \bbD ^s _T \to \bbA ^{r} _T \times _{T} \bbA ^s _T$.
Since the composition of 
$\bbD ^{r+s} _T 
\to 
\bbD ^{r} _T \times _{T} \bbD ^s _T 
\to 
\bbA ^{r} _T \times _{T} \bbA ^s _T$
is also relatively perfect, then this yields that so is 
the first morphism
$\bbD ^{r+s} _T \to \bbD ^{r} _T \times _{T} \bbD ^s _T$.
Hence, the composition morphism below is relatively perfect:
$$
\left ( X \times _{\bbD ^r _T} \bbD ^{r+s} _T \right ) 
\times _{\bbD ^{r+s} _T}
\left ( \bbD ^{r+s} _T \times _{\bbD ^s _T} Y \right ) 
\riso 
X \times _{\bbD ^r _T} \bbD ^{r+s} _T 
 \times _{\bbD ^s _T} Y 
\to 
X \times _{\bbD ^r _T} ( \bbD ^{r} _T \times _{T} \bbD ^s _T)
 \times _{\bbD ^s _T} Y 
\riso
X \times _T Y .
$$

2) The projection
$pr _2 $
is canonically isomorphic to the projection
$\left ( X \times _{\bbD ^r _T} \bbD ^{r+s} _T \right ) 
\times _{\bbD ^{s} _T}
Y
\to Y$, which is the base change 
of $X \times _{\bbD ^r _T} \bbD ^{r+s} _T \to \bbD ^{s} _T$
via $Y \to \bbD ^{s} _T$.
Thanks to \ref{Xr,s2s-flat}, we obtain the flatness of the projection $pr _2$. 
By symmetry, we get the flatness of $pr _1$.

Since $pr _1 $ (resp. $pr _2 $)
is the composition of \ref{XtimeCY2XtimesY}
with the projection 
$X \times _T Y \to X$
(resp. $X \times _T Y \to Y$), 
we get from 1) that 
$pr _1 $ (resp. $pr _2 $) have locally finite $p$-bases.

3) Set $R = \cV / \pi ^{j+1} \cV$. 
Since this is local, we can suppose $X= \Spec A$ and $Y = \Spec B$.
We have $\bbD ^r _T= \Spec R [[ \underline{u}]]$, 
$\bbD ^s _T= \Spec R [[ \underline{v}]] $,
$\bbD ^{r+s} _T= \Spec R [[\underline{u}, \underline{v}]] $,
with
$R [[ \underline{u}]] := R [[ u_1, \dots, u _{r}]]$,
$R [[ \underline{v}]] := R [[ v_1, \dots, v _{s}]]$
and
$R [[ \underline{u},  \underline{v}]] 
:=
R [[ u_1, \dots, u _{r}, v _1, \dots, v _{s}]]$. 
We compute 
\begin{gather}
C:=\Gamma (Z, \cO _Z)
=
\left (A \otimes _{R [[ \underline{u}]]} 
R [[ \underline{u},  \underline{v}]] \right) 
\otimes _{R [[ \underline{u},  \underline{v}]]}
\left (B \otimes _{R [[ \underline{v}]]} 
R [[ \underline{u},  \underline{v}]] \right) 
\notag
\\
\label{CisoABuv}
\liso
\left (
A \otimes _{R [[ \underline{u}]]} 
R [[ \underline{u},  \underline{v}]] 
\right) 
\otimes _{R [[ \underline{v}]]}
B.
\end{gather}

Set $E := \Gamma (X, \cE)$ 
and 
$F := \Gamma (X, \cF)$.
Since $\cE$ and $\cF$ are quasi-coherent, then 
we have
$\Gamma (Z, pr _1 ^* (\cE))
\riso
E \otimes _{A} C
$
and
$\Gamma (Z, pr _2 ^* (\cF))
\riso
F \otimes _{B} C$.
Following the part 2), 
since $C/B$ is flat then we get the following last isomorphism
\begin{gather*}
\Gamma (Z, pr _1 ^* (\cE)) 
\otimes _{C} ^{\bbL}
\Gamma (Z, pr _2 ^* (\cF))
\riso 
(E \otimes _{A} C) 
\otimes _{C} ^{\bbL}
(F \otimes _{B} C) 
\riso
(E \otimes _{A} C) 
\otimes _{B} ^{\bbL}
F.
\end{gather*}
We have the isomorphisms 
$$E \otimes _{A} C
\underset{\ref{CisoABuv}}{\riso} 
\left (E \otimes _{R [[ \underline{u}]]} 
R [[ \underline{u},  \underline{v}]]  \right) 
\otimes _{R [[ \underline{v}]]} 
B
\underset{\ref{prop-flat-u2uv}}{\riso} 
\left (E \otimes _{R [[ \underline{u}]]} 
R [[ \underline{u},  \underline{v}]]  \right) 
\otimes _{R [[ \underline{v}]]} ^{\bbL}
B.
$$ 
Hence: 
\begin{gather*}
(E \otimes _{A} C) 
\otimes _{B} ^{\bbL}
F
\riso
\left ( 
\left (E \otimes _{R [[ \underline{u}]]} 
R [[ \underline{u},  \underline{v}]]  \right) 
\otimes _{R [[ \underline{v}]]} ^{\bbL}
B
\right )
\otimes _{B} ^{\bbL}
F 
\riso 
( E \otimes _{R [[ \underline{u}]]} 
R [[ \underline{u},  \underline{v}]] )
\otimes _{R [[ \underline{v}]]} ^{\bbL}
F 
\end{gather*}
By composition, this yields
$$\Gamma (Z, pr _1 ^* (\cE)) 
\otimes _{C} ^{\bbL}
\Gamma (Z, pr _2 ^* (\cF)).
\riso 
( E \otimes _{R [[ \underline{u}]]} 
R [[ \underline{u},  \underline{v}]] )
\otimes _{R [[ \underline{v}]]} ^{\bbL}
F.$$
Similarly we get 
$$\Gamma (Z, pr _1 ^* (\cE)) 
\otimes _{C} 
\Gamma (Z, pr _2 ^* (\cF))
\riso 
( E \otimes _{R [[ \underline{u}]]} 
R [[ \underline{u},  \underline{v}]] )
\otimes _{R [[ \underline{v}]]} 
F.$$
We conclude via the canonical isomorphism
$$
( E \otimes _{R [[ \underline{u}]]} 
R [[ \underline{u},  \underline{v}]] )
\otimes _{R [[ \underline{v}]]} ^{\bbL}
F
\underset{\ref{prop-flat-u2uv}}{\riso} 
( E \otimes _{R [[ \underline{u}]]} 
R [[ \underline{u},  \underline{v}]] )
\otimes _{R [[ \underline{v}]]} 
F.$$
\end{proof}

\begin{ex}
\label{cor-flat-u2uv}
We set 
$k [[ \underline{u}]] := k [[ u_1, \dots, u _{r}]]$,
$k [[ \underline{v}]] := k [[ v_1, \dots, v _{s}]]$
and
$k [[ \underline{u},  \underline{v}]] 
:=
k [[ u_1, \dots, u _{r}, v _1, \dots, v _{s}]]$. 
Let $M$ be a $k [[ \underline{u}]]$-module.
and 
$N$ be a $k [[ \underline{v}]]$-module.
Let 
$pr _1 \colon 
k [[ \underline{u}]]
\to 
k [[ \underline{u},  \underline{v}]]$
and 
$pr _2 \colon 
k [[ \underline{u}]]
\to 
k [[ \underline{u},  \underline{v}]]$
be the canonical monomorphisms. 
We get 
$pr _1 ^* M
= 
M \otimes _{k [[ \underline{u}]]}k [[ \underline{u},  \underline{v}]] $
and
$pr _2 ^* N = N \otimes _{k [[ \underline{v}]]}k [[ \underline{u},  \underline{v}]] $. 
Then we have the canonical isomorphism
$$pr _1 ^* M 
\otimes ^{\bbL} _{k [[ \underline{u},  \underline{v}]] }
pr _2 ^* N
\riso 
pr _1 ^* M 
\otimes _{k [[ \underline{u},  \underline{v}]] }
pr _2 ^* N.$$
\end{ex}

\begin{empt}
[Local $p$-basis of the product in $\scr{C} _{T}$]
\label{XtimesY-pbases}
We keep notation \ref{graphCT}. 
Suppose we have two relatively perfect $T$-morphisms
$\alpha \colon X \to \bbA ^{n} _T$
and
$\beta \colon Y \to \bbA ^{m} _T$.
Then we get the commutative diagram of $T$-morphisms
\begin{equation}
\label{XtimesY-pbases-diag1}
\xymatrix{
{} 
& 
{X }
\ar[r] ^-{\alpha}
&
{\bbA ^{n} _T}
\\
{X \times _{\scr{C} _{T}} Y} 
\ar[r] ^-{\ref{XtimeCY2XtimesY}}
\ar[ru] ^-{pr _1}
\ar[rd] ^-{pr _2}
& 
{X \times _{T} Y}
\ar[r] ^-{\alpha \times \beta}
\ar[u] ^-{pr _1}
\ar[d] ^-{pr _2}
&
{\bbA ^{n+m} _T}
\ar[u] ^-{pr _1}
\ar[d] ^-{pr _2}
 \\ 
& 
\ar[r] ^-{ \beta}
{Y}
&
{\bbA ^{m} _T}
}
\end{equation}
whose horizontal arrows are relatively perfect. 
\end{empt}

\begin{empt}
\label{XtimesY-boxtimesD}
We keep notation \ref{graphCT} and we
suppose 
$X/T$ and $Y/T$ have locally finite $p$-bases.
We set $Z: =X \times _{\scr{C} _{T}} Y$.
By computation using finite $p$-bases (see \ref{XtimesY-pbases}), 
we check the canonical morphism of $\cO _Z$-algebras (for both structure)
$\cP ^{n} _{Z/T, (m)}
\to 
\cP ^{n} _{Z/X, (m)}$
is surjective. 
By duality, this yields the canonical morphism
of left $\cD ^{(m)} _{Z /X}$-modules
\begin{equation}
\label{XtimesY-boxtimesD-mor1}
\cD ^{(m)} _{Z /X}\to \cD ^{(m)} _{Z /T}
\end{equation}
is injective.
The sheaf $\cO _Z$ has a canonical structure of  left 
$\cD ^{(m)} _{Z /X}$-module 
and canonical structure of  left 
$\cD ^{(m)} _{Z /T}$-module.
Both structures are compatible with the morphism
$\cD ^{(m)} _{Z /X}\to \cD ^{(m)} _{Z /T}$.
Hence, via a local computation using finite $p$-bases, 
we can check that $\cD ^{(m)} _{Z /X}$ is an $\cO _Z$-submodule 
of $\cD ^{(m)} _{Z /T}$ for both (the left or the right) structures.
Via a local computation using finite $p$-bases, 
we can also check that 
$\cD ^{(m)} _{Z /X}\to \cD ^{(m)} _{Z /T}$ is in fact a morphism of rings.

Similarly, we check that the canonical morphism
of left $\cD ^{(m)} _{Z /Y}$-modules
\begin{equation}
\label{XtimesY-boxtimesD-mor2}
\cD ^{(m)} _{Z /Y}\to \cD ^{(m)} _{Z /T}
\end{equation}
is injective, is a morphism of rings and that 
$\cD ^{(m)} _{Z /Y}$ is an $\cO _Z$-submodule 
of $\cD ^{(m)} _{Z /T}$ for both (the left or the right) structures.

The morphisms \ref{XtimesY-boxtimesD-mor1} and \ref{XtimesY-boxtimesD-mor2}
induce the homomorphism
\begin{equation}
\label{XtimesY-boxtimesD-iso1}
\cD ^{(m)} _{Z /X} 
\otimes _{\cO _{Z}}
\cD ^{(m)} _{Z /Y}
\to 
\cD ^{(m)} _{Z /T},
\end{equation}
where to define the tensor product
we use 
the left structure of $\cO_Z$-module of $\cD ^{(m)} _{Z /Y}$
and 
the right structure of $\cO_Z$-module of $\cD ^{(m)} _{Z /X}$.
By local computations with finite $p$-bases (see \ref{XtimesY-pbases-diag1}),
we compute that the morphism \ref{XtimesY-boxtimesD-iso1} is an isomorphism.

We have by functoriality the canonical morphisms of
left $\cD ^{(m)} _{Z /T}$-modules:
\begin{equation}
\label{XtimeCY-D2pr*D}
\cD ^{(m)} _{Z /T}
\to 
pr _1 ^* \cD ^{(m)} _{X/T}
\
\text{~and~} 
\cD ^{(m)} _{Z /T}
\to 
pr _2 ^* \cD ^{(m)} _{Y/T}.
\end{equation}
By local computations with finite $p$-bases (see \ref{XtimesY-pbases-diag1}),
we check that the composition morphisms
\begin{gather}
\notag
\cD ^{(m)} _{Z /X}
\underset{\ref{XtimesY-boxtimesD-mor1}}{\longrightarrow}
\cD ^{(m)} _{Z /T}
\underset{\ref{XtimeCY-D2pr*D}}{\longrightarrow}
pr _2 ^* \cD ^{(m)} _{Y/T},
\\
\label{XtimesY-boxtimesD-iso2}
\cD ^{(m)} _{Z /Y}
\underset{\ref{XtimesY-boxtimesD-mor2}}{\longrightarrow}
\cD ^{(m)} _{Z /T}
\underset{\ref{XtimeCY-D2pr*D}}{\longrightarrow}
pr _1 ^* \cD ^{(m)} _{X/T}
\end{gather}
are isomorphisms.

We have the natural morphism
$pr _1 ^{-1} \cD ^{(m)} _{X/T} \to 
pr _1 ^{*} \cD ^{(m)} _{X/T} 
\underset{\ref{XtimesY-boxtimesD-iso2}}{\liso} 
\cD ^{(m)} _{Z /Y}$.
By a local computation with finite $p$-bases, we can check that this is a morphism
of $\cO _T$-algebras.
By composition, this yields the homomorphism of $\cO _T$-algebras 
\begin{equation}
\label{XtimesY-boxtimesD-mor3} 
pr _1 ^{-1} \cO _{X} 
\to
pr _1 ^{-1} \cD ^{(m)} _{X/T} 
\to 
\cD ^{(m)} _{Z /Y}
\underset{\ref{XtimesY-boxtimesD-mor1}}{\longrightarrow} 
\cD ^{(m)} _{Z /T}.
\end{equation}
Similarly, we get the homomorphism of $\cO _T$-algebras 
\begin{equation}
\label{XtimesY-boxtimesD-mor3bis} 
pr _2 ^{-1} \cO _{Y} 
\to
pr _2 ^{-1} \cD ^{(m)} _{Y/T} \to 
\cD ^{(m)} _{Z /X}
\underset{\ref{XtimesY-boxtimesD-mor2}}{\longrightarrow} 
\cD ^{(m)} _{Z /T}.
\end{equation}

Consider the following $\cO _T$-algebras :
$\cO  _{X} \boxtimes ^{\mathrm{top}} _{T} \cO _{Y}
:=
pr _1 ^{-1} \cO _{X}
\otimes _{\O _T}
pr _2 ^{-1} \cO _{Y}$
and
$\cD ^{(m)} _{X/T} \boxtimes ^{\mathrm{top}} _{T} \cD ^{(m)} _{Y/T}
:=
pr _1 ^{-1} \cD ^{(m)} _{X/T}
\otimes _{\O _T}
pr _2 ^{-1} \cD ^{(m)} _{Y/T}$.
From \ref{XtimesY-boxtimesD-mor3} and \ref{XtimesY-boxtimesD-mor3bis},
we get the homomorphism of $\cO _T$-algebras
\begin{equation}
\label{XtimesY-boxtimesD-mor4}
\cD ^{(m)} _{X/T} \boxtimes ^{\mathrm{top}} _{T} \cD ^{(m)} _{Y/T}
\to 
\cD ^{(m)} _{Z/T} .
\end{equation}
This yields by extension the homomorphism 
of $(\cO _X , \cD ^{(m)} _{X/T} \boxtimes ^{\mathrm{top}} _{T} \cD ^{(m)} _{Y/T})$-bimodules:
\begin{equation}
\label{XtimesY-boxtimesD-mor5}
\cO _Z
\otimes _{\cO  _{X} \boxtimes ^{\mathrm{top}} _{T} \cO _{Y}}
( \cD ^{(m)} _{X/T} \boxtimes ^{\mathrm{top}} _{T} \cD ^{(m)} _{Y/T}) 
\to 
\cD ^{(m)} _{Z/T} .
\end{equation}

We have the isomorphism
\begin{equation}
\label{boxtimestop-dfn2-pre}
\left ( 
(\cO  _{X} \boxtimes ^{\mathrm{top}} _{T} \cO _{Y} )
\otimes _{pr _1 ^{-1} \O _{X}}
pr _1 ^{-1} \cD ^{(m)} _{X/T}
 \right)
\otimes _{\cO  _{X} \boxtimes ^{\mathrm{top}} _{T} \cO _{Y} }
\left ( 
(\cO  _{X} \boxtimes ^{\mathrm{top}} _{T} \cO _{Y} )
\otimes _{pr _2 ^{-1} \O _{Y}}
pr _2 ^{-1} \cD ^{(m)} _{Y/T}
 \right)
 \riso 
 \cD ^{(m)} _{X/T} \boxtimes ^{\mathrm{top}} _{T} \cD ^{(m)} _{Y/T}.
\end{equation}
By applying the functor $\cO _Z
\otimes _{\cO  _{X} \boxtimes ^{\mathrm{top}} _{T} \cO _{Y}} -$ 
to \ref{boxtimestop-dfn2-pre}, 
this yields 
\begin{equation}
\label{boxtimestop-dfn2-prebis}
\left ( 
\cO _Z
\otimes _{pr _1 ^{-1} \O _{X}}
pr _1 ^{-1} \cD ^{(m)} _{X/T}
 \right)
\otimes _{\cO _Z}
\left ( 
\cO _Z
\otimes _{pr _2 ^{-1} \O _{Y}}
pr _2 ^{-1} \cD ^{(m)} _{Y/T}
 \right)
 \riso 
\cO _Z
\otimes _{\cO  _{X} \boxtimes ^{\mathrm{top}} _{T} \cO _{Y}}
( \cD ^{(m)} _{X/T} \boxtimes ^{\mathrm{top}} _{T} \cD ^{(m)} _{Y/T}).
\end{equation}
By composing 
\ref{XtimesY-boxtimesD-mor5} with \ref{boxtimestop-dfn2-prebis}, 
we get the morphism 
\begin{equation}
\label{boxtimestop-dfn2-preter}
pr _1 ^{*} \cD ^{(m)} _{X/T}
\otimes _{\cO _Z}
pr _2 ^{*} \cD ^{(m)} _{Y/T}
=
\left ( 
\cO _Z
\otimes _{pr _1 ^{-1} \O _{X}}
pr _1 ^{-1} \cD ^{(m)} _{X/T}
 \right)
\otimes _{\cO _Z}
\left ( 
\cO _Z
\otimes _{pr _2 ^{-1} \O _{Y}}
pr _2 ^{-1} \cD ^{(m)} _{Y/T}
 \right)
\to 
\cD ^{(m)} _{Z/T}.
\end{equation}
By a local computation with finite $p$-bases, we can check that the map \ref{boxtimestop-dfn2-preter}
is an isomorphism of left $\cD ^{(m)} _{Z/T}$-modules.
This implies that 
\ref{XtimesY-boxtimesD-mor5} is an isomorphism  
of $(\cO _X , \cD ^{(m)} _{X/T} \boxtimes ^{\mathrm{top}} _{T} \cD ^{(m)} _{Y/T})$-bimodules.

\end{empt}

\subsection{Exterior tensor products on schemes}
\label{ntn-ext-tens-prd}
Fix some integer $j\geq 0$ and put $S = S _j$. 
Since the base scheme $S$ is fixed, so we can remove it in the notation.
If $\phi \colon S'\to S$ is a morphism, 
by abuse of notation, we sometimes denote $\phi ^{-1} \O _S$ simply by
$\O _S$. Moreover, $S$-schemes will be supposed to be quasi-compact and separated.

For any $i\in \{ 1,\dots, n\}$, 
let 
$ p _i \colon X _i \to \bbD ^{r _i} _S$ 
and 
$q _i \colon Y _i \to \bbD ^{r _i} _S$ 
be two finite type morphisms,
let $(f _i, id) \colon p _i \to q _i$ be a morphism of 
$\fC _S$ (see notation \ref{ntn-CT}). 
We suppose 
$X _i /S$ and $Y _i/S$ have locally finite $p$-bases.

Set $X := X _1 \times _{\fC _S} X _2 \times _{\fC _S} \dots \times _{\fC _S} X _n$,
$Y := Y _1 \times _{\fC _S} Y _2 \times _{\fC _S} \dots \times _{\fC _S} Y _n$ (see notation \ref{ntnXCTY}),
and 
$f: = f _1 \times _{\fC _S}  \dots \times  _{\fC _S} f _n \colon X \to Y$.
For $i= 1,\dots, n$, let $pr _i  \colon X \to X _i$,
$pr '_i  \colon Y \to Y _i$ be the projections.
Following \ref{corXtimeCY2XtimesY}, the projections $pr _i$ and $pr ' _i$ are flat and have locally finite $p$-bases.
We denote by 
 $\varpi \colon X \to S$, $\varpi _i \colon X _i \to S$,
 $\varpi ' \colon Y \to S$ and $\varpi _i ' \colon Y _i \to S$ the structural morphisms.

Notice that since 
$X _i /S$ and $Y _i /S$ are formally smooth and of formal finite type, then 
$X _i /S$ and $Y _i /S$ are flat (see \ref{f0formétale-fforméta0}). 
Remark also that  $f \colon X \to Y$ 
is a morphism of $S$-schemes of formal finite type and having
locally finite $p$-bases.
Moreover, $f$ and $f _i$ are morphisms of finite type of noetherian schemes of finite Krull dimension.

\begin{rem}
\label{n=2,f2=id-diag-parag}
Suppose $n = 2$ and $f _2 $ is the identity. 
In that case, denoting by $Z _2:= X _2=Y_2$, we get the cartesian square
\begin{equation}
\label{n=2,f2=id-diag}
\xymatrix@ C=2cm {
{X= X _1 \times _{\fC _S} Z _2} 
\ar[d] ^-{f= f _1 \times _{\fC _S} id}
\ar@{}[rd] ^-{}| \square
\ar@/^{0,5cm}/[rr] ^-{pr _1}
\ar[r] ^-{\sim}
&
{ X _1 \times _{\bbD ^{r _1} _S} 
 \bbD ^{r _1+ r _2} _S \times _{\bbD ^{r _2} _S} Z _2} 
\ar[d] ^-{f= f _1 \times id \times id}
\ar[r] ^-{}
\ar@{}[rd] ^-{}| \square
&
{X _1} 
\ar[d] ^-{f _1}
\\ 
{Y= Y _1 \times _{\fC _S} Z _2} 
\ar@/_{0,5cm}/[rr] _-{pr '_1}
\ar[r] ^-{\sim}
&
{Y _1 \times _{\bbD ^{r _1} _S} 
 \bbD ^{r _1+ r _2} _S \times _{\bbD ^{r _2} _S} Z _2} 
\ar[r] ^-{}
 & 
{Y _1.} 
}
\end{equation}
Since $f _1$ is a morphism of finite type of noetherian schemes, 
since $pr _1$ and $pr ' _1$ are flat (see \ref{corXtimeCY2XtimesY}.\ref{corXtimeCY2XtimesY-item2}),
then 
for any $\E _1 \in D _{\mathrm{qc}} ( \cO _{X _1})$, 
the canonical base change morphism
\begin{equation}
\label{n=2,f2=id-diag-parag-iso}
pr ^{\prime *} _1 \R f _{1*}( \E _1)
\to
\R f _{*} pr ^{*} _1  (\E _1)
\end{equation}
 is an isomorphism (see \cite[II.5.12]{HaRD}).
\end{rem}

\begin{empt}
We define below similar to \cite[7.1.2]{caro-6operations} definitions.
\begin{enumerate}[(a)]
\item For $i= 1,\dots, n$, let 
$\E _i $ be a sheaf of $\varpi _i ^{-1}\O _{S}$-module.
We get the $\varpi ^{-1}\O _{S}$-module by setting
$$\underset{i}{\boxtimes} ^{\mathrm{top}}  \E _i
:=
pr _1 ^{-1} \E _1 
\otimes _{\O _S}
pr _2 ^{-1} \E _2
\otimes _{\O _S}
\cdots
\otimes _{\O _S}
pr _n ^{-1} \E _n.$$

\item For $i= 1,\dots, n$, let 
$\E _i $ be an $\O _{X _i}$-module.
The sheaf 
$\underset{i}{\boxtimes} ^{\mathrm{top}}  \E _i$
has a canonical structure of 
$\underset{i}{\boxtimes} ^{\mathrm{top}} 
\O _{X _i}$-module.
We put
$\underset{i}{\boxtimes}  \E _i
:=
\O _{X}
\otimes ^{\bbL} _{\underset{i}{\boxtimes} ^{\mathrm{top}}  \O _{X _i}}
\underset{i}{\boxtimes} ^{\mathrm{top}}  \E _i$. 
Moreover, by commutativity and associativity of tensor products,
we get the canonical isomorphism of 
$\underset{i}{\boxtimes} ^{\mathrm{top}} 
\O _{X _i}$-modules
\begin{equation}
\label{boxtimestop-dfn2}
\underset{i}{\boxtimes} ^{\mathrm{top}}  \E _i
\riso 
\left ( pr _1 ^{-1} \E _1 
\otimes _{pr _1 ^{-1} \O _{X _1}}
\underset{i}{\boxtimes} ^{\mathrm{top}}  \O _{X _i} \right)
\otimes _{\underset{i}{\boxtimes} ^{\mathrm{top}}  \O _{X _i} }
\cdots
\otimes _{\underset{i}{\boxtimes} ^{\mathrm{top}}  \O _{X _i} }
\left ( pr _n ^{-1} \E _n 
\otimes _{pr _n ^{-1} \O _{X _n}}
\underset{i}{\boxtimes} ^{\mathrm{top}}  \O _{X _i} \right).
\end{equation}
 Using the isomorphism \ref{boxtimestop-dfn2}, 
we get the isomorphism of $\O _X$-modules
\begin{equation}
\label{boxtimes-dfn2}
\underset{i}{\boxtimes}  \E _i 
\riso 
pr _1 ^{*} \E _1 
\otimes _{\O _{X} }
\cdots
\otimes _{\O _{X} }
 pr _n ^{*} \E _n .
\end{equation}

Since $pr _i ^{-1} \D ^{(m)} _{X _i}$ are $\O _S$-algebras,
we get a canonical structure of $\O _S$-algebra on
$\underset{i}{\boxtimes} ^{\mathrm{top}} 
\D ^{(m)} _{X _i}$.

\item For $i= 1,\dots, n$, 
let $\cF _i$ be a left  $\D ^{(m)} _{X _i}$-module 
(resp. $\cG _i$ be a right  $\D ^{(m)} _{X _i}$-module).
Then 
$\underset{i}{\boxtimes} ^{\mathrm{top}}  \cF _i$
(resp. $\underset{i}{\boxtimes} ^{\mathrm{top}}  \cG _i$)
has a canonical structure of left (resp. right)
$\underset{i}{\boxtimes} ^{\mathrm{top}} 
\D ^{(m)} _{X _i}$-module.
The canonical homomorphism 
of $\O _S$-algebras
$\underset{i}{\boxtimes} ^{\mathrm{top}} 
\D ^{(m)} _{X _i}
\to \D ^{(m)} _{X }$ induces the canonical isomorphism 
of $\O _X$-modules
$\underset{i}{\boxtimes} 
\D ^{(m)} _{X _i}
\riso \D ^{(m)} _{X }$ (this was checked in \ref{XtimesY-boxtimesD-mor5} in the case where $n=2$,
but the proof is identical).
This yields the isomorphism of $\O _{X }$-modules
$\underset{i}{\boxtimes}  \cF _i 
\riso 
\D ^{(m)} _{X }
\otimes _{\underset{i}{\boxtimes} ^{\mathrm{top}}  \D ^{(m)}  _{X _i}}
\underset{i}{\boxtimes}  ^{\mathrm{top}} \cF _i $
(resp. 
$\underset{i}{\boxtimes}  \cG _i 
\riso 
\underset{i}{\boxtimes}  ^{\mathrm{top}} \cF _i 
\otimes _{\underset{i}{\boxtimes} ^{\mathrm{top}}  \D ^{(m)}  _{X _i}}
\D ^{(m)} _{X }$).
Via this isomorphism, we endowed 
$\underset{i}{\boxtimes}  \cF _i $
(resp. $\underset{i}{\boxtimes}  \cG _i $)
with a structure of left (resp. right) $\D ^{(m)} _{X }$-module.

\item For $i= 1,\dots, n$, 
let $\cF _i$ be a left  $\D ^{(m)} _{X _i}$-module.
Then 
$pr _1 ^{*} \cF _1 
\otimes _{\O _{X} }
\cdots
\otimes _{\O _{X} }
pr _n ^{*} \cF _n$
has a canonical structure of left $\D ^{(m)} _{X}$-module (see \cite[2.3.3]{Be1}).  
By making a local computation with finite $p$-bases, 
we can check that the isomorphism 
\ref{boxtimes-dfn2} is in fact an isomorphism of 
left $\D ^{(m)} _{X}$-modules.

\end{enumerate}
\end{empt}

\begin{empt}
\label{rem-ext-prod}
\begin{enumerate}[(a)]
\item When $S$ is the spectrum of a field, 
the multi-functor $\underset{i}{\boxtimes} ^{\mathrm{top}}  $ is exact. 
This is not clear if the extensions
$\underset{i}{\boxtimes} ^{\mathrm{top}}  \O _{X _i}
\to \O _{X}$ 
and 
$\underset{i}{\boxtimes} ^{\mathrm{top}}  \D ^{(m)}  _{X _i}
\to 
\D ^{(m)} _{X }$ are right and left flat. 
But, following \ref{corXtimeCY2XtimesY}.\ref{corXtimeCY2XtimesY-item3}, 
the multi-functor $\underset{i}{\boxtimes}  $ is also exact when $S$ is the spectrum of a field.

\item When $S$ is not the spectrum of a field, 
the multi-functor $\underset{i}{\boxtimes} ^{\mathrm{top}}  $ is not necessarily exact. 
We get the multi-functor 
$\underset{i}{\overset{\L}{\boxtimes}} {}^{\mathrm{top}}  
\colon 
D ^{-} (\varpi _1 ^{-1}\O _{S})
\times \cdots \times  
D ^{-} (\varpi _n ^{-1}\O _{S})
\to 
D ^{-} (\varpi ^{-1}\O _{S})$ by setting 
for any 
$\E _i \in D ^{-} (\varpi _i ^{-1}\O _{S})$
$$\underset{i}{\overset{\L}{\boxtimes}} {}^{\mathrm{top}}   \E _i
:=
pr _1 ^{-1} \E _1 
\otimes _{\O _S} ^{\L}
pr _2 ^{-1} \E _2
\otimes _{\O _S} ^{\L}
\cdots
\otimes _{\O _S} ^{\L}
pr _n ^{-1} \E _n.$$

\item We have the multi-functor 
$\underset{i}{\overset{\L}{\boxtimes}}
\colon 
D ^{-} (\O _{X _1})
\times \cdots \times  
D ^{-} (\O _{X _n})
\to 
D ^{-} (\O _X)$ 
by setting for any 
$\E _i \in D ^{-} (\O _{X _i})$
\begin{equation}
\label{boxtimes-dfn2L}
\underset{i}{\overset{\L}{\boxtimes}}
 \E _i
:=
\O _{X}
\otimes ^{\bbL} _{\underset{i}{\boxtimes} ^{\mathrm{top}}  \O _{X _i}}
\underset{i}{\overset{\L}{\boxtimes}} {}^{\mathrm{top}} 
\cE _i
\riso 
pr _1 ^{*} \E _1 
\otimes _{\O _{X} } ^{\L}
\cdots
\otimes _{\O _{X} } ^{\L}
 pr _n ^{*} \E _n ,
\end{equation}
where the last isomorphism is, after using flat resolutions,
a consequence of \ref{boxtimes-dfn2}.

\item For any $i = 1,\dots, n$, let 
$\cF _i \in D ^{-} ({}^{l} \D ^{(m)} _{X _i})$,
$\cM _i \in D ^{-} ({}^{r} \D ^{(m)} _{X _i})$.
Since 
we have the canonical isomorphisms 
$\underset{i}{\overset{\L}{\boxtimes}}
\D ^{(m)} _{X _i}
\riso 
\underset{i}{\boxtimes} 
\D ^{(m)} _{X _i}
\riso \D ^{(m)} _{X }$,
then the canonical morphisms
\begin{gather}
\notag
\O _{X}
\otimes ^{\bbL} _{\underset{i}{\boxtimes} ^{\mathrm{top}}  \O _{X _i}}
\underset{i}{\overset{\L}{\boxtimes}} {}^{\mathrm{top}} 
\cF _i
\to 
\D ^{(m)} _{X }
\otimes _{\underset{i}{\boxtimes} ^{\mathrm{top}}  \D ^{(m)}  _{X _i}}
^{\bbL}
\underset{i}{\overset{\L}{\boxtimes}} {}^{\mathrm{top}} 
\cF _i,
\\
\notag
\O _{X}
\otimes ^{\bbL} _{\underset{i}{\boxtimes} ^{\mathrm{top}}  \O _{X _i}}
\underset{i}{\overset{\L}{\boxtimes}} {}^{\mathrm{top}} 
\cM _i
\riso 
\underset{i}{\overset{\L}{\boxtimes}} {}^{\mathrm{top}} 
\cM _i
\otimes _{\underset{i}{\boxtimes} ^{\mathrm{top}}  \O _{X _i}} ^{\bbL}
\O _{X}
\to 
\underset{i}{\overset{\L}{\boxtimes}} {}^{\mathrm{top}} 
\cM _i
\otimes _{\underset{i}{\boxtimes} ^{\mathrm{top}}  \D ^{(m)}  _{X _i}} ^{\bbL} 
\D ^{(m)} _{X }
\end{gather}
are isomorphisms. 
Hence, there is no problem (up to canonical isomorphism) with respect to \ref{boxtimes-dfn2L}  to set 
$\underset{i}{\overset{\L}{\boxtimes}}
 \cF _i
:=
\D ^{(m)} _{X }
\otimes _{\underset{i}{\boxtimes} ^{\mathrm{top}}  \D ^{(m)}  _{X _i}}
\underset{i}{\overset{\L}{\boxtimes}} {}^{\mathrm{top}} 
\cF _i$
and
$\underset{i}{\overset{\L}{\boxtimes}}
 \cM _i
:=
\underset{i}{\overset{\L}{\boxtimes}} {}^{\mathrm{top}} 
\cM _i
\otimes _{\underset{i}{\boxtimes} ^{\mathrm{top}}  \D ^{(m)}  _{X _i}}
\D ^{(m)} _{X }$.
For $ * \in \{ l,r\}$, we get the multi-functor
$$\underset{i}{\overset{\L}{\boxtimes}}
\colon 
D ^{-} ({} ^* \D ^{(m)} _{X _1})
\times \cdots \times  
D ^{-} ({} ^* \D ^{(m)} _{X _n})
\to 
D ^{-} ({} ^* \D ^{(m)} _{X }).$$ 

\item If we would like to clarify the basis $S$, 
we may add it in the notation. For instance,
we write 
$\underset{S ,i}{\overset{\L}{\boxtimes}} {}^{\mathrm{top}} $
and 
$\underset{S,i}{\overset{\L}{\boxtimes}}$
(or $\underset{\O _S ,i}{\overset{\L}{\boxtimes}} {}^{\mathrm{top}} $
and 
$\underset{\O _S,i}{\overset{\L}{\boxtimes}}$)
instead of
$\underset{i}{\overset{\L}{\boxtimes}} {}^{\mathrm{top}} $
and 
$\underset{i}{\overset{\L}{\boxtimes}}$.

\end{enumerate}

\end{empt}

\begin{lem}
\label{com-botimestop}
For $i= 1,\dots, n$, 
let $\D _i$ be a sheaf of $\varpi _i ^{-1} \O _S$-algebras,
$\cM _i\in D ^{-} (\D _i, \O _{X _i}) $,
$\E _i \in D ^- (\O _{X _i})$,
$\cN _i \in D ^{-} (\D _i,  \D ^{(m)} _{X _i})$,
$\cF _i \in D ^-(\D ^{(m)}  _{X _i})$.
\begin{enumerate}[(a)]
\item We have the canonical isomorphism of 
$D ^- (\underset{i}{\boxtimes} ^{\mathrm{top}} \D _i, 
\underset{i}{\boxtimes} ^{\mathrm{top}}  \O _{X _i} )$
\begin{equation}
\label{com-botimestop-iso1}
\underset{i}{\overset{\L}{\boxtimes}} {}^{\mathrm{top}}   
(\cM _i \otimes ^{\L} _{\O _{X _i}} \E _i)
\riso 
\underset{i}{\overset{\L}{\boxtimes}} {}^{\mathrm{top}}   
\cM _i \otimes ^{\L} _{\underset{i}{\boxtimes} ^{\mathrm{top}}  \O _{X _i}} 
\underset{i}{\overset{\L}{\boxtimes}} {}^{\mathrm{top}}   \E _i.
\end{equation}

\item We have the canonical isomorphism of $\underset{i}{\boxtimes} ^{\mathrm{top}} \D _i$-modules
\begin{equation}
\label{com-botimestop-iso2}
\underset{i}{\overset{\L}{\boxtimes}} {}^{\mathrm{top}}   
(\cN _i \otimes ^{\L} _{\D ^{(m)}  _{X _i}} \cF _i)
\riso 
\underset{i}{\overset{\L}{\boxtimes}} {}^{\mathrm{top}}   
\cN _i \otimes ^{\L} _{\underset{i}{\boxtimes} ^{\mathrm{top}}  \D ^{(m)}  _{X _i}} \underset{i}{\overset{\L}{\boxtimes}} {}^{\mathrm{top}}   \cF _i.
\end{equation}

\end{enumerate}
\end{lem}

\begin{proof}
We can copy the proof of \cite[7.1.4]{caro-6operations}.
\end{proof}

\begin{lemm}
\label{com-botimes}
For $i= 1,\dots, n$, 
let $\D _i$ be a sheaf of $\varpi _i ^{-1} \O _S$-algebras.

\begin{enumerate}[(i)]
\item 
\label{com-botimes1} 
For $i= 1,\dots, n$, 
for $ * \in \{ l,r\}$, let 
$\cM _i\in D ^{-} ({}^*\D _i, \O _{X _i}) $,
$\E _i \in D ^- (\O _{X _i})$.
We have the canonical isomorphism of the form
$\underset{i}{\overset{\L}{\boxtimes}} 
(\cM _i \otimes ^{\L} _{\O _{X _i}} \E _i)
\riso
\underset{i}{\overset{\L}{\boxtimes}}  
\cM _i \otimes ^{\L} _{\O _X} \underset{i}{\overset{\L}{\boxtimes}}   \E _i$ 
of 
$D ^- ({}^*\underset{i}{\boxtimes} ^{\mathrm{top}}\D _i, \O _{X} )$.
Moreover, this isomorphism is compatible with that of \ref{com-botimestop-iso1}, i.e. 
the following diagram 
of $D ^- ({}^*\underset{i}{\boxtimes} ^{\mathrm{top}} \D _i, 
\underset{i}{\boxtimes} ^{\mathrm{top}} \O _{X _i})$
\begin{equation}
\label{com-botimes-diag2}
\xymatrix {
{\underset{i}{\overset{\L}{\boxtimes}} {}^{\mathrm{top}}   
(\cM _i \otimes ^{\L} _{\O _{X _i}} \E _i)}
\ar[r] ^-{\sim}  _-{\ref{com-botimestop-iso1}}
\ar[d] ^-{}
&
{\underset{i}{\overset{\L}{\boxtimes}} {}^{\mathrm{top}}   
\cM _i \otimes ^{\L} _{\underset{i}{\boxtimes} ^{\mathrm{top}}  \O _{X _i}} 
\underset{i}{\overset{\L}{\boxtimes}} {}^{\mathrm{top}}   \E _i} 
\ar[d] 
\\ 
{\underset{i}{\overset{\L}{\boxtimes}} 
(\cM _i \otimes ^{\L} _{\O _{X _i}} \E _i)} 
\ar[r] ^-{\sim}  
&
{\underset{i}{\overset{\L}{\boxtimes}}  
\cM _i \otimes ^{\L} _{\O _X} \underset{i}{\overset{\L}{\boxtimes}}   \E _i}
}
\end{equation}
is commutative.
\item 
\label{com-botimes2} 
For $i= 1,\dots, n$, 
for $ * \in \{ l,r\}$, let
$\cM _i \in D ^{-} ({}^*  \D _i,  {}^l \D ^{(m)} _{X _i})$,
$\cE _i \in D ^-({}^l \D ^{(m)}  _{X _i})$.
Then, 
the isomorphism
$\underset{i}{\overset{\L}{\boxtimes}} 
(\cM _i \otimes ^{\L} _{\O _{X _i}} \cE _i)
\riso
\underset{i}{\overset{\L}{\boxtimes}}  
\cM _i \otimes ^{\L} _{\O _X} \underset{i}{\overset{\L}{\boxtimes}}   \cE _i$ 
constructed in  \ref{com-botimes}.(\ref{com-botimes1})
is in fact an isomorphism
of 
$D ^- ({}^*  \underset{i}{\boxtimes} ^{\mathrm{top}} \D _i,{}^l \D ^{(m)}  _{X} )$.

\end{enumerate}

\end{lemm}

\begin{proof}
We can copy the proof of \cite[7.1.5]{caro-6operations}.
\end{proof}

\begin{empt}
\label{boxtimes-left-right}
It follows from \ref{XtimeCY2XtimesY} that the canonical morphism
$X \to 
X _1 \times _{S} X _2 \times _{S} \dots \times _{S} X _n$ is relatively perfect
and in particular is formally étale. 
This implies that 
the canonical morphism 
$\oplus _{i=1} ^n pr _i ^*  \Omega ^1 _{X _i} \to \Omega ^1 _{X}$
is an isomorphism. By applying determinants,
this yields the isomorphism of $\O _X$-modules
$\underset{i}{\boxtimes}~  \omega _{X _i} \riso \omega  _{X}$.
Using the canonical structure of 
right $\D  ^{(m)} _{X _i}$-module on $\omega _{X _i}$, we get a structure of 
right $\D  ^{(m)} _{X}$-module on $\underset{i}{\boxtimes} \omega _{X _i}$.
By local computations, we can check the canonical isomorphism
$\underset{i}{\boxtimes} \omega _{X _i} \riso \omega  _{X}$ 
is in fact an isomorphism of right $\D  ^{(m)} _{X}$-modules.

For $i= 1,\dots, n$, 
$\E _i $ be a left $\D ^{(m)}  _{X _i}$-module, and 
$\cF _i $ be a right $\D ^{(m)}  _{X _i}$-module.
Then we have the canonical morphism of right $\D  ^{(m)} _{X}$-modules (resp. left $\D  ^{(m)} _{X}$-modules)
$\underset{i}{\boxtimes} ( \omega _{X _i} \otimes _{\O _{X _i}} \cE _i)
\riso \omega _{X } \otimes _{\O _{X }} \underset{i}{\boxtimes}  \cE _i$
(resp. $\underset{i}{\boxtimes}  ( \cF _i \otimes _{\O _{X _i}} \omega _{X _i} ^{-1} )
\riso \underset{i}{\boxtimes}  \cF _i \otimes _{\O _{X}} \omega _{X} ^{-1} $).
Taking flat resolutions, we have similar isomorphisms in derived categories.

\end{empt}

\subsection{Commutation with pullbacks and push forwards}
We keep the notation of the section \ref{ntn-ext-tens-prd}.
\begin{prop}
\label{comm-boxtimes-f*}
For $i= 1,\dots, n$, 
let $\D _i$ be a sheaf 
of $\varpi _i ^{\prime -1} \O _S$-algebras,
$\cN _i \in D ^- (\D _i,  \D ^{(m)} _{Y _i})$.
We have the canonical isomorphism of 
$D ^- (\underset{i}{\boxtimes} ^{\mathrm{top}} \D _i,\D ^{(m)}  _{X _i} )$:
$$\L f ^* ( \underset{i}{\overset{\L}{\boxtimes}}  \cN _i )
\riso 
\underset{i}{\overset{\L}{\boxtimes}}  \,
\L f _i ^* (  \cN _i ) .$$
\end{prop}

\begin{proof}
We can copy the proof of \cite[7.2.4]{caro-6operations}.
\end{proof}

\begin{thm}
\label{sch-prop-boxtimes-v*}
For $i= 1,\dots, n$, let
$\cE _i\in D ^{\mathrm{b}} _{\mathrm{qc}} (\O _{X _i})$.
The canonical morphism
\begin{equation}
\label{sch-prop-boxtimes-v*iso}
\underset{i}{\overset{\L}{\boxtimes}}  
\R f _{i*} (\cE _i )
\to
\R f _* (\underset{i}{\overset{\L}{\boxtimes}}   \cE _i ).
\end{equation}
is an isomorphism.
\end{thm}

\begin{proof}
By copying word by word the proof of \cite[7.2.8]{caro-6operations}, the theorem is 
a consequence of the base change isomorphism 
\ref{n=2,f2=id-diag-parag-iso}.
\end{proof}

\begin{ntn}
\label{ntn-Tf}
Let 
$\cE _i\in D ^{\mathrm{b}} _{\mathrm{qc}} (\D ^{(m)} _{X _i})$.
We denote by 
$\mathrm{T} _{f _i}
\colon 
\colon \R f _{i,*}  (\cE _i)
\to
\R f _{i,*}  (
 \cD _{Y _i  \leftarrow X _i} ^{(m)} \otimes _{  \cD _{X _i} ^{(m)}} ^{\L}
\cE _i)
=
f _{i,+} ^{ (m)}(\cE_i)$,
the canonical morphism induced by
the homomorphism 
$\cD _{X _i} ^{(m)}
\to 
\cD _{Y _i  \leftarrow X _i} ^{(m)}  $
given by the left $\cD _{X _i} ^{(m)}$-module structure of 
$\cD _{Y _i  \leftarrow X _i} ^{(m)}  $.
Similarly for $f$. 
\end{ntn}

\begin{thm}
\label{sch-prop-boxtimes-v+}
For $i= 1,\dots, n$, let 
$\cE _i\in D ^{\mathrm{b}} _{\mathrm{qc}} (\D ^{(m)} _{X _i})$.
We have the canonical isomorphism
\begin{equation}
\label{sch-prop-boxtimes-v+iso}
\underset{i}{\overset{\L}{\boxtimes}}  f _{i+} ^{ (m)} (\cE _i )
\riso
f _+ ^{ (m)} (\underset{i}{\overset{\L}{\boxtimes}}   \cE _i )
\end{equation}
making commutative the canonical diagram
\begin{equation}
\label{sch-prop-boxtimes-v+iso2}
\xymatrix @ R=0,3cm{
{\underset{i}{\overset{\L}{\boxtimes}}  
\R f _{i*} (\cE _i )} 
\ar[r] ^-{\ref{sch-prop-boxtimes-v*iso}} _-{\sim}
\ar[d] ^-{\underset{i}{\overset{\L}{\boxtimes}} \mathrm{T} _{f_i}}
_-{\ref{ntn-Tf}}
& 
{\R f _* (\underset{i}{\overset{\L}{\boxtimes}}   \cE _i )} 
\ar[d] ^-{\mathrm{T} _{f}}
_-{\ref{ntn-Tf}}
\\ 
{\underset{i}{\overset{\L}{\boxtimes}}  f _{i+} ^{ (m)} (\cE _i )} 
\ar[r] ^-{\sim}
& 
{f _+ ^{ (m)} (\underset{i}{\overset{\L}{\boxtimes}}   \cE _i ). } 
}
\end{equation}

\end{thm}

\begin{proof}
We can copy the proof of \cite[7.2.10]{caro-6operations}, i.e. this is a consequence of \ref{sch-prop-boxtimes-v*}.
\end{proof}

\subsection{Application : base change in the projection case}

We keep notation \ref{ntn-ext-tens-prd} and we 
suppose $n = 2$ and $f _2 $ is the identity.
\begin{prop}
\label{theo-iso-chgtbase2-pre} 
For any $\cE _1\in D ^{\mathrm{b}} _{\mathrm{qc}} (\D ^{(m)} _{X _1})$, we have the canonical isomorphism
$pr  _1 ^{\prime ! (m)}  \circ f ^{ (m)} _{1,+}( \E _1)
\riso 
f ^{ (m)} _{+}  \circ pr _1 ^{ ! (m)} (\E _1)$
of 
$D ^{\mathrm{b}} _{\mathrm{qc}} (\D ^{(m)} _{Y})$
making commutative the diagram 
\begin{equation}
\label{theo-iso-chgtbase2-pre-iso1}
\xymatrix{
{pr  _1 ^{\prime*}  \circ \R f  _{1,*}( \E _1)} 
\ar[r] ^-{\sim}
\ar[d] ^-{}
& 
{\R f _{*}  \circ pr _1 ^{*} (\E _1)} 
\ar[d] ^-{}
\\ 
{pr  _1 ^{\prime *}  \circ f ^{ (m)} _{1,+}( \E _1)} 
\ar[r] ^-{\sim}
& 
{f ^{ (m)} _{+}  \circ pr _1 ^{ *} (\E _1),} 
}
\end{equation}
where the top isomorphism is the usual base change isomorphism (see \ref{n=2,f2=id-diag-parag-iso}).
\end{prop}

\begin{proof}
We can copy the proof of \cite[7.3.1]{caro-6operations}: this is an easy consequence of 
Theorem \ref{sch-prop-boxtimes-v+}.
\end{proof}

\begin{ntn}
\label{ntn-psharp}
Let $g \colon Z \to T$ be a 
flat morphism 
having locally finite $p$-bases 
of $S$-schemes of formal finite type
having locally finite $p$-bases  over $S$.

\begin{enumerate}[(a)]
\item Similarly to \cite[III.2]{HaRD},
we define a functor
$g ^\sharp \colon D ( \O _T) \to D ( \O _Z)$
by setting 
\begin{equation}
\label{ntn-psharp-dfng^sharp}
g ^\sharp (\cM ) : = g ^* ( \cM) \otimes _{\O _T} \omega _{Z/T} [ \delta _{Z/T}],
\end{equation}
where
$\delta  _{Z/T}:= \delta ^{S}_{Z} - \delta ^{S} _{T} \circ g$.

\item Let
$\cM \in D ^{\mathrm{b}} _{\mathrm{qc}} ({} ^{r}\D ^{(m)} _{T})$.
We have the isomorphisms
\begin{gather}
\notag
g ^{ ! (m)} (\cM)
\riso 
g ^{ ! (m)} (
\omega  _{T / S}
\otimes _{\O _{T }}
\cM
\otimes _{\O _{T }} \omega ^{-1} _{T / S})
\underset{\ref{fund-isom2bisprepre}}{\riso} 
\omega  _{Z / S}
\otimes _{\O _{Z }}
g ^{ ! (m)} ( \cM \otimes _{\O _{T }} \omega ^{-1} _{T / S})
\\
\label{ntn-psharp-dfng^sharp-iso}
\underset{\ref{inv-image-def0-iso1}}{\riso} 
\omega  _{Z / S}
\otimes _{\O _{Z }}
g ^{*} ( \cM \otimes _{\O _{T }} \omega ^{-1} _{T / S})
[\delta _{Z/T}]
\riso 
g ^* ( \cM) \otimes _{\O _T} \omega _{Z/T} [ \delta _{Z/T}]
=
g ^\sharp (\cM ) ,
\end{gather}
where the last isomorphism is a consequence of 
$\omega _{Z/T}
\riso 
\omega _{Z/S}
\otimes _{\cO _Z}
g ^* (\omega _{T/S} ^{-1})$.
\end{enumerate}

\end{ntn}

\begin{prop}
\label{theo-iso-chgtbase2-pre2}
We keep notation \ref{ntn-psharp}.
\begin{enumerate}[(a)]
\item For any $\cM _1\in D ^{\mathrm{b}} _{\mathrm{qc}} (\O  _{X _1})$, 
we have the  isomorphism
\begin{equation}
\label{theo-iso-chgtbase2-pre-isopre}
pr  _1 ^{\prime \sharp}  \circ \R f  _{1,*}( \cM _1) \riso \R f _{*}  \circ pr _1 ^{\sharp} (\cM _1)
\end{equation}
of 
$D ^{\mathrm{b}} _{\mathrm{qc}} (\O _{Y})$
canonically induced by the usual base change isomorphism.

\item For any $\cM _1\in D ^{\mathrm{b}} _{\mathrm{qc}} ({} ^r \D ^{(m)} _{X _1})$, 
we have the  isomorphism
 the canonical
$pr  _1 ^{\prime ! (m)}  \circ f ^{ (m)} _{1,+}( \cM _1)
\riso 
f ^{ (m)} _{+}  \circ pr _1 ^{ ! (m)} (\cM _1)$
of 
$D ^{\mathrm{b}} _{\mathrm{qc}} ({} ^r \D ^{(m)} _{Y})$
making commutative the diagram 
\begin{equation}
\label{theo-iso-chgtbase2-pre-iso}
\xymatrix{
{pr  _1 ^{\prime \sharp}  \circ \R f  _{1,*}( \cM _1)} 
\ar[r] ^-{\sim} _-{\ref{theo-iso-chgtbase2-pre-isopre}}
\ar[d] ^-{}
& 
{\R f _{*}  \circ pr _1 ^{\sharp} (\cM _1)} 
\ar[d] ^-{}
\\ 
{pr  _1 ^{\prime  ! (m)}  \circ f ^{ (m)} _{1,+}( \cM _1)} 
\ar[r] ^-{\sim}
& 
{f ^{ (m)} _{+}  \circ pr _1 ^{ ! (m)} (\cM _1).} 
}
\end{equation}

\end{enumerate}

\end{prop}

\begin{proof}
We can copy the proof of \cite[7.3.3]{caro-6operations}
(we have only to replace ``$[ d _{Z/T}]$'' by ``$[ \delta _{Z/T}]$'' and sometimes 
``smooth'' 
by ``having locally finite $p$-bases, flat and of formal finite type'').
\end{proof}

\subsection{Application : relative duality isomorphism and adjunction for projective morphisms}
We keep notation \ref{ntn-ext-tens-prd}, we 
suppose $n = 2$, $f _2 $ is the identity, 
$X _1 = \bbP ^{d} _{Y _1}$,
$f _1 \colon \bbP ^{d} _{Y _1} \to Y _1$ is the canonical projection.
We set $T:= X _2=Y_2$.

\begin{lem}
With notation \ref{ntn-Tf} and \ref{ntn-psharp}, 
for any $\cN _1\in D ^{\mathrm{b}} _{\mathrm{qc}} (\O  _{Y _1})$, 
we have the commutative diagram
\begin{equation}
\label{diag-comm-Trf1-Trf}
\xymatrix{
{pr  _1 ^{\prime \sharp}  \circ \R f  _{1,*} \circ f _1 ^{\sharp} ( \cN _1)} 
\ar[r] ^-{\sim} _-{\ref{theo-iso-chgtbase2-pre-isopre}}
\ar[d] ^-{\mathrm{Tr} _{f _1}}
& 
{\R f  _{*} \circ  pr  _1 ^{\sharp}  \circ f _1 ^{\sharp} ( \cN _1)} 
\ar[r] ^-{\sim}
& 
{\R f  _{*} \circ  f  ^{\sharp}   \circ pr  _1 ^{\prime \sharp} ( \cN _1)} 
\ar[d] ^-{\mathrm{Tr} _{f}}
\\ 
{pr  _1 ^{\prime \sharp}  ( \cN _1)} 
\ar@{=}[rr] ^-{}
&& 
{ pr  _1 ^{\prime \sharp} ( \cN _1),} 
}
\end{equation}
where $\mathrm{Tr} _{f}$ and 
$\mathrm{Tr} _{f _1}$ are the trace map isomorphisms
(see \cite[III.4.3]{HaRD}).
\end{lem}

\begin{proof}
We can copy the proof of \cite[7.4.1]{caro-6operations}.
\end{proof}

\begin{prop}
\label{thm-tracemap-Virrion}
Let $\cN _1\in D ^{\mathrm{b}} _{\mathrm{qc}} ({} ^r \D ^{(m)} _{Y _1})$.
Suppose we have 
 the canonical morphism 
$\mathrm{Tr} _{+,f _1}
\colon 
f ^{ (m)} _{1,+}\circ f _1 ^{!(m)} ( \cN _1)
\to 
\cN _1$
of 
$D ^{\mathrm{b}} _{\mathrm{qc}} ({} ^r \D ^{(m)} _{Y_1})$
making commutative the diagram 
\begin{equation}
\label{theo-iso-chgtbase2-pre-iso3pre}
\xymatrix{
{ \R f  _{1,*}\circ f _1 ^{\sharp} ( \cN _1)} 
\ar[d] ^-{}
\ar[r] ^-{\mathrm{Tr} _{f _1}}
&
{ \cN _1}
\\
{ f ^{ (m)} _{1,+}\circ f _1 ^{!(m)} ( \cN _1).} 
\ar[ur] _-{\mathrm{Tr} _{+,f _1}}
}
\end{equation}

Then, there exists a  canonical morphism 
$\mathrm{Tr} _{+,f }
\colon 
f ^{ (m)} _{+}  \circ f  ^{! (m)} \circ pr _1 ^{\prime  ! (m)}  ( \cN _1)
\to 
pr _1 ^{\prime  ! (m)}  ( \cN _1)$
of 
$D ^{\mathrm{b}} _{\mathrm{qc}} ({} ^r \D ^{(m)} _{Y})$
making commutative the diagram 
\begin{equation}
\label{thm-tracemap-Virrion-diag}
\xymatrix{
{\R f  _{*} \circ  f  ^{\sharp}   \circ pr  _1 ^{\prime \sharp} ( \cN _1)} 
\ar[r] ^-{\mathrm{Tr} _{f}}
\ar[d] ^-{}
&
{pr  _1 ^{\prime \sharp} ( \cN _1)} 
\ar[d] ^-{\sim}
\\
{f ^{ (m)} _{+}  \circ f  ^{! (m)} \circ pr _1 ^{\prime  ! (m)}  ( \cN _1)} 
\ar[r] ^-{\mathrm{Tr} _{+,f }}
&
{pr  _1 ^{\prime !(m)} ( \cN _1).} 
}
\end{equation}
\end{prop}

\begin{proof}
We can copy the proof of \cite[7.4.2]{caro-6operations}
(for instance, our schemes are noetherian, projections are flat, 
$f _1$ and $f$ are projective so we can apply 
\cite[III.10.5.Tra 4)]{HaRD} ; moreover, since
$D ^{\mathrm{b}} _{\mathrm{qc}} (\O  _{Y _1})
= 
D ^{\mathrm{b}} _{\mathrm{qc}, \mathrm{tdf}} (\O  _{Y _1})$
then we can apply \cite[III.4.4]{HaRD}).
\end{proof}

\begin{empt}
Suppose 
$Y _1 =S$, $X _1 = \bbP ^{d} _S$,
$f _1 \colon \bbP ^{d} _S \to S$ is the canonical projection
and 
$\cN _1= \O _S \in D ^{\mathrm{b}} _{\mathrm{qc}} ({} ^r \D ^{(m)} _{Y_1/S})
=
D ^{\mathrm{b}} _{\mathrm{qc}} (\O _S)$.
We have 
$f _1 ^{!(m)} ( \cO _S)=f _1 ^\sharp (\O _S) =  \omega _{\bbP ^{d} _S/S}  [d]$
and 
the trace map 
$Tr _{f _1}
\colon 
\R f _* ( \omega _{\bbP ^{d} _S/S} ) [d]
\to \O _S$ 
is an isomorphism of 
$D ^{\mathrm{b}} _{\mathrm{qc}} (\O _S)$.
Since the canonical morphism 
$\R f _* ( \omega _{\bbP ^{d} _S/S} ) [d]
\to 
f ^{ (m)} _{1,+}
( \omega _{\bbP ^{d} _S/S} ) [d]$ is an isomorphism after applying the trunctation
functor 
$\tau _{\geq 0}$, 
we get the morphism
$\mathrm{Tr} _{+,f _1}
\colon 
f ^{ (m)} _{1,+}
( \omega _{\bbP ^{d} _S/S} ) [d]
\to 
 \O _S$ making commutative the diagram
\begin{equation}
\label{theo-iso-chgtbase2-ex}
\xymatrix{
{ \R f  _{1,*}( \omega _{\bbP ^{d} _S/S} ) [d]} 
\ar[d] ^-{}
\ar[r] ^-{\mathrm{Tr} _{f _1}}
&
{ \cO _S}
\\
{ f ^{ (m)} _{1,+} ( \omega _{\bbP ^{d} _S/S} ) [d].} 
\ar[ur] _-{\mathrm{Tr} _{+,f _1}}
}
\end{equation}
Hence, following Proposition \ref{thm-tracemap-Virrion},
there exists a  canonical morphism 
$\mathrm{Tr} _{+,f }
\colon 
f ^{ (m)} _{+}  ( \omega _{\bbP ^{d} _T/S} ) [d]
\to 
( \omega _{T/S} )$
of 
$D ^{\mathrm{b}} _{\mathrm{qc}} ({} ^r \D ^{(m)} _{T/S})$
making commutative the diagram 
\begin{equation}
\label{theo-iso-chgtbase2-pre-iso3ex}
\xymatrix{
{\R f  _{*} \circ  ( \omega _{\bbP ^{d} _T/S} ) [d]} 
\ar[r] ^-{\mathrm{Tr} _{f}}
\ar[d] ^-{}
&
{\omega _{T/S}.} 
\\
{f ^{ (m)} _{+}  ( \omega _{\bbP ^{d} _T/S} ) [d].} 
\ar[ur] ^-{\mathrm{Tr} _{+,f }}
}
\end{equation}

\end{empt}

\begin{thm}
\label{rel-dual-isom-proj}
Let $f\colon X \to Y$ be a morphism 
of $S$-schemes of formal finite type and having locally finite $p$-bases.
We suppose $f$ is projective (in the strong sense), i.e. $f$ is 
the composition of a closed immersion of the form
$X \hookrightarrow \bbP ^d _Y$ with the projection $\bbP ^d _Y \to Y$. 

\begin{enumerate}[(a)]
\item Let 
$\cE \in D ^\mathrm{b}  _{\mathrm{coh}}
({} ^l \D _X ^{(m)} )$.
We have the isomorphism of $D ^\mathrm{b}  _{\mathrm{coh}}
({} ^l \D _Y ^{(m)} )$:
\begin{equation}
\label{rel-dual-isom-proj1}
\DD ^{(m)} \circ f _+ (\cE)
\riso 
 f _+ \circ \DD ^{(m)} (\cE).
\end{equation}

\item Let $\cE \in D ^\mathrm{b}  _{\mathrm{coh}}
({} ^l \D _X ^{(m)} )$,
and
$\cF \in D ^\mathrm{b}  _{\mathrm{coh}}
({} ^l \D _Y ^{(m)} )$.
We have 
the isomorphisms
\begin{gather}
\label{cor-adj-formul-proj-bij1}
\R \mathcal{H} om _{\D _Y ^{(m)}}
( f _{+} ( \E ) , \cF) 
\riso 
\R f _* 
\R \mathcal{H} om _{\D _X ^{(m)}}
( \E  ,f ^!   ( \cF)),
\\
\label{cor-adj-formul-proj-bij2}
\R \mathrm{Hom}  _{\D _Y ^{(m)}}
( f _{+} ( \E ) , \cF) 
\riso 
\R \mathrm{Hom}  _{\D _X ^{(m)}}
( \E  ,f ^!   ( \cF)). 
\end{gather}
\end{enumerate}

\end{thm}

\begin{proof}
1) Let us check \ref{rel-dual-isom-proj1}. 
Following \ref{rel-dual-isom-imm}, 
the case of a closed immersion is already checked.
Hence, we reduce to the case where $f$ is the projection $\bbP ^d _Y \to Y$.
Using \ref{theo-iso-chgtbase2-pre-iso3ex}, to check such an isomorphism,
we can copy Virrion's proof (more precisely :
a) the construction is given in \cite[IV.1.3]{Vir04},
b) for induced modules, using Grothedieck's duality isomorphism for coherent 
$\O$-modules, we construct in another way such an isomorphism : see \cite[IV.2.2.4]{Vir04},
c) the equality between both constructions is a consequence of 
the commutativity of \ref{theo-iso-chgtbase2-pre-iso3ex}: see \cite[IV.2.2.5]{Vir04}).

2) The second statement
is an easy consequence of \ref{rel-dual-isom-proj1}
(see the proof of  \cite[4.4.2]{caro-6operations}).
\end{proof}

\subsection{Going to formal $\fS$-schemes}

\begin{empt}
\label{dfn-CfS}
We give below a formal version of \ref{ntn-CT} :
\begin{enumerate}[(a)]
\item We denote by $\scr{F} _{\fS}$ the full subcategory of the category
of formal $\fS$-schemes consisting in objects $\fP/\fS$ of formal finite type, i.e. such that there exists a finite type morphism of the form 
$\fP \to \bbD ^r _\fS$ for some integer $r$.

\item We denote by $\scr{C} _{\fS}$ the category whose objects are
finite type morphisms of formal $\fS$-schemes of the form 
$\fP \to \bbD ^r _{\fS}$ for some integer $r$.
A morphism $f \to g$ of $\scr{C} _{\fS}$, where 
$f\colon \fP \to \bbD ^r _{\fS}$ and
$g\colon \fQ \to \bbD ^s _{\fS}$, consists in 
two morphisms 
$\alpha \colon 
\fP \to \fQ$ and 
$\beta \colon \bbD ^r _{\fS} \to \bbD ^s _{\fS}$
of $\scr{F} _{\fS}$ 
making commutative the diagram
\begin{equation}
\xymatrix{
{\fP} 
\ar[r] ^-{\alpha}
\ar[d] ^-{f}
& 
{\fQ} 
\ar[d] ^-{g}
\\ 
{\bbD ^r _{\fS}} 
\ar[r] ^-{\beta}
& 
{\bbD ^s _{\fS}.} 
}
\end{equation}
We notice that $\beta$ is associated to 
a local homomorphism of local rings (see \ref{rem-loc-hom-fft}). 
We denote by 
$(\alpha, \beta) $ such a morphism.

\item 
\label{dfn-CfS-item3}
We get a functor 
$\mathscr{S} _{\fS}
\colon 
\scr{C} _{\fS} \to \scr{F} _{\fS}$ defined by setting 
$\mathscr{S} _{\fS}
(
\fP \to \bbD ^r _\fS
)
= \fP$.

\end{enumerate}

\end{empt}

\begin{ntn}
\label{dfn-CfS2}
Let $f\colon \fP \to \bbD ^r _{\fS}$ and
$g\colon \fQ \to \bbD ^s _{\fS}$ be two objects of 
$\scr{C} _{\fS}$ (see notation \ref{dfn-CfS}).
We put 
\begin{equation}
\label{ntnPCTQ}
\fP \times _{\scr{C} _{\fS}} \fQ : = 
\left ( \fP \times _{\bbD ^r _{\fS}} \bbD ^{r+s} _{\fS} \right ) 
\times _{\bbD ^{r+s} _{\fS}}
\left ( \bbD ^{r+s} _{\fS} \times _{\bbD ^s _{\fS}} \fQ \right )
\riso 
 \fP \times _{\bbD ^r _{\fS}} 
 \bbD ^{r+s} _{\fS} \times _{\bbD ^s _{\fS}} \fQ.
\end{equation}
We denote by 
$f \times  _{\scr{C} _{\fS}} g 
\colon \fP \times _{\scr{C} _{\fS}} \fQ
\to 
 \bbD ^{r+s} _{\fS} $, 
$pr _1 \colon \fP  \times _{\scr{C} _{\fS}} \fQ \to \fP$
and by
$pr_2 \colon \fP \times _{\scr{C} _{\fS}} \fQ \to \fQ$
the canonical projections, 
by 
$pr _{1} ^{r ,s}\colon 
\bbD ^{r+s} _{\fS}
\to 
\bbD ^r _{\fS}$
and
$pr _{2} ^{r ,s}\colon 
\bbD ^{r+s} _{\fS}
\to \bbD ^s _{\fS}$ the canonical morphisms.
Such morphisms of the form 
$pr_1$ or $pr_2$ are called ``projection morphism''. 
Similarly to \ref{priX,Y}, 
we check that $f \times  _{\scr{C} _{\fS}}g $ equipped with the morphisms
$(pr _1 , pr _{1} ^{r ,s}) $
and 
$(pr _2 , pr _{2} ^{r  ,s}) $
of $\scr{C} _{\fS}$
satisfies the universal property of the product
in $\scr{C} _{\fS}$.
As for \ref{graphCT}, 
for any morphism 
$(\alpha, \beta) 
\colon
f \to g$, this yields the graph morphism
$(\gamma _{(\alpha, \beta)}, \gamma _\beta)$ of $\scr{C} _{\fS}$,
which is such that 
$\gamma _{(\alpha, \beta)}
\colon 
\fP 
\hookrightarrow 
\fP 
\times _{\scr{C} _{\fS}}
\fQ$
is an immersion.

\end{ntn}

\begin{prop}
\label{corXtimeCY2XtimesYform}
We keep notation \ref{dfn-CfS2}.
\begin{enumerate}[(a)]
\item The canonical  morphism 
\begin{equation}
\label{XtimeCY2XtimesYform}
\fP \times _{\scr{C} _{\fS}} \fQ
\to 
\fP \times _{\fS} \fQ
\end{equation}
is relatively perfect, i.e. the reductions modulo $\pi ^{i+1}$ are relatively perfect for any integer $i$.

\item 
\label{corXtimeCY2XtimesYform-item2}
Suppose 
$\fP/\fS$ and $\fQ/\fS$ have locally finite $p$-bases.
Then $\fP \times _{\scr{C} _{T}} \fQ/\fS$ has locally finite $p$-bases and 
the projections
\begin{gather}
\notag
pr _1 
\colon 
\fP \times _{\scr{C} _{\fS}} \fQ
\to \fP,
\\
\label{XtimeCY2XtimesY-corpr}
pr _2 
\colon 
\fP \times _{\scr{C} _{\fS}} \fQ
\to 
\fQ
\end{gather}
are flat and have locally finite $p$-bases. 

\end{enumerate}
\end{prop}

\begin{proof}
We get the relative perfectness of \ref{XtimeCY2XtimesYform} from that of
\ref{XtimeCY2XtimesY}.
Suppose 
$\fP/\fS$ and $\fQ/\fS$ have locally finite $p$-bases.
This implies that $\fP \times _{\scr{C} _{\fS}} \fQ \to \fS$
is formally smooth. Hence, following \ref{f0formétale-fforméta0-cor}, 
$\fP \times _{\scr{C} _{\fS}} \fQ \to \fS$ is flat
(because $\fP \times _{\scr{C} _{\fS}} \fQ$ is of finite type over $\bbD ^{r+s} _{\fS} $). 
Finally, by using \ref{rem-dfn-pbasispadicbis} and \ref{corXtimeCY2XtimesY}.\ref{corXtimeCY2XtimesY-item2},
this yields that $\fQ \times _{\scr{C} _{T}} \fQ/\fS$ have locally finite $p$-bases.
\end{proof}

\begin{ntn}
From now, we keep notation \ref{dfn-CfS2} and 
we suppose moreover that 
$\fP/\fS$ and $\fQ/\fS$ have locally finite $p$-bases.

\end{ntn}

\begin{empt}
Using the tensor product defined in \ref{def-otimes-coh1qc}, 
we get the bifunctor
\begin{equation}
\label{dfnboxtimes}
 \smash{\widehat{\boxtimes}}
^\L _{\O _{\fS }}
\colon 
\smash{\underrightarrow{LD}}  ^\mathrm{b} _{\Q, \mathrm{qc}}
( \smash{\widehat{\D}} _{\PP ^{ }/\fS }  ^{(\bullet)} )
\times 
\smash{\underrightarrow{LD}}  ^\mathrm{b} _{\Q, \mathrm{qc}}
( \smash{\widehat{\D}} _{\fQ ^{ }/\fS }  ^{(\bullet)} )
\to 
\smash{\underrightarrow{LD}}  ^\mathrm{b} _{\Q, \mathrm{qc}}
( \smash{\widehat{\D}} _{\PP \times _{\scr{C} _{\fS}} \fQ/\fS }  ^{(\bullet)} )
\end{equation}
defined as follows: 
for any $\E ^{ (\bullet)} 
\in 
\smash{\underrightarrow{LD}}  ^\mathrm{b} _{\Q, \mathrm{qc}}
( \smash{\widehat{\D}} _{\PP ^{ }/\fS }  ^{(\bullet)} )$,
$\FF  ^{ (\bullet)} 
\in 
\smash{\underrightarrow{LD}}  ^\mathrm{b} _{\Q, \mathrm{qc}}
( \smash{\widehat{\D}} _{\fQ ^{ }/\fS }  ^{(\bullet)} )$, 
we set
$$ \E ^{ (\bullet)}
 \smash{\widehat{\boxtimes}}
^\L _{\O _{\fS }}
\FF ^{ (\bullet)}
:= 
pr _1 ^{(\bullet) *} \E ^{ (\bullet)}
 \smash{\widehat{\otimes}}
^\L _{\O ^{(\bullet)}  _{\PP \times _{\scr{C} _{\fS}} \fQ} }
pr _2 ^{(\bullet) *} \FF ^{ (\bullet)}.$$
As for \cite[4.3.5]{Beintro2}, this functor induces 
the following one
\begin{equation}
\label{boxtimesLDcoh}
 \smash{\widehat{\boxtimes}} 
^\L _{\O _{\fS }}
\colon 
\smash{\underrightarrow{LD}}  ^\mathrm{b} _{\Q, \mathrm{coh}}
( \smash{\widehat{\D}} _{\PP ^{ }/\fS }  ^{(\bullet)} )
\times 
\smash{\underrightarrow{LD}}  ^\mathrm{b} _{\Q, \mathrm{coh}}
( \smash{\widehat{\D}} _{\fQ ^{ }/\fS }  ^{(\bullet)} )
\to 
\smash{\underrightarrow{LD}}  ^\mathrm{b} _{\Q, \mathrm{coh}}
( \smash{\widehat{\D}} _{\PP \times _{\scr{C} _{\fS}} \fQ/\fS }  ^{(\bullet)} ).
\end{equation}

\end{empt}

\begin{empt}
For any $\E ^{ (\bullet)} 
\in 
\smash{\underrightarrow{LD}}  ^\mathrm{b} _{\Q, \mathrm{qc}}
( \smash{\widehat{\D}} _{\fP ^{ }/\fS }  ^{(\bullet)} )$,
$\FF  ^{ (\bullet)} 
\in 
\smash{\underrightarrow{LD}}  ^\mathrm{b} _{\Q, \mathrm{qc}}
( \smash{\widehat{\D}} _{\fQ ^{ }/\fS }  ^{(\bullet)} )$, 
we have the isomorphism
\begin{equation}
\label{boxtimesalg-formal}
\E ^{ (\bullet)}
 \smash{\widehat{\boxtimes}}
^\L _{\O _{\fS }}
\FF ^{ (\bullet)}
\riso 
\R \underleftarrow{\lim}_i \,  
\left (
\E _i ^{ (\bullet)}
 \smash{\widehat{\boxtimes}}
^\L _{\O _{S _i}}
\FF _i ^{ (\bullet)}
\right ),
\end{equation}
where as usual we set 
$\E ^{ (\bullet)} _i := \smash{\widehat{\D}} _{P  _i/ S  _i} ^{(\bullet)}  \otimes ^\L _{\smash{\widehat{\D}} _{\fP /\fS } ^{(\bullet)} } \E ^{ (\bullet)}$,
and
$\cF ^{ (\bullet)} _i := \smash{\widehat{\D}} _{P  _i/ S  _i} ^{(\bullet)}  \otimes ^\L _{\smash{\widehat{\D}} _{\fP /\fS } ^{(\bullet)} } \cF ^{ (\bullet)}$.

\end{empt}

\begin{lem}
\label{exact-boxtimes}
The bifunctor \ref{boxtimesLDcoh} induces the
 exact bifunctor
$$ \smash{\widehat{\boxtimes}} ^\L _{\O _{\fS }}
\colon 
\smash{\underrightarrow{LM}}  _{\Q, \mathrm{coh}}
( \smash{\widehat{\D}} _{\fP ^{ }/\fS }  ^{(\bullet)} )
\times 
\smash{\underrightarrow{LM}}   _{\Q, \mathrm{coh}}
( \smash{\widehat{\D}} _{\fQ ^{ }/\fS }  ^{(\bullet)} )
\to 
\smash{\underrightarrow{LM}}  _{\Q, \mathrm{coh}}
( \smash{\widehat{\D}} _{\fP \times _{\scr{C} _{\fS}} \fQ/\fS }  ^{(\bullet)} ).$$ 
\end{lem}

\begin{proof}
Let 
$\E ^{ (\bullet)} 
\in 
\smash{\underrightarrow{LM}}  _{\Q, \mathrm{coh}}
( \smash{\widehat{\D}} _{\fP ^{ }/\fS }  ^{(\bullet)} )$,
$\FF  ^{ (\bullet)} 
\in 
\smash{\underrightarrow{LM}}   _{\Q, \mathrm{coh}}
( \smash{\widehat{\D}} _{\fQ ^{ }/\fS }  ^{(\bullet)} )$.
Let 
$\E:= \underrightarrow{\lim} 
\,
\E ^{ (\bullet)}$,
$\FF:= \underrightarrow{\lim} 
\,
\FF ^{ (\bullet)}$,
where $\underrightarrow{\lim} $ is the equivalence of categories
of \ref{M-eq-coh-lim}.
Choose $m _0$ large enough so that there exists
a coherent 
$ \smash{\widehat{\D}} _{\fP ^{ }/\fS }  ^{(m _0)} $-module 
$\mathscr{E} ^{ (m _0)} $ without $p$-torsion such that 
$\D ^\dag _{\fP/\fS, \Q}
\otimes _{\smash{\widehat{\D}} _{\fP ^{ }/\fS }  ^{(m _0)} }
\mathscr{E} ^{ (m _0)}  
\riso 
\E $, 
and 
a coherent 
$ \smash{\widehat{\D}} _{\fQ ^{ }/\fS }  ^{(m _0)} $-module 
$\mathscr{F} ^{ (m _0)} $ without $p$-torsion such that 
$\D ^\dag _{\fQ/\fS, \Q}
\otimes _{\smash{\widehat{\D}} _{\fQ ^{ }/\fS }  ^{(m _0)} }
\mathscr{F} ^{ (m _0)}  
\riso 
\FF $. 
For any $m \geq m _0$, 
let 
$\mathscr{E} ^{ (m )} $ and 
(resp. $\mathscr{F} ^{ (m )} $) 
 be 
the quotient of 
$\smash{\widehat{\D}} _{\fP ^{ }/\fS }  ^{(m )} 
\otimes _{\smash{\widehat{\D}} _{\fP ^{ }/\fS }  ^{(m _0)} }
\mathscr{E} ^{ (m _0)} $
(resp. 
$\smash{\widehat{\D}} _{\fQ ^{ }/\fS }  ^{(m )} 
\otimes _{\smash{\widehat{\D}} _{\fQ ^{ }/\fS }  ^{(m _0)} }
\mathscr{F} ^{ (m _0)} $)
by its torsion part. 
 We get 
 $\mathscr{E} ^{ (\bullet + m _0)} 
\in 
\smash{\underrightarrow{LM}}   _{\Q, \mathrm{coh}}
( \smash{\widehat{\D}} _{\fP ^{ }/\fS }  ^{(\bullet)} )$,
$\mathscr{F}  ^{ (\bullet+ m _0)} 
\in 
\smash{\underrightarrow{LM}}  _{\Q, \mathrm{coh}}
( \smash{\widehat{\D}} _{\fQ ^{ }/\fS }  ^{(\bullet)} )$
such that 
$\underrightarrow{\lim} 
\,
\mathscr{E} ^{ (\bullet + m _0)} 
\riso 
\E$,
and 
$\underrightarrow{\lim} 
\,
\mathscr{F} ^{ (\bullet + m _0)} 
\riso 
\FF$.
Hence, we obtain the isomorphisms 
$\E ^{ (\bullet)}  \riso \mathscr{E} ^{ (\bullet + m _0)} $
and 
$\FF ^{ (\bullet)}  \riso \mathscr{F} ^{ (\bullet + m _0)} $.
Since 
$\mathscr{E} ^{ (m)}$ and 
$\mathscr{F} ^{ (m)}$ have no $p$-torsion, then 
$Tor _1 ^{\cO _{S _i}} ( \cO _{S _0}, \cE _i) = 0$ 
and
$Tor _1 ^{\cO _{S _i}} ( \cO _{S _0}, \cF _i) = 0$.
Hence, 
following \ref{corXtimeCY2XtimesY}.\ref{corXtimeCY2XtimesY-item3},
the canonical morphism
\begin{equation}
\label{exact-boxtimes-iso1}
pr _1 ^* (\cE _i ^{ (m)})
\otimes ^{\bbL} _{\cO _{P _i \times _{\fC _{S _i}} Q _i}}
pr _2 ^* (\cF _i^{ (m)})
\to
pr _1 ^* (\cE _i ^{ (m)})
\otimes  _{\cO _{P _i \times _{\fC _{S _i}} Q _i}}
pr _2 ^* (\cF _i ^{ (m)})
\end{equation}
is an isomorphism.
Hence, 
\begin{gather}
\notag
pr _1 ^{*} \mathscr{E} ^{ (m)}
 \smash{\widehat{\otimes}}
^\L _{\O   _{\fP \times \fQ} }
pr _2 ^{*} \mathscr{F} ^{ (m)}
\riso 
\R \underleftarrow{\lim}_i ~
pr _1 ^* (\cE _i ^{ (m)})
\otimes ^{\bbL} _{\cO _{P _i \times _{\fC _{S _i}} Q _i}}
pr _2 ^* (\cF _i ^{ (m)})
\\
\notag
\underset{\ref{exact-boxtimes-iso1}}{\riso} 
\R \underleftarrow{\lim}_i ~
pr _1 ^* (\cE _i ^{ (m)})
\otimes  _{\cO _{P _i \times _{\fC _{S _i}} Q _i}}
pr _2 ^* (\cF _i ^{ (m)})
\riso 
\underleftarrow{\lim}_i ~
pr _1 ^* (\cE _i ^{ (m)})
\otimes  _{\cO _{P _i \times _{\fC _{S _i}} Q _i}}
pr _2 ^* (\cF _i ^{ (m)})
\\
\notag 
\riso 
pr _1 ^{*} \mathscr{E} ^{ (m)}
 \smash{\widehat{\otimes}}
 _{\O   _{\fP \times \fQ} }
pr _2 ^{*} \mathscr{F} ^{ (m)},
\end{gather}
where the third isomorphism is checked using 
Mittag-Leffler.
\end{proof}

\begin{coro}
We get the t-exact bifunctor
\begin{equation}
\label{dfn-boxtimesDLMcoh}
 \smash{\widehat{\boxtimes}} 
^\L _{\O _{\fS }}
\colon 
D ^{\mathrm{b}}  (\underrightarrow{LM} _{\Q,\mathrm{coh}} (\smash{\widehat{\D}} _{\fP /\fS } ^{(\bullet)} ))
\times 
D ^{\mathrm{b}}  (\underrightarrow{LM} _{\Q,\mathrm{coh}} (\smash{\widehat{\D}} _{\fQ /\fS } ^{(\bullet)} ))
\to 
D ^{\mathrm{b}}  (\underrightarrow{LM} _{\Q,\mathrm{coh}} (\smash{\widehat{\D}} _{\cP \times _{\scr{C} _{\fS}} \fQ /\fS } ^{(\bullet)})).
\end{equation}

\end{coro}

\begin{prop}
\label{prop-boxtimes}
\begin{enumerate}[(a)]
\item 
\label{prop-boxtimes1}
Let 
$\E ^{ (\bullet)} 
\in 
D ^{\mathrm{b}}  (\underrightarrow{LM} _{\Q,\mathrm{coh}} (\smash{\widehat{\D}} _{\fP /\fS } ^{(\bullet)} ))$,
$\FF ^{ (\bullet)} 
\in 
D ^{\mathrm{b}}  (\underrightarrow{LM} _{\Q,\mathrm{coh}} (\smash{\widehat{\D}} _{\fQ /\fS } ^{(\bullet)} ))$. 
We get
the spectral sequence in 
$\underrightarrow{LM} _{\Q,\mathrm{coh}} (\smash{\widehat{\D}} _{\cP \times _{\scr{C} _{\fS}}\fQ /\fS } ^{(\bullet)})$
of the form
$$\H ^r (\E  ^{ (\bullet)} )
\smash{\widehat{\boxtimes}} 
^\L _{\O _{\fS }}
\H ^s
 ( 
\FF  ^{ (\bullet)}
)
=:
E _{2} ^{r,s}
\Rightarrow
E ^n :=
\H ^n 
\left ( 
\E  ^{ (\bullet)} 
\smash{\widehat{\boxtimes}} 
^\L _{\O _{\fS }}
\FF  ^{ (\bullet)}
\right ) .$$
In particular,
when
$\E ^{ (\bullet)} 
\in 
\smash{\underrightarrow{LM}}  _{\Q, \mathrm{coh}}
( \smash{\widehat{\D}} _{\fP ^{ }/\fS }  ^{(\bullet)} )$,
this yields 
$\H ^n 
\left ( 
\E  ^{ (\bullet)} 
\smash{\widehat{\boxtimes}} 
^\L _{\O _{\fS }}
\FF  ^{ (\bullet)}
\right ) 
\riso 
\E  ^{ (\bullet)} 
\smash{\widehat{\boxtimes}} 
^\L _{\O _{\fS }}
\H ^n 
 ( 
\FF  ^{ (\bullet)}
) $.

\item 
\label{prop-boxtimes2}
Suppose $\fQ$ affine.
Let 
$\E ^{ (\bullet)} 
\in 
\smash{\underrightarrow{LM}}  _{\Q, \mathrm{coh}}
( \smash{\widehat{\D}} _{\fP ^{ }/\fS }  ^{(\bullet)} )$,
$\FF  ^{ (\bullet)} 
\in 
\smash{\underrightarrow{LD}}  ^\mathrm{b} _{\Q, \mathrm{coh}}
( \smash{\widehat{\D}} _{\fQ ^{ }/\fS }  ^{(\bullet)} )$.
We have 
$\H ^n 
\left ( 
\E  ^{ (\bullet)} 
\smash{\widehat{\boxtimes}} 
^\L _{\O _{\fS }}
\FF  ^{ (\bullet)}
\right ) 
\riso 
\E  ^{ (\bullet)} 
\smash{\widehat{\boxtimes}} 
^\L _{\O _{\fS }}
\H ^n 
 ( 
\FF  ^{ (\bullet)}
) $.

\end{enumerate}

\end{prop}

\begin{proof}
We can copy the proof of \cite[7.5.5]{caro-6operations}.
\end{proof}

\begin{prop}
\label{prop-boxtimes-v+}
Let $f \colon \fP  \to \bbD ^{r} _{\fS}$ 
and 
$f'  \colon \fP ' \to \bbD ^{r} _{\fS}$ 
be two finite type morphisms,
let $(u, id) \colon f '  \to f$ be a morphism of 
$\fC _S$ 
(see notation \ref{dfn-CfS}). 
Let $g \colon \fQ  \to \bbD ^{s} _{\fS}$ 
and 
$g ' \colon \fQ ' \to \bbD ^{s} _{\fS}$ 
be two finite type morphisms,
let $(v, id) \colon g'  \to g$ be a morphism of 
$\fC _S$.
We suppose moreover that 
$\fP/\fS$,
$\fP'/\fS$,
$\fQ/\fS$
 and $\fQ '/\fS$ have locally finite $p$-bases.
Let
$\ZZ: =\fP \times _{\scr{C} _{\fS}} \fQ $,
$\ZZ ' := \fP '\times _{\scr{C} _{\fS}} \fQ'$,
and $w:= (u, v) \colon \ZZ ' \to \ZZ$ be the induced morphism.
\begin{enumerate}[(a)]
\item For any
$\E ^{(\bullet)}
\in \underrightarrow{LD}  ^\mathrm{b} _{\Q, \mathrm{qc}}
(\overset{^\mathrm{g}}{} \smash{\widehat{\D}} _{\fP} ^{(\bullet)} )$
and 
$\FF ^{ (\bullet)}
\in \underrightarrow{LD}  ^\mathrm{b} _{\Q, \mathrm{qc}}
(\overset{^\mathrm{g}}{} \smash{\widehat{\D}} _{\fQ} ^{(\bullet)} )$, 
with notation \ref{ntn-Lf!+*}, we have in 
 $\underrightarrow{LD}  ^\mathrm{b} _{\Q, \mathrm{qc}}
(\overset{^\mathrm{g}}{} \smash{\widehat{\D}} _{\ZZ} ^{(\bullet)})$
the isomorphism:
\begin{equation}
\label{boxtimes-v!formal}
\L w ^{*(\bullet)}  (\E ^{(\bullet)}
\smash{\widehat{\boxtimes}} ^\L _{\O _{\fS }}
\FF ^{(\bullet)})
\riso
\L u ^{*(\bullet)} (\E ^{(\bullet)})
\smash{\widehat{\boxtimes}} ^\L _{\O _{\fS }}
\L v ^{*(\bullet)} (\FF ^{(\bullet)}).
\end{equation}

\item For any
$\E ^{\prime(\bullet)}
\in \underrightarrow{LD}  ^\mathrm{b} _{\Q, \mathrm{qc}}
(\overset{^\mathrm{g}}{} \smash{\widehat{\D}} _{\fP'} ^{(\bullet)} )$
and 
$\FF ^{\prime (\bullet)}
\in \underrightarrow{LD}  ^\mathrm{b} _{\Q, \mathrm{qc}}
(\overset{^\mathrm{g}}{} \smash{\widehat{\D}} _{\fQ'} ^{(\bullet)} )$, 
we have in 
 $\underrightarrow{LD}  ^\mathrm{b} _{\Q, \mathrm{qc}}
(\overset{^\mathrm{g}}{} \smash{\widehat{\D}} _{\ZZ} ^{(\bullet)})$
the isomorphism:
\begin{equation}
\label{boxtimes-v+}
w ^{(\bullet)} _+ (\E ^{\prime(\bullet)}
\smash{\widehat{\boxtimes}} ^\L _{\O _{\fS }}
\FF ^{\prime(\bullet)})
\riso
u ^{(\bullet)}_+ (\E ^{\prime(\bullet)})
\smash{\widehat{\boxtimes}} ^\L _{\O _{\fS }}
v ^{(\bullet)}_+ (\FF ^{\prime(\bullet)}).
\end{equation}

\end{enumerate}

\end{prop}

\begin{proof}
The first statement is a consequence of \ref{comm-boxtimes-f*} and \ref{boxtimesalg-formal}.
The second one is a consequence of \ref{sch-prop-boxtimes-v+} and \ref{boxtimesalg-formal}.\end{proof}

\begin{coro}
\label{theo-iso-chgtbase2}
We keep notation \ref{prop-boxtimes-v+} and we suppose
$v $ is the identity. 
Let $\pi \colon  \fZ  \to \fP $, 
and $\pi '\colon  \fZ ^{ \prime } \to \fP ^{ \prime }$
be the projections. 
Let 
$\E ^{\prime (\bullet)}
\in  \smash{\underrightarrow{LD}} ^\mathrm{b} _{\Q, \mathrm{qc}}
(\overset{^\mathrm{l}}{} \smash{\widehat{\D}} _{\fP ^{\prime }} ^{(\bullet)})$. 
There exists a canonical isomorphism in 
$\smash{\underrightarrow{LD}} ^\mathrm{b} _{\Q, \mathrm{qc}}
(\overset{^\mathrm{l}}{} \smash{\widehat{\D}} _{\fZ } ^{(\bullet)})$ of the form:
\begin{equation}
\label{iso-chgtbase2}
\pi ^{ !(\bullet)} \circ u ^{(\bullet)}_{ +} (\E ^{\prime (\bullet)})
\riso
w  ^{(\bullet)}_{+}  \circ \pi  ^{\prime (\bullet) !} (\E ^{\prime (\bullet)}). 
\end{equation}
\end{coro}

\begin{proof}
This is a consequence of \ref{theo-iso-chgtbase2-pre} (or we can deduce it from \ref{prop-boxtimes-v+}).\end{proof}

\begin{rem}
\label{remoftheo-iso-chgtbase2}
We will prove later (see \ref{theo-iso-chgtbase}) a coherent version of Corollary \ref{theo-iso-chgtbase2}.
In this version, we can use for instance Berthelot-Kashiwara theorem which allow us to extend geometrically 
the context. 
\end{rem}

\begin{dfn}
\label{projectivefscheme}
Let $f\colon \fX \to \fY$ be a morphism 
of formal $\fS$-schemes of formal finite type and having locally finite $p$-bases.
We say that $f$ is projective (resp. quasi-projective) if $f$ is the composition of a closed immersion (resp. immersion) of the form
$\fX \hookrightarrow \widehat{\bbP} ^d \times _\fS \fY$ with the projection 
$\widehat{\bbP} ^d \times _\fS \fY \to \fY$. 
Beware that this notion is stronger than that appearing in \cite{EGAII}.
\end{dfn}

\begin{prop}
\label{rel-dual-isom-proj-formal}
Let $f\colon \fX \to \fY$ be a projective morphism 
of formal $\fS$-schemes of formal finite type and having locally finite $p$-bases.

\begin{enumerate}[(a)]
\item For any
$\cE ^{(\bullet)} \in  \smash{\underrightarrow{LD}} ^\mathrm{b} _{\Q, \mathrm{coh}}
(\overset{^\mathrm{l}}{} \smash{\widehat{\D}} _{\fX } ^{(\bullet)})$,
we have a canonical isomorphism of 
$\smash{\underrightarrow{LD}} ^{\mathrm{b}} _{\Q,\mathrm{coh}} ( \smash{\widehat{\D}} _{\fY} ^{(\bullet)})$
of the form 
\begin{equation}
\label{rel-dual-isom-proj-formal-iso}
\DD ^{(\bullet)} \circ f ^{(\bullet)} _+ (\cE ^{(\bullet)})
\riso 
 f ^{(\bullet)} _+ \circ \DD ^{(\bullet)} (\cE ^{(\bullet)}).
\end{equation}

\item Let $\E  \in D ^\mathrm{b} _{\mathrm{coh}}
(\D ^{\dag} _{\fX,\Q})$,
and
$\cF 
\in 
D ^\mathrm{b} _{\mathrm{coh}}
(\D ^{\dag} _{\fY ,\Q})$.
We have 
the isomorphisms
\begin{gather}
\label{cor-adj-formul-proj-formal-bij1}
\R \mathcal{H} om _{\D ^{\dag} _{\fY,\Q}}
( f _{+} ( \E ) , \cF) 
\riso 
\R f _* 
\R \mathcal{H} om _{\D ^{\dag} _{\fX ,\Q}}
( \E  ,f ^!   ( \cF)),
\\
\label{cor-adj-formul-proj-formal-bij2}
\R \mathrm{Hom}  _{\D ^{\dag} _{\fY ,\Q}}
( f _{+} ( \E) , \cF) 
\riso 
\R \mathrm{Hom}  _{\D ^{\dag} _{\fX,\Q}}
( \E  ,f ^!   ( \cF)). 
\end{gather}
\end{enumerate}

\end{prop}

\begin{proof}
The first statement is a consequence of \ref{rel-dual-isom-proj}.
Similarly to  \cite[4.4.2]{caro-6operations}, we check that \ref{rel-dual-isom-proj-formal-iso}
implies the second statement.
\end{proof}

\subsection{Relative duality isomorphism and adjunction for relatively proper complexes and quasi-projective morphisms}

\begin{dfn}
[Proper support with respect to a morphism]
\label{dfn-prop-support}
Let $g \colon \X '\to \X$ be a morphism 
of formal $\fS$-schemes of formal finite type and having locally finite $p$-bases.
Let 
$\E ^{\prime (\bullet)} \in \smash{\underrightarrow{LD}} ^{\mathrm{b}} _{\Q,\mathrm{coh}} ( \smash{\widehat{\D}} _{\X '/\V} ^{(\bullet)})$.
We say that $\E ^{\prime (\bullet)}$
has a proper support over $X$
if there exist a closed subscheme $Z'$ of $X'$ such that 
$\E ^{\prime (\bullet)}$ has his support in $Z$' (i.e. 
$\E ^{\prime (\bullet)} | \U ' = 0$ with $\fU ': =\fX ' \setminus Z'$) 
and such that the composite morphism 
$Z ' \hookrightarrow  X' \overset{g}{\to} X$ is proper. 
\end{dfn}

\begin{prop}
\label{stab-propersupp}
Let $g \colon \X '\to \X$ be a quasi-projective (in the sense of \ref{projectivefscheme}) morphism of formal $\fS$-schemes of formal finite type and having locally finite $p$-bases.
For any 
$\E ^{\prime (\bullet)} \in \smash{\underrightarrow{LD}} ^{\mathrm{b}} _{\Q,\mathrm{coh}} ( \smash{\widehat{\D}} _{\X '/\V} ^{(\bullet)})$
with proper support over $X$ (see \ref{dfn-prop-support}),
 the object 
$g _{+} (\E ^{\prime (\bullet)} ) $
belongs to 
$\smash{\underrightarrow{LD}} ^{\mathrm{b}} _{\Q,\mathrm{coh}} ( \smash{\widehat{\D}} _{\X/\V} ^{(\bullet)})$.
\end{prop}

\begin{proof}
We can copy the proof of \cite[10.3.2]{caro-6operations}.
\end{proof}

\begin{thm}
[Relative duality isomorphism]
\label{rel-dual-isom}
Let $g \colon \fP '\to \fP$ be a quasi-projective morphism of formal $\fS$-schemes of formal finite type and having locally finite $p$-bases.
For any 
$\E ^{\prime (\bullet)} 
\in 
\smash{\underrightarrow{LD}} ^{\mathrm{b}} _{\Q,\mathrm{coh}} ( \smash{\widehat{\D}} _{\fP '} ^{(\bullet)})$
with proper support over $P$, 
we have the isomorphism of $\smash{\underrightarrow{LD}} ^{\mathrm{b}} _{\Q,\mathrm{coh}} ( \smash{\widehat{\D}} _{\fP} ^{(\bullet)})$
of the form 
$$g _{+} \circ \DD 
(\E ^{\prime (\bullet)} ) 
\riso 
\DD \circ g _{+} 
(\E ^{\prime (\bullet)} ) .
$$
\end{thm}

\begin{proof}
By copying the proof of \cite[10.4.1]{caro-6operations},
we check that this is a consequence of \ref{rel-dual-isom-proj-formal}.
\end{proof}

\begin{cor}
\label{cor-adj-formulbis}
Let $g \colon \fP '\to \fP$ be a quasi-projective morphism of formal $\fS$-schemes of formal finite type and having locally finite $p$-bases.
Let $\E ' \in D ^\mathrm{b} _{\mathrm{coh}}
(\D ^{\dag} _{\PP ^{\prime },\Q})$
with proper support over $P$,
and
$\E 
\in 
D ^\mathrm{b} _{\mathrm{coh}}
(\D ^{\dag} _{\PP  ,\Q})$.
We have 
the isomorphisms
\begin{gather}
\label{cor-adj-formulbis-bij1}
\R \mathcal{H} om _{\D ^{\dag} _{\PP ,\Q}}
( g _{+} ( \E ') , \E) 
\riso 
\R g _* 
\R \mathcal{H} om _{\D ^{\dag} _{\PP ' ,\Q}}
( \E ' ,g ^! ( \E)). 
\\
\label{cor-adj-formulbis-bij2}
\R \mathrm{Hom}  _{\D ^{\dag} _{\PP ,\Q}}
( g _{+} ( \E ') , \E) 
\riso 
\R \mathrm{Hom}  _{\D ^{\dag} _{\PP ' ,\Q}}
( \E ' ,g ^!   ( \E)). 
\end{gather}
\end{cor}

\begin{proof}
By copying the proof of \cite[4.4.2]{caro-6operations},
we check that this is a consequence of 
\ref{rel-dual-isom}.
\end{proof}

\section{On the differential coherence of $\O _{\X} ({} ^\dag Z) _\Q$}

\subsection{Descent of coherence via finite base change}

\begin{lem}
\label{lem-desc-coh-chgbase}
Let $\V \to \V'$ be a finite morphism of complete discrete valuation rings of mixed characteristics $(0,p)$.
We get the finite morphism $\fS ' := \Spf \V ' \to \fS$.
Let $\X $ be a formal  $\fS$-scheme
of formal finite type
and having locally finite $p$-bases  over $\fS$.
Let
$\X ' := \X \times _{\fS} \fS '$, 
and $f \colon \X' \to \X$ be the canonical projection.
Let $Z$ be a divisor of $X$ and $Z':= f ^{-1} (Z)$.
\begin{enumerate}[(a)]
\item  
\label{lem-desc-coh-chgbase0}
The canonical 
homomorphism 
$\D ^\dag _{\X'/\fS '} (\hdag Z') _{\Q} 
\to 
\D ^\dag _{\X'\to \X/\fS '\to \fS} (\hdag Z') _{\Q} $
is an isomorphism. 
The composite morphism
$f ^{-1}\D ^\dag _{\X/\fS} (\hdag Z) _{\Q} 
\to 
\D ^\dag _{\X'\to \X/\fS '\to \fS} (\hdag Z') _{\Q} 
\liso
\D ^\dag _{\X'/\fS '} (\hdag Z') _{\Q} $
is a homomorphism of rings. 
Hence, if $\E$ is a coherent $\D ^\dag _{\X/\fS} (\hdag Z) _{\Q}$-module, then
$f _Z ^! (\E) \riso \D ^\dag _{\X'/\fS '} (\hdag Z') _{\Q} \otimes _{f ^{-1}\D ^\dag _{\X/\fS} (\hdag Z) _{\Q}}
f ^{-1}\E$, where $f ^! _Z$ is the extraordinary inverse image of $\X' \to \X$ above $\fS' \to \fS$
 with overconvergent singularities along $Z$, i.e. $f ^! _Z$ is the base change inverse image.

\item 
\label{lem-desc-coh-chgbase1}
Suppose $\X$ is affine.
Let $\E$ be a coherent $\D ^\dag _{\X/\fS} (\hdag Z) _{\Q}$-module.
 Then the canonical morphisms
$$ \V ' \otimes _\V \Gamma (\X, \E) 
\to 
D ^\dag _{\X'/\fS '} (\hdag Z') _{\Q} 
\otimes _{D ^\dag _{\X/\fS} (\hdag Z) _{\Q}}
\Gamma (\X, \E) 
\to 
\Gamma (\X', f _Z ^! (\E)) $$
 are isomorphisms.
 Moreover, 
$D ^\dag _{\X'/\fS'} (\hdag Z') _{\Q}$ is a faithfully flat 
$D ^\dag _{\X/\fS} (\hdag Z) _{\Q}$-module for both left or right structure.

\item 
\label{lem-desc-coh-chgbase2}
For any  $\D ^\dag _{\X/\fS} (\hdag Z) _{\Q}$-module $\E$,
the canonical morphisms
\begin{equation}
\notag
f ^* (\E) := \O _{\X'} \otimes _{f ^{-1}\O _{\X}}
f ^{-1}\E 
\to
\O _{\X'} (\hdag Z') _{\Q} \otimes _{f ^{-1}\O _{\X} (\hdag Z) _{\Q}}
f ^{-1}\E 
\to
\D ^\dag _{\X'/\fS '} (\hdag Z') _{\Q} \otimes _{f ^{-1}\D ^\dag _{\X/\fS} (\hdag Z) _{\Q}}
f ^{-1}\E 
\end{equation}
are isomorphisms.

 \item 
 \label{lem-desc-coh-chgbase3}
 Let $\phi \colon \E' \to \E$ be a morphism of $\O _{\X}$-modules.
 Then $\phi$ is an isomorphism if and only if 
$f ^* (\phi)$ 
is an isomorphism.
\end{enumerate}

\end{lem}

\begin{proof}
We can copy the proof of \cite[8.3.1]{caro-6operations}.
\end{proof}

\begin{prop}
\label{desc-coh-chgbase}
With notation \ref{lem-desc-coh-chgbase}, let
 $\E$ be a   $\D ^\dag _{\X/\fS} (\hdag Z) _{\Q}$-coherent module.
Then $\E$ is a coherent $\D ^\dag _{\X/\fS, \Q}$-module if and only if 
$f _Z ^! (\E) $ is a coherent $\D ^\dag _{\X'/\fS ', \Q}$-module.  
\end{prop}

\begin{proof}
We can copy the proof of \cite[8.3.2]{caro-6operations}.
\end{proof}

For completeness, we add  Proposition \ref{desc-coh-chgbase2}, 
which is useless in this paper but  which 
extends somehow Lemma \ref{desc-coh-chgbase}.

\begin{rem}
\label{gonflement-cohen}
Let $k \to l$ be an extension of perfect field of characteristic $p>0$.
Since $k \to l$ is separable, 
following \cite[$0.19.8.2.(ii)$]{EGAIV1}, 
there exists a unique up to (non unique) isomorphism $\V$-algebra of Cohen  $\W$ (in the sense of 
\cite[$0.19.8.1$]{EGAIV1})
which is a lifting of $k \to l$. 
\end{rem}

\begin{prop}
\label{desc-coh-chgbase2}
With notation \ref{gonflement-cohen}, 
suppose $l$ is algebraic over $k$.
Let $\T := \Spf \W \to \fS$ be the corresponding morphism of formal $p$-adic schemes.
Let $\X $ be a 
 formal  $\fS$-schemes of formal finite type
and having locally finite $p$-bases  over $\fS$,
$\Y := \X \times _{\fS} \T$, 
and $f \colon \Y \to \X$ be the canonical projection.
Let $Z _X$ be a divisor of $X$ and 
$Z _Y := f ^{-1} (Z _X)$ be the corresponding divisor of $Y$.

The homomorphisms 
$\widehat{\D} ^{(m)} _{\X/\fS} ( Z _X) 
\to 
f _* \widehat{\D} ^{(m)} _{\Y/\T} ( Z _Y) $
and 
$\D ^\dag _{\X/\fS} (\hdag Z _X) _{\Q}
\to 
f _* \D ^\dag _{\Y/\T} (\hdag Z _Y ) _{\Q}$
are right and left faithfully flat (in the sense of the definition after \cite[Lemma 4.3.8]{Be1}). 
\end{prop}

\begin{proof}
We can copy the proof of \cite[8.4.7]{caro-6operations}.
\end{proof}

\subsection{Extraordinary pullbacks by a projective morphism: comparison between $\O$-modules and
$\cD$-modules}
We prove in this subsection the isomorphism
\ref{fund-isom2-corpre-bis-Gen-dag}
which will be useful in \ref{256Be2-iso1}.
\begin{empt}
Let $u \colon \ZZ \hookrightarrow \X$ be a closed immersion of 
formal  $\fS$-schemes of formal finite type and having locally finite $p$-bases  over $\fS$.
Following \ref{fflat=f!-immer-f}, when we get the isomorphism
\begin{equation}
\label{fflat=f!-immerbis}
u ^{\flat} \riso u ^!
\end{equation}
of functors
$D  ({} ^r \D ^{\dag} _{\X/\fS,\Q} )
\to 
D  ({} ^r \D ^{\dag} _{\fZ/\fS,\Q} )$
(resp. 
$D  ({} ^r\D ^{\dag} _{\X/\fS,\Q} ,{} ^* \D ^{\dag} _{\fX/\fS,\Q} )
\to 
D  ({} ^r \D ^{\dag} _{\fZ/\fS,\Q} ,{} ^* u ^{-1}\D ^{\dag} _{\fX/\fS,\Q} )$).

\end{empt}

\begin{empt}
\label{ntn-PX2X}
Let $\fX$ be a formal  $\fS$-scheme
of formal finite type
and having locally finite $p$-bases  over $\fS$.
Let $\fY := \widehat{\bbP} ^{d} _{\fS}\times _{\fS} \fX$
and $f \colon \fY \to \fX$ be the canonical projection.
We get $f _i \colon Y _i \to X _i$. 
Following \cite[III.2]{HaRD},
since $f$ is smooth, then we have the functor
$f _i ^{\sharp}
\colon 
D ( \O _{X _i})
\to 
D ( \O _{Y _i})$ defined by setting
for any $\cM _i  \in D ( \O _{X _i})$, 
$$f _i ^{\sharp} (\cM _i)
:=
f _i ^{*} (\cM _i ) \otimes _{\O _{Y _i}} \omega _{Y _i /Y _i}[d]
\riso
f _i ^{*} (\cM _i \otimes _{\O_{X _i}}\omega ^{-1}_{X _i /S _i}) \otimes _{\O _{Y _i}} \omega _{Y _i /S _i}[d],$$
where the isomorphism comes from the fact $f _i $ is smooth 
(recall also we have defined it in a wider context in \ref{ntn-psharp-dfng^sharp}).

Let 
$\cM _i$ be a right $\cD ^{(m)} _{X _i/S _i}$-module.
Using $m$-PD-costratification, 
since the functors of the form $p _{ij} ^{\flat}$ are exact, 
since $f _i ^{\sharp} [-n]$ is acyclic, then
by using \cite[III.8.7]{HaRD} we get 
a canonical $m$-PD-costratification on 
$f _i ^{\sharp} [-n] (\cM _i)$.
Hence, the functor 
$f _i ^{\sharp}$
induces the functor
$f _i ^{\sharp}
\colon 
D ( \cD ^{(m)} _{X _i/S _i})
\to 
D ( \cD ^{(m)} _{Y _i/S _i})$.

\end{empt}

\begin{prop}
\label{thm !=sharp}
We keep notation \ref{ntn-PX2X}. 
Let $* \in \{\mathrm{l}, \mathrm{r}\}$.
\begin{enumerate}[(a)]
\item We have an isomorphism 
\begin{equation}
\label{thm !=sharpiso2}
f ^{!} \riso f ^{\sharp} 
\end{equation}
of functors 
$D  ({} ^r \D _{X _i/S _i} ^{(m)}) \to 
D  ({} ^r \D _{Y _i/S _i} ^{(m)})$
(resp. $D  ( {} ^r \D _{X _i/S _i} ^{(m)}, {} ^* \D _{X _i/S _i} ^{(m)})
\to 
D  ({} ^r \D _{Y _i/S _i} ^{(m)},{} ^* f ^{-1}\D _{X _i/S _i} ^{(m)})$).

\item Let $\E \in D  ( {} ^l \D _{X _i/S _i} ^{(m)})$ (resp. $\E \in D  ( {} ^l \D _{X _i/S _i} ^{(m)}, {} ^* \D _{X _i/S _i} ^{(m)})$). 
We have the canonical isomorphism of 
$D  ({} ^r \D _{Y _i/S _i} ^{(m)})$
(resp. $D  ({} ^r \D _{Y _i/S _i} ^{(m)},{} ^* f ^{-1}\D _{X _i/S _i} ^{(m)})$)
\begin{equation}
\label{thm !=sharpiso}
\omega _{Y _i/S _i} \otimes _{\O _{X _i}} f ^{!} (\E )
\riso
f ^{\sharp} ( \omega _{X _i/S _i} \otimes _{\O _{X _i}} \E).
\end{equation}

\item We have the canonical isomorphism of $( u ^{-1}\D _{X} ^{(m)}, \D _{Z} ^{(m)}) $-bimodules of the form
\begin{equation}
\label{fund-isom2-corpre-bis-smooth}
 \D _{X _i\leftarrow Y _i} ^{(m)} [d]
\riso
f _i ^{\sharp}  (  \D _{X _i}  ^{(m)}) .
\end{equation}

\end{enumerate}
\end{prop}

\begin{proof}
1) Let us check the first statement.
The canonical isomorphism \ref{thm !=sharpiso2}
is already known (see \ref{ntn-psharp-dfng^sharp-iso}).
To check the $\D _{Y _i/S _i} ^{(m)}$-linearity, 
we reduce to the case where $X _i$ has a finite $p$-basis. 
Then this is an easy computation.

2) Since we have also the isomorphism
$\omega _{Y _i/S _i} \otimes _{\O _{X _i}} f ^{!} (\E )
\riso
f ^{!} ( \omega _{X _i/S _i} \otimes _{\O _{X _i}} \E)$,
then 
\ref{thm !=sharpiso} is a straightforward consequence of \ref{thm !=sharpiso2}.

3) The third statement is a consequence of \ref{thm !=sharpiso}.
\end{proof}

\begin{empt}
\label{proj-gen-flat*def}
Let $g \colon \fZ \to \fX$ be a projective morphism of 
formal  $\fS$-schemes of formal finite type and having locally finite $p$-bases  over $\fS$ in the following strong sense : 
there exists a closed immersion 
$u\colon \fZ \hookrightarrow \fY := \widehat{\bbP} ^{d} _{\fS}\times _{\fS} \fX$
such that 
$g = f \circ u$ where 
$f \colon \fY \to \fX$ be the canonical projection.

\begin{enumerate}[(a)]
\item We set 
$g _i ^{\flat \sharp} := u _i ^{\flat} f _i ^{\sharp}
\colon 
D ({} ^r  \cD ^{(m)} _{X _i/S _i})
\to 
D ( {} ^r \cD ^{(m)} _{Z _i/S _i})$.
We have another functor
$g _i ^{!} 
\colon 
D ( \cD ^{(m)} _{X _i/S _i})
\to 
D ( \cD ^{(m)} _{Z _i/S _i})$.
Following \ref{fflat=f!-immerbis} and \ref{thm !=sharpiso2}, we have in fact the isomorphism
\begin{equation}
\label{proj-gen-flat*def-iso1}
g _i ^{\flat \sharp}
\riso 
g _i ^{!}
\end{equation}
of functors
$D  ({} ^r \D _{X _i/S _i} ^{(m)}) \to 
D  ({} ^r \D _{Z _i/S _i} ^{(m)})$.

\item When $g$ is a finite morphism, 
we have $g _i ^{\flat \sharp} \liso g _i ^{\flat}$ 
as functors 
of 
$D ( \cO _{X _i})
\to 
D ( \cO _{Z _i})$
(see \cite[III.8.7]{HaRD}).
In fact, by construct of both functors, 
this isomorphism is horizontal, i.e. commutes with the $m$-PD-costratification. 
Hence, 
we get the first isomorphism 
\begin{equation}
\label{flat=flatsharp=!}
g _i ^{\flat}
\riso 
g _i ^{\flat \sharp}
\underset{\ref{proj-gen-flat*def-iso1}}{\riso} 
g _i ^{!} 
\end{equation}
of functors
$D  ({} ^r \D _{X _i/S _i} ^{(m)}) \to 
D  ({} ^r \D _{Z _i/S _i} ^{(m)})$.
This yields the isomorphism of right $\D _{Z _i/S _i} ^{(m)}$-modules
\begin{equation}
\label{flat=flatsharp=!2}
g _i ^{\flat} ( \omega _{X _i/S _i} ) 
\riso 
g _i ^{!} ( \omega _{X _i/S _i} ) 
\riso 
\omega _{Y _i/S _i}.
\end{equation}
Using 
\ref{fund-isom2-corpre-bis}
and 
\ref{fund-isom2-corpre-bis-smooth}, 
we construct the canonical isomorphism of $( g _i ^{-1}\D _{X _i} ^{(m)}, \D _{Z _i} ^{(m)}) $-bimodules of the form
\begin{equation}
\label{flat=flatsharp=!3}
 \D _{X _i\leftarrow Z _i} ^{(m)}
\riso
g _i ^{\flat}  (  \D _{X _i}  ^{(m)}) .
\end{equation}

\item We still suppose $g$ is a finite morphism. 
Taking projective and inductive limits, the isomorphism \ref{flat=flatsharp=!} induces the isomorphism 
$g  ^{\flat}
\riso 
 g  ^{!} $
of functors
$D  ({} ^r \D ^{\dag} _{\fX/\fS,\Q}) \to 
D  ({} ^r \D ^{\dag} _{\fZ/\fS,\Q})$.
Again, taking projective and inductive limits, we get from \ref{flat=flatsharp=!2} and \ref{flat=flatsharp=!3} the isomorphisms
\begin{gather}
\label{fund-isom2-corpre-bis-Gen-dagpre}
g ^{\flat} ( \omega _{\fX /\fS } ) 
\riso 
g  ^{!} ( \omega _{\fX/\fS} ) 
\riso 
\omega _{\fY /\fS },
\\
\label{fund-isom2-corpre-bis-Gen-dag}
 \D ^\dag _{\fX \leftarrow \fZ} 
\riso
g  ^{\flat}  (  \D ^\dag _{\fX}) .
\end{gather}

\end{enumerate}

\end{empt}

\begin{rem}
\label{rem-proj-gen-flat*def}
With notation \ref{proj-gen-flat*def},
the induced functor 
$g _i ^{\flat \sharp} 
\colon 
D ( \cO _{Z _i})
\to 
D ( \cO _{Y _i})$
is denoted by 
$g _i ^{!}$ in \cite[III.8.7]{HaRD}.
Since 
$g _i ^{!} 
\riso 
g _i ^{\flat \sharp} $ then these notations are compatible. 
But we have written $g _i ^{\flat \sharp}$ in order to avoid confusion with 
the functor $g _i ^{!} 
\colon 
D ( \cD ^{(m)} _{Z _i/S _i})
\to 
D ( \cD ^{(m)} _{Y _i/S _i})$.
\end{rem}

\subsection{Descent of coherence via universal homeomorphisms}

\begin{empt}
[Universal homeomorphism]

\label{univ-homeo}
Let $f\colon X \to Y$ be a morphism of schemes. 

\begin{enumerate}[(a)]
\item 
\label{univ-homeo-1}
Following Definitions \cite[3.5.4]{EGAI} (and Remark \cite[3.5.11]{EGAI})
or \cite[2.4.2]{EGAIV2}, 
$f$ is by definition a universal homeomorphism (resp. is universally injective) if
for any morphism of schemes $g \colon Y ' \to Y$, the morphism 
$f _{Y'}\colon X \times _{Y} Y' \to Y'$ is a homeomorphism (resp. is injective). 

\item 
\label{univ-homeo-2}
Some authors use the name of ``purely inseparable'' (e.g. \cite[5.3.13]{Liu-livre-02}) or ``radicial'' (e.g. \cite[3.5.4]{EGAI}) 
instead of ``universally injective''. From Definition \cite[3.5.4]{EGAI}, Proposition \cite[3.5.8]{EGAI} and Remark 
\cite[3.5.11]{EGAI}, the following conditions are equivalent : 

\begin{enumerate}[(a)]
\item $f$ is universally injective ; 
\item for any field $K$, the map $X (K) \to Y (K)$ is injective ; 
\item $f$ is injective and for any point $x$ of $X$ the monomorphism of the 
residue fields $k (f(x)) \to k (x)$ induced by $f$ is purely inseparable (some authors say ``radicial'' instead of ``purely inseparable''). 
	
\end{enumerate}

\item 
\label{univ-homeo-3}
Suppose $f\colon X \to Y$ is a morphism of finite type such that $Y$ is locally noetherian.
Following Proposition \cite[2.4.5]{EGAIV2}, $f$ is a universal homeomorphism if and only if 
$f$ is finite, surjectif and universally injective. 

\end{enumerate}

\end{empt}

\begin{lem}
\label{Ex5.3.9Liu}
Let $f\colon X \to Y$ is a finite, surjective morphism such that $Y$ is a normal noetherian scheme and $X$ is integral.
The morphism $f$ is  a universal homeomorphism
if and only if $k(X)/k(Y)$ is radicial.
\end{lem}

\begin{proof}
Suppose $f$ is  a universal homeomorphism. 
Since the generic point of $X$ is sent to the generic point of 
$Y$ then 
the hypothesis that $f$ is universally injective 
implies that $k(X)/k(Y)$ is radicial. 

Conversely, suppose $k(X)/k(Y)$ is radicial.
We have to check that $f$ is universally injective.
We can suppose $X = \Spec A$ and $Y = \Spec B$.
Since $f$ is surjective then 
$B \to A$ is injective. 
Since  $\Frac (A)/\Frac (B)$ is radicial, then for any $a \in A$, 
there exists $s$ large enough such that
$a ^{p ^s} \in \Frac (B)$.
Since 
$B$ is normal and $A/B$ is finite, this yields
$a ^{p ^s} \in B$.
This implies that 
$f $ is injective and that 
for any point $x$ of $X$ the monomorphism of the 
residue fields $k (f(x)) \to k (x)$ induced by $f$ is radicial. 
Following \ref{univ-homeo}.\ref{univ-homeo-2}, this yields that 
$f$ is universally injective.
By using \ref{univ-homeo}.\ref{univ-homeo-3}, 
this implies that $f$ is a universal homeomorphism.\end{proof}
\begin{lem}
\label{F-locfreeft2}
Let $f\colon \Y  \to \X  $ be 
a morphism of 
formal  $\fS$-schemes
of formal finite type
and having locally finite $p$-bases  over $\fS$.
We suppose  that the induced morphism 
$f _0 \colon Y \to X$ is a finite, surjective and radicial morphism. 
Then the morphism $f$ is locally free of finite type,
i.e. $f _* \O _\Y$ is a locally free of finite type
$\O _\X$-module. 
\end{lem}

\begin{proof}
Since this is local, we can suppose $\fX$ and $\fY$ affine. 
Since $X$ and $Y$ are regular (see \ref{regularity/formalsm}), then 
following \cite[4.3.11]{Liu-livre-02}
(or see \cite[15.4.2]{EGAIV4}, \cite[5.4.2]{EGAIV4}
and
\cite[0.17.3.5]{EGAIV1})
$f _0$ is flat.
Since $\X$ and $\fY$ are noetherian,
since $X$ and $Y$ are regular then $\X$ and $\fY$ are regular
(use \cite[Lemma 6.1]{MonskyWashnitzer}), i.e. 
the formal spectrum of a regular ring.
Using \cite[11.3.10]{EGAIV3}, this yields that $f $ is flat. 
\end{proof}

\begin{empt}
\label{tau-dagm}
Let $\X $, $\Y $ be two 
formal  $\fS$-schemes
of formal finite type
and having locally finite $p$-bases  over $\fS$.
Let $f,f'\colon \Y \to \X$ be two finite  morphisms of formal $\fS$-schemes
such that $f _{0} =f' _{0}$.
Using \ref{215Be2},
we get the isomorphism of respectively $(\widehat{\D} ^{(m)} _{\Y /\fS}, f _0 ^{-1} \widehat{\D} ^{(m)} _{\X /\fS})$-bimodules
and
$(f _0 ^{-1} \widehat{\D} ^{(m)} _{\X /\fS},\widehat{\D} ^{(m)} _{\Y /\fS})$-bimodules
$$\tau _{f,f'}
\colon 
f ^{\prime *}  \widehat{\D} ^{(m)} _{\X /\fS} 
\riso 
f ^*  \widehat{\D} ^{(m)} _{\X /\fS} 
,
\
\sigma _{f,f'}
\colon 
f ^{\prime \flat}  \widehat{\D} ^{(m)} _{\X /\fS} 
\riso 
f ^{\flat}  \widehat{\D} ^{(m)} _{\X /\fS} .$$
Looking at the construction of the  isomorphism $\tau _{f,f'}$, 
we get the following explicit local description of the isomorphism: 
suppose $\X$ has  the finite $p$-basis
$t  _1, \dots, t  _d \in \Gamma (\X, \O _{\X })$.
Then the image of 
$1 \otimes 1$ 
is 
$$
\tau _{f,f'} (  1 \otimes 1)
=
\sum _{\underline{i} \in \N ^d}
(f ^{\prime *} (\underline{t} )- f ^{*} (\underline{t}) )^{\{\underline{i} \} _{(m)}}
\otimes \underline{\partial} ^{<\underline{i}> _{(m)}}.$$

Taking inductive limits on the level,
this yields
 the isomorphism of respectively 
 $(\D ^\dag _{\Y /\fS}, f _0 ^{-1} \D ^\dag _{\X /\fS})$-bimodules
and 
$(f _0 ^{-1} \D ^\dag _{\X /\fS}, \D ^\dag _{\Y /\fS})$-bimodules
$$\tau _{f,f'}
\colon 
f ^{\prime *}  \D ^\dag _{\X /\fS} 
\riso 
f ^*  \D ^\dag _{\X /\fS} ,
\ 
\sigma _{f,f'}
\colon 
f ^{\prime \flat}  \D ^\dag _{\X /\fS} 
\riso 
f ^{\flat}  \D ^\dag _{\X /\fS} .$$
Hence, for any $\cE \in D ^{\mathrm{b}} _{\mathrm{coh}} ({} ^{\mathrm{l}} \D ^\dag _{\X /\fS}) $
and
$\cM \in D ^{\mathrm{b}} _{\mathrm{coh}} ({} ^{\mathrm{r}}\D ^\dag _{\X /\fS}) $, we get the isomorphism
\begin{gather}
\label{glueff'-dag-l}
\tau _{f,f'} 
\colon 
f ^{\prime *}  \cE
=
f ^{\prime *}  \D ^\dag _{\X /\fS}  \otimes _{f ^{-1} \D ^\dag _{\X /\fS} } f ^{-1}\cE
\underset{\tau _{f,f'}}{\riso} 
f ^* \D ^\dag _{\X /\fS}  \otimes _{f ^{-1} \D ^\dag _{\X /\fS} } f ^{-1} \cE
=
f ^* \cE;
\\
\label{glueff'-dag-r}
\sigma _{f,f'} 
\colon 
f ^{\prime \flat}  \cM
=
 f ^{-1}\cM
 \otimes _{f ^{-1} \D ^\dag _{\X /\fS} }
 f ^{\prime \flat}  \D ^\dag _{\X /\fS}
\underset{\sigma _{f,f'}}{\riso} 
f ^{-1} \cM
\otimes _{f ^{-1} \D ^\dag _{\X /\fS} } 
f ^\flat \D ^\dag _{\X /\fS}  
=
f ^\flat \cM
\end{gather}
Moreover, 
for any $\cE \in D ^{\mathrm{b}} _{\mathrm{coh}} ({} ^{\mathrm{l}} \D ^\dag _{\X /\fS}) $
and
for any $\cM \in D ^{\mathrm{b}} _{\mathrm{coh}} ({} ^{r} \D ^\dag _{\Y /\fS}) $, we get the isomorphisms
\begin{gather}
\label{glueff'-dag-r+}
\tau _{f,f'} 
\colon 
f ^{\prime } _+ \cM
=
\cM \otimes _{\D ^\dag _{\Y /\fS} }  f ^{\prime *}  \D ^\dag _{\X /\fS} 
\underset{\tau _{f,f'}}{\riso} 
\cM \otimes _{\D ^\dag _{\Y /\fS} }  f ^{ *}  \D ^\dag _{\X /\fS} 
=
f _+ \cM;
\\
\label{glueff'-dag-l+}
f ^{\prime } _+ \cE
:=
\left (\omega _{\Y /\fS} 
\otimes _{\O _{\fY}}
f ^{\prime *} ( \D ^\dag _{\X /\fS} \otimes _{\O _{\fX}} \omega ^{-1} _{\X /\fS}) \right )
\otimes _{\D ^\dag _{\Y /\fS} }  
\cE 
\underset{\tau _{f,f'}}{\riso} 
\left (\omega _{\Y /\fS} 
\otimes _{\O _{\fY}}
f ^{ *} ( \D ^\dag _{\X /\fS} \otimes _{\O _{\fX}} \omega ^{-1} _{\X /\fS}) \right )
\otimes _{\D ^\dag _{\Y /\fS} }  
\cE 
=
f _+ \cE.
\end{gather}

\end{empt}

\begin{lem}
\label{incl-Xs0pre}
Let $T _0$ be a noetherian $S _0$-scheme. 
Let $X _0$ be a integral, noetherian $T _0$-scheme
having locally finite $p$-bases  over $T _0$.
Then for any integer $s$, the morphism
$F ^s _{X _0/T _0}\colon X _0 \to X _0 ^{(s)}$ is finite, radicial, surjective. 
\end{lem}

\begin{proof}
We can suppose we have a relatively perfect morphism of the form $g _0\colon X _0 \to \A ^d _{T _0}$. 
Then we get the cartesian square
\begin{equation}
\label{relperf-Frob-ter}
\xymatrix@ R= 0,4cm{
{X _0} 
\ar[r] ^-{F ^s _{X _0/T _0}}
\ar[d] ^-{g _0}
\ar@{}[dr] ^-{} |\square
& 
{X _0 ^{(s)}} 
\ar[d] ^-{g _0 ^{(s)}}
\\ 
{ \A ^d _{T _0}} 
\ar[r] _-{F ^s _{ \A ^d _{T _0}/T _0}}
& 
{\A ^d _{T _0} } 
}
\end{equation}
Since $F ^s _{ \A ^d _{T _0}/T _0}$
is a finite, radicial surjective morphism,
we can conclude.
\end{proof}

\begin{lem}
\label{incl-Xs0}
Let $T _0:= \Spec A _0$ be an $S _0$-scheme such that the absolute Frobenius
$F _{T _0} \colon T _0 \to  T _0$ is finite. 
Let $Y _0:= \Spec A _0 [[t _1, \dots, t _r]] /I$ be an affine $T _0$-scheme of formal finite type.
For any integer $s$, let
$Y _0  ^{(s)}$ be the base change of $Y _0$ by the $s$-th power of the absolute Frobenius of $T _0$.

\begin{enumerate}[(a)]

\item Then 
$Y _0  ^{(s)} = \Spec A _0 [[u _1, \dots, u _r]] /I ^{(p^s)}$,
where 
$I ^{(p^s)}$ is the ideal of $A _0 [[u _1, \dots, u _r]] $ generated by the elements 
of the form 
$\sum _{\nu \in \bbN ^{r}} a _{\nu} ^{p ^s} u ^{\nu}$,
with 
$\sum _{\nu \in \bbN ^{r}} a _{\nu} u ^{\nu} \in I$.

\item The relative Frobenius $F ^s _{Y _0/T _0}\colon Y _0 \to Y _0 ^{(s)}$
is induced by the $A _0$-algebra homomorphism
$u _i \mapsto t _i ^{p^s}$.

\item If $Y _0$ and $Y _0 ^{(s)}$ are integral and $A _0 =k$, then 
$k ( Y _0^{(s)}) 
\subset
k \cdot k ( Y_0)  ^{p ^s}$,
where
$k ( Y _0^{(s)}) 
:=
\Frac \left (k [[u _1, \dots, u _r]] /I ^{(p^s)}\right) $
and 
$k ( Y _0) 
:=
\Frac (k [[t _1, \dots, t _r]] /I )$.

\end{enumerate}

\end{lem}

\begin{proof}
Set 
$\bbD := \Spec A _0 [[t _1, \dots, t _r]] $. 
Since 
$F ^s _{T _0} \colon T _0 \to  T _0$ is finite, then 
$\bbD ^{(s)} 
= 
\bbD \times _{T _0, F ^r _{T _0}} T _0
=
\Spec A _0 [[u _1, \dots, u _r]] $.
The relative Frobenius morphism 
$\bbD 
\to 
\bbD^{(s)} $
corresponds to 
the $A _0$-algebra homomorphism
$A _0 [[u _1, \dots, u _r]] 
\to 
A _0 [[t _1, \dots, t _r]] $
sending 
$u _i $ to $t ^{p ^s} _i$.
We conclude via the commutative diagram
\begin{equation}
\notag
\xymatrix @C=2cm{
{\bbD } 
\ar[r] ^-{F ^s _{\bbD /T _0}}
&
{\bbD ^{(s)}}
\ar@{}[rd] ^-{} |\square
\ar[r] ^-{}
&
{\bbD } 
\\ 
{Y _0} 
\ar@{^{(}->}[u] ^-{}
\ar[r] ^-{F ^s _{Y _0/T _0}}
& 
{Y _0 ^{(s)}} 
\ar@{^{(}->}[u] ^-{}
\ar[r] ^-{}
& 
{Y _0,} 
\ar@{^{(}->}[u] ^-{}
}
\end{equation}
and noticing that 
$I A _0 [[u _1, \dots, u _r]] 
= 
I ^{(p ^s)}$. 

Since the absolute Frobenius $k \to k$ is finite, 
then the first two statements are satisfied in the case where $A _0 =k$. 
Moreover, following 
\ref{incl-Xs0pre},
$F ^s _{Y _0/T _0}$ is finite, radicial and surjective.
Hence, 
the induced $k$-homomorphism
$k [[u _1, \dots, u _r]] /I ^{(p^s)}
\to 
k [[t _1, \dots, t _r]] /I $
is injective. 
Let $P \in k [[u _1, \dots, u _r]] $.
If $\lambda _1,\dots, \lambda _N$ is a basis of 
$k/k ^{p^s}$, then 
we can write 
$P
= 
\sum _{i =1} ^N
\lambda _i P _i$, 
with 
$P _i \in k ^{p^s} [[u _1, \dots, u _r]]$.
Let us denote by 
$\phi \colon k [[u _1, \dots, u _r]]
\to k [[t _1, \dots, t _r]]$ 
the $k$-homomorphism given by 
$u _i \mapsto t _i ^{p ^s}$. 
Then $\phi (P _i) \in \left ( k [[t _1, \dots, t _r]] \right ) ^{p ^s}$.
This implies that $\phi (P) $ belongs to the $k$-subalgebra of $ k [[t _1, \dots, t _r]]$ generated by 
$\left ( k [[t _1, \dots, t _r]] \right ) ^{p ^s}$.
Hence, 
the image of 
$k [[u _1, \dots, u _r]] /I ^{(p^s)}
\to 
k [[t _1, \dots, t _r]] /I $
is included in the 
$k$-subalgebra of $ k [[t _1, \dots, t _r]]/I$ generated by 
$\left ( k [[t _1, \dots, t _r]] /I \right ) ^{p ^s}$, 
and in particular 
in $k \cdot k ( Y_0)  ^{p ^s}$
Hence, we are done.
 \end{proof}

\begin{thm}
\label{univhomeo-eqcat}
Let $f\colon \Y  \to \X  $ be 
a morphism of 
formal  $\fS$-schemes
of formal finite type
and having locally finite $p$-bases  over $\fS$.
We suppose  that the induced morphism 
$f _0 \colon Y \to X$ is a finite, surjective and radicial morphism. 
The functor
$\E \mapsto f ^* \E$ (resp. $\cM \mapsto f ^\flat \cM$)
from the category of left (resp. right) $\D ^{\dag} _{\X /\fS }$-modules
to that of left  (resp. right) 
$\D ^{\dag} _{\Y /\fS }$-modules is an exact equivalence of categories.

\end{thm}

\begin{proof}

a) The fact that $f ^*$ (resp. $f ^\flat$) is exact is a consequence of \ref{F-locfreeft2}.

b) Since this is local, we can suppose $X$ affine
(resp. $Y$ affine of the form $Y:= \Spec k [[t _1, \dots, t _r]] /I$
with $I$ an ideal of $k [[t _1, \dots, t _r]] $).

Following \ref{regularity/formalsm}, $X$, $Y$, $X ^{(s)}$ and 
$Y ^{(s)}$ are regular.
Since $\X$ is affine and noetherian and 
since $X$ is regular then $\X$ is regular
(see \cite[Lemma 6.1]{MonskyWashnitzer}).
Since $X$ and $Y$ are normal and $f _0$ is a universal homeomorphism, then 
$k (X) \subset k(Y)$ is radicial
(see \ref{Ex5.3.9Liu}).
Hence, 
for $s$ large enough, 
we have the inclusion 
$k (Y) ^{p ^s} \subset k (X)$
and then 
$k \cdot k (Y) ^{p ^s} \subset k (X)$. 
Moreover, using 
\ref{incl-Xs0}, we get
$k ( Y ^{(s)}) 
\subset
k \cdot k ( Y)  ^{p ^s}$.
Hence, 
$k ( Y ^{(s)}) 
\subset
 k (X)$.
 From
$\Gamma (Y ^{(s)}, \O _{Y ^{(s)}})
\subset
k ( Y ^{(s)})
\cap 
\Gamma (Y, \O _Y)$,
this yields
$\Gamma (Y ^{(s)}, \O _{Y ^{(s)}})
\subset
k (X)
\cap 
\Gamma (Y, \O _Y)
$.
Since $X$ is normal and $f$ is finite, we get
$\Gamma (Y, \O _Y) \cap k (X)
=
\Gamma (X, \O _X)$.
Hence, we have checked 
$\Gamma (Y ^{(s)}, \O _{Y ^{(s)}})
\subset
\Gamma (X, \O _X)$.
In other words, 
there exists a morphism 
$g _0  \colon X \to Y ^{(s)}$ making commutative the diagram of $S$-schemes
\begin{equation}
\label{univhomeo-eqcat-diag}
\xymatrix{
{Y}
\ar[d] _-{F ^s _{Y/S}}
\ar[r] ^-{f _0 } 
& 
{X} 
\ar[ld] _-{g _0 }
\ar[d] ^-{F ^s _{X/S}}
\\ 
{Y ^{(s)}}
\ar[r] _-{f _0 ^{ (s)}} 
& 
{X ^{(s)}.} 
}
\end{equation}

Following corollary \ref{lifting-pbasis}, 
there exists a 
formal $\fS$-scheme of formal finite type 
$\X ^{\prime}$
(resp. $\Y ^{\prime}$) having finite $p$-basis
and such that 
$\X ^{\prime} \times _{\fS} S  \riso X ^{(s)} $
(resp. $\Y ^{\prime} \times _{\fS} S  \riso Y ^{(s)} $).
Since $\Y ^{\prime}$ is formally smooth and $\X$ is affine, 
there exists a lifting
$g \colon \X  \to \Y ^{\prime}$ of $g _0 $.
Similarly, there exists 
a lifting 
$f ' \colon \X  ^{\prime} \to \Y ^{\prime}$ of $f _0 ^{(s)}$.
We get the lifting
$F _\fX := f' \circ g \colon \fX  \to \X ^{\prime}$
(resp. $F _\fY:=g \circ f \colon  \Y  \to \Y ^{\prime}$)
of $F ^s _{X/S}$
(resp. $F ^s _{Y/S}$).

c) Following \ref{thm236Be2-cordag} (resp. \ref{thm236Be2-cordag-cohrightpre}),
$F ^{*} _{\X}$ 
(resp. $F ^{\flat} _{\X}$)
induces an equivalence between the category of 
left $\D ^\dag _{\fX ^{\prime},\Q} $-modules 
and that of 
left $\D ^\dag _{\fX,\Q} $-modules,
and similarly for 
$F ^{*} _{\Y}$ 
(resp. $F ^{\flat} _{\Y}$)
Hence, 
using the transitivity with respect to the composition of morphisms of the functors
$\E \mapsto f ^* \E$
(resp. $\E \mapsto f ^{\flat} \E$)
we conclude.
\end{proof}

\subsection{Quasi-inverse functor for universal homeomorphisms}

\begin{empt}
\label{empt-flat**flat}
Let 
$m,s\geq 0$
be two integers, $T$ be an ${S_i}$-scheme of finite type endowed with a quasi-coherent $m$-PD-ideal
$(\mathfrak{a},\mathfrak{b},\alpha)$ such that $p \in \mathfrak{a}$.

Let $f\colon Y  \to X  $ be 
a morphism of $T$-schemes
of formal finite type
and having locally finite $p$-bases  over $T$.
Following \ref{rem-Dnn'Ynoeth}.\ref{rem-Dnn'Ynoeth-i)} and \ref{f0formétale-fforméta0}, $X/T$ and $Y/T$ are flat,
$X$ and $Y$ are noetherian. 
We suppose $f$ is a homeomorphism and is a finite and locally free morphism.
Let $\cM$ be a $\cD ^{(m)} _{X/T}$-bimodule.
Since
$f$ is a finite, locally free homeomorphism, then 
$f ^\flat (\cM)
=
f ^{-1} \cH om _{\O _X} ( f _{*} \O _{Y}, \cM)$. 
By functoriality, 
$f ^\flat (\cM)$ is 
$(f ^{-1}\cD ^{(m)} _{X/T}, \cD ^{(m)} _{Y/T})$-bimodule.
This yields a structure of 
$(\cD ^{(m)} _{X/T}, f _* \cD ^{(m)} _{Y/T})$-bimodule
on  $f _* f ^\flat (\cM)$.
By functoriality, 
we get a structure of 
$(\cD ^{(m)} _{Y/T},  \cD ^{(m)} _{Y/T})$-bimodule
on  
$f ^{*\flat}(\cM )
:= 
f ^* f _* f ^\flat (\cM)$.
Similarly we have a functorial
structure of 
$(\cD ^{(m)} _{Y/T},  \cD ^{(m)} _{Y/T})$-bimodule
on  
$f ^{\flat *}(\cM )
:= 
f ^\flat f _* f ^*  (\cM)$.
We have the canonical isomorphism
\begin{gather}
\notag
f ^{*\flat}
(\cM )
=
\O _{Y} \otimes _{f ^{-1} \O _{X}}
f ^{-1} \cH om _{\O _X} ( f _{*} \O _{Y}, \cM)
\riso
\O _{Y} \otimes _{f ^{-1} \O _{X}}
 \cH om _{f ^{-1} \O _X} (  \O _{Y}, f ^{-1} \cM)
 \\
\notag
\riso
\O _{Y} 
\otimes _{f ^{-1} \O _{X}}
f ^{-1} \cM
\otimes _{f ^{-1} \O _{X}}
 \cH om _{f ^{-1} \O _X} (  \O _{Y}, f ^{-1} \O _{X})
 \\
\label{flat**flat}
\riso
 \cH om _{f ^{-1} \O _X} (  \O _{Y}, \O _{Y} 
\otimes _{f ^{-1} \O _{X}}
f ^{-1} \cM
)
\riso
f ^{-1} \cH om _{ \O _X} (  f_* \O _{Y}, \
f _* f ^*  (\cM))
=
f ^{\flat *}(\cM ).
\end{gather}
By functoriality, we can check that the 
isomorphism
$f ^{*\flat}(\cM )
\riso f ^{\flat *}(\cM )$
of \ref{flat**flat}
is an isomorphism of 
$(\cD ^{(m)} _{Y/T},  \cD ^{(m)} _{Y/T})$-bimodules.

Similarly to the proof of \cite[2.5.2]{Be2},
we construct a morphism of 
$(\cD ^{(m)} _{Y/T},  \cD ^{(m)} _{Y/T})$-bimodules
of the form
\begin{equation}
\label{252Be2-morp-uh}
\rho _f 
\colon 
\cD ^{(m)} _{Y/T}
\to
f ^{*\flat}(\cD ^{(m)} _{X/T}).
\end{equation}

\end{empt}

\begin{empt}
\label{empt-flat**flat-trans}
We keep notation \ref{empt-flat**flat}.
The homomorphism \ref{252Be2-morp-uh} is transitive, i.e. we have the following properties.
Let $g\colon Z  \to Y  $ be 
a second morphism of 
noetherian flat  $T$-schemes
of formal finite type
and having locally finite $p$-bases  over $T$.
We suppose $g$ is a homeomorphism and is a finite and locally free morphism.
By transitivity of the functor $f ^*$ and $f ^\flat$ 
we get the canonical isomorphism
$g ^{*\flat} f ^{*\flat}(\cD ^{(m)} _{X/T})
\riso
(f \circ g)  ^{*\flat}(\cD ^{(m)} _{X/T})$.
By construction, we have the commutative diagram:
\begin{equation}
\label{empt-flat**flat-trans-diag}
\xymatrix{
{g ^{*\flat} \cD ^{(m)} _{Y/T}} 
\ar[r] ^-{g ^{*\flat}(\rho _f) }
& 
{g ^{*\flat} f ^{*\flat}(\cD ^{(m)} _{X/T})} 
\ar[d] ^-{\sim}
\\ 
{\cD ^{(m)} _{Z/T}} 
\ar[u] ^-{\rho _g }
\ar[r] ^-{\rho _{f\circ g} }
& 
{(f \circ g)  ^{*\flat}(\cD ^{(m)} _{X/T}).} 
}
\end{equation}

\end{empt}

\begin{empt}
With notation and hypotheses of \ref{subsecElevFrob}, 
we have the commutative diagram 
\begin{equation}
\label{comp-252uh-iso}
\xymatrix{
{\D ^{(m)} _{X/T}} 
\ar[r] ^-{}
\ar[d] _-{\ref{252Be2-morp-uh}}  ^-{\rho _F}
& 
{\D ^{(m+s)} _{X/T}} 
\ar[d] ^-{\ref{252Be2-iso}} _-{\sim}
\\ 
{F ^{*\flat}\D ^{(m)} _{X'/T}} 
\ar@{=}[r] ^-{}
& 
{F ^* F ^{\flat} \D ^{(m)} _{X'/T}.} 
}
\end{equation}
\end{empt}

\begin{empt}
\label{empt-flat**flat-form}
Let $f\colon \Y  \to \X  $ be 
a morphism of 
formal  $\fS$-schemes
of formal finite type
and having locally finite $p$-bases  over $\fS$.
We suppose $f$ is a homeomorphism which is a finite and locally free morphism.

\begin{enumerate}[(a)]
\item Let $\cM$ be a $\D ^{\dag} _{\X /\fS }$-bimodule.
It follows from \ref{empt-flat**flat} that we get a structure of 
$\D ^{\dag} _{\Y /\fS }$-bimodule on 
$f ^{*\flat} \cM
:=
\O _{\fY} \otimes _{f ^{-1} \O _{\X}}
f ^{-1} \cH om _{\O _{\X}} ( f _{*} \O _{\Y}, \cM)
$.

Taking projective limits of some morphisms of the form
\ref{252Be2-morp-uh} and next taking inductive limits on the level, we get the morphism
of $\D ^{\dag} _{\Y /\fS }$-bimodules:
\begin{equation}
\label{252Be2-morp-uh-form}
\rho _f \colon 
\D ^{\dag} _{\fY /\fS}
\to 
f ^{*\flat} \D ^{\dag} _{\X /\fS}.
\end{equation}

\item From \ref{empt-flat**flat-trans}, we can check that 
the homomorphism \ref{252Be2-morp-uh-form} is transitive, i.e. we have the following property.
Let $g\colon \Y  \to \X  $ be 
a morphism of 
formal  $\fS$-schemes
of formal finite type
and having locally finite $p$-bases  over $\fS$.
We suppose $g$ is a homeomorphism and is a finite and locally free morphism.
Then we have 
the canonical isomorphism
$g ^{*\flat} f ^{*\flat}(\D ^{\dag} _{\fX /\fS})
\riso
(f \circ g)  ^{*\flat}(\D ^{\dag} _{\fX /\fS})$
and the commutative diagram:
\begin{equation}
\label{empt-flat**flat-trans-diag-form}
\xymatrix{
{g ^{*\flat} \D ^{\dag} _{\fY /\fS}} 
\ar[r] ^-{g ^{*\flat} (\rho _f )}
& 
{g ^{*\flat} f ^{*\flat}(\D ^{\dag} _{\fX /\fS})} 
\ar[d] ^-{\sim}
\\ 
{\D ^{\dag} _{\fZ /\fS}} 
\ar[u] ^-{\rho _g }
\ar[r] ^-{\rho _{f\circ g} }
& 
{(f \circ g)  ^{*\flat}(\D ^{\dag} _{\fX /\fS}).} 
}
\end{equation}

\end{enumerate}

\end{empt}

\begin{empt}
\label{comp-252uh-iso-emptform}
Let $\X$ be a formal $\fS$-scheme
of formal finite type
and having locally finite $p$-bases  over $\fS$.
Let $X _0$ be its special fiber and 
$X _0  ^{(s)}$ be the base change of $X _0$ by the $s$-th power of the absolute Frobenius of $S _0$.
Suppose there exists 
$F \colon \X \to \X ^{\prime}$ a morphism of 
formal $\fS$-schemes
of formal finite type
and having locally finite $p$-bases  over $\fS$
which is a lifting of the relative Frobenius
$F ^s _{X _0/S _0}\colon X _0 \to X _0 ^{(s)}$.
By taking projective limits and next inductive limits on the level of some diagrams of the form
\ref{comp-252uh-iso}, 
we get that 
$\rho _F \colon 
\D ^{\dag} _{\fX /\fS}
\to 
F ^{*\flat} \D ^{\dag} _{\fX' /\fS}
=
F ^* F ^{\flat} 
\D ^{\dag} _{\fX' /\fS}$ is the canonical isomorphism.

\end{empt}

\begin{prop}
\label{univhomeo-252Be2}
Let $f\colon \Y  \to \X  $ be 
a morphism of 
formal  $\fS$-schemes
of formal finite type
and having locally finite $p$-bases  over $\fS$.
We suppose  that the induced morphism 
$f _0 \colon Y \to X$ is a finite, surjective and radicial morphism. 
Then the morphism of \ref{252Be2-morp-uh-form}
\begin{equation}
\label{thm-252Be2-morp-uh-form}
\rho _f \colon 
\D ^{\dag} _{\fY /\fS}
\to 
f ^{*\flat} \D ^{\dag} _{\X /\fS}.
\end{equation}
is an isomorphism.
\end{prop}

\begin{proof}
Since this is local, we can suppose $X$ affine
(resp. $Y$ affine of the form $Y:= \Spec k [[t _1, \dots, t _r]] /I$
with $I$ an ideal of $k [[t _1, \dots, t _r]] $).
Then, we can use the constructions and notation of the part b) of the proof of \ref{univhomeo-eqcat}. 
Following \ref{empt-flat**flat-trans-diag-form}; we have the following commutative diagram
\begin{equation}
\label{empt-flat**flat-trans-diag-form-appl1}
\xymatrix{
{f ^{*\flat} \D ^{\dag} _{\fX /\fS}} 
\ar[r] ^-{\rho _g }
& 
{f ^{*\flat} g ^{*\flat}(\D ^{\dag} _{\fY ^{\prime}/\fS})} 
\ar[d] ^-{\sim}
\\ 
{\D ^{\dag} _{\fY /\fS}} 
\ar[u] ^-{\rho _f }
\ar[r] ^-{\rho _{f\circ g} }
& 
{(g \circ f)  ^{*\flat}(\D ^{\dag} _{\fY ^{\prime} /\fS}).} 
}
\end{equation}
Since $g _0 \circ f _0  = F ^s _{X/S}$, then following \ref{comp-252uh-iso-emptform}, 
$\rho _{f\circ g} $ is an isomorphism.
Hence, $\rho _{f}$ is injective.
Since $g _0$ is also  a finite, surjective and radicial morphism,
then  $\rho _{g}$ is also injective.
Hence, $\rho _{f}$ is bijective.
\end{proof}

\begin{coro}
\label{253Be2-coro}
We keep notation \ref{univhomeo-252Be2}.

\begin{enumerate}[(a)]
\item The $\D ^{\dag} _{\Y /\fS }$-module 
$f ^{*} \D ^{\dag} _{\X /\fS}$
(resp. 
$f ^{\flat} \D ^{\dag} _{\X /\fS}$) is locally projective of finite type. 

\item Let $\E $ (resp. $\cM$) be a left (resp. right) $\D ^{\dag} _{\X /\fS }$-module.
Then $\E$ is $\D ^{\dag} _{\X /\fS }$-coherent if and only if 
$f ^* \E$ (resp. $f ^{\flat} \cM$) is $\D ^{\dag} _{\Y /\fS }$-coherent. 

\end{enumerate}

\end{coro}

\begin{proof}
By copying the proof of 
\cite[2.5.3]{Be2}, this is a consequence of 
Proposition \ref{univhomeo-252Be2}.
\end{proof}

\begin{coro}
\label{256Be2}
We keep notation \ref{univhomeo-252Be2}.
Let $\E $ (resp. $\cM$) be a left (resp. right) $\D ^{\dag} _{\X /\fS }$-module.
\begin{enumerate}[(a)]
 
\item There exist canonical $\D ^{\dag} _{\fX /\fS }$-linear isomorphisms  
\begin{equation}
\label{256Be2-iso1}
f _+ (f ^{*} \cE) 
\underset{\ref{fund-isom2-corpre-bis-Gen-dag}}{\riso}
\R f _{*} \left (
f ^{\flat} \D ^{\dag} _{\X /\fS}
\otimes ^{\bbL}_{\D ^{\dag} _{\Y /\fS }}
f ^{*} \cE
\right ) 
\riso
f _{*} \left (
f ^{\flat} \D ^{\dag} _{\X /\fS}
\otimes _{\D ^{\dag} _{\Y /\fS }}
f ^{*} \cE
\right ) 
\riso 
\cE. 
\end{equation}

\item There exist canonical $\D ^{\dag} _{\fX /\fS }$-linear isomorphisms  
\begin{equation}
\label{256Be2-iso2}
f _+ (f ^{\flat} \cM) 
\riso 
\R f _{*} \left (
f ^{\flat} \cM
\otimes ^{\bbL}_{\D ^{\dag} _{\Y /\fS }}
f ^{*} \D ^{\dag} _{\X /\fS}
\right ) 
\riso
\R f _{*} \left (
f ^{\flat} \cM
\otimes _{\D ^{\dag} _{\Y /\fS }}
f ^{*} \D ^{\dag} _{\X /\fS}
\right ) 
\riso 
\cM. 
\end{equation}
 \item 
 \label{256Be2-iso3} 
 Let $\cF $  be a left (resp. right) $\D ^{\dag} _{\fY /\fS }$-module.
Then $\cF$ is $\D ^{\dag} _{\fY/\fS }$-coherent if and only if 
$f _+ \cF$ is $\D ^{\dag} _{\fX /\fS }$-coherent. 

 \item 
 \label{256Be2-iso4}
The functor $f _+$
from the category of 
coherent left (resp. right) $\D ^{\dag} _{\Y /\fS ,\Q}$-modules 
to that of 
coherent left  (resp. right) $\D ^{\dag} _{\X /\fS ,\Q}$-modules
is an exact  quasi-inverse equivalence of categories of $f ^* $ (resp. $f ^\flat$).

\end{enumerate}
\end{coro}

\begin{proof}
By copying the proof of 
\cite[2.5.6]{Be2}, the first two assertions are a consequence of 
Proposition \ref{univhomeo-252Be2}.
Since $f$ is finite then $f$ is proper. 
Hence, if $\cF$ is $\D ^{\dag} _{\fY/\fS }$-coherent then
$f _+ \cF$ is $\D ^{\dag} _{\fX /\fS }$-coherent. 
Conversely,  following \ref{253Be2-coro},  if $f _+ \cF$ is $\D ^{\dag} _{\fX /\fS }$-coherent, 
then  $f ^* f _+ \cF$ (resp. $f ^\flat f _+ \cF$) is $\D ^{\dag} _{\fX /\fS }$-coherent.
This yields that 
$f _+ f ^* f _+ \cF$ (resp. $f _+ f ^\flat f _+ \cF$) is $\D ^{\dag} _{\fY /\fS }$-coherent.
Hence, we get the statement \ref{256Be2-iso3} by using respectively \ref{256Be2-iso1} and \ref{256Be2-iso2}.
Finally, 
using Theorem \ref{univhomeo-eqcat} and
Corollary \ref{253Be2-coro}, 
this yield the last statement.\end{proof}

\begin{lem}
\label{2.2.13-carocourbe}
Let $\X$ be an affine 
formal  $\fS$-scheme
of formal finite type
and having finite $p$-bases  over $\fS$.
Let $\E$ be a coherent left $\D ^\dag _{\X/\fS,\Q}$-module. 
The following conditions are equivalent
\begin{enumerate}[(a)]
\item The sheaf $\E$ is coherent $\O _{\X,\Q}$-module. 
\item $\Gamma ( \X, \E)$ is a $\Gamma (\X, \O _{\X,\Q})$-module of finite type.
\end{enumerate}
\end{lem}

\begin{proof}
We can copy word by word the proof of 
\cite[2.2.13]{caro_courbe-nouveau}.
\end{proof}

\begin{lem}
\label{cvisoc-desc-finite}
Let $f \colon \X \to \Y$ be a finite morphism of formal  $\fS$-schemes
of formal finite type
and having locally finite $p$-bases  over $\fS$. 
Let $\E$ be a left $\D ^\dag _{\X/\fS,\Q}$-module. 
The following conditions are equivalent
\begin{enumerate}[(a)]
\item $\E$ is a coherent $\D ^\dag _{\X/\fS,\Q}$-module which is also 
$\O _{\X,\Q}$-coherent, i.e. $\E$ is an object of $\mathrm{MIC} ^{\dag \dag} (\fX/K)$ 
(see notation \ref{ntnMICdag2fs}). 
\item $\E$ is $\O _{\X,\Q}$-coherent. 
\item $f ^* (\E)$ is $\O _{\Y,\Q}$-coherent. 
\item $f ^* (\E)$ is coherent $\D ^\dag _{\Y/\fS,\Q}$-module which is also 
$\O _{\Y,\Q}$-coherent, i.e. $f ^* (\E)$ is an object of $\mathrm{MIC} ^{\dag \dag} (\fY/K)$. 
\end{enumerate}

\end{lem}

\begin{proof}
Following \ref{thm-eqcat-cvisoc},
we get that the equivalence between 1 and 2 and between 3 and 4. 
Using \ref{2.2.13-carocourbe} and theorem of type $A$ for coherent $\O _{\X,\Q}$-modules, 
we get the equivalence between 2 and 3.
\end{proof}

\begin{coro}
\label{univhomeo-eqcat-isoc}
Let $f\colon \Y  \to \X  $ be 
a morphism of 
formal  $\fS$-schemes
of formal finite type
and having locally finite $p$-bases  over $\fS$.
We suppose  that the induced morphism 
$f _0 \colon Y \to X$ is a finite, surjective and radicial morphism. 

The functors 
$f _+$ and $f ^!$ are quasi-inverse equivalences of categories between 
\ref{ntnMICdag2fs}
$\mathrm{MIC} ^{\dag \dag} (\fX/K)$
and 
$\mathrm{MIC} ^{\dag \dag} (\Y/K)$ 
(see notation \ref{ntnMICdag2fs}).
 \end{coro}

\begin{proof}
This is a consequence of \ref{256Be2} and \ref{cvisoc-desc-finite}.
\end{proof}

\begin{coro}
\label{directsummandCVIsoc}
Let $f \colon \X \to \Y$ be a projective (in the sense of Definition \ref{projectivefscheme}), finite, surjective morphism of formal  $\fS$-scheme
of formal finite type and having finite $p$-bases  over $\fS$.
Let $\cE$ be an object of $\mathrm{MIC} ^{\dag \dag} (\Y/K)$ (see notation \ref{ntnMICdag2fs}).
Then $f ^! (\cE) \in \mathrm{MIC} ^{\dag \dag} (\fX/K)$ and 
$\cE$ is a direct summand of 
$f _+ f ^! (\cE)$.
\end{coro}

\begin{proof}
1) Since $f$ is projective, then we have the adjunction morphisms
$f _+ f ^! (\cE) \to \cE$ 
and 
$f _+ f ^! (\DD (\cE)) \to \DD (\cE)$ (see \ref{cor-adj-formul-proj-formal-bij2}).
By applying the dual functor to this latter morphism,
we get 
$$\cE \riso 
\DD \circ \DD (\cE) 
\riso 
\DD f _+ f ^! (\DD (\cE))
\underset{\ref{rel-dual-isom-proj-formal-iso}}{\riso}
 f _+ \DD f  ^! (\DD (\cE)).$$
By using \ref{dualisoscvdag} and \ref{sp*f*com}, we get
$\DD f  ^! (\DD (\cE))
\riso 
\left (f  ^* (\E ^\vee) \right) ^\vee
\riso 
f  ^* (\E ) 
\riso 
f  ^! (\E ) $.
Hence, we get 
\begin{equation}
\label{comp-f+f!+-drsumm}
\cE \to f _+ f ^! (\cE) \to \cE.
\end{equation}

2) We check in this step that the composition \ref{comp-f+f!+-drsumm} is an isomorphism.
Since this is local, we can suppose $\fX$ and $\fY$ affine and integral, and
there exists 
a relatively perfect morphism
$\varpi \colon \fY \to \widehat{\A}  ^d _\fS$. 
Let $L$ be the subextension of $k (X) /k(Y)$ such that
$L/k(Y)$ is separable and 
$k(X) / L$ is radicial (i.e. purely inseparable).
Let $\widetilde{X}$ be the normalisation of $Y$ in $L$.
Since $L/k(Y)$ is separable, then the canonical morphism
$h _0\colon \widetilde{X} \to Y$ is finite. 
Since the result that we have to check is local in $Y$, 
since $L/k(Y)$ is separable, then we can suppose 
$h _0$ is finite and étale.

Following \ref{EGAIV18.1.2} there exists 
a formally étale morphism
$\alpha \colon \widetilde{\fX} \to \widehat{\A}  ^d _\fS$ such that 
$\widetilde{\fX}$ is a formal $\fS$-scheme of formal finite type and 
the reduction of $\alpha$ modulo $\pi$
is $\varpi _0\circ h _0$.
Since $\varpi$ is formally étale, there exists (a unique)
$h \colon \widetilde{\fX} \to \fY$ making commutative the diagram
$$
\xymatrix{
{\widetilde{\fX}} 
\ar[rr] ^-{\alpha}
\ar@{.>}[rrd] ^-{h}
&& 
{ \widehat{\A}  ^d _\fS} 
\\ 
{\widetilde{X}} 
\ar@{^{(}->}[u] ^-{}
\ar[r] ^-{h _0}
& 
{Y}
\ar@{^{(}->}[r] ^-{}
&
{\fY.} 
\ar[u] ^-{\varpi}
}
$$

Since $X$ is normal, we get a morphism
$g _0 \colon X \to \widetilde{X} $
whose composition with $\widetilde{X} \to Y$
is $f _0$. 
Since $h$ is formally étale, 
then there exists a lifting $g \colon \fX \to \widetilde{\fX}$ of $g _0$
such that $h\circ g = f$.
Since $f$ is finite, then so is $g$. Following \ref{Ex5.3.9Liu},
this yields that $g _0$ is a universal homeomorphism. 
Hence, by using \ref{univhomeo-eqcat-isoc}, we reduce to the case where $g = id$, i.e. 
to the case where $f$ is finite and étale. Then, this is well known.
\end{proof}

\subsection{Differential coherence of the constant coefficient : the case of formal $\fS$-schemes
of finite type over $\V [[t]]$}
In this subsection, we suppose $k$ perfect (we need the perfectness in \ref{split-basechange}). 

\begin{dfn}
\label{dfn-morph-trait}
A morphism of complete discrete valuation rings $R \to  R'$ will refer
to a local ring homomorphism such that a uniformizer of $R$ is not mapped to zero (which is equivalent to saying that $R \to R'$ is injective
or that $\Spec R ' \to \Spec R $ is surjective).
A scheme $S$ is called a trait 
if it is isomorphic to a scheme of the form
$\Spec R$, where $R$ is a complete discrete valuation ring. 
A morphism of traits is a morphism $S'\to  S$ corresponding to a morphism of complete discrete valuation rings 
$R \to  R'$ as above. Such a morphism is said to be a finite extension of traits if the extension $S'\to S$ is finite.
Remark that in the case, since $R$ and $R'$ are regular, 
then 
the fact that $S'\to S$ is finite implies that 
$S'\to S$ is flat
\cite[4.3.11]{Liu-livre-02}
(or see \cite[15.4.2]{EGAIV4}, \cite[5.4.2]{EGAIV4}
and
\cite[0.17.3.5]{EGAIV1}).
\end{dfn}

\begin{empt}
\label{split-basechange}
Let $f \colon R \to  R'$ be a finite morphism of complete discrete valuation rings of equal characteristic $p>0$.
We denote by $\fm$ et $\fm'$ the maximal ideals of $R$ and $R'$.
We suppose $R /\fm$ is a perfect field. 
Let $l$ be the (unique) field of representative of $R$.
Since $l \to R '/\fm'$ is separable, then following
\cite[IX.\S 3, Proposition 1]{Bourbaki-AC89},
 there exists 
a field of representative $l'$ of $R'$ such that
$f (l) \subset l'$, i.e. $l'$ is a field extension of $l$.
(Beware that 
when $l \to R '/\fm'$ is not separable then this is not necessarily possible (see the counter-example
of \cite[VIII, Exercice 29, p. 98]{Bourbaki-AC89}).) 
Hence we can split $R \to R'$ in 
$R \to R\otimes _{l} l' \to R'$.
Since $l \to l'$ is a finite extension, then 
$R\otimes _{l} l'$ is a complete discrete valuation rings of equal characteristic $p>0$
and 
$R\otimes _{l} l' \to R'$ is a morphism of complete discrete valuation rings whose induced morphism of residue fields is an isomorphism.

The $l$-algebra $R$ is isomorphic to $l[[t]]$ 
and $l'$-algebra $R'$ is isomorphic to $l'[[u]]$. Hence $R \to R'$ corresponds to a morphism of $l$-algebras 
of the form $l [[t]] \to l' [[u]]$ which is decomposed in 
$l [[t]] \to l' [[t]] \to  l' [[u]]$.

Let $\W$ be complete discrete valuation ring of unequal characteristic and residue field $l$. 
Let $\W'$ be the unramified extension of $\W$ whose special fiber is 
$l'$. 
The canonical morphism of $\cW$-algebra
$\cW [[t]] \to \cW' [[t]]$ is a lifting of 
$l [[t]] \to l' [[t]]$.

\end{empt}

\begin{lemm}
\label{coh-ssp}
Let  $(X,Z)$ be a strict semi-stable pair over $\Spec R$ 
where $R = k [[t]]$ (see \ref{dfn-ssp}). 
Then 
$\O _{\fX} (\hdag T ) _\Q$ a coherent left 
$\cD ^{\dag} _{\fX,\bbQ}$-module.
\end{lemm}

\begin{proof}
Following \ref{proplocdesc-sstp}, this is an application of 
\ref{NCDgencoh}.
\end{proof}

\begin{thm}
[Berthelot]
\label{coh-cst-div}
Let $\X $ be a
formal  $\Spf \, \V [[t]]$-scheme of finite type
and having locally finite $p$-bases  over $\fS$.
Let $Z$ be a divisor of $X$.
Then 
$\O _{\X} (\hdag Z) _\Q$ is a coherent 
$\D ^\dag _{\X, \Q}$-module.
\end{thm}

\begin{proof}
We can adapt the proof of Berthelot of \cite{Becohdiff} as follows. 

0) Following Theorem \cite[6.5]{dejong}, 
there exist a trait 
$\mathbb{D} ^1 _{S'}=\Spec k '[[u]]$ (with $S ': =\Spec k'$) 
finite over 
$\mathbb{D} ^1 _{S}=\Spec k [[t]]$ such that the corresponding morphism
$k [[t]] \to k' [[u]] $ is a morphism of traits, 
a separated $\mathbb{D} ^1 _{S'}$-scheme of finite type $X '$, an alteration of schemes over $\mathbb{D} ^1 _{S}$ (in the sense of \cite[2.20]{dejong})
$g _0 \colon X ' \to X$ 
and an open immersion $j ' \colon X ' \to \overline{X} '$ of $\mathbb{D} ^1 _{S'}$-schemes, with the following properties: 
\begin{enumerate}[(a)]
\item $\overline{X} '$ is an integral projective $\mathbb{D} ^1 _{S'}$-scheme with geometrically irreducible generic fibre, and 
\item  the pair $(\overline{X} ', g _0 ^{-1} (Z ) _\mathrm{red} \cup (\overline{X} ' \setminus j ' (X '))$ is strict semi-stable.
\end{enumerate}
In particular, we get that 
$(X ', g _0  ^{-1} (Z ) _\mathrm{red})$ is a strict semi-stable pair
and there exists a closed immersion of the form
$u _0\colon X ' \hookrightarrow \P ^n _{X}$
whose composition with the projection 
$\P ^n _{X} \to X$ is $g _0 $.

1) Since $k$ is perfect, then 
following \ref{split-basechange},
replacing $k'$ by another field of representative if necessary,
we can suppose that the canonical diagram
$$
\xymatrix{
{\mathbb{D} ^1 _{S'}} 
\ar[d] ^-{}
\ar[r] ^-{}
&
 {\mathbb{D} ^1 _{S}} 
\ar[d] ^-{}
 \\ 
 {S'} 
\ar[r] ^-{}
&
 {S} 
 }
$$
is commutative. 
Hence, using 
\ref{split-basechange}
and 
\ref{desc-coh-chgbase},
we reduce to the case where $S' \to S$ is the identity.

2) i) 
Let 
$\fP := \widehat{\bbP} ^n _{\X}$, 
$f \colon \fP \to \X$ be the projection. 
Since $f$ is projective, we have the adjoint morphism
$ f _{+}  \circ f ^{!} (\O _{\X,\Q})
 \to 
\O _{\X,\Q}$
in 
$D ^{\mathrm{b}} _{\mathrm{coh}}( \smash{\D} ^\dag _{\X ,\Q} )$ (see 
\ref{cor-adj-formul-proj-formal-bij2}).
Following \ref{extriangleloc} and \ref{coh-smoothsubsch} (see also \ref{locdesc-sstvar}
and \ref{cor-closed-immer-local}), 
we have in $D ^{\mathrm{b}} _{\mathrm{coh}}( \smash{\D} ^\dag _{\fP /\fS ,\Q} )$
the morphism
$\R \underline{\Gamma} ^{\dag} _{X '}  (\O _{\fP,\Q})
\to 
\O _{\fP,\Q}$.
Since 
$f ^{!} (\O _{\X,\Q}) 
\riso 
\O _{\fP,\Q} [n]$, then we get 
the morphism
in 
$D ^{\mathrm{b}} _{\mathrm{coh}}( \smash{\D} ^\dag _{\X ,\Q} )$
\begin{equation}
\label{cstrcf+GammaO->O}
 f _{+} ( \R \underline{\Gamma} ^{\dag} _{X '} \O _{\fP,\Q} [n] )
 \to 
\O _{\X,\Q}.
\end{equation}

ii) In this step, 
we construct 
the morphism
$\O _{\X,\Q}
\to 
f _{+} ( \R \underline{\Gamma} ^{\dag} _{X '} \O _{\fP,\Q} [n])$ as follows:
 we have 
\begin{equation}
\label{cstrcf+GammaO->O-duality}
\O _{\X,\Q} 
\underset{\ref{dualisoscvdag}}{\riso}
\DD (\O _{\X,\Q})
\underset{\ref{cstrcf+GammaO->O}}{\longrightarrow} 
\DD f _{+}  ( \R \underline{\Gamma} ^{\dag} _{X '} \O _{\fP,\Q} [n])
\underset{\ref{rel-dual-isom-proj-formal-iso}}{\riso} 
f _{+} \DD  (  \R \underline{\Gamma} ^{\dag} _{X '}  \O _{\fP,\Q} [n])
\underset{\ref{propspetdualsansfrob}}{\riso} 
f _{+} ( \R \underline{\Gamma} ^{\dag} _{X '} \O _{\fP,\Q} [n]).
\end{equation}

iii) The composite morphism
$\O _{\X,\Q}
\to 
f _{+} ( \R \underline{\Gamma} ^{\dag} _{X '} \O _{\fP,\Q} [n])
 \to 
\O _{\X,\Q}$
in 
$D ^{\mathrm{b}} _{\mathrm{coh}}( \smash{\D} ^\dag _{\X ,\Q} )$
is an isomorphism.
Indeed, using Proposition \ref{lem-projff}.\ref{lem-projff-it3}, since this composition is a morphism of 
the abelian category $\mathrm{MIC} ^{\dag\dag} (\X /K) $,
we reduce to check that  its restriction to an open dense subset is an isomorphism.
Hence, we can suppose that 
$X$ affine,  
$X ' \to X$ is a finite surjective morphism,
and finally that $X$ and $X'$ have a finite $p$-basis over $S$.
Hence, by using \ref{lifting-pbasis}, 
there exists  a 
formal scheme $\fX '$ of formal finite type and having a finite $p$-basis over $\fS $
such that $\fX ' \times _{\fS} S \riso X ' $.
Since $X'$ is affine and $\fX'$ is noetherian, 
then for any integer $i$ the $S _i$-scheme 
$\fX ' \times _{\fS} S _i$ is affine (see  \cite[5.1.9]{EGAI}).
Hence 
$\fX' 
\riso 
\underrightarrow{\lim}
X' _i$, then $\fX'$ is affine. 
Since $\fP  /\fS$ is formally smooth and $\fX'$ is affine, then there exists a morphism
$u \colon \X '  \to \fP  $ which is a lifting of 
$u _0 \colon X '  \to P $.
We denote by 
$g := f \circ u
\colon 
\fX' \to \fX$, which is a lifting of . 
By using \ref{coro-trace-upre-BK}, 
we get that the morphism \ref{cstrcf+GammaO->O} (resp. 
\ref{cstrcf+GammaO->O-duality})
corresponds to  the trace map 
$g _+ g ^! (\O _{\X,\Q}) 
\to 
\O _{\X,\Q}$
(resp. 
to the map 
$\O _{\X,\Q}
\to 
g _+ g ^! (\O _{\X,\Q}) $
which is induced by duality from the trace map, 
i.e. is the adjunction morphism of $g _! = g_+$ and $g ^! = g ^+$).

Since $X' \to X$ is finite, since $X'$ is normal (and even regular),
then $X' \to X$ is the normalisation of $X$ in $k(X')$.
We can split 
the extension 
$k (  X ') / k (X)$ into 
a finite separable extension
and a finite radicial extension 
$k(X')/L$.
Let $X''$ be the normalization of $X$ in $L$. 
Then $X'' \to X$ is surjective and finite (see \cite[4.1.25]{Liu-livre-02}).
Shrinking $X$ if necessary, we can suppose
that $X'' \to X$ is finite, étale, surjective. 
Hence, $X' \to X$ is the composition of a universal homeomorphism
$X ' \to X''$ (use \ref{Ex5.3.9Liu}) with a finite, etale surjective morphism 
$X '' \to X  $.
Hence, by using \ref{lifting-pbasis} (and also  \cite[5.1.9]{EGAI}), 
there exists  an affine formal scheme $\fX ''$ of formal finite type and having a $p$-basis over $\fS $
such that $\fX '' \times _{\fS} S \riso X '' $.
As above, 
we get the liftings 
$\fX' \to \fX''$ and 
$\fX'' \to \fX$ of 
$X' \to X''$ 
and 
$X'' \to X$. 
Using \ref{univhomeo-eqcat-isoc}, we reduce to the case where 
$X ' \to X$ is finite, étale and surjective. 
In that case, 
$g _+ = g _*$ and $g ^! = g ^*$
and the trace maps are the trap maps computed in the categories of coherent $\cO _{\fX,\bbQ}$-modules,
which is well known.

3) 
Following the step 2), 
$\O _{\X,\Q}$ 
is a direct summand of
$f _{+} ( \R \underline{\Gamma} ^{\dag} _{X '} \O _{\fP,\Q} [n])$
in the category 
$D ^{\mathrm{b}} _{\mathrm{coh}}( \smash{\D} ^\dag _{\X /\fS ,\Q} )$.
This yields that 
$\O _{\X} (\hdag Z) _\Q$ 
is a direct summand of 
$(\hdag Z) f _{+} ( \R \underline{\Gamma} ^{\dag} _{X '} \O _{\fP,\Q} [n])$
in the category 
$D ^{\mathrm{b}} _{\mathrm{coh}}( \smash{\D} ^\dag _{\X /\fS } (\hdag Z) _{\Q} )$.
Using \ref{rema-fct-qcoh2coh} and \ref{surcoh2.1.4-cor},
we get in 
$D ^{\mathrm{b}} _{\mathrm{coh}}( \smash{\D} ^\dag _{\X /\fS } (\hdag Z) _{\Q} )$
the morphism
$$(\hdag Z) f _{+} ( \R \underline{\Gamma} ^{\dag} _{X '} \O _{\fP,\Q} [n]) \riso 
f _{Z,+} \circ (\hdag f ^{-1} (Z))  ( \R \underline{\Gamma} ^{\dag} _{X '} \O _{\fP,\Q} [n]).$$
Hence,  it is sufficient to
check that this latter object is 
$ \smash{\D} ^\dag _{\X /\fS ,\Q} $-coherent.
Since $f$ is proper and since 
$(\hdag f ^{-1} (Z))  ( \R \underline{\Gamma} ^{\dag} _{X '} \O _{\fP,\Q} [n])$ is 
already known to be $ \smash{\D} ^\dag _{\fP /\fS} (\hdag f ^{-1} (Z)) _{\Q} $-coherent,
using the remark \ref{oub-div-opcoh}.\ref{oub-div-opcoha)},
we reduce to check that 
$(\hdag f ^{-1} (Z))  ( \R \underline{\Gamma} ^{\dag} _{X '} \O _{\fP,\Q} [n])$ is 
$ \smash{\D} ^\dag _{\fP /\fS ,\Q} $-coherent.
Since this is local in $\fP $, we can suppose $\fP $ affine. 
Hence, there exists  a morphism 
$u \colon \X '  \to \fP  $
of formal schemes having locally finite $p$-bases over $\fS $
 which is $u _0 \colon X '  \to P $ modulo $\pi$.
 We get 
\begin{gather}
\notag
(\hdag f ^{-1} (Z)) (\R \underline{\Gamma} ^\dag _{X'} \O _{\fP,\Q}   [n])
\underset{\ref{coro-trace-upre}}{\riso} 
(\hdag f ^{-1} (Z)) (u _+ (\O _{\X',\Q}))
\underset{\ref{surcoh2.1.4-cor}}{\riso} 
 u _{f ^{-1} (Z),+} (\O _{\X'}(\hdag g _0 ^{-1} (Z)) _{\Q}).
\end{gather}
Since $(X ',g _0  ^{-1} (Z ))$ is a strict semi-stable pair, 
then following  \ref{coh-ssp},  
the left $ \smash{\D} ^\dag _{\X' /\fS ,\Q} $-module 
$\O _{\X'}(\hdag g _0 ^{-1} (Z)) _{\Q})$ is 
coherent.
Hence, 
using the remark \ref{oub-div-opcoh}.\ref{oub-div-opcoha)},
$ u _{f ^{-1} (Z),+} (\O _{\X'}(\hdag g _0 ^{-1} (Z)) _{\Q})
\riso
 u _{+} (\O _{\X'}(\hdag g _0 ^{-1} (Z)) _{\Q})$
 is 
$ \smash{\D} ^\dag _{\fP/\fS ,\Q} $-coherent. 
\end{proof}
\begin{dfn}
\label{nice-div}
Let $\fP \to \bbD ^r _{\fS}$ be an object of $\scr{C} _{\fS}$ (see \ref{dfn-CfS}) such that 
$\fP /\fS$ has locally finite $p$-bases.
Let $X$ be a closed subscheme of the special fiber of $\fP$
and having locally finite $p$-bases  over $S$.
Let $T$ be a divisor of $X$.
We say that $T$ is a ``nice divisor of $X/S$'' if Zariski locally in $T$ there exists 
\begin{enumerate}[(a)]
\item a finite morphism $\V \to \V'$ of complete discrete valuation rings of mixed characteristics $(0,p)$, 

\item a finite morphism of formal schemes $\bbD ^r _{\fS '}\to \bbD ^r _{\fS}$ with $\fS ' := \Spf~\cV '$ making commutative the diagram 
$$
\xymatrix{
{\mathbb{D} ^r _{\fS'}} 
\ar[d] ^-{}
\ar[r] ^-{}
&
 {\mathbb{D} ^r _{\fS }} 
\ar[d] ^-{}
 \\ 
 {\fS'} 
\ar[r] ^-{}
&
 {\fS.} 
 }
$$

\item a projective morphism $g\colon \fP' \to \fP$ of formal schemes
such that $\fP ' /\fS '$ has locally finite $p$-bases, 
a closed subscheme $X' \subset g _0 ^{-1} (X)$ of the special fiber $P'$
such that $X' /S$ has locally finite $p$-base
and the induced morphism $\phi \colon X ' \to X$ is 
an alteration of schemes (in the sense of \cite[2.20]{dejong})
and 
$\phi  ^{-1} (T ) _\mathrm{red}$ is a strictly nice divisor of $X'/S'$
(see Definition \ref{st-nice-div}).
\end{enumerate}

\end{dfn}

\begin{ex}
\label{ex-nice-div}
Following de Jong desingularisation Theorem \cite[6.5]{dejong} (see the proof \ref{coh-cst-div})
in the case where $r =1$
any divisor is an nice divisor. 
\end{ex}

\begin{thm}
[Berthelot]
\label{coh-cst-nicediv}
Let $\X \to \bbD ^r _{\fS}$ be an object of $\scr{C} _{\fS}$ (see \ref{dfn-CfS}) such that 
$\fX /\fS$ has locally finite $p$-bases.
Let $Z$ be a nice divisor of $X/S$.
Then 
$\O _{\X} (\hdag Z) _\Q$ is a coherent 
$\D ^\dag _{\X, \Q}$-module.
\end{thm}

\begin{proof}
We copy word by word the proof of \ref{coh-cst-div}.
\end{proof}

\begin{empt}
[Warning]
\label{warning}
Let $\X \to \bbD ^r _{\fS}$ be an object of $\scr{C} _{\fS}$ (see \ref{dfn-CfS}) such that 
$\fX /\fS$ has locally finite $p$-bases.
This is not clear that any divisor 
$Z$ of $X$ is a nice divisor of $X/S$.
I do not know if 
$\O _{\X} (\hdag Z) _\Q$ is a coherent 
$\D ^\dag _{\X, \Q}$-module, even if we can hope so.
\end{empt}

\begin{coro}
\label{coh-Bbullet}
With notation \ref{coh-cst-nicediv}, we have
$\smash{\widetilde{\B}} _{\X} ^{(\bullet)} (Z ) \in 
\underrightarrow{LM}  _{\Q, \mathrm{coh}} (\smash{\widehat{\D}} _{\X } ^{(\bullet)})
\cap 
\underrightarrow{LM}  _{\Q, \mathrm{coh}} (\smash{\widetilde{\D}} _{\X } ^{(\bullet)} (Z))$.
\end{coro}

\begin{proof}
We already know that
$\smash{\widetilde{\B}} _{\X} ^{(\bullet)} (Z ) 
\in \underrightarrow{LM}  _{\Q, \mathrm{coh}} (\smash{\widetilde{\D}} _{\X } ^{(\bullet)} (Z))$.
Following \ref{coh-cst-div}, 
$\O _{\X} (\hdag Z) _{\Q}=
\underrightarrow{\lim}
\smash{\widetilde{\B}} _{\X} ^{(\bullet)} (Z ) $
is a coherent
$\smash{\D} ^\dag _{\X,\Q}$-module.
Using  \ref{coro1limTouD}, we can conclude.
\end{proof}

\begin{thm}
\label{coh-ss-div-bis}
Let  $\fP $ be a 
formal  $\fS$-scheme
of formal finite type
and having locally finite $p$-bases  over $\fS$.
Let $X$ be a 
closed subscheme of $P$
and having locally finite $p$-bases  over $\Spec k$.
We suppose also that 
$X$ is of finite type over 
$\Spec ~k[[t]]$. 
Let 
$\E ^{(\bullet)}$ be an object of 
$\mathrm{MIC} ^{(\bullet)} (X, \fP/K) $.
Then 
$\E ^{(\bullet)}
\in 
\smash{\underrightarrow{LM}}  _{\Q, \mathrm{ovcoh}}
(\smash{\widehat{\D}} _{\fP /\fS } ^{(\bullet)})$.
\end{thm}

\begin{proof}
By stability under inverse images of 
$\mathrm{MIC} ^{(\bullet)} (X, \fP/K) $,
we reduce to check that for any
divisor $T$ of $X$,
we have 
$(\hdag T ) (\E ^{(\bullet)})
\in 
\smash{\underrightarrow{LM}}  _{\Q, \mathrm{coh}}
(\smash{\widehat{\D}} _{\fP /\fS } ^{(\bullet)})$.
Using the inductive system version of Berthelot-Kashiwara's theorem
(see \ref{u!u+=id}),
we reduce to the case where $X = P$. 
In this case, we write $\X$ (resp. $X$) instead of $\fP$
(resp. $T$) and we will use the notation of the proof of \ref{coh-cst-div}.
Now, 
similarly to the part 1) of the proof of 
\ref{coh-cst-div},
we reduce to the case where $S =S'$.
Finally we can copy the proof
of \cite[10.2.3]{caro-6operations}.
\end{proof}

\begin{thm}
\label{coh-ss-div-bis-nice}
Let $\fP \to \bbD ^r _{\fS}$ be an object of $\scr{C} _{\fS}$ (see \ref{dfn-CfS}) such that 
$\fP /\fS$ has locally finite $p$-bases.
Let $X$ be a 
closed subscheme of $P$
and having locally finite $p$-bases  over $\Spec k$.
Let 
$\E ^{(\bullet)}$ be an object of 
$\mathrm{MIC} ^{(\bullet)} (X, \fP/K) $.
Then for any nice divisor 
$T$ of $X$ (see definition \ref{nice-div}), 
$(\hdag T) (\E ^{(\bullet)})
\in 
\smash{\underrightarrow{LM}}  _{\Q, \mathrm{coh}}
(\smash{\widehat{\D}} _{\fP /\fS } ^{(\bullet)})$.
\end{thm}

\begin{proof}
This is checked similarly to \ref{coh-cst-nicediv}.
\end{proof}

We will need later the following proposition.

\begin{prop}
\label{ind-desc-coh-chgbase}
Let $\V \to \V'$ be a finite morphism of complete discrete valuation rings of mixed characteristics $(0,p)$.
We get the finite morphism $\fS ' := \Spf \V ' \to \fS$.
Let  $\X \to \bbD ^r _{\fS}$ be an object of $\scr{C} _{\fS}$ 
such that $\fX /\fS$ has locally finite $p$-bases.
Let
$\X ' := \X \times _{\fS} \fS '$, 
and $f \colon \X' \to \X$ be the canonical projection.
Let $Z$ be a divisor of $X$ and $Z':= f ^{-1} (Z)$.
Let 
$\E ^{(\bullet)}
\in 
\smash{\underrightarrow{LD}} ^{\mathrm{b}} _{\Q,\mathrm{coh}} ( \smash{\widehat{\D}} _{\X} ^{(\bullet)})$.
Let 
$ \E ^{\prime (\bullet)}
 :=
  \V' \otimes 
_{\V}  \E ^{(\bullet)}$.
If $(\hdag Z ' )  ( \E ^{\prime (\bullet)})
\in 
\smash{\underrightarrow{LD}} ^{\mathrm{b}} _{\Q,\mathrm{coh}} ( \smash{\widehat{\D}} _{\X '} ^{(\bullet)})$,
then
$(\hdag Z )  ( \E ^{(\bullet)})
\in 
\smash{\underrightarrow{LD}} ^{\mathrm{b}} _{\Q,\mathrm{coh}} ( \smash{\widehat{\D}} _{\X} ^{(\bullet)})$.
\end{prop}

\begin{proof}
Using  \ref{coro1limTouD}, this is a consequence of 
Lemma \ref{desc-coh-chgbase}.
\end{proof}

\section{Local cohomological functors}
\label{chapter12}
\subsection{Local cohomological functor with strict support over a divisor}

Let  $\fP $ be a 
formal  $\fS$-scheme
of formal finite type
and having locally finite $p$-bases  over $\fS$.
Let $T$ be a divisor of $P$.
We have already defined in \ref{hdagT-nota}
the localisation functor $(\hdag T)$ outside $T$. 
In this subsection, we define and study 
the local cohomological functor with support in $T$, which we denote by
$\R \underline{\Gamma} ^\dag _{T} $.

\begin{lemm}
\label{annulationHom-hdag}
\begin{enumerate}[(a)]
\item 
\label{annulationHom-hdag-item1}
Let 
$\FF ^{(\bullet)}
\to 
\E ^{(\bullet)} 
\to 
(\hdag T) (\E ^{(\bullet)} )  
 \to 
\FF ^{(\bullet)} [1] 
$
be a distinguished triangle of 
$\smash{\underrightarrow{LD}} ^\mathrm{b} _{\Q,\mathrm{qc}} ( \smash{\widehat{\D}} _{\fP /\fS } ^{(\bullet)})$
where the second arrow is the canonical morphism. 
For any divisor $T \subset T'$, 
we have the isomorphism 
$(\hdag T') (\FF ^{(\bullet)} )\riso 0$
of 
$\smash{\underrightarrow{LD}} ^\mathrm{b} _{\Q,\mathrm{qc}} ( \smash{\widehat{\D}} _{\fP /\fS } ^{(\bullet)})$.

\item 
\label{annulationHom-hdag-item2}
Let 
$\E ^{(\bullet)}
\in \smash{\underrightarrow{LD}} ^\mathrm{b} _{\Q,\mathrm{qc}} ( \smash{\widehat{\D}} _{\fP /\fS } ^{(\bullet)})$
et 
$\FF ^{(\bullet)}
\in \smash{\underrightarrow{LD}} ^\mathrm{b} _{\Q,\mathrm{qc}} ( \smash{\widetilde{\D}} _{\fP /\fS } ^{(\bullet) } (T))$.
We suppose we have in 
$\smash{\underrightarrow{LD}} ^\mathrm{b} _{\Q,\mathrm{qc}} ( \smash{\widehat{\D}} _{\fP /\fS } ^{(\bullet)})$
the isomorphism
$(\hdag T) (\E ^{(\bullet)} )\riso 0$.
Then
$\mathrm{Hom} _{\smash{\underrightarrow{LD}}  _{\Q} ( \smash{\widehat{\D}} _{\fP /\fS } ^{(\bullet)})}
(\E ^{(\bullet)} , \FF ^{(\bullet)}) =0.$

\end{enumerate}

\end{lemm}

\begin{proof}
Using \ref{hdagT'T=cup}, this is checked similarly to
\cite[4.1.2 and 4.1.3]{caro-stab-sys-ind-surcoh}.
\end{proof}

\begin{empt}
\label{exist-RrmHom}
Let  $\mathrm{Ab}$ be the category of abelian groups. 
Similarly to 
\cite[1.4.2]{caro-stab-sys-ind-surcoh}, we construct the bifunctor
(which is the standard construction of the homomorphim bifunctor of the abelian category
$\underrightarrow{LM} _{\Q} (\smash{\widetilde{\D}} _{\fP /\fS } ^{(\bullet)} (T))$):
$$\mathrm{Hom} ^{\bullet} _{\underrightarrow{LM} _{\Q} (\smash{\widetilde{\D}} _{\fP /\fS } ^{(\bullet)} (T))}
(-,-)
\colon 
K 
(\underrightarrow{LM} _{\Q} (\smash{\widetilde{\D}} _{\fP /\fS } ^{(\bullet)} (T))) ^\circ 
\times
K 
(\underrightarrow{LM} _{\Q} (\smash{\widetilde{\D}} _{\fP /\fS } ^{(\bullet)} (T))) 
\to
K (\mathrm{Ab}).$$
Similarly to 
\cite[1.4.7]{caro-stab-sys-ind-surcoh},
we check that the bifunctor 
$\mathrm{Hom} ^{\bullet} _{\underrightarrow{LM} _{\Q} (\smash{\widetilde{\D}} _{\fP /\fS } ^{(\bullet)} (T))} (-,-)$ 
is right localizable. We get the bifunctor 
$$\R \mathrm{Hom} _{D (\underrightarrow{LM} _{\Q} (\smash{\widetilde{\D}} _{\fP /\fS } ^{(\bullet)} (T)))} (-,-)
\colon 
D ^{\mathrm{b}}(\underrightarrow{LM} _{\Q} (\smash{\widetilde{\D}} _{\fP /\fS } ^{(\bullet)} (T))) ^\circ 
\times
D ^{\mathrm{b}} (\underrightarrow{LM} _{\Q} (\smash{\widetilde{\D}} _{\fP /\fS } ^{(\bullet)} (T)))
\to 
D (\mathrm{Ab}).$$
Moreover, we have the isomorphism of bifunctors
$D ^{\mathrm{b}}(\underrightarrow{LM} _{\Q} (\smash{\widetilde{\D}} _{\fP /\fS } ^{(\bullet)} (T))) ^\circ 
\times
D ^{\mathrm{b}} (\underrightarrow{LM} _{\Q} (\smash{\widetilde{\D}} _{\fP /\fS } ^{(\bullet)} (T)))
\to 
\mathrm{Ab} $ of the form:
\begin{equation}
\label{H0Homrm-DLM}
\mathcal{H} ^{0} (\R \mathrm{Hom} _{D (\underrightarrow{LM} _{\Q} (\smash{\widetilde{\D}} _{\fP /\fS } ^{(\bullet)} (T)))} (-,-))
\riso 
\mathrm{Hom} _{D (\underrightarrow{LM} _{\Q} (\smash{\widetilde{\D}} _{\fP /\fS } ^{(\bullet)} (T)))} (-,-).
\end{equation}
\end{empt}

\begin{empt}
\label{GammaT}
Let  $T \subset T'$ be a second divisor.
Suppose we have the commutative diagram in 
$\smash{\underrightarrow{LD}} ^\mathrm{b} _{\Q,\mathrm{qc}} ( \smash{\widehat{\D}} _{\fP /\fS } ^{(\bullet)})$
of the form
\begin{equation}
\label{prefonct-GammaT}
\xymatrix @=0,4cm{
{\FF ^{(\bullet)}  } 
\ar[r] ^-{}
& 
{\E ^{(\bullet)}  } 
\ar[r] ^-{}
\ar[d] ^-{\phi}
& 
{(\hdag T) (\E ^{(\bullet)} )} 
\ar[d] ^-{(\hdag T)(\phi)}
\ar[r] ^-{}
&
{\FF ^{(\bullet)}  [1]} 
\\ 
{\FF ^{\prime (\bullet)}  } 
\ar[r] ^-{}
& 
{\E ^{\prime (\bullet)}  } 
\ar[r] ^-{}
& 
{(\hdag T) (\E ^{\prime (\bullet)} )} 
\ar[r] ^-{}
&
{\FF ^{\prime (\bullet)}   [1]} 
}
\end{equation}
where middle horizontal morphisms are the canonical ones and
where both horizontal triangles are distinguished. 
Modulo  the equivalence of categories
$\underrightarrow{LD} ^{\mathrm{b}} _{\Q} (\smash{\widetilde{\D}} _{\fP /\fS } ^{(\bullet)} (T))
\cong 
D ^{\mathrm{b}}
(\underrightarrow{LM} _{\Q} (\smash{\widetilde{\D}} _{\fP /\fS } ^{(\bullet)} (T)))$
(see \ref{eqcatLD=DSM-fonct})
which allows us to see \ref{prefonct-GammaT} as a diagram of 
$D ^{\mathrm{b}}
(\underrightarrow{LM} _{\Q} (\smash{\widetilde{\D}} _{\fP /\fS } ^{(\bullet)} (T)))$,
we have
$$H ^{-1} (\R \mathrm{Hom} _{D (\underrightarrow{LM} _{\Q} (\smash{\widetilde{\D}} _{\fP /\fS } ^{(\bullet)} (T)))}
( \FF ^{(\bullet)}  ,(\hdag T) (\E ^{\prime (\bullet)} )))
\underset{\ref{H0Homrm-DLM}}{\riso}
\mathrm{Hom} _{D (\underrightarrow{LM} _{\Q} (\smash{\widetilde{\D}} _{\fP /\fS } ^{(\bullet)} (T)))}
( \FF ^{(\bullet)}  ,(\hdag T) (\E ^{\prime (\bullet)}) [-1]) 
\underset{\ref{annulationHom-hdag}}{=}0.$$
Following 
\cite[1.1.9]{BBD},
this implies there exists a unique morphism
$\FF ^{(\bullet)}\to \FF ^{\prime (\bullet)} $
making commutative in 
$\smash{\underrightarrow{LD}} ^\mathrm{b} _{\Q,\mathrm{qc}} ( \smash{\widehat{\D}} _{\fP /\fS } ^{(\bullet)})$
the diagram:
\begin{equation}
\label{fonct-hdagT}
\xymatrix @=0,4cm{
{\FF ^{(\bullet)}  } 
\ar[r] ^-{}
\ar@{.>}[d] ^-{\exists !}
& 
{\E ^{(\bullet)}  } 
\ar[r] ^-{}
\ar[d] ^-{\phi}
& 
{(\hdag T) (\E ^{(\bullet)} )} 
\ar[d] ^-{(\hdag T)(\phi)}
\ar[r] ^-{}
& 
{\FF ^{(\bullet)}  [1]} 
\ar@{.>}[d] ^-{\exists !}
\\ 
{\FF ^{\prime (\bullet)}  } 
\ar[r] ^-{}
& 
{\E ^{\prime (\bullet)}  } 
\ar[r] ^-{}
& 
{(\hdag T) (\E ^{\prime (\bullet)} )} 
\ar[r] ^-{}
& 
{\FF ^{\prime (\bullet)} [1].} 
}
\end{equation}
Similarly to \cite[1.1.10]{BBD}, 
this implies that the cone of 
$\E ^{(\bullet)}
\to  
(\hdag T) (\E ^{(\bullet)} ) $
is unique up to canonical isomorphism. 
Hence, such a complex $\FF ^{ (\bullet)}$ is unique up to canonical isomorphism.
We denote it by 
$ \R \underline{\Gamma} ^\dag _{T} (\E ^{(\bullet)})$.
Moreover, 
the complex 
$\R \underline{\Gamma} ^\dag _{T} (\E ^{(\bullet)}) $ is functorial in
$\E ^{(\bullet)}$.
\end{empt}

\begin{dfn}
With notation \ref{GammaT}, 
the functor
$\R \underline{\Gamma} ^\dag _{T} 
\colon 
\smash{\underrightarrow{LD}} ^\mathrm{b} _{\Q,\mathrm{qc}} ( \smash{\widehat{\D}} _{\fP /\fS } ^{(\bullet)})
\to 
\smash{\underrightarrow{LD}} ^\mathrm{b} _{\Q,\mathrm{qc}} ( \smash{\widehat{\D}} _{\fP /\fS } ^{(\bullet)})
$
is 
the ``local cohomological functor with strict support over the divisor $T$''. 
For $\E ^{(\bullet)} \in 
\smash{\underrightarrow{LD}} ^\mathrm{b} _{\Q,\mathrm{qc}} ( \smash{\widehat{\D}} _{\fP /\fS } ^{(\bullet)})$, 
we denote by $\Delta _{T} (\E ^{(\bullet)})$  the canonical exact triangle 
\begin{equation}
\label{tri-loc-berthelot}
 \R \underline{\Gamma} ^\dag _{T} (\E ^{(\bullet)})
\to 
\E ^{(\bullet)}
\to 
(\hdag T) (\E ^{(\bullet)})
\to 
 \R \underline{\Gamma} ^\dag _{T} (\E ^{(\bullet)})
 [1].
\end{equation}

Let $\fU : = \fP \setminus T$. Remark that 
since the restriction to $\fU$ of the canonical morphism
$\E ^{(\bullet)}
\to 
(\hdag T) (\E ^{(\bullet)})$
is an isomorphism, then 
$ \R \underline{\Gamma} ^\dag _{T} (\E ^{(\bullet)}) |\fU = 0$, which justifies the terminology.
\end{dfn}

\begin{lemm}
Let  $T \subset T'$ be a second divisor,
and 
$\E ^{(\bullet)} 
\in 
\smash{\underrightarrow{LD}} ^\mathrm{b} _{\Q,\mathrm{qc}} ( \smash{\widehat{\D}} _{\fP /\fS } ^{(\bullet)})$.
There exists a unique morphism
$\R \underline{\Gamma} ^\dag _{T} (\E ^{(\bullet)}) 
\to
\R \underline{\Gamma} ^\dag _{T'} (\E ^{(\bullet)}) $
making commutative the following diagram
\begin{equation}
\label{fonct-hdagX2}
\xymatrix @=0,4cm{
{\R \underline{\Gamma} ^\dag _{T} (\E ^{(\bullet)}) } 
\ar[r] ^-{}
\ar@{.>}[d] ^-{\exists !}
& 
{\E ^{(\bullet)}  } 
\ar[r] ^-{}
\ar@{=}[d] ^-{}
& 
{(\hdag T) (\E ^{(\bullet)} )} 
\ar[d] ^-{}
\ar[r] ^-{}
& 
{\R \underline{\Gamma} ^\dag _{T} (\E ^{(\bullet)})[1]} 
\ar@{.>}[d] ^-{\exists !}
\\ 
{\R \underline{\Gamma} ^\dag _{T'} (\E ^{ (\bullet)})  } 
\ar[r] ^-{}
& 
{\E ^{ (\bullet)}  } 
\ar[r] ^-{}
& 
{(\hdag T') (\E ^{ (\bullet)} )} 
\ar[r] ^-{}
& 
{\R \underline{\Gamma} ^\dag _{T'} (\E ^{ (\bullet)})[1].} 
}
\end{equation}
In other words, 
$\R \underline{\Gamma} ^\dag _{T} (\E ^{(\bullet)}) $
is functorial in $T$.
\end{lemm}

\begin{proof}
We can copy 
\cite[4.1.4.3]{caro-stab-sys-ind-surcoh}.
\end{proof}

\begin{empt}
[Commutation with tensor products]
\label{iso-comm-locaux-prod-tens}
Let  $\E ^{(\bullet)},~\FF ^{(\bullet)}
\in \smash{\underrightarrow{LD}} ^\mathrm{b} _{\Q,\mathrm{qc}} ( \smash{\widehat{\D}} _{\fP /\fS } ^{(\bullet)})$.
By commutativity and associativity of tensor products, 
we have the canonical isomorphisms
$$(\hdag T) (\E ^{(\bullet)})
\smash{\widehat{\otimes}}^\L  _{\O ^{(\bullet)}  _{\fP} } 
\FF ^{(\bullet)} 
\riso 
(\hdag T) (\E ^{(\bullet)}
\smash{\widehat{\otimes}}^\L  _{\O ^{(\bullet)}  _{\fP} } 
\FF ^{(\bullet)}) 
\riso 
 \E ^{(\bullet)}
\smash{\widehat{\otimes}}^\L  _{\O ^{(\bullet)}  _{\fP} } 
(\hdag T) (\FF ^{(\bullet)}) .$$
Hence, there exists a unique isomorphism  of the form
$\R \underline{\Gamma} ^\dag _{T} 
(\E ^{(\bullet)}
\smash{\widehat{\otimes}}^\L  _{\O ^{(\bullet)}  _{\fP} } 
\FF ^{(\bullet)})
\riso 
\R \underline{\Gamma} ^\dag _{T}( \E ^{(\bullet)})
\smash{\widehat{\otimes}}^\L  _{\O ^{(\bullet)}  _{\fP} } 
\FF ^{(\bullet)}$
(resp. $\R \underline{\Gamma} ^\dag _{T} 
(\E ^{(\bullet)}
\smash{\widehat{\otimes}}^\L  _{\O ^{(\bullet)}  _{\fP} } 
\FF ^{(\bullet)})
\riso 
\E ^{(\bullet)}
\smash{\widehat{\otimes}}^\L  _{\O ^{(\bullet)}  _{\fP} } 
\R \underline{\Gamma} ^\dag _{T}
(\FF ^{(\bullet)})$)
making commutative the following diagram
\begin{equation}
\label{fonct-hdagXbis}
\xymatrix @=0,4cm{
{\R \underline{\Gamma} ^\dag _{T} 
(\E ^{(\bullet)})
\smash{\widehat{\otimes}}^\L  _{\O ^{(\bullet)}  _{\fP} } 
\FF ^{(\bullet)}} 
\ar[r] ^-{}
& 
{\E ^{(\bullet)}
\smash{\widehat{\otimes}}^\L  _{\O ^{(\bullet)}  _{\fP} } 
\FF ^{(\bullet)} } 
\ar@{=}[d] ^-{}
\ar[r] ^-{}
& 
{(\hdag T) (\E ^{(\bullet)})
\smash{\widehat{\otimes}}^\L  _{\O ^{(\bullet)}  _{\fP} } 
\FF ^{(\bullet)} } 
\ar[r] ^-{}
& 
{\R \underline{\Gamma} ^\dag _{T} 
(\E ^{(\bullet)} )
\smash{\widehat{\otimes}}^\L  _{\O ^{(\bullet)}  _{\fP} } 
\FF ^{(\bullet)} [1]}
\\
{\R \underline{\Gamma} ^\dag _{T} 
(\E ^{(\bullet)}
\smash{\widehat{\otimes}}^\L  _{\O ^{(\bullet)}  _{\fP} } 
\FF ^{(\bullet)})} 
\ar[r] ^-{}
\ar@{.>}[d] ^-{\exists !}
\ar@{.>}[u] ^-{\exists !}
& 
{\E ^{(\bullet)}
\smash{\widehat{\otimes}}^\L  _{\O ^{(\bullet)}  _{\fP} } 
\FF ^{(\bullet)} } 
\ar@{=}[d] ^-{}
\ar[r] ^-{}
& 
{(\hdag T) (\E ^{(\bullet)}
\smash{\widehat{\otimes}}^\L  _{\O ^{(\bullet)}  _{\fP} } 
\FF ^{(\bullet)})  } 
\ar[d] ^-{\sim}
\ar[r] ^-{}
\ar[u] ^-{\sim}
& 
{\R \underline{\Gamma} ^\dag _{T} 
(\E ^{(\bullet)}
\smash{\widehat{\otimes}}^\L  _{\O ^{(\bullet)}  _{\fP} } 
\FF ^{(\bullet)}) [1]}
\ar@{.>}[d] ^-{\exists !}
\ar@{.>}[u] ^-{\exists !}
\\ 
 {\E ^{(\bullet)}
\smash{\widehat{\otimes}}^\L  _{\O ^{(\bullet)}  _{\fP} } 
\R \underline{\Gamma} ^\dag _{T}
(\FF ^{(\bullet)})}
\ar[r] ^-{}
& 
{\E ^{(\bullet)}
\smash{\widehat{\otimes}}^\L  _{\O ^{(\bullet)}  _{\fP} } 
\FF ^{(\bullet)} } 
\ar[r] ^-{}
& 
{ \E ^{(\bullet)}
\smash{\widehat{\otimes}}^\L  _{\O ^{(\bullet)}  _{\fP} } 
(\hdag T) (\FF ^{(\bullet)})  }\ar[r] ^-{}
& 
{\E ^{(\bullet)}
\smash{\widehat{\otimes}}^\L  _{\O ^{(\bullet)}  _{\fP} } 
\R \underline{\Gamma} ^\dag _{T}
(\FF ^{(\bullet)})[1].} 
}
\end{equation}
Theses isomorphisms are functorial in 
$\E ^{(\bullet)},~\FF ^{(\bullet)},~ T$ (for the meaning of the functoriality in $T$, see  \ref{fonct-hdagX2}).

\end{empt}

\begin{empt}
[Commutation between local cohomological functors and localization functors]
\label{iso-comm-locaux}
Let $T _1, T _2$ be two divisors of  $P$, 
$\E ^{(\bullet)}
\in \smash{\underrightarrow{LD}} ^\mathrm{b} _{\Q,\mathrm{qc}} ( \smash{\widehat{\D}} _{\fP /\fS } ^{(\bullet)})$.

\begin{enumerate}[(a)]
\item By commutativity of the tensor product, 
we have the functorial in  $T _1$, $T _2$ and $\E ^{(\bullet)}$ canonical isomorphism
\begin{equation}
\label{commhdagT1T2}
(\hdag T _2) \circ (\hdag T _1) (\E ^{(\bullet)})
\riso 
(\hdag T _1) \circ (\hdag T _2) (\E ^{(\bullet)}).
\end{equation}

\item  There exists a unique  isomorphism
$(\hdag T _2) \circ \R \underline{\Gamma} ^\dag _{T _1}(\E ^{(\bullet)} )
\riso 
\R \underline{\Gamma} ^\dag _{T _1}
\circ
(\hdag T _2)(\E ^{(\bullet)} )$
inducing the canonical morphism of triangles
$ (\hdag T _2) ( \Delta _{T _1} (\E ^{(\bullet)} ))
\to 
\Delta _{T _1} ((\hdag T _2) (\E ^{(\bullet)} ))$ (see \cite[4.2.2.2]{caro-stab-sys-ind-surcoh}).
This isomorphism is functorial in  $T _1$, $T _2$, $\E ^{(\bullet)}$.

\item Similarly there exists a unique isomorphism 
$\R \underline{\Gamma} ^\dag _{T _2} \circ \R \underline{\Gamma} ^\dag _{T _1} (\E ^{(\bullet)})
\riso 
\R \underline{\Gamma} ^\dag _{T _1}\circ \R \underline{\Gamma} ^\dag _{T _2} (\E ^{(\bullet)})$
functorial in $T _1$, $T _2$, $\E ^{(\bullet)}$
and inducing the canonical morphism of triangles 
$\Delta _{T _2} (\R \underline{\Gamma} ^\dag _{T _1} (\E ^{(\bullet)} ))
\to 
\R \underline{\Gamma} ^\dag _{T _1} ( \Delta _{T _2} (\E ^{(\bullet)} ))$.

\end{enumerate}

\end{empt}

\begin{empt}
The three isomorphisms of \ref{iso-comm-locaux} are compatible with that 
of \ref{iso-comm-locaux-prod-tens} (for more precision, see  \cite[4.2.3]{caro-stab-sys-ind-surcoh}).
\end{empt}

We will need the following Lemmas (e.g. see the construction of \ref{dfn-4.3.4} or Proposition \ref{prop2.2.9}) in the next section.
\begin{lemm}
\label{lemme2.2.3}
Let  $D,~T$ be two divisors of $P$, 
$\E ^{(\bullet)}
\in 
\smash{\underrightarrow{LD}} ^{\mathrm{b}} _{\Q,\mathrm{coh}} ( \smash{\widetilde{\D}} _{\fP /\fS } ^{(\bullet)} (D))$,
$\U$ be the open subset of  $\fP$ complementary to the support of $T$.
The following assertions are equivalent :
\begin{enumerate}[(a)]
\item We have in  $\smash{\underrightarrow{LD}} ^{\mathrm{b}} _{\Q,\mathrm{coh}} ( \smash{\widetilde{\D}} _{\U /\fS } ^{(\bullet)} (D\cap U))$
the isomorphism $\E ^{(\bullet)}|\U \riso 0$.
\item The canonical morphism 
$\R \underline{\Gamma} ^\dag _{T} (\E ^{(\bullet)})
\to 
\E ^{(\bullet)}$ of $\smash{\underrightarrow{LD}} ^{\mathrm{b}} _{\Q,\mathrm{qc}} ( \smash{\widetilde{\D}} _{\fP /\fS } ^{(\bullet)} (D))$
is an isomorphism.
\item We have in $\smash{\underrightarrow{LD}} ^{\mathrm{b}} _{\Q,\mathrm{qc}} ( \smash{\widetilde{\D}} _{\fP /\fS } ^{(\bullet)} (D))$
the isomorphism $(\hdag T) (\E ^{(\bullet)} )\riso 0$.
\end{enumerate}
\end{lemm}

\begin{proof}
We can copy the proof of \cite[4.3.2]{caro-stab-sys-ind-surcoh}.
\end{proof}

\begin{cor}
\label{cor-induction-div-coh}
Let $\fP \to \bbD ^1 _{\fS}$ be a finite type morphism of 
formal  $\Spf \, \V$-schemes
having locally finite $p$-bases  over $\fS$.
Let $T _1, \dots, T _{r}$ be some divisors of $P$.
Let $T$ be a divisor of $P$.
Then
$\R \underline{\Gamma} ^\dag _{T _r} \circ \cdots \circ 
\R \underline{\Gamma} ^\dag _{T _1} (\smash{\widetilde{\B}} _{\fP} ^{(\bullet)} (T  ))
\in 
\smash{\underrightarrow{LD}} ^{\mathrm{b}} _{\Q,\mathrm{coh}} ( \smash{\widehat{\D}} _{\fP /\fS } ^{(\bullet)} )$.
Moreover, there exists 
a canonical isomorphism
$$\R \sp _* 
\left (
\underline{\Gamma} ^\dag _{T _r} 
\circ \cdots \circ 
\underline{\Gamma} ^\dag _{T _1} ( j ^\dag _{T}\O _{\fP _K})
\right )
\riso 
\underrightarrow{\lim}\,
\R \underline{\Gamma} ^\dag _{T _r} 
\circ \cdots \circ 
\R \underline{\Gamma} ^\dag _{T _1} (\smash{\widetilde{\B}} _{\fP} ^{(\bullet)} (T  ))$$
which are functorial in $T _i$ and $T$, i.e. making commutative the following diagram of 
$D ^{\mathrm{b}} _{\mathrm{coh}}( \smash{\D} ^\dag _{\fP,\Q} )$
\begin{gather}
\notag
\xymatrix{
{\R \sp _* 
\left (
\underline{\Gamma} ^\dag _{T _r} 
\circ \cdots \circ 
\underline{\Gamma} ^\dag _{T _1} ( \O _{\fP _K})
\right )} 
\ar[r] ^-{\sim}
\ar[d] ^-{}
& 
{\underrightarrow{\lim}\,
\R \underline{\Gamma} ^\dag _{T _r} 
\circ \cdots \circ 
\R \underline{\Gamma} ^\dag _{T _1} (\O _{\fP} ^{(\bullet)} )} 
\ar[d] ^-{}
\\ 
{\R \sp _* 
\left (
\underline{\Gamma} ^\dag _{T _r} 
\circ \cdots \circ 
\underline{\Gamma} ^\dag _{T _1} ( j ^\dag _{T}\O _{\fP _K})
\right )} 
\ar[r] ^-{\sim}
\ar[d] ^-{}
& 
{\underrightarrow{\lim}\,
\R \underline{\Gamma} ^\dag _{T _r} 
\circ \cdots \circ 
\R \underline{\Gamma} ^\dag _{T _1} (\smash{\widetilde{\B}} _{\fP} ^{(\bullet)} (T  ))} 
\ar[d] ^-{}
\\
{\R \sp _* 
\left (
\underline{\Gamma} ^\dag _{T _{r-1}} 
\circ \cdots \circ 
\underline{\Gamma} ^\dag _{T _1} ( j ^\dag _{T}\O _{\fP _K})
\right )} 
\ar[r] ^-{\sim}
& 
{\underrightarrow{\lim}\,
\R \underline{\Gamma} ^\dag _{T _{r-1}} 
\circ \cdots \circ 
\R \underline{\Gamma} ^\dag _{T _1} (\smash{\widetilde{\B}} _{\fP} ^{(\bullet)} (T  )),} 
}
\end{gather}
where the vertical arrows are the canonical ones induced by 
$\O _{\fP _K} \to j ^\dag _{T}\O _{\fP _K}$,
$\O _{\fP} ^{(\bullet)} 
\to \smash{\widetilde{\B}} _{\fP} ^{(\bullet)} (T  )$,
$\underline{\Gamma} ^\dag _{T _r}  \to id$, 
$\R \underline{\Gamma} ^\dag _{T _r} \to id$, 
and where $\underrightarrow{\lim} $
is the equivalence of categories
$\underrightarrow{\lim} 
\colon
\underrightarrow{LD} ^{\mathrm{b}} _{\Q, \mathrm{coh}} (\smash{\widehat{\D}} _{\fP /\fS } ^{(\bullet)})
\cong
D ^{\mathrm{b}} _{\mathrm{coh}}( \smash{\D} ^\dag _{\fP,\Q} )$
(see \ref{eqcat-limcoh}).
\end{cor}

\begin{proof}
Thanks to \ref{coh-cst-div}, 
we can copy the proof of \cite[9.1.11]{caro-6operations}.
\end{proof}

\subsection{Local cohomological functor with strict support over closed subscheme
for overconvergent complexes}
Let  $\fP $ be a 
formal  $\fS$-scheme
of formal finite type
and having locally finite $p$-bases  over $\fS$.

\begin{dfn}
[Overcoherent complexes]
\label{dfn-ovch}
Let 
$\cE ^{(\bullet)} 
\in 
\smash{\underrightarrow{LD}} ^{\mathrm{b}} _{\Q,\mathrm{coh}} 
(\overset{^\mathrm{l}}{} \smash{\widehat{\D}} _{\fP /\fS  } ^{(\bullet)} )$.
We say that $\cE ^{(\bullet)} $ is overcoherent 
if for any smooth formal $\fS$-scheme $\fX$, 
for any divisor $T$ of $P\times _S X$, denoting 
by $\varpi \colon \fP \times _{\fS} \fX \to \fP$ the projection, 
we have
$$(\hdag T)\circ \varpi ^{ !(\bullet)} (\cE ^{(\bullet)}  ) 
\in 
\smash{\underrightarrow{LD}} ^{\mathrm{b}} _{\Q,\mathrm{coh}} 
(\overset{^\mathrm{l}}{} \smash{\widehat{\D}} _{\fP \times _{\fS} \fX /\fS  } ^{(\bullet)} ).$$
We denote by 
$\smash{\underrightarrow{LD}} ^{\mathrm{b}} _{\Q,\mathrm{ovcoh}} 
(\overset{^\mathrm{l}}{} \smash{\widehat{\D}} _{\fP /\fS  } ^{(\bullet)} )$ 
the full subcategory of 
$\smash{\underrightarrow{LD}} ^{\mathrm{b}} _{\Q,\mathrm{coh}} 
(\overset{^\mathrm{l}}{} \smash{\widehat{\D}} _{\fP /\fS  } ^{(\bullet)} )$
consisting of 
overcoherent complexes.
This notion is an analogue of 
that of overcoherence as defined in
\cite[5.4]{caro-stab-sys-ind-surcoh} (see also \ref{pre-loc-tri-B-t1Tovcoh-cor2}).
Moreover, if $T$ is a divisor of $P$, with the notation of chapter \ref{ntn-tildeD(Z)}, 
for simplicity we set 
$\smash{\underrightarrow{LD}} ^{\mathrm{b}} _{\Q,\mathrm{ovcoh}} 
(\overset{^\mathrm{l}}{} \smash{\widetilde{\D}} _{\fP /\fS  } ^{(\bullet)} (T)):=
\smash{\underrightarrow{LD}} ^{\mathrm{b}} _{\Q,\mathrm{ovcoh}} 
(\overset{^\mathrm{l}}{} \smash{\widehat{\D}} _{\fP /\fS  } ^{(\bullet)})
\cap 
\smash{\underrightarrow{LD}} ^{\mathrm{b}} _{\Q,\mathrm{coh}} 
(\overset{^\mathrm{l}}{} \smash{\widetilde{\D}} _{\fP /\fS  } ^{(\bullet)} (T))$.

\end{dfn}

\begin{ex}
\label{O-ovcoh}
Suppose $\fP$ is of finite type over $\bbD ^1 _{\fS}$. 
It follows from \ref{coh-cst-div} that
$\cO _{\fP} ^{(\bullet) } \in \smash{\underrightarrow{LD}} ^{\mathrm{b}} _{\Q,\mathrm{ovcoh}} 
(\overset{^\mathrm{l}}{} \smash{\widehat{\D}} _{\fP /\fS  } ^{(\bullet)} )$.
More generally, this is not clear 
that 
$\cO _{\fP} ^{(\bullet) } $
is overcoherent
(see \ref{warning}).

\end{ex}

\begin{lemm}
\label{induction-div-coh}
Let $T _1, \dots, T _{r}$ be some divisors of $P$.
Let $T$ be a divisor of $P$.
Let 
$\cE  ^{(\bullet) } \in \smash{\underrightarrow{LD}} ^{\mathrm{b}} _{\Q,\mathrm{ovcoh}} 
(\overset{^\mathrm{l}}{} \smash{\widehat{\D}} _{\fP /\fS  } ^{(\bullet)} )$.
Then
$\R \underline{\Gamma} ^\dag _{T _r} \circ \cdots \circ 
\R \underline{\Gamma} ^\dag _{T _1} \circ (\hdag T  ) (\cE  ^{(\bullet) }) 
\in 
\smash{\underrightarrow{LD}} ^{\mathrm{b}} _{\Q,\mathrm{ovcoh}} ( \smash{\widehat{\D}} _{\fP /\fS } ^{(\bullet)} )$.
\end{lemm}

\begin{proof}
This can be easily checked by devissage.
\end{proof}

\begin{dfn}
\label{dfn-4.3.4}
Let $X$ be a closed subscheme of $P$.
Similarly to \cite[2.2]{caro_surcoherent}, we define
the local cohomological functor  
$\R \underline{\Gamma} ^\dag _{X}
\colon 
\smash{\underrightarrow{LD}} ^{\mathrm{b}} _{\Q,\mathrm{ovcoh}} 
(\overset{^\mathrm{l}}{} \smash{\widehat{\D}} _{\fP /\fS  } ^{(\bullet)} )
\to 
\smash{\underrightarrow{LD}} ^{\mathrm{b}} _{\Q,\mathrm{ovcoh}} 
(\overset{^\mathrm{l}}{} \smash{\widehat{\D}} _{\fP /\fS  } ^{(\bullet)} )$
with strict support in $X$ as follows.
Since $P$ is the sum of its irreducible components $U _i$,
then we reduce to the case where $P$ is integral.

\begin{enumerate}[(a)]

\item When $X= P$, the functor $\R \underline{\Gamma} ^\dag _{X}$ 
is by definition the identity.

\item Suppose now $X \not = P$. 
Similarly to \cite[2.2.5]{caro_surcoherent} (there was a typo: we need to add the hypothesis {`` $P$ is integral''})
the underlying space of $X$ is equal to a finite intersection of (the support of some) divisors of $P$.
Choose some divisors
$T _1, \dots, T _{r}$ of $P$ such that
$X = \cap _{i =1} ^{r} T _i$.
For $\E ^{(\bullet)}
\in \smash{\underrightarrow{LD}} ^\mathrm{b} _{\Q,\mathrm{qc}} ( \smash{\widehat{\D}} _{\fP /\fS } ^{(\bullet)} )$,
the complex
$\R \underline{\Gamma} ^\dag _{X} (\E ^{(\bullet)} ):=
\R \underline{\Gamma} ^\dag _{T _r} \circ \cdots \circ 
\R \underline{\Gamma} ^\dag _{T _1} (\E ^{(\bullet)} )$
does not depend canonically on the choice of the divisors  $T _1,\dots, T _r$ of $P$
satisfying $X = \cap _{i =1} ^{r} T _i$.
(Indeed, 
thanks to Lemmas \ref{lemme2.2.3} and \ref{induction-div-coh},
it is useless to add divisors containing $ X$.)

\end{enumerate}

\end{dfn}

\begin{prop}
\label{cor-induction-div-coh2}
Suppose there exists a morphism $\fP\to \bbD ^1 _{\fS}$ of  finite type.
Let $X$ be a closed subscheme of $P$ 
having locally finite $p$-bases  over $\Spec k$.
The complex
$\R \underline{\Gamma} ^\dag _X \O _{\fP,\Q} 
:=\R \sp _* \underline{\Gamma} ^\dag _X ( \O _{\fP _K})$
defined at \ref{ntn-GammaZO-rig} 
is canonically isomorphic to 
$\underrightarrow{\lim}\,
\R \underline{\Gamma} ^\dag _{X} 
 (\O _{\fP} ^{(\bullet)} )$, which confirms the compatibility of our notation.

\end{prop}

\begin{proof}
By the construction explained in \ref{dfn-4.3.4}, 
this is a consequence of \ref{cor-induction-div-coh}.\end{proof}

\begin{lem}
\label{lem-u*XT=0}
Let $u \colon \fX \hookrightarrow \fP$ be a closed immersion
of formal  $\fS$-schemes
of formal finite type
and having locally finite $p$-bases  over $\fS$.
Let $T$ be a divisor of $P$ containing $X$.
Then for any 
$\E ^{(\bullet)} 
\in 
\smash{\underrightarrow{LD}} ^{\mathrm{b}} _{\Q,\mathrm{qc}}(\overset{^\mathrm{l}}{} \smash{\widehat{\D}} _{\fP /\fS  } ^{(\bullet)} )$,
we have the isomorphism 
$u ^{! (\bullet)} \circ (\hdag T) ( \cE ^{(\bullet)} )\riso 0$
in $\smash{\underrightarrow{LD}} ^{\mathrm{b}} _{\Q,\mathrm{qc}}(\overset{^\mathrm{l}}{} \smash{\widehat{\D}} _{\fP /\fS  } ^{(\bullet)} )$.
\end{lem}

\begin{proof}
Following \ref{iso-comm-locaux-prod-tens}, 
we reduce to check 
$u ^{ !(\bullet)} ((\hdag T) (\cO _{\fP} ^{(\bullet)}))= 0$.
Since this is local, we can suppose $\fP$ affine and there exist $f \in \Gamma (\fP, \cO _{\fP})$ such that
$T = V ( \overline{f})$, where $\overline{f}$ is the image of $f$ in $\Gamma (P, \cO _{P})$.
Then 
$\widehat{\B} ^{(m)} _{\fP} (T)$ is the $p$-adic completion of 
$ \O _{\fP} [t] / (f ^{p ^{m +1}} T-p)$, denoted by 
$ \O _{\fP} \{ t \} / (f ^{p ^{m +1}} t-p)$.
Hence, 
$u ^* (\widehat{\B} ^{(m)} _{\fP} (T) )
=
 \O _{\fX} \{ t \} / (p)
 =
 \O _{\fX} [t ] / (p)$.
Hence, 
$p u ^*( \widehat{\B} ^{(m)} _{\fP} (T)) =0$.
This yields,
$u ^* (\widehat{\B} ^{(\bullet)} _{\fP} (T)) \riso 0$.
Since, 
$u ^* (\widehat{\B} ^{(\bullet)} _{\fP} (T)) [\delta _{\fX/\fP}]
\riso 
u ^{!(\bullet)} \circ (\hdag T) ( \cO _{\fP }^{(\bullet)})$, 
we are done.
\end{proof}

\begin{prop}
\label{pre-loc-tri-B-t1Tovcoh}
Let $u \colon \fX  \to \fP $ be a closed immersion of 
 formal  $\fS$-schemes of formal finite type
and having locally finite $p$-bases  over $\fS$.
Let 
$\E ^{(\bullet)} 
\in 
\smash{\underrightarrow{LD}} ^{\mathrm{b}} _{\Q,\mathrm{ovcoh}}(\overset{^\mathrm{l}}{} \smash{\widehat{\D}} _{\fP /\fS  } ^{(\bullet)} )$.
\begin{enumerate}[(a)]
\item Then
$u ^{ !(\bullet)} (\E ^{(\bullet)})
\in 
\smash{\underrightarrow{LD}} ^{\mathrm{b}} _{\Q,\mathrm{ovcoh}}(\overset{^\mathrm{l}}{} \smash{\widehat{\D}} _{\fX /\fS  } ^{(\bullet)} )$

\item and we have the canonical isomorphism
of 
$\smash{\underrightarrow{LD}} ^{\mathrm{b}} _{\Q,\mathrm{ovcoh}}(\overset{^\mathrm{l}}{} \smash{\widehat{\D}} _{\fP /\fS  } ^{(\bullet)} )$ 
of the form
\begin{equation}
\label{pre-loc-tri-B-t1T-isoovcoh}
\R \underline{\Gamma} ^\dag _{X} (\E ^{(\bullet)}) 
\riso 
u _{+} ^{ (\bullet)} \circ  u ^{ !(\bullet)} (\E ^{(\bullet)}).
\end{equation}

\end{enumerate}
\end{prop}

\begin{proof}
I) First, suppose that $X$ is a divisor of $P$.
i) We prove that 
$u ^{ !(\bullet)} (\E ^{(\bullet)})
\in 
\smash{\underrightarrow{LD}} ^{\mathrm{b}} _{\Q,\mathrm{coh}}(\overset{^\mathrm{l}}{} \smash{\widehat{\D}} _{\fX /\fS  } ^{(\bullet)} )$.

Following 
\ref{tri-loc-berthelot} we get 
 the exact triangle of 
 $\smash{\underrightarrow{LD}} ^{\mathrm{b}} _{\Q,\mathrm{ovcoh}}(\overset{^\mathrm{l}}{} \smash{\widehat{\D}} _{\fP /\fS  } ^{(\bullet)} )$
\begin{equation}
\label{pre-loc-tri-B-t1Tovcoh-Tri1}
\R \underline{\Gamma} ^\dag _{X} (\E ^{(\bullet)})
\to 
\E ^{(\bullet)}
\to 
(\hdag X) (\E ^{(\bullet)})
\to 
 \R \underline{\Gamma} ^\dag _{X} (\E ^{(\bullet)})
 [1].
\end{equation}
It follows from \ref{lem-u*XT=0} that we have 
$u ^{ !(\bullet)} ((\hdag X) (\E ^{(\bullet)}))= 0$.
This yields that the canonical morphism
$u ^{ !(\bullet)} \circ \R \underline{\Gamma} ^\dag _{X} (\E ^{(\bullet)})
\to 
u ^{ !(\bullet)} (\E ^{(\bullet)})
$
is an isomorphism.
Since $\R \underline{\Gamma} ^\dag _{X} (\E ^{(\bullet)})$ is a coherent complex with support in $X$, 
then 
by using Berthelot-Kashiwara's theorem (\ref{u!u+=id-isobis}), 
$u ^{ !(\bullet)} \circ \R \underline{\Gamma} ^\dag _{X} (\E ^{(\bullet)})$
is coherent. 
Hence, we are done. 

ii) For any divisor $T _X$ of $X$, 
we prove in this step that 
$(\hdag T _X) \circ u ^{ !(\bullet)} (\E ^{(\bullet)})
\in 
\smash{\underrightarrow{LD}} ^{\mathrm{b}} _{\Q,\mathrm{coh}}(\overset{^\mathrm{l}}{} \smash{\widehat{\D}} _{\fX /\fS  } ^{(\bullet)} )$.
Since this is local, 
we can suppose there exists a divisor $T$ of $P$ such that 
$T \cap X = T _X$.
Following \ref{f!commoub},
we have 
$(\hdag T _X) \circ u ^{ !(\bullet)} (\E ^{(\bullet)})
\riso 
u ^{ !(\bullet)} \circ (\hdag T )  (\E ^{(\bullet)})$.
Since 
$\E ^{(\bullet)}$ is overcoherent, then 
$(\hdag T )  (\E ^{(\bullet)})$ is also overcoherent. 
Hence, it follows from I.i) that 
$u ^{ !(\bullet)} \circ (\hdag T )  (\E ^{(\bullet)})$ is coherent. 

iii) By using the base change isomorphism of the form \ref{theo-iso-chgtbase2}, 
it follows from  I.ii) that 
$u ^{ !(\bullet)} (\E ^{(\bullet)})
\in 
\smash{\underrightarrow{LD}} ^{\mathrm{b}} _{\Q,\mathrm{ovcoh}}(\overset{^\mathrm{l}}{} \smash{\widehat{\D}} _{\fX /\fS  } ^{(\bullet)} )$.

II) Let us go back to the general case. 
i) Let $\I$ be the ideal given by $u$. 
Since this is local, it follows from \ref{cor-closed-immer-local} that we can suppose
there exist
$t _{r +1},\dots , t _{d}\in \Gamma (\fP,\I)$ generating 
$\Gamma (\fP,\I)$,
$t _{1},\dots , t _{r}  \in \Gamma ( \fP,\O _{\fP})$
such that,
denoting by 
$\overline{t} _1, \dots, \overline{t} _d$ 
the image of 
$t _1, \dots, t _d$ on $\Gamma ( \fX,\O _{\fX})$,
 the following (the third one is useless here) properties hold :
\begin{enumerate}[(1)]
\item $t _{1},\dots ,t  _{d}$ form a finite $p$-basis of $\fP $ over $\fS$ ;
\item  $\overline{t} _{1},\dots ,\overline{t} _{r}$ form a finite $p$-basis of $\fX$ over $\fS$.
\end{enumerate}
By induction in $d -r$ from the part I), we can check that 
$u ^{ !(\bullet)} (\E ^{(\bullet)})
\in 
\smash{\underrightarrow{LD}} ^{\mathrm{b}} _{\Q,\mathrm{ovcoh}}(\overset{^\mathrm{l}}{} \smash{\widehat{\D}} _{\fX /\fS  } ^{(\bullet)} )$.

ii) 
It follows from \ref{lem-u*XT=0} that for any divisor $T$ containing $X$, the canonical morphism
$u ^{ !(\bullet)} \circ \R \underline{\Gamma} ^\dag _{T} (\E ^{(\bullet)})
\to 
u ^{ !(\bullet)} (\E ^{(\bullet)})
$
is an isomorphism.
Hence, 
the canonical morphism
$u ^{ !(\bullet)} \circ \R \underline{\Gamma} ^\dag _{X} (\E ^{(\bullet)})
\to 
u ^{ !(\bullet)} (\E ^{(\bullet)})
$
is an isomorphism.
This yields that the canonical morphism
$u ^{(\bullet)} _+ \circ u ^{ !(\bullet)} \circ \R \underline{\Gamma} ^\dag _{X} (\E ^{(\bullet)})
\to 
u ^{(\bullet)} _+ \circ  u ^{ !(\bullet)} (\E ^{(\bullet)})
$
is an isomorphism. 
Since $\R \underline{\Gamma} ^\dag _{X} (\E ^{(\bullet)})$ is a coherent complex with support in $X$, 
then 
by using Berthelot-Kashiwara's theorem (\ref{u!u+=id-isobis}), 
$u ^{(\bullet)} _+ \circ u ^{ !(\bullet)} \circ \R \underline{\Gamma} ^\dag _{X} (\E ^{(\bullet)})
\riso 
 \R \underline{\Gamma} ^\dag _{X} (\E ^{(\bullet)})$.
Hence, we get the isomorphism \ref{pre-loc-tri-B-t1T-isoovcoh}.
\end{proof}

\begin{cor}
\label{pre-loc-tri-B-t1Tovcoh-cor}
Let 
$\fX$ be a smooth formal $\fS$-scheme, 
$\fP'$ and $\fP$ be two formal  $\fS$-schemes of formal finite type
and having locally finite $p$-bases  over $\fS$.
Let $\fP' \hookrightarrow \fX \times _{\fS} \fP$ be an immersion
and $f \colon \fP '  \to \fP $ be the induced morphism by composition with the canonical projection.
For any 
$\E ^{(\bullet)} 
\in 
\smash{\underrightarrow{LD}} ^{\mathrm{b}} _{\Q,\mathrm{ovcoh}}(\overset{^\mathrm{l}}{} \smash{\widehat{\D}} _{\fP /\fS  } ^{(\bullet)} )$,
we have 
$f ^{ !(\bullet)} (\E ^{(\bullet)})
\in 
\smash{\underrightarrow{LD}} ^{\mathrm{b}} _{\Q,\mathrm{ovcoh}}(\overset{^\mathrm{l}}{} \smash{\widehat{\D}} _{\fP' /\fS  } ^{(\bullet)} )$.
\end{cor}

\begin{proof}
Let $\varpi \colon \fX \times _{\fS} \fP \to \fP$ be the canonical projection.
From the definition of overcoherence, since a product of smooth formal $\fS$-schemes is a smooth formal $\fS$-scheme, 
then 
$\varpi ^{ !(\bullet)} (\E ^{(\bullet)})
\in 
\smash{\underrightarrow{LD}} ^{\mathrm{b}} _{\Q,\mathrm{ovcoh}}(\overset{^\mathrm{l}}{} \smash{\widehat{\D}} _{\fX \times _{\fS} \fP} ^{(\bullet)} )$.
By using \ref{pre-loc-tri-B-t1Tovcoh}, we check that the overcoherence is stable under 
$u ^{ !(\bullet)}$ when $u$ is an immersion (indeed, the case of an open immersion is easy).
\end{proof}

\begin{cor}
\label{pre-loc-tri-B-t1Tovcoh-cor2}
Let $f \colon \fP'\to   \fP$ be a finite type morphism of 
formal  $\fS$-schemes of formal finite type
and having locally finite $p$-bases  over $\fS$.
For any 
$\E ^{(\bullet)} 
\in 
\smash{\underrightarrow{LD}} ^{\mathrm{b}} _{\Q,\mathrm{ovcoh}}(\overset{^\mathrm{l}}{} \smash{\widehat{\D}} _{\fP /\fS  } ^{(\bullet)} )$,
we have 
$f ^{ !(\bullet)} (\E ^{(\bullet)})
\in 
\smash{\underrightarrow{LD}} ^{\mathrm{b}} _{\Q,\mathrm{ovcoh}}(\overset{^\mathrm{l}}{} \smash{\widehat{\D}} _{\fP' /\fS  } ^{(\bullet)} )$.
\end{cor}

\begin{proof}
Since the overcoherence is local, 
this is a straightforward consequence of \ref{pre-loc-tri-B-t1Tovcoh-cor}.
\end{proof}

\begin{prop}
Let $X,~X'$ be two closed subschemes of  $P$, 
$\E ^{(\bullet)},~\FF ^{(\bullet)}
\in \smash{\underrightarrow{LD}} ^\mathrm{b} _{\Q,\mathrm{ovcoh}} ( \smash{\widehat{\D}} _{\fP /\fS } ^{(\bullet)})$.

\begin{enumerate}[(a)]
\item 
We have the canonical isomorphism functorial in $\E ^{(\bullet)},~ X$, and $X'$ :
\begin{equation}
\label{theo2.2.8}
\R \underline{\Gamma} ^\dag _{X} \circ \R \underline{\Gamma} ^\dag _{X'} 
(\E ^{(\bullet)} )
\riso
\R \underline{\Gamma} ^\dag _{X\cap X'}  (\E ^{(\bullet)} ).
\end{equation}

\item We have the canonical isomorphism functorial in $\E ^{(\bullet)},~\FF ^{(\bullet)},~ X$, and $X'$ :
\begin{equation}
\label{fonctXX'Gamma-iso}
\R \underline{\Gamma} ^\dag _{X \cap X'} (\E ^{(\bullet)}
\smash{\widehat{\otimes}}^\L  _{\O ^{(\bullet)}  _{\fP} }  \FF ^{(\bullet)} )
\riso 
\R \underline{\Gamma} ^\dag _{X} 
(\E ^{(\bullet)})
\smash{\widehat{\otimes}}^\L  _{\O ^{(\bullet)}  _{\fP} }  
\R \underline{\Gamma} ^\dag _{X'}
(\FF ^{(\bullet)}).
\end{equation}
\end{enumerate}

\end{prop}

\begin{proof}
The first statement is obvious by construction of the local cohomological functor with strict support.
We can copy \cite[4.3.6]{caro-stab-sys-ind-surcoh} for the last one.
\end{proof}

\subsection{Localisation outside a closed subscheme functor for overconvergent complexes}
Let  $\fP $ be a 
formal  $\fS$-scheme
of formal finite type
and having locally finite $p$-bases  over $\fS$.
\begin{dfn}
Let $\E ^{(\bullet)}
\in \smash{\underrightarrow{LD}} ^\mathrm{b} _{\Q,\mathrm{ovcoh}} ( \smash{\widehat{\D}} _{\fP /\fS } ^{(\bullet)})$.
Let $X$ be a closed subscheme of  $P$.
Using \cite[1.1.10]{BBD} and Lemma \ref{annulationHom-hdag}, 
we check that the cone of the morphism
$\R \underline{\Gamma} ^\dag _{X} (\E ^{(\bullet)} )
\to  
\E ^{(\bullet)} $
is unique up to canonical isomorphism (for more details, 
see \cite[4.4.3]{caro-stab-sys-ind-surcoh}).
We will denote it by 
$(\hdag X) (\E ^{(\bullet)} )$.
We check that 
$(\hdag X) (\E ^{(\bullet)} )$
is functorial in $X$, and $\E ^{(\bullet)} $.
We get by construction the distinguished triangle
\begin{equation}
\label{caro-stab-sys-ind-surcoh4.4.3}
\R \underline{\Gamma} ^\dag _{X} (\E ^{(\bullet)} )
\to 
\E ^{(\bullet)} 
\to 
(\hdag X) (\E ^{(\bullet)} )
 \to 
\R \underline{\Gamma} ^\dag _{X} (\E ^{(\bullet)} )[1] .
\end{equation}
\end{dfn}

\begin{empt}
\label{coh-OGammaY}
Let  $X, X'$ be two closed subschemes of $P$.

\begin{enumerate}[(a)]
\item 
For any 
$\cE ^{(\bullet)} 
\in 
\smash{\underrightarrow{LD}} ^{\mathrm{b}}  
_{\Q, \mathrm{ovcoh}}
(\overset{^\mathrm{l}}{} \smash{\widehat{\D}} _{\fP /\fS  } ^{(\bullet)} )$, 
we have 
$(\hdag X ') \circ \R \underline{\Gamma} ^\dag _{X}  (\cE ^{(\bullet)}  )
\in 
\smash{\underrightarrow{LD}} ^{\mathrm{b}}  
_{\Q, \mathrm{ovcoh}}
(\overset{^\mathrm{l}}{} \smash{\widehat{\D}} _{\fP /\fS  } ^{(\bullet)} )$.

\item Suppose there exist a morphism $\fP\to \bbD ^1 _{\fS}$ of  finite type.
Then following, \ref{O-ovcoh}, we get 
$(\hdag X ') \circ \R \underline{\Gamma} ^\dag _{X}  (\O _{\fP}^{(\bullet)}  )
\in 
\smash{\underrightarrow{LD}} ^{\mathrm{b}}  
_{\Q, \mathrm{ovcoh}}
(\overset{^\mathrm{l}}{} \smash{\widehat{\D}} _{\fP /\fS  } ^{(\bullet)} )$.
\end{enumerate}

\end{empt}

\begin{empt}
\label{(hdagX)otimes}
For a closed subscheme  $X$ of $P$, 
for
$\E ^{(\bullet)},~\FF ^{(\bullet)}
\in \smash{\underrightarrow{LD}} ^\mathrm{b} _{\Q,\mathrm{ovcoh}} ( \smash{\widehat{\D}} _{\fP /\fS } ^{(\bullet)})$,
there exists a canonical isomorphism 
$(\hdag X) (\E ^{(\bullet)}
\smash{\widehat{\otimes}}^\L
_{\O ^{(\bullet)}  _{\fP}} 
\FF ^{(\bullet)}) 
\riso
\E ^{(\bullet)}
\smash{\widehat{\otimes}}^\L
_{\O ^{(\bullet)}  _{\fP}} 
(\hdag X) (\FF ^{(\bullet)}) $, 
which is moreover 
functorial in $X,~\E ^{(\bullet)},~\FF ^{(\bullet)}$ (for more details see 
\cite[4.4.4]{caro-stab-sys-ind-surcoh}).
\end{empt}

\begin{empt}
Let  $X, X'$ be two closed subschemes of $P$,
$\E ^{(\bullet)}
\in \smash{\underrightarrow{LD}} ^\mathrm{b} _{\Q,\mathrm{ovcoh}} ( \smash{\widehat{\D}} _{\fP /\fS } ^{(\bullet)})$.
There exists a canonical isomorphism 
$(\hdag X') \circ \R \underline{\Gamma} ^\dag _{X}(\E ^{(\bullet)} )
\riso
\R \underline{\Gamma} ^\dag _{X} \circ (\hdag X') (\E ^{(\bullet)})$
functorial in  $X,~X ',~\E ^{(\bullet)}$.

Similarly to  \cite[2.2.14]{caro_surcoherent}, 
we get the canonical isomorphism
\begin{equation}
\label{2.2.14-surcoh}
(\hdag X) \circ (\hdag X') (\E ^{(\bullet)})
\riso
(\hdag X \cup X')(\E ^{(\bullet)}),
\end{equation}
functorial in $X,~X ',~\E ^{(\bullet)}$.
Similarly to \cite[2.2.16]{caro_surcoherent},
we get the Mayer-Vietoris distinguished triangles :
\begin{gather}
\notag
  \R \underline{\Gamma} ^\dag _{X \cap X'}(\E ^{(\bullet)}) \rightarrow
  \R \underline{\Gamma} ^\dag _{X }(\E ^{(\bullet)}) \oplus
\R \underline{\Gamma} ^\dag _{X ' }(\E ^{(\bullet)})  \rightarrow
\R \underline{\Gamma} ^\dag _{X \cup X '}(\E^{(\bullet)} ) \rightarrow
\R \underline{\Gamma} ^\dag _{X \cap X'}(\E^{(\bullet)} )[1],\\
\label{eq1mayer-vietoris}
 (\hdag X \cap X')(\E ) \rightarrow  (\hdag X )(\E^{(\bullet)} ) \oplus
 (\hdag X ') (\E^{(\bullet)} )  \rightarrow   (\hdag X \cup X ')(\E^{(\bullet)} ) \rightarrow (\hdag X \cap X')(\E ^{(\bullet)})[1].
\end{gather}

\end{empt}

 \begin{prop}
\label{prop2.2.9}
Let  $D$ be a divisor of  $P$, 
$X$ be a closed subscheme of $P$,
$\U$ be the open subset of  $\fP$ complementary to the support of $X$.
Let $\E ^{(\bullet)}
\in 
\smash{\underrightarrow{LD}} ^{\mathrm{b}} _{\Q,\mathrm{ovcoh}} ( \smash{\widetilde{\D}} _{\fP /\fS } ^{(\bullet)} (D))$.
The following assertions are equivalent :
\begin{enumerate}[(a)]
\item We have in  $\smash{\underrightarrow{LD}} ^{\mathrm{b}} _{\Q,\mathrm{ovcoh}} ( \smash{\widetilde{\D}} _{\U /\fS } ^{(\bullet)} (D\cap U))$
the isomorphism $\E ^{(\bullet)}|\U \riso 0$.
\item The canonical morphism 
$\R \underline{\Gamma} ^\dag _{X} (\E ^{(\bullet)})
\to 
\E ^{(\bullet)}$ is an isomorphism in 
$\smash{\underrightarrow{LD}} ^{\mathrm{b}} _{\Q} ( \smash{\widetilde{\D}} _{\fP /\fS } ^{(\bullet)} (D))$.
\item We have in $\smash{\underrightarrow{LD}} ^{\mathrm{b}} _{\Q,\mathrm{ovcoh}} ( \smash{\widetilde{\D}} _{\fP /\fS } ^{(\bullet)} (D))$
the isomorphism $(\hdag X) (\E ^{(\bullet)} )\riso 0$.
\end{enumerate}

\end{prop}

\begin{proof}
Using Lemma \ref{lemme2.2.3}, 
we can copy the proof of 
\cite[4.4.6]{caro-stab-sys-ind-surcoh}.
\end{proof}

\begin{empt}
[Support]
\label{dfn-support}
Let  $D$ be a divisor of  $P$, 
$\E ^{(\bullet)}
\in 
\smash{\underrightarrow{LD}} ^{\mathrm{b}} _{\Q,\mathrm{ovcoh}} ( \smash{\widetilde{\D}} _{\fP /\fS } ^{(\bullet)} (D))$.
The support of $\E ^{(\bullet)} $ 
is by definition 
the biggest closed subscheme $X$ of $P$ such that 
$(\hdag X) (\E ^{(\bullet)} )\riso 0$ (one of the equivalent conditions of \ref{prop2.2.9}).

Remark if 
$\E ^{(\bullet)}
\in 
\smash{\underrightarrow{LM}} _{\Q,\mathrm{ovcoh}} ( \smash{\widetilde{\D}} _{\fP /\fS } ^{(\bullet)} (D))$, 
then 
this is equal to the  support (for the usual definition) of the coherent 
$\smash{\D} ^\dag _{\fP} (\hdag D) _{\Q} $-module
$\underrightarrow{\lim}\, \E ^{(\bullet)} $, which justifies the terminology.

\end{empt}

\subsection{Local cohomological functor with strict support over a subscheme for overconvergent complexes}
Let  $\fP $ be a 
formal  $\fS$-scheme
of formal finite type
and having locally finite $p$-bases  over $\fS$.

\begin{empt}
\label{3.2.1caro-2006-surcoh-surcv}
Let  $X$, $X'$, $T$, $T'$ be closed subschemes of  $P$ such that
  $X \setminus T = X' \setminus T'$.
For any $\E  ^{(\bullet)} \in \underrightarrow{LD} ^{\mathrm{b}} _{\Q ,\mathrm{ovcoh}}
(\smash{\widehat{\D}} _{\fP} ^{(\bullet)})$,
we have the canonical isomorphism:
\begin{equation}
  \label{xtx't'}
\R\underline{\Gamma} ^\dag _{X} (\hdag T ) (\E  ^{(\bullet)})
\riso
\R\underline{\Gamma} ^\dag _{X'} (\hdag T ') (\E  ^{(\bullet)}).
\end{equation}
Indeed, using \ref{coh-OGammaY},
\ref{theo2.2.8},
\ref{2.2.14-surcoh},
\ref{prop2.2.9},
we get the isomorphism
$\R\underline{\Gamma} ^\dag _{X} (\hdag T ) (\cE ^{(\bullet)})
\riso
\R\underline{\Gamma} ^\dag _{X \cap X'} (\hdag T \cup T ') (\cE   ^{(\bullet)})$.
We conclude by symmetry. 

Setting $ Y := X \setminus T$, we denote by 
$\R\underline{\Gamma} ^\dag _{Y} (\E  ^{(\bullet)}) $ one of both complexes
of \ref{xtx't'}.

\end{empt}

\begin{empt}
\label{propertiesGammaYY'}
Let  $Y $ and $Y'$ be two subschemes of  $P$. 
Let $\E  ^{(\bullet)} ,\FF ^{(\bullet)}\in \underrightarrow{LD} ^{\mathrm{b}} _{\Q ,\mathrm{ovcoh}}
(\smash{\widehat{\D}} _{\fP} ^{(\bullet)})$.
\begin{enumerate}[(a)]
\item 
Using \ref{theo2.2.8},
\ref{2.2.14-surcoh},
we get the canonical isomorphism functorial in $\E ^{(\bullet)},~ Y$, and $Y'$ :
\begin{equation}
  \label{gammayY'}
  \R\underline{\Gamma} ^\dag _{Y} \circ \R\underline{\Gamma} ^\dag _{Y'} (\E  ^{(\bullet)})
  \riso
  \R\underline{\Gamma} ^\dag _{Y \cap Y'} (\E  ^{(\bullet)}).
\end{equation}

\item  Using \ref{fonctXX'Gamma-iso} and \ref{(hdagX)otimes}
we get the canonical isomorphism functorial in $\E ^{(\bullet)},~\FF ^{(\bullet)},~ Y$, and $Y'$ :
\begin{equation}
\label{fonctYY'Gamma-iso}
\R \underline{\Gamma} ^\dag _{Y \cap Y'} (\E ^{(\bullet)}
\smash{\widehat{\otimes}}^\L  _{\O ^{(\bullet)}  _{\fP} }  \FF ^{(\bullet)} )
\riso 
\R \underline{\Gamma} ^\dag _{Y} 
(\E ^{(\bullet)})
\smash{\widehat{\otimes}}^\L  _{\O ^{(\bullet)}  _{\fP} }  
\R \underline{\Gamma} ^\dag _{Y'}
(\FF ^{(\bullet)}).
\end{equation}
\item  If $Y '$ is an open (resp. a closed) subscheme of  $Y$, we have the canonical homomorphism
$\R\underline{\Gamma} ^\dag _{Y} (\E  ^{(\bullet)}) \rightarrow \R\underline{\Gamma} ^\dag _{Y'} (\E  ^{(\bullet)})$
(resp. $\R\underline{\Gamma} ^\dag _{Y'} (\E  ^{(\bullet)}) \rightarrow \R\underline{\Gamma} ^\dag _{Y} (\E  ^{(\bullet)})$).
If $Y '$ is a closed subscheme of  $Y$, we have the localization distinguished triangle
$\R\underline{\Gamma} ^\dag _{Y'} (\E  ^{(\bullet)}) \rightarrow \R\underline{\Gamma} ^\dag _{Y} (\E  ^{(\bullet)})
\rightarrow \R\underline{\Gamma} ^\dag _{Y \setminus Y'} (\E  ^{(\bullet)}) \rightarrow +1.$

\end{enumerate}

\end{empt}

\subsection{Commutation with local cohomological functors for overconvergent complexes}

\begin{theo}
\label{2.2.18}
Let $f\colon \X '\to \X$ 
be a quasi-projective (in the sense of Definition \ref{projectivefscheme}) 
morphism of formal $\fS$-schemes of formal finite type and having locally finite $p$-bases.
Let $Y$ be a subscheme of $X$, $Y':= f ^{-1} (Y)$.
\begin{enumerate}[(a)]
\item 
\label{2.2.18-item1}
Let  $\E ^{(\bullet)} \in 
\smash{\underrightarrow{LD}} ^{\mathrm{b}} _{\Q,\mathrm{ovcoh}} 
(\overset{^\mathrm{l}}{} \smash{\widehat{\D}} _{\X /\fS  } ^{(\bullet)} )$. 
We have the functorial in $Y$ isomorphism of 
$\smash{\underrightarrow{LD}} ^{\mathrm{b}} _{\Q,\mathrm{ovcoh}} ( \smash{\widehat{\D}} _{\X '/\V} ^{(\bullet)})$:
\begin{gather}
\label{commutfonctcohlocal1}
  f ^{ !(\bullet)}  \circ\R \underline{\Gamma} ^\dag _{Y}(\E ^{(\bullet)}) 
  \riso
   \R \underline{\Gamma} ^\dag _{Y' }\circ f ^{ !(\bullet)}  (\E ^{(\bullet)}).
\end{gather}
   
\item 
\label{2.2.18-item2}
Let 
$\E ^{\prime (\bullet)} \in \smash{\underrightarrow{LD}} ^{\mathrm{b}} _{\Q,\mathrm{ovcoh}} ( \smash{\widehat{\D}} _{\X '/\V} ^{(\bullet)})$
with proper support over $X$ (see \ref{dfn-prop-support}).
Then the object 
$f _{+} (\E ^{\prime (\bullet)} ) $
belongs to 
$\smash{\underrightarrow{LD}} ^{\mathrm{b}} _{\Q,\mathrm{ovcoh}} ( \smash{\widehat{\D}} _{\X/\V} ^{(\bullet)})$.
Moreover, we have the functorial in $Y$ isomorphism:
\begin{gather}
\label{commutfonctcohlocal2}
\R \underline{\Gamma} ^\dag _{Y}\circ f ^{ (\bullet)} _{+} (\E ^{\prime (\bullet)})
\riso
f ^{ (\bullet)}  _{+} \circ \R \underline{\Gamma} ^\dag _{Y'}(\E ^{\prime (\bullet)}).
\end{gather}

\end{enumerate}

\end{theo}

\begin{proof}
a) Let us check \ref{commutfonctcohlocal1}.
Following \ref{pre-loc-tri-B-t1Tovcoh-cor}, 
the complexes are indeed overcoherent. 
By devissage and construction, we reduce to the case where $Y$ is the complement of a divisor $T$.
By definition, 
$f$ is the composition of an immersion of the form
$\fX' \hookrightarrow \widehat{\bbP} ^d \times _\fS \fX$ with the projection 
$\widehat{\bbP} ^d \times _\fS \fX \to \fX$. 
 Since the case where $f$ is a flat morphism is known (see \ref{f!commoub}),
we reduce to the case where $f$ is a closed immersion. 
We conclude by using again \ref{f!commoub} (indeed, either $T \cap X'$ is a divisor and we can use \ref{f!commoub}, or 
$T\cap X'= X'$ and then the isomorphism \ref{commutfonctcohlocal1} is $0\riso 0$). 
  
b) Let us check that $f _{+} (\E ^{\prime (\bullet)} ) $ is overcoherent. 
Let $\fZ$ be a smooth formal $\fS$-scheme, 
$\varpi \colon \fX \times _{\fS} \fZ \to \fX$
and
$\varpi '\colon \fX ' \times _{\fS} \fZ \to \fX'$
be the canonical projections.
Set $g = f \times id \colon  
\fX ' \times _{\fS} \fZ
\to 
\fX  \times _{\fS} \fZ$.
Let $T$ be a divisor of $\fX \times _{\fS} \fZ$.
We have to check that 
$(\hdag T) \circ  \varpi ^{ ! (\bullet)} f ^{ (\bullet)} _{+} (\E ^{\prime (\bullet)})$
is coherent. 
It follows from \ref{theo-iso-chgtbase2} that we have the isomorphism
$$(\hdag T) \circ  \varpi ^{ ! (\bullet)} f ^{ (\bullet)} _{+} (\E ^{\prime (\bullet)})
\riso 
(\hdag T) \circ   g ^{ (\bullet)} _{+}\circ \varpi ^{\prime ! (\bullet)}  (\E ^{\prime (\bullet)}).$$

i) First suppose $T ' = g ^{-1} (T)$ is a divisor of $X ' \times _{S} Z$.
It follows from 
\ref{surcoh2.1.4-cor}
that we have
$$(\hdag T) \circ   g ^{ (\bullet)} _{+}\circ \varpi ^{\prime ! (\bullet)}  (\E ^{\prime (\bullet)})
\riso 
  g ^{ (\bullet)} _{+} \circ (\hdag T')\circ \varpi ^{\prime ! (\bullet)}  (\E ^{\prime (\bullet)}).
$$
Since 
$\E ^{\prime (\bullet)}$
is overcoherent with proper support over $X$, 
then 
$(\hdag T')\circ \varpi ^{\prime ! (\bullet)}  (\E ^{\prime (\bullet)})$ is coherent with proper support over $X \times _S Z$. 
It follows from \ref{stab-propersupp} that 
$g ^{ (\bullet)} _{+} \circ (\hdag T')\circ \varpi ^{\prime ! (\bullet)}  (\E ^{\prime (\bullet)})$ is coherent. 

ii) In general, since $f$ is quasi-projective, then using part i) we reduce to the case $f$ is an immersion. 
Hence, we come down to treat two cases : either $g ^{-1}(T)$ is a divisor of $X ' \times _{S} Z$,
or $g ^{-1}(T)$ contains $X '\times _{S} Z$.
In the first case, we have already checked that $(\hdag T) \circ  \varpi ^{ ! (\bullet)} f ^{ (\bullet)} _{+} (\E ^{\prime (\bullet)})$
is coherent.
When $g ^{-1}(T)$ contains $X '\times _{S} Z$, 
since $\varpi ^{ ! (\bullet)} f ^{ (\bullet)} _{+} (\E ^{\prime (\bullet)})$ is coherent, 
then it follows from \ref{lemme2.2.3}
that 
$(\hdag T) \circ   g ^{ (\bullet)} _{+}\circ \varpi ^{\prime ! (\bullet)}  (\E ^{\prime (\bullet)})
=0$, which is coherent.

c) It remains to prove the isomorphism \ref{commutfonctcohlocal2}.
By devissage and construction, we reduce to the case where $Y$ is the complement of a divisor $T$.
Since the case where $f$ is smooth is already known 
(see \ref{surcoh2.1.4-cor}), it remains to check the case where $f$ is a closed immersion.
When $T$ contains $X$, then the isomorphism 
\ref{commutfonctcohlocal2} is $0 \riso 0$.
When $T \cap X$ is a divisor, this is \ref{surcoh2.1.4-cor}.
   Hence we are done.
\end{proof}

\begin{prop}
\label{bc-com-locfunc}
With notation \ref{subsect-comm-bc}, 
let $U$ be a subscheme of 
$X$ and $V : = \varpi ^{-1} (U)$ be the reduced subscheme of $Y$.
Let 
$\E ^{ (\bullet)}
\in  \smash{\underrightarrow{LD}} ^\mathrm{b} _{\Q, \mathrm{ovcoh}}
(\overset{^\mathrm{l}}{} \smash{\widehat{\D}} _{\X /\fS } ^{(\bullet)})$. 
We have the canonical isomorphism of $\smash{\underrightarrow{LD}} ^\mathrm{b} _{\Q, \mathrm{ovcoh}}
(\overset{^\mathrm{l}}{} \smash{\widehat{\D}} _{\fY /\fT } ^{(\bullet)})$
\begin{equation}
\label{bc-com-locfunc-iso}
\varpi  ^{* (\bullet) }( \R\underline{\Gamma} ^\dag _{U} (\E  ^{(\bullet)}) )
\riso 
\R\underline{\Gamma} ^\dag _{V} (\varpi  ^{* (\bullet) } (\E  ^{(\bullet)}) ).
\end{equation}

\end{prop}

\begin{proof}
By construction (see \ref{3.2.1caro-2006-surcoh-surcv}), we reduce to the case where there exists a divisor $D$ of $X$ such that 
$U = X \setminus D$. Since $\varpi$ is flat, then $E:= \varpi ^{-1} (D) $ is (the support of) a divisor of $Y$
such that 
$V= Y \setminus E$.
Since $\varpi ^{*}
(\widehat{\B} ^{(m)} _{\fX} ( D))
\riso 
\widehat{\B} ^{(m)} _{\fY} ( E)$, since the functor $\varpi ^*$ commutes with tensor products,
we are done.\end{proof}

\subsection{Local cohomological functors for quasi-coherent complexes over weak admissible subschemes}
Let  $\fP $ be a 
formal  $\fS$-scheme
of formal finite type
and having locally finite $p$-bases  over $\fS$.

\begin{empt}
\label{div-closed-sub}
Let $X$ be an integral closed subscheme of $P$ having locally $p$-bases over $S$.

We have the map from the set of (reduced) divisors of $P$ not containing $X$
to the set of (reduced) divisors of $X$ defined by $T \mapsto T \cap X$.
Recall that since $P$ and $X$ are regular (see \ref{regularity/formalsm}), then the notions of Cartier divisors 
or Weil divisors are similar. 
Then this map is ``locally surjective'' in the following sense.
Let $D$ be an integral divisor of $X$. Then 
there exists an open covering 
$(\fP _i ) _i$ of $\fP$ such that 
$ D \cap P _i $ is a principal divisor of $X$. 
Hence, $D \cap P _i$ is the intersection of a principal divisor of $P _i$ with $X$.

\end{empt}

\begin{empt}
\label{XpbasesGamma-desc}
Let $u \colon \fX \hookrightarrow \fP$ be a closed immersion
of formal  $\fS$-schemes
of formal finite type
and having locally finite $p$-bases  over $\fS$.
Let $\I$ be the ideal given by $u$. 

i) Following \ref{cor-closed-immer-local}, there exists a covering 
$(\fP _i) _{i =1, \dots , n}$ of $\fX$ by open affine subsets of $\fP$ such that 
there exist
$t _{i, r +1},\dots , t _{i, d}\in \Gamma (\fP _i,\I)$ generating 
$\Gamma (\fP _i,\I)$,
$t _{i, 1},\dots , t _{i, r}  \in \Gamma (\fP _i,\O _{\fP _i})$
such that,
denoting by 
$\overline{t} _{i,1}, \dots, \overline{t} _{i,d}$ 
the image of 
$t _{i,1}, \dots, t _{i,d}$ on $\Gamma ( \fX \cap \fP _i,\O _{\fX})$,
the following two (the third one is useless here) properties hold :
\begin{enumerate}[(a)]
\item $t _{i,1},\dots ,t  _{i, d}$ form a finite $p$-basis of $\fP _i$ over $\fS$ ;
\item  $\overline{t} _{i, 1},\dots ,\overline{t} _{i,r}$ form a finite $p$-basis of $\fX \cap \fP _i$ over $\fS$.
\end{enumerate}
Let $T _{i,j}$ be the divisor of $P$ equal to the closure in $P$ of the principal
divisor defined by $\overline{t} _{i,j}$ for $i = 1,\dots , n$ and $j = r +1, \dots d$. 
Since $X$ is irreducible, then $X \subset T _{i,j}$ for any $i = 1,\dots , n$ and $j = r +1, \dots d$.
Hence, $X \subset  \cap _{i,j}T _{i,j}$ (the intersection is over $i = 1,\dots , n$ and $j = r +1, \dots d$). 
Conversely, for any integer $i _0 \in \{ 1, \dots, n\}$,
we have the inclusion $ P _{i _0} \cap ( \cap _{i,j}T _{i,j} ) \subset P _{i _0} \cap ( \cap _{j}T _{i _0,j} ) = 
\cap _{j} ( P _{i _0} \cap T _{i _0,j} ) 
=P _{i _0} \cap X $ (recall $t _{i _{0}, r +1},\dots , t _{i _{0}, d}\in \Gamma (\fP _i,\I)$ generate 
$\Gamma (\fP _{i _{0}},\I)$).
Hence, 
$X = \cap _{i,j}T _{i,j}$.
Set 
$\R \underline{\Gamma} ^\dag _{X _i} :=
\R \underline{\Gamma} ^\dag _{T _{i,1}} \circ \cdots \circ 
\R \underline{\Gamma} ^\dag _{T _{i,d}}$.
Consider the functor
$\phi :=
\R \underline{\Gamma} ^\dag _{X _{1}} \circ \cdots \circ 
\R \underline{\Gamma} ^\dag _{X _{n}} $.

ii) Since $X \subset T _{i,j}$, then following \ref{lem-u*XT=0} we have
$u  ^{! (\bullet )}  (
(\hdag T _{i,j}) (\cO _{\fP } ^{(\bullet)}))
=0$. Hence, the canonical morphism
$u  ^{! (\bullet )}  \circ \R \underline{\Gamma} ^\dag _{T _{i,1}} 
(\cO _{\fP } ^{(\bullet)})
\to
u  ^{! (\bullet )}   (\cO _{\fP } ^{(\bullet)})$
is an isomorphism.
This yields the canonical isomorphism
$u  ^{! (\bullet )}  (\phi (\cO _{\fP } ^{(\bullet)}))
\riso 
u  ^{! (\bullet )}   (\cO _{\fP } ^{(\bullet)})$.
Hence, 
\begin{equation}
\label{XpbasesGamma-desc-iso0}
u _{+} ^{(\bullet)} \circ u  ^{! (\bullet )}  (\phi (\cO _{\fP } ^{(\bullet)}))
\riso 
u _{+} ^{(\bullet)} \circ  u  ^{! (\bullet )}   (\cO _{\fP } ^{(\bullet)}).
\end{equation}

iii) Let $u _i \colon \fX \cap \fP _i \hookrightarrow \fP _i$ be the closed immersion induced by $u$. 
In this step, we check that we have 
\begin{equation}
\label{XpbasesGamma-desc-iso1}
\R \underline{\Gamma} ^\dag _{X _i} (\cO _{\fP} ^{(\bullet)}) | \fP _i 
\riso 
u _{i +} ^{ (\bullet )} \circ u _i ^{! (\bullet )} (\cO _{\fP _i } ^{(\bullet)}).
\end{equation}
Let $u _{i,d} \colon V (t _{i,d}) \hookrightarrow \fP _i$ be the canonical closed immersion. 
Since 
$V (t _{i,d})=T _{i, d} \cap P _i$ is a smooth divisor of $P _i$ then 
$\cO _{\fP _i } ^{(\bullet)}  (\hdag T _{i, d} \cap P _i)
\in 
\smash{\underrightarrow{LD}} ^{\mathrm{b}} _{\Q,\mathrm{coh}} 
(\overset{^\mathrm{l}}{} \smash{\widehat{\D}} _{\fP /\fS  } ^{(\bullet)} )$.
Thanks to the inductive system version of Berthelot-Kashiwara's theorem
(see \ref{u!u+=id}) and since 
$u _{i,d} ^{! (\bullet)} (\cO _{\fP _i } ^{(\bullet)}  (\hdag T _{i, d} \cap P _i))=0$,
then 
\begin{equation}
\label{XpbasesGamma-desc-iso2}
\R \underline{\Gamma} ^\dag _{T _{i,d}\cap P _i} (\cO _{\fP _i } ^{(\bullet)}) 
\riso 
u _{i,d +} ^{ (\bullet )} \circ u _{i,d} ^{! (\bullet )} (\cO _{\fP _i } ^{(\bullet)}).
\end{equation}
Let $v _{i,d} \colon \fX \cap \fP _i \hookrightarrow V (t _{i,d}) $ be the canonical closed immersion. 
Set
$\R \underline{\Gamma} ^\dag _{X ' _i} :=
\R \underline{\Gamma} ^\dag _{T _{i,1} \cap P _i} \circ \cdots \circ 
\R \underline{\Gamma} ^\dag _{T _{i,d-1} \cap P _i}$
and
$\R \underline{\Gamma} ^\dag _{X '' _i} :=
\R \underline{\Gamma} ^\dag _{T _{i,1}\cap V (t _{i,d})} \circ \cdots \circ 
\R \underline{\Gamma} ^\dag _{T _{i,d-1}\cap V (t _{i,d})}$.
It follows from \ref{surcoh2.1.4-cor} that we have
$u _{i,d +} ^{ (\bullet )} \R \underline{\Gamma} ^\dag _{X ''_i} 
\riso 
\R \underline{\Gamma} ^\dag _{X ' _i} u _{i,d +} ^{ (\bullet )}  $.
Proceeding by  induction in $d -r$, 
we get 
\begin{equation}
\label{XpbasesGamma-desc-iso2-bis}
\R \underline{\Gamma} ^\dag _{X ''  _i} (\cO _{ V (t _{i,d})} ^{(\bullet)}) 
\riso 
v _{i,d +} ^{ (\bullet )} \circ v _{i,d} ^{! (\bullet )} (\cO _{  V (t _{i,d})} ^{(\bullet)}).
\end{equation}
Hence, we get
\begin{gather}
\notag
\R \underline{\Gamma} ^\dag _{X _i} (\cO _{\fP} ^{(\bullet)}) | \fP _i 
\riso 
\R \underline{\Gamma} ^\dag _{X ' _i} 
\R \underline{\Gamma} ^\dag _{T _{i,d}\cap P _i} (\cO _{\fP _i } ^{(\bullet)}) 
\underset{\ref{XpbasesGamma-desc-iso2}}{\riso} 
\R \underline{\Gamma} ^\dag _{X ' _i} u _{i,d +} ^{ (\bullet )}  
(\cO _{ V (t _{i,d})} ^{(\bullet)}) [-1]
\\
\notag 
\underset{\ref{surcoh2.1.4-cor}}{\riso} 
u _{i,d +} ^{ (\bullet )} \R \underline{\Gamma} ^\dag _{X ''_i} (\cO _{ V (t _{i,d})} ^{(\bullet)}) 
[-1]
\underset{\ref{XpbasesGamma-desc-iso2-bis}}{\riso} 
u _{i,d+} ^{ (\bullet )} \circ v _{i,d +} ^{ (\bullet )}  \circ v _{i,d} ^{! (\bullet )} (\cO _{  V (t _{i,d})} ^{(\bullet)})
[-1]
\riso 
u _{ i +} ^{ (\bullet )} u _{i} ^{! (\bullet )} (\cO _{\fP _i } ^{(\bullet)})  
.
\end{gather}

iv) Since $u _{i +} ^{ (\bullet )} \circ u _i ^{! (\bullet )} (\cO _{\fP _i } ^{(\bullet)})$
is coherent complex with support in $X \cap P _i$,
since $X \subset T _{i,j}$ 
for any $i = 1,\dots , n$ and $j = r +1, \dots d$,
then by using 
the inductive system version of Berthelot-Kashiwara's theorem
and iii), we get 
\begin{equation}
\label{XpbasesGamma-desc-iso3}
( \phi (\cO _{\fP} ^{(\bullet)}) )| \fP _i
\riso 
u _{i +} ^{ (\bullet )} \circ u _i ^{! (\bullet )} (\cO _{\fP _i } ^{(\bullet)}).
\end{equation}

v) It follows from \ref{XpbasesGamma-desc-iso3} that 
$\phi (\cO _{\fP} ^{(\bullet)})$ is coherent with support in $X$. 
Hence by using the
inductive system version of Berthelot-Kashiwara's theorem,
we get the canonical isomorphism
\begin{equation}
\label{XpbasesGamma-desc-iso4}
u _{+} ^{(\bullet)} \circ u  ^{! (\bullet )}  (\phi (\cO _{\fP } ^{(\bullet)}))
\riso 
\phi (\cO _{\fP } ^{(\bullet)}).
\end{equation}
Hence, from \ref{XpbasesGamma-desc-iso0} and 
\ref{XpbasesGamma-desc-iso4} we get by composition the canonical isomorphism
\begin{equation}
\label{XpbasesGamma-desc-iso5}
u _{+} ^{(\bullet)} \circ u  ^{! (\bullet )}  (\cO _{\fP } ^{(\bullet)})
\riso 
\phi (\cO _{\fP } ^{(\bullet)}).
\end{equation}

\end{empt}

\begin{dfn}
\label{dfn-4.3.4bis}
Let $X$ be a closed subscheme of $P$.

\begin{enumerate}[(a)]
\item Suppose $P$ integral.
\begin{enumerate}[(a)]
\item We say that $(P \subset \fP)$ is weak admissible, 
if $\cO _{\fP} ^{(\bullet)} \in 
\smash{\underrightarrow{LD}} ^{\mathrm{b}} _{\Q,\mathrm{ovcoh}} 
(\overset{^\mathrm{l}}{} \smash{\widehat{\D}} _{\fP /\fS  } ^{(\bullet)} )$
and then we put
$\R \underline{\Gamma} ^\dag _{P} ( \cO _{\fP} ^{(\bullet)})
:=
\cO _{\fP} ^{(\bullet)}$.

\item Suppose $X \not = P$.
We say that $(X \subset \fP)$ is weak admissible
if there exist some divisors 
$T _1, \dots, T _{r}$ of $P$ such that
 1) $X = \cap _{i =1} ^{r} T _i$
and such that 2) 
$\R \underline{\Gamma} ^\dag _{T _r} \circ \cdots \circ 
\R \underline{\Gamma} ^\dag _{T _1} (\cO _{\fP} ^{(\bullet)} )
\in 
\smash{\underrightarrow{LD}} ^{\mathrm{b}} _{\Q,\mathrm{ovcoh}} 
(\overset{^\mathrm{l}}{} \smash{\widehat{\D}} _{\fP /\fS  } ^{(\bullet)} )$.
Then we put 
$\R \underline{\Gamma} ^\dag _{X}  (\cO _{\fP} ^{(\bullet)} ):=
\R \underline{\Gamma} ^\dag _{T _r} \circ \cdots \circ 
\R \underline{\Gamma} ^\dag _{T _1} (\cO _{\fP} ^{(\bullet)} )$.
We remark that this does not depend on the choice of the divisors 
satisfying properties 1) and 2) above 
(Indeed, 
using Lemmas \ref{lemme2.2.3} and \ref{induction-div-coh},
it is useless to add divisors containing $ X$),
which justifies the notation.

\end{enumerate}

\item In general, $\fP$ is the sum of its irreducible components $\fP _i$.
We say that 
$(X \subset \fP)$ is weak admissible
if $(X \cap P _i, \fP _i)$ is weak admissible 
for any $i$.
In that case 
$\R \underline{\Gamma} ^\dag _{X} ( \cO _{\fP} ^{(\bullet)})$ is the object of 
$\smash{\underrightarrow{LD}} ^{\mathrm{b}} _{\Q,\mathrm{ovcoh}} 
(\overset{^\mathrm{l}}{} \smash{\widehat{\D}} _{\fP /\fS  } ^{(\bullet)} )$
so that
$\R \underline{\Gamma} ^\dag _{X} ( \cO _{\fP} ^{(\bullet)}) | \fP _i = 
\R \underline{\Gamma} ^\dag _{X \cap P _i} ( \cO _{\fP _i} ^{(\bullet)})$.

\end{enumerate}
\end{dfn}

\begin{rem}
\label{GammaannulationHom-hdag2-rem}
Let $X$ be a closed subscheme of $P$.
Contrary to the standard context of arithmetic $\cD$-modules,
this is not clear that $(X \subset \fP)$ is always weak admissible.
Suppose $(X \subset \fP)$ is weak admissible. We have the following further problems, which explains why we need to be careful. 
\begin{enumerate}[(a)]
\item Let $T _1, \dots, T _{r}$ be some divisors of $P$ such that
$X = \cap _{i =1} ^{r} T _i$.
Then this is not clear that 
$\R \underline{\Gamma} ^\dag _{T _r} \circ \cdots \circ 
\R \underline{\Gamma} ^\dag _{T _1} (\cO _{\fP} ^{(\bullet)} )
\in 
\smash{\underrightarrow{LD}} ^{\mathrm{b}} _{\Q,\mathrm{ovcoh}} 
(\overset{^\mathrm{l}}{} \smash{\widehat{\D}} _{\fP /\fS  } ^{(\bullet)} )$.

\item Let  $\fP'$ be a 
formal  $\fS$-scheme
of formal finite type
and having locally finite $p$-bases  over $\fS$.
Let $X'$ be a closed subscheme of $P'$ 
such that 
$X $ and $X'$ are isomorphic $S$-schemes. Then this is not  
clear that $(X ' \subset \fP')$ is weak admissible, 
even if $\fP$ is of finite type over $\bbD ^1 _\fS$.
This is an open question.

\end{enumerate}

\end{rem}

\begin{lemm}
\label{lemm-X'XGammadfn}
Let $X' \subset X$ be two closed subschemes of $P$.
If $(X '\subset \fP)$ is weak admissible
then we have the isomorphism 
\begin{equation}
\R \underline{\Gamma} ^\dag _{X } 
( 
\R \underline{\Gamma} ^\dag _{X'} ( \cO _{\fP} ^{(\bullet)})
)
\riso
\R \underline{\Gamma} ^\dag _{X'} ( \cO _{\fP} ^{(\bullet)}),
\end{equation}
where the functor 
$\R \underline{\Gamma} ^\dag _{X } $
is defined at \ref{dfn-4.3.4}.
\end{lemm}

\begin{proof}
We can suppose $P$ integral and $X ' \not = P$.
Let $T$ be a  divisor  containing $X$.
Since $\R \underline{\Gamma} ^\dag _{X'} ( \cO _{\fP} ^{(\bullet)})$ is a coherent complex with support in $X$,
since $T$ contains $X'$, since 
$(\hdag T) ( \R \underline{\Gamma} ^\dag _{X'} ( \cO _{\fP} ^{(\bullet)})) | (\fP \setminus T) =0$,
then it follows from \ref{Be1-4310&12}.\ref{Be1-4310&12-item3}
that 
$(\hdag T) ( \R \underline{\Gamma} ^\dag _{X'} ( \cO _{\fP} ^{(\bullet)})) =0$.
Hence, 
we have the canonical isomorphism
$
\R \underline{\Gamma} ^\dag _{T} ( 
\R \underline{\Gamma} ^\dag _{X'} ( \cO _{\fP} ^{(\bullet)}))
\riso
\R \underline{\Gamma} ^\dag _{X'} ( \cO _{\fP} ^{(\bullet)})$.
By definition of the functor $\R \underline{\Gamma} ^\dag _{X } $,
we conclude.\end{proof}

\begin{dfn}
\label{dfn-4.3.4bisY}
\begin{enumerate}[(a)]
\item
 \label{dfn-4.3.4bisY-item1}
Let $Y$ be a  subscheme of $P$.
We say that 
$(Y\subset \fP)$ is ``weak admissible''
if $(\overline{Y}\subset \fP)$ is weak admissible,
where 
$\overline{Y}$ is the closure of $Y$ in $P$.
In that case, 
we set
$$\R \underline{\Gamma} ^\dag _{Y} ( \cO _{\fP} ^{(\bullet)}):=
\R \underline{\Gamma} ^\dag _{Y} \left (
\R \underline{\Gamma} ^\dag _{\overline{Y}} ( \cO _{\fP} ^{(\bullet)})\right) \in \smash{\underrightarrow{LD}} ^{\mathrm{b}} _{\Q,\mathrm{ovcoh}} 
(\overset{^\mathrm{l}}{} \smash{\widehat{\D}} _{\fP /\fS  } ^{(\bullet)} ),$$
where 
$\R \underline{\Gamma} ^\dag _{Y}$ is the functor defined at \ref{3.2.1caro-2006-surcoh-surcv}
and 
$\R \underline{\Gamma} ^\dag _{\overline{Y}} ( \cO _{\fP} ^{(\bullet)})\in \smash{\underrightarrow{LD}} ^{\mathrm{b}} _{\Q,\mathrm{ovcoh}} 
(\overset{^\mathrm{l}}{} \smash{\widehat{\D}} _{\fP /\fS  } ^{(\bullet)} )$ is the object defined at 
\ref{dfn-4.3.4bis}.
It follows from \ref{lemm-X'XGammadfn} that this definition is compatible with 
\ref{dfn-4.3.4bis}.

\item 
\label{dfn-4.3.4bisY-item2}
Let $u \colon \fQ \hookrightarrow \fP$ be an (closed) immersion.
We say that $u$ is a ``weak admissible''  (closed) immersion 
if $(Q \subset \fP)$ is weak admissible.
\end{enumerate}
\end{dfn}

\begin{dfn}
\label{dfn-4.3.4bisY-bc}
Suppose there exists an integer  $r\geq 0$ such that 
$\fP $ is 
a formal  $\bbD ^r _{\fS}$-scheme of finite type
having locally finite $p$-bases  over $\fS$.
Let $Y$ be a  subscheme of $P$.
We say $(Y \subset \fP)$ is ``weak admissible after any base change''
if for any morphism of local algebras
$\alpha \colon \V \to \W$ 
 such that 
 $\cV$ and $\W$ are  complete discrete valued ring  of mixed characteristic $(0,p)$ with respective perfect residue fields $k$ and $l$, 
$(Y \times _{\bbD ^r _{\Spec k}} \bbD ^r _{\Spec l} 
\subset 
\fP \times _{\bbD ^r _{\Spf (\cV) }} \bbD ^r _{\Spf(\cW)} )$
is weak admissible. 
\end{dfn}

\begin{lemm}
\label{lemm-Y'YGammadfn}
Let $Y' \subset Y$ be two subschemes of $P$.
\begin{enumerate}[(a)]
\item 
\label{lemm-Y'YGammadfn-item1}
If $(Y'\subset \fP)$ is weak admissible
then we have the isomorphism of 
$\smash{\underrightarrow{LD}} ^{\mathrm{b}} _{\Q,\mathrm{ovcoh}} 
(\overset{^\mathrm{l}}{} \smash{\widehat{\D}} _{\fP /\fS  } ^{(\bullet)} )$:
\begin{equation}
\label{lemm-Y'YGammadfn-iso1}
\R \underline{\Gamma} ^\dag _{Y } ( 
\R \underline{\Gamma} ^\dag _{Y'} ( \cO _{\fP} ^{(\bullet)}))
\riso 
\R \underline{\Gamma} ^\dag _{Y'} ( \cO _{\fP} ^{(\bullet)}).
\end{equation}

\item 
\label{lemm-Y'YGammadfn-item2}
If $(Y\subset \fP)$ is weak admissible
then 
so is 
$(Y'\subset \fP)$ is weak admissible
and we have the isomorphism of
$\smash{\underrightarrow{LD}} ^{\mathrm{b}} _{\Q,\mathrm{ovcoh}} 
(\overset{^\mathrm{l}}{} \smash{\widehat{\D}} _{\fP /\fS  } ^{(\bullet)} )$
\begin{equation}
\label{lemm-Y'YGammadfn-iso2}
\R \underline{\Gamma} ^\dag _{Y '} ( 
\R \underline{\Gamma} ^\dag _{Y} ( \cO _{\fP} ^{(\bullet)}))
\riso 
\R \underline{\Gamma} ^\dag _{Y'} ( \cO _{\fP} ^{(\bullet)}).
\end{equation}

\end{enumerate}

\end{lemm}

\begin{proof}
We can suppose $P$ integral. 
Let $X'$ (resp. $X$) be the closure of $Y'$ (resp. $Y$) in $P$
The first statement is obvious.
a) Suppose that $(Y'\subset \fP)$ is weak admissible, i.e. 
that 
$(X'\subset \fP)$ is weak admissible.
Then the complex 
$\cE ^{\prime (\bullet) }:= \R \underline{\Gamma} ^\dag _{X'}  (\cO _{\fP} ^{(\bullet)} )$
is well defined as in \ref{dfn-4.3.4bis}.
The isomorphism \ref{lemm-Y'YGammadfn-iso1}, 
is 
$$\R \underline{\Gamma} ^\dag _{Y } ( 
\R \underline{\Gamma} ^\dag _{Y'} ( \cE ^{\prime (\bullet)}))
\underset{\ref{gammayY'}}{\riso} 
\R \underline{\Gamma} ^\dag _{Y'} ( \cE ^{\prime (\bullet)}).$$

b) Suppose that $(Y\subset \fP)$ is weak admissible, i.e. 
that 
$(X\subset \fP)$ is weak admissible.
Then the complex 
$\cE ^{ (\bullet) }:= \R \underline{\Gamma} ^\dag _{X}  (\cO _{\fP} ^{(\bullet)} )$
is well defined as in \ref{dfn-4.3.4bis} and is overcoherent.
This yields that 
$\cE ^{\prime (\bullet) }:= \R \underline{\Gamma} ^\dag _{X '}  (\cO _{\fP} ^{(\bullet)} )$
is well defined as in \ref{dfn-4.3.4bis}, i.e. 
$(X'\subset \fP)$ is weak admissible.
Moreover by construction we have
$\cE ^{\prime (\bullet) }\riso 
\R \underline{\Gamma} ^\dag _{X'} \cE ^{ (\bullet) }$.
Hence, we get the isomorphisms 
$$\R \underline{\Gamma} ^\dag _{Y '} ( 
\R \underline{\Gamma} ^\dag _{Y} ( \cO _{\fP} ^{(\bullet)}))
=
\R \underline{\Gamma} ^\dag _{Y '} ( 
\R \underline{\Gamma} ^\dag _{Y} ( \cE ^{ (\bullet) }))
\underset{\ref{gammayY'}}{\riso} 
\R \underline{\Gamma} ^\dag _{Y'} (\cE ^{ (\bullet) })
\underset{\ref{gammayY'}}{\riso} 
\R \underline{\Gamma} ^\dag _{Y'} \R \underline{\Gamma} ^\dag _{X'} (\cE ^{ (\bullet) })
\riso
\R \underline{\Gamma} ^\dag _{Y'} (\cE ^{\prime (\bullet) })
=
\R \underline{\Gamma} ^\dag _{Y'} ( \cO _{\fP} ^{(\bullet)}),
$$
whose composition is \ref{lemm-Y'YGammadfn-iso2}.
\end{proof}

\begin{lemm}
\label{u+closedimm-YGamma-pre}
Let $u \colon \fX \hookrightarrow \fP$ be a closed immersion
of formal  $\fS$-schemes
of formal finite type
and having locally finite $p$-bases  over $\fS$.
Let $T$ be divisor of $X$.
Then there exists 
an open covering 
$(\fP _{n}) _{n=1,\dots ,N}$ 
of $T$ by open subsets of $\fP$
(i.e. $T \subset \cup _n \fP _n$),
there exist divisors 
$D _1, \dots, D _N$ of $P$ 
such that 
for any $n = 1,\dots, N$ we have
\begin{enumerate}[(a)]
\item $(D _n \cap X )\cap P _n = T \cap P _n$ ;
\item $D _n \supset T $.
\end{enumerate}

\end{lemm}

\begin{proof}
Let $T _1, \dots, T _r$ be the irreducible components of $T$. 
We proceed by induction on $r$.

I) Suppose $r =1$.
There exist 
an open covering 
$(\fP _{n}) _{n=1,\dots ,N}$
of $T$ by open subsets of $\fP$ such that 
$T \cap P _n$ is non empty and is a principal divisor of $X \cap P _n$ for any $n$. 
Hence there exists a principal divisor $\overset{_\circ}{D} _n$ of $P _n$ such that 
$\overset{_\circ}{D} _n \cap X = T \cap P _n$.
Denoting by $D _n$ the closure of $\overset{_\circ}{D} _n$ in $P$,
we get 
$(D _n \cap X )\cap P _n = T \cap P _n$.
Since $T$ is irreducible, then the closure of $T \cap P _n$ is $T$.
This yields $T \subset D _n$.

II) We suppose now that $r \geq 2$.

1) Let $x \in T _1 \cap \dots \cap T _r$.
Let $\fP _x$ be an open subset of $\fP$ containing $x$ such that 
$T \cap P _x$ is a principal divisor of $X \cap P _x$.
Hence there exists a principal divisor $\overset{_\circ}{D} _x$ of $P _x$ such that 
$\overset{_\circ}{D} _x \cap X = T \cap P _x$.
Denoting by $D _x$ the closure of $\overset{_\circ}{D} _x$ in $P$,
we get 
$(D _x \cap X )\cap P _x = T \cap P _x$.
Since $P _x$ is an open subset of $P$ containing $x \in  T _1 \cap \dots \cap T _r$,
then $P _x$ contains the generic points of $T _1, \dots, T _r$.
Hence, the closure of $T \cap P _x$ is $T$.
This yields 
$T \subset D _x$.

2) We set $\fU _n = \fP \setminus T _n$ for $n =1,\dots, r$.

i) For any $n =1,\dots, r$, 
by using the induction hypothesis, 
there exist 
$(\fU  _{ni}) _{i=1,\dots, i _n}$ an open covering of $T \cap U _n=T \setminus T _n$ 
by open subsets of $\fU _n$,
$(\overset{_\circ}{D} _{ni}) _{i=1,\dots, i _n}$ some divisors of $U _{n}$
such that 
for any $i=1,\dots, i _n$ we have 
\begin{enumerate}[(a)]
\item $(\overset{_\circ}{D} _{ni} \cap X )\cap U  _{ni} = T \cap U  _{ni}$ ;
\item $\overset{_\circ}{D} _{ni}  \supset T \cap U  _{n}$.
\end{enumerate}
We denote by 
$D _{ni} $
the closure of 
$\overset{_\circ}{D} _{ni} $ in $P$.
Then $D _{ni} $ is a divisor of $P$ such that 
\begin{enumerate}[(a)]
\item $(D _{ni} \cap X )\cap U  _{ni} = T \cap U  _{ni}$ ;
\item $D _{ni}  \supset T \cap U  _{n}$.
\end{enumerate}

ii) From the part I), 
there exists 
an open covering 
$(\fP _{nj}) _{j=1,\dots ,j _n}$ 
of $T _n$ by open subsets of $\fP$,
there exist divisors 
$E _{n1}, \dots, E _{n j _n}$ of $P$ 
such that 
for any $j = 1,\dots, j _n$ we have
\begin{enumerate}[(a)]
\item $(E _{nj} \cap X )\cap P _{nj} = T _n \cap P _{nj}$ ;
\item $E _{nj} \supset T _n$.
\end{enumerate}

iii) Put $\fU _{nij}: = \fU _{ni} \cap \fP _{n j}$.
Then $ \cup _n T \setminus T _n \subset  \cup _n  \cup _{i} \cup _{j} U _{nij}$.
We get the divisor $F _{nij}:= D _{ni} \cup E _{nj}$ of $P$.
Since $D _{ni} \supset T \setminus T _n$ and 
$E _{nj} \supset T _n$, then 
$F _{nij}\supset T$.
Moreover, 
$(F _{nij} \cap X ) \cap U _{nij}
=
(D _{ni} \cap X \cap U _{ni} \cap P _{n j})
\cup 
(E _{nj} \cap X \cap U _{ni} \cap P _{n j})$.
We have 
$D _{ni} \cap X \cap U _{ni} \cap P _{n j}
=
T \cap U _{ni} \cap P _{n j}$
and
$E _{nj} \cap X \cap U _{ni} \cap P _{n j}
=
T  \cap U _{ni} \cap P _{n j}$.
Hence, 
$(F _{nij} \cap X ) \cap U _{nij}
=
T\cap U _{nij}$.

3) Since $T = (T _1 \cap \dots \cap T _r)
\cup \left ( \cup _{n=1} ^{r} T\setminus T _n\right )$,
then using II.1) and II.2.iii) we get an open covering 
$(\fP _{m}) _{m=1,\dots ,M}$ 
of $T$ by open subsets of $\fP$,
there exist divisors 
$D _1, \dots, D _M$ of $P$ 
such that 
for any $m = 1,\dots, M$ we have
\begin{enumerate}[(a)]
\item $(D _m \cap X )\cap P _m = T \cap P _m$ ;
\item $D _m \supset T $.
\end{enumerate}
More precisely, 
$\fP _m$ is either of the form $\fP _x$ and $D _m$ of the form $D _x$ (case of 1)
or $\fP _m$ is of the form $\fU _{nij}$ and $D _m$ of the form $F _{nij}$ (case of 2.iii).
Hence, we are done.\end{proof}

\begin{prop}
\label{u+closedimm-YGamma}
Let $u \colon \fX \hookrightarrow \fP$ be a closed immersion
of formal  $\fS$-schemes
of formal finite type
and having locally finite $p$-bases  over $\fS$.
Let $Y$ be a subscheme of $X$
such that 
$(Y\subset \fX)$ is weak admissible.
Then 
$(u(Y), \fP)$ is weak admissible
and we have the isomorphism of
$\smash{\underrightarrow{LD}} ^{\mathrm{b}} _{\Q,\mathrm{ovcoh}} 
(\overset{^\mathrm{l}}{} \smash{\widehat{\D}} _{\fP /\fS  } ^{(\bullet)} )$:
\begin{equation}
\label{u+closedimm-YGamma-iso1}
u _+ ^{(\bullet) }
\left ( 
\R \underline{\Gamma} ^\dag _{Y} ( \cO _{\fX} ^{(\bullet)})
\right)
[  \delta _{\fX /\fP}]
\riso 
\R \underline{\Gamma} ^\dag _{Y} ( \cO _{\fP} ^{(\bullet)}).
\end{equation}

\end{prop}

\begin{proof}
We can suppose $P$ and $X$ are integral. 
By using \ref{commutfonctcohlocal2} and by definition (recall \ref{dfn-4.3.4bisY}), 
we reduce to the case where $Y$ is a closed subscheme of $P$.
We keep notation \ref{XpbasesGamma-desc}.

a) Suppose $Y= X$. 
Since 
$u  ^{! (\bullet )}  (\cO _{\fP } ^{(\bullet)})
\riso 
\cO _{\fX } ^{(\bullet)}
[  \delta _{\fP /\fX}]$
is overcoherent, 
then 
$u _{+} ^{(\bullet)} \circ u  ^{! (\bullet )}  (\cO _{\fP } ^{(\bullet)})$
is overcoherent
(use \ref{2.2.18}.\ref{2.2.18-item2}).
Following \ref{XpbasesGamma-desc-iso5},
we have an isomorphism of the form 
$u _{+} ^{(\bullet)} \circ u  ^{! (\bullet )}  (\cO _{\fP } ^{(\bullet)})
\riso 
\phi (\cO _{\fP } ^{(\bullet)})$.
Hence, 
$\phi (\cO _{\fP } ^{(\bullet)})$
is overcoherent.
Since $X = \cap _{i,j}T _{i,j}$, we are done.

b) Suppose $Y\not = X$. 
Hence, 
there exist some divisors 
$T _1, \dots, T _{N}$ of $X$ such that
 1) $Y = \cap _{n =1} ^{N} T _n$
and such that 2) 
$\R \underline{\Gamma} ^\dag _{T _N} \circ \cdots \circ 
\R \underline{\Gamma} ^\dag _{T _1} (\cO _{\fX} ^{(\bullet)} )
\in 
\smash{\underrightarrow{LD}} ^{\mathrm{b}} _{\Q,\mathrm{ovcoh}} 
(\overset{^\mathrm{l}}{} \smash{\widehat{\D}} _{\fX /\fS  } ^{(\bullet)} )$.

Following \ref{u+closedimm-YGamma-pre}, 
for any $n = 1,\dots, N$,
there exists 
an open covering 
$(\fP _{nl}) _{l=1,\dots ,L_n}$ 
of $T _{n}$ by open subsets of $\fP$,
there exist divisors 
$D _{n1}, \dots, D _{nL_n}$ of $P$ 
such that 
for any $l = 1,\dots, L _n$ we have
\begin{enumerate}[(a)]
\item $(D _{nl} \cap X )\cap P _{nl} = T _n \cap P _{nl}$ ;
\item $D _{nl} \supset T _n$.
\end{enumerate}
Since 
$\cap _l D _{nl}  \cap X = T _n$, 
we have
$(\cap _n (\cap _l D _{nl} ) \cap X = \cap _n T _n= Y$.
Let us denote 
$\R \underline{\Gamma} ^\dag _{D _{n \bullet}} 
:= 
\R \underline{\Gamma} ^\dag _{D _{n 1}} 
\circ 
\dots
\circ 
\R \underline{\Gamma} ^\dag _{D _{n L _n}}$
and
$\R \underline{\Gamma} ^\dag _{D _{n \bullet}\cap X} 
:= 
\R \underline{\Gamma} ^\dag _{D _{n 1}\cap X} 
\circ 
\dots
\circ 
\R \underline{\Gamma} ^\dag _{D _{n L _n}\cap X}$,
for any $n$.
Following \ref{surcoh2.1.4-cor}, 
since $D _{nl}  \cap X$ is a divisor of $X$ and $D _{nl}$ is a divisor of $P$, 
we have the isomorphism
\begin{equation}
\label{u+closedimm-YGamma-iso2}
\R \underline{\Gamma} ^\dag _{D _{1 \bullet}} 
\circ \dots \circ 
\R \underline{\Gamma} ^\dag _{D _{N \bullet}} 
u _+ ^{(\bullet) } \cO _{\fX} ^{(\bullet)}
\riso 
u _+ ^{(\bullet) } 
\left ( 
\R \underline{\Gamma} ^\dag _{D _{1 \bullet}\cap X} 
\circ \dots \circ 
\R \underline{\Gamma} ^\dag _{D _{N \bullet}\cap X} 
\cO _{\fX} ^{(\bullet)}
\right).
\end{equation}

For any 
$\underline{l} = (l _1,\dots, l _N) 
\in 
\{1,\dots, L _1 \}
\times \dots \times 
\{N,\dots, L _N\}$,
we set
$\fP _{\underline{l}}:=
\fP _{1 l _1}\cap \dots \cap \fP _{Nl _N}$.
For any 
$\underline{l} = (l _1,\dots, l _N) 
\in 
\{1,\dots, L _1 \}
\times \dots \times 
\{N,\dots, L _N\}$,
we set
$$\label{u+closedimm-YGamma-iso3}
\R \underline{\Gamma} ^\dag _{D _{n \underline{l}}} 
:= 
\R \underline{\Gamma} ^\dag _{D _{n 1}} 
\circ 
\dots
\circ 
\R \underline{\Gamma} ^\dag _{D _{n (l _{n-1})}}
\circ 
\R \underline{\Gamma} ^\dag _{D _{n (l _{n+1})}}
\circ 
\dots \circ 
\R \underline{\Gamma} ^\dag _{D _{n L _n}}.$$
Hence, 
$\R \underline{\Gamma} ^\dag _{D _{n \bullet}} 
=
\R \underline{\Gamma} ^\dag _{D _{n \underline{l}}} 
\circ 
\R \underline{\Gamma} ^\dag _{D _{n l _n}} $.
We define similarly 
$\R \underline{\Gamma} ^\dag _{D _{n \underline{l}}\cap X} $ so that 
$\R \underline{\Gamma} ^\dag _{D _{n \bullet}\cap X} 
=
\R \underline{\Gamma} ^\dag _{D _{n \underline{l}}\cap X} 
\circ 
\R \underline{\Gamma} ^\dag _{D _{n l _n}\cap X} $.

Since 
$D _{n l _n }\cap X \cap P  _{nl _n} 
= 
T _n \cap P  _{nl _n}$,
then 
$D _{n l _n }\cap X \cap P  _{\underline{l}} 
= 
T _n \cap P  _{\underline{l}}$.
By setting $\fX _{\underline{l}}:=
\fX \cap \fP _{\underline{l}}$, we get
$$
\label{u+closedimm-YGamma-iso4}
\left ( 
\R \underline{\Gamma} ^\dag _{D _{n \bullet}\cap X} 
\cO _{\fX} ^{(\bullet)} 
\right) 
|\fX _{\underline{l}}
=
\R \underline{\Gamma} ^\dag _{D _{n \underline{l}}\cap X _{\underline{l}}} 
\circ 
\R \underline{\Gamma} ^\dag _{T _n\cap X _{\underline{l}}} 
\cO _{\fX _{\underline{l}}} ^{(\bullet)} .$$
This yields
\begin{equation}
\label{u+closedimm-YGamma-iso5}
\left (
\R \underline{\Gamma} ^\dag _{D _{1 \bullet}\cap X} 
\circ \dots \circ 
\R \underline{\Gamma} ^\dag _{D _{N \bullet}\cap X} 
\cO _{\fX} ^{(\bullet)} 
\right) 
|\fX _{\underline{l}}
\riso 
\R \underline{\Gamma} ^\dag _{D _{1 \underline{l}}\cap X _{\underline{l}}}
\circ \dots \circ 
\R \underline{\Gamma} ^\dag _{D _{N \underline{l}}\cap X _{\underline{l}}}
\circ 
\R \underline{\Gamma} ^\dag _{T _1 \cap X _{\underline{l}}} 
\circ \dots \circ 
\R \underline{\Gamma} ^\dag _{T _N\cap X _{\underline{l}}} 
\cO _{\fX _{\underline{l}}} ^{(\bullet)} .
\end{equation}
Since
$\R \underline{\Gamma} ^\dag _{T _1 \cap X _{\underline{l}}} 
\circ \dots \circ 
\R \underline{\Gamma} ^\dag _{T _N\cap X _{\underline{l}}} 
\cO _{\fX _{\underline{l}}} ^{(\bullet)} 
=
\left (\R \underline{\Gamma} ^\dag _{T _1} 
\circ \dots \circ 
\R \underline{\Gamma} ^\dag _{T _N} 
\cO _{\fX} ^{(\bullet)}  \right) 
| \fX _{\underline{l}}
=
\R \underline{\Gamma} ^\dag _{Y} 
\cO _{\fX} ^{(\bullet)}  
| \fX _{\underline{l}}$, 
then it follows from 
\ref{lemm-Y'YGammadfn}.\ref{lemm-Y'YGammadfn-item1}
and \ref{u+closedimm-YGamma-iso5}, 
\begin{equation}
\label{u+closedimm-YGamma-iso6}
\left (
\R \underline{\Gamma} ^\dag _{D _{1 \bullet}\cap X} 
\circ \dots \circ 
\R \underline{\Gamma} ^\dag _{D _{N \bullet}\cap X} 
\cO _{\fX} ^{(\bullet)} 
\right) 
|\fX _{\underline{l}}
\riso 
\R \underline{\Gamma} ^\dag _{Y} 
\cO _{\fX} ^{(\bullet)}  
| \fX _{\underline{l}}
.
\end{equation}
Since $\cup _{l} \fP _{nl} \supset T _n$, then 
$\cup _{\underline{l}} \fP _{\underline{l}} \supset \cap _{n} T _n=Y$.
Since 
$\R \underline{\Gamma} ^\dag _{D _{1 \bullet}\cap X} 
\circ \dots \circ 
\R \underline{\Gamma} ^\dag _{D _{N \bullet}\cap X} 
(\cO _{\fX} ^{(\bullet)} )|\fP \setminus Y \riso 0$,
then it follows from 
\ref{u+closedimm-YGamma-iso6} that 
$\R \underline{\Gamma} ^\dag _{D _{1 \bullet}\cap X} 
\circ \dots \circ 
\R \underline{\Gamma} ^\dag _{D _{N \bullet}\cap X} 
\cO _{\fX} ^{(\bullet)} $ is overcoherent.
Since
$\cap _{n,l}D _{nl}  \cap X = \cap _n T _n= Y$,
then we have the canonical isomorphism
$$\R \underline{\Gamma} ^\dag _{D _{1 \bullet}\cap X} 
\circ \dots \circ 
\R \underline{\Gamma} ^\dag _{D _{N \bullet}\cap X} 
\cO _{\fX} ^{(\bullet)} 
\riso 
\R \underline{\Gamma} ^\dag _{Y} 
\cO _{\fX} ^{(\bullet)} .$$
Hence, we get the isomorphism of overcoherent complexes
\begin{equation}\label{u+closedimm-YGamma-iso7}
u _{+} ^{(\bullet)} \R \underline{\Gamma} ^\dag _{Y} 
\cO _{\fX} ^{(\bullet)}
[  \delta _{\fX /\fP}]
\underset{\ref{u+closedimm-YGamma-iso2}}{\riso} 
\R \underline{\Gamma} ^\dag _{D _{1 \bullet}} 
\circ \dots \circ 
\R \underline{\Gamma} ^\dag _{D _{N \bullet}} 
u _{+} ^{(\bullet)} \circ u  ^{! (\bullet )}  (\cO _{\fP } ^{(\bullet)}).
\end{equation}
Following \ref{XpbasesGamma-desc-iso5}
$u _{+} ^{(\bullet)} \circ u  ^{! (\bullet )}  (\cO _{\fP } ^{(\bullet)})
\riso 
\phi (\cO _{\fP } ^{(\bullet)})$, 
where 
$\phi :=
\R \underline{\Gamma} ^\dag _{X _{1}} \circ \cdots \circ 
\R \underline{\Gamma} ^\dag _{X _{n}} $
and 
$\R \underline{\Gamma} ^\dag _{X _i} :=
\R \underline{\Gamma} ^\dag _{T _{i,1}} \circ \cdots \circ 
\R \underline{\Gamma} ^\dag _{T _{i,d}}$
where $T _{ij}$ are some divisors of $P$.
Hence, 
$(u(Y), \fP)$ is weak admissible and we have 
\begin{equation}
\label{u+closedimm-YGamma-iso8}
\R \underline{\Gamma} ^\dag _{D _{1 \bullet}} 
\circ \dots \circ 
\R \underline{\Gamma} ^\dag _{D _{N \bullet}} 
u _{+} ^{(\bullet)} \circ u  ^{! (\bullet )}  (\cO _{\fP } ^{(\bullet)})
\riso 
\R \underline{\Gamma} ^\dag _{Y} ( \cO _{\fP} ^{(\bullet)}).
\end{equation}
Finally, by composing \ref{u+closedimm-YGamma-iso7} with \ref{u+closedimm-YGamma-iso8} we are done.
\end{proof}

\begin{cor}
Let $f \colon \fP \to \bbD ^r _\fS$ be a finite type morphism of formal $\fS$-schemes.
Let $\Delta \colon \fP \hookrightarrow \fP \times _{\fC _S} \fP$ be the diagonal morphism. 
Let $Y$ be a subscheme of $P$ such that 
$(Y\subset \fP)$ is weak admissible. 
Then 
$(\Delta ( Y) \subset \fP \times _{\fC _S} \fP)$ 
is weak admissible. 
\end{cor}

\begin{lem}
\label{u!closedimm-YGamma}
Let $u \colon \fX \hookrightarrow \fP$ be an immersion
of formal  $\fS$-schemes
of formal finite type
and having locally finite $p$-bases  over $\fS$.
Let $Y$ be a subscheme of $P$
such that 
$(Y\subset \fP)$ is weak admissible.
Then 
$(u ^{-1}(Y), \fX)$ is weak admissible
and we have the isomorphism of
$\smash{\underrightarrow{LD}} ^{\mathrm{b}} _{\Q,\mathrm{ovcoh}} 
(\overset{^\mathrm{l}}{} \smash{\widehat{\D}} _{\fX /\fS  } ^{(\bullet)} )$:
$$
\label{u!closedimm-YGamma-iso1pre}
\R \underline{\Gamma} ^\dag _{u ^{-1}(Y)} ( \cO _{\fX} ^{(\bullet)})
[  \delta _{\fX /\fP}]
\riso 
u  ^{!(\bullet) }
\R \underline{\Gamma} ^\dag _{Y} ( \cO _{\fP} ^{(\bullet)}).$$
\end{lem}

\begin{proof}
Following \ref{lemm-Y'YGammadfn}.\ref{lemm-Y'YGammadfn-item2}, 
$(u (u ^{-1}(Y)), \fP)$ is weak admissible. Hence, we reduce to the case where 
$Y= u (u ^{-1}(Y))$.
We can suppose $P$ and $X$ are integral and $Y$ is a closed subscheme of $P$.
We can suppose $X$ is a subscheme of $P$.

There exist some divisors 
$T _1, \dots, T _{r}$ of $P$ such that
 1) $Y = \cap _{i =1} ^{r} T _i$
and such that 2) 
$\R \underline{\Gamma} ^\dag _{T _r} \circ \cdots \circ 
\R \underline{\Gamma} ^\dag _{T _1} (\cO _{\fP} ^{(\bullet)} )
\in 
\smash{\underrightarrow{LD}} ^{\mathrm{b}} _{\Q,\mathrm{ovcoh}} 
(\overset{^\mathrm{l}}{} \smash{\widehat{\D}} _{\fP /\fS  } ^{(\bullet)} )$.
For any $i =1,\dots, r$, we have two cases : either the divisor  $T _i$ contains $X$, 
or $T _i \cap X$ a divisor of $X$.
Reordering the divisors if necessary, 
we can suppose that for any 
$i = 1, \dots , s$, 
$T _i \cap X$ is a divisor of $X$
and for 
$i = s+1,\dots, r$
the divisor $T _i $ contains $X$.
It follows from \ref{lem-u*XT=0}, 
that 
for any 
$i = s+1,\dots, r$
the canonical morphism
$u  ^{!(\bullet) } \R \underline{\Gamma} ^\dag _{T _i} (\cO _{\fP} ^{(\bullet)} )
\to 
u  ^{!(\bullet) } (\cO _{\fP} ^{(\bullet)} )$
is an isomorphism.
It follows from 
\ref{f!commoub}.\ref{f!commoub-item1}
that for any $i=1,\dots, s$, we have the isomorphism
$u  ^{!(\bullet) } 
\circ \R \underline{\Gamma} ^\dag _{T _i} 
\riso 
\R \underline{\Gamma} ^\dag _{T _i\cap X} 
\circ 
u  ^{!(\bullet) } $.
Since 
$u  ^{!(\bullet) } (\cO _{\fP} ^{(\bullet)} )
\riso 
\cO _{\fX} ^{(\bullet)}[  \delta _{\fX /\fP}]$,
this yields 
the isomorphism of
$\smash{\underrightarrow{LD}} ^{\mathrm{b}} _{\Q,\mathrm{qc}} 
(\overset{^\mathrm{l}}{} \smash{\widehat{\D}} _{\fX /\fS  } ^{(\bullet)} )$
\begin{equation}
\label{u!closedimm-YGamma-proofiso1}
u  ^{!(\bullet) } 
\R \underline{\Gamma} ^\dag _{T _r} \circ \cdots \circ 
\R \underline{\Gamma} ^\dag _{T _1} (\cO _{\fP} ^{(\bullet)} )
\riso 
\R \underline{\Gamma} ^\dag _{T _s \cap X} \circ \cdots \circ 
\R \underline{\Gamma} ^\dag _{T _1\cap X} (\cO _{\fX} ^{(\bullet)} )
[  \delta _{\fP /\fX}].
\end{equation}
Since $\R \underline{\Gamma} ^\dag _{T _r} \circ \cdots \circ 
\R \underline{\Gamma} ^\dag _{T _1} (\cO _{\fP} ^{(\bullet)} )
$ is overcoherent then so is 
$u  ^{!(\bullet) } 
\R \underline{\Gamma} ^\dag _{T _r} \circ \cdots \circ 
\R \underline{\Gamma} ^\dag _{T _1} (\cO _{\fP} ^{(\bullet)} )$ 
(see \ref{2.2.18}) and then so is the right term of \ref{u!closedimm-YGamma-proofiso1}.
Moreover, since 
$Y = \cap _{i =1} ^{r} T _i \subset X$,
since $T _i \supset X$ for any $i \geq s+1$, then 
$Y = \cap _{i =1} ^{s} (T _i \cap X)$.
Hence, 
we are done.
\end{proof}

\begin{lem}
\label{u!proj-ft-YGamma}
Let $\fX$ be a smooth formal  $\fS$-scheme
of finite type and let 
$\fP ': = \fP \times _{\fS} \fX$.
Let $Y$ be a subscheme of $P$
such that 
$(Y\subset \fP )$ is weak admissible.
Let $\varpi \colon 
\fP '\to \fP $ 
be the canonical projection and 
let
$Y': =\varpi ^{-1}(Y)$. 
Then $(Y', \fP')$ is weak admissible
and we have the isomorphism of
$\smash{\underrightarrow{LD}} ^{\mathrm{b}} _{\Q,\mathrm{ovcoh}} 
(\overset{^\mathrm{l}}{} \smash{\widehat{\D}} _{\fP '/\fS  } ^{(\bullet)} )$:
$$
\label{u!proj-ft-YGamma-YGamma-iso1}
\R \underline{\Gamma} ^\dag _{Y'} ( \cO _{\fP'} ^{(\bullet)})
[ d _X]
\riso 
\varpi  ^{!(\bullet) }
\R \underline{\Gamma} ^\dag _{Y} ( \cO _{\fP} ^{(\bullet)}).$$
\end{lem}

\begin{proof}
We can suppose $P$ and $X$ are integral and $Y$ is a closed subscheme of $P$.
We can suppose $X$ is a subscheme of $P$.
There exist some divisors 
$T _1, \dots, T _{r}$ of $P$ such that
 1) $Y = \cap _{i =1} ^{r} T _i$
and such that 2) 
$\R \underline{\Gamma} ^\dag _{T _r} \circ \cdots \circ 
\R \underline{\Gamma} ^\dag _{T _1} (\cO _{\fP} ^{(\bullet)} )
\in 
\smash{\underrightarrow{LD}} ^{\mathrm{b}} _{\Q,\mathrm{ovcoh}} 
(\overset{^\mathrm{l}}{} \smash{\widehat{\D}} _{\fP /\fS  } ^{(\bullet)} )$.
Set
$T ' _i := \varpi ^{-1}(T _i)$. 
Since $T ' _i$ is a divisor of $P'$, 
then it follows from 
\ref{f!commoub}.\ref{f!commoub-item1}
that for any $i=1,\dots, s$, we have the isomorphism
$\varpi  ^{!(\bullet) } 
\circ \R \underline{\Gamma} ^\dag _{T _i} 
\riso 
\R \underline{\Gamma} ^\dag _{T ' _i } 
\circ 
\varpi  ^{!(\bullet) } $.
Since 
$\varpi  ^{!(\bullet) } (\cO _{\fP} ^{(\bullet)} )
\riso 
\cO _{\fP '} ^{(\bullet)}[  d _X]$,
this yields 
$$\varpi  ^{!(\bullet) }
\R \underline{\Gamma} ^\dag _{Y} ( \cO _{\fP} ^{(\bullet)})
\riso 
\R \underline{\Gamma} ^\dag _{T '_r} \circ \cdots \circ 
\R \underline{\Gamma} ^\dag _{T '_1} (\cO _{\fP '} ^{(\bullet)} )[ d _X].$$
Hence, we are done.
\end{proof}

\begin{prop}
\label{u!-YGamma}
Let $f\colon \fP '\to \fP$ 
be a quasi-projective 
(in the sense of Definition \ref{projectivefscheme}) morphism of formal $\fS$-schemes of formal finite type and having locally finite $p$-bases.
Let $Y$ be a subscheme of $P$, $Y':= f ^{-1} (Y)$.
If $(Y\subset \fP)$ 
is weak admissible then so is $(Y '\subset \fP')$
and we have the isomorphism of
$\smash{\underrightarrow{LD}} ^{\mathrm{b}} _{\Q,\mathrm{ovcoh}} 
(\overset{^\mathrm{l}}{} \smash{\widehat{\D}} _{\fP ' /\fS  } ^{(\bullet)} )$:
$$
\label{u!closedimm-YGamma-iso1N2}
\R \underline{\Gamma} ^\dag _{Y'} ( \cO _{\fP '} ^{(\bullet)})
[  \delta _{\fP ' /\fP}]
\riso 
f  ^{!(\bullet) }
\R \underline{\Gamma} ^\dag _{Y} ( \cO _{\fP} ^{(\bullet)}).$$
\end{prop}

\begin{proof}
This is a consequence of 
Lemmas \ref{u!closedimm-YGamma} and \ref{u!proj-ft-YGamma}.
\end{proof}

\begin{dfn}
\label{dfn-Gamma-adm}
Let $Y$ be a subscheme of $P$
such that 
$(Y\subset \fP)$ is weak admissible.
Then we define the functor 
$\R\underline{\Gamma} ^\dag _{Y}
\colon 
\smash{\underrightarrow{LD}} ^{\mathrm{b}} _{\Q,\mathrm{qc}} 
(\overset{^\mathrm{l}}{} \smash{\widehat{\D}} _{\fP /\fS  } ^{(\bullet)} )
\to 
\smash{\underrightarrow{LD}} ^{\mathrm{b}} _{\Q,\mathrm{qc}} 
(\overset{^\mathrm{l}}{} \smash{\widehat{\D}} _{\fP /\fS  } ^{(\bullet)} )
$ by setting for any 
$\E  ^{(\bullet)} 
\in 
\smash{\underrightarrow{LD}} ^{\mathrm{b}} _{\Q,\mathrm{qc}} 
(\overset{^\mathrm{l}}{} \smash{\widehat{\D}} _{\fP /\fS  } ^{(\bullet)} )$
\begin{equation}
\notag
\R \underline{\Gamma} ^\dag _{Y} ( \E  ^{(\bullet)})
:=
\R \underline{\Gamma} ^\dag _{Y} ( \cO _{\fP} ^{(\bullet)})
 \smash{\widehat{\otimes}}
^\L _{\cO ^{(\bullet)}  _{\fP}}
\E  ^{(\bullet)}.
\end{equation}
We retrieve the functor of
\ref{3.2.1caro-2006-surcoh-surcv} for overcoherent complexes (use \ref{fonctYY'Gamma-iso}).
\end{dfn}

We can extend Corollary \ref{coro-trace-upre} for quasi-coherent complexes :

\begin{coro}
\label{pre-loc-tri-B-t1T}
Let $u \colon \ZZ  \to \X $ be a closed immersion of 
 formal  $\fS$-schemes of formal finite type
and having locally finite $p$-bases  over $\fS$.
We suppose $(Z \subset \fX)$ weak admissible.
For any 
$\E ^{(\bullet)} \in \smash{\underrightarrow{LD}} ^{\mathrm{b}} _{\Q,\mathrm{qc}}(\overset{^\mathrm{l}}{} \smash{\widehat{\D}} _{\X /\fS  } ^{(\bullet)} )$,
we have the isomorphism
\begin{equation}
\label{pre-loc-tri-B-t1T-iso}
\R \underline{\Gamma} ^\dag _{Z} (\E ^{(\bullet)}) 
\riso 
u _{+} ^{ (\bullet)} \circ  u ^{ !(\bullet)} (\E ^{(\bullet)}),
\end{equation}
where by abuse of notation we denote $u (Z)$ by $Z$.
\end{coro}

\begin{proof}
Using \ref{surcoh2.1.4-cor1},
we reduce to the case where
$\E ^{(\bullet) }= \O _{\X}^{(\bullet)} $.
From Berthelot-Kashiwara's theorem \ref{u!u+=id},
since 
$\R \underline{\Gamma} ^\dag _{Z } (\O _{\X}^{(\bullet)}) $
is coherent with support in $Z $,
we get 
$$u _{+} ^{ (\bullet)}u ^{ !(\bullet)}
\R \underline{\Gamma} ^\dag _{Z } (\O _{\X}^{(\bullet)}) 
\riso 
\R \underline{\Gamma} ^\dag _{Z } (\O _{\X}^{(\bullet)})  .$$
On the other hand,
$$u ^{ !(\bullet)}
\R \underline{\Gamma} ^\dag _{Z } (\O _{\X}^{(\bullet)}) 
\underset{\ref{u!closedimm-YGamma}}{\riso} 
\R \underline{\Gamma} ^\dag _{Z } 
u ^{ !(\bullet)} (\O _{\X}^{(\bullet)} ) 
\riso 
u ^{ !(\bullet)} (\O _{\X}^{(\bullet)} ) .
$$
Hence 
$u _{+} ^{ (\bullet)}u ^{ !(\bullet)}
\R \underline{\Gamma} ^\dag _{Z } (\O _{\X}^{(\bullet)}) 
\riso 
u _{+} ^{ (\bullet)}u ^{ !(\bullet)}
 (\O _{\X}^{(\bullet)} ) $,
 and we are done.\end{proof}
\begin{empt}
\label{propertiesGammaYY'qc}
Let  $Y $ and $Y'$ be two subschemes of  $P$ such that 
$(Y \subset \fP)$ and $(Y '\subset \fP)$ are weak admissible. 
Then, $(Y \cap Y' \subset \fP)$ is also weak admissible.
Let $\E  ^{(\bullet)} ,\FF ^{(\bullet)}\in \underrightarrow{LD} ^{\mathrm{b}} _{\Q ,\mathrm{qc}}
(\smash{\widehat{\D}} _{\fP} ^{(\bullet)})$.
The following properties are obvious consequences of \ref{propertiesGammaYY'}.
\begin{enumerate}[(a)]
\item 
We have the canonical isomorphism functorial in $\E ^{(\bullet)},~ Y$, and $Y'$ :
\begin{equation}
  \label{gammayY'qc}
  \R\underline{\Gamma} ^\dag _{Y} \circ \R\underline{\Gamma} ^\dag _{Y'} (\E  ^{(\bullet)})
  \riso
  \R\underline{\Gamma} ^\dag _{Y \cap Y'} (\E  ^{(\bullet)}).
\end{equation}

\item 
We have the canonical isomorphism functorial in $\E ^{(\bullet)},~\FF ^{(\bullet)},~ Y$, and $Y'$ :
\begin{equation}
\label{fonctYY'Gamma-isoqc}
\R \underline{\Gamma} ^\dag _{Y \cap Y'} (\E ^{(\bullet)}
\smash{\widehat{\otimes}}^\L  _{\O ^{(\bullet)}  _{\fP} }  \FF ^{(\bullet)} )
\riso 
\R \underline{\Gamma} ^\dag _{Y} 
(\E ^{(\bullet)})
\smash{\widehat{\otimes}}^\L  _{\O ^{(\bullet)}  _{\fP} }  
\R \underline{\Gamma} ^\dag _{Y'}
(\FF ^{(\bullet)}).
\end{equation}
\item If $Y '$ is an open (resp. a closed) subscheme of  $Y$, we have the canonical homomorphism
$\R\underline{\Gamma} ^\dag _{Y} (\E  ^{(\bullet)}) \rightarrow \R\underline{\Gamma} ^\dag _{Y'} (\E  ^{(\bullet)})$
(resp. $\R\underline{\Gamma} ^\dag _{Y'} (\E  ^{(\bullet)}) \rightarrow \R\underline{\Gamma} ^\dag _{Y} (\E  ^{(\bullet)})$).
If $Y '$ is a closed subscheme of  $Y$, we have the localization distinguished triangle
$\R\underline{\Gamma} ^\dag _{Y'} (\E  ^{(\bullet)}) \rightarrow \R\underline{\Gamma} ^\dag _{Y} (\E  ^{(\bullet)})
\rightarrow \R\underline{\Gamma} ^\dag _{Y \setminus Y'} (\E  ^{(\bullet)}) \rightarrow +1.$

\end{enumerate}

\end{empt}

\begin{theo}
\label{2.2.18qc}
Let $f\colon \fP '\to \fP$ 
be a quasi-projective (in the sense of Definition \ref{projectivefscheme}) 
morphism of formal $\fS$-schemes of formal finite type and having locally finite $p$-bases.
Let $Y$ be a subscheme of $P$, $Y':= f ^{-1} (Y)$.
We suppose 
$(Y \subset \fP)$ is weak  admissible. 
Then $(Y '\subset \fP')$ is weak admissible. Moreover, we have the following properties.
\begin{enumerate}[(a)]
\item 
\label{2.2.18qc-item1}
Let  $\E ^{(\bullet)} \in 
\smash{\underrightarrow{LD}} ^{\mathrm{b}} _{\Q,\mathrm{qc}} 
(\overset{^\mathrm{l}}{} \smash{\widehat{\D}} _{\fP /\fS  } ^{(\bullet)} )$. 
We have the functorial in $Y$ isomorphism of 
$\smash{\underrightarrow{LD}} ^{\mathrm{b}} _{\Q,\mathrm{qc}} ( \smash{\widehat{\D}} _{\fP '/\fS} ^{(\bullet)})$:
\begin{gather}
\label{commutfonctcohlocal1qc}
  f ^{ !(\bullet)}  \circ\R \underline{\Gamma} ^\dag _{Y}(\E ^{(\bullet)}) 
  \riso
   \R \underline{\Gamma} ^\dag _{Y' }\circ f ^{ !(\bullet)}  (\E ^{(\bullet)}).
\end{gather}
   
\item 
\label{2.2.18qc-item2}
Let 
$\E ^{\prime (\bullet)} \in \smash{\underrightarrow{LD}} ^{\mathrm{b}} _{\Q,\mathrm{qc}} ( \smash{\widehat{\D}} _{\fP '/\fS} ^{(\bullet)})$.
Then 
we have the functorial in $Y$ isomorphism of 
$\smash{\underrightarrow{LD}} ^{\mathrm{b}} _{\Q,\mathrm{qc}} 
(\overset{^\mathrm{l}}{} \smash{\widehat{\D}} _{\fP /\fS  } ^{(\bullet)} )$ :
\begin{gather}
\label{commutfonctcohlocal2qc}
\R \underline{\Gamma} ^\dag _{Y}\circ f ^{ (\bullet)} _{+} (\E ^{\prime (\bullet)})
\riso
f ^{ (\bullet)}  _{+} \circ \R \underline{\Gamma} ^\dag _{Y'}(\E ^{\prime (\bullet)}).
\end{gather}

\end{enumerate}

\end{theo}

\begin{proof}
By definition of the local cohomological functor (see \ref{dfn-Gamma-adm})
and by commutation of tensor products with extraordinary inverse images,
to check the first statement we reduce to the case where
$\E ^{(\bullet)} = \cO _{\fP} ^{(\bullet)}$.
Then, this is \ref{u!-YGamma}.
Now, let us check that \ref{commutfonctcohlocal2qc} is a consequence of
\ref{commutfonctcohlocal1qc}.
\begin{equation}
\R \underline{\Gamma} ^\dag _{Y}\circ f ^{ (\bullet)} _{+} (\E ^{\prime (\bullet)})
=
\R \underline{\Gamma} ^\dag _{Y}(\O ^{(\bullet)}  _{\fP})
\smash{\widehat{\otimes}}^\L  _{\O ^{(\bullet)}  _{\fP} }
 f ^{ (\bullet)} _{+} (\E ^{\prime (\bullet)})
\underset{\ref{surcoh2.1.4}}{\riso} 
 f ^{ (\bullet)} _{+} 
 (
f  ^{!(\bullet)} ( \R \underline{\Gamma} ^\dag _{Y}(\O ^{(\bullet)}  _{\fP}) )
\smash{\widehat{\otimes}}^\L  _{\O ^{(\bullet)}  _{\fP'} }
\E ^{\prime (\bullet)}))[-\delta _{P'/P}]
\end{equation}
Using \ref{commutfonctcohlocal1qc}, we get
$f  ^{!(\bullet)} ( \R \underline{\Gamma} ^\dag _{Y}(\O ^{(\bullet)}  _{\fP}) )
[-\delta _{P'/P}]
\riso 
   \R \underline{\Gamma} ^\dag _{Y' } \O ^{(\bullet)}  _{\fP'}$.
   Hence we are done.
\end{proof}

We will need (see \ref{rem-S+wa}, \ref{lemdfnqupropbis+mdiv} etc.)
to extend the previous theorem when $f$ is not necessarily of finite type as follows. 
\begin{theo}
\label{2.2.18qcgen}
Let $f\colon \fP '\to \fP$ 
be a  
morphism of formal $\fS$-schemes of formal finite type and having locally finite $p$-bases.
We suppose that there exists 
a closed immersion of the form
$u \colon \fP '\hookrightarrow \fQ$
and a flat morphism of the form
$ \varpi \colon \fQ \to \fP$
such that 
$f = \varpi \circ u$. 
Let $Y$ be a subscheme of $P$, 
$Y'$ be a subscheme of $P'$. 
We suppose 
$(Y '\subset \fP')$
is weak admissible.
We suppose moreover either
$(Y \subset \fP)$ is weak admissible,
either $Y$ is the support of a divisor of $P$ or
$P \setminus Y$ is the support of a divisor of $P$. 
Let  $\E ^{(\bullet)} \in 
\smash{\underrightarrow{LD}} ^{\mathrm{b}} _{\Q,\mathrm{qc}} 
(\overset{^\mathrm{l}}{} \smash{\widehat{\D}} _{\fP /\fS  } ^{(\bullet)} )$. 
We have the canonical isomorphism of 
$\smash{\underrightarrow{LD}} ^{\mathrm{b}} _{\Q,\mathrm{qc}} ( \smash{\widehat{\D}} _{\fP '/\fS} ^{(\bullet)})$:
\begin{gather}
\label{commutfonctcohlocal1qc-pre}
\R \underline{\Gamma} ^\dag _{Y' } \circ  f ^{ !(\bullet)}  \circ\R \underline{\Gamma} ^\dag _{Y}(\E ^{(\bullet)}) 
\riso
\R \underline{\Gamma} ^\dag _{Y' \cap f ^{-1} (Y)}\circ f ^{ !(\bullet)}  (\E ^{(\bullet)}).
\end{gather}

\end{theo}

\begin{proof}
By definition of the local cohomological functor (see \ref{dfn-Gamma-adm})
and by commutation of tensor products with extraordinary inverse images,
to check the first statement we reduce to the case where
$\E ^{(\bullet)} = \cO _{\fP} ^{(\bullet)}$.
Since the other cases are easier, let us suppose
$(Y \subset \fP)$ is weak admissible.
We can suppose $P$ is integral and $Y$ is a closed subscheme of $P$.
By hypothesis, there exist some divisors 
$T _1, \dots, T _{N}$ of $P$ such that
 1) $Y = \cap _{n =1} ^{N} T _n$
and such that 2) 
$\R \underline{\Gamma} ^\dag _{Y}
(\cO _{\fP} ^{(\bullet)} ) 
:=
R \underline{\Gamma} ^\dag _{T _N} \circ \cdots \circ 
\R \underline{\Gamma} ^\dag _{T _1} 
(\cO _{\fP} ^{(\bullet)} )
\in 
\smash{\underrightarrow{LD}} ^{\mathrm{b}} _{\Q,\mathrm{ovcoh}} 
(\overset{^\mathrm{l}}{} \smash{\widehat{\D}} _{\fP /\fS  } ^{(\bullet)} )$.
Since $\varpi$ is flat, 
then the functor 
$\varpi ^{!(\bullet) }$ commutes with 
functors of the form 
$(\hdag T)$, where $T$ is a divisor of $P$.
Hence, 
the functor 
$\varpi ^{!(\bullet) }$ commutes with 
functors of the form 
$\R \underline{\Gamma} ^\dag _{T} $, where $T$ is a divisor of $P$.
Since 
$\varpi ^{!(\bullet) } (\cO _{\fP} ^{(\bullet)} )
\riso
\cO _{\fQ} ^{(\bullet)} [\delta _{Q/P}]$, 
this yields the last isomorphism
\begin{equation}
\label{commutfonctcohlocal1qc-pre-iso1}
f ^{ !(\bullet)} \circ \R \underline{\Gamma} ^\dag _{Y}
(\cO _{\fP} ^{(\bullet)} ) 
\riso
u ^{(\bullet) !} \circ \varpi ^{!(\bullet) }\circ \R \underline{\Gamma} ^\dag _{Y}
(\cO _{\fP} ^{(\bullet)} ) 
\riso
u ^{(\bullet) !} \circ   \R \underline{\Gamma} ^\dag _{\varpi ^{-1}(T _N)} \circ \cdots \circ 
\R \underline{\Gamma} ^\dag _{\varpi ^{-1}(T _1)}
(\cO _{\fQ} ^{(\bullet)}) [\delta _{Q/P}].
\end{equation}
Since $(Y '\subset \fP')$ is weak admissible, then so is 
$(Y '\subset \fQ)$ (see \ref{u+closedimm-YGamma}).

By definition of the local cohomological functor (see \ref{dfn-Gamma-adm})
and by commutation of tensor products with extraordinary inverse images,
we get the isomorphisms
\begin{gather}
\notag
\R \underline{\Gamma} ^\dag _{Y' } u ^{(\bullet) !} \circ   \R \underline{\Gamma} ^\dag _{\varpi ^{-1}(T _N)} \circ \cdots \circ 
\R \underline{\Gamma} ^\dag _{\varpi ^{-1}(T _1)}
(\cO _{\fQ} ^{(\bullet)}) 
\riso 
\R \underline{\Gamma} ^\dag _{Y '} ( \cO _{\fP'} ^{(\bullet)})
 \smash{\widehat{\otimes}}
^\L _{\cO ^{(\bullet)}  _{\fP '}  }
u ^{(\bullet) !} \circ   \R \underline{\Gamma} ^\dag _{\varpi ^{-1}(T _N)} \circ \cdots \circ 
\R \underline{\Gamma} ^\dag _{\varpi ^{-1}(T _1)}
(\cO _{\fQ} ^{(\bullet)})
\\
\notag
\underset{\ref{u!closedimm-YGamma}}{\riso}
u ^{(\bullet) !} \R \underline{\Gamma} ^\dag _{Y '} ( \cO _{\fQ} ^{(\bullet)})
 \smash{\widehat{\otimes}}
^\L _{\cO ^{(\bullet)}  _{\fP'}  }
u ^{(\bullet) !} \circ   \R \underline{\Gamma} ^\dag _{\varpi ^{-1}(T _N)} \circ \cdots \circ 
\R \underline{\Gamma} ^\dag _{\varpi ^{-1}(T _1)}
(\cO _{\fQ} ^{(\bullet)})
[-\delta _{P'/Q}]
\\
\notag
\riso 
u ^{(\bullet) !} 
\left (
\R \underline{\Gamma} ^\dag _{Y '} ( \cO _{\fQ} ^{(\bullet)})
 \smash{\widehat{\otimes}}
^\L _{\cO ^{(\bullet)}  _{\fQ}  }
\R \underline{\Gamma} ^\dag _{\varpi ^{-1}(T _N)} \circ \cdots \circ 
\R \underline{\Gamma} ^\dag _{\varpi ^{-1}(T _1)}
(\cO _{\fQ} ^{(\bullet)})
\right) 
\\
\label{commutfonctcohlocal1qc-pre-iso2}
\riso 
u ^{(\bullet) !} 
\R \underline{\Gamma} ^\dag _{Y '} \circ
\R \underline{\Gamma} ^\dag _{\varpi ^{-1}(T _N)} \circ \cdots \circ 
\R \underline{\Gamma} ^\dag _{\varpi ^{-1}(T _1)}
(\cO _{\fQ} ^{(\bullet)})
\riso 
u ^{(\bullet) !} 
\R \underline{\Gamma} ^\dag _{Y' \cap \varpi  ^{-1} (Y)} 
(\cO _{\fQ} ^{(\bullet)})
\underset{\ref{u!closedimm-YGamma}}{\riso}
\R \underline{\Gamma} ^\dag _{Y' \cap f ^{-1} (Y)} ( \cO _{\fP'} ^{(\bullet)})
[\delta _{P'/Q}].
\end{gather}
By applying the functor 
$\R \underline{\Gamma} ^\dag _{Y' }$ to the isomorphism \ref{commutfonctcohlocal1qc-pre-iso1}
and by composing it with 
\ref{commutfonctcohlocal1qc-pre-iso2},
we get the first isomorphism of the isomorphisms of 
$\smash{\underrightarrow{LD}} ^{\mathrm{b}} _{\Q,\mathrm{ovcoh}} ( \smash{\widehat{\D}} _{\fP '/\fS} ^{(\bullet)})$:
$$
\R \underline{\Gamma} ^\dag _{Y' } \circ  f ^{ !(\bullet)}  \circ\R \underline{\Gamma} ^\dag _{Y}
 (\cO _{\fP} ^{(\bullet)} )
\riso 
\R \underline{\Gamma} ^\dag _{Y '\cap f ^{-1} (Y)} 
(\cO _{\fP'} ^{(\bullet)}) [\delta _{P'/P}]
\riso 
\R \underline{\Gamma} ^\dag _{Y '\cap f ^{-1} (Y)} 
\circ 
f ^{ !(\bullet)}  (\cO _{\fP} ^{(\bullet)}) .
$$
\end{proof}

\begin{ex}
\label{2.2.18qcgen-ex}
Let $a\colon \fP' \to \bbD ^r _{\fS}$ and
$b\colon \fP \to \bbD ^s _{\fS}$ be two objects of 
$\scr{C} _{\fS}$ (see notation \ref{dfn-CfS}).
We suppose moreover that 
$\fP'/\fS$ and $\fP/\fS$ have locally finite $p$-bases.
Let $(f,g) 
\colon
a \to b$ be a morphism of $\scr{C} _{\fS}$.
Then 
$f $ 
is the composition of the graph morphism
$\gamma _{(f,g)}
\colon 
\fP' 
\hookrightarrow 
\fP' \times _{\scr{C} _{\fS}}
\fP$ which is a closed immersion
(see \ref{dfn-CfS2})
with the projection
$\fP' \times _{\scr{C} _{\fS}}
\fP
\to 
\fP$
which is flat (see \ref{corXtimeCY2XtimesYform}).
Hence we are in the situation to apply 
\ref{2.2.18qcgen}.

\end{ex}

\subsection{Base change isomorphism for relatively proper coherent complexes with respect to quasi-projective morphisms}

\begin{theo}
\label{theo-iso-chgtbase}

Let $a\colon \fY \to \bbD ^r _{\fS}$ and
$b\colon \fX \to \bbD ^s _{\fS}$ be two objects of 
$\scr{C} _{\fS}$ (see notation \ref{dfn-CfS}).
We suppose moreover that 
$\fY/\fS$ and $\fX/\fS$ have locally finite $p$-bases.
Let $(g,h) 
\colon
a \to b$ be a morphism of $\scr{C} _{\fS}$.
We suppose that 
$(Y \subset \fY)$ 
is weak admissible and that $g$ 
has locally finite $p$-bases.

Let  $f \colon \X ^{\prime } \to \X $ be a quasi-projective  morphism 
of formal $\fS$-schemes
having locally finite $p$-bases. 
Set $\Y ^{\prime } := \X ^{\prime } \times _{\X} \Y$, 
$f ' \colon \Y ^{\prime } \to \Y$, $g ' \colon \Y ^{\prime }\to \X ^{\prime }$ 
be the canonical projections. 
Let
$\E ^{\prime (\bullet)}
\in  \smash{\underrightarrow{LD}} ^\mathrm{b} _{\Q, \mathrm{coh}}
(\overset{^\mathrm{l}}{} \smash{\widehat{\D}} _{\X ^{\prime }/\fS  } ^{(\bullet)} )$ with proper support over $X$. 
There exists a canonical isomorphism in 
$\smash{\underrightarrow{LD}} ^\mathrm{b} _{\Q, \mathrm{coh}}
(\overset{^\mathrm{l}}{} \smash{\widehat{\D}} _{\X ^{\prime }/\fS  } ^{(\bullet)} )$:
\begin{equation}
\label{basechange}
g ^{ !(\bullet)} \circ f ^{(\bullet)}_{+} (\E ^{\prime (\bullet)})
\riso
f  ^{\prime (\bullet)}_{+}  \circ g ^{\prime(\bullet) !} (\E ^{\prime (\bullet)}). 
\end{equation}
\end{theo}

\begin{proof}
By copying the proof of \cite[10.3.4]{caro-6operations}, we can check that this is a corollary of 
Theorem \ref{u!u+=id}, 
\ref{theo-iso-chgtbase2}, 
\ref{commutfonctcohlocal2qc}
and \ref{pre-loc-tri-B-t1T-iso}.
\end{proof}

\section{Stability under Grothendieck's six operations}

\subsection{Data of absolute coefficients}

\begin{dfn}
\label{DVR}
We define the category 
$\mathrm{DVR}  (\V)$ 
as follows : 
an object is a
morphism of local algebras
$\V \to \W$ such that 
$\W$ is a complete discrete valued ring  of mixed characteristic $(0,p)$
with perfect residue field.
A morphism 
$\W
\to 
\W '$ 
is the data of 
a morphism of local $\V$-algebras
$\W \to \W'$.

\end{dfn}

\begin{empt}
[Convention]
\label{ftV[[t]]}
Let $\W$ be an object of $\mathrm{DVR}  (\V)$.
For simplicity, for any integer $r \geq 0$, we set
$\bbD ^r _{\cW}:= \bbD ^r _{\Spf \cW}$
(we hope this can not be confused with 
$\bbD ^r _{\Spec \cW}$)
and $\scr{C} _{\cW} := \scr{C} _{\Spf \cW}$ (see notation \ref{dfn-CfS}).
In this section, we work exclusively 
with 
the full subcategory of 
$\scr{C} _{\cW}$ consisting of 
formal  $\bbD ^r _{\Spf \cW}$-schemes of finite type for some integer $r$ (not fixed), 
having locally finite $p$-bases  over $\Spf \W$. 
By abuse of notation, 
an object 
$\fP \to \bbD ^r _{\cW}$ of $\scr{C} _{\cW}$ will simply be denoted by 
$\fP$
and a morphism 
$(f,g) \colon 
(\fP ' \to \bbD ^{r '} _{\cW})
\to 
(\fP \to \bbD ^r _{\cW})$
of $\scr{C} _{\cW}$ will simply be denoted by 
$\fP' \to \fP$. 
Moreover a morphism between formal $\Spf \cW$-schemes 
will  mean a morphism in $\scr{C} _{\Spf \cW}$.
\end{empt}

\begin{empt}
\label{t-structure-coh}
Let $\W$ be an object of $\mathrm{DVR}  (\V)$, and 
$\X$ be a formal  $\Spf(\cW)$-scheme of formal finite type, having locally finite $p$-bases  over $\Spf \cW$.
If there is no possible confusion (some confusion might arise specially when the homomorphism $\V \to \W$ is not finite and etale), 
for any integer $m \in \N$, 
we denote 
$\smash{\widehat{\D}} _{\X/\Spf (\W)} ^{(m)}$
(resp. $\smash{\D} ^\dag _{\X/\Spf (\W), \Q}$)
simply by 
$\smash{\widehat{\D}} _{\X} ^{(m)}$
(resp. $\smash{\D} ^\dag _{\X\Q}$).
Berthelot checked the following equivalence of categories (see \cite[4.2.4]{Beintro2}, or \ref{eq-catLDBer-LD-D}):
\begin{equation}
\label{limeqcat}
\underrightarrow{\lim}
\colon 
\smash{\underrightarrow{LD}} ^{\mathrm{b}} _{\Q,\mathrm{coh}} ( \smash{\widehat{\D}} _{\X} ^{(\bullet)})
\cong 
D ^{\mathrm{b}}  _{\mathrm{coh}} (\smash{\D} ^\dag _{\X\Q}).
\end{equation}

The category 
$D ^{\mathrm{b}}  _{\mathrm{coh}} (\smash{\D} ^\dag _{\X\Q})$ 
is endowed with its usual t-structure.
Via \ref{limeqcat}, we get a t-structure on 
$\smash{\underrightarrow{LD}} ^{\mathrm{b}} _{\Q,\mathrm{coh}} ( \smash{\widehat{\D}} _{\X} ^{(\bullet)})$
whose heart is 
$\smash{\underrightarrow{LM}}  _{\Q,\mathrm{coh}} ( \smash{\widehat{\D}} _{\X} ^{(\bullet)})$
(see Notation \ref{nota-(L)Mcoh}).
Recall, following \ref{empt-diag-Hn-comp},
we have canonical explicit cohomological functors
$\mathcal{H} ^n 
\colon 
\smash{\underrightarrow{LD}} ^{\mathrm{b}} _{\Q,\mathrm{coh}} ( \smash{\widehat{\D}} _{\X} ^{(\bullet)})
\to 
\smash{\underrightarrow{LM}} _{\Q,\mathrm{coh}} ( \smash{\widehat{\D}} _{\X} ^{(\bullet)})$.
The equivalence of categories 
\ref{limeqcat} commutes with the 
cohomogical functors $\mathcal{H} ^n $
(where the cohomogical functors $\mathcal{H} ^n $ on 
$D ^{\mathrm{b}}  _{\mathrm{coh}} (\smash{\D} ^\dag _{\X\Q})$
are the obvious ones),
 i.e. 
$\underrightarrow{\lim}
\mathcal{H} ^n (\E ^{(\bullet)})$
is  canonically isomorphic
to 
$\mathcal{H} ^n (\underrightarrow{\lim} \,\E ^{(\bullet)})$.

Last but not least, 
following \ref{eqcat-limcoh} 
we have the equivalence of categories 
$\smash{\underrightarrow{LD}} ^{\mathrm{b}} _{\Q,\mathrm{coh}} ( \smash{\widehat{\D}} _{\X} ^{(\bullet)})
\cong 
D ^{\mathrm{b}} _{\mathrm{coh}}
(
\smash{\underrightarrow{LM}} _{\Q} ( \smash{\widehat{\D}} _{\X} ^{(\bullet)})
)
$
which is also compatible with t-structures, 
where the t-structure on 
$D ^{\mathrm{b}} _{\mathrm{coh}}
(\smash{\underrightarrow{LM}} _{\Q} ( \smash{\widehat{\D}} _{\X} ^{(\bullet)}))$
is the canonical one as the derived category of an abelian category.
\end{empt}

\begin{dfn}
\label{dfn-datacoef}
\begin{enumerate}[(a)]
\item A {\it data of absolute coefficients $\mathfrak{C}$ over $\fS$}
(resp. a {\it weak data of absolute coefficients $\mathfrak{C}$ over $\fS$})
will be the data for any object
$\W$ of $\mathrm{DVR}  (\V)$ (see notation \ref{DVR}), 
for any formal  $\Spf(\cW)$-scheme of formal finite type, having locally finite $p$-bases  over $\Spf \W$
of a strictly full subcategory of 
$\smash{\underrightarrow{LD}} ^{\mathrm{b}} _{\Q,\mathrm{coh}} ( \smash{\widehat{\D}} _{\X} ^{(\bullet)})$
(resp. $\smash{\underrightarrow{LD}} ^{\mathrm{b}} _{\Q,\mathrm{qc}} ( \smash{\widehat{\D}} _{\X} ^{(\bullet)})$),
which will be denoted by $\mathfrak{C} (\X/\W)$,
or simply $\mathfrak{C} (\X)$ if there is no ambiguity with the base $\W$.
If there is no ambiguity with $\V$, 
we simply say a {\it data of absolute coefficients}.

\item A {\it restricted data of absolute coefficients $\mathfrak{C}$ over $\fS$}
(resp. a {\it restricted weak data of absolute coefficients $\mathfrak{C}$ over $\fS$})
will be the data for any object
$\W$ of $\mathrm{DVR}  (\V)$,
for any formal $\bbD ^1 _{\cW}$-scheme of finite type, having locally finite $p$-bases  over $\Spf \W$
of a strictly full subcategory of 
$\smash{\underrightarrow{LD}} ^{\mathrm{b}} _{\Q,\mathrm{coh}} ( \smash{\widehat{\D}} _{\X} ^{(\bullet)})$
(resp. $\smash{\underrightarrow{LD}} ^{\mathrm{b}} _{\Q,\mathrm{qc}} ( \smash{\widehat{\D}} _{\X} ^{(\bullet)})$),
which will be denoted by $\mathfrak{C} (\X/\W)$,
or simply $\mathfrak{C} (\X)$ if there is no ambiguity with the base $\W$.
If there is no ambiguity with $\V$, 
we simply say a {\it restricted data of absolute coefficients}.

\item Let $\mathfrak{C}$ be a (weak) data of coefficients over $\fS$. 
By restriction, we get a restricted (weak) data of coefficients over $\fS$,
that we will denote by 
$\fC ^{(1)}$. 
\end{enumerate}

\end{dfn}

\begin{exs}
\label{ex-Dcst}

\begin{enumerate}[(a)]
\item We define the data of absolute coefficients $\fB _\emptyset$ as follows: 
for any object $\cW$ of $\mathrm{DVR}  (\cV)$,  
for any  formal  $\Spf(\cW)$-scheme of formal finite type $\fX$ having locally finite $p$-bases  over $\Spf \cW$,  
the category $\fB _\emptyset (\X)$ is the full subcategory of 
$\smash{\underrightarrow{LD}} ^{\mathrm{b}} _{\Q,\mathrm{coh}} ( \smash{\widehat{\D}} _{\X} ^{(\bullet)})$
whose unique object is 
$\O _{\X} ^{(\bullet)}$ 
(where $\O _{\X} ^{(\bullet)}$ is the constant object $\O _{\X} ^{(m)} = \O _{\X}$ for any $m\in \N$ with the identity as transition maps).

\item We define the {\it weak} data of absolute coefficients $\fB _\mathrm{div}$ 
as follows: 
for any object $\cW$ of $\mathrm{DVR}  (\cV)$,  
for any  formal  $\Spf(\cW)$-scheme of formal finite type $\fX$ having locally finite $p$-bases  over $\Spf \cW$,  
the category $\fB _\mathrm{div}(\X)$ is the full subcategory of 
$\smash{\underrightarrow{LD}} ^{\mathrm{b}} _{\Q,\mathrm{qc}} ( \smash{\widehat{\D}} _{\X} ^{(\bullet)})$
whose objects are of the form  
$\widehat{\B} ^{(\bullet)} _{\X} (T)$, 
where $T$ is any divisor of the special fiber of $\X$.

Following \ref{O-ovcoh}, 
the restricted weak data 
$\fB _\mathrm{div} ^{(1)}$ is in fact 
a restricted data of absolute coefficients.

\item We define the 
restricted data of absolute coefficients $\fB _\mathrm{cst} ^{(1)}$
as follows: 
for any object $\cW$ of $\mathrm{DVR}  (\cV)$, 
for any  formal  $\bbD ^1 _{\cW}$-scheme of finite type $\fX$ having locally finite $p$-bases  over $\Spf \cW$,  
the category $\fB _\mathrm{cst}(\X)$ is the full subcategory of 
$\smash{\underrightarrow{LD}} ^{\mathrm{b}} _{\Q,\mathrm{coh}} ( \smash{\widehat{\D}} _{\X} ^{(\bullet)})$
whose objects are of the form  
$\R \underline{\Gamma} ^\dag _{Y} \O _\X ^{(\bullet)} $, 
where $Y$ is a subscheme of the special fiber of $\X$
and the functor 
$\R \underline{\Gamma} ^\dag _{Y}$ is defined in 
\ref{3.2.1caro-2006-surcoh-surcv} (use also \ref{O-ovcoh}).

\item We define the data of absolute coefficients $\fB _\mathrm{wa}$ 
as follows: 
for any object $\cW$ of $\mathrm{DVR}  (\cV)$,  
for any  formal  $\Spf(\cW)$-scheme of formal finite type $\fX$ having locally finite $p$-bases  over $\Spf \cW$,  
the category $\fB _\mathrm{wa}(\X)$ is the full subcategory of 
$\smash{\underrightarrow{LD}} ^{\mathrm{b}} _{\Q,\mathrm{coh}} ( \smash{\widehat{\D}} _{\X} ^{(\bullet)})$
whose objects are of the form  
$\R \underline{\Gamma} ^\dag _{Y} \O _\X ^{(\bullet)} $,
where $Y$ is a subscheme of the special fiber of $\X$
is such that 
$(Y \subset \fX)$
is weak admissible after any base change (see \ref{dfn-4.3.4bisY-bc}).
Remark that following 
\ref{O-ovcoh},
we have 
$\fB _\mathrm{wa} ^{(1)}
=
\fB _\mathrm{cst} ^{(1)}$.

\item We define weak data (resp. data) 
$\fM _\mathrm{div}$,
(resp. $\fM _\emptyset$
resp. $\fM _\mathrm{sn}$,
resp. $\fM _\mathrm{n}$) of  absolute coefficients over $\fS$ 
as follows: 
for any object $\W$ of $\mathrm{DVR}  (\V)$ with special fiber $l$, 
for any  formal  $\Spf(\cW)$-scheme of formal finite type $\fP$ having locally finite $p$-bases  over $\Spf \cW$,  
the category 
$\fM _\mathrm{div} (\fP)$
(resp. $\fM _\emptyset(\fP)$
resp. $\fM _\mathrm{sn}(\fP)$)
is the full subcategory of 
$\smash{\underrightarrow{LD}} ^{\mathrm{b}} _{\Q,\mathrm{qc}} ( \smash{\widehat{\D}} _{\fP} ^{(\bullet)})$
(resp. 
$\smash{\underrightarrow{LD}} ^{\mathrm{b}} _{\Q,\mathrm{coh}} ( \smash{\widehat{\D}} _{\fP} ^{(\bullet)})$)
consisting of objects 
of the form  
$(\hdag T)  (\E ^{(\bullet)} )$, 
where 
$\E ^{(\bullet)}
\in 
\mathrm{MIC} ^{(\bullet)} (X, \fP/K) $ (see notation
\ref{ntnMICdag2fs3}),
where $X$ is a closed subscheme of $P$ having locally finite $p$-bases  over $\Spec l$,
$T$ is a divisor of $X$
(resp. $T$ is the empty set,
resp. $T$ is a strictly nice divisor of $X/\Spec l$ in the sense of \ref{st-nice-div},
resp. $T$ is a nice divisor of $X$ in the sense of \ref{nice-div}).
Recall that following 
\ref{coh-ss-div-bis-nice},
these respective objects are indeed coherent.
Following \ref{coh-ss-div-bis},
$\fM _\mathrm{div} ^{(1)}$ is a restricted data of  absolute coefficients over $\fS$.

\end{enumerate}

\end{exs}

\begin{dfn}
\label{dfn-stable-data}
In order to be precise, let us fix some terminology.
Let $\mathfrak{C}$ and $\mathfrak{D}$ be two data of absolute coefficients over $\fS$. 
\begin{enumerate}[(a)]
\item We will say that the  data of absolute coefficients $\mathfrak{C}$  is
stable under 
pushforwards 
if for any object $\cW$ of $\mathrm{DVR}  (\cV)$,  
for any 
{\it quasi-projective} morphism $g \colon \X ' \to \X$ (in the sense of Definition \ref{projectivefscheme}) of  
formal  $\Spf (\cW)$-schemes
of formal finite type having locally finite $p$-bases  over $\Spf \cW$, 
for any object $\E ^{\prime (\bullet)}$ of $\mathfrak{C} (\X')$ with proper support over $X$ via $g$,
the complex $g _{+} (\E ^{\prime (\bullet)})$ is an object of  $\mathfrak{C} (\X)$.

\item We will say that the  data of absolute coefficients $\mathfrak{C}$  is stable under extraordinary pullbacks 
(resp. extraordinary pullbacks by smooth projections,
resp. extraordinary pullbacks by projections,
resp. quasi-projective extraordinary pullbacks,
resp. extraordinary pullbacks by closed immersions,
resp. extraordinary pullbacks by  weak admissible closed immersions)
if for any object $\cW$ of $\mathrm{DVR}  (\cV)$, 
for any morphism (resp. projection morphism in the sense of \ref{dfn-CfS2} which is smooth, 
resp. projection morphism,
resp. quasi-projective morphism,
resp. closed immersion,
resp. weak admissible closed immersion)
$f \colon \Y \to \X$ 
of  formal  $\Spf \cW$-schemes, having locally finite $p$-bases  over $\Spf \W$
(and in  the essential image of the functor 
$\mathscr{S} _{\Spf \cW}$), 
for any object $\E ^{(\bullet)}$ of $\mathfrak{C} (\X)$, 
we have $f ^{!} (\E ^{(\bullet)})\in \mathfrak{C} (\Y)$ (see \ref{dfn-4.3.4bisY}.\ref{dfn-4.3.4bisY-item2}).

\item We will say that the  data of absolute coefficients $\mathfrak{C}$  satisfies 
the first property (resp. the second property) of Berthelot-Kashiwara theorem
or satisfies $BK ^!$ (resp. $BK _+$) for short if the following property is satisfied:
for any object $\cW$ of $\mathrm{DVR}  (\cV)$,  
for any closed immersion  $u \colon \ZZ \hookrightarrow \X$ of  formal  $\Spf (\cW)$-schemes of formal finite type, having locally finite $p$-bases  over $\Spf \W$, 
for any object $\E ^{(\bullet)}$ of $\mathfrak{C} (\X)$ with support in $\ZZ$, 
we have $u ^{!} (\E ^{(\bullet)})\in \mathfrak{C} (\ZZ)$
(resp. 
for any object $\G ^{(\bullet)}$ of $\mathfrak{C} (\ZZ)$, 
we have $u _{+} (\G ^{(\bullet)})\in \mathfrak{C} (\X)$).
Remark that $BK ^!$ and $BK _+$ hold if and only if the  data of absolute coefficients $\mathfrak{C}$ satisfies 
(an analogue of) Berthelot-Kashiwara theorem, which justifies the terminology. 

\item We will say that the  data of absolute coefficients $\mathfrak{C}$ is stable under base change 
if for any morphism 
$\W
\to 
\W '$ 
of $\mathrm{DVR}  (\V)$,
for any integer $r \geq 0$, 
for any  formal  $\bbD ^r _{\cW}$-scheme of finite type $\fX$ having locally finite $p$-bases  over $\Spf \cW$, 
for any object $\E ^{(\bullet)}$ of $\mathfrak{C} (\X)$, 
we have 
$\bbD ^r _{\cW '} 
\smash{\widehat{\otimes}}  ^\L   _{\bbD ^r _{\cW}}  
\E^{(\bullet) }\in \mathfrak{C} (\X \times _{\bbD ^r _{\cW}} \bbD ^r _{\cW'})$.

\item We will say that the  data of absolute coefficients $\mathfrak{C}$ is stable under tensor products (resp. weak admissible tensor products) if
for any object $\cW$ of $\mathrm{DVR}  (\cV)$,  
for any  formal  $\Spf(\cW)$-scheme of formal finite type $\fX$ having locally finite $p$-bases  over $\Spf \cW$, 
for any objects $\E ^{(\bullet)}$ and $\FF ^{(\bullet)}$ of $\mathfrak{C} (\X)$
(resp. and for any weak admissible inclusion $(Y \subset \fX)$)
we have 
$\R \underline{\Gamma} ^\dag _{Y}   \FF ^{(\bullet)}
\smash{\widehat{\otimes}}^\L
_{\O  _{\X}} \E^{(\bullet) }
\in \mathfrak{C} (\X)$.

\item We will say that the  data of absolute coefficients $\mathfrak{C}$ is stable under duality (resp. weak admissible duality) if
for any object $\cW$ of $\mathrm{DVR}  (\cV)$,  
for any  formal  $\Spf(\cW)$-scheme of formal finite type $\fX$ having locally finite $p$-bases  over $\Spf \cW$, 
for any object $\E ^{(\bullet)}$ of $\mathfrak{C} (\X)$
we have 
$\DD  _{\X}(\E^{(\bullet) }) \in  \mathfrak{C} (\X)$
(resp. 
$\DD  _{\X}(\R \underline{\Gamma} ^\dag _{Y}  \E^{(\bullet) }) \in  \mathfrak{C} (\X)$).

\item We will say that the  data of absolute coefficients $\mathfrak{C}$ is stable under weak admissible external tensor products if
for any object $\cW$ of $\mathrm{DVR}  (\cV)$, 
for any  formal  $\Spf (\cW)$-scheme
of formal finite type $\fP$ and having locally finite $p$-bases  over $\Spf \cW$, 
for any weak admissible inclusion $(Y \subset \fP)$,
for any  formal  $\Spf (\cW)$-scheme $\fQ$ of formal finite type and having locally finite $p$-bases  over $\Spf \W$, 
for any objects $\E ^{(\bullet)}\in \mathfrak{C} (\fP)$, 
$\FF ^{ (\bullet)}
\in
\mathfrak{C} (\fQ)$,
we have 
$\left ( \R \underline{\Gamma} ^\dag _{Y}  \E ^{(\bullet)} \right) 
\smash{\widehat{\boxtimes}}^\L
_{\O  _{\Spf \W }} 
\FF^{ (\bullet) } 
\in
\mathfrak{C} (\fP\times _{\scr{C} _{\W}} \fQ)$.

\item \label{dfn-stable-data-h}
We will say that the  data of absolute coefficients $\mathfrak{C}$ is stable under weak admissible local cohomological functors
(resp. under localizations outside a divisor,
resp. localizations outside a  weak admissible  divisor), if 
for any object $\cW$ of $\mathrm{DVR}  (\cV)$,  
for any  formal  $\Spf(\cW)$-scheme of formal finite type $\fP$ having locally finite $p$-bases  over $\Spf \cW$, 
for any object $\E ^{(\bullet)}$ of $\mathfrak{C} (\fP)$, 
for any weak admissible inclusion $(Y \subset \fP)$ 
(resp.  for any divisor $T$ of the special fiber of $\fP$, 
resp. for any divisor $T$ of the special fiber of $\fP$ such that $(T \subset \fP)$ is admissible , 
we have
$\R \underline{\Gamma} ^\dag _{Y} \E ^{(\bullet)} \in \mathfrak{C} (\fP)$ 
(resp. $(\hdag T) (\E ^{(\bullet)} )\in \mathfrak{C} (\fP)$).

\item We will say that the  data of absolute coefficients $\mathfrak{C}$ is stable under cohomology if,
for any object $\cW$ of $\mathrm{DVR}  (\cV)$,  
for any  formal  $\Spf(\cW)$-scheme of formal finite type $\fX$ having locally finite $p$-bases  over $\Spf \cW$,
for any object  $\E ^{(\bullet)}$ of $\smash{\underrightarrow{LD}} ^{\mathrm{b}} _{\Q,\mathrm{coh}} ( \smash{\widehat{\D}} _{\X} ^{(\bullet)})$, the property 
$\E ^{(\bullet)}$ is an object of 
$\mathfrak{C} (\X)$ is equivalent to the fact that, for any integer $n$, 
$\H ^n (\E ^{(\bullet)})$ is an object of $\mathfrak{C} (\X)$. 

\item We will say that the  data of absolute coefficients $\mathfrak{C}$ is stable under shifts if,
for any object $\cW$ of $\mathrm{DVR}  (\cV)$,  
for any  formal  $\Spf(\cW)$-scheme of formal finite type $\fX$ having locally finite $p$-bases  over $\Spf \cW$,
for any object  $\E ^{(\bullet)}$ of $\mathfrak{C} (\X)$, for any integer $n$, 
$\E ^{(\bullet)} [n]$ is an object of $\mathfrak{C} (\X)$.

\item We will say that the  data of absolute coefficients $\mathfrak{C}$ is stable by devissages if
 $\mathfrak{C}$ is stable by shifts and if
for any object $\cW$ of $\mathrm{DVR}  (\cV)$,  
for any  formal  $\Spf(\cW)$-scheme of formal finite type $\fX$ having locally finite $p$-bases  over $\Spf \cW$, 
for any exact triangle 
$\E ^{(\bullet)} _1
\to 
\E ^{(\bullet)} _2
\to 
\E ^{(\bullet)} _3
\to 
\E ^{(\bullet)} _1 [1]$
of 
$\smash{\underrightarrow{LD}} ^{\mathrm{b}} _{\Q,\mathrm{coh}} ( \smash{\widehat{\D}} _{\X} ^{(\bullet)})$, 
if two objects are in $\mathfrak{C} (\X)$, then so is the third one. 

\item We will say that the  data of absolute coefficients $\mathfrak{C}$ is stable under direct summands if,
for any object $\cW$ of $\mathrm{DVR}  (\cV)$,  
for any  formal  $\Spf(\cW)$-scheme of formal finite type, having locally finite $p$-bases  over $\Spf \cW$ we have the following property: 
any direct summand in $\smash{\underrightarrow{LD}} ^{\mathrm{b}} _{\Q,\mathrm{coh}} ( \smash{\widehat{\D}} _{\X} ^{(\bullet)})$
of an object of $\mathfrak{C} (\X)$ is an object of
$\mathfrak{C} (\X)$.

\item We say that $\mathfrak{C}$ contains $\mathfrak{D}$ 
(or $\mathfrak{D}$ is contained in $\mathfrak{C}$) 
if for any object $\cW$ of $\mathrm{DVR}  (\cV)$,  
for any  formal  $\Spf(\cW)$-scheme of formal finite type, having locally finite $p$-bases  over $\Spf \cW$
the category $\mathfrak{D} (\X)$ is a full subcategory of $\mathfrak{C} (\X)$.

\item We say that the  data of absolute coefficients $\mathfrak{C}$ is local 
if for any object $\cW$ of $\mathrm{DVR}  (\cV)$,  
for any  formal  $\Spf(\cW)$-scheme of formal finite type $\fX$ having locally finite $p$-bases  over $\Spf \cW$, 
for any open covering $(\X _i) _{i\in I}$ of $\X$, 
for any object $\E ^{(\bullet)}$ of
$\smash{\underrightarrow{LD}} ^{\mathrm{b}} _{\Q,\mathrm{qc}} ( \smash{\widehat{\D}} _{\X} ^{(\bullet)})$, 
we have 
$\E ^{(\bullet)}\in \mathrm{Ob} \mathfrak{C} (\X)$ if and only if 
$\E ^{(\bullet)}| \X _i \in \mathrm{Ob} \mathfrak{C} (\X _i)$ for any $i \in I$. 
For instance, it follows from 
\ref{thick-subcat}.\ref{thick-subcat2}
that 
the  data of absolute coefficients 
$\smash{\underrightarrow{LD}} ^{\mathrm{b}} _{\Q,\mathrm{coh}}$ is local. 

\item We say that the  data of absolute coefficients $\mathfrak{C}$ is quasi-local 
if for any object $\cW$ of $\mathrm{DVR}  (\cV)$,  
for any  formal  $\Spf(\cW)$-scheme of formal finite type $\fX$ having locally finite $p$-bases  over $\Spf \cW$, 
for any open immersion $j \colon \Y \hookrightarrow \X$ 
for any object 
$\E ^{(\bullet)}\in  \mathfrak{C} (\X)$, 
we have 
$ j ^{ !(\bullet)}\E ^{(\bullet)} \in  \mathfrak{C} (\Y)$.
\end{enumerate}

\end{dfn}

\begin{dfn}
\label{dfn-stable-data-restricted}
Let $\mathfrak{C}$ and $\mathfrak{D}$ be two restricted data of absolute coefficients over $\fS$. 
Then, we have the similar definition than in \ref{dfn-stable-data} : we have only to restrict to
formal  $\bbD ^1 _{\cW}$-scheme of finite type, having locally finite $p$-bases  over $\Spf \W$ 
and to morphisms of formal  $\bbD ^1 _{\cW}$-scheme of finite type, having locally finite $p$-bases  over $\Spf \W$.
In the definition \ref{dfn-stable-data}.\ref{dfn-stable-data-h}, since in the restricted context every subschemes are weak admissible,
then we can remove ``weak admissible'' in the definitions.
\end{dfn}

We finish the subsection with some notation.

\begin{empt}
[Duality]
\label{ntn-dual}
Let $\mathfrak{C}$ be a data (resp. a restricted data) of absolute coefficients. We define its dual (restricted) data of absolute coefficients 
$\mathfrak{C}  ^{\vee}$ as follows: 
for any object $\cW$ of $\mathrm{DVR}  (\cV)$, for any integer $r\geq 0$ (resp. $r =1$), 
for any   formal  $\bbD ^r _{\cW}$-scheme of finite type, having locally finite $p$-bases  over $\Spf \W$, 
the category $\mathfrak{C}  ^{\vee} (\X)$ is the subcategory of 
$\smash{\underrightarrow{LD}} ^{\mathrm{b}} _{\Q,\mathrm{coh}} ( \smash{\widehat{\D}} _{\X} ^{(\bullet)})$
of objects $\E^{(\bullet) }$
such that $\DD _{\X} (\E^{(\bullet) }) \in \mathfrak{C} (\X)$.
\end{empt}

\begin{ntn}
\label{dfn-Delta(C)}
Let $\mathfrak{C}$ be a  (restricted) data of absolute coefficients.
We denote by $\mathfrak{C} ^+$ 
the smallest (restricted) data of absolute coefficients containing $\fC$ and stable under shifts.
We define by induction on $n\in \N$ the (restricted) data of absolute coefficients
$\Delta _n (\mathfrak{C} )$  as follows: 
for $n=0$, we put $\Delta _0 (\mathfrak{C} )= \mathfrak{C} ^+$.
Suppose $\Delta _n (\mathfrak{C} )$ constructed for $n\in\N$. 
for any object $\cW$ of $\mathrm{DVR}  (\cV)$, for any integer $r\geq 0$ (resp. for $r =1$), 
for any  formal  $\bbD ^r _{\cW}$-scheme of finite type $\fX$ having locally finite $p$-bases  over $\Spf \cW$,  
the category 
$\Delta _{n+1} (\mathfrak{C} ) (\X)$
is the full subcategory of 
$\smash{\underrightarrow{LD}} ^{\mathrm{b}} _{\Q,\mathrm{coh}} ( \smash{\widehat{\D}} _{\X} ^{(\bullet)})$
of objects $\E^{(\bullet) }$
such that 
there exists an exact triangle of the form 
$\E^{(\bullet) } \to \FF^{(\bullet) } \to \G ^{(\bullet) } \to \E^{(\bullet) } [1]$ 
such that 
$\FF^{(\bullet) } $ 
and
$\G ^{(\bullet) } $ 
are objects of 
$\Delta _{n} (\mathfrak{C} ) (\X)$.
Finally, we put 
$\Delta  (\mathfrak{C} ): = \cup _{n\in\N} \Delta _n (\mathfrak{C} )$. 
The (restricted) data of absolute coefficients 
$\Delta  (\mathfrak{C} )$ is the smallest (restricted) data of absolute coefficients
containing $\mathfrak{C}$ and stable under devissage. 
\end{ntn}

\begin{ex}
\label{stab-cst}
\begin{enumerate}[(a)]
\item 
\label{stab-cst-1}
Thanks to \ref{O-ovcoh}, 
using the isomorphisms \ref{fonctYY'Gamma-iso}, \ref{pre-loc-tri-B-t1T-isoovcoh} and Theorem \ref{2.2.18}, we check 
that $\fB _\mathrm{cst} ^{(1)+}$ satisfies $BK _+$, and is stable under local cohomological functors, 
extraordinary pullbacks 
and tensor products.

\item 
\label{stab-cst-2}
Following 
\ref{u+closedimm-YGamma},
\ref{u!-YGamma},
\ref{fonctYY'Gamma-isoqc},
we check 
that $\fB _\mathrm{wa} ^+$ 
satisfies $BK _+$, and is stable under weak admissible local cohomological functors, 
quasi-projective extraordinary pullbacks 
and tensor products.

\end{enumerate}
\end{ex}

The following lemma is obvious.
\begin{lemm}
\label{Delta-lemm-stab}
Let $\mathfrak{D}$ be a (restricted) data of absolute coefficients over $\fS$.

\begin{enumerate}
\item Let $P$ be one of the  stability property of \ref{dfn-stable-data} which is neither
the stability under cohomology, nor the stability under direct summands, nor the $BK ^!$ property, 
nor the localness. 
If $\fC$  satisfies the property $P$
then so is 
$\Delta(\mathfrak{D})$.

\item 
If $\mathfrak{D}$ satisfies $BK ^!$ and is stable under
local cohomological functors then so is 
$\Delta(\mathfrak{D})$.
\end{enumerate}

\end{lemm}

\begin{empt}
Beware also that if $\mathfrak{D}$ is local (resp. stable under cohomology, 
resp. stable under direct summands, 
resp. satisfies $BK ^!$), then it is not clear that 
so is $\Delta(\mathfrak{D})$.

\end{empt}

Since the converse of \ref{Delta-lemm-stab} is not true, let us introduce the following definition.
\begin{dfn}
\label{Delta-lemm-stab-bis}
Let $\mathfrak{D}$ be a (restricted) data of absolute coefficients over $\fS$. 
Let $P$ be one of the stability property of 
\ref{dfn-stable-data}.
We say that $\mathfrak{D}$ is 
$\Delta$-stable under $P$ (or satisfies the property $\Delta$-P) if there exists
a (restricted) data of absolute coefficients $\fD'$ over $\fS$ 
such that 
$\Delta(\fD ') =\Delta(\mathfrak{D})$  and $\fD '$ is stable under $P$.

Suppose $P$ is one of the  stability property of \ref{dfn-stable-data} which is neither
the stability under cohomology, nor the stability under direct summands, nor the $BK ^!$ property, 
nor the localness. 
A translation of Lemma \ref{Delta-lemm-stab} is the following : 
the data of coefficients $\mathfrak{D}$ is 
$\Delta$-stable under $P$
if and only if 
$\Delta(\mathfrak{D})$
is stable under $P$. 

Beware, it is not clear that if $\mathfrak{D}$ satisfies $\Delta$-$BK ^!$ and is $\Delta$-stable under
local cohomological functors then 
$\Delta(\mathfrak{D})$ satisfies $BK ^!$.

\end{dfn}

\subsection{Overcoherence, (over)holonomicity (after any base change)}

In this subsection, we explain how to get stable restricted 
data of absolute coefficients (see \ref{ovcoh-invim-prop}) which will be later the main ingredient of Theorem
\ref{dfnquprop}. For completeness, we extend some preliminary lemmas in the case 
of data of absolute coefficients when it is valid.

\begin{dfn}
\label{dfnS(D,C)}
Let $\mathfrak{C}$ and $\mathfrak{D}$ be two (restricted) data of absolute coefficients.

\begin{enumerate}[(a)]
\item We denote by 
$S _0 (\mathfrak{D}, \mathfrak{C})$
the (restricted) data of absolute coefficients defined as follows: 
for any object $\cW$ of $\mathrm{DVR}  (\cV)$, for any integer $r\geq 0$ (resp. for $r =1$), 
for any  formal  $\bbD ^r _{\cW}$-scheme of finite type $\fX$ having locally finite $p$-bases  over $\Spf \cW$,  
the category 
$S _0 (\mathfrak{D}, \mathfrak{C}) (\X)$
is the full subcategory of 
$\smash{\underrightarrow{LD}} ^{\mathrm{b}} _{\Q,\mathrm{coh}} ( \smash{\widehat{\D}} _{\X} ^{(\bullet)})$
of objects  $\E ^{(\bullet)}$
satisfying the following properties :
\begin{enumerate}[(a)]
\item [($\star$)] 
if for any smooth formal $\fS$-scheme $\fZ$, 
denoting by $\fY := \fX \times _{\fS} \fZ$ and
by $f \colon \fY \to \fX$ the projection, 
for any object 
$\FF ^{(\bullet)}
\in
\mathfrak{D} (\Y)$,
we have 
$\FF ^{(\bullet)}
\smash{\widehat{\otimes}}^\L
_{\O  _{\Y}} f ^{!} (\E^{(\bullet) })
\in
\mathfrak{C} (\Y)$.
\end{enumerate}

\item We denote by 
$S(\mathfrak{D}, \mathfrak{C})$ 
the (restricted) data of absolute coefficients defined as follows: 
for any object $\cW$ of $\mathrm{DVR}  (\cV)$, for any integer $r\geq 0$ (resp. for $r =1$), 
for any  formal  $\bbD ^r _{\cW}$-scheme of finite type $\fX$ having locally finite $p$-bases  over $\Spf \cW$,  
the category 
$S  (\mathfrak{D}, \mathfrak{C}) (\X)$ 
is the full subcategory of 
$\smash{\underrightarrow{LD}} ^{\mathrm{b}} _{\Q,\mathrm{coh}} ( \smash{\widehat{\D}} _{\X} ^{(\bullet)})$
of objects  $\E ^{(\bullet)}$
satisfying the following 
property :
\begin{enumerate}[(a)]
\item [($\star \star$)] for any morphism 
$\W
\to 
\W '$ 
of $\mathrm{DVR}  (\V)$, 
with notation \ref{chg-base}, 
we have 
$$ \cO _{\bbD ^r _{\cW'}}  \smash{\widehat{\otimes}}^\L_{\cO _{\bbD ^r _{\cW}} }  \E^{(\bullet) }\in 
S _0 (\mathfrak{D}, \mathfrak{C}) (\X \times _{\bbD ^r _{\cW}} \bbD ^r _{\cW '}).$$

\end{enumerate}
 \item Let $\sharp$ be a symbol so that 
either $S _\sharp = S _0$
or
$S _\sharp = S $.

\end{enumerate}

\end{dfn}

\begin{exs}
\label{ex-cst-surcoh}
\begin{enumerate}[(a)]

\item 
We have 
$ \smash{\underrightarrow{LD}} ^{\mathrm{b}} _{\Q,\mathrm{ovcoh}}= 
S _0 (\fB _\mathrm{div}, \smash{\underrightarrow{LD}} ^{\mathrm{b}} _{\Q,\mathrm{coh}})$
(see the second example  of \ref{ex-Dcst}).
We get again the notion of overcoherence of \ref{dfn-ovch}.

We denote by 
$ \smash{\underrightarrow{LD}} ^{\mathrm{b}} _{\Q,\mathrm{oc}}= 
S (\fB _\mathrm{div}, \smash{\underrightarrow{LD}} ^{\mathrm{b}} _{\Q,\mathrm{coh}})$.
This notion is an analogue of that of overcoherence after any base change as defined in
\cite{surcoh-hol}.

\item \label{hstab} 
We put 
$\mathfrak{H} _0 :=S (\fB _\mathrm{div}, \smash{\underrightarrow{LD}} ^{\mathrm{b}} _{\Q,\mathrm{coh}})$ 
and 
by induction on $i \in \N$, 
we put $\mathfrak{H} _{i+1} :=
\mathfrak{H} _{i} \cap S(\fB _\mathrm{div}, \mathfrak{H} _{i} ^{\vee})$
(see Notation \ref{ntn-dual}).
The absolute coefficients of 
$\mathfrak{H} _{i}$ are called 
{\it $i$-overholonomic after any base change}.
We get the (restricted) data of absolute coefficients 
$\smash{\underrightarrow{LD}} ^{\mathrm{b}} _{\Q,\mathrm{h}}:= \mathfrak{H} _{\infty}
:=  \cap _{i\in \N} \mathfrak{H} _{i}$
whose objects are called {\it overholonomic after any base change}.

\item \label{ovholstab}
Replacing $S $ by $S _0$ in the definition of $\smash{\underrightarrow{LD}} ^{\mathrm{b}} _{\Q,\mathrm{h}}$, 
we get a (restricted) data of absolute coefficients that we will denote by 
$\smash{\underrightarrow{LD}} ^{\mathrm{b}} _{\Q,\mathrm{ovhol}}$.

\item Finally, 
we set
$ \smash{\underrightarrow{LM}}  _{\Q,\star}:= 
 \smash{\underrightarrow{LD}} ^{\mathrm{b}} _{\Q,\star} \cap 
  \smash{\underrightarrow{LM}} _{\Q,\mathrm{coh}}$,
  for $\star \in \{\mathrm{ovcoh}, \mathrm{oc}, \mathrm{h}, \mathrm{ovhol} \}$.
\end{enumerate}

\end{exs}

\begin{rem}
\label{rem-overhol}
\begin{enumerate}[(a)]

\item Let $\mathfrak{C}$ be a (restricted) data of absolute coefficients.
The (restricted) data of absolute coefficients $\mathfrak{C}$ is stable under extraordinary pullbacks by  smooth projections, localizations outside a divisor
(resp. under extraordinary pullbacks by  smooth projections, localizations outside a divisor, and base change) if and only if 
$S _0 (\fB _\mathrm{div}, \mathfrak{C})=\mathfrak{C}$
(resp. $S  (\fB _\mathrm{div}, \mathfrak{C})=\mathfrak{C}$).

\item Let $\mathfrak{C}$ be a (restricted) data of absolute coefficients.
The (restricted) data of absolute coefficients $\mathfrak{C}$ is stable under extraordinary pullbacks by  smooth projections, 
weak admissible localizations
(resp. under extraordinary pullbacks by  smooth projections, weak admissible localizations, and base change) if and only if 
$S _0 (\fB _\mathrm{wa}, \mathfrak{C})\supset \mathfrak{C}$
(resp. $S  (\fB _\mathrm{wa}, \mathfrak{C})\supset \mathfrak{C}$).
Beware that the inclusion 
$S _0 (\fB _\mathrm{wa}, \mathfrak{C})\subset \mathfrak{C}$
(resp. $S  (\fB _\mathrm{wa}, \mathfrak{C})\subset \mathfrak{C}$)
is not clear (see \ref{preS(D,C)stability}.\ref{S(D,C)stability1}).

\item By construction, 
we remark that 
$\smash{\underrightarrow{LD}} ^{\mathrm{b}} _{\Q,\mathrm{ovhol}}$
is the biggest (restricted) data of absolute coefficients 
which contains 
$\fB _\mathrm{div}$, 
is stable by devissage, 
duality
and the operation
$S _{0} (\fB _\mathrm{div}, -)$. 
Moreover, 
$\smash{\underrightarrow{LD}} ^{\mathrm{b}} _{\Q,\mathrm{h}}$
is the biggest (restricted) data of absolute coefficients 
which contains 
$\fB _\mathrm{div}$, 
is stable by devissage, 
duality
and the operation
$S  (\fB _\mathrm{div}, -)$.

\end{enumerate}
\end{rem}

\

We will need later the following Lemmas.

\begin{lem}
\label{rem-div-cst}

We have the following properties. 

\begin{enumerate}[(a)]
\item 
\label{rem-div-cst1}
We have the equality 
$\Delta (\fB _\mathrm{div} ^{(1)})=\Delta (\fB _\mathrm{cst} ^{(1)})$ (see Notation \ref{ex-Dcst}).

\item \label{S(D,C)stability3-pre}
Let $\mathfrak{C}$ be a restricted data of absolute coefficients stable under devissage.
The following properties are equivalent :  
\begin{enumerate}[(a)]
\item $\mathfrak{C}$ is stable under local cohomological functors

\item $\mathfrak{C}$ is stable under localizations outside a divisor.

\end{enumerate}

\end{enumerate}
 
\end{lem}

\begin{proof}
Both  statements are checked by using exact triangles of localisation 
\ref{caro-stab-sys-ind-surcoh4.4.3}
and Mayer-Vietoris exact triangles \ref{eq1mayer-vietoris}.
\end{proof}

\begin{lem}
\label{lem-stabextpullback}
Let $\mathfrak{C}$ be a (restricted) data of absolute coefficients which is 
stable under local cohomological functors
(resp. weak admissible local cohomological functors)
and satisfies $BK ^!$.
Then  
$\mathfrak{C}$ is stable under 
extraordinary pullbacks by closed immersions
(resp. extraordinary pullbacks by  weak admissible closed immersions).

\end{lem}

\begin{proof}
Since the proof is the same, let us only check the non respective case. 
Let $\cW$ be an object of $\mathrm{DVR}  (\V)$, $r \geq 0$ be an integer, 
$f \colon \Y \to \X$ be a 
closed immersion 
of  formal  $\Spf (\cW)$-schemes
of formal finite type having locally finite $p$-bases  over $\Spf \cW$, 
and $\E ^{(\bullet)}$ be an object of $\mathfrak{C} (\X)$. 
We have to check 
$f ^{!(\bullet)} ( \E ^{(\bullet)}) \in \mathfrak{C} (\Y)$.
From the stability under local cohomological functors, 
$\R \underline{\Gamma} ^\dag _{Y} \E ^{(\bullet)} \in \mathfrak{C} (\X)$. 
Since 
 $\mathfrak{C}$ satisfies $BK ^!$,
 then 
$f ^{!(\bullet)}  \R \underline{\Gamma} ^\dag _{Y} \E ^{(\bullet)}
\in 
\mathfrak{C} (\Y)$.
We conclude using 
 the isomorphism 
$f ^{!(\bullet)}  \R \underline{\Gamma} ^\dag _{Y} \E ^{(\bullet)} \riso 
f ^{!(\bullet)}  ( \E ^{(\bullet)})$ 
(use \ref{2.2.18}).\end{proof}
\begin{rem}
\label{rem-qprojpullbacks}
The above lemma \ref{lem-stabextpullback} is important because of the following remark. 
Let $\mathfrak{C}$ be a (restricted) data of absolute coefficients which 
is quasi-local, stable under extraordinary pullbacks by closed immersions and by smooth projections. 
Then 
$\mathfrak{C}$ is a (restricted) data of absolute coefficients 
is stable under quasi-projective extraordinary pullbacks.
\end{rem}

\begin{lem}
\label{S(D,C)stability3bis}
Let $\mathfrak{D}$ be a (restricted) data of absolute coefficients over $\fS$.
If $\mathfrak{D}$ contains $\fB _\mathrm{div}$
(resp. $\fB _\mathrm{wa}$)
and if $\mathfrak{D}$ is stable under 
tensor products (resp. weak admissible tensor products),  
then 
$\mathfrak{D}$ is stable 
under localizations outside a divisor
(resp. weak admissible local cohomological functors).
\end{lem}

\begin{proof}
The non respective case is obvious.
The respective case 
is a consequence of the isomorphism 
\ref{fonctYY'Gamma-iso} (we use the case where 
$\E ^{ (\bullet)}= \O _{\X} ^{(\bullet)}$).
\end{proof}

\begin{lem}
\label{stab-Dvee-3prop}
Let $\mathfrak{C}$ be a (restricted) data of absolute coefficients. 
If the (restricted) data of absolute coefficients $\mathfrak{C}$ is local (resp. is stable under devissages, resp. is stable under direct summands,
resp. is stable under pushforwards, 
resp. is stable under base change,
resp. satisfies $BK ^!$),
then so is $\mathfrak{C} ^{\vee}$  (see Notation \ref{ntn-dual}).

\end{lem}

\begin{proof}
Thanks to Berthelot-Kashiwara theorem (see \ref{u!u+=id})
and to the relative duality isomorphism in the form of \ref{rel-dual-isom},
we can copy the proof \cite[11.2.7]{caro-6operations}.
\end{proof}

\begin{lem}
\label{stab-Dvee-3propbis}
Let $\mathfrak{C}$ and $\mathfrak{D}$ be two (restricted) data of absolute coefficients.
\begin{enumerate}[(a)]
\item If $\mathfrak{D} \subset \mathfrak{C} $ then 
$\mathfrak{D} ^\vee \subset \mathfrak{C} ^\vee$.
\item We have the equality
$\Delta (\fC) ^\vee =\Delta (\fC ^\vee)$. 
\end{enumerate}
\end{lem}

\begin{proof}
We can copy the proof \cite[11.2.8]{caro-6operations}.
\end{proof}

\begin{lem}
\label{preS(D,C)stability}
Let $\mathfrak{C}$ and $\mathfrak{D}$ be two data (resp. two restricted data) of absolute coefficients.
With the notation of \ref{dfnS(D,C)}, 
we have the following properties.

\begin{enumerate}[(a)]
\item 
\label{S(D,C)stability1}
With Notation \ref{ex-Dcst}, 
if $\mathfrak{D}$ contains $\fB _\emptyset$
(resp. if $\mathfrak{D}$ contains $\fB _\emptyset ^{(1)}$)
then  $S _\sharp (\mathfrak{D}, \mathfrak{C})$
 is contained in  $ \mathfrak{C}$.
 In the respective case, 
if $\mathfrak{D}$ contains $\fB _\mathrm{div} ^{(1)}$, then 
$S  _0(\mathfrak{D} , \mathfrak{C})$ 
is included in 
$\smash{\underrightarrow{LD}} ^{\mathrm{b} ~(1)} _{\Q,\mathrm{ovcoh}}$
and 
$S(\mathfrak{D}, \mathfrak{C})$
is included in 
$\smash{\underrightarrow{LD}} ^{\mathrm{b} ~(1)} _{\Q,\mathrm{oc}}$.

 \item 
\label{S(D,C)stability1bis}
If $\fC \subset \fC '$ and $\fD ' \subset \fD$, then 
$S _\sharp (\mathfrak{D}, \mathfrak{C}) \subset
S _\sharp (\mathfrak{D}', \mathfrak{C}') $.

 \item 
 \label{S(D,C)stabilitynew3}
 If either $\mathfrak{C}$ or $\mathfrak{D}$ is stable under devissages (resp. shifts),
then so is
$S _\sharp (\mathfrak{D}, \mathfrak{C})$
and we have the equality
$S _\sharp (\Delta  (\mathfrak{D} ), \mathfrak{C})
=
S _\sharp (\mathfrak{D}, \mathfrak{C})$
(resp.
$S _\sharp (\mathfrak{D}^+, \mathfrak{C})
=
S _\sharp (\mathfrak{D}, \mathfrak{C})$).

\item \label{S(D,C)stability2}
Suppose that $\mathfrak{D} $ is stable under extraordinary pullbacks by smooth projections, 
  tensor products (resp. and base change),
 and that $\mathfrak{C} $ contains $\mathfrak{D} $.
\begin{enumerate}[(a)]
\item The (restricted) data of absolute coefficients
$S _0  (\mathfrak{D}, \mathfrak{C})$
contains $\mathfrak{D}$
(resp. $S   (\mathfrak{D}, \mathfrak{C})$
contains $\mathfrak{D}$).

\item If $\mathfrak{D}$ contains $\fB _\emptyset$ (resp. $\fB _\emptyset ^{(1)}$), 
if either $\mathfrak{C} $ or $\mathfrak{D} $ is stable under shifts,
then 
$S _0 (\mathfrak{D}, \mathfrak{C})
=
S _0  \left (\mathfrak{D}, S _0 (\mathfrak{D}, \mathfrak{C}) \right)$ 
(resp. 
$S(\mathfrak{D}, \mathfrak{C})
=S \left (\mathfrak{D}, S(\mathfrak{D}, \mathfrak{C}) \right)$ ).

\item If either $\mathfrak{C} $ or $\mathfrak{D} $ is stable under shifts then
$S _0  \left (S _0 (\mathfrak{D}, \mathfrak{C}), S _0 (\mathfrak{D}, \mathfrak{C}) \right)$ 
(resp. $S \left (S(\mathfrak{D}, \mathfrak{C}), S(\mathfrak{D}, \mathfrak{C}) \right)$ )
contains $\mathfrak{D}$.

\end{enumerate}

\end{enumerate}
\end{lem}

\begin{proof}
Thanks to \ref{f!T'Totimes}, we can copy the proof of 
 \cite[11.2.9]{caro-6operations}.
\end{proof}

\begin{rem}
\label{rem-div-cst21bis} 
Let $\mathfrak{C}$,
$\mathfrak{D} $ be two (restricted) data of absolute coefficients.
Let $\fD'$ be a (restricted) data of absolute coefficients
such that 
$\Delta(\fD ') =\Delta(\mathfrak{D})$.
If $\fC$ is stable under devissages, then 
$S _\sharp (\mathfrak{D}', \mathfrak{C})
=
S _\sharp (\mathfrak{D}, \mathfrak{C})$.
Hence, in the case of stable properties appearing in Lemma \ref{Delta-lemm-stab-bis}
and when $\fC$ is stable under devissages,
to study 
$S _\sharp (\mathfrak{D}, \mathfrak{C})$ 
 it is enough to consider 
$\Delta$-stable properties  instead of stable properties satisfied by 
$\mathfrak{D}$ (e.g. see the beginning of the proof of \ref{ovcoh-invim-prop}).

\end{rem}

\begin{rem}
\label{rem-div-cst2}
Let $\mathfrak{C}$,
$\mathfrak{D} $ be two restricted data of absolute coefficients.

\begin{enumerate}[(a)]
\item 
\label{rem-div-cst21}
If $\fC$ is stable under devissages, then 
using \ref{preS(D,C)stability}.\ref{S(D,C)stabilitynew3}
and \ref{rem-div-cst}
we get
$S _\sharp (\fB _\mathrm{div} ^{(1)},\mathfrak{C}) 
=
S _\sharp (\fB _\mathrm{cst} ^{(1)+},\mathfrak{C}) 
$.

\item \label{rem-div-cst22} If
$\mathfrak{D} $ is stable under extraordinary pullbacks by smooth projections, 
  tensor products,
 and that $\mathfrak{D} $ contains $\fB _{\mathrm{div}} ^{(1)}$ and is contained in $\mathfrak{C} $,
if moreover 
either $\mathfrak{C} $ or $\mathfrak{D} $ is stable under shifts, 
then using \ref{preS(D,C)stability} (1, 2 and 4.b), we get
\begin{equation}
\label{rem-div-cst22eq1}S _0 (\mathfrak{D}, \mathfrak{C})
=
S _0  \left (\mathfrak{D}, S _0 (\fB _{\mathrm{div}} ^{(1)}, \mathfrak{C}) \right)
=
S _0  \left (\mathfrak{D}, S _0 (\mathfrak{D}, \mathfrak{C}) \right).
\end{equation}

 If moreover $\fD$ is stable under base change,
 then 
\begin{equation}
\label{rem-div-cst22eq2}
S  (\mathfrak{D}, \mathfrak{C})
=
S   \left (\mathfrak{D}, S  (\fB _{\mathrm{div}} ^{(1)}, \mathfrak{C}) \right)
=S   \left (\mathfrak{D}, S  (\mathfrak{D}, \mathfrak{C}) \right).
\end{equation}

\end{enumerate}
\end{rem}

\begin{lem}
\label{S(D,C)stability}
Let $\mathfrak{C}$ and $\mathfrak{D}$ be two (restricted) data of absolute coefficients.
We have the following properties.

\begin{enumerate}[(a)]

\item 
\label{S(D,C)stability3}
If  $\mathfrak{C}$ is local
and if $\mathfrak{D}$ is quasi-local
then 
$S _\sharp  (\mathfrak{D}, \mathfrak{C})$ is local. 
If  $\mathfrak{C}$ is stable under direct summands,
then so is
$S _\sharp  (\mathfrak{D}, \mathfrak{C})$.

\item \label{S(D,C)stability4}
The (restricted) data of absolute coefficients 
$S _0 (\mathfrak{D}, \mathfrak{C})$ 
(resp. $S(\mathfrak{D}, \mathfrak{C})$)
is stable under extraordinary pullbacks by smooth projections 
(resp. and under base change).

\item 
\label{S(D,C)stability5}
If 
$\mathfrak{D}$ is stable under weak admissible local cohomological functors (resp. localizations outside a weak admissible divisor), 
then so is
$S _\sharp (\mathfrak{D}, \mathfrak{C})$.

\item \label{S(D,C)stability6}
Suppose that $\mathfrak{C}$ 
is stable under pushforwards and shifts. Suppose that
$\mathfrak{D}$ 
is stable under quasi-projective extraordinary pullbacks.   
Then the (restricted) data of absolute coefficients 
$S _\sharp  (\mathfrak{D}, \mathfrak{C})$
 are stable under pushforwards.

\item \label{S(D,C)stability7}
Suppose that $\mathfrak{C}$ 
stable under shifts, and satisfies $BK ^!$.
Moreover, suppose that $\mathfrak{D}$ 
satisfies $BK _+$.
Then 
the (restricted) data of absolute coefficients 
$S _\sharp (\mathfrak{D}, \mathfrak{C})$
satisfies $BK ^!$.

\end{enumerate}
\end{lem}

\begin{proof}
Using \ref{thick-subcat}, 
\ref{fonctYY'Gamma-isoqc},
\ref{2.2.18qc},
\ref{subsect-comm-bc},
\ref{iso-chgtbase2},
\ref{surcoh2.1.4-iso}
we can copy the proof of 
 \cite[11.2.11.1--4]{caro-6operations}
 to check the first four statements. 
The check of the fifth one is very similar to that of  \cite[11.2.11.5]{caro-6operations}. 
For the reader, 
let us clarify it as follows. 
Since extraordinary pullbacks commute with base change, 
we reduce to check that 
$S _0(\mathfrak{D}, \mathfrak{C})$
satisfies $BK ^!$.
Let $\cW$ be an object of $\mathrm{DVR}  (\V)$, $r \geq 0$ be an integer, 
and $u\colon \X \hookrightarrow \fP$ be a closed immersion of  formal  $\Spf (\cW)$-schemes
of formal finite type having locally finite $p$-bases  over $\Spf \cW$. 
Let $\E ^{(\bullet)} \in S _0  (\mathfrak{D}, \mathfrak{C}) (\fP)$ with support in $\X$.
We have to check that $u ^! (\E ^{(\bullet)} ) \in S _0  (\mathfrak{D}, \mathfrak{C}) (\X)$.
We already know that
$u ^! (\E ^{(\bullet)})
\in 
\smash{\underrightarrow{LD}} ^{\mathrm{b}} _{\Q,\mathrm{coh}} ( \smash{\widehat{\D}} _{\X} ^{(\bullet)})$
(thanks to Berthelot-Kashiwara theorem
\ref{u!u+=id}).
Let $\fZ$ be a smooth formal $\fS$-scheme, 
$\fY := \fX \times _{\fS} \fZ$ and
$f \colon \fY \to \fX$ be the projection.
Let $\FF ^{(\bullet)}
\in
\mathfrak{D} (\Y)$.
We have to check 
$\FF ^{(\bullet)}
\smash{\widehat{\otimes}}^\L
_{\O  _{\Y}} f ^{!} (u ^! \E^{(\bullet) })
\in
\mathfrak{C} (\Y)$.
We denote by $v=id \times u 
\colon 
 \fX \times _{\fS} \fZ
 \hookrightarrow 
  \fP \times _{\fS} \fZ$.
Let $\fQ:=  \fP \times _{\fS} \fZ$
and $g \colon \fQ \to \fP$ be the projection. 
Since 
$\mathfrak{D} $ 
satisfies $BK _+$, 
then
$v _+ (\FF ^{(\bullet)}) \in \mathfrak{D} ( \fQ)$.
Since 
$\E ^{(\bullet)} \in S _0  (\mathfrak{D}, \mathfrak{C}) (\fP)$
and $g$ is a smooth projection morphism,
this yields 
$v _+ (\FF ^{(\bullet)})
\smash{\widehat{\otimes}}^\L
_{\O  _{\fQ}}  g ^{!} (\E^{(\bullet) })
\in \mathfrak{C} (\fQ)$.
Since $\mathfrak{C} $  satisfies $BK ^!$, 
this implies 
$v ^!\left (v _+ (\FF ^{(\bullet)})
\smash{\widehat{\otimes}}^\L
_{\O  _{\fQ}} g ^{!} (\E^{(\bullet) }) \right ) 
\in \mathfrak{C} (\Y)$.
Since 
$v ^!\left (v _+ (\FF ^{(\bullet)})
\smash{\widehat{\otimes}}^\L
_{\O  _{\fQ}}   g ^{!} (\E^{(\bullet) }) \right ) 
\riso 
v ^! v _+ (\FF ^{(\bullet)})
\smash{\widehat{\otimes}}^\L
_{\O  _{\Y}} 
v ^! g ^{!} (\E^{(\bullet) }) [r]$ with $r$ an integer (see \ref{f!T'Totimes}), 
since 
$v ^! v _+ (\FF ^{(\bullet)}) \riso 
\FF ^{(\bullet)}$
(see Berthelot-Kashiwara theorem
\ref{u!u+=id}),
since $\mathfrak{C}$ is stable under shifts, 
since by transitivity
$v ^! g ^{!}  \riso f ^! u ^!$,
we get
$\FF ^{(\bullet)}
\smash{\widehat{\otimes}}^\L
_{\O  _{\Y}} 
f ^! u ^! (\E^{(\bullet) })
\in 
\mathfrak{C} (\Y)$.
\end{proof}

\begin{prop}
\label{ovcoh-invim-prop}
Let $\mathfrak{C}$ and $\mathfrak{D}$ be two restricted data of absolute coefficients satisfying the 
following properties:
\begin{enumerate}[(a)]
\item We suppose either 
\begin{enumerate}[(i)]
\item $\fD$ contains $\fB _\mathrm{div} ^{(1)}$, 
satisfies $\Delta$-$BK _+$, and is $\Delta$-stable under 
quasi-projective extraordinary pullbacks
and tensor products

\item or $\fD$ contains $\fB _\emptyset ^{(1)}$, 
satisfies $\Delta$-$BK _+$, and is $\Delta$-stable under 
extraordinary pullbacks 
and
local cohomological functors.
\end{enumerate}

 \item  We suppose $\fC$ is local,
satisfies $BK ^!$, is stable under devissages, pushforwards, and direct summands.

\end{enumerate}

In both cases, the restricted data of absolute coefficients 
$S _0 (\mathfrak{D}, \mathfrak{C})$ 
(resp. $S(\mathfrak{D}, \mathfrak{C})$)
is local, stable under devissages, direct summands, 
local cohomological functors,
quasi-projective extraordinary pullbacks, pushforwards (resp. and base change).

\end{prop}

\begin{proof}
We can copy the proof of 
 \cite[11.2.12]{caro-6operations} (use also the remark \ref{rem-qprojpullbacks}).
\end{proof}

\begin{coro}
\label{ovcoh-invim}
Let $i \in \N \cup \{ \infty\}$.
The restricted data of absolute coefficients
$ \smash{\underrightarrow{LD}} ^{\mathrm{b} (1)} _{\Q,\mathrm{ovcoh}}$
(resp. $ \smash{\underrightarrow{LD}} ^{\mathrm{b} (1)} _{\Q,\mathrm{oc}}$,
resp. 
$\mathfrak{H} _{i}$)
contains
$\fB _\mathrm{cst} ^{(1)}$,
is local, stable under devissages, direct summands, 
local cohomological functors,
quasi-projective extraordinary pullbacks, pushforwards (resp. and base change).
Moreover, 
$\smash{\underrightarrow{LD}} ^{\mathrm{b} (1)} _{\Q,\mathrm{h}}$
is stable under duality.
\end{coro}

\begin{proof}
We can copy the proof of 
 \cite[11.2.13]{caro-6operations}.
\end{proof}

\subsection{On the stability under cohomology}
\begin{ntn}
\label{ntnC0}
Let $\mathfrak{C}$ be a (restricted) data of absolute coefficients.
We denote by $\mathfrak{C} ^0 $ the (restricted) data of absolute coefficients
defined as follows. 
Let $\cW$ be an object of $\mathrm{DVR}  (\V)$, $\fX$ be a formal  $\Spf (\cW)$-scheme of formal finite type, having locally finite $p$-bases  over $\Spf \cW$.
Then we set
$\mathfrak{C} ^0 (\X)
:=
\mathfrak{C}(\X)
\cap
\smash{\underrightarrow{LM}}  _{\Q,\mathrm{coh}} 
( \smash{\widehat{\D}} _{\X} ^{(\bullet)})
$.
\end{ntn}

\begin{lem}
\label{ntnC0-lemm}
Let $\mathfrak{C}$ be a (restricted) data of absolute coefficients.
Let $\cW$ be an object of $\mathrm{DVR}  (\V)$, 
$\X$ be a formal  $\Spf (\cW)$-scheme of formal finite type, having locally finite $p$-bases  over $\Spf \cW$.
\begin{enumerate}[(a)]
\item 
\label{ntnC0-lemm1}
If $\mathfrak{C}$ is stable under cohomology, then 
$\Delta (\mathfrak{C})
= 
\Delta (\mathfrak{C} ^0) $.
\item 
\label{ntnC0-lemm2}
If $\mathfrak{C}$ is stable under devissages and cohomology,
then the category $\mathfrak{C} ^0 (\X)$
is an abelian strictly full subcategory of 
$\smash{\underrightarrow{LM}}  _{\Q,\mathrm{coh}} 
( \smash{\widehat{\D}} _{\X} ^{(\bullet)})$ which is 
stable under extensions.

\end{enumerate}

\end{lem}

\begin{proof}
We can copy the proof of 
 \cite[11.2.15]{caro-6operations}.
\end{proof}

\begin{prop}
\label{Bcst-st-cohom}
Let $\mathfrak{C}$ 
be a data (resp. restricted data) of absolute coefficients
which is stable under cohomology, and devissage.
Then $S _\sharp (\fB _\mathrm{div} , \mathfrak{C})$ 
(resp. $S _\sharp (\fB _\mathrm{cst}^{(1)+} , \mathfrak{C})$) is stable under devissages  and cohomology.
\end{prop}

\begin{proof}
Since localizations outside a divisor and the functor $f ^{(\bullet)*}$ when $f$ is any smooth morphism
are t-exact (for the canonical t-structure of $\smash{\underrightarrow{LD}}  ^\mathrm{b} _{\Q, \mathrm{coh}}$), 
then the non respective case is straightforward.
Following \ref{rem-div-cst2},
$S _\sharp (\fB _\mathrm{cst} ^{(1)+}, \mathfrak{C})
=
S _\sharp (\fB _\mathrm{div} ^{(1)}, \mathfrak{C})$,
which yields the respective case.
\end{proof}

\begin{coro}
\label{coro-ovcoh-oc-tstr}
The restricted data of absolute coefficients 
$\smash{\underrightarrow{LD}} ^{\mathrm{b}(1)} _{\Q,\mathrm{ovcoh}}$,
and 
$\smash{\underrightarrow{LD}} ^{\mathrm{b}(1)} _{\Q,\mathrm{oc}}$ 
are stable under cohomology.
\end{coro}

\begin{empt}
Let $\cW$ be an object of $\mathrm{DVR}  (\V)$, $\fX$ be a formal  $\Spf (\cW)$-scheme of formal finite type, having locally finite $p$-bases  over $\Spf \cW$, 
$\E ^{(\bullet)}\in 
\smash{\underrightarrow{LM}}  _{\Q,\mathrm{coh}} 
( \smash{\widehat{\D}} _{\X} ^{(\bullet)})$.
Following \ref{ntn-dualfunctor}, 
we have the dual functor
$\DD ^{(\bullet)}
\colon 
\smash{\underrightarrow{LD}} ^{\mathrm{b}} _{\Q,\mathrm{coh}}
( \widehat{\D} _{\X /\fS } ^{(\bullet)})
\to 
\smash{\underrightarrow{LD}} ^{\mathrm{b}} _{\Q,\mathrm{coh}}
( \widehat{\D} _{\X /\fS } ^{(\bullet)})$.
Similarly to \cite[2.8]{caro-holo-sansFrob},
we say that 
$\E ^{(\bullet)}$ is holonomic if 
for any $i \not = 0$, 
$\H ^i (\DD ^{(\bullet)}(\E ^{(\bullet)})) =0$.
We denote by 
$\smash{\underrightarrow{LM}}  _{\Q,\mathrm{hol}} ( \widehat{\D} _{\X /\fS } ^{(\bullet)})$ the strictly subcategory 
of 
$\smash{\underrightarrow{LM}}  _{\Q,\mathrm{coh}} ( \widehat{\D} _{\X /\fS } ^{(\bullet)})$
of holonomic 
$\smash{\widehat{\D}} _{\X} ^{(\bullet)}$-modules.
By copying \cite[2.14]{caro-holo-sansFrob}, we check 
$\smash{\underrightarrow{LM}}  _{\Q,\mathrm{hol}} ( \widehat{\D} _{\X /\fS } ^{(\bullet)})$ 
is in fact a Serre subcategory 
of 
$\smash{\underrightarrow{LM}}  _{\Q,\mathrm{coh}} ( \widehat{\D} _{\X /\fS } ^{(\bullet)})$.

We denote by 
$\smash{\underrightarrow{LD}} ^{\mathrm{b}} _{\Q,\mathrm{hol}}
( \widehat{\D} _{\X /\fS } ^{(\bullet)})$ the strictly full subcategory of
$\smash{\underrightarrow{LD}} ^{\mathrm{b}} _{\Q,\mathrm{coh}}
( \widehat{\D} _{\X /\fS } ^{(\bullet)})  $ 
consisting of complexes 
$\E ^{(\bullet)}$ such that 
$\H ^n \E ^{(\bullet)} \in 
\smash{\underrightarrow{LM}}  _{\Q,\mathrm{hol}} ( \widehat{\D} _{\X /\fS } ^{(\bullet)})$
for any 
$n\in\Z$.
This yields the t-exact equivalence of categories
$\DD ^{(\bullet)}
\colon 
\smash{\underrightarrow{LD}} ^{\mathrm{b}} _{\Q,\mathrm{hol}}
( \widehat{\D} _{\X /\fS } ^{(\bullet)})
\cong 
\smash{\underrightarrow{LD}} ^{\mathrm{b}} _{\Q,\mathrm{hol}}
( \widehat{\D} _{\X /\fS } ^{(\bullet)})$.

Warning:  since in  the proof of 
\cite[3.3.5]{surcoh-hol} we have used Kedlaya's paper \cite{Kedlaya-coveraffinebis},
then this is not clear we have the inclusion
$\smash{\underrightarrow{LM}}  _{\Q,\mathrm{oc}} 
\subset
\smash{\underrightarrow{LM}}  _{\Q,\mathrm{hol}}$.
In particular the following inclusion is not clear
\begin{equation}
\label{ocinchol}
\smash{\underrightarrow{LD}} ^{\mathrm{b}} _{\Q,\mathrm{oc}} 
\subset
\smash{\underrightarrow{LD}} ^{\mathrm{b}}  _{\Q,\mathrm{hol}}.
\end{equation}
\end{empt}

\begin{empt}
\label{t-structure-ovcoh-oc-hol-h}
Let $\mathfrak{C}$ be a (restricted) data of absolute coefficients stable under devissages  and cohomology. 
Let $\cW$ be an object of $\mathrm{DVR}  (\V)$, $\fX$ be a formal  $\Spf (\cW)$-scheme of formal finite type, having locally finite $p$-bases  over $\Spf \cW$.
Recall that following 
\ref{t-structure-coh}
we have a canonical t-structure on 
$ \smash{\underrightarrow{LD}} ^{\mathrm{b}} _{\Q,\mathrm{coh}} 
( \smash{\widehat{\D}} _{\X} ^{(\bullet)})$.
We get a canonical t-structure on 
$\fC (\X/\W)$ whose heart is $\fC  ^0  (\X/\W)$ 
and so that the t-structure of 
$\fC (\X/\W)$ is induced 
by that of 
$ \smash{\underrightarrow{LD}} ^{\mathrm{b}} _{\Q,\mathrm{coh}} 
( \smash{\widehat{\D}} _{\X} ^{(\bullet)})$,
i.e. the truncation functors are the same
and 
$\fC  ^{\geq n}  (\X/\W):=
\smash{\underrightarrow{LD}} ^{\geq n} _{\Q,\mathrm{coh}} 
( \smash{\widehat{\D}} _{\X} ^{(\bullet)})
\cap 
\fC  (\X/\W)$,
$\fC  ^{\leq n}  (\X/\W):=
\smash{\underrightarrow{LD}} ^{\leq n} _{\Q,\mathrm{coh}} 
( \smash{\widehat{\D}} _{\X} ^{(\bullet)})
\cap 
\fC  (\X/\W)$.

For instance, using \ref{coro-ovcoh-oc-tstr}
we get 
  for $\star \in \{\mathrm{ovcoh}, \mathrm{oc}\}$
  a canonical t-structure on 
$\smash{\underrightarrow{LD}} ^{\mathrm{b}} _{\Q,\star} $.
The heart of 
$\smash{\underrightarrow{LD}} ^{\mathrm{b}} _{\Q,\star} $
is  $\smash{\underrightarrow{LM}} ^{\mathrm{b}} _{\Q,\star}$.
\end{empt}

\subsection{Constructions of stable restricted data of absolute coefficients}

\begin{dfn}
\label{dfn-almostdual}
Let $\mathfrak{D}$ be a (restricted) data of absolute coefficients over $\fS$.
We say that $\mathfrak{D}$ is ``almost stable under duality''
if the following property holds:
for any (restricted) data of absolute coefficients $\mathfrak{C}$ over $\fS$ which is local, 
stable under 
devissages, 
direct summands and 
pushforwards,
if $\mathfrak{D} \subset \mathfrak{C}$
then 
$\mathfrak{D} ^{\vee} \subset \mathfrak{C}$.
Remark from the biduality isomorphism that 
the inclusion 
$\mathfrak{D} ^{\vee} \subset \mathfrak{C}$
is equivalent to 
the following one 
$\mathfrak{D} \subset  \mathfrak{C} ^\vee$.
\end{dfn}

\begin{rem}
Compared to a  previous version of ``almost stability'' (see \cite{caro-6operations}),
we have added ``local'' in the hypotheses that $\fC$ have to satisfy.
This is 
because of the Zariski localness of 
the notion of ``nice divisor''. 
Hence, we can check that 
$\fM _\mathrm{n}$ is almost stable under stability (see
\ref{prop=div-almostst}) with our new notion.
\end{rem}

\begin{lem}
\label{almostdual-delta}
Let $\mathfrak{D}$ be a (restricted) data of absolute coefficients over $\fS$.
The (restricted) data $\mathfrak{D}$ is almost stable under duality
if and only if $\Delta (\mathfrak{D})$ is almost stable under duality.
\end{lem}
\begin{proof}
This is a consequence of \ref{stab-Dvee-3propbis}.
\end{proof}

\begin{lem}
\label{lem-prop=div-almostst}
With notation \ref{ex-Dcst},
we have the equalities
$\fM _\emptyset ^\vee
=
\fM _\emptyset $,
$(\Delta (\fM _\emptyset) ) ^\vee
=
\Delta (\fM _\emptyset)$
and 
$\Delta (\fM _\mathrm{sn})
= 
\Delta (\fM _\emptyset)$.
\end{lem}

\begin{proof}
The first equality is a consequence of
\ref{cor-com-sp+-f*}.
The second one follows from \ref{stab-Dvee-3propbis}.
It remains to  check the inclusion
$\fM _\mathrm{sn}
\subset \Delta (\fM _\emptyset)$.
Let $\cW$ be an object of $\mathrm{DVR}  (\V)$, $l$ be its residue field, 
let $\fX$ be a formal  $\Spf (\cW)$-scheme
of formal finite type having locally finite $p$-bases  over $\Spf \cW$,  
$Z$ be a closed subscheme of $X$ having locally finite $p$-bases over $\Spec l$, 
$T$ be a nice divisor of $Z/\Spec~ l$ 
and 
$\E ^{(\bullet)}
\in 
\mathrm{MIC} ^{(\bullet)} (Z, \X/K) $.
We have to prove that 
$(\hdag T ) (\E ^{(\bullet)})
\in 
\Delta (\fM _\emptyset) (\X)$.
We proceed by induction on the dimension of $T$ and next 
on the number of irreducible components of $T$.

Let $Z'$ be one irreducible component of $T$
and $T'$ be the union of the other irreducible components 
(hence $T  = Z' \cup T' $).
Then 
$T'\cap Z'$ is a strictly nice divisor of $Z'/\Spec l$. 
We have the localisation triangle
\begin{equation}
\label{coh-ss-div-bis-extri2}
(\hdag ~ T'\cap Z')   \R \underline{\Gamma} ^{\dag} _{Z'}   (  \E ^{(\bullet)})
\to 
 (\hdag T' )   (\E ^{(\bullet)})
\to 
 (\hdag T )   (\E ^{(\bullet)})
 \to +1.
 \end{equation}
Following \ref{cor-com-sp+-f*}, we have
$\R \underline{\Gamma} ^{\dag} _{Z'}   (  \E ^{(\bullet)}) [1]
\in 
\mathrm{MIC} ^{(\bullet)} (Z', \X/K) $.
Hence, since 
$T'\cap Z'$ is a strictly nice divisor of $Z'/\Spec l$,
by induction hypothesis we get
$(\hdag ~ T'\cap Z')   \R \underline{\Gamma} ^{\dag} _{Z'}   (  \E ^{(\bullet)})
\in
\Delta (\fM _\emptyset) (\X)$.
By induction hypothesis, we have also
$ (\hdag T')   (\E ^{(\bullet)}) \in \Delta (\fM _\emptyset) (\X)$.
Hence, by devissage, we get 
$ (\hdag T)   (\E ^{(\bullet)}) \in \Delta (\fM _\emptyset) (\X)$.
\end{proof}

\begin{prop}
\label{prop=div-almostst}
The (restricted) data of absolute coefficients $\fB _\mathrm{div} ^{(1)}$,
$\fB _\mathrm{cst} ^{(1)}$
and $\fM _\mathrm{n}$
are almost stable under duality.
\end{prop}

\begin{proof}
I) Since $\Delta (\fB ^{(1)} _\mathrm{cst} )= \Delta (\fB ^{(1)}_\mathrm{div})$ (see \ref{rem-div-cst}.\ref{rem-div-cst1}) and using \ref{almostdual-delta}, 
since the case $\fB _\mathrm{div} ^{(1)}$ is checked similarly, 
we reduce to prove the almost dual stability of 
$\fM _\mathrm{n}$.

II) Let $\mathfrak{C}$
be a restricted data of absolute coefficients  over $\fS$ which 
contains $\fM _\mathrm{n}$, and which 
is stable under 
devissages, 
direct summands and 
pushforwards.
Let $\cW$ be an object of $\mathrm{DVR}  (\V)$, 
$\fP$ be a formal  $\Spf (\cW)$-scheme
of formal finite type having locally finite $p$-bases  over $\Spf \cW$,  
$X$ be a closed subscheme of the special fiber of $\fP$
and having locally finite $p$-bases  over $\Spec l$, 
$T$ be a nice (see \ref{nice-div}) divisor of $X/S$,
and 
$\E ^{(\bullet)}
\in 
\mathrm{MIC} ^{(\bullet)} (X, \fP/K) $
be such that
$
(\hdag T ) (\E ^{(\bullet)})
\in 
\fC  (\X)$.
We have to check that 
$
(\hdag T ) (\E ^{(\bullet)})
\in 
\fC ^\vee (\X)$.
Since $\fC$ is local and is stable under pushforwards,
we can suppose that $X= P$
and we prefer to write $\fX$ instead of $\fP$. 
We can suppose there exists 
\begin{enumerate}[(a)]
\item a finite morphism $\V \to \V'$ of complete discrete valuation rings of mixed characteristics $(0,p)$, 

\item a finite morphism of formal schemes $\bbD ^r _{\fS '}\to \bbD ^r _{\fS}$ with $\fS ' := \Spf~\cV '$ making commutative the diagram 
$$
\xymatrix{
{\mathbb{D} ^r _{\fS'}} 
\ar[d] ^-{}
\ar[r] ^-{}
&
 {\mathbb{D} ^r _{\fS '}} 
\ar[d] ^-{}
 \\ 
 {\fS'} 
\ar[r] ^-{}
&
 {\fS,} 
 }
$$

\item a projective morphism $g\colon \fP' \to \fP$ of formal schemes
such that $\fP ' /\fS '$ has locally finite $p$-bases, 
a closed subscheme $X' \subset g _0 ^{-1} (X)$ of the special fiber $P'$
such that $X' /S$ has locally finite $p$-base
and the induced morphism $\phi \colon X ' \to X$ is 
an alteration of schemes (in the sense of \cite[2.20]{dejong})
and 
$T':= \phi  ^{-1} (T ) _\mathrm{red}$ is a strictly nice divisor of $X'/S'$
(see Definition \ref{st-nice-div}).
\end{enumerate}

1) Since $\mathfrak{C}$
is stable under 
devissages, 
direct summands and 
pushforwards, then using all the arguments of 
the step 1) of the proof of Proposition \ref{coh-cst-div},
we can suppose $S' = S$.

2) i) By copying 
the step 2 ) i) of the proof of Proposition \ref{coh-cst-div},
we get both morphisms by adjunction
$ f  ^{(\bullet)} _{+}  \R \underline{\Gamma} ^{\dag} _{X '}  f ^{!(\bullet)}(\E ^{(\bullet)} )
\overset{\rho _1}{\longrightarrow}
\E ^{(\bullet)}$
and 
$ f  ^{(\bullet)} _{+}  \R \underline{\Gamma} ^{\dag} _{X '}  f ^{!(\bullet)}(\DD ^{(\bullet)} \E ^{(\bullet)}  )
\overset{\rho '_2}{\longrightarrow}
\DD ^{(\bullet)} (\E ^{(\bullet)} )$.

ii) Next, 
we construct 
$\E ^{(\bullet)}
\overset{\rho _2}{\longrightarrow}
 f  ^{(\bullet)} _{+}  \R \underline{\Gamma} ^{\dag} _{X '}  f ^{!(\bullet)}(\E ^{(\bullet)} )$ by duality as follows:
$$\E ^{(\bullet)} 
\riso 
\DD ^{(\bullet)} \DD ^{(\bullet)} \E ^{(\bullet)} 
\overset{\DD (\rho '_2)}{\longrightarrow}
\DD  ^{(\bullet)}  f  ^{(\bullet)} _{+}  \R \underline{\Gamma} ^{\dag} _{X '}  f ^{!(\bullet)}
(\DD ^{(\bullet)} \E ^{(\bullet)} )
\underset{\ref{rel-dual-isom-proj-formal-iso}}{\riso}
f  ^{(\bullet)} _{+} 
\DD  ^{(\bullet)}   \R \underline{\Gamma} ^{\dag} _{X '}  f ^{!(\bullet)}(\DD ^{(\bullet)} \E ^{(\bullet)} ).$$
Following \ref{cor-com-sp+-f*2}, we have the following isomorphism
$\DD ^{(\bullet)} \R \underline{\Gamma} ^\dag _{X'} f ^{!(\bullet)} (\DD ^{(\bullet)} \E ^{(\bullet)} )
\riso 
\R \underline{\Gamma} ^\dag _{X'} f ^{!(\bullet)} \E ^{(\bullet)} $,
and we are done by composition.

3) By copying 
the step 2 ) iii) of the proof of Proposition \ref{coh-cst-div},
we check that $\rho _1 \circ \rho _2$ is an isomorphism. 
In particular, 
$\E ^{(\bullet)}$
is a direct summand of 
$ f  ^{(\bullet)} _{+}  \R \underline{\Gamma} ^{\dag} _{X '}  f ^{!(\bullet)}(\E ^{(\bullet)} )$.
Hence, 
$(\hdag T ) ( \E ^{(\bullet)})$
is a direct summand of 
$(\hdag T )  f  ^{(\bullet)} _{+}  \R \underline{\Gamma} ^{\dag} _{X '}  f ^{!(\bullet)}(\E ^{(\bullet)} )$
Using the commutation of localisation functor with pushforwards,
this yields
$(\hdag T )  (\E ^{(\bullet)})$
is a direct summand of 
$  f  ^{(\bullet)} _{+} (\hdag T ') \R \underline{\Gamma} ^{\dag} _{X '}  f ^{!(\bullet)}(\E ^{(\bullet)} )$.

4)
Since
$\E ^{\prime (\bullet)} := 
\R \underline{\Gamma} ^{\dag} _{X '}  f ^{!(\bullet)}(\E ^{(\bullet)} )
\in 
\mathrm{MIC} ^{(\bullet)} (X', \fP'/K) $ (use \ref{cor-com-sp+-f*}), 
$T'$ is a strictly nice divisor of $X'/S'$,
then
$(\hdag T ' )(\E ^{\prime (\bullet)} )
\in 
\fM _\mathrm{sn} (\X')$.
Since 
$\mathfrak{C}$ contains 
$\fM _\mathrm{sn}$ and 
is stable under 
devissages, 
then 
using \ref{lem-prop=div-almostst}
we get 
$\fM _\mathrm{sn} \subset 
\fC ^\vee$.
Hence, 
$(\hdag T ' )(\E ^{\prime (\bullet)} )
\in 
\fC ^\vee ( \X')$.
Since 
$\mathfrak{C}$ is stable under 
direct summands and 
pushforwards, we are done.
\end{proof}

\begin{ntn}
\label{dfnqupre}
Let $\mathfrak{C}, \mathfrak{D}$ be two restricted data of absolute coefficients.
We put 
$T _0 (\mathfrak{D} ,\mathfrak{C}) := 
S  (\mathfrak{D}  ,\mathfrak{C})$.
By induction on $i \in \N$, we set 
$U _i (\mathfrak{D}  ,\mathfrak{C}):= T _i (\mathfrak{D}  ,\mathfrak{C}) \cap T _i (\mathfrak{D}  ,\mathfrak{C} ) ^{\vee}$, 
$\widetilde{T} _{i} (\mathfrak{D}  ,\mathfrak{C}) := 
S (\mathfrak{D} , U _i (\mathfrak{D}  ,\mathfrak{C}))$
and
$T _{i+1} (\mathfrak{D}  ,\mathfrak{C}):= 
S  (\widetilde{T} _{i} (\mathfrak{D}  ,\mathfrak{C}), \widetilde{T} _{i} (\mathfrak{D}  ,\mathfrak{C}) )$.
We put $T (\mathfrak{D}  ,\mathfrak{C}) := \cap _{i\in \N} T _{i} (\mathfrak{D}  ,\mathfrak{C}) $.
\end{ntn}

\begin{thm}
\label{dfnquprop}
Let $\fB _\mathrm{div} ^{(1)} \subset \mathfrak{D}\subset \fC$ be two restricted data of absolute coefficients.
We suppose 
\begin{enumerate}[a)]
\item The restricted data $\fD$ is stable under 
extraordinary pullbacks by smooth projections ;
\item The restricted data $\Delta (\fD)$ satisfies $BK _+$, is stable under quasi-projective extraordinary pullbacks, base change, 
tensor products 
and  is almost stable under duality ;
\item The restricted data $\mathfrak{C}$ 
satisfies $BK ^!$, 
is local and stable under devissages, direct summands, 
pushforwards.
\end{enumerate}
Then, 
the restricted data of absolute coefficients $T(\mathfrak{D}, \mathfrak{C})$ 
(see Definition \ref{dfnqupre})
 is included in 
$\mathfrak{C} $,
contains 
$\mathfrak{D} $,
is 
local,
stable by devissages, direct summands, 
local cohomological functors, 
pushforwards, quasi-projective extraordinary pullbacks, base change, tensor products, duality.

\end{thm}

\begin{proof}
We can copy the proof of 
 \cite[11.6.6]{caro-6operations}.
\end{proof}

\begin{ex}
\label{ex-datastableevery}
We can choose 
$\fD = \fB _\mathrm{div} ^{(1)+}$ (or 
$\fD =\fM _\mathrm{n}$)
and 
$\fC = \smash{\underrightarrow{LD}} ^{\mathrm{b}(1)} _{\Q,\mathrm{coh}}$.
\end{ex}

\subsection{On the stability under external products}

In order to get some data stable under external products,
we need to extend the constructions 
of \ref{dfnS(D,C)} as follows. 
\begin{dfn}
\label{dfnS(D,C)+}
Let $\mathfrak{C}$ and $\mathfrak{D}$ be two  data of absolute coefficients.

\begin{enumerate}[(a)]
\item We denote by 
$S ^+ _0 (\mathfrak{D}, \mathfrak{C})$
the data of absolute coefficients defined as follows: 
for any object $\cW$ of $\mathrm{DVR}  (\cV)$, 
for any formal  $\Spf (\cW)$-scheme $\fX$ of formal finite type, having locally finite $p$-bases over $\Spf \cW$,
the category 
$S _0 (\mathfrak{D}, \mathfrak{C}) (\X)$
is the full subcategory of 
$\smash{\underrightarrow{LD}} ^{\mathrm{b}} _{\Q,\mathrm{coh}} ( \smash{\widehat{\D}} _{\X} ^{(\bullet)})$
of objects  $\E ^{(\bullet)}$
satisfying the following properties :
\begin{enumerate}[(a)]
\item [($\star$)] 
for any formal  $\Spf (\cW)$-scheme $\fY$ of formal finite type, having locally finite $p$-bases over $\Spf \cW$,
for any object 
$\FF ^{(\bullet)}
\in
\mathfrak{D} (\fX \times _{\scr{C} _\cW}\fY)$,
we have 
$\FF ^{(\bullet)}
\smash{\widehat{\otimes}}^\L
_{\O  _{\fX \times _{\scr{C} _\cW}\fY}} \varpi ^{!} (\E^{(\bullet) })
\in
\mathfrak{C} (\fX \times _{\scr{C} _\cW}\fY)$,
\end{enumerate}
where $\varpi \colon \fX \times _{\scr{C} _\cW}\fY \to \fX$ is the projection.

\item We denote by 
$S ^+(\mathfrak{D}, \mathfrak{C})$ 
the data of absolute coefficients defined as follows: 
for any object $\cW$ of $\mathrm{DVR}  (\cV)$, for any integer $r\geq 0$,
for any  formal  $\bbD ^r _{\cW}$-scheme of finite type $\fX$ having locally finite $p$-bases  over $\Spf \cW$,  
the category 
$S  (\mathfrak{D}, \mathfrak{C}) (\X)$ 
is the full subcategory of 
$\smash{\underrightarrow{LD}} ^{\mathrm{b}} _{\Q,\mathrm{coh}} ( \smash{\widehat{\D}} _{\X} ^{(\bullet)})$
of objects  $\E ^{(\bullet)}$
satisfying the following 
property :
\begin{enumerate}[(a)]
\item [($\star \star$)] for any morphism 
$\W
\to 
\W '$ 
of $\mathrm{DVR}  (\V)$, 
with notation \ref{chg-base}, 
we have 
$$ \cO _{\bbD ^r _{\cW'}}  \smash{\widehat{\otimes}}^\L_{\cO _{\bbD ^r _{\cW}} }  \E^{(\bullet) }\in 
S _0 (\mathfrak{D}, \mathfrak{C}) (\X \times _{\bbD ^r _{\cW}} \bbD ^r _{\cW '}).$$

\end{enumerate}
 \item Let $\sharp$ be a symbol so that 
either $S _\sharp = S _0$
or
$S _\sharp = S $.

\end{enumerate}
\end{dfn}

\begin{lem}
\label{S(D,C)stability+}
Let $\mathfrak{C}$ and $\mathfrak{D}$ be two  data of absolute coefficients.
We have the following properties.

\begin{enumerate}[(a)]

\item 
\label{S(D,C)stability3+}
If  $\mathfrak{C}$ is local
and if $\mathfrak{D}$ is quasi-local
then 
$S ^{+} _\sharp  (\mathfrak{D}, \mathfrak{C})$ is local. 
If  $\mathfrak{C}$ is stable under direct summands (resp. devissage),
then so is
$ S ^{+}  _\sharp  (\mathfrak{D}, \mathfrak{C})$.

\item \label{S(D,C)stability4+}
The  data of absolute coefficients 
$ S ^{+}  _0 (\mathfrak{D}, \mathfrak{C})$ 
(resp. $S ^{+}(\mathfrak{D}, \mathfrak{C})$)
is stable under extraordinary pullbacks by projections 
(resp. and under base change).

\item \label{S(D,C)stability6+}
Suppose that $\mathfrak{C}$ 
is stable under pushforwards and shifts. Suppose that
$\mathfrak{D}$ 
is stable under quasi-projective extraordinary pullbacks.   
Then the  data of absolute coefficients 
$ S ^{+}  _\sharp  (\mathfrak{D}, \mathfrak{C})$
 are stable under pushforwards.

\item \label{S(D,C)stability7+}
Suppose that $\mathfrak{C}$ 
is stable under shifts, and satisfies $BK ^!$.
Moreover, suppose that $\mathfrak{D}$ 
satisfies $BK _+$.
Then 
the  data of absolute coefficients 
$ S ^{+}  _\sharp (\mathfrak{D}, \mathfrak{C})$
satisfies $BK ^!$.

\end{enumerate}
\end{lem}

\begin{proof}
By using \ref{u+-com-otimes}, we check \ref{S(D,C)stability7+}.
We check the other assertions similarly to \ref{S(D,C)stability}.
\end{proof}

\begin{rem}
\label{rem-S+wa}
Let $\mathfrak{C}$ be a data of absolute coefficients which is stable under shifts.
It is not clear that the analogue of
\ref{S(D,C)stability}.\ref{S(D,C)stability5} is true. 
But we have the following remarks.
 
\begin{enumerate}[(a)]
\item 
\label{rem-S+wa-a)}
Since this is not clear 
that $\mathfrak{B} _{\emptyset}  \subset \mathfrak{B} _{\mathrm{wa}}$, then 
beware that the inclusion 
$ S ^+ (\mathfrak{B} _{\mathrm{wa}}  ,\fC) \subset \fC $
might be wrong.
But we have this property ``on weak admissible support'' (i.e. see \ref{rem-S+wa-d)}).

\item 
\label{rem-S+wa-b)}
The data $S ^+ (\mathfrak{B} _{\mathrm{wa}}  ,\fC) $ is stable under 
weak admissible cohomological functors
and under extraordinary pullbacks by projections.
Indeed, since a projection morphism is flat, 
this is a consequence of \ref{2.2.18qcgen}.

\item 
\label{rem-S+wa-c)}
If $\fC$ is stable under weak admissible cohomological functors
and under extraordinary pullbacks by projections,
then
we have the inclusion $\fC \subset S ^+ (\mathfrak{B} _{\mathrm{wa}}  ,\fC) $.
Beware that this is not clear in this case that the inclusion is an equality.

\item 
\label{rem-S+wa-d)}
For any object $\cW$ of $\mathrm{DVR}  (\cV)$, 
for any  formal  $\Spf (\cW)$-scheme of formal finite type $\fP$ having locally finite $p$-bases  over $\Spf \cW$,  
for any weak admissible inclusion $(Y \subset \fP)$,
for any object $\E ^{(\bullet)}$ of $S ^+ (\mathfrak{B} _{\mathrm{wa}}  ,\fC) (\fP)$,  we have
$\R \underline{\Gamma} ^\dag _{Y} \E ^{(\bullet)} \in \mathfrak{C} (\fP)$.

\item 
\label{rem-S+wa-e)}
We have the equality
\begin{equation}
\label{rem-S+wa-c)=}
S ^+ (\mathfrak{B} _{\mathrm{wa}}  ,S ^+ (\mathfrak{B} _{\mathrm{wa}}  ,\fC) ) 
=S ^+ (\mathfrak{B} _{\mathrm{wa}}  ,\fC) .
\end{equation}
Indeed, by using the above remarks \ref{rem-S+wa-b)} and \ref{rem-S+wa-c)}, 
we can check that the left term contain the right term. 
The reverse inclusion is a consequence of \ref{rem-S+wa-d)}.

\end{enumerate}

\end{rem}

\begin{lem}
\label{lemdfnqupropbis+mdiv}
Let $\mathfrak{C}$ be a data of absolute coefficients which contains $\mathfrak{B} _{\mathrm{wa}} $ and is stable under shifts.
We have the inclusions :
$\mathfrak{B} _{\emptyset} 
\subset
S ^+ (\mathfrak{B} _{\mathrm{wa}}  ,\fC)$
and 
$\mathfrak{B} _{wa} 
\subset
S ^+ (\mathfrak{B} _{\mathrm{wa}}  ,\fC)$.
\end{lem}

\begin{proof}
This is a consequence of 
Theorem \ref{2.2.18qcgen}.
\end{proof}

\begin{prop}
\label{waS(D,C)stability+}
Let $\mathfrak{C}$ be a data of absolute coefficients which is local, satisfies $BK ^!$, is stable under 
devissage, direct summands, 
pushforwards and which
contains $\mathfrak{B} _{\mathrm{wa}} $.
Then 
$ S ^{+}  (\mathfrak{B} _{\mathrm{wa}}  ,\fC)$
is local, 
satisfies $BK ^!$, is stable under base change, 
devissage, direct summands, 
weak admissible local cohomological functors, 
pushforwards, extraordinary pullbacks by projections and 
contains $\mathfrak{B} _{\emptyset} $.
\end{prop}

\begin{proof}
Following \ref{lemdfnqupropbis+mdiv}, 
$S ^+ (\mathfrak{B} _{\mathrm{wa}}  ,\fC)$
contains $\mathfrak{B} _{\emptyset} $.
The other properties follows from \ref{stab-cst}.\ref{stab-cst-2}, \ref{S(D,C)stability+} and \ref{rem-S+wa}.\ref{rem-S+wa-b)}.
\end{proof}

\begin{dfn}
\label{dfnSboxtimes(D,C)}
Let $\mathfrak{C}$ and $\mathfrak{D}$ be two  data of absolute coefficients.

\begin{enumerate}[(a)]
\item We denote by 
$\boxtimes _0 (\mathfrak{D}, \mathfrak{C})$
the  data of absolute coefficients defined as follows: 
for any object $\cW$ of $\mathrm{DVR}  (\cV)$, for any integer $r\geq 0$, 
for any  a formal  $\Spf (\cW)$-scheme of formal finite type $\fX$ having locally finite $p$-bases  over $\Spf \cW$,
the category 
$\boxtimes _0 (\mathfrak{D}, \mathfrak{C}) (\X)$
is the full subcategory of 
$\smash{\underrightarrow{LD}} ^{\mathrm{b}} _{\Q,\mathrm{coh}} ( \smash{\widehat{\D}} _{\X} ^{(\bullet)})$
consisting of objects  $\E ^{(\bullet)}$
satisfying the following property :
\begin{enumerate}[(a)]
\item [($\star$)]
for any formal  $\Spf (\cW)$-scheme $\Y$ of formal finite type, having locally finite $p$-bases over $\Spf \W$, 
for any object 
$\FF ^{ (\bullet)}
\in
\mathfrak{D} (\Y)$,
we have 
$\E ^{(\bullet)}
\smash{\widehat{\boxtimes}}^\L
_{\O  _{\Spf \W }} 
\FF^{ (\bullet) }
\in
\mathfrak{C} (\X\times _{\scr{C} _{\W}} \Y)$.
\end{enumerate}

\item We denote by 
$\boxtimes (\mathfrak{D}, \mathfrak{C})$
the  data of absolute coefficients defined as follows: 
for any object $\cW$ of $\mathrm{DVR}  (\cV)$, for any integer $r\geq 0$, 
for any  formal  $\bbD ^r _{\cW}$-scheme of finite type $\fX$ having locally finite $p$-bases  over $\Spf \cW$,  
the category 
$\boxtimes (\mathfrak{D}, \mathfrak{C}) (\X)$
is the full subcategory of 
$\smash{\underrightarrow{LD}} ^{\mathrm{b}} _{\Q,\mathrm{coh}} ( \smash{\widehat{\D}} _{\X} ^{(\bullet)})$
consisting of objects  $\E ^{(\bullet)}$
satisfying the following property :
\begin{enumerate}[(a)]
\item 
[($\star \star$)]
for any morphism 
$\W
\to 
\W '$ 
of $\mathrm{DVR}  (\V)$, 
$\cO _{\bbD ^r _{\cW'}}  \smash{\widehat{\otimes}}^\L_{\cO _{\bbD ^r _{\cW}} }  \E^{(\bullet) }
\in 
\boxtimes _0 (\mathfrak{D}, \mathfrak{C})
(\X \times _{\bbD ^r _{\cW}} \bbD ^r _{\cW '})$.
\end{enumerate}
 \item Let $\sharp$ be a symbol so that 
either $\boxtimes _\sharp = \boxtimes _0$
or
$\boxtimes _\sharp = \boxtimes $.

\end{enumerate}

\end{dfn}

\begin{lem}
\label{lem-boxtimesDC}
Let $\mathfrak{C}$ and $\mathfrak{D}$ be two  data of absolute coefficients.

\begin{enumerate}[(a)]
\item 
\label{lem-boxtimesDC0}
Suppose $\fD\subset \fC$.
If $\fD$ or $\fC$ is stable under extraordinary pullbacks by projections and shift, 
then $\fB _{\emptyset} \subset \boxtimes _\sharp(\mathfrak{D}, \mathfrak{C}) $.

\item 
\label{lem-boxtimesDC1}
Suppose 
for any object $\cW$ of $\mathrm{DVR}  (\cV)$, 
$\fB _\emptyset(\Spf (\cW) /\cW) \subset \mathfrak{D}(\Spf (\cW) /\cW)$. 
Then 
$\boxtimes _\sharp(\mathfrak{D}, \mathfrak{C})\subset \fC$.

 \item 
\label{lem-boxtimesDC2}
If $\fC \subset \fC '$ and $\fD ' \subset \fD$, then 
$\boxtimes _\sharp(\mathfrak{D}, \mathfrak{C}) \subset
\boxtimes _\sharp(\mathfrak{D}', \mathfrak{C}')$.

\item 
\label{lem-boxtimesDC3} 
If $\fC$ is stable under devissage
then so is 
$\boxtimes _\sharp(\mathfrak{D}, \mathfrak{C})$.
Moreover, 
$\boxtimes _\sharp(\mathfrak{D}, \mathfrak{C})=
\boxtimes _\sharp(\Delta(\mathfrak{D}), \mathfrak{C})$.

\item 
\label{lem-boxtimesDC4} 
If $\fC$ is stable under extraordinary pullbacks by projections, then so is 
$\boxtimes _\sharp(\mathfrak{D}, \mathfrak{C}) $.

\item The  data $\boxtimes (\mathfrak{D}, \mathfrak{C})$ is stable under base change.

\item If $\fC$ is stable under pushforwards (resp. satisfies  $BK ^!$, resp. is local, 
resp. is stable under direct summands), 
then so is 
$\boxtimes _\sharp(\mathfrak{D}, \mathfrak{C})$.

\end{enumerate}

\end{lem}

\begin{proof}
To check \ref{lem-boxtimesDC4}, we remark that 
for any object $\cW$ of $\mathrm{DVR}  (\cV)$, for any integers $r,s,y\geq 0$, 
for any  formal  $\bbD ^r _{\cW}$-scheme (resp. $\bbD ^s _{\cW}$-scheme,
resp. $\bbD ^u _{\cW}$-scheme) of finite type $\fX$ (resp. $\fY$, resp. $\fT$) having locally finite $p$-bases  over $\Spf \cW$,  
for any  $\E ^{(\bullet)}
\in \smash{\underrightarrow{LD}} ^{\mathrm{b}} _{\Q,\mathrm{coh}} ( \smash{\widehat{\D}} _{\X} ^{(\bullet)})$,
$\FF ^{ (\bullet)}
\in \smash{\underrightarrow{LD}} ^{\mathrm{b}} _{\Q,\mathrm{coh}} ( \smash{\widehat{\D}} _{\Y} ^{(\bullet)})$,
we have the formula
$$\varpi ^{! (\bullet)} 
\left ( 
\E ^{(\bullet)}
\smash{\widehat{\boxtimes}}^\L
_{\O  _{\Spf \W }} 
\FF^{ (\bullet) }
\right) 
\riso
\varpi ^{! (\bullet)}  (\E ^{(\bullet)}) 
\smash{\widehat{\boxtimes}}^\L
_{\O  _{\Spf \W }} 
\FF^{ (\bullet) }$$
where $\varpi $ is either the projection 
$
\X\times _{\scr{C} _{\W}} \Y \times _{\scr{C} _{\W}} \fT
\to 
\X\times _{\scr{C} _{\W}} \Y$ 
or 
$\X\times _{\scr{C} _{\W}}\fT
\to 
\X$.
Thanks to \ref{prop-boxtimes-v+}, we can copy the proof of 
 \cite[11.2.25]{caro-6operations} to check the other properties.
\end{proof}

\begin{lem}
\label{lemdfnqupropbis}
Let $\mathfrak{C}$ and $\mathfrak{D}$ be two  data of absolute coefficients.
Suppose 
for any object $\cW$ of $\mathrm{DVR}  (\cV)$, 
$\fB _\emptyset(\Spf (\cW) /\cW) \subset \mathfrak{D}(\Spf (\cW) /\cW)$. 
Then $
\boxtimes \left (
\fD,
S ^+ ( \mathfrak{B} _{\mathrm{wa}} , \fC) 
\right)$ is stable under weak admissible local cohomological functors.
\end{lem}

\begin{proof}
Let  $\cW$ be an object of $\mathrm{DVR}  (\cV)$, 
let $\fP$ be a formal  $\Spf (\cW)$-scheme
of formal finite type having locally finite $p$-bases  over $\Spf \cW$,  
$\E ^{(\bullet)}
\in 
\boxtimes \left (
\fD,
S ^+ ( \mathfrak{B} _{\mathrm{wa}} , \fC) 
\right) (\fP)$.
Let 
$(Y \subset \fP)$ be a weak admissible inclusion. 
We have to check that 
$\R \underline{\Gamma} ^\dag _{Y} (\E ^{(\bullet)})
\in 
\boxtimes 
\left (
\fD,
S ^+ ( \mathfrak{B} _{\mathrm{wa}} , \fC) 
\right )(\fP)$.
From \ref{lem-boxtimesDC}.\ref{lem-boxtimesDC1}, 
we have 
$\boxtimes \left (
\fD,
S ^+ ( \mathfrak{B} _{\mathrm{wa}} , \fC) 
\right)
\subset 
S ^+ ( \mathfrak{B} _{\mathrm{wa}} , \fC)$.
Hence, 
$\R \underline{\Gamma} ^\dag _{Y} (\E ^{(\bullet)})
\in 
\fC (\fP)$, and is coherent in particular.
Let $\fQ$ be a formal  $\Spf (\cW)$-scheme
of formal finite type having locally finite $p$-bases  over $\Spf \cW$,
$\cF ^{ (\bullet)}
\in 
\fD (\fQ)$.
We have to check that 
$$\R \underline{\Gamma} ^\dag _{Y} (\E ^{ (\bullet)})
\smash{\widehat{\boxtimes}}^\L
\cF ^{ (\bullet)}
\in 
S ^+ ( \mathfrak{B} _{\mathrm{wa}} , \fC) 
(\fP \times _{\scr{C} _\cW} \fQ).$$
 
Let $\fR$ be a formal  $\Spf (\cW)$-scheme
of formal finite type having locally finite $p$-bases  over $\Spf \cW$,
let $(U \subset \fP \times _{\scr{C} _\cW} \fQ \times _{\scr{C} _\cW} \fR )$ 
be a weak admissible inclusion. 
We have to prove that 
$$\R \underline{\Gamma} ^\dag _{U}
\circ  \varpi ^{!} 
\left ( 
\R \underline{\Gamma} ^\dag _{Y} (\E ^{ (\bullet)})
\smash{\widehat{\boxtimes}}^\L
\cF ^{ (\bullet)}
\right) 
\in
\mathfrak{C} (\fP \times _{\scr{C} _\cW} \fQ \times _{\scr{C} _\cW} \fR ),$$
where $\varpi \colon\fP \times _{\scr{C} _\cW} \fQ \times _{\scr{C} _\cW} \fR \to \fP \times _{\scr{C} _\cW} \fQ $ is the projection.
Recall by definition we have
$$\R \underline{\Gamma} ^\dag _{Y} (\E ^{ (\bullet)})
\smash{\widehat{\boxtimes}}^\L
\cF ^{ (\bullet)}
\riso 
\varpi _1 ^{* (\bullet) } (\R \underline{\Gamma} ^\dag _{Y} (\E ^{ (\bullet)}))
\smash{\widehat{\otimes}} ^\L _{\cO _{\fP \times _{\scr{C} _\cW} \fQ}}
\varpi _2 ^{* (\bullet) }\cF ^{ (\bullet)},$$
where 
$\varpi _1 \colon 
\fP \times _{\scr{C} _\cW} \fQ
\to \fP$
and 
$\varpi _1 \colon 
\fP \times _{\scr{C} _\cW} \fQ
\to \fP$
are the projections.
Since the functor 
$\R \underline{\Gamma} ^\dag _{U}
\circ  \varpi ^{!} $ commutes with tensor products, 
then it follows from 
Theorem \ref{2.2.18qcgen} 
that we have the isomorphism
$$\R \underline{\Gamma} ^\dag _{U}
\circ  \varpi ^{!} 
\left ( 
\R \underline{\Gamma} ^\dag _{Y} (\E ^{ (\bullet)})
\smash{\widehat{\boxtimes}}^\L
\cF ^{ (\bullet)}
\right) 
\riso
\R \underline{\Gamma} ^\dag _{U \cap (\varpi _1 \circ \varpi) ^{-1} (Y)}
\circ  \varpi ^{!} 
\left ( 
\E ^{ (\bullet)}
\smash{\widehat{\boxtimes}}^\L
\cF ^{ (\bullet)}
\right) 
.$$
Since 
$\E ^{(\bullet)}
\in 
\boxtimes \left (
\fD,
S ^+ ( \mathfrak{B} _{\mathrm{wa}} , \fC) 
\right) (\fP)$,
then
$\R \underline{\Gamma} ^\dag _{U \cap (\varpi _1 \circ \varpi) ^{-1} (Y)}
\circ  \varpi ^{!} 
\left ( 
\E ^{ (\bullet)}
\smash{\widehat{\boxtimes}}^\L
\cF ^{ (\bullet)}
\right) 
 \in
\mathfrak{C} (\fP \times _{\scr{C} _\cW} \fQ \times _{\scr{C} _\cW} \fR )$.
Hence, we are done.\end{proof}

\begin{ntn}
\label{dfnquprebis}
For any data of absolute coefficients $\fC$,
which contains $\fB _{\mathrm{wa}}$,
we set
$T _0 ( \fC) := S ^+ ( \mathfrak{B} _{\mathrm{wa}} , \fC )$.
By induction on the integer $n \geq 0$, we define 
$U _{n} ( \fC) := 
\boxtimes \left (
T _n (\fC),
T _n (\fC) \right)$
and 
$T _{n+1} ( \fC) := 
S ^+ ( \mathfrak{B} _{\mathrm{wa}} , U _{n} ( \fC) )$.
We set 
$T (\fC) := \cap _{n\geq 0}T _{n} (\fC)$.
\end{ntn}

\begin{prop}
\label{dfnqupropbis}
Let $\mathfrak{C}$ be a data of absolute coefficients which is local, satisfies $BK ^!$, is stable under 
devissage, direct summands, 
pushforwards and which
contains $\mathfrak{B} _{\mathrm{wa}} $.

\begin{enumerate}[(a)]
\item We have the inclusions 
$T _{n+1} ( \fC) 
\subset
T _{n} ( \fC) $
and the equality
$T  ( \fC) =
S ^+ ( \mathfrak{B} _{\mathrm{wa}} , T( \fC) )$.
\item 
The data of absolute coefficients $T ( \fC )$ 
contains 
$\mathfrak{B} _{\emptyset}$,
satisfies $BK ^!$,
is 
local,
is stable by devissages, direct summands, 
weak admissible local cohomological functors, 
pushforwards, extraordinary pullbacks by projections, 
base change, 
weak admissible external tensor products.

\end{enumerate}

\end{prop}

\begin{proof}
i) Following \ref{lemdfnqupropbis+mdiv}, 
$T _{0} ( \fC)$
contains 
$\fB _{\emptyset}$.
Since $T _{0} ( \fC)$ is stable under extraordinary pullbacks by projections and shift
(see \ref{S(D,C)stability+}),
then, 
from \ref{lem-boxtimesDC}.\ref{lem-boxtimesDC0}, 
$U _{0} ( \fC) $ contains 
$\fB _{\emptyset}$.
Since 
$T _{0} ( \fC) := 
S ^+ ( \mathfrak{B} _{\mathrm{wa}} , \fC )$,
then it follows from \ref{lemdfnqupropbis}
that 
$U _{0} ( \fC) $ is stable under weak admissible 
local cohomological functors.
Hence, 
$U _{0} ( \fC) $ contains 
$\fB _{\mathrm{wa}}$.
Similarly, we check by induction in $n\geq0$ that
$T _{n} ( \fC) $
and
$U _{n} ( \fC) $
are stable under weak admissible 
local cohomological functors,
contain both 
$\fB _{\emptyset}$
and 
$\fB _{wa}$.

ii) Since $T _n (\fC)$ contains $\fB _{\emptyset}$,
then 
$U _{n} ( \fC) 
\subset 
T _n (\fC)$
(use \ref{lem-boxtimesDC}.\ref{lem-boxtimesDC1}).
Hence,
$T_{n+1} ( \fC) =
S ^+ ( \mathfrak{B} _{\mathrm{wa}} , U _{n} ( \fC) )
\subset
S ^+ ( \mathfrak{B} _{\mathrm{wa}} , T _{n} ( \fC) )$.
Since, $T _{n} ( \fC)  = S ^+ ( \mathfrak{B} _{\mathrm{wa}} , U _{n-1} ( \fC) )$,
then 
$S ^+ ( \mathfrak{B} _{\mathrm{wa}} , T _{n} ( \fC) )
=T _{n} ( \fC) $
(use \ref{rem-S+wa-c)=}).
This yields 
$T _{n+1} ( \fC) 
\subset
T _{n} ( \fC) $.
This yields 
$S ^+ ( \mathfrak{B} _{\mathrm{wa}} , T( \fC) )
\subset 
S ^+ ( \mathfrak{B} _{\mathrm{wa}} , T _n ( \fC) )
= 
T _n ( \fC)$. 
Hence, 
$S ^+ ( \mathfrak{B} _{\mathrm{wa}} , T( \fC) )
\subset 
T( \fC) $. 

iii) 
Let us check now that $T (\fC)$ is stable under
weak admissible external tensor products.
Let $\cW$ be an object  of $\mathrm{DVR}  (\cV)$, 
let $\fP$ be a  formal  $\Spf (\cW)$-scheme
of formal finite type  and having locally finite $p$-bases  over $\Spf \cW$, 
let $(Y \subset \fP)$ be a weak admissible inclusion ,
let $\fQ$ be a  formal $\Spf (\cW)$-scheme of formal finite type and having locally finite $p$-bases  over $\Spf \W$, 
let $\E ^{(\bullet)}\in T (\fC)(\fP)$, 
$\FF ^{ (\bullet)}
\in
T (\fC) (\fQ)$.
We have to check
$\left ( \R \underline{\Gamma} ^\dag _{Y}  \E ^{(\bullet)} \right) 
\smash{\widehat{\boxtimes}}^\L
_{\O  _{\Spf \W }} 
\FF^{ (\bullet) } 
\in
T (\fC) (\fP\times _{\scr{C} _{\W}} \fQ)$.
Let $n \geq 0$ be an integer.
Since 
$\cE ^{ (\bullet)}
\in
T _{n+1}(\fC) (\fQ)$,
since 
$T_{n+1} ( \fC) =
S ^+ ( \mathfrak{B} _{\mathrm{wa}} , U _{n} ( \fC) )$,
then it follows from 
the remark 
\ref{rem-S+wa}.\ref{rem-S+wa-d)}, 
that 
$\R \underline{\Gamma} ^\dag _{Y}  \E ^{(\bullet)}
\in
U _{n}(\fC) (\fQ)$.
Since $\cF ^{ (\bullet)}
\in
T _{n}(\fC) (\fQ)$,
then 
$\left ( \R \underline{\Gamma} ^\dag _{Y}  \E ^{(\bullet)} \right) 
\smash{\widehat{\boxtimes}}^\L
_{\O  _{\Spf \W }} 
\FF^{ (\bullet) } 
\in
T _n (\fC) (\fP\times _{\scr{C} _{\W}} \fQ)$
Hence, we are done.

iv) By induction on $n\geq 0$, we check that 
by using \ref{waS(D,C)stability+}
and \ref{lem-boxtimesDC},
that 
$U _{n} (\fC)$
and 
$T _{n} (\fC)$
are local, satisfy $BK ^!$, are stable under 
devissage, direct summands, 
pushforwards, extraordinary pullbacks by projections, 
base change.
Then so is $T  (\fC)$.
It follows from 
\ref{rem-S+wa}.\ref{rem-S+wa-c)}
that we get 
$
T( \fC) 
\subset
S ^+ ( \mathfrak{B} _{\mathrm{wa}} , T( \fC) )$. 
\end{proof}

\subsection{Stability under duality, admissible subschemes}

\begin{ntn}
\label{ntn-DeltaYinP}
Let $\cW$ be an object of $\mathrm{DVR}  (\cV)$,   
$\fP$ be a formal  $\Spf (\cW)$-scheme of formal finite type having locally finite $p$-bases  over $\Spf \cW$.
Let $Y$ be a subscheme of $P$ such that $(Y \subset \fP)$ is weak admissible. 
We denote by $\Delta  (Y \subset \fP)$ the  stable under devissage category generated by 
the objects of the form 
$\R \underline{\Gamma} ^\dag _{Y'}  \cO _{\fP} ^{(\bullet)}$
where $Y'$ is a subscheme of $Y$.
\end{ntn}

\begin{dfn}
\label{dfn-n-admissible}
Let $\cW$ be an object of $\mathrm{DVR}  (\cV)$,   $r\geq 0$ be an integer, 
$\fP$ be a  formal  $\bbD ^r _{\cW}$-scheme of finite type  having locally finite $p$-bases  over $\Spf \cW$.
Let $Y$ be a subscheme of $P$. 
In order to get a data of coefficient almost stable under duality (see \ref{dfn-almostdual}), 
which is not a priori the case of $\fB ^+ _{wa}$, we need to introduce the notion of admissibility.

\begin{enumerate}[(a)]
\item We define by induction on $n \in \bbN$ the notion of $n$-admissibility as follows.  
We say that $(Y \subset \fP)$ is $0$-admissible if 
$(Y \subset \fP)$  is weak admissible. 
We say that $(Y \subset \fP)$ is $(n+1)$-admissible if 
$(Y \subset \fP)$ is $n$-admissible
and if 
for any projective smooth formal $\cW$-scheme $\fQ$, for any 
subscheme $U$ of $Q \times Y$, 
we have the following property : 
\begin{enumerate}[(i)]
\item there exists a projective smooth formal $\cW$-scheme
$\fQ '$,
\item there exists a subscheme 
$U'$ of $Q ' \times U$,
\item there exists an object 
$\cE ^{(\bullet)}$ 
of 
$\Delta ( U ' \subset \fQ ' \times \fQ \times \fP )$ 
\end{enumerate}
such that
$\bbD 
\left ( 
\R \underline{\Gamma} ^\dag _{U}  \cO _{ \fQ \times \fP} ^{(\bullet)}
\right) $ 
is a direct summand of 
$f ^{(\bullet)} _{+} (\cE ^{(\bullet)})$
where $f \colon \fQ ' \times \fQ \times \fP \to  \fQ \times \fP $ is the canonical projection.

\item We say that $(Y \subset \fP)$ is admissible if 
 $(Y \subset \fP)$ is $n$-admissible for any integer $n$.

\item 
We say $(Y \subset \fP)$ is ``admissible after any base change''
if for any morphism 
$\alpha \colon \cW \to \cW '$ 
of $\mathrm{DVR}  (\V)$
$(Y \times _{\bbD ^r _{\Spec l}} \bbD ^r _{\Spec l'} 
\subset 
\fP \times _{\bbD ^r _{\Spf (\cW) }} \bbD ^r _{\Spf(\cW')} )$
is admissible,
where $l$ and $l'$ are the residue fields of $\cW$ and $\cW'$.

\item Replacing ``weak admissible inclusions'' by ``admissible inclusions''
we get the notion of data of coefficients stable under admissible local cohomological functors,
under admissible tensor products and admissible duality.

\end{enumerate}

\end{dfn}

\begin{empt}
\label{empt-admi1}
Let $\cW$ be an object of $\mathrm{DVR}  (\cV)$,   
$\fP$ be a formal  $\Spf (\cW)$-scheme of formal finite type having locally finite $p$-bases  over $\Spf \cW$.
Let $Y$ be a subscheme of $P$ such that $(Y \subset \fP)$ is $n$-admissible.
Then, for any projective smooth formal $\cW$-scheme $\fQ$, for any 
subscheme $U$ of $Q \times Y$,
$(U \subset  \fQ \times \fP )$ is $n$-admissible.
\end{empt}

\begin{lemm}
\label{stab-n-admiss}
Let $\cW$ be an object of $\mathrm{DVR}  (\cV)$,   
$\fP' \hookrightarrow \fP$ be a closed immersion of   formal  
formal  $\Spf (\cW)$-schemes of formal finite type having locally finite $p$-bases  over $\Spf \cW$.
Let $Y'$ be a subscheme of $P'$. Then the following conditions are equivalent. 
\begin{enumerate}[(a)]
\item The inclusion $(Y' \subset \fP')$ is $n$-admissible. 
\item The inclusion $(Y' \subset \fP)$ is $n$-admissible.
\end{enumerate}
\end{lemm}

\begin{proof}
0) When $n =0$, the Lemma is already known (see \ref{u+closedimm-YGamma} and \ref{u!closedimm-YGamma}).
We prove the lemma by induction in $n$ as follows.

1) Suppose that $(Y' \subset \fP')$ is $n+1$-admissible.
Then by induction hypothesis, 
$(Y' \subset \fP)$ is $n$-admissible.
Let $\fQ$ be a projective smooth formal $\cW$-scheme,
$U$ be a subscheme  of $Q \times Y'$. 
By hypothesis, 
\begin{enumerate}[(i)]
\item there exists a projective smooth formal $\cW$-scheme $\fQ '$,
\item there exists a subscheme 
$U'$ of $Q ' \times U$,
\item there exists an object 
$\cE ^{\prime (\bullet)}$ 
of 
$\Delta ( U ' \subset \fQ ' \times \fQ \times \fP ')$ 
\end{enumerate}
such that
$\bbD 
\left ( 
\R \underline{\Gamma} ^\dag _{U}  \cO _{ \fQ \times \fP '} ^{(\bullet)}
\right) $ 
is a direct summand of 
$f ^{\prime (\bullet)} _{+} (\cE ^{\prime (\bullet)})$
where $f' \colon \fQ ' \times \fQ \times \fP' \to  \fQ \times \fP' $ is the canonical projection.
Let us denote by 
$f \colon \fQ ' \times \fQ \times \fP \to  \fQ \times \fP$
the canonical projection,
by 
$u\colon \fQ ' \times \fQ \times \fP ' \hookrightarrow 
\fQ ' \times \fQ \times \fP $,
$v\colon  \fQ \times \fP ' \hookrightarrow 
\fQ \times \fP $ 
the closed immersions induced by $\fP' \hookrightarrow \fP$.
Then 
$v ^{(\bullet)} _{+} \bbD 
\left ( 
\R \underline{\Gamma} ^\dag _{U}  \cO _{ \fQ \times \fP '} ^{(\bullet)}
\right)$ 
is a direct summand of 
$v ^{(\bullet)} _{+} f ^{\prime (\bullet)} _{+} (\cE ^{\prime (\bullet)})$.

It follows from \ref{rel-dual-isom} and \ref{u+closedimm-YGamma-iso1}
that we have the isomorphism
$$
v ^{(\bullet)} _{+} \bbD 
\left ( 
\R \underline{\Gamma} ^\dag _{U}  \cO _{ \fQ \times \fP'} ^{(\bullet)}
\right)
\riso 
\bbD 
v ^{(\bullet)} _{+} 
\left ( 
\R \underline{\Gamma} ^\dag _{U}  \cO _{ \fQ \times \fP'} ^{(\bullet)}
\right)
\riso
\bbD 
\R \underline{\Gamma} ^\dag _{U}  
\cO _{ \fQ \times \fP} ^{(\bullet)}
[-\delta _{\fP' /\fP}] .
$$ 
By using again \ref{u+closedimm-YGamma-iso1},
we can check 
$u ^{(\bullet)} _{+} (\cE ^{\prime (\bullet)}) 
\in 
\Delta ( U ' \subset \fQ ' \times \fQ \times \fP )$.
Since $v ^{(\bullet)} _{+} f ^{\prime (\bullet)} _{+} (\cE ^{\prime (\bullet)})
\riso 
f ^{ (\bullet)} _{+} ( u ^{(\bullet)} _{+}  (\cE ^{\prime (\bullet)}))$, 
then we conclude that 
 $(Y' \subset \fP)$ is $n+1$-admissible.

2) Conversely, suppose  $(Y' \subset \fP)$ is $n+1$-admissible.
Then by induction hypothesis, 
$(Y' \subset \fP')$ is $n$-admissible.
Let $\fQ$ be a projective smooth formal $\cW$-scheme,
$U$ be a subscheme  of $Q \times Y'$. 
By hypothesis, 
\begin{enumerate}[(i)]
\item there exists a projective smooth formal $\cW$-scheme $\fQ '$,
\item there exists a subscheme 
$U'$ of $Q ' \times U$,
\item there exists an object 
$\cE ^{(\bullet)}$ 
of 
$\Delta ( U ' \subset \fQ ' \times \fQ \times \fP )$ 
\end{enumerate}
such that
$\bbD 
\left ( 
\R \underline{\Gamma} ^\dag _{U}  \cO _{ \fQ \times \fP } ^{(\bullet)}
\right) $ 
is a direct summand of 
$f ^{ (\bullet)} _{+} (\cE ^{(\bullet)})$
where 
$f \colon \fQ ' \times \fQ \times \fP \to  \fQ \times \fP$ is the canonical projection.
Let us denote by 
$f' \colon \fQ ' \times \fQ \times \fP' \to  \fQ \times \fP'$
the canonical projection,
by 
$u\colon \fQ ' \times \fQ \times \fP ' \hookrightarrow 
\fQ ' \times \fQ \times \fP $,
$v\colon  \fQ \times \fP ' \hookrightarrow 
\fQ \times \fP $ 
the closed immersions induced by $\fP' \hookrightarrow \fP$.
Then 
$v ^{!(\bullet)}\bbD 
\left ( 
\R \underline{\Gamma} ^\dag _{U}  \cO _{ \fQ \times \fP } ^{(\bullet)}
\right)$ 
is a direct summand of 
$v ^{!(\bullet)}  f ^{ (\bullet)} _{+} (\cE ^{(\bullet)})$.
Since $\R \underline{\Gamma} ^\dag _{U}  \cO _{ \fQ \times \fP } ^{(\bullet)}$ has his support in $\fQ \times \fP'$, 
then it follows from Berthelot-Kashiwara's theorem 
\ref{u!u+=id} and the relative duality isomorphism (see \ref{rel-dual-isom}) that we have the first isomorphism: 
\begin{equation}
\notag
v ^{!(\bullet)}\bbD 
\left ( 
\R \underline{\Gamma} ^\dag _{U}  \cO _{ \fQ \times \fP } ^{(\bullet)}
\right)
\riso 
\bbD  v ^{!(\bullet)}
\left ( 
\R \underline{\Gamma} ^\dag _{U}  \cO _{ \fQ \times \fP } ^{(\bullet)}
\right)
\underset{\ref{u!closedimm-YGamma}}{\riso}
\bbD  
\left ( 
\R \underline{\Gamma} ^\dag _{U}  \cO _{ \fQ \times \fP '} ^{(\bullet)}
\right)
[\delta _{\fP' /\fP}] .
\end{equation}
By using \ref{u!closedimm-YGamma}, 
we can check that 
$u ^{!(\bullet)}  (\cE ^{ (\bullet)}) 
\in 
\Delta ( U ' \subset \fQ ' \times \fQ \times \fP ')$.
By using twice 
Berthelot-Kashiwara theorem, we get the isomorphism
$$v ^{!(\bullet)} f ^{ (\bullet)} _{+} (\cE ^{ (\bullet)})
\riso 
v ^{!(\bullet)} f ^{ (\bullet)} _{+} 
u ^{(\bullet)} _+ u ^{!(\bullet)} (\cE ^{ (\bullet)})
\riso 
v ^{!(\bullet)} 
v ^{(\bullet)} _+ 
f ^{ \prime (\bullet)} _{+} 
u ^{!(\bullet)} (\cE ^{ (\bullet)})
\riso 
f ^{ \prime (\bullet)} _{+}  u ^{!(\bullet)}  (\cE ^{(\bullet)}).$$ 
then we conclude that 
 $(Y' \subset \fP')$ is $n+1$-admissible.
\end{proof}

\begin{lemm}
\label{u!-YGamma-a}
Let $f\colon \fP '\to \fP$ 
be a quasi-projective 
(in the sense of Definition \ref{projectivefscheme}) morphism of formal $\fS$-schemes of formal finite type and having locally finite $p$-bases.
Let $Y$ be a subscheme of $P$, $Y':= f ^{-1} (Y)$.
If $(Y\subset \fP)$ 
is admissible then so is $(Y '\subset \fP')$
and we have the isomorphism of
$\smash{\underrightarrow{LD}} ^{\mathrm{b}} _{\Q,\mathrm{ovcoh}} 
(\overset{^\mathrm{l}}{} \smash{\widehat{\D}} _{\fP ' /\fS  } ^{(\bullet)} )$:
$$
\label{u!closedimm-YGamma-iso1}
\R \underline{\Gamma} ^\dag _{Y'} ( \cO _{\fP '} ^{(\bullet)})
[  \delta _{\fP ' /\fP}]
\riso 
f  ^{!(\bullet) }
\R \underline{\Gamma} ^\dag _{Y} ( \cO _{\fP} ^{(\bullet)}).$$
\end{lemm}

\begin{proof}
The fact that $(Y '\subset \fP')$ is admissible 
is a consequence of 
\ref{empt-admi1} and 
\ref{stab-n-admiss}.
The isomorphism is a consequence of \ref{u!-YGamma}.
\end{proof}

\begin{ntn}
We define the data of absolute coefficients $\fB _\mathrm{a}$ 
as follows: 
for any object $\cW$ of $\mathrm{DVR}  (\cV)$, for any integer $r\geq 0$, 
for any  formal  $\Spf(\cW)$-scheme of formal finite type $\fX$ having locally finite $p$-bases  over $\Spf \cW$,  
the category $\fB _\mathrm{a}(\X)$ is the full subcategory of 
$\smash{\underrightarrow{LD}} ^{\mathrm{b}} _{\Q,\mathrm{coh}} ( \smash{\widehat{\D}} _{\X} ^{(\bullet)})$
whose objects are of the form  
$\R \underline{\Gamma} ^\dag _{Y} \O _\X ^{(\bullet)} $,
where $Y$ is a subscheme of the special fiber of $\X$
is such that 
$(Y \subset \fX)$
is admissible after any base change.
Remark that following 
\ref{O-ovcoh} and the proof of \ref{prop=div-almostst},
we have 
$\fB _\mathrm{a} ^{(1)}
=
\fB _\mathrm{cst} ^{(1)}$.
\end{ntn}

\begin{prop}
The data of absolute coefficients $\fB _\mathrm{a} ^{+}$
satisfies $BK _+$, and is stable under admissible local cohomological functors, 
quasi-projective extraordinary pullbacks
and is almost stable under duality. 

\end{prop}

\begin{proof}
It follows from  \ref{stab-n-admiss}
(resp. \ref{u!-YGamma-a})
that $\fB _\mathrm{a} ^{+}$
satisfies $BK _+$, (resp. is stable under 
quasi-projective extraordinary pullbacks).
It follows from \ref{empt-admi1} and 
\ref{gammayY'qc} that 
 $\fB _\mathrm{a} ^{+}$ is stable under 
 admissible local cohomological functors.
Finally, by construction, 
 $\Delta (\fB _\mathrm{a})$
is almost stable under duality.
\end{proof}

Similarly to \ref{rem-S+wa}, we have the following remark.

\begin{rem}
\label{rem-S+a}
Let $\mathfrak{C}$ be a data of absolute coefficients which is stable under shifts.
It is not clear that the analogue of
\ref{S(D,C)stability}.\ref{S(D,C)stability5} is true. 
But we have the following remarks.
 
\begin{enumerate}[(a)]
\item 
\label{rem-S+a-a)}
Since this is not clear 
that $\mathfrak{B} _{\emptyset}  \subset \mathfrak{B} _{\mathrm{a}}$, then 
beware that the inclusion 
$ S ^+ (\mathfrak{B} _{\mathrm{a}}  ,\fC) \subset \fC $
might be wrong.
But we have this property ``on admissible support'' (i.e. see \ref{rem-S+a-d)}).

\item 
\label{rem-S+a-b)}
The data $S ^+ (\mathfrak{B} _{\mathrm{a}}  ,\fC) $ is stable under 
 admissible cohomological functors
and under extraordinary pullbacks by projections.

\item 
\label{rem-S+a-c)}
If $\fC$ is stable under  admissible cohomological functors
and under extraordinary pullbacks by projections,
then
we have the inclusion $\fC \subset S ^+ (\mathfrak{B} _{\mathrm{a}}  ,\fC) $.
Beware that this is not clear in this case that the inclusion is an equality.

\item 
\label{rem-S+a-d)}
For any object $\cW$ of $\mathrm{DVR}  (\cV)$, 
for any  formal  $\Spf (\cW)$-scheme of formal finite type $\fP$ having locally finite $p$-bases  over $\Spf \cW$,  
for any  admissible inclusion $(Y \subset \fP)$,
for any object $\E ^{(\bullet)}$ of $S ^+ (\mathfrak{B} _{\mathrm{a}}  ,\fC) (\fP)$,  we have
$\R \underline{\Gamma} ^\dag _{Y} \E ^{(\bullet)} \in \mathfrak{C} (\fP)$.

\item 
\label{rem-S+a-e)}
We have the equality
\begin{equation}
\label{rem-S+a-c)=}
S ^+ (\mathfrak{B} _{\mathrm{a}}  ,S ^+ (\mathfrak{B} _{\mathrm{a}}  ,\fC) ) 
=S ^+ (\mathfrak{B} _{\mathrm{a}}  ,\fC) .
\end{equation}

\end{enumerate}

\end{rem}

Similarly to \ref{lemdfnqupropbis+mdiv}, \ref{waS(D,C)stability+}
and \ref{lemdfnqupropbis}, we have the following  proposition.

\begin{prop}
\label{aS(D,C)stability+}
Let $\mathfrak{C}$ and $\mathfrak{D}$ be two  data of absolute coefficients.
\begin{enumerate}[(a)]
\item If $\mathfrak{C}$ 
contains 
$\mathfrak{B} _{\mathrm{a}} $ and is stable under shifts, then 
we have the inclusions :
$\mathfrak{B} _{\emptyset} 
\subset
S ^+ (\mathfrak{B} _{\mathrm{a}}  ,\fC)$
and 
$\mathfrak{B} _{a} 
\subset
S ^+ (\mathfrak{B} _{\mathrm{a}}  ,\fC)$.

\item If $\mathfrak{C}$ is local, satisfies $BK ^!$, is stable under 
devissage, direct summands, 
pushforwards and
contains $\mathfrak{B} _{\mathrm{a}} $,
then 
$ S ^{+}  (\mathfrak{B} _{\mathrm{a}}  ,\fC)$
is local, 
satisfies $BK ^!$, is stable under base change, 
devissage, direct summands, 
admissible local cohomological functors, 
pushforwards, extraordinary pullbacks by projections and 
contains $\mathfrak{B} _{\emptyset} $.

\item 
Suppose 
for any object $\cW$ of $\mathrm{DVR}  (\cV)$, 
$\fB _\emptyset(\Spf (\cW) /\cW) \subset \mathfrak{D}(\Spf (\cW) /\cW)$. 
Then, the data of coefficients 
$\boxtimes \left (
\fD,
S ^+ ( \mathfrak{B} _{\mathrm{a}} , \fC) 
\right)$ is stable under admissible local cohomological functors.

\end{enumerate}

\end{prop}

\begin{ntn}
\label{dfnqupreter}
For any data of absolute coefficients $\fC$,
which contains $\fB _{\emptyset}$,
we set
$T ^{a} _0 ( \fC) := S ^+ ( \mathfrak{B} _{\mathrm{a}} , \fC )$.
By induction on the integer $n \geq 0$, we define 
$U ^{a} _{n} ( \fC) := 
\boxtimes \left (
T ^{a} _n (\fC),
T ^{a} _n (\fC) \right)$
and 
$T ^{a} _{n+1} ( \fC) := 
S ^+ ( \mathfrak{B} _{\mathrm{a}} , U ^{a} _{n} ( \fC) )$.
We set 
$T ^{a} (\fC) := \cap _{n\geq 0}T ^{a} _{n} (\fC)$.
\end{ntn}

Similarly to \ref{dfnqupropbis}, we have the following proposition. 
\begin{prop}
\label{dfnqupropter}
Let $\mathfrak{C}$ be a data of absolute coefficients which is local, satisfies $BK ^!$, is stable under 
devissage, direct summands, 
pushforwards and which
contains $\mathfrak{B} _{\mathrm{a}} $.

\begin{enumerate}[(a)]
\item 
\label{dfnqupropter-a}
We have the inclusions 
$T ^{a} _{n+1} ( \fC) 
\subset
T ^{a} _{n} ( \fC) $
and the equality
$T ^{a}  ( \fC) =
S ^+ ( \mathfrak{B} _{\mathrm{a}} , T ^{a}( \fC) )$.
\item 
\label{dfnqupropter-b}
The data of absolute coefficients $T ^{a} ( \fC )$ 
satisfies $BK ^!$,
is 
local, is stable by devissages, direct summands, 
admissible local cohomological functors, 
pushforwards, extraordinary pullbacks by projections, 
base change, 
admissible external tensor products
and contains 
$\mathfrak{B} _{\emptyset}$.
\end{enumerate}
\end{prop}

\begin{theo}
\label{theo-V-6operations}
Let $\mathfrak{C}$ be a data of absolute coefficients which is local, satisfies $BK ^!$, is stable under 
devissage, direct summands, 
pushforwards and which
contains $\mathfrak{B} _{\mathrm{a}} $.
We set $V _{0} (\fC)
:=
T ^{a}  (\fC)$.
By induction on the integer $n \geq 0$, we define 
$V _{n+1} (\fC)
:= 
T ^{a}  (V _{n} (\fC ) \cap V _{n} (\fC ) ^{\vee} ) $.

\begin{enumerate}[(a)]
\item We have the inclusions 
$V _{n+1} ( \fC) 
\subset
V _{n} ( \fC) $
and the equality
$V  ( \fC) =
S ^+ ( \mathfrak{B} _{\mathrm{a}} , V( \fC) )$.
\item 
The data of absolute coefficients $V ( \fC )$ 
contains 
$\mathfrak{B} _{\emptyset}$,
satisfies $BK ^!$,
is 
local,
is stable by devissages, direct summands, 
admissible local cohomological functors, 
pushforwards, extraordinary pullbacks by projections, 
base change, 
admissible external tensor products,
admissible duality.
\end{enumerate}

\end{theo}

\begin{proof}
Let $n \geq 0$ be an integer. 
Following \ref{dfnqupropter}, we get 
$S ^+ ( \mathfrak{B} _{\mathrm{a}} , V  _n( \fC) )
=
V  _n( \fC) $.
Moreover, 
$T ^{a}  (V _{n} (\fC ) )
\subset
T ^{a} _0 (V  _n( \fC) )=
S ^+ ( \mathfrak{B} _{\mathrm{a}} , V  _n( \fC) )$.
Hence,
$T ^{a}  (V _{n} (\fC ) )
\subset
V _{n} (\fC ) $.
Hence,
we have the inclusions
$V _{n+1} ( \fC) 
=
T ^{a}  (V _{n} (\fC ) \cap V _{n} (\fC ) ^{\vee} ) 
\subset 
T ^{a}  (V _{n} (\fC ) )
\subset
V _{n} (\fC ) $.
Moreover, 
$S ^+ ( \mathfrak{B} _{\mathrm{a}} , V( \fC) )
\subset 
S ^+ ( \mathfrak{B} _{\mathrm{a}} , V _n( \fC) )
=
V _n( \fC)$. 
Hence  
$S ^+ ( \mathfrak{B} _{\mathrm{a}} , V ( \fC) )
\subset 
V ( \fC)$. 

Thanks to \ref{stab-Dvee-3prop} and \ref{dfnqupropter},
we can check by induction in $n$ that
$V _n( \fC)$  satisfies $BK ^!$,
is 
local, is stable by devissages, direct summands, 
admissible local cohomological functors, 
pushforwards, extraordinary pullbacks by projections, 
base change, 
admissible external tensor products
and contains 
$\mathfrak{B} _{\emptyset}$
 (and then 
$\mathfrak{B} _{a}$).

Let $\cW$ be an object  of $\mathrm{DVR}  (\cV)$, 
$\fP$ be a formal  $\Spf(\cW)$-scheme of formal finite type  having locally finite $p$-bases  over $\Spf \cW$, 
$(Y\subset \fP)$ be an admissible inclusion, 
and $\E ^{(\bullet)}$ be an  object  of $V ( \fC) (\X)$.
Then $\E^{(\bullet) }
\in V _{n+1} ( \fC) (\X)$.
Since 
$V _{n+1} ( \fC) (\X)
=
T ^{a}  (V _{n} (\fC ) \cap V _{n} (\fC ) ^{\vee} ) (\X)
\subset 
S ^+ (  \mathfrak{B} _{\mathrm{a}} , V _{n} (\fC ) \cap V _{n} (\fC ) ^{\vee}  )(\X)
\subset
S ^+ (  \mathfrak{B} _{\mathrm{a}} , V _{n} (\fC ) ^{\vee}  )(\X)$, 
then by using the remark \ref{rem-S+a}.\ref{rem-S+a-b)} we get
$\R \underline{\Gamma} ^\dag _{Y}  \E^{(\bullet) }
\in
V _{n} (\fC ) ^{\vee} (\X)$.
This means
$\DD  _{\X}(\R \underline{\Gamma} ^\dag _{Y}  \E^{(\bullet) }) \in V _{n} (\fC )(\X)$.
Hence, 
$\DD  _{\X}(\R \underline{\Gamma} ^\dag _{Y}  \E^{(\bullet) }) \in V  (\fC )(\X)$.
\end{proof}

\begin{ex}
Take $\fC:= \smash{\underrightarrow{LD}} ^{\mathrm{b}} _{\Q,\mathrm{coh}}$.
\end{ex}

\section{Formalism of Grothendieck six operations for arithmetic $\D$-modules over couples}

\subsection{Data of absolute coefficients over frames}

\begin{dfn}
We define the category of admissible frames over $\cV$ as follows. 
\begin{enumerate}[(a)]
\item 
An {\it admissible frame} $(Y,X,\fP,\bbD ^{r} _{\fS})$ 
over $\fS$ means that $r$ is an integer,
$\fP$ is a  
quasi-projective smooth  formal $\bbD ^r _{\fS}$-scheme, 
$X$ is a reduced closed subscheme of the special fiber $P$ of $\fP$ 
such that $(X \subset \fP)$ is admissible and $Y$ is an open subscheme of $X$. 
Let 
$(Y',X',\fP',\bbD ^{r'} _{\fS})$ 
and $(Y,X,\fP,\bbD ^{r} _{\fS})$ be two admissible frames over $\cV$. 

A morphism 
$\theta= (b,a,f,\alpha) 
\colon (Y', X', \fP',\bbD ^{r'} _{\fS})\to (Y,X,\fP,\bbD ^{r} _{\fS})$ 
of admissible frames over $\cV$ 
is the data of a  
morphism $f\colon \fP' \to \fP$ of formal $\fS$-schemes 
a morphism $a\colon X' \to X$ of 
$S$-schemes,
a morphism $b \colon Y' \to Y$ of 
schemes, 
and a morphism 
$\alpha \colon 
\bbD ^{r'} _{\fS}
\to 
\bbD ^{r} _{\fS}$ 
making commutative the following diagram
$$\xymatrix{
{Y'} 
\ar[d] ^-{b}
\ar@{^{(}->}[r] ^-{}
& 
{X'}
\ar[d] ^-{a}
\ar@{^{(}->}[r] ^-{}
& 
{\fP'} 
\ar[r] ^-{}
\ar[d] ^-{f}
& 
{\bbD ^{r'} _{\fS}} 
\ar[d] ^-{\alpha}
\\
{Y} 
\ar@{^{(}->}[r] ^-{}
& 
{X}
\ar@{^{(}->}[r] ^-{}
& 
{\fP}
\ar[r] ^-{}
& 
{\bbD ^{r} _{\fS}.} 
}$$
If there is no ambiguity with $\cV$, we simply say admissible frame or morphism of admissible frames. 

\item A morphism 
$\theta= (b,a,f,\alpha) 
\colon (Y', X', \fP',\bbD ^{r'} _{\fS})\to (Y,X,\fP,\bbD ^{r} _{\fS})$  of admissible frames over $\cV$ 
is said to be {\it complete} 
(resp. {\it strictly complete})
if $\alpha=id$ and $a$ is proper (resp. $\alpha=id$,
$f$ and $a$ are proper). 
We can also call such morphisms, 
morphisms of admissible frames over $\bbD ^{r} _{\fS}$
and write them 
$\theta= (b,a,f) 
\colon (Y', X', \fP')\to (Y,X,\fP)$.

\end{enumerate}

\end{dfn}

\begin{dfn}
\begin{enumerate}[(a)]
\item We define the category of {\it admissible couples} over $\cV$ as follow. 
A couple $(Y, X,\bbD ^{r} _{\fS})$ over $\cV$ 
is the data of a quasi-projective $\bbD ^{r} _S$-scheme $X$  (for some integer $r$) together with 
an open subscheme $Y$.

A morphism of admissible couples $u=(b,a,\alpha)\colon (Y', X',\bbD ^{r'} _\fS) \to (Y, X,\bbD ^{r} _\fS)$ over $\cV$
is the data of 
is the data of 
a morphism $a\colon X' \to X$ of 
$S$-schemes,
a morphism $b \colon Y' \to Y$ of 
schemes, 
and a morphism 
$\alpha \colon 
\bbD ^{r'} _{\fS}
\to 
\bbD ^{r} _{\fS}$ 
making commutative the following diagram
$$\xymatrix{
{Y'} 
\ar[d] ^-{b}
\ar@{^{(}->}[r] ^-{}
& 
{X'}
\ar[d] ^-{a}
\ar[r] ^-{}
& 
{\bbD ^{r'} _{\fS}} 
\ar[d] ^-{\alpha}
\\
{Y} 
\ar@{^{(}->}[r] ^-{}
& 
{X}
\ar[r] ^-{}
& 
{\bbD ^{r} _{\fS}.} 
}$$

\item A morphism of couples $u=(b,a,\alpha)\colon (Y', X',\bbD ^{r'} _\fS) \to (Y, X,\bbD ^{r} _\fS)$ over $\cV$
is said to be {\it complete} if $a$ is proper and $\alpha = id$.

\end{enumerate}
 
\end{dfn}

\begin{lem}
\label{lem-complete-frame-coup}
Let $u=(b,a,\alpha )\colon (Y', X',\bbD ^{r'} _\fS) \to (Y, X,\bbD ^{r} _\fS)$ 
be a morphism of admissible couples over $\cV$.
\begin{enumerate}[(a)]

\item There exists a morphism  of admissible frames over $\cV$ of the form
$\theta= (b,a,f,\alpha) 
\colon (Y', X', \fP',\bbD ^{r'} _{\fS})\to (Y,X,\fP,\bbD ^{r} _{\fS})$
such that $f$ is a projection morphism in the sense of \ref{dfn-CfS2}.

\item When $u$ is complete, such a morphism $\theta$ can be chosen strictly complete. 

\end{enumerate}

\end{lem}

\begin{proof}
There exist an immersion of the form 
$\iota \colon X \hookrightarrow \widehat{\bbP} ^{n} _{\bbD ^r _{\fS}}$
for some integer $n$.
This yields the commutative diagram
$$\xymatrix{
{Y'} 
\ar[d] ^-{b}
\ar@{^{(}->}[r] ^-{}
& 
{X'}
\ar[d] ^-{a}
\ar[r] ^-{u'}
& 
{\widehat{\bbP} ^{n} _{\bbD ^{r'} _{\fS}}} 
\ar[r] ^-{}
\ar@{}[rd] ^-{}|\square
\ar[d] ^-{\varpi}
& 
{\bbD ^{r'} _{\fS}} 
\ar[d] ^-{\alpha}
\\
{Y} 
\ar@{^{(}->}[r] ^-{}
& 
{X}
\ar@{^{(}->}[r] ^-{\iota}
& 
{\widehat{\bbP} ^{n} _{\bbD ^r _{\fS}}}
\ar[r] ^-{}
& 
{\bbD ^{r} _{\fS},} 
}$$
where $u'$ is the morphism making the diagram commutative.
There exists an immersion of the form 
$\iota ' \colon X '\hookrightarrow \widehat{\bbP} ^{n'} _{\bbD ^{r'} _{\fS}}$
for some integer $n'$.
We get the morphism
$\iota''=( \iota' , u') \colon X \to 
\widehat{\bbP} ^{n'} _{\bbD ^{r'} _{\fS}}
\times _{\bbD ^{r'} _{\fS}}
\widehat{\bbP} ^{n} _{\bbD ^{r'} _{\fS}}$.
Since $\iota'$ is an immersion, then so is $\iota''$.
Let 
$\varpi _1  
\colon 
\widehat{\bbP} ^{n'} _{\bbD ^{r'} _{\fS}}
\times _{\bbD ^{r'} _{\fS}}
\widehat{\bbP} ^{n} _{\bbD ^{r'} _{\fS}}
\to 
\widehat{\bbP} ^{n} _{\bbD ^{r'} _{\fS}}$
be the canonical projection.
We get the morphism 
$g:= \varpi \circ \varpi _1
\colon \widehat{\bbP} ^{n'} _{\bbD ^{r'} _{\fS}}
\times _{\bbD ^{r'} _{\fS}}
\widehat{\bbP} ^{n} _{\bbD ^{r'} _{\fS}}
\to 
\widehat{\bbP} ^{n} _{\bbD ^{r} _{\fS}}$.
Let 
$\fP$ be an open subscheme of  
$\widehat{\bbP} ^{n} _{\bbD ^r _{\fS}}$
containing $X$ and such that the factorization
$v \colon X \hookrightarrow \fP$ of $\iota $ is a closed immersion.

a)  We can choose an open formal subscheme
$\fP'$ of  $g ^{-1} (\fP)$ such that 
the factorization
$v '\colon X ' \hookrightarrow \fP'$ of $\iota''$ is a closed immersion.
This yields the morphism
$\theta= (b,a,f,\alpha) 
\colon (Y', X', \fP',\bbD ^{r'} _{\fS})\to (Y,X,\fP,\bbD ^{r} _{\fS})$,
where $f$ is the morphism induced by $g$.

b) Suppose now $a$ is proper and $\alpha = id$.
We get $\varpi = id$ and then $g=\varpi _1$ is proper. 
Set  $\fP ' : = g ^{-1} (\fP)$.
Since $f \colon \fP ' \to \fP$ is proper,
since the morphism $v \circ a \colon X ' \hookrightarrow \fP'$ is proper, 
then the immersion $X ' \hookrightarrow \fP'$ (induced by $\iota ''$) is proper, i.e. 
is a closed immersion.
Hence, we get the strictly complete morphism $\theta= (b,a,f,id) 
\colon 
(Y', X', \fP',\bbD ^{r'} _{\fS})\to (Y,X,\fP,\bbD ^{r} _{\fS})$. 
\end{proof}

\begin{dfn}
\label{dfn-framW[[t]]}
\begin{enumerate}[(a)]
\item We define the category of ``frames over $\cV[[t]]$'' 
whose objects are 
the frames of the form
 $(Y,X,\fP,\bbD ^{1} _{\fS})$  
 and whose morphisms are morphisms of 
 frames of the form
 $\theta= (b,a,f,\mathrm{id}) 
\colon (Y', X', \fP',\bbD ^{1} _{\fS})\to (Y,X,\fP,\bbD ^{1} _{\fS})$.
Since the morphism on $\bbD ^{1} _{\fS}$ are the identity, 
we denote a frame over $\cV[[t]]$ simply by 
 $(Y,X,\fP)$ and morphisms of frames over $\cV [[t]]$ are denoted by
 $\theta= (b,a,f) 
\colon (Y', X', \fP')\to (Y,X,\fP)$.

\item  We define the category of ``couples over $\cV[[t]]$'' 
whose objects are 
the couples of the form
 $(Y,X,\bbD ^{1} _{\fS})$  
 and whose morphisms are morphisms of 
 couples of the form
 $\theta= (b,a,\mathrm{id}) 
\colon (Y', X', \bbD ^{1} _{\fS})\to (Y,X,\bbD ^{1} _{\fS})$.
Since the morphism on $\bbD ^{1} _{\fS}$ are the identity, 
we denote a frame over $\cV[[t]]$ simply by 
 $(Y,X)$ and morphisms of couples over $\cV [[t]]$ are denoted by
 $\theta= (b,a) 
\colon (Y', X')\to (Y,X)$.

\end{enumerate}

\end{dfn}

\begin{ntn}
\label{ntn-6operations}
Let $\mathfrak{C} $ be a  data of absolute coefficients over $\cV$.

\begin{enumerate}[(a)]
\item Let  $(Y,X,\fP,\bbD ^{r} _{\fS})$ be  an admissible frame over $\cV$. 
We denote by $\mathfrak{C} (Y,\fP,\bbD ^{r} _{\fS}/\cV)$
the full subcategory of 
$\mathfrak{C}  (\fP)$ 
of objects $\E$ such that there exists an isomorphism of the form 
$\E \riso \R \underline{\Gamma} ^\dag _{Y} (\E)$.
We remark that $\mathfrak{C} (Y,\fP,\bbD ^{r} _{\fS}/\cV)$ only depend on the immersion and the structural map
$Y \hookrightarrow \fP\to \bbD ^{r} _{\fS}$ which explains the notation.
We might choose $X$ equal to the closure of $Y$ in $P$.

\item Let  $(Y,X,\fP)$ be  a frame over $\cV[[t]]$. 
Similarly, we denote by $\mathfrak{C} (Y,\fP/\cV[[t]])$
the full subcategory of 
$\mathfrak{C}  (\fP)$ 
consisting of objects $\E$ such that there exists an isomorphism of the form 
$\E \riso \R \underline{\Gamma} ^\dag _{Y} (\E)$.
\end{enumerate}

\end{ntn}

\begin{ntn}
\label{ntn-t-structureovcoh}
Let $\mathfrak{C}$ be a  data of absolute coefficients stable under devissages  and cohomology. 
Let $(Y,X,\fP,\bbD ^{r} _{\fS})$ be an admissible frame over $\cV$
(resp. let $(Y,X,\fP)$ be a frame over $\bbD ^{1} _{\fS}$).
Choose $\U$ an open set of $\fP$ such that 
$Y$ is closed in $\U$.
We introduce the following notation (in the respective case, we remove the indication 
$\bbD ^{1} _{\fS}$).
\begin{enumerate}[(a)]
\item Similarly to \cite[1.2.1-5]{Abe-Caro-weights} , 
we define a canonical t-structure on 
$\mathfrak{C} (Y , \fP ,\bbD ^{r} _{\fS}/\cV)$
as follows.
We denote by 
$\fC  ^{\leq n}    (Y,\fP,\bbD ^{r} _{\fS}/\cV)$
(resp. 
$\fC   ^{\geq n}  (Y,\fP,\bbD ^{r} _{\fS}/\cV)$)
the full subcategory of 
$\mathfrak{C} (Y , \fP ,\bbD ^{r} _{\fS}/\cV)$
of complexes 
$\E$ such that
$\E |\U 
\in 
\fC  ^{\leq n}  (Y,\U,\bbD ^{r} _{\fS}/\cV):=
\fC   (Y,\U,\bbD ^{r} _{\fS}/\cV)
\cap 
\fC  ^{\leq n}  (\U,\bbD ^{r} _{\fS}/\cV)$
(resp. 
$\E |\U 
\in 
\fC  ^{\geq n}  (Y,\U,\bbD ^{r} _{\fS}/\cV)
:=
\fC   (Y,\U,\bbD ^{r} _{\fS}/\cV)
\cap 
\fC  ^{\geq n}  (\U,\bbD ^{r} _{\fS}/\cV)$),
where the t-structure on 
$\fC  (\U,\bbD ^{r} _{\fS}/\cV)$
is the canonical one (see \ref{t-structure-ovcoh-oc-hol-h}).
The heart of this t-structure
will be denoted by
$\mathfrak{C} ^0   (Y,\fP,\bbD ^{r} _{\fS}/\cV)$.
Finally, we denote by 
$\mathcal{H} ^i _{\mathrm{t}}$
the $i$th space of cohomology with respect to this canonical t-structure. 

\item Suppose $Y/S$ has locally finite $p$-bases. 
Then, we denote by
$\fC _{\mathrm{isoc}}   (Y,\fP,\bbD ^{r} _{\fS}/\cV)$
(resp. $\fC ^{\geq n}  _{\mathrm{isoc}}   (Y,\fP,\bbD ^{r} _{\fS}/\cV)$,
resp. $\fC ^{\leq n}  _{\mathrm{isoc}}   (Y,\fP,\bbD ^{r} _{\fS}/\cV)$,
resp. $\fC ^0 _{\mathrm{isoc}}   (Y,\fP,\bbD ^{r} _{\fS}/\cV)$)
the full subcategory of 
(resp. $\fC ^{\geq n}  (Y,\fP,\bbD ^{r} _{\fS}/\cV)$,
resp. $\fC ^{\leq n}   (Y,\fP,\bbD ^{r} _{\fS}/\cV)$,
resp. $\fC ^0  (Y,\fP,\bbD ^{r} _{\fS}/\cV)$)
consisting of complexes
$\E ^{(\bullet)}$ such that 
$\mathcal{H} ^i  (\E ^{(\bullet)} |\U ) 
\in 
\mathrm{MIC} ^{(\bullet)} (Y, \U/K)$.
We refer ``$\mathrm{isoc}$'' as isocrystals. 
The reason is 
the equivalence of categories
of \cite[5.4.6.1]{caro-pleine-fidelite} in the context of smooth formal schemes.
In this paper, we avoid trying to check such equivalence of categories
(other than the easier case where the partial compactification is smooth).

\end{enumerate}

\end{ntn}

\begin{rem}
\label{rem-tstructure-exact}
Let $\mathfrak{C}$ be a  data of absolute coefficients stable under devissages  and cohomology. 
Let $\fP$ be a quasi-projective smooth formal $\bbD ^r _\fS$-scheme,
$Y$ be a subscheme of $P$,
$Z$ be a closed subscheme of $Y$,
and 
$Y':= Y \setminus Z$.

\begin{enumerate}[(a)]
\item We get the t-exact functor 
$(\hdag Z) 
\colon
\mathfrak{C} (Y,\fP,\bbD ^{r} _{\fS}/\cV)
\to 
\mathfrak{C} (Y', \fP,\bbD ^{r} _{\fS}/\cV)$.
Beware the functor 
$(\hdag Z) 
\colon
\mathfrak{C} (Y,\fP,\bbD ^{r} _{\fS}/\cV)
\to 
\mathfrak{C} (Y,\fP,\bbD ^{r} _{\fS}/\cV)$
is not always t-exact.

\item 
\label{rem-tstructure-exact2}
We say that $Z$ locally comes from a divisor of $P$ if 
locally in $P$, there exists a divisor $T$ of $P$ such that 
$Z = Y \cap T$ (this is equivalent to saying that locally in $P$, 
the ideal defining $Z \hookrightarrow Y$ is generated by one element).
In that case, 
we get the t-exact functor 
$(\hdag Z) 
\colon
\mathfrak{C} (Y,\fP,\bbD ^{r} _{\fS}/\cV)
\to 
\mathfrak{C} (Y,\fP,\bbD ^{r} _{\fS}/\cV)$.
Indeed, 
by construction of our t-structures, 
we can suppose $Y$ is closed in $\fP$ (and then we reduce to the case where the t-structure on 
$\mathfrak{C} (Y,\fP,\bbD ^{r} _{\fS}/\cV)$ is induced by 
the standard t-structure
of $\smash{\underrightarrow{LD}} ^{\mathrm{b}} _{\Q,\mathrm{coh}} ( \smash{\widehat{\D}} _{\fP} ^{(\bullet)})$).
Since the property is local, 
we can suppose there exists a divisor
 $T$ such that $Z = T \cap Y$. 
Then both functors 
$(\hdag Z) $ and $(\hdag T)$ 
of 
$\mathfrak{C} (Y,\fP,\bbD ^{r} _{\fS}/\cV)
\to 
\mathfrak{C} (Y,\fP,\bbD ^{r} _{\fS}/\cV)$
are isomorphic. Since 
$(\hdag T)$ is exact,
we are done.
\end{enumerate}

\end{rem}

\begin{rem}
\label{rem-sp-desc-123}
Let $\Lambda ^{\flat}:= k  (( t ^{p ^{-\infty}})) $ be a perfect closure of $\Lambda := k ((t))$.
Let $Y$ be a reduced 
$\Spec \Lambda $-scheme of finite type. 
Let $Y^{\flat}:= Y \times _{\Spec \Lambda} \Spec \Lambda ^{\flat}$,
and 
$\widetilde{Y}  ^{\flat}:= 
Y^{\flat}_{\mathrm{red}}:= (Y \times _{\Spec \Lambda} \Spec \Lambda ^{\flat}) _{\mathrm{red}}$
be the corresponding reduced scheme.
Let $ \Lambda'$  be a finite radicial extension of $ \Lambda$ included in $ \Lambda ^\flat$
(i.e.  $ \Lambda' =k  (( t ^{p ^{-n}}))$ for some integer $n$).
We put $Y ': =Y \times _{\Spec \Lambda} \Spec ( \Lambda')$.

\begin{enumerate}[(a)]

\item \label{rem-sp-desc-123-a}
 By using
 \cite[8.7.2]{EGAIV3}, 
 \cite[8.8.2.(ii)]{EGAIV3} and \cite[8.10.5.(v)]{EGAIV3}, 
for $ \Lambda'$ large enough,
there exist a reduced $ \Lambda'$-scheme $\widetilde{Y} '$ of finite type
satisfying 
$\widetilde{Y} ^\flat \riso \widetilde{Y} ' \times _{\Spec ( \Lambda')} \Spec ( \Lambda ^\flat)$.
For $ \Lambda'$ large enough, 
it follows from \cite[8.8.2.(i)]{EGAIV3} 
that 
there exists a morphism 
$\widetilde{Y} ' \to Y '$
inducing the closed immersion 
$\widetilde{Y} ^\flat \hookrightarrow Y ^\flat$.
By using \cite[8.10.5]{EGAIV3}, 
for $ \Lambda'$ large enough, 
we can suppose that $\widetilde{Y} ' \to Y '$ is a surjective closed immersion.
Since $\widetilde{Y} '$ is reduced, this yields
$\widetilde{Y} ' = Y' _{\mathrm{red}}$, for $ \Lambda'$ large enough.

\item 
\label{rem-sp-desc-123-3b}
Hence, by using \cite[17.7.8]{EGAIV4}, 
we check that if 
$Y^{\flat}_{\mathrm{red}}$ is smooth 
(resp. étale) over $ \Lambda ^\flat$, then so is 
$Y' _{\mathrm{red}}$ over $ \Lambda'$
for $ \Lambda'$ large enough.

\item 
\label{rem-sp-desc-123-3c}
Hence, if $Y$ is of dimension $0$, then
$Y' _{\mathrm{red}}$ is a finite and étale $ \Lambda'$-scheme
for $ \Lambda'$ large enough.
Indeed, 
since $\Lambda ^{\flat}$ is perfect and 
since
$Y^{\flat}_{\mathrm{red}}$
is a reduced $\Lambda ^{\flat}$-scheme of finite type of dimension $0$, then 
$Y^{\flat}_{\mathrm{red}}$
is a finite and étale $\Lambda ^{\flat}$-scheme.
We conclude using the previous remark.

\end{enumerate}

\end{rem}

\begin{dfn}
 \label{recallofstrat}
Let $(Y,X,\fP,\bbD ^{r} _{\fS})$ be an admissible frame over $\cV$.
 An ordered set of
 subschemes $\{Y_i\}_{i=1,\dots,r}$ of $Y$ is said to be a {\em stratification having locally finite $p$-bases} 
 if the following
 holds: 1.\ $\{Y_i\}$ is a stratification, namely putting
 $Y_0:=\emptyset$, $Y_k$ is an open subscheme of
 $Y\setminus\bigcup_{i<k}Y_i$ and $Y=\bigcup_{1\leq i\leq r} Y_i$. 2.\
 $Y_i/S$ has locally finite $p$-bases. 
 
\end{dfn}

\begin{lem}
\label{rem-sp-desc-123-3d}
Let $Y$ be a reduced 
$\Spec k[[t]]$-scheme of finite type. 
Then for $n$ large enough, 
there exists an open dense subscheme $U'$ of  
$( Y \times _{\Spec k[[t]] } \Spec (k[[t ^{p ^{-n}}]]) ) _{\mathrm{red}}$ such that $U' /S$ has locally finite $p$-bases.
Moreover, we can choose $U'$ so that
the irreducible components of $U' $ are either smooth over $S$ or smooth over 
$\Spec (k((t ^{p ^{-n}}))) $.

\end{lem}

\begin{proof}
Replacing $Y$ by an open dense subset, we can suppose
$Y$ is a direct sum of its irreducible components. Hence, 
we can suppose $Y$ irreducible. 
When $Y _\eta: = Y \times _{\Spec k[[t]] } \Spec k ((t))  $ is empty, then $Y$ is a reduced scheme of finite type over $S$.
Since $k$ is perfect, there exists an open dense subscheme $U$ of $Y$ such that $U/S$ is smooth and in particular has locally finite $p$-bases. 
Suppose now $Y _\eta$ is not empty, i.e. is dense in $Y$.
We can suppose $Y _\eta = Y$, i.e. $Y$ is in fact a scheme of finite type over $\Spec k (( t))$.
With notation \ref{rem-sp-desc-123},
since $\widetilde{Y}  ^{\flat}$ is a reduced scheme of finite type over the perfect field $\Lambda ^{\flat}$, 
there exists an open dense subscheme $ U ^{\flat}$ of $\widetilde{Y}  ^{\flat}$ 
such that $ U ^{\flat}/\Spec \Lambda ^\flat$ is smooth.
 By using
 \cite[8.7.2]{EGAIV3}, 
 \cite[8.8.2.(ii)]{EGAIV3} and \cite[8.10.5.(v)]{EGAIV3}, 
for $n$ large enough,
there exist a reduced $\Spec (k((t ^{p ^{-n}}))$-scheme $U'$ of finite type
satisfying 
$U ^\flat \riso U ' \times _{\Spec (k((t))} 
\Spec (k((t ^{p ^{-n}}))$.
By using \cite[17.7.8]{EGAIV4}, 
we can suppose 
$U'/\Spec (k((t ^{p ^{-n}}))$ is smooth.
In particular, 
$U' /S$ has locally finite $p$-bases.
\end{proof}

\begin{prop}
[Devissage in isocrystals]
\label{ovcoharedev}
Let $\mathfrak{C}$ be a restricted data of absolute coefficients over $ \cV$ 
stable under devissages, cohomology, local cohomological functors,
quasi-projective extraordinary pullbacks.
Let $(Y,X,\fP)$ be a frame over $\cV[[t]]$.
Let $\E ^{(\bullet)} \in \fC  (Y,\fP/\cV)$.
For any integer $n $, 
let $\fP _{(n)} := \fP \times _{\Spf \cV [[t ]]} \Spf \cV [[t ^{p ^{-n}}]]$,
and
$f _{(n)}
\colon 
\fP _{(n)}
\to \fP $ be the projection.
For $n $ large enough
there exists a stratification having locally finite $p$-bases
$(U ' _{i}) _{i=1, \dots , r}$ of $( Y \times _{\Spec k[[t]] } \Spec (k[[t ^{p ^{-n}}]]) ) _{\mathrm{red}}$ 
(see Definition \ref{recallofstrat})
such that we have
$ \R \underline{\Gamma} ^\dag _{U ' _i}  (f _{(n)} ^! (\E ^{(\bullet)}) )
\in 
\mathfrak{C} _{\mathrm{isoc}} (U ' _i, \fP _{(n)})$
 for any $i =1,\dots, r$.
Moreover, we can choose such a stratification $(U ' _{i}) _{i=1, \dots , r}$ so that
the irreducible components of $U' _i $ are either smooth over $S$ or smooth over 
$\Spec (k((t ^{p ^{-n}}))) $,
for any $i =1,\dots, r$.

\end{prop}

\begin{proof}
Since $\fP _{(n)}$ is smooth over $\Spf \cV [[t ^{p ^{-n}}]]$ 
and $\Spf \cV [[t ^{p ^{-n}}]]$ has a finite $p$-basis over 
$\Spf \cV$, then $\fP _{(n)}$ has locally finite $p$-basis and 
then the theorem is well defined. It remains to prove it. 
Following \ref{rem-sp-desc-123-3d}, 
for $n$ large enough, 
there exists an open dense subscheme $U' _1$ of  
$( Y \times _{\Spec k[[t]] } \Spec ( k [[t ^{p ^{-n}}]]) ) _{\mathrm{red}}$ such that 
the irreducible components of $U' _1 /S$ are either smooth over $S$ or smooth over 
$\Spec (k[[t ^{p ^{-n}}]]) $.
We have
$ \R \underline{\Gamma} ^\dag _{U ' _1}  f _{(n)} ^! (\E ^{(\bullet)}) 
\in 
\fC  (U ' _1,\fP _{(n)}/\cV)$.
Hence, shrinking $U ' _1$ if necessary, 
using Theorem \cite[3.4]{caro-holo-sansFrob} (which is still valid in our context), 
we get 
$ \R \underline{\Gamma} ^\dag _{U ' _1}  f _{(n)} ^! (\E ^{(\bullet)}) 
\in  \mathfrak{C}  _{\mathrm{isoc}}  (U ' _1,\fP _{(n)}/\cV)$.
By iterating this processus, we can conclude. 
\end{proof}

\begin{ex}
Following \ref{ovcoh-invim} and \ref{coro-ovcoh-oc-tstr},
we can apply the Proposition \ref{ovcoharedev} in the case where
the restricted data of absolute coefficients $\fC$ is either
$\smash{\underrightarrow{LD}} ^{\mathrm{b}(1)} _{\Q,\mathrm{ovcoh}}$, or
$\smash{\underrightarrow{LD}} ^{\mathrm{b}(1)} _{\Q,\mathrm{oc}}$. 

\end{ex}

\subsection{Formalism of Grothendieck six operations over couples over $\cV[[t]]$}

\begin{lem}
\label{ind-CYWpre}
Let $\mathfrak{C}$ be a  data of absolute coefficients over $\fS$
which contains $\fB _\mathrm{\emptyset}$, 
which is 
stable by devissages, 
pushforwards, extraordinary pullbacks by projections,
and under admissible local cohomological functors.

Let
$\theta= (b,a,f) 
\colon (Y', X', \fP')\to (Y,X,\fP)$ 
be  a morphism of admissible frames over $\bbD ^{r} _{\fS}$
such that the morphisms $a$ and $b$ are proper.
For any 
$\E ^{(\bullet)} 
\in 
\fC (Y , \fP /\cV)$, 
for any $\E ^{\prime (\bullet)} \in \fC (Y ', \fP '/\cV)$  (recall notation \ref{ntn-6operations}),
we have 
$$\mathrm{Hom}  _{\fC (Y , \fP /\cV)}
( f ^{(\bullet) } _{+} ( \E ^{\prime (\bullet) }) , \E ^{(\bullet) }) 
\riso 
 \mathrm{Hom}  _{\fC (Y ', \fP '/\cV)}
( \E ^{\prime (\bullet) } ,
\R \underline{\Gamma} ^{\dag} _{Y'}
f ^{!(\bullet )}  ( \E ^{ (\bullet) })).$$
\end{lem}

\begin{proof}
Let us check the first statement. 
Replacing $X$ and $X'$ by the closure of $Y$ in $P$ and $Y'$ in $P'$ if necessary,
we can suppose $Y$ is dense in $X$ and $Y'$ is dense in $X'$.
Let 
$\E ^{(\bullet)} 
\in 
\fC (Y , \fP /\cV)$, 
and $\E ^{\prime (\bullet)} \in \fC (Y ', \fP '/\cV)$.
Since $a$ is proper,
using \ref{cor-adj-formulbis-bij2}, 
the stability of $\fC$  under extraordinary pullbacks,
and the equivalence of categories
\ref{eqcat-limcoh}, we get the bijection
$$\mathrm{Hom}  _{\smash{\underrightarrow{LD}} ^{\mathrm{b}} _{\Q,\mathrm{coh}}
(\smash{\widehat{\D}} _{\fP} ^{(\bullet)})}
( f ^{(\bullet) } _{+} ( \E ^{\prime (\bullet) }) , \E ^{(\bullet) }) 
\riso 
 \mathrm{Hom}  _{\smash{\underrightarrow{LD}} ^{\mathrm{b}} _{\Q,\mathrm{coh}}
(\smash{\widehat{\D}} _{\fP'} ^{(\bullet)})}
( \E ^{\prime (\bullet) } ,
f ^{!(\bullet )}  ( \E ^{ (\bullet) })).$$
Since $\E ^{\prime (\bullet) } \in \fC (X ', \fP '/\cV)$,
then the functor $\R \underline{\Gamma} ^{\dag} _{X'}$ induces the bijection
$$ \mathrm{Hom}  _{\smash{\underrightarrow{LD}} ^{\mathrm{b}} _{\Q,\mathrm{coh}}
(\smash{\widehat{\D}} _{\fP'} ^{(\bullet)})}
( \E ^{\prime (\bullet) } ,
f ^{!(\bullet )}  ( \E ^{ (\bullet) }))
\riso
 \mathrm{Hom}  _{\fC (X ', \fP '/\cV)}
( \E ^{\prime (\bullet) } ,
\R \underline{\Gamma} ^{\dag} _{X'}
f ^{!(\bullet )}  ( \E ^{ (\bullet) })).$$
Since $a$ and $b$ are proper, then 
the open immersion 
$Y' \subset a ^{-1} (Y)$ is proper. 
Since $Y'$ is dense in $X'$, then 
$Y'=a ^{-1} (Y)$.
Hence,
the functors 
$\R \underline{\Gamma} ^{\dag} _{X'}
f ^{!(\bullet )}$
and
$\R \underline{\Gamma} ^{\dag} _{Y'}
f ^{!(\bullet )}$
(resp. 
$f ^{(\bullet) } _{+} $ 
and 
$\R \underline{\Gamma} ^{\dag} _{Y} f ^{(\bullet) } _{+} $)
are isomorphic over 
$\fC(Y , \fP /\cV)$
(resp. $\fC(Y ', \fP '/\cV)$).
This yields,
the functor 
$\R \underline{\Gamma} ^{\dag} _{X'}
f ^{!(\bullet )}$
(resp. $f ^{(\bullet) } _{+} $)
induces
$\R \underline{\Gamma} ^{\dag} _{X'}
f ^{!(\bullet )}
\colon
\fC (Y , \fP /\cV)
\to 
\fC (Y ', \fP '/\cV)$
(resp. 
$f ^{(\bullet) } _{+} 
\colon 
\fC (Y ', \fP '/\cV)
\to 
\fC (Y , \fP /\cV)$).
Hence, we are done.
\end{proof}

\begin{thm}
[Independence]
\label{ind-CYW}
Let $\mathfrak{C}$ be a restricted data of absolute coefficients over $\fS$
which contains $\fB _\mathrm{\emptyset}$, 
which is 
stable  under 
devissages, 
pushforwards, quasi-projective extraordinary pullbacks,
and under local cohomological functors.

Let
$\theta= (id,a,f) 
\colon (Y, X', \fP')\to (Y,X,\fP)$ 
be  a morphism of frames over $\bbD ^{1} _{\fS}$
such that $a$ is proper.

\begin{enumerate}[(a)]
\item 
Suppose moreover that 
$\mathfrak{C}$ is stable under cohomology.
 Then, 
for any 
$\E ^{(\bullet)} \in \fC ^0 (Y , \fP /\cV)$,
for any 
$\E ^{\prime(\bullet)} \in 
\fC ^0 (Y , \fP '/\cV)$, 
for any $n \in \Z \setminus \{ 0\}$, we have
$$\mathcal{H} _t ^n \R \underline{\Gamma} ^{\dag} _{Y} f ^{!(\bullet)} ( \E ^{(\bullet)}) =0,
\hspace{1 cm}
\mathcal{H} ^n _t f _{+} ^{(\bullet)} ( \E ^{\prime (\bullet)}) =0.$$

\item 
For any 
$\E ^{(\bullet)} 
\in 
\mathfrak{C} (Y , \fP /\cV)$, 
for any $\E ^{\prime (\bullet)} \in \mathfrak{C} (Y , \fP '/\cV)$, 
the adjunction morphisms 
$\R \underline{\Gamma} ^{\dag} _{Y} f ^{!(\bullet)}  f _{+} ^{(\bullet)} (\E ^{\prime (\bullet)})
\to 
\E ^{\prime (\bullet)}$ 
and 
$f _{+} ^{(\bullet)}  \R \underline{\Gamma} ^{\dag} _{Y} f ^{!(\bullet)}  
(\E ^{(\bullet)} )
\to 
\E ^{(\bullet)} $
are isomorphisms.
In particular, the functors
$\R \underline{\Gamma} ^{\dag} _{Y} f ^{!(\bullet)}$
and 
$f _{+} ^{(\bullet)}$ induce quasi-inverse equivalences of categories between 
$\mathfrak{C} (Y , \fP /\cV)$
and 
$\mathfrak{C} (Y , \fP '/\cV)$.

\end{enumerate}

\end{thm}

\begin{proof}
Using the stability properties that $\mathfrak{C}$ satisfies, 
we check that the functors 
$f _+ ^{(\bullet )}\colon 
\mathfrak{C} (Y , \fP' /\cV)
\to 
\mathfrak{C} (Y , \fP /\cV)$
and
$\R \underline{\Gamma} ^{\dag} _{Y} f ^{!(\bullet )}
\colon 
\mathfrak{C} (Y , \fP /\cV)
\to 
\mathfrak{C} (Y , \fP '/\cV)$
are well defined. 
Since $\mathfrak{C}$ is included in $\smash{\underrightarrow{LD}} ^{\mathrm{b}(1)} _{\Q,\mathrm{ovcoh}} $, 
we reduce to check the case where 
$\mathfrak{C}=\smash{\underrightarrow{LD}} ^{\mathrm{b}} _{\Q,\mathrm{ovcoh}} $.
We proceed similarly to \cite[3.2.6]{caro_surcoherent}:
Choose 
$\U$ (resp. $\U'$) an open set of $\fP$ (resp. $\fP '$) such that 
$Y$ is closed in $\U$ (resp. $Y$ is closed in $\U'$),
and such that $f (\U') \subset \U$.
The functor 
$| \U \colon 
\smash{\underrightarrow{LD}} ^{\mathrm{b}} _{\Q,\mathrm{ovcoh}} (Y , \fP /\cV)
\to 
\smash{\underrightarrow{LD}} ^{\mathrm{b}} _{\Q,\mathrm{ovcoh}} (Y , \U /\cV)$
is t-exact, and the same with some primes.
Moreover, 
for any $\E ^{(\bullet)}
\in 
\smash{\underrightarrow{LM}}  _{\Q,\mathrm{ovcoh}} (Y , \fP /\cV)$
(or 
$\E ^{(\bullet)}
\in \smash{\underrightarrow{LD}} ^{\mathrm{b}} _{\Q,\mathrm{ovcoh}} (Y , \fP /\cV)$), 
the property 
$\E ^{(\bullet)}= 0$ is equivalent to 
$\E ^{(\bullet)}|\U= 0$. Hence, we can suppose 
$\U= \fP$ and $\U'= \fP'$, i.e.
$Y \hookrightarrow P$
and 
$Y \hookrightarrow P'$ are closed immersions.
For any integer $n $, 
let $\fP _{(n)} := \fP \times _{\Spf \cV [[t ]]} \Spf \cV [[t ^{p ^{-n}}]]$,
$\fP ' _{(n)} := \fP ' \times _{\Spf \cV [[t ]]} \Spf \cV [[t ^{p ^{-n}}]]$,
$f _{(n)}
\colon 
\fP _{(n)}
\to \fP $ 
and
$f '_{(n)}
\colon 
\fP '_{(n)}
\to \fP '$
be the canonical projections, and
$Y _{(n)}:= ( Y \times _{\Spec k[[t]] } \Spec (k[[t ^{p ^{-n}}]]) ) _{\mathrm{red}}$.
As a topological space,
$Y _{(n)} = f _{(n)} ^{-1} (Y) = f _{(n)} ^{\prime -1} (Y)$.
Since $f _{(n)}$ (resp. $f ' _{(n)}$) 
is surjective, finite and radicial
then it is a universal homeomorphism (see \ref{univ-homeo}.\ref{univ-homeo-3}
and 
then the functors $f _{(n)+}$ and $f _{(n)} ^!$ induce exact quasi-inverse equivalence of categories (\ref{256Be2})
between 
$\mathfrak{C} (Y , \fP /\cV)$
and
$\mathfrak{C} (Y _{(n)}  , \fP _{(n)}/\cV)$
(resp. $\mathfrak{C} (Y , \fP ' /\cV)$
and
$\mathfrak{C} (Y _{(n)}  , \fP '_{(n)}/\cV)$).

1) In this step we make the following assumption : 
for $n$ large enough, we suppose that 
$(Y _{(n)} )  _{\mathrm{red}}/S$ has locally finite $p$-bases and moreover that 
the irreducible components of $(Y _{(n)} )  _{\mathrm{red}}$ are either  smooth over $\Spec (k((t ^{p ^{-n}}))$
or smooth over $S$. 

Since $f _{(n)+}$ and $f _{(n)} ^!$ induce exact quasi-inverse equivalences of categories (\ref{256Be2}),
then we reduce to the case where $n = 0$, i.e. we suppose 
$Y/S$ has locally finite $p$-bases and moreover 
the irreducible components of $Y$ are either smooth 
over $\Spec k((t))$ or over $\Spec k$. Hence we reduce to check the following two cases.

a) Suppose that $Y$ is integral and smooth over $\Spec k((t))$.
Remark that $\cV [[t]]\{ \frac{1}{t}\}$ is a complete local ring whose residue field is $k((t))$ and
whose maximal ideal is generated by a uniformizer of $\cV$ (which is also non nilpotent in $\cV [[t]]\{ \frac{1}{t}\}$). 
Hence, $\cV [[t]]\{ \frac{1}{t}\}$ is a complete  discrete valuation whose residue field is $k((t))$.
There exists a affine smooth formal $\Spf \cV [[t]]\{ \frac{1}{t}\}$-scheme
$\Y$ which is a lifting of $Y/\Spec k((t))$. 
We get $\fY \hookrightarrow \fP'$ a lifting of $Y \hookrightarrow P'$. 
This yields by composition with $f$ the lifting  $\fY \hookrightarrow \fP$ of $Y \hookrightarrow P$. 
Since $\fY/\fS$, $\fP/\fS$ and $\fP '/\fS$ have locally finite $p$-basis, then 
this is an obvious consequence of Berthelot-Kashiwara theorem 
\ref{u!u+=id}. 

b) Suppose that $Y$ is integral and smooth over $S$. We proceed similarly to 
the step 1.a).

2) Let us go back to the general case. We proceed by induction on the lexicographic order 
$(\dim Y, N _Y)$ where $\dim Y$ is dimension of $Y$ and $N _Y$ is the number of irreducible components of maximal dimension
of $Y$ (i.e. of dimension $\dim Y$). 

a) Suppose $\dim Y = 0$. Since $Y$ is noetherian, then $Y$ is artinian (see  \cite[6.2.2]{EGAI}).
Hence we can suppose $Y$ irreducible.
Since $Y$ is reduced, then $Y$ is an affine integral scheme of dimension $0$, i.e. 
$Y= \Spec L$ where $L$ is a field. When $t  =0$ in $L$
then $L/k$ is a finite extension. Since $k$ is perfect, $Y/S$ is finite and étale and we are done. 
When $t \not =0$ in $L$, we get that $L/k((t))$ is a finite extension.
Let $\Lambda ^{\flat}:= k  (( t ^{p ^{-\infty}})) $ be a perfect closure of $\Lambda := k ((t))$.
Let 
$\widetilde{Y}  ^{\flat}:= 
Y^{\flat}_{\mathrm{red}}:= (Y \times _{\Spec \Lambda} \Spec \Lambda ^{\flat}) _{\mathrm{red}}$.
Since $\widetilde{Y}  ^{\flat}$ is a reduced scheme of dimension $0$ of finite type over $\Lambda ^{\flat}$, 
then $\widetilde{Y}  ^{\flat}/\Spec (\Lambda  ^{\flat})$ is finite and étale. 
For some integer $n$, 
let $ \Lambda _{(n)}:=k  (( t ^{p ^{-n}}))$ and  
$Y _{(n)}: =Y \times _{\Spec (\Lambda)} \Spec ( \Lambda_{(n)})$.
Then, following \ref{rem-sp-desc-123}.\ref{rem-sp-desc-123-3b}, 
$(Y _{(n)} )  _{\mathrm{red}}$ is smooth over $ \Lambda _{(n)}$
for $n$ large enough. 
Hence, this is a consequence of the step 1). 

b)  Now, suppose $(\dim Y, N _Y)$ is such that $\dim Y\geq 1$
 and suppose the theorem holds
for $(\dim Y, N _Y)$ strictly lower.

Following \ref{rem-sp-desc-123-3d},
for $n$ large enough, there exists a dense open $U ' _1$ of 
$(Y _{(n)} )  _{\mathrm{red}}$ having locally finite $p$-bases and moreover such that 
the irreducible components of $U ' _1$ are either smooth over $\Spec (k[[t ^{p ^{-n}}]]$
or smooth over $\Spec k$. 
Since $f _{(n)+}$ and $f _{(n)} ^!$ induce exact quasi-inverse equivalences of categories (\ref{256Be2}),
then we reduce to the case where $n = 0$, i.e. 
we can suppose there exists a dense open $U ' _1$ of 
$Y$ having locally finite $p$-bases and moreover such that 
the irreducible components of $U ' _1$ are either  smooth over $\Spec (k((t)))$
or smooth over $\Spec k$. 
Let $V $ be one irreducible component of $U' _1$ of dimension $\dim Y$.
Shrinking $V$ if necessary, we can suppose there exists a divisor 
$D$ of $P$ such that $V = Y \setminus D$.
Set $Z := D \cap Y$. 
We get 
$(\dim Z, N _{Z})
<
(\dim Y, N _Y)$ and we can use the induction hypothesis on $Z$.

3) We check in this step that 
for any 
    $\E ^{\prime (\bullet)} \in \smash{\underrightarrow{LM}}  _{\Q,\mathrm{ovcoh}} (Y , \fP '/\cV)$, 
for any integer $r \neq 0$,
 $\mathcal{H} ^r f _{+} ^{(\bullet)} ( \E ^{\prime(\bullet)}) = 0$.
\medskip

The localisation triangle in $Z$ of  $\E ^{\prime(\bullet)}$ induces the exact sequence in 
$\smash{\underrightarrow{LM}} _{\Q,\mathrm{ovcoh}} (Y , \fP' /\cV)$:
\begin{equation}\label{s.e.loca}
0 \rightarrow \mathcal{H} ^{\dag , 0} _{Z} (\E ^{\prime(\bullet)} ) \rightarrow
\E ^{\prime(\bullet)} \rightarrow (\hdag Z ) ( \E ^{\prime(\bullet)}) \rightarrow \mathcal{H}
^{\dag , 1} _{Z} (\E ^{\prime(\bullet)} ) \rightarrow 0 .
\end{equation}
Since $Z$ locally comes from a divisor of $P'$, then the functor
$(\hdag Z )\colon 
\smash{\underrightarrow{LD}} ^{\mathrm{b}} _{\Q,\mathrm{ovcoh}} (Y , \fP '/\cV)
\to 
\smash{\underrightarrow{LD}} ^{\mathrm{b}} _{\Q,\mathrm{ovcoh}} (Y , \fP '/\cV)$
is exact (see \ref{rem-tstructure-exact}.\ref{rem-tstructure-exact2}). 
Let  $\FF ^{\prime(\bullet)}$ be the kernel of the epimorphism
$(\hdag Z ) ( \E ^{\prime(\bullet)}) 
\rightarrow  \mathcal{H} ^{\dag , 1} _{Z} (\E ^{\prime(\bullet)} )$.
We get the exact sequence in $\smash{\underrightarrow{LM}} _{\Q,\mathrm{ovcoh}} (Y , \fP' /\cV)$
$$0 \rightarrow \FF ^{\prime(\bullet)} \rightarrow (\hdag Z ) ( \E ^{\prime(\bullet)})
\rightarrow  \mathcal{H} ^{\dag , 1} _{Z} (\E ^{\prime(\bullet)} ) \rightarrow
0.$$ 
By applying the functor 
$f ^{(\bullet)} _{+}$ to this latter exact sequence, we get a long exact sequence.
We have
$(\hdag Z)(\E ^{\prime(\bullet)}) 
\in 
\smash{\underrightarrow{LM}} _{\Q,\mathrm{ovcoh}} (Y \setminus Z, \fP' /\cV)$ 
and
$\mathcal{H} ^{\dag , 1} _{Z} (\E ^{\prime(\bullet)}) \in 
\smash{\underrightarrow{LM}} _{\Q,\mathrm{ovcoh}} (Z , \fP' /\cV)$.
Hence, following the step 1), using the induction hypothesis, 
using the long exact sequence, we check that 
for any integer 
$r \not \in \{0,1\}$, 
we have 
$\mathcal{H} ^r (f _{+} ^{(\bullet)} ) ( \FF ^{\prime(\bullet)} ) = 0$. 
Moreover, 
 $\mathcal{H} ^1 (f _{+} ^{(\bullet)} ) ( \FF ^{\prime(\bullet)} ) = 0$ if and only if 
the morphism
$s\colon \mathcal{H} ^0 (f _{+} ^{(\bullet)} )
((\hdag Z ) ( \E ^{\prime(\bullet)}))  \rightarrow \mathcal{H} ^0 (f _{+} ^{(\bullet)} )(
\mathcal{H} ^{\dag , 1} _{Z} (\E ^{\prime(\bullet)}))$ 
is an epimorphism. 
We split the check of this latter property in the following two steps a) and b).

3.a) In this step, we check that
the morphism $s ':=
    \mathcal{H} ^0(\R \underline{\Gamma} ^\dag _{Y}
    \circ f ^{!(\bullet)})  (s ) $
is an epimorphism.
Since  $(\hdag Z)(\E ^{\prime(\bullet)}) \in 
\smash{\underrightarrow{LM}} _{\Q,\mathrm{ovcoh}} (Y \setminus Z, \fP' /\cV)$, 
since the functors
$ \R \underline{\Gamma} ^\dag _{Y}
    \circ f ^{!(\bullet)} $
and
$\R \underline{\Gamma} ^\dag _{Y\setminus Z}
    \circ f ^{!(\bullet)}$
    are canonically isomorphic  
    over 
    $\smash{\underrightarrow{LD}} ^{\mathrm{b}} _{\Q,\mathrm{ovcoh}} (Y \setminus Z, \fP /\cV)$
then following the step 1), 
the canonical morphism
$$(\hdag Z)(\E ^{\prime(\bullet)}) \rightarrow 
\mathcal{H} ^0( \R \underline{\Gamma} ^\dag _{Y}
    \circ f ^{!(\bullet)} )\circ \mathcal{H} ^0(f _{+} ^{(\bullet)}) ((\hdag Z)(\E ^{\prime(\bullet)}) )$$ 
    is an isomorphism.
Since 
$\mathcal{H} ^{\dag , 1} _{Z} (\E ^{\prime(\bullet)}) \in 
\smash{\underrightarrow{LM}} _{\Q,\mathrm{ovcoh}} (Z , \fP' /\cV)$,
since the functors
$ \R \underline{\Gamma} ^\dag _{Y}
    \circ f ^{!(\bullet)} $
and
$\R \underline{\Gamma} ^\dag _{Z}
    \circ f ^{!(\bullet)}$
    are canonically isomorphic  
    over 
    $\smash{\underrightarrow{LD}} ^{\mathrm{b}} _{\Q,\mathrm{ovcoh}} (Z , \fP /\cV)$
then by induction hypothesis the canonical morphism
$$\mathcal{H} ^{\dag ,1} _{Z} (\E ^{\prime(\bullet)}) 
\rightarrow
    \mathcal{H} ^0(\R \underline{\Gamma} ^\dag _{Y}
    \circ f ^{!(\bullet)} ) \circ \mathcal{H} ^0(f _{+} ^{(\bullet)} )(\mathcal{H} ^{\dag ,1} _{Z} (\E ^{\prime(\bullet)}))$$ 
    is an isomorphism.
Since $(\hdag Z)(\E ^{\prime(\bullet)}) \rightarrow \mathcal{H} ^{\dag ,1} _{Z} (\E ^{\prime(\bullet)})$ 
is an epimorphism,
this yields that so is $s '$.

\bigskip

3.b)
Let us check that  $s $ is an epimorphism.
Let  
$\FF ^{(\bullet)} 
\in
\smash{\underrightarrow{LM}} _{\Q,\mathrm{ovcoh}} (Y, \fP /\cV)$ 
be the image of  $s $, and $i $ be the canonical monomorphism 
$\FF ^{(\bullet)}  
\hookrightarrow 
\mathcal{H} ^0 (f _{+} ^{(\bullet)} )( \mathcal{H} ^{\dag , 1} _{Z} (\E ^{\prime(\bullet)}))$.
Since  
$\mathcal{H} ^0 (f _{+} ^{(\bullet)} )( \mathcal{H} ^{\dag , 1} _{Z} (\E ^{\prime(\bullet)}))$
has his support in $Z$, then $i $ is in fact a monomorphism of 
$\smash{\underrightarrow{LM}} _{\Q,\mathrm{ovcoh}} (Z, \fP /\cV)$.
Using the induction hypothesis,
since the functors
$ \R \underline{\Gamma} ^\dag _{Y}
    \circ f ^{!(\bullet)} $
and
$\R \underline{\Gamma} ^\dag _{Z}
    \circ f ^{!(\bullet)}$
    are canonically isomorphic  
    over 
    $\smash{\underrightarrow{LD}} ^{\mathrm{b}} _{\Q,\mathrm{ovcoh}} (Z , \fP /\cV)$
this yields that 
$i ' : =\mathcal{H} ^0( \R \underline{\Gamma} ^\dag _{Y}
    \circ f ^{!(\bullet)} ) (i )$
    is a monomorphism.
Since  $s '$ is an epimorphism, then so is $i '$. 
Hence, the morphism $i'$ is an isomorphism.
Using the induction hypothesis,
this implies that $i $ is an isomorphism.
This yields that $s$ is an epimorphism.
    
\medskip

3.c) Hence, we have checked that 
for any integer 
$r \not= 0$, 
we have 
$\mathcal{H} ^r (f _{+} ^{(\bullet)} ) ( \FF ^{\prime(\bullet)} ) = 0$. 
From \ref{s.e.loca}, we get the exact sequence 
    $0 \rightarrow \mathcal{H} ^{\dag , 0} _{Z} (\E ^{\prime(\bullet)} ) \rightarrow \E ^{\prime(\bullet)} \rightarrow \FF ^{\prime(\bullet)} \rightarrow 0 .$
By applying the functor  $f _{+} ^{(\bullet)}$ 
to this latter sequence, we get a long exact sequence.
Looking at this later one, we remark that the property 
``for any $r \neq 0$, $\mathcal{H}
    ^r (f _{+} ^{(\bullet)} ) ( \FF ^{\prime(\bullet)} ) = 0$ 
    and $\mathcal{H} ^r ( f _{+} ^{(\bullet)} ) (\mathcal{H} ^{\dag , 0} _{Z} (\E ^{\prime(\bullet)} ))=0$'',
implies that 
    "for any $r \neq 0$, $\mathcal{H} ^r (f _{+} ^{(\bullet)} ) (\E ^{\prime(\bullet)} ) = 0$". 
    \medskip

    4) Similarly to the step 3), we check that for any  $r\neq 0$, for any 
    $\E ^{(\bullet)} \in \smash{\underrightarrow{LM}}  _{\Q,\mathrm{ovcoh}} (Y , \fP /\cV)$, 
    we have 
    $\mathcal{H} ^r (\R \underline{\Gamma} ^{\dag} _{Y}\circ  f ^{!(\bullet)}  )  ( \E ^{(\bullet)}) = 0$. 
    \medskip

5)    
It remains to check the last statement of the theorem. 
Let $\E ^{(\bullet)}
\in \smash{\underrightarrow{LD}} ^{\mathrm{b}} _{\Q,\mathrm{ovcoh}} (Y , \fP /\cV)$.
    Using the localisation triangle with respect to $Z$, 
     to check that the morphism
     $f _{+} ^{(\bullet)}\circ \R \underline{\Gamma} ^{\dag} _{{Y} }\circ  f ^{!(\bullet)}( \E   ^{(\bullet)} ) \rightarrow  \E   ^{(\bullet)}  $
     is an isomorphism, we reduce to check we get an isomorphism after applying 
     $\R \underline{\Gamma} ^{\dag} _{{Z} }$ and $(\hdag Z)$.
     Using \ref{2.2.18} and \ref{gammayY'}, 
     after applying $\R \underline{\Gamma} ^{\dag} _{{Z} }$, we get a morphism canonically
     isomorphic to the canonical morphism 
     $f _{+} ^{(\bullet)}\circ \R \underline{\Gamma} ^{\dag} _{{Z} }\circ  f ^{!(\bullet)}
     ( \R \underline{\Gamma} ^{\dag} _{{Z} } \E   ^{(\bullet)} ) \rightarrow  
     \R \underline{\Gamma} ^{\dag} _{{Z} } \E   ^{(\bullet)}  $. By induction hypothesis, this latter is an isomorphism.
     Moreover, after applying $(\hdag Z)$, we get the morphism
     $f _{+} ^{(\bullet)}\circ \R \underline{\Gamma} ^{\dag} _{{Y  \setminus Z} }\circ  f ^{!(\bullet)}
     ( \R \underline{\Gamma} ^{\dag} _{{Y  \setminus Z} } \E   ^{(\bullet)} ) \rightarrow  
     \R \underline{\Gamma} ^{\dag} _{{Y  \setminus Z} } \E   ^{(\bullet)}  $, which is an isomorphism following the step 1).

We proceed similarly to check 
that the canonical morphism
     $\E   ^{\prime (\bullet)}  
     \rightarrow  
\R \underline{\Gamma} ^{\dag} _{{Y} }\circ  f ^{!(\bullet)} \circ f _{+} ^{(\bullet)}( \E   ^{\prime (\bullet)} )      $
     is an isomorphism for any 
     $\E ^{\prime (\bullet)}
\in \smash{\underrightarrow{LD}} ^{\mathrm{b}} _{\Q,\mathrm{ovcoh}} (Y , \fP '/\cV)$.     
\end{proof}

\begin{cor}
\label{ind-cat-overcouples}
Let $\mathfrak{C}$ be a restricted data of absolute coefficients over $\fS$
which contains $\fB _\mathrm{div}$, 
which is stable  under devissages, 
pushforwards, quasi-projective  extraordinary pullbacks,
and local cohomological functors.
Let $\mathbb{Y}:= (Y,X)$ be a couple over $\cV[[t]]$. 

\begin{enumerate}[(a)]
\item Choose a frame over $\cV[[t]]$ of the form $(Y,X,\fP)$. 
The category 
$\mathfrak{C} (Y,\fP/\cV[[t]])$
 does not depend, up to a canonical equivalence of categories, 
on the choice of the frame $(Y,X,\fP)$ over $\cV[[t]]$ enclosing $(Y,X/\cV[[t]])$.
Hence, we can simply write 
$\mathfrak{C} (\mathbb{Y}/\cV[[t]])$ 
instead of 
$\mathfrak{C} (Y,\fP/\cV[[t]])$
without ambiguity (up to equivalences of categories).

\item If moreover $\fC$ is stable under cohomology, 
then we get a canonical t-structure on 
$\mathfrak{C} (\mathbb{Y}/\cV[[t]])$. 
\end{enumerate}

\end{cor}

\begin{proof}
We can copy word by word the proof of \cite[12.2.2]{caro-6operations}.
\end{proof}

\begin{lem}
\label{ind-dual}
Let $\mathfrak{C}$ be a restricted data of absolute coefficients over $\fS$
which contains $\fB _\mathrm{div}$, 
which is stable  under devissages, 
pushforwards, quasi-projective  extraordinary pullbacks, 
local cohomological functors, and duals.
Let $\mathbb{Y}:= (Y,X)$ be a couple over $\cV[[t]]$. 
 Choose a frame of the form $(Y,X/\bbD ^{1} _{\fS})$. 
The functor 
$\R \underline{\Gamma} ^\dag _{Y} \DD _{\fP} 
\colon 
\mathfrak{C} (Y,\fP/\cV) \to \mathfrak{C} (Y,\fP/\cV)$
does not depend,
up to the canonical equivalences of categories of \ref{ind-cat-overcouples}, 
on the choice of the frame enclosing $(Y,X)$. 
Hence, 
we will denote by 
$\DD _{\mathbb{Y}}
\colon \mathfrak{C} (\mathbb{Y}/\cV) \to \mathfrak{C} (\mathbb{Y}/\cV)$
the functor 
$\R \underline{\Gamma} ^\dag _{Y} \DD _{\fP}$. 
\end{lem}

\begin{proof}
We can copy the proof of \cite[12.2.3]{caro-6operations}.
\end{proof}

\begin{lem}
\label{ind-pushforward-extinv}
Let $\mathfrak{C}$ be a restricted data of absolute coefficients over $\fS$
which contains $\fB _\mathrm{div}$, 
which is stable  under devissages, 
pushforwards, quasi-projective  extraordinary pullbacks,
and local cohomological functors.
Let  $u=(b,a)\colon (Y', X') \to (Y, X)$ be a morphism of couples over $\cV[[t]]$.
Put $\mathbb{Y}:= (Y,X)$ and $\mathbb{Y}':= (Y',X')$.
Let us choose a morphism of frames $\theta= (b,a,f) \colon (Y', X', \fP')\to (Y,X,\fP)$ over $\cV[[t]]$ enclosing $u$.

\begin{enumerate}[(a)]
\item 
The functor $\theta ^{!(\bullet )} := \R \underline{\Gamma} ^\dag _{Y'} \circ f ^{!(\bullet)}
\colon
\mathfrak{C} (Y,\fP/\cV[[t]]) \to \mathfrak{C} (Y',\fP'/\cV [[t]])$ 
does not depend on the choice of such $\theta$ enclosing $u$
(up to canonical equivalences
of categories).
Hence, it will be denoted by $u ^{!}
\colon \mathfrak{C} (\mathbb{Y}/\cV[[t]]) \to \mathfrak{C} (\mathbb{Y}'/\cV[[t]])$. 

\item Suppose that $u$ is complete, i.e. that $a \colon X' \to X$ is proper. 
The functor $\theta _{+}:= f _+ ^{(\bullet )}
\colon 
\mathfrak{C} (Y',\fP'/\cV[[t]]) \to \mathfrak{C} (Y,\fP/\cV[[t]])$
does not depend on the choice of such $\theta$ enclosing $u$
(up to canonical equivalences
of categories).
Hence, it will be denoted by 
$u _{+} \colon \mathfrak{C} (\mathbb{Y}'/\cV) \to \mathfrak{C} (\mathbb{Y}/\cV)$.
\end{enumerate}

\end{lem}

\begin{proof}
We can copy the proof of \cite[12.2.4]{caro-6operations}.
\end{proof}

\begin{lem}
\label{ind-prod-tensor}
Let $\mathfrak{C}$ be a restricted data of absolute coefficients over $\fS$
which contains $\fB _\mathrm{div}$, 
which is stable under devissages, 
pushforwards, quasi-projective extraordinary pullbacks, 
and tensor products.
Let $\mathbb{Y}:= (Y,X)$ be a couple over $\cV[[t]]$. 
 Choose a frame of the form $(Y,X,\fP/\bbD ^{1} _{\fS})$. 
The bifunctor $-\smash{\widehat{\otimes}}^\L
_{\O  _{\fP}} - [-\dim P]
\colon
\mathfrak{C} (Y,\fP/\cV[[t]]) \times \mathfrak{C} (Y,\fP/\cV[[t]]) \to \mathfrak{C} (Y,\fP/\cV[[t]])$
does not depend, up to the canonical equivalences of categories of \ref{ind-cat-overcouples}, 
on the choice of the frame enclosing $(Y,X)$.
It will be denoted by 
$\widetilde{\otimes} _{\mathbb{Y}}
\colon 
\mathfrak{C} (\mathbb{Y}/\cV[[t]]) \times \mathfrak{C} (\mathbb{Y}/\cV[[t]]) \to \mathfrak{C} (\mathbb{Y}/\cV[[t]])$. 
\end{lem}

\begin{proof}
We can copy the proof of \cite[12.2.5]{caro-6operations}.
\end{proof}

\begin{empt}
[Formalism of Grothendieck six operations]
\label{6operations}
Let $\mathfrak{C}$ be a restricted data of absolute coefficients over $\fS$
which contains $\fB _\mathrm{div}$, 
which is 
stable  under devissages, 
pushforwards, quasi-projective extraordinary pullbacks, 
duals,
and tensor products.
To sum-up the above Lemmas
we can define a formalism of Grothendieck six operations on couples 
as follows.
Let  $u=(b,a)\colon (Y', X') \to (Y, X)$ be a morphism of couples over $\cV[[t]]$.
Put $\mathbb{Y}:= (Y,X)$ and $\mathbb{Y}':= (Y',X')$.
\begin{enumerate}[(a)]

\item We have the dual functor $\DD _{\mathbb{Y}}
\colon \mathfrak{C} (\mathbb{Y}/\cV[[t]]) \to \mathfrak{C} (\mathbb{Y}/\cV[[t]])$ (see \ref{ind-dual}).

\item 
We have the extraordinary pullback $u ^{!}
\colon
\mathfrak{C} (\mathbb{Y}/\cV[[t]]) \to \mathfrak{C} (\mathbb{Y}'/\cV[[t]])$ (see \ref{ind-pushforward-extinv}).
We get the pullbacks $u ^{+}:= \DD _{\mathbb{Y}'} \circ u ^{!} \circ \DD _{\mathbb{Y}}$.

\item Suppose that $u$ is complete.
Then, we have the functor $u _+
\colon \mathfrak{C} (\mathbb{Y}'/\cV[[t]]) \to \mathfrak{C} (\mathbb{Y}/\cV[[t]])$ (see \ref{ind-pushforward-extinv}). 
We denote by  $u _{!}:= \DD _{\mathbb{Y}} \circ u _{+} \circ \DD _{\mathbb{Y}'}$, 
the extraordinary pushforward by $u$.

\item 
We have the tensor product 
$-\widetilde{\otimes} _{\mathbb{Y}}-
\colon 
\mathfrak{C} (\mathbb{Y}/\cV[[t]]) \times \mathfrak{C} (\mathbb{Y}/\cV[[t]]) \to \mathfrak{C} (\mathbb{Y}/\cV[[t]])$
(see \ref{ind-prod-tensor})
\end{enumerate}

\end{empt}

\begin{exs}
\label{nota-h-ovhol}
\begin{enumerate}[(a)]
\item  We recall
the restricted data of absolute coefficients $\smash{\underrightarrow{LD}} ^{\mathrm{b}(1)} _{\Q,\mathrm{ovhol}}$
and $\smash{\underrightarrow{LD}} ^{\mathrm{b}(1)} _{\Q,\mathrm{h}}$
are defined respectively in \ref{ex-cst-surcoh}.\ref{hstab} and \ref{ex-cst-surcoh}.\ref{ovholstab}.
Using Lemmas \ref{lem-stabextpullback} and \ref{S(D,C)stability} (and  \ref{rem-div-cst2}),
they 
are stable 
under local cohomological functors, 
pushforwards, quasi-projective extraordinary pullbacks, and duals.
Hence, 
with the notation 
\ref{ind-cat-overcouples},
using Lemmas \ref{ind-pushforward-extinv}, \ref{ind-prod-tensor}, and \ref{ind-dual},
for any frame 
$(Y,X,\fP)$ over $\cV[[t]]$,
we get the categories of the forms 
$\smash{\underrightarrow{LD}} ^{\mathrm{b}(1)} _{\Q,\mathrm{h}} (Y, \fP /\cV[[t]])$,
$\smash{\underrightarrow{LD}} ^{\mathrm{b}(1)} _{\Q,\mathrm{h}}(\mathbb{Y}/\cV[[t]])$,
$\smash{\underrightarrow{LD}} ^{\mathrm{b}(1)} _{\Q,\mathrm{ovhol}} (Y, \fP /\cV[[t]])$
or
$\smash{\underrightarrow{LD}} ^{\mathrm{b}(1)} _{\Q,\mathrm{ovhol}}(\mathbb{Y}/\cV[[t]])$
endowed with five of Grothendieck cohomological operations (the tensor product is a priori missing).

\item Following theorem \ref{dfnquprop} and the example \ref{ex-datastableevery} 
(or this is a restricted consequence of Theorem \ref{theo-V-6operations}),
there exist a  data of absolute coefficients $T$ 
which contains 
$\mathfrak{B} _{\mathrm{div}}$,
local,
stable by devissages, direct summands, 
local cohomological functors, 
pushforwards, quasi-projective extraordinary pullbacks, base change, tensor products, duals.
Hence, 
for any frame 
$(Y,X,\fP/\cV [[t]])$,
we get the triangulated category 
$T (Y, \fP /\cV [[t]])$ or $T(\mathbb{Y}/\cV [[t]])$,
endowed with 
a formalism of Grothendieck six operations.

\end{enumerate}

\end{exs}

\subsection{Formalism of Grothendieck six operations over quasi-projective schemes over $\cV [[t]]$}

\begin{dfn}
[Projective compactification]
\begin{enumerate}[(a)]
\item A frame $(Y,X,\fP)$ over $\cV [[t]]$ is said to be {\it projective} if $\fP$ is projective over $\Spf \cV [[t]]$. 
The category of projective frames over $\cV[[t]]$ is the full subcategory of the category
of  frames over $\cV [[t]]$ whose objects are projective frames over $\cV[[t]]$.

\item The category of {\it projective couples} over $\cV [[t]]$ is the full subcategory of 
the category of couples over $\cV [[t]]$ whose objects $(Y, X)$ are such that 
$X$ is projective over $\Spec k [[t]]$. 
We remark that if $(Y, X)$ is a projective couple over $\cV[[t]]$  then 
there exists a projective frame over $\cV[[t]]$ of the form $(Y,X,\fP)$.

\item Let $Y$ be a {quasi-projective scheme} 
over $\cV[[t]]$.
Then there exists 
a projective frame over $\cV[[t]]$ of the form 
$(Y,X,\fP)$. For such frame $(Y,X,\fP)$, we say that 
the projective frame $(Y,X,\fP)$ encloses $Y$ or that the projective couple $(Y,X)$ encloses $Y$. 
\end{enumerate}

\end{dfn}

\begin{empt}
[Formalism of Grothendieck six operations]
\label{6operations-variety}
Let $\mathfrak{C}$ be a restricted data of absolute coefficients over $\fS$
which contains $\fB _\mathrm{div}$, 
which is 
stable  under devissages, 
pushforwards, quasi-projective extraordinary pullbacks, 
duals,
and tensor products.
Similarly to Lemma \ref{ind-cat-overcouples}, we check using Theorem \ref{ind-CYW}
that the category 
$\mathfrak{C} (Y,\fP/\cV[[t]])$
(resp. $\mathfrak{C} (Y,X/\cV[[t]])$)
 does not depend, up to a canonical equivalence of categories, 
on the choice of the projective frame $(Y,X/\cV[[t]])$ (resp. the projective couple $(Y,X)$) over $\cV[[t]]$ enclosing $Y$.
Hence, we simply denote it by
$\mathfrak{C} (Y/\cV[[t]])$.
As for \ref{6operations}, we can define a formalism of Grothendieck six operations on quasi-projective schemes over $\cV[[t]]$ 
as follows.
Let  $u\colon Y'\to Y$ be a morphism of quasi-projective schemes over $\cV[[t]]$.
\begin{enumerate}[(a)]

\item We have the dual functor $\DD _{Y}
\colon \mathfrak{C} (Y/\cV[[t]]) \to \mathfrak{C} (Y/\cV[[t]])$ (see \ref{ind-dual}).

\item 
We have the extraordinary pullback $u ^{!}
\colon
\mathfrak{C} (Y/\cV[[t]] \to \mathfrak{C} (Y'/\cV[[t]])$ (see \ref{ind-pushforward-extinv}).
We get the pullbacks $u ^{+}:= \DD _{Y'} \circ u ^{!} \circ \DD _{Y}$.

\item 
We have the functor $u _+
\colon \mathfrak{C} (Y'/\cV[[t]]) \to \mathfrak{C} (Y/\cV[[t]])$ (see \ref{ind-pushforward-extinv}). 
We denote by  $u _{!}:= \DD _{Y} \circ u _{+} \circ \DD _{Y'}$, 
the extraordinary pushforward by $u$.

\item 
We have the tensor product 
$-\widetilde{\otimes} _{Y}-
\colon 
\mathfrak{C} (Y/\cV[[t]]) \times \mathfrak{C} (Y/\cV[[t]]) \to \mathfrak{C} (Y/\cV[[t]])$
(see \ref{ind-prod-tensor})
\end{enumerate}

\end{empt}

\subsection{Constructible t-structure}

For completeness,
we introduce the notion of constructibility.
Let $\mathfrak{C}$ be a restricted data of absolute coefficients over $\fS$
which contains $\fB _\mathrm{div}$, 
which is 
stable  under devissages, 
pushforwards, quasi-projective extraordinary pullbacks, 
duals,
tensor products, and cohomology.

\begin{empt}
[Constructible t-structure]
\label{t-structure}

Let $\mathbb{Y}:= (Y,X)$ be a couple over $\cV [[t]]$. Choose a  frame $(Y,X/\cV [[t]])$. 
If $Y' \rightarrow Y$ is an immersion, then 
we denote by 
$i _{Y'} \colon (Y', X', \fP) \to  (Y, X, \fP)$ the induced morphism
where $X'$ is the closure of $Y'$ in $X$. 
For any integer $n $, 
let $\fP _{(n)} := \fP \times _{\Spf \cV [[t ]]} \Spf \cV [[t ^{p ^{-n}}]]$,
and
$f _{(n)}
\colon 
\fP _{(n)}
\to \fP $ be the projection.
We define on $\fC(\mathbb{Y}/\cV)$
the constructible t-structure as follows.

An object $\E\in \fC(\mathbb{Y}/\cV) $
belongs to 
$\fC^{c, \geq 0}(\mathbb{Y}/\cV)$ 
(resp. $\fC^{c, \leq 0}(\mathbb{Y}/\cV)$)
if there exists 
for $n $ large enough
a stratification having locally finite $p$-bases
$(Y _{i}) _{i=1, \dots , r}$ of $( Y \times _{\Spec k[[t]] } \Spec (k[[t ^{p ^{-n}}]]) ) _{\mathrm{red}}$ 
(see Definition 
\ref{ntn-t-structureovcoh})
such that we have
$ i _{Y _i} ^{+}   (f _{(n)} ^+ (\E ^{(\bullet)}) ) [ \delta _{Y _i}]
\in 
\fC ^{\geq 0}_{\textrm{isoc}} (Y _i, \fP _{(n)})$
(resp. 
$ i _{Y _i} ^{+}   (f _{(n)} ^+ (\E ^{(\bullet)}) ) [ \delta _{Y _i}]
\in 
\fC ^{\leq 0}_{\textrm{isoc}} (Y _i, \fP _{(n)})$)
 for any $i =1,\dots, r$.

\end{empt}

\begin{prop}
Let $\mathbb{Y}:= (Y,X)$ be a couple.
\begin{enumerate}[(a)]

\item Let $\E ^{\prime(\bullet)}\to \E \to \E ^{\prime \prime(\bullet)} \to \E ^{\prime(\bullet)} [1]$ be an exact triangle in $\fC(\mathbb{Y}/\cV) $. 
If $\E ^{\prime(\bullet)}$ and $\E ^{\prime \prime(\bullet)}$are in $\fC^{c, \geq 0}(\mathbb{Y}/\cV)$
(resp. $\fC^{c, \leq 0}(\mathbb{Y}/\cV)$)
then so is $\E$.

\item 
Suppose that $Y$ has locally finite $p$-bases. 
Let $\E\in \fC _{\textrm{isoc}}(\mathbb{Y}/\cV) $.
Then $\E \in \fC^{c, \geq 0}(\mathbb{Y}/\cV)$
(resp. $\E \in \fC^{c, \leq 0}(\mathbb{Y}/\cV)$)
if and only if 
$\E \in \fC^{\geq \delta _X}_{\textrm{isoc}}(\mathbb{Y}/\cV)$
(resp. $\E \in \fC^{\leq \delta _X}_{\textrm{isoc}}(\mathbb{Y}/\cV)$). 

\end{enumerate}
\end{prop}

\begin{proof}
This is left to the reader.
\end{proof}

\small
\bibliographystyle{alpha}

\begin{thebibliography}{Gro61b}

\bibitem[Abe18]{Abe-Langlands}
Tomoyuki Abe.
\newblock Langlands correspondence for isocrystals and the existence of
  crystalline companions for curves.
\newblock {\em J. Amer. Math. Soc.}, 31(4):921--1057, 2018.

\bibitem[AC18]{Abe-Caro-weights}
Tomoyuki Abe and Daniel Caro.
\newblock Theory of weights in {$p$}-adic cohomology.
\newblock {\em Amer. J. Math.}, 140(4):879--975, 2018.

\bibitem[BBD82]{BBD}
A.~A. Be{\u\i}linson, J.~Bernstein, and P.~Deligne.
\newblock Faisceaux pervers.
\newblock In {\em Analysis and topology on singular spaces, {I} ({L}uminy,
  1981)}, volume 100 of {\em Ast{\'e}risque}, pages 5--171. Soc. Math. France,
  Paris, 1982.

\bibitem[Ber90]{Be0}
Pierre Berthelot.
\newblock {Cohomologie rigide et th\'eorie des $\mathcal{D}$-modules}.
\newblock In {\em $p$-adic analysis (Trento, 1989)}, pages 80--124. Springer,
  Berlin, 1990.

\bibitem[Ber96a]{Becohdiff}
Pierre Berthelot.
\newblock Coh\'erence diff\'erentielle des alg\`ebres de fonctions
  surconvergentes.
\newblock {\em C. R. Acad. Sci. Paris S\'er. I Math.}, 323(1):35--40, 1996.

\bibitem[Ber96b]{Be1}
Pierre Berthelot.
\newblock ${\mathcal{d}}$-modules arithm\'etiques. {I}. {O}p\'erateurs
  diff\'erentiels de niveau fini.
\newblock {\em Ann. Sci. \'Ecole Norm. Sup. (4)}, 29(2):185--272, 1996.

\bibitem[Ber00]{Be2}
Pierre Berthelot.
\newblock {$\mathcal{D}$}-modules arithm\'etiques. {I}{I}. {D}escente par
  {F}robenius.
\newblock {\em M\'em. Soc. Math. Fr. (N.S.)}, (81):vi+136, 2000.

\bibitem[Ber02]{Beintro2}
Pierre Berthelot.
\newblock {Introduction \`a la th\'eorie arithm\'etique des
  {$\mathcal{D}$}-modules}.
\newblock {\em Ast\'erisque}, (279):1--80, 2002.
\newblock Cohomologies {$p$}-adiques et applications arithm\'etiques, {II}.

\bibitem[Bou61]{bourbaki}
N.~Bourbaki.
\newblock {\em {\'{E}l\'ements de math\'ematique. {F}ascicule {X}{X}{V}{I}{I}.
  {A}lg\`ebre commutative. {C}hapitre 1: {M}odules plats.{C}hapitre 2:
  {L}ocalisation}}.
\newblock Herman, Paris, 1961.

\bibitem[Bou06]{Bourbaki-AC89}
N.~Bourbaki.
\newblock {\em \'{E}l\'{e}ments de math\'{e}matique. {A}lg\`ebre commutative.
  {C}hapitres 8 et 9}.
\newblock Springer, Berlin, 2006.
\newblock Reprint of the 1983 original.

\bibitem[Car04]{caro_surcoherent}
Daniel Caro.
\newblock {$\mathcal{D}$}-modules arithm{\'e}tiques surcoh{\'e}rents.
  {A}pplication aux fonctions {L}.
\newblock {\em Ann. Inst. Fourier, Grenoble}, 54(6):1943--1996, 2004.

\bibitem[Car05]{caro_comparaison}
Daniel Caro.
\newblock Comparaison des foncteurs duaux des isocristaux surconvergents.
\newblock {\em Rend. Sem. Mat. Univ. Padova}, 114:131--211, 2005.

\bibitem[Car06]{caro_courbe-nouveau}
Daniel Caro.
\newblock Fonctions {L} associ{\'e}es aux {$\mathcal{D}$}-modules
  arithm{\'e}tiques. {C}as des courbes.
\newblock {\em Compositio Mathematica}, 142(01):169--206, 2006.

\bibitem[Car09a]{caro-construction}
Daniel Caro.
\newblock {Arithmetic $\mathcal D$-modules associated with overconvergent
  isocrystals. Smooth case. ($\mathcal D$-modules arithm{\'e}tiques
  associ{\'e}s aux isocristaux surconvergents. Cas lisse.)}.
\newblock {\em Bull. Soc. Math. Fr.}, 137(4):453--543, 2009.

\bibitem[Car09b]{caro_log-iso-hol}
Daniel Caro.
\newblock {Overconvergent log-isocrystals and holonomy. (Log-isocristaux
  surconvergents et holonomie.)}.
\newblock {\em Compos. Math.}, 145(6):1465--1503, 2009.

\bibitem[Car09c]{caro-frobdualrel}
Daniel Caro.
\newblock Sur la compatibilit\'{e} \`a {F}robenius de l'isomorphisme de
  dualit\'{e} relative.
\newblock {\em Rend. Semin. Mat. Univ. Padova}, 122:235--286, 2009.

\bibitem[Car11a]{caro-holo-sansFrob}
Daniel Caro.
\newblock Holonomie sans structure de {F}robenius et crit\`eres d'holonomie.
\newblock {\em Ann. Inst. Fourier (Grenoble)}, 61(4):1437--1454 (2012), 2011.

\bibitem[Car11b]{caro-pleine-fidelite}
Daniel Caro.
\newblock {Pleine fid{\'e}lit{\'e} sans structure de Frobenius et isocristaux
  partiellement surconvergents}.
\newblock {\em Math. Ann.}, 349:747--805, 2011.

\bibitem[Car11c]{caro-stab-holo}
Daniel Caro.
\newblock Stabilit\'e de l'holonomie sur les vari\'et\'es quasi-projectives.
\newblock {\em Compos. Math.}, 147(6):1772--1792, 2011.

\bibitem[Car15]{caro-stab-prod-tens}
Daniel Caro.
\newblock Sur la stabilit\'e par produit tensoriel de complexes de
  {$\mathcal{D}$}-modules arithm\'etiques.
\newblock {\em Manuscripta Math.}, 147(1-2):1--41, 2015.

\bibitem[Car16a]{surcoh-hol}
Daniel Caro.
\newblock La surcoh\'{e}rence entra\^{\i}ne l'holonomie.
\newblock {\em Bull. Soc. Math. France}, 144(3):429--475, 2016.

\bibitem[Car16b]{caro-stab-sys-ind-surcoh}
Daniel Caro.
\newblock Syst\`emes inductifs coh\'{e}rents de {$\mathcal D$}-modules
  arithm\'{e}tiques logarithmiques, stabilit\'{e} par op\'{e}rations
  cohomologiques.
\newblock {\em Doc. Math.}, 21:1515--1606, 2016.

\bibitem[Car18]{caro-unip}
Daniel Caro.
\newblock Unipotent monodromy and arithmetic {$\mathcal {D}$}-modules.
\newblock {\em Manuscripta Math.}, 156(1-2):81--115, 2018.

\bibitem[Car19]{caro-6operations}
Daniel Caro.
\newblock Arithmetic $\mathcal{D}$-modules over algebraic varieties of
  characteristic $p>0$.
\newblock {\em ArXiv Mathematics e-prints}, 2019.

\bibitem[CT12]{caro-Tsuzuki}
Daniel Caro and Nobuo Tsuzuki.
\newblock Overholonomicity of overconvergent {$F$}-isocrystals over smooth
  varieties.
\newblock {\em Ann. of Math. (2)}, 176(2):747--813, 2012.

\bibitem[CV17]{Caro-Vauclair}
Daniel Caro and David Vauclair.
\newblock Logarithmic {$p$}-bases and arithmetical differential modules.
\newblock 2017.

\bibitem[dJ96]{dejong}
A.~J. de~Jong.
\newblock Smoothness, semi-stability and alterations.
\newblock {\em Inst. Hautes \'Etudes Sci. Publ. Math.}, (83):51--93, 1996.

\bibitem[FK18]{FujiwaraKatoBookI}
Kazuhiro Fujiwara and Fumiharu Kato.
\newblock {\em Foundations of rigid geometry. {I}}.
\newblock EMS Monographs in Mathematics. European Mathematical Society (EMS),
  Z\"{u}rich, 2018.

\bibitem[Gro57]{Tohoku}
Alexander Grothendieck.
\newblock Sur quelques points d'alg\`ebre homologique.
\newblock {\em T\^{o}hoku Math. J. (2)}, 9:119--221, 1957.

\bibitem[Gro60]{EGAI}
A.~Grothendieck.
\newblock \'{E}l\'ements de g\'eom\'etrie alg\'ebrique. {I}. {L}e langage des
  sch\'emas.
\newblock {\em Inst. Hautes \'Etudes Sci. Publ. Math.}, (4):228, 1960.

\bibitem[Gro61a]{EGAII}
A.~Grothendieck.
\newblock \'{E}l\'ements de g\'eom\'etrie alg\'ebrique. {II}. \'{E}tude globale
  \'el\'ementaire de quelques classes de morphismes.
\newblock {\em Inst. Hautes \'Etudes Sci. Publ. Math.}, (8):222, 1961.

\bibitem[Gro61b]{EGAIII1}
A.~Grothendieck.
\newblock \'{E}l\'ements de g\'eom\'etrie alg\'ebrique. {III}. \'{E}tude
  cohomologique des faisceaux coh\'erents. {I}.
\newblock {\em Inst. Hautes \'Etudes Sci. Publ. Math.}, (11):167, 1961.

\bibitem[Gro64]{EGAIV1}
A.~Grothendieck.
\newblock \'{E}l\'ements de g\'eom\'etrie alg\'ebrique. {IV}. \'{E}tude locale
  des sch\'emas et des morphismes de sch\'emas. {I}.
\newblock {\em Inst. Hautes \'Etudes Sci. Publ. Math.}, (20):259, 1964.

\bibitem[Gro65]{EGAIV2}
A.~Grothendieck.
\newblock \'{E}l\'ements de g\'eom\'etrie alg\'ebrique. {IV}. \'{E}tude locale
  des sch\'emas et des morphismes de sch\'emas. {II}.
\newblock {\em Inst. Hautes \'Etudes Sci. Publ. Math.}, (24):231, 1965.

\bibitem[Gro66]{EGAIV3}
A.~Grothendieck.
\newblock \'{E}l\'ements de g\'eom\'etrie alg\'ebrique. {IV}. \'{E}tude locale
  des sch\'emas et des morphismes de sch\'emas. {III}.
\newblock {\em Inst. Hautes \'Etudes Sci. Publ. Math.}, (28):255, 1966.

\bibitem[Gro67]{EGAIV4}
A.~Grothendieck.
\newblock \'{E}l\'ements de g\'eom\'etrie alg\'ebrique. {IV}. \'{E}tude locale
  des sch\'emas et des morphismes de sch\'emas {IV}.
\newblock {\em Inst. Hautes \'Etudes Sci. Publ. Math.}, (32):361, 1967.

\bibitem[Har66]{HaRD}
Robin Hartshorne.
\newblock {\em Residues and duality}.
\newblock Springer-Verlag, Berlin, 1966.

\bibitem[Hub94]{Huber-gen-rig-an-var}
R.~Huber.
\newblock A generalization of formal schemes and rigid analytic varieties.
\newblock {\em Math. Z.}, 217(4):513--551, 1994.

\bibitem[Kas95]{kashiwarathesis}
Masaki Kashiwara.
\newblock Algebraic study of systems of partial differential equations.
\newblock {\em M\'em. Soc. Math. France (N.S.)}, (63):xiv+72, 1995.

\bibitem[Kat91]{Kato-explicity-recip91}
Kazuya Kato.
\newblock The explicit reciprocity law and the cohomology of
  {F}ontaine-{M}essing.
\newblock {\em Bull. Soc. Math. France}, 119(4):397--441, 1991.

\bibitem[Ked05]{Kedlaya-coveraffinebis}
Kiran~S. Kedlaya.
\newblock More \'etale covers of affine spaces in positive characteristic.
\newblock {\em J. Algebraic Geom.}, 14(1):187--192, 2005.

\bibitem[KS06]{Kashiwara-schapira-book}
Masaki Kashiwara and Pierre Schapira.
\newblock {\em Categories and sheaves}, volume 332 of {\em Grundlehren der
  Mathematischen Wissenschaften [Fundamental Principles of Mathematical
  Sciences]}.
\newblock Springer-Verlag, Berlin, 2006.

\bibitem[Liu02]{Liu-livre-02}
Qing Liu.
\newblock {\em Algebraic geometry and arithmetic curves}, volume~6 of {\em
  Oxford Graduate Texts in Mathematics}.
\newblock Oxford University Press, Oxford, 2002.
\newblock Translated from the French by Reinie Ern{{\'e}}, Oxford Science
  Publications.

\bibitem[LP16]{Lazda-Pal-Book}
Christopher {Lazda} and Ambrus {P\'al}.
\newblock {\em {Rigid cohomology over Laurent series fields.}}
\newblock Cham: Springer, 2016.

\bibitem[Mon02]{these_montagnon}
Claude Montagnon.
\newblock {\em {G{\'e}n{\'e}ralisation de la th{\'e}orie arithm{\'e}tique des
  $\mathcal{D}$-modules {\`a} la g{\'e}om{\'e}trie logarithmique}}.
\newblock PhD thesis, Universit{\'e} de {R}ennes {I}, 2002.

\bibitem[MW68]{MonskyWashnitzer}
P.~Monsky and G.~Washnitzer.
\newblock Formal cohomology. {I}.
\newblock {\em Ann. of Math. (2)}, 88:181--217, 1968.

\bibitem[SGA6]{sga6}
{\em Th\'eorie des intersections et th\'eor\`eme de {R}iemann-{R}och}.
\newblock Springer-Verlag, Berlin, 1971.
\newblock S\'eminaire de G\'eom\'etrie Alg\'egrique du Bois-Marie 1966--1967
  (SGA 6), Dirig\'e par P. Berthelot, A. Grothendieck et L. Illusie. Avec la
  collaboration de D. Ferrand, J. P. Jouanolou, O. Jussila, S. Kleiman, M.
  Raynaud et J. P. Serre, Lecture Notes in Mathematics, Vol. 225.

\bibitem[SGA4]{sga4-2}
{\em Th\'eorie des topos et cohomologie \'etale des sch\'emas. {T}ome 2}.
\newblock Springer-Verlag, Berlin, 1972.
\newblock S\'eminaire de G\'eom\'etrie Alg\'ebrique du Bois-Marie 1963--1964
  (SGA 4), Dirig\'e par M. Artin, A. Grothendieck et J. L. Verdier. Avec la
  collaboration de N. Bourbaki, P. Deligne et B. Saint-Donat, Lecture Notes in
  Mathematics, Vol. 270.

\bibitem[Vir04]{Vir04}
Anne Virrion.
\newblock Trace et dualit\'e relative pour les {$\mathcal{D}$}-modules
  arithm\'etiques.
\newblock In {\em Geometric aspects of Dwork theory. Vol. I, II}, pages
  1039--1112. Walter de Gruyter GmbH \& Co. KG, Berlin, 2004.

\end{thebibliography}
\def\cprime{$'$}

\bigskip
\noindent Daniel Caro\\
Laboratoire de Mathématiques Nicolas Oresme\\
Université de Caen
Campus 2\\
14032 Caen Cedex\\
France.\\
email: daniel.caro@unicaen.fr

\end{document}